\documentclass{amsart}
\usepackage[utf8]{inputenc}
\usepackage{imakeidx}
\usepackage{makeidx}
\makeindex[columns=2, options= -s example_style.ist]
\usepackage{eurosym}
\usepackage[utf8]{inputenc}
\usepackage{amscd}
\usepackage{amsthm}
\usepackage{extarrows}
\usepackage{latexsym}
\usepackage{mathrsfs}
\usepackage{pst-all,multido,ifthen}
\usepackage{amssymb,amsmath,amsfonts}
\usepackage{tikz}
\usepackage{tikz-cd}
\usepackage[all]{xy}
\usepackage{xcolor}
\usepackage[unicode=true]{hyperref}
\hypersetup{
    colorlinks=true,      
    citecolor=blue!70!black, 
    linkcolor=red!70!black,  
    urlcolor=blue!90!black   
}
\usepackage{stmaryrd}

\newtheorem{theorem}{Theorem}[section]%
\newtheorem{corollary}[theorem]{Corollary}%
\newtheorem{definition}[theorem]{Definition}%
\newtheorem{example}[theorem]{Example}%
\newtheorem{lemma}[theorem]{Lemma}%
\newtheorem{remark}[theorem]{Remark}%
\newtheorem{notation}[theorem]{Notation}%
\newtheorem{proposition}[theorem]{Proposition}
\newtheorem{fact}[theorem]{Fact}%
\newtheorem{property}[theorem]{Property}%

\def\Aut{\textup{Aut}}
\def\Alt{\textup{Alt}}

\def\a{\textup{a}}
\def\triv{\textup{triv}}

\def\d=1{\textup{coll}}

\def\e{\textup{e}}
\def\E{\textup{E}}

\def\shift{\textup{sh}}

\def\Perm{\textup{TPerm}}
\def\perm{\textup{Perm}}

\def\B{\textup{B}}
\def\r{\textup{r}}
\def\sup{\textup{sup}}
\def\Span{\textup{Span}}

\def\symm{\textup{s}}
\def\nsymm{\textup{ns}}
\def\C{\textup{C}}

\def\char{\textup{char}}

\def\dir{\textup{dir}}

\def\End{\textup{End}}
\def\Res{\textup{Res}}
\def\R{\textup{R}}

\def\NAT{\textup{NAT}}
\def\RNAT{\textup{RNAT}}
\def\TIS{\textup{TIS}}

\def\supp{\textup{supp}}

\def\S{\textup{S}}
\def\GA{\textup{GA}}
\def\GGA{\textup{GGA}}
\def\SGA{\textup{SGA}}
\def\Stab{\textup{Stab}}
\def\Fix{\textup{Fix}}

\def\TGA{\textup{TGA}}
\def\STGA{\textup{STGA}}
\def\Gal{\textup{Gal}}
\def\GL{\textup{GL}}
\def\AGL{\textup{AGL}}

\def\c{\textup{c}}
\def\d{\textup{d}}
\def\l{\textup{l}}
\def\s{\textup{s}}
\def\z{\textup{z}}

\def\L{\textup{L}}

\def\Im{\textup{Im}}

\def\aff{\textup{aff}}

\def\j{\textup{j}}
\def\Ker{\textup{Ker}}
\def\L{\textup{L}}

\def\m{\textup{m}}

\def\O{\textup{O}}

\def\rank{\textup{rank}}

\def\Reg{\textup{Reg}}

\def\i{\textup{i}}
\def\l{\textup{l}}
\def\q{\textup{q}}
\def\r{\textup{r}}
\def\n{\textup{n}}
\def\N{\textup{N}}

\def\w{\textup{w}}
\def\sc{\textup{sc}}
\def\der{\textup{der}}
\def\ASL{\textup{ASL}}
\def\SL{\textup{SL}}

\def\Spec{\textup{Spec\,}}

\def\t{\textup{t}}
\def\v{\textup{v}}

\def\cycle{\textup{simple}}
\def\s{\textup{s}}

\def\t{\textup{t}}

\def\mq{\textup{mq}}

\begin{document}
\title[transitivity of tame groups of automorphisms of affine spaces]{Infinite transitivity of tame groups of automorphisms of affine spaces}
\author{Alexander Borisov, Ofer Gabber, and Adrian Vasiu}
\maketitle
\centerline{\today}
\vskip 0.7cm 

\noindent
{\bf Abstract.} For positive integers $n$ and $m$, we study the actions of the groups of tame automorphisms of the $n$-dimensional affine spaces over finite fields on ordered subsets of $m$ points. Our primary interest lies in constructing tame automorphisms that take one ordered sequence to another and in proving upper and lower bounds on the maximal complexity of such automorphisms. Our preferred measure of complexity of an automorphism is the maximum of the degrees of the polynomials that define it and its inverse. Using methods and results from various branches of mathematics, including the theory of symmetric groups, affine geometry over finite fields, polynomial interpolation, combinatorics of projective spaces over fields, and polynomial automorphisms, we obtain a wide variety of qualitative and quantitative results.



\newpage\tableofcontents

\newpage\section*{Preface}\label{S0}

Let $K$ be a field\index{field!finite field} and $n$ a positive integer. The group of automorphisms of the affine space $\mathbb A^n_K$ over $\Spec K$ acts naturally on $K^n.$ The action is given by affine transformations for $n=1$ and by tame transformations for $n=2$ by a result of van der Kulk (\cite{vdK}). For $n\geq 3,$ tame automorphisms form an important subgroup of the group of all automorphisms of $\mathbb A^n_K$. Our goal is to study the action of this subgroup on $K^n$ for $n\ge 2$. We are especially interested in the action on ordered finite subsets of $K^n.$ For $K=\mathbb C$ (and other infinite fields) it was studied previously by Jelonek (\cite{Je}), Kaliman (\cite{Kali}), Srinivas (\cite{Sr}), and Kaliman--Zaidenberg (\cite{KZ}). So, to begin with, we first assume that $K$ is finite and denote by $|K|$ its number of elements and by $\char(K)$ its characteristic.

At the qualitative level, we show that the image of the group of tame automorphisms of $\mathbb A^n_K$ in $\perm(K^n)$ is $\perm(K^n)$ if $4\nmid |K|$ and is $\Alt(K^n)$ if $4\mid |K|$. Thus for each integer $m$ between $1$ and $|K|^n$ if $4\nmid |K| $ or between $1$ and $|K|^n-2$ if $4\mid |K|$, for every two sequences $P_1,\ldots,P_m$ and $Q_1,\ldots,Q_m$ of distinct elements of $K^n$ there exists a tame automorphism $a$ of $\mathbb A^n_K$ such that $a(P_i)=Q_i$ for each $i\in\{1,\ldots,m\}.$ It is natural to ask for a quantitative version of this result, by looking for the worst minimal complexity of such an automorphism $a$. While various definitions of ``complexity'' are possible, we focus on the maximum of the total degrees of $a$ and $a^{-1}.$ So we define $\pi_{n,m}(K)$ to be the smallest number, such that, regardless of what $P_1,\ldots,P_m$ and $Q_1,\ldots,Q_m$ are, we can choose $a$ with the property that all polynomials defining it and its inverse have total degrees at most $\pi_{n,m}(K)$.

We use many different techniques to get upper and lower bounds for $\pi_{n,m}(K)$ for $n\geq 2$ and various $m$ and $K$. In particular, we use new and old invariants and notions to bring together methods and results from the theory of symmetric groups, affine geometry over finite fields, polynomial interpolation, combinatorics of projective spaces over fields, and polynomial automorphisms. For instance, fourteen invariants of finite subsets of $K^n$ are introduced, studied, and applied, implicitly providing refined Lagrange interpolations in multiple indeterminates. Also, we use extensively product decompositions of permutations in $\perm(K^n)$ into $i$ disjoint $s$-cycles, where $s$ is mainly $2$ if $|K|=2$ and is an odd integer in the interval $[3,2\char(K)-1]$ if $|K|\ge 3$ and where the natural number $i$ varies in an interval $[0,r]$ with $r$ well chosen in terms of $m$, $n$, $K$, and $s$. 

Some of our results, especially for small $m,$ are optimal or close to optimal. For instance, if $|K|\ge\frac{m^2-m}{4}$ and $s:=\min(m,|K|)-1$, then $\pi_{n,m}(K)\le s^2.$ In contrast, for large $m$s, our tightest upper bounds are proven only for $m$ greater than some explicit large constant or are asymptotic in $m.$ 

The lower bounds when $|K|>m$ rely on works of Abhyankar--Moh, Richman, and Kanga (see \cite{AM}, \cite{Ri}, and \cite{Kan}); if $m$ is $3r$ (or $3r+1$ or $3r+2$), then for generic collinear sequences $P_1,\ldots,P_m$ and $Q_1,\ldots,Q_m$ of distinct elements of $K^2$, each such automorphism $a$ is defined by two polynomials among which one has degree at least $4r^2$ (or $4r^2+4r$ or $4r^2+6r+2$), and hence degree at least $\frac{4m^2}{9}$.
 
For $n=2$ our lower and upper bounds are quite robust for large $m$ as well. For example, as $|K|$ goes to infinity along prime powers, we have
$$2\leq \liminf \frac{\ln \bigl( \pi_{2 , |K|^2-2}(K)\bigr)}{|K|\ln |K|} \leq \limsup  \frac{\ln\bigl( \pi_{2 , |K|^2-2}(K)\bigr)}{|K|\ln |K|} \le 18$$
and the same holds with $|K|^2-2$ replaced by $|K|^2$ if we restrict to only $|K|$s with $4\nmid |K|$. Moreover, if we restrict to only $|K|$s that are primes or satisfy $\char(K)=3$ (or satisfy $\char(K)=p$ with $p\ge 5$ a fixed prime), then $18$ can be replaced by $13.5$ (or by $\frac{27p}{2(p-1)}$). The sequence $(\frac{27i}{2(i-1)})_{i\ge 5}$ has initial value $\frac{135}{8}=16.875$, is strictly decreasing, and converges to $13.5$.

In some cases we are able to get the precise values of $\pi_{n,m}(K)$. 

We also study several variations of $\pi_{n,m}(K)$. For instance, $\upsilon_{n,m}(K)\in\mathbb N^{\ast}$ is the smallest such that for every two subsets $Y$ and $Z$ of $K^n$ with $m$ elements, there exists a tame automorphism $a$ such that $a(Y)=Z$ and all polynomials defining it and its inverse have total degrees at most $\upsilon_{n,m}(K)$. As a key difference between the $\pi_{n,m}(K)$s and the $\upsilon_{n,m}(K)$s we have $\upsilon_{n,m}(K)=\upsilon_{n,|K|^n-m}(K)$ and we do not know if the sequence $\bigl(\upsilon_{n,m}(K)\bigr)_{m=1}^{\lfloor\frac{|K|^n}{2}\rfloor}$ is non-decreasing. For $n=2$, as $|K|$ goes to infinity along prime powers, we have
$$2\leq \liminf \frac{\ln \bigl( \upsilon_{2,\lfloor\frac{|K|^2}{2}\rfloor}(K)\bigr)}{|K|\ln |K|} \leq \limsup \frac{\ln\bigl(\upsilon_{2,\lfloor\frac{|K|^2}{2}\rfloor}(K)\bigr)}{|K|\ln |K|} \le 14$$
and $14$ can be replaced by $9$ (or by $13.5$ or $\frac{27p}{2(p-1)}$) if we restrict to only $|K|$s that are primes or even (or satisfy $\char(K)=3$ or $\char(K)=p$ with $p\ge 29$ a fixed prime).

If $K$ is an infinite field, then $\pi_{n,m}(K)$ is defined similarly and we show that $\pi_{n,m}(K)\le (m-1)^2$ and that the equality holds if $n=2$ and either $m\le 4$ or $m\in\{5,6\}$ and $K$ has a primitive $m$-th root of unity. The inequality $\pi_{2,m}(K)\ge\frac{4m^2}{9}$ we have for $m\ge 5$ (see above) is improved to $\pi_{2,m}(K)\ge \lfloor\frac{m+2}{2}\rfloor(m-1)$ if $m\ge 7$ and $K$ has a primitive $m$-th root of unity.

We propose a number of problems (resp.\ questions) that we found important or tried to solve but obtained only partial results (resp.\ we did not have sufficient numerical data to state as conjectures). 

While non-tame (i.e., wild) automorphisms are not referred to explicitly in what follows except in the section on Nagata's automorphism introduced in \cite{N}, many of our problems pertain indirectly to them and we hope that such problems and several of our results will be used in the future in connection to them.

We refer to Section \ref{S1} for precise definitions, notation, and a detailed exposition of our main results. Additional notation is introduced throughout the text, with the most commonly used math symbols listed and briefly described after the bibliography.

\medskip\noindent
{\bf Acknowledgement.} The first author would like to thank Inna Sysoeva for numerous helpful
conversations regarding this work. The third author would like to thank SUNY Binghamton and IH\'ES for good working conditions and Inna Sysoeva for sharing the result of a Fortran code mentioned in Footnote \ref{foot13}. 

\medskip\noindent
{\bf Interests and AI statement.} Competing interests: authors declare none. No AI was used in this work, and no results are computer-based or inspired by computer-generated data. Standard calculators were used to compute constants such as $\binom{27}{10}$ and values of the natural logarithm and exponential functions. Footnote \ref{foot13} mentions that computer programming can obtain a slightly better result in one particular case.

\newpage\section{Introduction}\label{S1}

\phantomsection{For a set $A$, let $\perm(A)$ be the group of permutations of it, and let $|A|\in\mathbb N\cup\{\infty\}$ be its cardinality if it is finite and be $\infty$ if it is infinite. If $A$ is finite, let $\Alt(A)$ be the subgroup of $\perm(K^n)$ formed by even permutations of $A$. For $(r,s)\in\mathbb N^2$, let}\label{PH1} 
$$\llbracket r,s\rrbracket:=\{i\in\mathbb N|r\le i\le s\}$$ 
(so $\llbracket r,s\rrbracket=\emptyset$ if $r>s$). 

\phantomsection{Let $K$ be a field.\index{field} Let $K^{\ast}$ be the multiplicative group of units of $K$. Let $\mathbb A^r_K$ be the affine space of dimension $r$ over $\Spec K$.}\label{PH2} 

\phantomsection{For a scheme $X$ of finite type over $\Spec K$ let $\Reg(X)$ be its regular locus. For $m\in\mathbb N^{\ast}$ with $m\le|X(K)|$, we consider the set
$$\mathbb D_m(X):=\{(P_1,\ldots,P_m)\in X(K)^m|P_i\neq P_j\;\forall (i,j)\in \llbracket1,m\rrbracket^2, i<j\}$$
of $m$-tuples formed by distinct elements of $X(K)$; its elements are denoted by $\underline{\star}=(\star_1,\ldots,\star_m)$ with $\star\in\{P,Q,O\}$.}\label{PH3} 

\phantomsection{Let $\GA(X)$ be the {\it group of automorphisms}\index{automorphism} of $X$; the left action 
$$\mathbb T_m(X): \GA(X)\times\mathbb D_m(X)\rightarrow \mathbb D_m(X)$$
is defined by the rule $(a,\underline{P})\mapsto a(\underline{P}):=\bigl(a(P_1),\ldots,a(P_m)\bigr)$.}\label{PH3+} 

\phantomsection{Recall that $\mathbb G_{a,K}$ and $\mathbb G_{m,K}$ are the affine smooth connected group schemes over $\Spec K$ of dimension $1$ such that for each $K$-algebra $S$, $\mathbb G_{\a,K}(S)$ is the additive group of $S$ and  $\mathbb G_{\m,K}(S)$ is the multiplicative group of units of $S$.}\label{PH3a}

\phantomsection{Identifying $\mathbb G_{\m,K}(K)=K^{\ast}$ and $\mathbb A^n_K(K)=K^n$ for $n\in\mathbb N^{\ast}$, let}\label{PH4} 
$$\mathbb D_{1,m}(K^{\ast}):=\mathbb D_m(\mathbb G_{\m,K})\subset (K^{\ast})^m,$$
$$\mathbb D_{n,m}(K):=\mathbb D_m(\mathbb A^n_K)\subset (K^n)^m,$$
$$\GA_n(K):=\GA(\mathbb A^n_K),$$
$$\mathbb T_{n,m}(K):=\mathbb T_m(\mathbb A^n_K): \GA_n(K)\times\mathbb D_{n,m}(K)\rightarrow \mathbb D_{n,m}(K).$$ 

\phantomsection{Let $\SGA(X)$ be the normal subgroup of $\GA(X)$ formed by {\it special automorphisms}\index{automorphism!special automorphism} of $X$, i.e., generated by $\mathbb G_{\a,K}(K)=K$ subgroups of $\GA(X)$ defined by faithful actions $\mathbb G_{\a,K}\times_{\Spec K} X\rightarrow X$. Let $\GGA(X)$ be the normal subgroup of $\GA(X)$ generated by $\SGA(X)$ and all $\mathbb G_{\m,K}(K)=K^{\ast}$ subgroups defined by faithful actions $\mathbb G_{\m,K}\times_{\Spec K} X\rightarrow X$. Let}\label{PH4a} 
$$\SGA_n(K):=\SGA(\mathbb A^n_K)$$
and
$$\GGA_n(K):=\GGA(\mathbb A^n_K).$$ 

\begin{definition}\label{D1}
{\bf (1)} If $X$ is affine and integral, then we call it a flexible (resp.\ an almost flexible) variety \index{flexible variety} \index{flexible variety!almost flexible variety} over $\Spec K$ if for each $m\in\mathbb N^{\ast}$ with $m\le|\Reg(X)(K)|$, $\mathbb T_m(X)$ restricts to a transitive action $\SGA(X)\times\mathbb D_m\bigl(\Reg(X)\bigr)\rightarrow \mathbb D_m\bigl(\Reg(X)\bigr)$ (resp.\ $\GGA(X)\times\mathbb D_m\bigl(\Reg(X)\bigr)\rightarrow \mathbb D_m\bigl(\Reg(X)\bigr)$). 

\smallskip
{\bf (2)} We say $K$ is infinitely transitive\index{field!infinitely transitive field} if $\mathbb A^n_K$ is an almost flexible variety over $\Spec K$ for each integer $n\ge 2$.
\end{definition}

The action $\mathbb T_{1,m}(K)$ is transitive iff $m\in\{1,2\}$ and this explains why in Definition \ref{D1}(2) we have $n\ge 2$. 

In this monograph we study restrictions of $\mathbb T_{n,m}(K)$ to different subgroups and subsets of $\GA_n(K)$ in order to classify all infinitely transitive fields in qualitative and quantitive ways.

The case when $K$ is algebraically closed has a long history, with the first major result going back at least to \cite{Je}, Thm.\ 1.2 whose proof implies that $\mathbb C$ is infinitely transitive. Similarly, the proof of \cite{Kali}, Thm.\ 1 implies that all algebraically closed fields of characteristic $0$ are infinitely transitive and \cite{Sr}, Thm.\ 1.2 implies that all infinite fields are infinitely transitive. Note that \cite{KZ}, Lem.\ 5.5 proves for a fixed $m$ using an induction on $n\ge 2$ that $\SGA_n(\mathbb C)$ acts transitively on $\mathbb D_{n,m}(\mathbb C)$, and hence that $\mathbb C$ is infinitely transitive. More recently, if $K$ is an algebraically closed field of characteristic $0$, \cite{AFKKZ}, Thm.\ 0.1 reobtains Kaliman's result by classifying all (affine integral) flexible varieties over $\Spec K$.

To describe our results, we introduce extra notation and recall basic notions. 

\phantomsection{Let $\End_n(K)$ be the multiplicative monoid of endomorphisms of $\mathbb A^n_K$ under composition; so $\GA_n(K)$ is its subgroup of invertible elements. The polynomial $K$-algebra $R:=K[x_1,\ldots,x_n]$ in indeterminates $x_1,\ldots,x_n$ is $\mathbb N$-graded by the {\it degree} function 
$$\deg:R\rightarrow\mathbb N\cup\{-\infty\};$$ so $\deg(0)=-\infty$, $\deg(1)=0$, and $\deg(x_i)=1$ for each $i\in \llbracket1,n\rrbracket$.}\label{PH5} 

\phantomsection{An $n$-tuple $(f_1,\ldots,f_n)\in R^n$ defines an endomorphism}\label{PH6} 
$$e:=\e(f_1,\ldots,f_n)\in\End_n(K)$$ 
as follows: it is the spectrum of the $K$-algebra endomorphism $e^{\#}:R\rightarrow R$ that maps $x_i$ to $f_i$ for each $i\in \llbracket1,m\rrbracket$. \phantomsection{We define the {\it degree} of $e$ by}\label{PH7}
$$\pi(e)=\pi(f_1,\ldots,f_n):=\max\bigl(\deg(f_i)|i\in \llbracket1,n\rrbracket\bigr)\in\mathbb N\cup\{-\infty\}$$
(cf.\ \cite{W}, Def.\ 1.8). If $a\in\GA_n(K)$, then $\bigl(\pi(a),\pi(a^{-1})\bigr)\in (\mathbb N^{\ast})^2$ and we define its {\it degree length}\index{length function!degree length}\footnote{The terminology `degree length' is used to avoid, for $n=2$, confusion with: (i) the usual length associated to the amalgamated product structure of $\GA_2(K)$ as in \cite{FM}, Sect.\ 2, Def.; (ii) Furter's length that was introduced in \cite{F1}, Sect.\ 4, Def.\ for $K=\mathbb C$ and in \cite{F2}, Sect.\ 1 c. for all fields and that is also used in Sections \ref{S12} and \ref{S27} for all (finite) fields; (iii) the affine length introduced in \cite{FP}, Def.\ 4.1 for all fields $K$.} by
\begin{equation}\label{EQ1}
\ell(a):=\max\bigl(\pi(a),\pi(a^{-1})\bigr).
\end{equation}
Hence 
\begin{equation}\label{EQ2}
\ell(a)=\ell(a^{-1}).
\end{equation} 
Note that $\ell(a)\in \llbracket\pi(a),\pi(a)^{n-1}\rrbracket$ by \cite{BCW}, Cor.\ (1.4); this result is optimal, for instance, by Example \ref{EX9}.
Thus for $n\in\{1,2\}$, for each $a\in\GA_n(K)$ we have identities
$$\ell(a)=\pi(a)=\pi(a^{-1}).$$

For $(a,b)\in\GA_n(K)^2$, we have an inequality $\pi(ab)\le\pi(a)\pi(b)$, and hence
\begin{equation}\label{EQ3}
\ell(ab)\le\ell(a)\ell(b).
\end{equation}
\phantomsection{Let $1_{\mathbb A^n_K}$ be the identity automorphism of $\mathbb A^n_K$. Clearly, $\ell(1_{\mathbb A^n_K})=\pi(1_{\mathbb A^n_K})=1$. Based on this and Equations (\ref{EQ2}) and (\ref{EQ3}), the rule $a\mapsto\ell(a)$ defines a function\index{length function!$\ell_{\GA_n(K)}$ degree length function}}\label{PH8}
$$\ell_{\GA_n(K)}:\GA_n(K)\rightarrow\mathbb N^{\ast}$$
which is a length function in the sense recalled in Definition \ref{D11}(1); it is the smallest length function such that for each $a\in\GA_n(K)$ we have $\ell(a)\ge\pi(a)$. 

\phantomsection{For $l\in\mathbb N^{\ast}$ and $\Gamma\subset \GA_n(K)$, let $\Gamma[l]:=\{a\in\Gamma|\ell(a)\le l\}$; we have inclusions}\label{PH8a}
$$\Gamma[1]\subset\Gamma[2]\subset\cdots\subset\Gamma[l]\subset\cdots\subset\Gamma.$$ If $\Gamma$ is a subgroup of $\GA_n(K)$, \phantomsection{then Equation (\ref{EQ3}) implies that for $(l,s)\in (\mathbb N^{\ast})^2$ we have an inclusion $\Gamma[l]\Gamma[s]\subset\Gamma[ls]$, hence $\Gamma[1]$ is a subgroup of $\Gamma$ but for $l\ge 2$, in general $\Gamma[l]$ is not a subgroup of $\Gamma$. In particular,}\label{PH9}
$$\AGL_n(K):=\GA_n(K)[1]$$
\phantomsection{is the subgroup of $\GA_n(K)$ of affine automorphisms\index{automorphism!affine automorphism} and}\label{PH10}
$$\ASL_n(K):=\SGA_n(K)[1]$$ 
is the subgroup of $\GA_n(K)$ of special affine automorphisms\index{automorphism!special affine automorphism}. 

\phantomsection{The affine (resp.\ special affine) automorphisms that fix $(0,\ldots,0)$ are called linear automorphisms\index{automorphism!linear automorphism} (resp.\ special linear automorphisms\index{automorphism!special linear automorphism}); they form a subgroup of $\GGA_n(K)$ (resp.\ $\SGA_n(K)$) naturally identified with $\GL_n(K)$ (resp.\ $\SL_n(K)$) as $\GL_n(K)$ (resp.\ $\SL_n(K)$ is generated by its $\mathbb G_{\a,K}$ and $\mathbb G_{\m,K}$ subgroups (resp.\ generated by its $\mathbb G_{\a,K}$ subgroups). The notation used is justified by the existence of a short exact sequence $1\rightarrow K^n\rightarrow\AGL_n(K)\rightarrow\GL_n(K)\rightarrow 1$ (resp.\ $1\rightarrow K^n\rightarrow\ASL_n(K)\rightarrow\SL_n(K)\rightarrow 1)$, with $K^n$ identified with the subgroup of translations of $\mathbb A^n_K$ or $K^n$.}\label{PH10a}

\phantomsection{Let $\TGA_n(K)$ be the subgroup of $\GA_n(K)$ of {\it tame automorphisms}\index{automorphism!tame automorphism} generated by affine automorphisms and, if $n\ge 2$, by automorphisms defined by $n$-tuples of the form $(x_1,\ldots,x_{n-1},x_n+f)$ with $f\in K[x_1,\ldots,x_{n-1}]$ a linear combination of primitive monomials of degree at least $2$.}\label{PH10b} 
\phantomsection{The intersection}\label{PH11}
$$\STGA_n(K):=\TGA_n(K)\cap\SGA_n(K)$$
is the subgroup of special tame automorphisms\index{automorphism!special tame automorphism} of $\mathbb A^n_K$.

\phantomsection{We have $\TGA_1(K)=\AGL_1(K)=\GA_1(K)$. Recall that $\TGA_2(K)=\GA_2(K)$ by \cite{vdK}, Thms.\ 1 and 2; thus $\STGA_2(K)=\SGA_2(K)$ and we only use the shorter $\SGA_2(K)$. Clearly, $\TGA_n(K)$ is a subgroup of $\GGA_n(K)$ and the Jacobian determinant surjective homomorphism $\det:\GA_n(K)\rightarrow K^{\ast}$ induces a short exact sequence}\label{PH90}
$$1\rightarrow \STGA_n(K)\rightarrow\TGA_n(K)\xrightarrow{\det} K^{\ast}\rightarrow 1.$$ 
As for $n>2$ we do not know when the inclusions $\TGA_n(K)\subset\GA_n(K)$ and $\STGA_n(K)\subset\SGA_n(K)$ are identities, we work primarily with the subgroup $\TGA_n(K)$ of $\GGA_n(K)$ but we point out the differences in the case of $\STGA_n(K)$.

If $K$ has positive characteristic, then the group of automorphisms $\GA(\mathbb G_{\a,K}^n)$ of the affine group variety $\mathbb G_{\a,K}^n$ over $\Spec K$, when viewed by forgetting the group scheme structure $\mathbb G_{\a,K}^n$ on $\mathbb A^n_K$ as a subgroup of $\GA_n(K)$, is a subgroup of $\TGA_n(K)$ by \cite{Kur}, Cor.\ 2.3 (the particular case when $K$ is algebraically closed was first proved in \cite{TK}, Thm.\ 1).

\begin{definition}\label{D2}
Let $(n,m)\in (\mathbb N^{\ast})^2$. Let $K$ be a field such that $|K|\ge\sqrt[n]{m}$. 

\medskip
{\bf (1)} We call $K$ a $T_{n,m}$ field\index{field!$T_{n,m}$ field} (resp.\ an $ST_{n,m}$ field)\index{field!$S_{n,m}$ field} if the action 
$$\TGA_n(K)\times \mathbb D_{n,m}(K)\rightarrow \mathbb D_{n,m}(K)\;\; \textup{(resp.}\;\; \STGA_n(K)\times \mathbb D_{n,m}(K)\rightarrow \mathbb D_{n,m}(K))$$ induced by $\mathbb T_{n,m}(K)$ is transitive.

\smallskip
{\bf (2)} For $(\underline{P},\underline{Q})\in\mathbb D_{n,m}(K)^2$, let $\pi_{\underline{P},\underline{Q}}$ (resp.\ $\pi^{\S}_{\underline{P},\underline{Q}}$) in $\mathbb N^{\ast}\cup\{\infty\}$ be the infimum of the set $\{\ell(a)|a\in\TGA_n(K), a(\underline{P})=\underline{Q}\}$ (resp.\ $\{\ell(a)|a\in\STGA_n(K), a(\underline{P})=\underline{Q}\}$).

\smallskip
{\bf (3)} Let $\pi_{n,m}(K)$\index{$\pi_{n,m}(K)$ main invariant} (resp.\ $\pi^{\S}_{n,m}(K)$\index{$\pi^{\S}_{n,m}(K)$ main invariant}) in $\mathbb N^{\ast}\cup\{\infty\}$ be the supremum of the set $\{\pi_{\underline{P},\underline{Q}}|(\underline{P},\underline{Q})\in\mathbb D_{n,m}(K)^2\}$ (resp.\ $\{\pi^{\S}_{\underline{P},\underline{Q}}|(\underline{P},\underline{Q})\in\mathbb D_{n,m}(K)^2\}$).
\end{definition}

The conventions are as follows: the infimum of the empty set is $\infty$ and the supremum of each subset of $\mathbb N\cup\{\infty\}$ that contains $\infty$ is $\infty$. 

Clearly, we have $\pi_{n,m}(K)\le\pi^{\S}_{n,m}(K)$. A field $K$ is a $T_{1,m}$ field iff $m\in\{1,2\}$ and is an $ST_{1,m}$ field iff either $m=1$ or $m=|K|=2$; moreover, we have identities $\pi^{\S}_{n,1}(K)=\pi^{\S}_{n,2}(K)=1$ for each $n\in N^{\ast}\setminus\{1\}$. Each infinite field $K$ is a $T_{n,m}$ field for every pair $(n,m)\in (\mathbb N^{\ast})^2$ with $n\ge 2$ by \cite{Sr}, Thm.\ 2 and the proof of \cite{KZ}, Lem.\ 5.5 worked over $\mathbb C$ can be adapted to give that it is also an $ST_{n,m}$ field for each $(n,m)\in (\mathbb N^{\ast})^2$ with $n\ge 2$. But the proofs in \cite{Je}, \cite{Kali}, \cite{Sr}, \cite{KZ}, or \cite{AFKKZ} do not give information on the $\pi_{n,m}(K)$s or $\pi^{\S}_{n,m}(K)$s.

\phantomsection{For a finite field $K$, we consider the homomorphism}\label{PH12} 
$$\varrho_{n,K}:\GA_n(K)\rightarrow\perm(K^n)$$
that defines $\mathbb T_{n,m}(K)$. 
Its restriction to $\AGL_n(K)$ is injective, hence, as in the literature, we identify $\AGL_n(K)$ and its subgroups with their images under $\varrho_{n,K}$.

\phantomsection{If $q\in\mathbb N^{\ast}$ and $p$ is a prime, then recall that $\mathbb F_{p^q}$ is a field with $p^q$ elements, uniquely determined up to isomorphisms. We do not know literature on finite fields $\mathbb F_{p^q}$ that are $T_{n,m}$ fields or on the $\pi_{n,m}(\mathbb F_{p^q})$s or $\pi^{\S}_{n,m}(\mathbb F_{p^q})$s, but various classifications of maximal subgroups of finite alternating or symmetric groups can be used to study finite $T_{n,m}$ fields. For instance, for $n\ge 2$ the group $\AGL_n(\mathbb F_p)$ is a maximal subgroup of $\perm(\mathbb F_p^n)$ (e.g., see \cite{LPS}, Thm.(I)) and based on this we get that $\varrho_{n,\mathbb F_p}\bigl(\TGA_n(\mathbb F_p)\bigr)=\perm(\mathbb F_p^n)$ by simply identifying an $a\in\TGA_n(\mathbb F_p)$ such that $\varrho_{n,\mathbb F_p}(a)\in\perm(\mathbb F_p^n)\setminus \AGL_n(\mathbb F_p)$. If $q\ge 2$, then each $\mathbb F_p$-linear isomorphism $\mathbb F_{p^q}^n\cong\mathbb F_p^{nq}$ allows us to view $\AGL_{nq}(\mathbb F_p)$ as a subgroup of $\perm(\mathbb F_{p^q}^n)$ that contains $\AGL_n(\mathbb F_{p^q})$, and thus $\AGL_n(\mathbb F_{p^q})$ is not a maximal subgroup of $\perm(\mathbb F_{p^q}^n)$; however, based on the fact that $\AGL_{nq}(\mathbb F_p)\cap\Alt(\mathbb F_p^{nq})$ is a maximal subgroup of $\Alt(\mathbb F_p^{nq})$ (e.g., see \cite{LPS}, Thm.(I)) one can check using \cite{Kur}, Cor.\ 2.3 or \cite{TK}, Thm.\ 1 that the subgroup of $\perm(\mathbb F_{p^q}^n)$ generated by $\AGL_n(\mathbb F_{p^q})$ and $\varrho_{n,\mathbb F_{p^q}}\bigl(\GA(\mathbb G_{\a,\mathbb F_{p^q}}^n)\bigr)$ is $\AGL_{nq}(\mathbb F_p)$ if $p$ is odd and is $\AGL_{nq}(\mathbb F_p)\cap\Alt(\mathbb F_p^{nq})$ if $p=2$.}\label{PH12a}

\phantomsection{If $K$ is finite (resp.\ infinite), then $\bigl(\pi_{n,m}(K)\bigr)_{m=1}^{|K|^n}$ and $\bigl(\pi^{\S}_{n,m}(K)\bigr)_{m=1}^{|K|^n-2}$ (resp.\ $\bigl(\pi_{n,m}(K)\bigr)_{m\ge 1}$ and $\bigl(\pi^{\S}_{n,m}(K)\bigr)_{m\ge 1}$) are non-decreasing sequences in $\mathbb N^{\ast}\cup\infty\}$ (resp.\ in $\mathbb N^{\ast}$). If $K$ is finite, then it is well-known that for each $\sigma\in\perm(K^n)$ there exists $e\in\End_n(K)$ such that $\deg(e)\le n(|K|-1)$ and $e(K)=\sigma$. But this $e$ is, in general, not an automorphism. In fact, there exists $c\in\GA_n(K)$ such that $c(K)=\sigma$ iff $4\nmid |K|$ or $\sigma$ is an even permutation (see Theorem \ref{T6}(2) and (3)). And even if such a $c\in \GA_n(K)$ exists, its smallest possible degree length is in general much larger than $n(|K|-1)$. Establishing upper and lower bounds on $\deg(c)$ when $c$ exists, by various methods, in general and in some natural subcases, is one of the main goals of this monograph (see Sections \ref{S17} through \ref{S28}).}\label{PH12b}

The monograph is organized as follows. 

Section \ref{S2} recalls, for a finite set $A$, basic notation and results on $\perm(A)$ and its alternating subgroup $\Alt(A)$. 

Section \ref{S3} contains decompositions results for permutations in $\perm(A)$ that are required in the study (from Section \ref{S11} onwards) of the degree length functions on $\varrho_{n,K}\bigl(\TGA_n(K)\bigr)$ and $\varrho_{n,K}\bigl(\STGA_n(K)\bigr)$ with $K$ a finite field.

Section \ref{S4} studies homomorphisms $\varrho_{\Lambda}:G(K)\rightarrow\perm\bigl(X(K)\bigr)$ attached to actions $\Lambda:G\times_{\Spec K} X\rightarrow X$ in the case when $K$ is a finite field and $G$ is a linear algebraic group over $\Spec K$. Theorem \ref{T3} is a criterion on when $\Im(\varrho_{\Lambda})\leqslant\Alt\bigl(X(K)\bigr)$; it implies that an affine integral variety over a finite field of odd characteristic is flexible iff it has at most one valued point in the field (see Corollary \ref{C5}). Proposition \ref{PR5} is a criterion on when two smooth closed subschemes of $G$ coincide.

Section \ref{S5} introduces first invariants of finite subsets of $K^n$ that are often used in what follows.

Section \ref{S6} studies arithmetic functions related to partitions of such finite subsets into linearly independent subsets.

Section \ref{S7} uses these invariants and arithmetic functions to represent functions $h:Y\rightarrow K$ on finite subsets $Y$ of $K^n$ by polynomials $f\in R$ of bounded degrees and thus provides refined Lagrange interpolations in multiple indeterminates.

Section \ref{S8} uses Section \ref{S5} to define tame automorphisms $a\in\TGA_n(K)$ for which, when $K$ is finite, the permutation $a(K)$ of $K^n$ has nice properties required in the subsequent proofs. In particular, Proposition \ref{PR10}(1) shows how to construct automorphisms $a\in\STGA_n(\mathbb F_{p^q})$ with $a(\mathbb F_{p^q})$ an even cycle of arbitrary length at most $2p-1$.

Section \ref{S9} contains transitivity properties related to the action $\mathbb T_{n,m}$. In particular, it computes $\varrho_{n,K}\bigl(\TGA_n(K)\bigr)$, $\varrho_{n,K}\bigl(\STGA_n(K)\bigr)$, and $\varrho_{n,K}\bigl(\SGA_n(K)\bigr)$ when $K$ is finite (see Theorem \ref{T6}(2) to (4)) using basic properties of alternating groups, Proposition \ref{PR10}(1), and Theorem \ref{T3}. The computations allow us to obtain in Subsection \ref{S9.1} the following classifications.

\begin{theorem}\label{T1} Let $K$ be a field. Let $(n,m)\in\mathbb N^{\ast}\setminus\{1\}\times\mathbb N^{\ast}$. Then the following properties hold.

\medskip
{\bf (1)} Suppose that $|K|\ge\sqrt[n]{m}$. Then $K$ is not an $ST_{n,m}$ field iff $K$ is finite with at least $3$ elements and $m\in\{|K|^n-1,|K|^n\}$. 

\smallskip
{\bf (2)} Suppose that $|K|\ge\sqrt[n]{m}$. Then $K$ is not a $T_{n,m}$ field iff there exists an integer $q\in\mathbb N^{\ast}\setminus\{1\}$ such that $K\cong\mathbb F_{2^q}$ and $m\in\{2^{qn}-1,2^{qn}\}$.

\smallskip
{\bf (3)} The field $K$ is not infinitely transitive iff $K\cong\mathbb F_{2^q}$ with $q\in\mathbb N^{\ast}\setminus\{1\}$.
\end{theorem}

We have the following consequence of Theorem \ref{T1}(1) or (2) (see Subsection \ref{S9.2}).

\begin{corollary}\label{C1}
Let $n\in\mathbb N^{\ast}\setminus\{1\}$ and $K$ a field. Let $m\in\mathbb N$ be such that $|K|\ge\sqrt[n]{m}$. Then the following properties hold.

\medskip
{\bf (1)} Given two subsets $Y$ and $Z$ of $K^n$ such that $|Y|=|Z|=m$, there exists $a\in\STGA_n(K)$ such that $a(Y)=Z$.

\smallskip
{\bf (2)} There exists a smallest $\upsilon^{\S}_{n,m}(K)\in\mathbb N^{\ast}$\index{$\upsilon^{\S}_{n,m}(K)$ main invariant} such that for each subsets $Y$ and $Z$ of $K^n$ of cardinality $m$, there exists $a\in\STGA_n(K)[\upsilon^{\S}_{n,m}(K)]$ with $a(Y)=Z$.

\smallskip
{\bf (3)} There exists a smallest $\upsilon_{n,m}(K)\in\mathbb N^{\ast}$\index{$\upsilon_{n,m}(K)$ main invariant} such that for each subsets $Y$ and $Z$ of $K^n$ of cardinality $m$, there exists $a\in\TGA_n(K)[\upsilon_{n,m}(K)]$ with $a(Y)=Z$.
\end{corollary}

Clearly, we have an inequality $\upsilon_{n,m}(K)\le\upsilon^{\S}_{n,m}(K)$.

With the notation of Corollary \ref{C1}, the isomorphism class of 
$$\mathbb B^n_{m,K}:=\mathbb A^n_{m,K}\setminus Y$$ depends only on $(n,m,K)$ and not on $Y$, hence we call it {\it the affine space of dimension $n$ over $\Spec K$ punctured in $m$ distinct $K$-valued points}. Note that $\mathbb B^n_{0,K}=\mathbb A^n_{K}$.

Restricting to finite fields $K$, as samples of what various quantitative classifications of $T_{n,m}$ fields when $4\nmid |K|$ (resp.\ $4\mid |K|^n$) and $n\ge 2$ would mean we mention only the following two problems considered to a certain extent in this monograph. 

\medskip\noindent
{\bf (QP1)} Compute or estimate the $\pi_{n,m}(K)$s and $\pi^{\S}_{n,m}(K)$s when $m\in \llbracket3,|K|^n\rrbracket$.

\smallskip\noindent
{\bf (QP2)} \phantomsection{For each $s\in \llbracket2,|K|^n\rrbracket$ (resp.\ odd $s\in \llbracket2,|K|^n\rrbracket$), compute the set}\label{PH14}
$$\Pi_{n,s\textup{-cycle}}(K):=\{\ell(a)|a\in\TGA_n(K),\, a(K)\in\perm(K^n)\,\textup{is an}\,s\textup{-cycle}\}.$$

\medskip
For all fields we mention the following third problem. 

\medskip\noindent
{\bf (QP3)} Compute the $\pi_{Y,s}(K)$s and $\pi^{\S}_{Y,s}(K)$s introduced in the following definition.

\begin{definition}\label{D2.5}
Let $n\in\mathbb N^{\ast}\setminus\{1\}$ and $K$ a field. Let $s\in\mathbb N^{\ast}$ and a finite non-empty subset $Y$ of $K^n$ be such that $s+|Y|\le |K|^n$. 

\medskip
{\bf (1)} For a subgroup $H$ of $\GA_n(K)$, let $\Fix_H(Y)$ be the subgroup of $H$ that fixes each element of $Y$ (thus $\Fix_H(Y)=H\cap\Fix_{\GA_n(K)}(Y)$). 

\smallskip
{\bf (2)} If $K$ is a finite field with $4\mid |K|$, then we assume that $s+|Y|\le |K|^n-2$. Let $\pi_{Y,s}(K)\in\mathbb N^{\ast}$\index{$\pi_{Y,s}$ main invariant} be the smallest such that for each pair $(\underline{P},\underline{Q})\in\mathbb D_s(\mathbb A^n_K\setminus Y)^2$ there exists $a\in\Fix_{\TGA_n(K)}(Y)[\pi_{Y,s}(K)]$ with $a(\underline{P})=\underline{Q}$.

\smallskip
{\bf (3)} If $K$ is a finite field, then we assume that $s+|Y|\le |K|^n-2$. We define $\pi^{\S}_{Y,s}(K)\in\mathbb N^{\ast}$ to be the smallest such that for each $(\underline{P},\underline{Q})\in\mathbb D_s(\mathbb A^n_K\setminus Y)^2$ there exists $a\in\Fix_{\STGA_n(K)}(Y)[\pi^{\S}_{Y,s}(K)]$ with $a(\underline{P})=\underline{Q}$.
\end{definition}

Section \ref{S10} reviews types of {\it length functions} and shows that \cite{vdK}, Thms.\ 1 and 2 imply that the length function $\ell_{\GA_2(K)}$ is regular in the sense of \cite{Pro}, Sect.\ 4, Def.\ recalled in Definition \ref{D11}(3). 

Section \ref{S11} translates the quantitative problems (QP1) to (QP3) in terms of the {\it length functions} $\ell_{n,K}:\varrho_{n,K}\bigl(\TGA_n(K)\bigr)\rightarrow\mathbb N^{\ast}$ and $\ell^{\S}_{n,K}:\varrho_{n,K}\bigl(\STGA_n(K)\bigr)\rightarrow\mathbb N^{\ast}$. For instance, for a permutation $\sigma\in\varrho_{n,K}\bigl(\TGA_n(K)\bigr)$, $\ell_{n,K}(\sigma)$ is the minimum of the set $\{\ell(a)|a\in\TGA_n(K),\, a(K)=\sigma\}$. We note that their real-valued natural logarithms $\ln\ell_{n,K}:\varrho_{n,K}\bigl(\TGA_n(K)\bigr)\rightarrow\mathbb R$ and $\ln\ell^{\S}_{n,K}:\varrho_{n,K}\bigl(\STGA_n(K)\bigr)\rightarrow\mathbb R$ satisfy Axioms $A_2$ and $A_3$ of \cite{Ha}, Sect.\ 1 and Axiom $A 1^{\prime}$ of \cite{Ch} that are variations of Axioms $A1$ to $A3$ of \cite{L}, Sect.\ 2. So $\ln\ell_{n,K}$ and $\ln\ell^{\S}_{n,K}$ are $\mathbb R$-valued semigauges\index{semigauge} in the terminology of \cite{Pro}, Sect.\ 2. 

Section \ref{S12} introduces only for $n=2$ the Furter's length functions $\ell_{2,K}^-$ and $\ell_{2,K}^+$ on $\Perm(K^2)$ for which we have inequalities $\ell_{2,K}^-\le\ell_{n,K}\le\ell_{2,K}^+$. Their properties are vital for different estimates and proofs as they are based not on degrees of polynomials but on the number of suitable composite factors of automorphisms in $\GA_2(K)$ that represent a suitable permutation in $\Perm(K^2)$. We have $\ell_{2,K}=\ell^-_{n,K}$ (or $\ell_{n,K}=\ell_{2,K}^+$) iff $|K|=3$ by Corollary \ref{C13}.

Section \ref{S13} provides estimates of values of the length functions $\ell_{n,K}$ and $\ell^{\S}_{n,K}$ via unions of conjugacy classes in $\perm(K^n)$ or $\Alt(K^n)$, such as conjugacy classes of $r$ disjoint cycles of length $s$ with $(r,s)\in \mathbb N^{\ast}\times (\mathbb N^{\ast}\setminus\{1\})$. For applications to (QP2), see Lemmas \ref{L8} and \ref{L9}, Proposition \ref{PR16}, and Example \ref{EX17}.

Section \ref{S14} recalls Nagata's automorphism over $\Spec K$ introduced in \cite{N}, Sect.\ 2.1 and lists properties of it and of some variations of it that relate to prior sections.

Section \ref{S15} studies additional invariants of finite subsets of $K^n$ that are related to directions parametrized by the $K$-valued points of the projective space $\mathbb P^{n-1}_K$ and that are essential in computing upper bounds of the $\pi_{n,m}(K)$s.

Section \ref{S16} uses Sections \ref{S13} and \ref{S15} to construct automorphisms $a\in\STGA_n(K)$ with $a(K)$ as products of a priori given numbers of disjoint $s$-cycles and with explicit upper bounds for $\ell(a)$ which often are optimal, where $s\in\mathbb N^{\ast}\setminus\{1\}$.

Sections \ref{S17} to \ref{S24} pertain to finite fields $K$ and provide explicit upper bounds for the $\pi_{n,m}(K)$s and $\pi^{\S}_{n,m}(K)$s in various ranges of $m$ in comparison to $|K|$ and $n$, by using a wide array of techniques to construct the required automorphisms. As the field $\mathbb F_2$ is often used and has a great potential in encryptions (for instance, see \cite{M}) and as the fields $\mathbb F_2$ and $\mathbb F_3$ are simpler for computations and refinements, many of our results in Sections \ref{S17} to \ref{S24} are stated separately for $|K|=2$ and $n\ge 3$ and for $|K|=3$ and $n\ge 2$.

Section \ref{S17} provides explicit upper bounds for the $\pi_{n,m}(K)$s and the $\pi^{\S}_{n,m}(K)$s for `small' $m$, by which we mean that $|K|\ge m$ or $n$ is greater than some suitable functions of $|K|$ and $m$ such as $n>m$ (see Theorems \ref{T7} and \ref{T8}). The upper bounds are functions of $m$ which in most cases are polynomials of degree $2$ or $3$. The methods used are basic combinatorics on linear projections and projective spaces over $\Spec K$. 

Section \ref{S18} refines the techniques of Section \ref{S17} by using Weil restriction of scalars (see Lemma \ref{L13}) and `relative' settings that involve $K^{n+n_1}$ with $n$ even and $n_1\in\mathbb N$ and selections (in pairs) of $m_O$-tuples of distinct points in $K^n\times\{O\}$ indexed by $O\in K^{n_1}$. The resulting upper bounds are polynomial in $m$ (see Proposition \ref{PR21} for $n$ even and Proposition \ref{PR23} for $n$ odd), with $m$ ranging from $|K|$ to $|K|^{\frac{n}{2}}$ if $n$ is even and to values slightly greater than or equal to $\bigl\lfloor\sqrt{\frac{|K|^n}{2}}\bigr\rfloor$ if $n$ is odd.

Section \ref{S19} presents the principles that allow us to push this to an upper range of $m$ up to $2\lceil\frac{|K|^{n-1}}{4}\rceil$ when $n\ge 3$, giving almost polynomial upper bounds that are polynomial in $m$ for fixed $n$ but with constants explicitly depending on $n$ (see Theorems \ref{T9} and \ref{T10}).

Section \ref{S20} translates the principles of Section \ref{S19} into concrete estimates of $\ell^{\S}_{n,K}(\sigma)$ for each $\sigma\in\Alt(K^n)$ and every $m\le 2\lceil\frac{|K|^{n-1}}{4}\rceil$ (see Corollary \ref{C20}).

Section \ref{S21} combines Sections \ref{S13} and \ref{S20} to explicitly estimate all values of $\ell^{\S}_{n,K}$. 

Section \ref{S22} deals with $m$ up to $|K|^n$ by combining results from Section \ref{S17} and the earlier Sections \ref{S2}, \ref{S7}, \ref{S9}, and \ref{S11} in the context of affine automorphisms that simplify the required permutations of $K^n$ by controlling the cardinalities of their supports. 

Sections \ref{S23} and \ref{S24} push this approach further for $n=2$ and $n\ge 3$ (respectively), by using the more sophisticated results from Sections \ref{S18} and \ref{S19} instead of Section \ref{S17} to get asymptotically better bounds of the form 
$$\pi_{n,m}(K)<e^{(C_1\sqrt{m}+C_2)[\ln(C_3m)]}\le e^{C\sqrt{m}\ln m}$$ 
with $C_1$, $C_2$, and $C_3$ (hence also $C$) as explicit universal positive constants, see Theorems \ref{T12} and \ref{T13}. In Section \ref{S24}, for $n\ge 3$ even better bounds 
$$\pi_{n,m}(K)<e^{C_4(2m)^{\frac{1}{n}}(\ln 2m)^2\bigl[C_5\ln\bigl(\ln(2m)\bigr)+C_6\bigr]}\le e^{Dm^{\frac{1}{n}}(\ln m)^2 \ln (\ln m)}$$ 
\phantomsection{are obtained with $C_4$ to $C_6$ (hence also $D$) explicit universal constants with $C_4$ and $C_5$ (hence also $D$) positive, see Theorem \ref{T14} and Corollary \ref{C24}. Note that the arguments in Sections \ref{S23} and \ref{S24} are mostly of analytic, rather than algebraic, nature, and between larger constants and larger $m$s we have chosen $m$ to be fairly large, i.e., $m\in\llbracket N,|K|^n\rrbracket$ if $4\nmid |K|$ and $m\in\llbracket N,|K|^n-2\rrbracket$ if $4\mid |K|$ with the lower bound $N\in\llbracket122,676\rrbracket$, and the constants to be small, closer to optimality.}\label{EXT5}

Section \ref{S25} pertains to (QP3): it provides inductively on $m$ upper bounds for the $\pi_{n,m}(K)$s and $\pi^{\S}_{n,m}(K)$ based on the study of stabilizers of the action $\mathbb T_{n,m}(K)$ and on upper bounds for the $\pi_{Y,s}$s and $\pi^{\S}_{Y,s}$s. 

Section \ref{S26} provides lower bounds for the $\pi_{n,m}(K)$s with $m$ small based on prior works of Abhyankar--Moh, Richman, and Kanga (see \cite{AM}, \cite{Ri}, and \cite{Kan}) and on a new notion of $\llbracket1,m-1\rrbracket$-generic $m$-tuples in $\mathbb D_{n,m}(K)$ (see Definition \ref{D23}(1); for $m\ge 3$, such $m$-tuples exist iff $|K|\ge m+1$ by Proposition \ref{PR28}(5). 

The following direct consequence of Theorems \ref{T7}(1) and \ref{T15}(4) and (5) refines the above mentioned consequence of \cite{Sr}, Thm.\ 1.2 for $n=2$.

\begin{corollary}\label{C2}
For $(r,\epsilon)\in\mathbb N^{\ast}\times\{0,1,2\}$, let $\phi(r,\epsilon):=4r^2+\epsilon(5-\epsilon)r+\epsilon(\epsilon-1)$ and $m:=3r+\epsilon$. Let $K$ be a field. Then the following inequalities hold
$$\frac{4m^2}{9}\le\phi(r,\epsilon)\le\pi_{2,m}(K)\le\pi^{\S}_{2,m}(K)\le 9r^2+6r(\epsilon-1)+(\epsilon-1)^2= (m-1)^2.$$
If moreover $K$ contains a primitive $m$-th root of unity, then we have
$$\frac{4m^2}{9}\le\phi(r,\epsilon)\le\Bigl\lfloor\frac{m+2}{2}\Bigr\rfloor(m-1)\le\pi_{2,m}(K)\le\pi^{\S}_{2,m}(K)\le (m-1)^2$$
for $m\ge 7$ and $\pi_{2,m}(K)=\pi^{\S}_{2,m}(K)=(m-1)^2$ for $m\in\{5,6\}$.\end{corollary}

From the case $l=1$ of Theorem \ref{T7}(1) and from Theorem \ref{T15}(1) we get directly the following practical result which also refines the above mentioned consequence of \cite{Sr}, Thm.\ 1.2 for $n\ge 3$.

\begin{corollary}\label{C3}
If $n\ge 3$, $m\ge 4$, and $|K|\ge \frac{m(m-1)}{2}$, then 
$$m+1\le\pi_{n,m}(K)\le\pi^{\S}_{n,m}(K)\le (m-1)^2.$$
\end{corollary}

\phantomsection{Section \ref{S27} uses a reformulation of van der Kulk's results that proved that $\TGA_2(K)=\GA_2(K)$ to get lower bounds for $\pi_{2,m}(K)$ and $\pi^{\S}_{2,m}(K)$ when $K$ is a finite field and $m\in \llbracket3,|K|^2\rrbracket$ (see Theorems \ref{T16} and \ref{T17}). For instance, as applications of Theorems \ref{T12}, \ref{T16}, and \ref{T17}, for each $\varepsilon\in (0,2)$ we have inequalities}\label{EXT3}
$$(2-\varepsilon) |K|\ln |K|<\ln\bigl(\pi_{2,|K|^2-2}(K)\bigr)<18|K|\ln |K|$$ 
if $|K|\gg 1$ and 
$$(2-\varepsilon) |K|\ln |K|<\ln\bigl(\pi_{2,|K|^2}(K)\bigr)<18|K|\ln |K|$$ 
if $|K|\gg 1$ is odd (see Corollary \ref{C27}(1) and (3)).

Section \ref{S28} uses the bounds obtained to compute exact values or small ranges of values for many of the $\pi_{n,m}(K)$s and $\pi^{\S}_{n,m}(K)$s with $m\le 7$. For example, $\pi_{n,3}(K)$ is $4$ if $|K|\ge 4$, is $2$ if $K\cong \mathbb F_3$, and is $1$ if $K\cong\mathbb F_2$ by Theorem \ref{T19}(1) and Lemma \ref{L11}(2). Similarly, $\pi_{2,4}(K)$ is $9$ if $|K|\ge 5$, is $6$ if $K\cong \mathbb F_4$, is $4$ if $K\cong \mathbb F_3$, and is $1$ if $K\cong\mathbb F_2$ by Theorem \ref{T19}(2), (3), and (5) and Lemma \ref{L11}(2); also, $\pi^{\S}_{2,4}(K)=9$ if $|K|\ge 7$ by Theorem \ref{T19}(2).

Section \ref{S29} gathered basic properties of the $\upsilon_{n,m}(K)$s and $\upsilon^{\S}_{n,m}(K)$s such as their symmetric property (see Lemma \ref{F14}(2)) and the classification of all cases when they are $1$ (see Lemma \ref{L28}).

Section \ref{S30} applies the methods of the prior sections to the $\upsilon_{n,m}(K)$s, $\upsilon^{\S}_{n,m}(K)$s, and analogs $\upsilon_{Y,s}(K)$s and $\upsilon^{\S}_{Y,s}(K)$s of the $\pi_{Y,s}(K)$s and $\pi^{\S}_{Y,s}(K)$s for sets. For instance, for each prime $p$ and $\varepsilon\in (0,2)$ we have inequalities 
$$(2-\varepsilon) |K|\ln |K|<\ln\bigl(\upsilon_{2,\lfloor\frac{|K|^2}{2}\rfloor}(K)\bigr)<\bigl[\min(14,\upsilon_{2,p})+\varepsilon\bigr]|K|\ln |K|$$ 
if $K$ has characteristic $p$ and $|K|\gg 1$, where we have $\upsilon_{2,2}:=9$, $\upsilon_{2,3}:=13.5$, and $\upsilon_{2,p}:=\frac{27p}{2(p-1)}$ if $p\ge 5$ (see Corollary \ref{C31}).

Section \ref{S31} considers sequences $\bigl((K_i,m_i)\bigr)_{i\in\mathbb N^{\ast}}$ of pairs with each $K_i$ a finite field and $m_i\in \llbracket2,|K_i|^2\rrbracket$ and shows that the existence of $C\in (0,\infty)$ such that we have $\star_{2,m_i}(K_i)<m_i^C$ for each $i\in\mathbb N^{\ast}$, where $\star$ is one of the main invariants such as $\pi$ and $\upsilon$, is equivalent to the boundedness of the set $\bigl\{\frac{m_i}{|K_i|}|i\in\mathbb N^{\ast}\bigr\}$ of rational numbers. 

Section \ref{S32} lists open problems. 

\section{Finite symmetric groups, part I: notation and basic properties}\label{S2}

We introduce notation and review properties of finite groups of permutations that are more or less well-known and that are often used in what follows.

\begin{notation}\normalfont\label{N1}
For $l\in\mathbb N^{\ast}$ let $S_l:=\perm(\llbracket1,l\rrbracket)$. Let $A_l$ be the alternating subgroup of $S_l$. For a subset $\nabla\subset S_l$ and $n\in\mathbb N^{\ast}$, let $\nabla^n$ be the set of products of $n$ elements in $\nabla$. 

For $\sigma\in S_l$, let $o(\sigma)\in\mathbb N^{\ast}$ be its order, let $\supp(\sigma):=\{i\in \llbracket1,l\rrbracket|\sigma(i)\neq i\}$ be its support, and let $\n(\sigma):=|\supp(\sigma)|\in \llbracket0,l\rrbracket$ be the number of its non-fixed elements. 

For $s\in\llbracket2,l\rrbracket$ let $\c_s(\sigma)\in \llbracket0,\lfloor\frac{l}{s}\rfloor\rrbracket$ be the number of $s$-cycles of $\sigma$. Also, let $\c_{\textup{even}}(\sigma):=\sum_{s=1}^{\lfloor\frac{l}{2}\rfloor} \c_{2s}(\sigma)$, $\c_{\textup{odd}}(\sigma):=\sum_{s=1}^{\lfloor\frac{l-1}{2}\rfloor} \c_{2s+1}(\sigma)$, and 
$$\c(\sigma):=\sum_{s=2}^l \c_s(\sigma)=\c_{\textup{even}}(\sigma)+\c_{\textup{odd}}(\sigma).$$ 

For $r\in\mathbb N^{\ast}$ and $\sigma\in A_l$ if $s$ is odd and $\sigma\in S_l$ if $s$ is even, let 
$$\nu_{s,r}(\sigma)\in\mathbb N$$ 
be the smallest such that we have a product decomposition $\sigma=\prod_{i=1}^{\nu_{s,r}(\sigma)} \theta_i$ with each $\theta_i\in S_l$ of order $s$ and a product of at most $r$ disjoint $s$-cycles; we have $\n(\theta_i)\le sr$ and if $s$ is odd, then $\theta_i\in A_l$. This makes sense even for $l=s=4$ as the $4$-cycles generate $S_4$. 

Similarly, for $(r,s)\in (2\mathbb N^{\ast})^2$ and $\sigma\in A_l$, let 
$$\nu^{\textup{even}}_{s,r}(\sigma)\in\mathbb N$$ be the smallest such that we have a product decomposition $\sigma=\prod_{i=1}^{\nu_{s,r}(\sigma)} \theta_i$ with each $\theta_i\in A_l$ a product of $2n_i$ disjoint $s$-cycles with $n_i\in\llbracket0,\frac{r}{2}\rrbracket$; we have $\n(\theta_i)\le sr$.

For $r\in\mathbb N^{\ast}\setminus\{1\}$ and $\sigma\in A_l$ with $s$ odd and $l\ge sr$, let 
$$\nu^+_{s,r}(\sigma)\in\mathbb N$$
be the smallest such that we have a product decomposition $\sigma=\prod_{i=1}^{\nu^+_{s,r}(\sigma)} \theta_i$ with each $\theta_i\in S_l$ of order $s$ and a product of exactly $r$ disjoint $s$-cycles; we have $\n(\theta_i)=sr$.

\phantomsection{For a conjugacy class $\mathcal C$ of $A_l$ or $S_l$, let $o(\mathcal C):=o(\sigma)$, $\n(\mathcal C)=\n(\sigma)$, $\c_s(\mathcal C):=\c_s(\sigma)$, $\c_{\textup{even}}(\mathcal C):=\c_{\textup{even}}(\sigma)$, $\c_{\textup{odd}}(\mathcal C):=\c_{\textup{odd}}(\sigma)$, and $\c(\mathcal C):=\c(\sigma)$ for any $\sigma\in\mathcal C$; also, if $\nu_{s,r}(\sigma)$ is defined let $\nu_{s,r}(\mathcal C):=\nu_{s,r}(\sigma)$ and if $\nu_{s,r}^+(\sigma)$ is defined let $\nu^+_{s,r}(\mathcal C):=\nu^+_{s,r}(\sigma)$.}\label{EXT2}

Similarly, if $A$ is a finite non-empty set, we use the above notation for each permutation $\sigma_A\in\perm(A)$ via its image under an isomorphism $\perm(A)\cong S_{|A|}$ defined by a bijection $A\rightarrow \llbracket1,|A|\rrbracket$; so $\n(\sigma_A)$ is the cardinality of the support $\supp(\sigma_A)$ of $\sigma_A$, i.e., is the number of elements of $A$ that are not fixed by $\sigma_A$, etc.
\end{notation}

If $s$ is odd, the existence of $\nu^+_{s,r}(\sigma)$ also follows from the existence of $\nu_{s,r}(\sigma)$ and the following lemma.

\begin{lemma}\label{F1}
Let $(l,s,r)\in (\mathbb N^{\ast}\setminus\{1\})^3$ and $t\in\llbracket1,r-1\rrbracket$. Suppose that $s$ is odd and $l\ge sr$. Let $\sigma\in A_l$ be a product of $t$ disjoint $s$-cycles. Then we have a product decomposition $\sigma=\sigma_1\sigma_2$ with $\sigma_1$ and $\sigma_2$ as products of $r$ disjoint $s$-cycles. In particular, $\nu^+_{s,r}(\sigma)\le 2$. 
\end{lemma}

\begin{proof}
We write $\sigma=\prod_{i=1}^t\varsigma_i$ as a product of $t$ disjoint $s$-cycles. As $l\ge sr$, there exists $s$-cycles $\varsigma_{t+1},\ldots,\varsigma_r\in S_l$ such that we have $\supp(\varsigma_i)\cap\supp(\varsigma_j)=\emptyset$ for each $(i,j)\in\llbracket1,r\rrbracket^2$ with $i<j$. The fact holds as we can take $\sigma_1:=\prod_{i=1}^t\varsigma_i^2\prod_{i=t+1}^r\varsigma_i$ and $\sigma_2:=\prod_{i=1}^t\varsigma_i^{-1}\prod_{i=t+1}^r\varsigma_i^{-1}$.
\end{proof}

\begin{lemma}\label{P1}
For each $l\in\mathbb N^{\ast}$ and $\sigma\in S_l$ the following properties hold.

\medskip
{\bf (1)} We have identities $\nu_{2,1}(\sigma)=\n(\sigma)-\c(\sigma)=\sum_{s=2}^l (s-1)\c_s(\sigma)$.

\smallskip
{\bf (2)}  We have inequalities $c_{\textup{even}}(\sigma)+2\c_{\textup{odd}}(\sigma)\le\nu_{2,1}(\sigma)\le l-1$.
\end{lemma}

\begin{proof}
The first identity of part (1) is well-known and the second one follows from the first one and the identity $\n(\sigma)=\sum_{s=2}^l s\c_s(\sigma)$. So part (1) holds.

From part (1) we get that the expression $\nu_{2,1}(\sigma)-c_{\textup{even}}(\sigma)-2\c_{\textup{odd}}(\sigma)$ is equal to $\sum_{t=1}^{\lfloor\frac{l}{2}\rfloor} (2t-2)\c_{2t}(\sigma)+\sum_{t=1}^{\lfloor\frac{l-1}{2}\rfloor} (2t-2)\c_{2t+1}(\sigma)$ and hence it is greater than or equal to $0$. From this and the fact that the inequality $\nu_{2,1}(\sigma)\le l-1$ is well-known, we get that part (2) holds.
\end{proof}

\begin{proposition}\label{P2}
For $l\in\mathbb N^{\ast}$ and $\sigma\in S_l$ the following properties hold.

\medskip
{\bf (1)  (Brenner)}\footnote{Proposer's solution to E 1118 in `The American Mathematical Monthly', Vol.\ {\bf 62}, No.\ 1 (Jan., 1955), p.\ 43.} If $l\ge 3$ and $\sigma$ is an $s$-cycle with $s\in\llbracket2,l\rrbracket$, then we have product decompositions $\sigma=\sigma_1\sigma_2=\varsigma_1\varsigma_2$ into elements of order at most $2$ for which the identities $\nu_{2,1}(\sigma_1)=\nu_{2,1}(\varsigma_2)=s-1-\lfloor\frac{s}{2}\rfloor$ and $\nu_{2,1}(\sigma_2)=\nu_{2,1}(\varsigma_1)=\lfloor\frac{s}{2}\rfloor$ hold. 

\smallskip
{\bf (2)} If $l\ge 4$ (resp.\ $l\ge 1$), then each $\sigma\in S_l$ is a product of two elements in $S_l$ of order $2$ (resp.\ of order at most $2$). 

\smallskip
{\bf (3) (Brenner and Riddell)} If $l\ge 4$ and $\sigma\in A_l$ is conjugate in $A_l$ to its inverse, then $\sigma$ is a product of two permutations in $A_l$ of order $2$.
\end{proposition}

\begin{proof}
For part (1) we can assume that $\sigma=(1\,\cdots\,s)$. So we have $\sigma=\sigma_1\sigma_2=\varsigma_1\varsigma_2$ with $\sigma_1:=\prod_{i=1}^{\lfloor\frac{s}{2}\rfloor} (s+1-i\, i+1)$, $\sigma_2=\varsigma_1:=\prod_{i=1}^{\lfloor\frac{s}{2}\rfloor} (s+1-i\,i)$, and $\varsigma_2:=\sigma_2\sigma_1\sigma_2$.\footnote{For instance, cf.\ the two colors in \url{https://i.sstatic.net/NbLub.png}.}

If $s$ is odd, then we have $\nu_{2,1}(\sigma_1)=\nu_{2,1}(\sigma_2)={\lfloor\frac{s}{2}\rfloor}$. If $s$ is even, then we have $\nu_{2,1}(\sigma_1)+1=\nu_{2,1}(\sigma_2)=\frac{s}{2}$. Thus $\nu_{2,1}(\sigma)=\nu_{2,1}(\sigma_1)+\nu_{2,1}(\sigma_2)=s-1$. So, as $\nu_{2,1}(\varsigma_1)=\nu_{2,1}(\sigma_2)$ and $\nu_{2,1}(\varsigma_2)=\nu_{2,1}(\sigma_1)$, part (2) holds.

Part (2) follows from part (1).

Part (3) is mentioned in \cite{BR2}, property (i), by citing the reference to part (1), but part (1) does not imply part (3) directly. So we include a detailed proof of part (3). We first consider five disjoint special cases as follows.

{\bf Case 1: $\sigma$ is a cycle and $4\mid \n(\sigma)-1$.} In this case part (3) holds by part (1). 

\phantomsection{{\bf Case 2: $\sigma=\vartheta_1\vartheta_2$ is a product of two disjoint cycles of even length.} By writing $\vartheta_i=\varsigma_{i,1}\varsigma_{i,2}$ as a product of two permutations of order at most $2$ with $\varsigma_{i,1}$ odd and $\supp(\vartheta_i)=\supp(\varsigma_{i,1})\cup\supp(\varsigma_{i,2})$ for each $i\in\{1,2\}$ by parts (1) and (2), with $\sigma_i:=\varsigma_{1,i}\varsigma_{2,i}$ for $i\in\{1,2\}$ we have $(\sigma_1,\sigma_2)\in A_l^2$ and $\sigma=\sigma_1\sigma_2$. We have $o(\sigma_1)=o(\sigma_2)=2$ except when $\vartheta_1$ and $\vartheta_2$ are transpositions. If $\vartheta_1$ and $\vartheta_2$ are transpositions, then $\sigma$ is a product of two elements of order $2$ as $A_4$ has a normal elementary abelian subgroup of order $4$.}\label{PH14i}

{\bf Case 3: $\sigma=\vartheta_0\vartheta_1\vartheta_2$ is a product of disjoint cycles with $4\mid \n(\vartheta_0)-3$ and both $\n(\vartheta_1)$ and $\n(\vartheta_2)$ even.} For $i\in\{1,2\}$ we write $\vartheta_i=\varsigma_{i,1}\varsigma_{i,2}$ as in the prior paragraph. We write $\vartheta_0=\varsigma_{0,1}\varsigma_{0,2}$ as a product of two permutations of order $2$ with $\nu_{2,1}(\varsigma_{0,1})=\nu_{2,1}(\varsigma_{0,1})=\frac{\n(\varsigma_0)-1}{2}$ odd by part (1). For $\sigma_1:=\varsigma_{0,1}\varsigma_{1,1}\varsigma_{2,2}$ and $\sigma_2:=\varsigma_{0,2}\varsigma_{1,2}\varsigma_{2,1}$ we have $(\sigma_1,\sigma_2)\in A_l^2$, $o(\sigma_1)=o(\sigma_2)=2$, and $\sigma=\sigma_1\sigma_2$.

{\bf Case 4: we have $\n(\sigma)\le l-2$, $\sigma$ is a cycle, and $4\mid\n(\sigma)-3$.} With the notation of the prior paragraph, we can assume that $\sigma=\vartheta_0$. As $\n(\vartheta_0)\le l-2$, let $i$ and $j$ be two distinct elements of $\llbracket1,l\rrbracket$ fixed by $\vartheta_0$. For $\sigma_1:=\varsigma_{0,1}(i\; j)$ and $\sigma_2:=\varsigma_{0,2}(i\;j)$ we have $(\sigma_1,\sigma_2)\in A_l^2$, $o(\sigma_1)=o(\sigma_2)=2$, and $\sigma=\vartheta_0=\sigma_1\sigma_2$.

{\bf Case 5: $\sigma=\vartheta_0\vartheta_1\vartheta_2$ is a product of disjoint cycles with $4\mid \n(\vartheta_0)-3$ and $\n(\vartheta_1)=\n(\vartheta_2)=2s+1$ with $s\in\mathbb N^{\ast}$.} We write $\vartheta_1=(i_1\;\cdots\;i_{2s+1})$ and $\break\vartheta_2=(j_1\;\cdots\;j_{2s+1})$. For $\vartheta_{0,3}:=\prod_{t=1}^{2s+1} (i_t\;j_t)$ and $\vartheta_{0,4}:=\prod_{t=1}^{2s+1} (i_t\;i_{t+1})$ with $i_{2t+2}:=i_{1}$ we have $\vartheta_1\vartheta_2=\vartheta_{0,3}\vartheta_{0,4}$. With $\vartheta_0=\vartheta_{0,1}\vartheta_{0,2}$ as above, for $\sigma_1:=\varsigma_{0,1}\varsigma_{0,3}$ and $\sigma_2:=\varsigma_{0,2}\varsigma_{0,4}$ we have $(\sigma_1,\sigma_2)\in A_l^2$, $o(\sigma_1)=o(\sigma_2)=2$, and $\sigma=\sigma_1\sigma_2$.

As $[(1\; 2)(3\;4)][(1\; 2)(3\;4)]$ is the identity permutation,  to prove part (3) in general, we can assume that $\c(\sigma)>0$. As $\sigma\in A_l$, we have $c_{\textup{even}}\in 2\mathbb N$. As $\sigma$ and $\sigma^{-1}$ are conjugate in $A_l$, one of the following properties holds: (i) $\n(\sigma)\le l-2$ or (ii) $\sum_{i=0}^{\lfloor\frac{l-3}{4}\rfloor} \c_{4i+3}\in 2\mathbb N^{\ast}$ and $\c_{\textup{even}}=0$ or (iii) $\c_{\textup{even}}\in 2\mathbb N^{\ast}$ or (iv) there exists $s\in 2\mathbb N^{\ast}+1$ with $\c_s(\sigma)\ge 2$. Based on this, we can write $\sigma$ as a product of disjoint permutations as in Cases 1 to 5 above. So part (3) holds by Cases 1 to 5.
\end{proof}

\begin{lemma}\label{P3}
For each $l\in\mathbb N^{\ast}$ and $\sigma\in A_l$ the following properties hold.

\medskip
{\bf (1)} We have identities $\nu_{3,1}(\sigma)=\frac{\n(\sigma)-\c_{\textup{odd}}(\sigma)}{2}=\frac{\nu_{2,1}(\sigma)+\c_{\textup{even}}(\sigma)}{2}.$

\smallskip
{\bf (2)} We have inequalities $\c(\sigma)\le\nu_{3,1}(\sigma)\le\min\left(\nu_{2,1}(\sigma),\lfloor\frac{l}{2}\rfloor\right)$.

\smallskip
{\bf (3)} The equality $\nu_{3,1}(\sigma)=\nu_{2,1}(\sigma)$ holds iff $\nu_{2,1}(\sigma)=\c_{\textup{even}}(\sigma)$ and iff $\sigma$ is a product of an even number of disjoint transpositions.
\end{lemma}

\begin{proof}
For the first identity of part (1) see \cite{HR1}, Cor.\ 2.4(i).\footnote{For Richard Stong's proof of this identity, see `The American Mathematical Monthly', Vol.\ {\bf 110}, No.\ 2 (Feb., 2003), p.\ 162.} The second identity of part (1) follows from the first identity of Lemma \ref{P1}(1). So part (1) holds.

The first inequality of part (2) follows from part (1) and Lemma \ref{P1}(2). As $\n(\sigma)\le l$, we have $\nu_{3,1}(\sigma)\le\bigl\lfloor\frac{l-\c_{\textup{odd}}(\sigma)}{2}\bigr\rfloor\le\lfloor\frac{l}{2}\rfloor$. As $c_{\textup{even}}(\sigma)\le\nu_{2,1}(\sigma)$ by Lemma \ref{P1}(2), from part (1) we get that $\nu_{3,1}(\sigma)\le\nu_{2,1}(\sigma)$. So part (2) holds.

The first equivalence of part (3) follows from part (1). As $\nu_{2,1}(\sigma)-\c_{\textup{even}}(\sigma)$ is equal to $\sum_{t=1}^{\lfloor\frac{l}{2}\rfloor} (2t-2)\c_{2t}(\sigma)+\sum_{t=1}^{\lfloor\frac{l-1}{2}\rfloor} 2t\c_{2t+1}(\sigma)$ by Lemma \ref{P1}(1), we have $\nu_{2,1}(\sigma)=\c_{\textup{even}}(\sigma)$ iff $\c_s(\sigma)=0$ for each $s\in\llbracket3,l\rrbracket$ and iff $\sigma$ is a product of disjoint transpositions, their number being even as $\sigma\in A_l$. So part (3) holds.\end{proof}

\begin{example}\normalfont\label{EX0}
To exemplify Lemma \ref{P3}(1), let $(s,t)\in (\mathbb N^{\ast})^2$. If $2t+1\le l$, then 
$$(1\,\cdots\,2t+1)=(1\,2\,3)(3\,4\,5)\cdots (2t-1\,2t\,2t+1)$$ 
is a product of $t$ permutations that are $3$-cycles. Similarly, if $2t+2s\le l$, then $(1\,\cdots\,2t)(2t+1\,\cdots\,2t+2s)$ is equal to
$$(1\,2\,\cdots\,2t-1)(2t-1\,2t+1\,2t+2)(2t-1\,2t\,2t+2)(2t+2\,2t+3\,\cdots\,2t+2s)$$
and hence is a product of $t-1+1+1+s-1=t+s$ permutations that are $3$-cycles.
\end{example}

\begin{theorem}\label{P4}
Let $(l,s)\in (\mathbb N^{\ast}\setminus\{1\})^2$ with $4\le s\le l$. Let $\epsilon\in\{0,1\}$ be such that $s\equiv\epsilon\pmod{2}$. Let $\sigma\in S_l$ if $\epsilon=0$ and $\sigma\in A_l$ if $\epsilon=1$. Then the following properties hold.

\medskip
{\bf (1) (Herzog and Reid)} We have lower and upper bounds
$$\frac{\n(\sigma)-\c(\sigma)}{s-1}\le\nu_{s,1}(\sigma)\le \frac{\n(\sigma)+\c(\sigma)+(2+\epsilon)s-6}{s-1}.$$

{\bf (2) (Herzog and Reid)} We have inequalities
$$\nu_{s,1}(\sigma)\in\begin{cases}\frac{3\n(\sigma)}{2s}+\frac{2s+15}{3}\quad\quad\, {\rm if}\; s\equiv 0\pmod{3}\\
\frac{3\n(\sigma)-2}{2(s-1)}+3\quad\quad\quad {\rm if}\; s\equiv 1\pmod{3}\\
\frac{3\n(\sigma)}{2s-1}+\frac{4s+9}{3}\quad\quad\;\,\, {\rm if}\; s\equiv 2\pmod{3}.\end{cases}$$

{\bf (3)} For $s=5$ we have inequalities
$$\frac{\n(\sigma)-\c_{\textup{odd}}(\sigma)}{4}\le\nu_{5,1}(\sigma)\le\min\Bigl(\frac{\n(\sigma)+29}{3},\frac{3\n(\sigma)+18}{8},\frac{\n(\sigma)+\c(\sigma)+9}{4}\Bigr).$$
\end{theorem}

\begin{proof}
See \cite{HR2}, Thm.\ 1 and its proof for part (1).

See \cite{HR2}, Props.\ 3 to 5 for part (2).

For part (3), the upper bound follows from parts (1) and (2) applied to $s=5$ and \cite{HR2}, Thm.\ 2 and the lower bound follows from Lemma \ref{P3}(1) and the fact that each $5$-cycle is a product of two $3$-cycles. So part (3) holds.
\end{proof}

For an exact formula for $\nu_{4,1}(\sigma)$ see \cite{HR1}, Cor.\ 2.4(ii).

For parts (1) to (3) of the next theorem see \cite{Ber}, Cor.\ 2.1, \cite{BH}, Thm.\ 2, and \cite{BH}, Thm.\ 3 (respectively).

\begin{theorem}\label{P5}
Let $s\in \llbracket2,l\rrbracket$. Then the following properties hold.

\medskip
{\bf (1) (Bertram)} Each $\sigma\in A_l$ is a product of two $s$-cycles iff either $s\in \llbracket\lfloor\frac{3l}{4}\rfloor,l\rrbracket$ or $(l,s)=(4,2)$. 

\smallskip
{\bf (2) (Bertram and Herzog)} Each $\sigma\in A_l$ is a product of three $s$-cycles iff either $s\in \llbracket\lceil\frac{l}{2}\rceil,l\rrbracket$ is odd or $(l,s)=(7,3)$. 

\smallskip
{\bf (3) (Bertram and Herzog)} Each $\sigma\in A_l$ is a product of four $s$-cycles iff either $s\in \llbracket\lceil\frac{3l}{8}\rceil,l\rrbracket$ with $l\not\equiv 1\pmod{8}$ or $s\in \llbracket\lfloor\frac{3l}{8}\rfloor,l\rrbracket$ with $l\equiv 1\pmod{8}$ or $(l,s)=(6,2)$. 
\end{theorem}

For other similar results and conjectures pertaining to five or more $s$-cycles see \cite{BH}, Conj.\ and \cite{HKL2}, Thms.\ 3.3 and 3.4 and Conjs.\ 1.1 and 1.2.

\begin{property}\normalfont\label{P6}
Suppose that $l\ge 3$. Let $s\in \llbracket3,l\rrbracket$ and $j\in \llbracket2,s\rrbracket$. For $r\in \llbracket1,\lfloor \frac{s}{j}\rfloor\rrbracket$, each $s$-cycle of $S_l$ is a product of $r$ disjoint $j$-cycles and an $(s-rj+r)$-cycle, e.g., $(1\,\cdots\,s)=(1\;\ldots\;j)(j+1\;\ldots\;2j)\cdots (rj-j+1\;\ldots\;rj)(j\;2j\;\cdots rj\;rj+1\;\cdots s)$.
\end{property}

\begin{property}\normalfont\label{P7}
Suppose that $l\ge 4$. Let $s\in \llbracket2,\lfloor\frac{l}{2}\rfloor\rrbracket$ and $t\in \llbracket1,s\rrbracket$. By denoting $r:=s-t\in \llbracket0,s-1\rrbracket$, we have $s+t=2t+r=2s-r$. Then the product of two disjoint $s$-cycles in $S_l$ is a product of two non-commuting $(s+t)$-cycles whose supports intersect into a set of cardinality $2t$, e.g., 
$$(1\;3\;\ldots\;2s-1)(2\;4\;\cdots\; 2s)=(2\; 4\cdots\; 2r\; 2r+1\;\cdots\;2s)(1\;3\;\cdots\;2r-1\;2r\;\cdots\;2s-1).$$ 
This property is also a particular case of \cite{HKL1}, Thm.\ 8 (applied to $l=2t$).
\end{property}

\begin{theorem}\label{P8}
For $l\in\mathbb N^{\ast}\setminus\{1,2\}$ and $\sigma\in A_l$ the following properties hold.

\medskip
{\bf (1) (Brenner and Herzog)} We have a product decomposition $\sigma=\theta_1\theta_2$ such that $o(\theta_1)=o(\theta_2)=3$. 

\smallskip
{\bf (2) (Brenner and Herzog)} If we have $l=4m+r$ with $(m,r)\in\mathbb N\times\{0,1,3\}$, then in part (1) we can assume moreover that $\n(\theta_1)=\n(\theta_2)=3m+3\lfloor\frac{r}{3}\rfloor$. 

\smallskip
{\bf (3)} Let the set $\mathcal B$ be either $\{0,3\}$ or $\{-3,0\}$. If $\n(\sigma)\ge 3$, then in part (1) we can assume moreover that $\max\bigl(\n(\theta_1),\n(\theta_2)\bigr)=3\lfloor\frac{\n(\sigma)+2}{4}\rfloor$ and $\n(\theta_1)-\n(\theta_2)\in\mathcal B$.
\end{theorem}

\begin{proof}
For part (1) see \cite{BR2}, property (iii).

For part (2), the case $l=3$ is trivial and thus we can assume that $l\ge 4$. For $l\ge 4$ see \cite{BR1}, Thms.\ 6.12 and 6.13. So part (2) holds.

For part (3), by taking inverses we can assume that $\mathcal B=\{0,3\}$. Based on part (2) and the fact that $\c_{\textup{even}}(\sigma)$ is even, we can assume that there exist $m\in\mathbb N^{\ast}$ and $(s_1,s_2)\in (\mathbb N^{\ast}\setminus\{1\})^2$ such that $\n(\sigma)=l=4m+2=2s_1+2s_2$ and we have a product decomposition $\sigma=\sigma_1\sigma_2$ into disjoint cycles of lengths $2s_1$ and $2s_2$. As $s_1+s_2=2m+1$ is odd and $\sigma_1\sigma_2=\sigma_2\sigma_1$ we can also assume that $s_1>s_2$. Based on Property \ref{P6} applied to $(j,r)=(3,1)$, we write $\sigma_1=\vartheta_0\sigma_3$ with $\vartheta_0$ a $3$-cycle, $\sigma_3$ a $(2s_1-2)$-cycle, $|\supp(\vartheta_0)\cap\supp(\sigma_3)|=1$, and $\supp(\sigma_1)=\supp(\vartheta_0)\cup\supp(\sigma_3)$. So for $\varsigma:=\sigma_3\sigma_2$ we have $\n(\varsigma)=l-2=4m$. Based of part (2), we write $\varsigma=\vartheta_1\theta_2$ with $o(\vartheta_1)=o(\theta_2)=3$ and $\n(\vartheta_1)=\n(\theta_2)=3m$; it follows that we have identities $\supp(\vartheta_1)\cup\supp(\theta_2)=\supp(\sigma_1)=4m$ and $|\supp(\vartheta_1)\cap\supp(\theta_2)|=2m$. So $2m$ elements that are not fixed by $\varsigma$ do not belong to $\supp(\vartheta_1)\cap\supp(\theta_2)$. Based on this and the circular property of the cycle $\varsigma$, we can assume that $\supp(\vartheta_0)\cap\supp(\vartheta_1)=\emptyset$, and hence by taking $\theta_1:=\vartheta_0\vartheta_1$ we have $\n(\theta_1)=3m+3=\n(\theta_2)+3$, $o(\theta_1)=3$, and $\sigma=\theta_1\theta_2$. So part (3) holds.
\end{proof}

For the next theorem see \cite{BE}, Thms.\ 1 and 2.

\begin{theorem}[{\bf Brenner and Evans}]
\label{P9}
If $l\ge 5$ and $l\notin\{8,9,14\}$, then each $\sigma\in A_l$ is a product of two permutations of order $5$.
\end{theorem}

\begin{property}\normalfont\label{P10}
If $(t,s)\in (\mathbb N^{\ast})^2$ and $l\ge 2t+s$, then we have an identity
$$(1\;2\;\cdots\; 2t+1)=(t+1\;t+2\;\ldots\; 2t+s)(2t+s\;2t+s-1\;\cdots\;2t+1\;1\;2\;\cdots\; t)$$
in $S_l$. In other words, each $(2t+1)$-cycle in $S_l$ is a product $\theta_1\theta_2$ with $\theta_1$ and $\theta_2$ as $(t+s)$-cycles in $S_l$ such that $|\supp(\theta_1)\cap\supp(\theta_2)|=s$.
\end{property}

\begin{theorem}\label{P11}
We consider integers $s\ge 4$, $r\ge 1$, $j\in \llbracket0,s-1\rrbracket$, and $l:=sr+j$. Then the following properties hold for $\sigma\in A_l$.

\medskip
{\bf (1) (Bardakov)} If $j=0$, then $\sigma$ is a product of two permutations in $S_l$ that are products of $r$ disjoint $s$-cycles.

\smallskip
{\bf (2)} Suppose that $j\in \llbracket1,s-1\rrbracket$. We define integers $\epsilon_1=\cdots=\epsilon_{\lfloor\frac{s}{2}\rfloor}:=1$ and $\epsilon_{\lfloor\frac{s}{2}\rfloor+1}=\cdots=\epsilon_{s-1}:=2$. Let $\jmath\in\{1,2\}$ be such that $\jmath=1$ iff $(s-1)\epsilon_j$ is even. Then $\sigma$ is a product of $2+\jmath\epsilon_j$ permutations in $S_l$, each one being an $s$-cycle or a product of $r$ disjoint $s$-cycles.
\end{theorem}

\begin{proof}
For part (1) see \cite{Bar}, Thm.

Part (2) follows from part (1) applied to $A_{sr}$ as there exist $s$-cycles $\theta_1,\ldots,\theta_{\jmath\epsilon_j}$ such that $\n(\sigma\prod_{i=1}^{\jmath\epsilon_j} \theta_i)\le sr=l-j$; we note that $\prod_{i=1}^{\jmath\epsilon_j} \theta_i\in A_l$ as $\jmath\epsilon_{\jmath}(s-1)$ is even. 

For instance, if $s$ is even and $\epsilon_{\jmath}$ is odd, then $\epsilon_{\jmath}=1$ and $\jmath=2$, and to check the statement we can assume that $j>0$ and $\n(\sigma)=sr+t$ with $t\in\llbracket1,j\rrbracket$. We first choose an $s$-cycle $\theta_1$ such that $\n(\sigma\theta_1)\le \n(\sigma)-\frac{s}{2}$ and then we choose an $s$-cycle $\theta_2$ such that $\n(\sigma\theta_1\theta_2)\le \n(\sigma\theta_1)-\frac{s}{2}$ if $\n(\sigma\theta_1)>sr$ and $\n(\sigma\theta_1\theta_2)\le sr$ if $\n(\sigma\theta_1)\le sr$. We have $\n(\sigma\theta_1\theta_2)\le\max\bigl(sr,\n(\sigma)-s\bigr)=sr$.
\end{proof}

Note that in Theorem \ref{P11}(2) we have $\jmath\epsilon_j\in\{1,2\}$ if $s$ is odd or $\epsilon_j$ is even and $\jmath\epsilon_j=2$ if $s(\epsilon_j+1)$ is even; in particular, $\jmath\epsilon_j=2$ if $s=4$.

\begin{theorem}\label{P12}
Suppose that $l\ge 5$. For $m\in \llbracket1,l\rrbracket$, we identity $S_m$ with the subgroup $\{\sigma\in S_l|\supp(\sigma)\subset \llbracket1,m\rrbracket\}$. Let $\mathcal D$ be a non-trivial conjugacy class of $A_l$. For each $t\in \llbracket\max\bigl(5-\n(\mathcal D),0\bigr),l-\n(\mathcal D)\rrbracket$, let $\mathcal D_t:=\mathcal D\cap A_{\n(\mathcal D)+t}$. Then the following properties hold.

\medskip
{\bf (1)} If $t\ge 2$, then $\mathcal D_t$ is a conjugacy class of $A_{\n(\mathcal D)+t}$. 

\smallskip
{\bf (2)} If $\mathcal D$ is also a conjugacy class of $S_l$, then $\mathcal D_t$ is a conjugacy class of $S_{\n(\mathcal D)+t}$. 

\smallskip
{\bf (3) (Dvir)} Let $\mu_t(\mathcal D)$ be the exponent of the conjugacy class $\mathcal D_t$ of $A_{\n(\mathcal D)+t}$, i.e., the smallest $\mu_t(\mathcal D)\in\mathbb N^{\ast}$ such that we have $A_{\n(\mathcal D)+t}=\mathcal D_t^{\mu_t(\mathcal D)}$. Then the inequality $\mu_t(\mathcal D)\le\max\bigl(3,\lfloor\frac{\n(\mathcal D)+t}{2}\rfloor\bigr)$ holds.

\smallskip
{\bf (4)} The number $\mu_t(\mathcal D)$ of part (4) does not depend on $l\ge\max\bigl(5,\n(\mathcal D)+t\bigr)$. 

\smallskip
{\bf (5) (Brenner)} If $\mathcal D$ is a conjugacy class of $S_l$ and $\nu_{2,1}(\mathcal D)\le\frac{\n(\mathcal D)+t-1}{2}$, then $\mu_t(\mathcal D)\le 4$.

\smallskip
{\bf (6) (Dvir)} If in part (5) we moreover have $o(\mathcal D)\ge 3$, then $\mu_t(\mathcal D)\le 3$.

\smallskip
{\bf (7)} Let $\mathcal U$ be a union of non-trivial conjugacy classes of $A_l$. We define $\break\n(\mathcal U):=\max\bigl(\n(\sigma)|\sigma\in\mathcal U\bigr)$. For each $t\in \llbracket\max\bigl(5-\n(\mathcal U),0\bigr),l-\n(\mathcal U)\rrbracket$ there exists a smallest 
$\mu_t(\mathcal U)\in \llbracket1,\lfloor\frac{\n(\mathcal U)+t}{2}\rfloor\rrbracket$
such that $A_{\n(\mathcal U)+t}=(\mathcal U\cap A_{\n(\mathcal U)+t}) ^{\mu_t(\mathcal U)}$.
\end{theorem}

\begin{proof}
Parts (1), (2), and (4) are clear.

For part (3) see \cite{Dv}, Thm.\ 9.1.

For part (5) see \cite{Bre}, Thm.\ 3.05.

For part (6) see \cite{Dv}, Thm.\ 10.2.

Part (7) follows from part (3) applied to a $\mathcal D$ contained in $\mathcal U$.
\end{proof}

\begin{lemma}\label{L1}
For $l\in\mathbb N^{\ast}\setminus\{1\}$, $s\in \llbracket1,l\rrbracket$, and $\sigma\in S_l$ the following properties hold.

\medskip
{\bf (1)} There exists $\theta\in S_l$ such that $\theta(i)=\sigma(i)$ for every $i\in \llbracket1,s\rrbracket$ and we have an inequality $\nu_{2,1}(\theta)\le\min(s,l-1)$. In particular, $\n(\theta)\le \min(2s,l)$. 

\smallskip
{\bf (2)} Suppose that $s\le l-2$ (so $l\ge 3$) or $\sigma\in A_l$. Then there exists $\vartheta\in A_l$ such that $\vartheta(i)=\sigma(i)$ for every $i\in \llbracket1,s\rrbracket$, $\nu_{2,1}(\vartheta)\le 2\lfloor\frac{s+1}{2}\rfloor$, for $s=1$ we have $\n(\vartheta)\le 3$ and for $s\ge 2$ we have $\n(\vartheta)\le\min(2s,l)$, and $\nu_{3,1}(\vartheta)\le\min\bigl(s,\lfloor\frac{l}{2}\rfloor\bigr)$.

\smallskip
{\bf (3)} Suppose that $2\le s\le l-2$ (so $l\ge 4$) or $\sigma\in A_l$. Then there exists $\vartheta\in A_l$ such that $\vartheta(i)=\sigma(i)$ for every $i\in \llbracket1,s\rrbracket$, $\nu_{2,1}(\vartheta)\le s+1$, $\n(\vartheta)\le \min(2s,l)$, and $\nu_{3,1}(\vartheta)\le\min\bigl(2\lfloor\frac{s}{2}\rfloor,\lfloor\frac{l}{2}\rfloor\bigr)$.\end{lemma}

\begin{proof}
Let $j:=l-s$. As the case $j=0$ is trivial we can assume that $j\in\mathbb N^{\ast}$; so $\min(s,l-1)=s$. 

We prove part (1) by induction on $j\in \llbracket1,l-1\rrbracket$. If $j=1$, then $l-1=s$, $\theta=\sigma$, and $\nu_{2,1}(\theta)\in \llbracket1,s\rrbracket$ by Lemma \ref{P1}(2). So the base of the induction holds. For $j\in \llbracket2,l-1\rrbracket$, for the passage from $j-1$ to $j$ we consider two disjoint cases.

\phantomsection{{\bf Case 1: there exists $t\in \llbracket s+1,l\rrbracket$ such that $\sigma(t)\in \llbracket s+1,l\rrbracket$.} Then, by reindexing, we can assume that $t=l$. Let $\tau\in S_l$ be the identity element of $S_l$ if $\sigma(l)=l$ and let $\tau:=\bigl(l\; \sigma(l)\bigr)$ if $\sigma(l)\neq l$. Then $\sigma_1:=\sigma\tau\in S_l$ fixes $\sigma(l)$ and we have $\sigma_1(i)=\sigma(i)$ for each $i\in \llbracket1,s\rrbracket$. By the inductive assumption applied to $\sigma_1$ and $j-1=l-1-s$, it follows that there exists $\theta\in S_l$ such that it fixes $l$, we have $\theta(i)=\sigma(i)$ for each $i\in \llbracket1,s\rrbracket$, and $\nu_{2,1}(\theta)\le s$.}\label{PH14f}

\phantomsection{{\bf Case 2: we have $\sigma(\llbracket s+1,l\rrbracket)\subset \llbracket1,s\rrbracket$ (so $j\le s$).} For $t\in\sigma(\llbracket s+1,l\rrbracket)$, let $\iota_t\in\mathbb N^{\ast}$ be the smallest such that $\sigma^{\iota_t}(t)\in \llbracket s+1,l\rrbracket$. The rule $t\mapsto \sigma^{\iota_t}(t)$ defines a bijection $\hbar:\sigma(\llbracket s+1,l\rrbracket)\rightarrow \llbracket s+1,l\rrbracket$. Let $\theta\in S_l$ be such that $\theta(i)=\sigma(i)$ for $i\in \llbracket1,s\rrbracket$ and $\theta\bigl(\hbar(t)\bigr)=t$ for $t\in\sigma(\llbracket s+1,l\rrbracket)$. In the cyclic decomposition of $\theta$, the elements $t\in\sigma(\llbracket s+1,l\rrbracket)$ belong to distinct cycles of length at least $2$, hence $\c(\theta)\ge l-s=j$. As $\n(\theta)\le l$, we have $\nu_{2,1}(\theta)\le l-(l-s)=s$ by Lemma \ref{P1}(1).}\label{PH14e}

This ends the inductive step and the induction and hence the proof of part (1). 

If $\theta$ is even, then we have $\nu_{3,1}(\theta)\le \frac{\n(\theta)}{2}$ by Lemma \ref{P3}(1) and from this and part (1) we get that part (2) holds. Thus to prove part (2) we can assume that each $\theta$ of part (1) is odd. If $s\in\{l-1,l\}$, then $\sigma$ is even by our hypotheses and we have $\sigma=\theta$, a contradiction. Thus $s\le l-2$. As $\theta\notin A_l$, $\nu_{2,1}(\theta)$ is odd. As $\nu_{2,1}(\theta)\le s$, it follows that $\nu_{2,1}(\theta)+1\le 2\lfloor\frac{s+1}{2}\rfloor$. As $s\le l-2$, there exists a transposition $\tau\in S_l$ with support in $\llbracket1,l\rrbracket\setminus\{\sigma(i)|i\in \llbracket1,s\rrbracket\}$. Let $\vartheta:=\tau\theta\in A_l$. Then $\vartheta$ is a product of at most $\nu_{2,1}(\theta)+1$ transpositions, hence $\nu_{2,1}(\vartheta)\le \nu_{2,1}(\theta)+1\le 2\lfloor\frac{s+1}{2}\rfloor$. We have $\vartheta(i)=\sigma(i)$ for every $i\in \llbracket1,s\rrbracket$. If $l\le 2s$, then $\n(\vartheta)\le l\le 2s$. If $l\ge 2s+1=3$ then we can choose $\vartheta\in A_l$ such that $\n(\vartheta)\le 2s+1=3$. If $l\ge 2s+1\ge 5$, then by reindexing the elements of $\llbracket1,l\rrbracket$ we can replace $l$ by $2s$ and hence we can assume that $\n(\vartheta)\le 2s$. In all situations it follows that $\nu_{3,1}(\vartheta)\le\min\bigl(s,\lfloor\frac{l}{2}\rfloor\bigr)$ by Lemma \ref{P3}(2) applied to $l\in\{l,\n(\vartheta)\}$ and that for $s=1$ we have $\n(\vartheta)\le 3$ and for $s\ge 2$ we have $\n(\vartheta)\le\min(2s,l)$. So part (2) holds. 

As $\n(\vartheta)\le 2\nu_{2,1}(\vartheta)\le 2\nu_{2,1}(\theta)+2$, part (3) follows from the proof of part (2) and Lemma \ref{P3}(2) except when each $\theta$ of part (1) is odd with $\nu_{2,1}(\theta)=s$ and in part (2) we compulsory have $l\ge 2s$; so $s$ is also odd and we must have $\theta(\llbracket1,s\rrbracket)\subset \llbracket s+1,l\rrbracket$. Thus $s\ge 3$. Taking $\tau:=(1\;2)$ and $\vartheta:=\tau\theta\in A_l$, part (3) holds.
\end{proof}

The following examples list situations when Lemma \ref{L1} cannot be improved. 

\begin{example}\normalfont\label{EX1}
{\bf (1)} Let $(l,s)\in(\mathbb N^{\ast})^2$ with $l\ge s+1$. For $\sigma=(1\; 2\;\cdots\; s\;s+1)$, if $\theta\in S_l$ is such that $\theta(i)=\sigma(i)$ for every $i\in \llbracket1,s\rrbracket$, then $\theta^i(1)\neq 1$ for each $i\in \llbracket1,s\rrbracket$. So the cycle of $\theta$ containing $1$ has length at least $s+1$. Thus we have $\nu_{2,1}(\theta)\ge s=\min(s,l-1)= \nu_{2,1}(\sigma)$. 

\smallskip
{\bf (2)} Let $s\in \mathbb N^{\ast}$ and $l:=2s$. If $\sigma=\prod_{i=1}^s (i\; s+i)$, then for each $\theta\in S_l$ such that $\theta(i)=\sigma(i)$ for every $i\in \llbracket1,s\rrbracket$ we have $\n(\theta)=2s=\min(2s+1,l)=\min(2s,l)$, $\c(\theta)\in \llbracket1,s\rrbracket$, and $\c_{\textup{odd}}(\theta)=0$, hence $\nu_{3,1}(\theta)=s=\min\bigl(s,\lfloor\frac{l}{2}\rfloor\bigr)$ by Lemma \ref{P3}(1) and $\nu_{2,1}(\theta)=2s-\c(\theta)\ge s$. In particular, if $s$ is odd and $\theta$ is even, then we cannot have $\c(\theta)=s$ and hence $\nu_{2,1}(\theta)\ge s+1$.

\smallskip
{\bf (3)} Let $s\in 1+2\mathbb N^{\ast}$ and $l:=2s+1$. If $\sigma(\llbracket1,s\rrbracket)=\llbracket s+1,2s\rrbracket$, then for $\vartheta\in A_l$ such that $\vartheta(i)=\sigma(i)$ for every $i\in \llbracket1,s\rrbracket$ we have $\n(\vartheta)-\c_{\textup{odd}}(\vartheta)=2s=\min(2s,l)$ and $\c_{\textup{even}}(\sigma)\in \llbracket0,s-1\rrbracket$ is even. Therefore $\nu_{2,1}(\vartheta)=2s-\c_{\textup{even}}(\sigma)\ge s+1$ by Lemma \ref{P1}(1); the equality holds when $\c_{\textup{odd}}(\vartheta)=1$ and $\c_{\textup{even}}(\vartheta)=s-1$. Note that $\nu_{3,1}(\vartheta)=s$ by Lemma \ref{P3}(1).\end{example}

\begin{lemma}\label{L2}
Let $(l,s,t)\in (\mathbb N^{\ast})^3$ be such that $2\le s< t$ and $s+t-1\le l$. Let $\theta$ and $\vartheta$ be two $s$-cycles in $S_l$. Then $\theta\vartheta$ is a product of two $t$-cycles.
\end{lemma}

\begin{proof}
If $\supp(\theta)=\supp(\vartheta)$, then there exists distinct elements $r_1,\ldots,r_{t-s}$ in $\llbracket1,l\rrbracket\setminus \supp(\theta)$ and we take $i\in\supp(\theta)$ and $j:=r_{t-s}$. If $\supp(\theta)\neq\supp(\vartheta)$, let $(i,j)\in\llbracket1,l\rrbracket^2$ be such that $i\in\supp(\theta)\setminus\supp(\vartheta)$ and $j\in\supp(\vartheta)\setminus\supp(\theta)$; as $l-2s\ge t-s-1\ge 0$, it follows that there exists distinct elements $r_1,\ldots,r_{t-s-1}$ in $\llbracket1,l\rrbracket\setminus [\supp(\theta)\cup\supp(\vartheta)]$. Then $\theta\vartheta=[\theta(i\;r_1\;\cdots\;r_{t-s-1}\;j)][(i\;r_1\;\cdots\;r_{t-s-1}\;j)^{-1}\vartheta]$ is a product of two $t$-cycles.
\end{proof}

\begin{lemma}\label{F2}
Let $j\in\mathbb N^{\ast}$ and let the triple $(l,s,r)\in (\mathbb N^{\ast}\setminus\{1\})^3$ be such that $l\ge (2r-1)s$. Let the pair $(\theta_{-1},\theta_1)\in S_l^2$ be such that for $\iota\in\{-1,1\}$ we have a product decomposition $\theta_{\iota}=\prod_{i\in\llbracket1,j\rrbracket} \vartheta_{\iota,i}$ with each $\vartheta_{\iota,i}$ as a product of $t_i$ disjoint $s$-cycles for some $t_i\in\llbracket1,r\rrbracket$. Then $\theta_{-1}\theta_1\in A_l$ and $\nu^+_{s,r}(\theta_{-1}\theta_1)\le 2j$.
\end{lemma}

\begin{proof}
The notation $\prod_{i\in\llbracket1,j\rrbracket}$ is used instead of $\prod_{i=1}^j$ as the order in the product decompositions of $\theta_1$ and $\theta_2$ is irrelevant. Conjugating multiple times, it suffices to show that for $\sigma:=\prod_{i\in\llbracket1,j\rrbracket} \vartheta_{-1,i}\vartheta_{1,i}$ we have $\sigma\in A_l$ and $\nu^+_{s,r}(\sigma)\le 2j$. To check this we can assume that $j=1$ and $t_1\le r-1$. As $\sigma$ is a product of $2t_1$ permutations that are $s$-cycles, we have $\sigma\in A_l$. As $l\ge (2r-1)s\ge (r+t_1)s=t_1s+t_1s+(r-t_1)s$, there exists $\vartheta\in S_l$ which is a product of $r-t_1$ disjoint $s$-cycles and such that $\supp(\vartheta)\cap [\supp(\vartheta_{-1,1})\cup\supp(\vartheta_{1,1})]=\emptyset$. We have $\sigma=\sigma_{-1}\sigma_1$, where $\sigma_{\iota}=\vartheta_{\iota}\vartheta^{\iota}$ is a product of $r$ disjoint $s$-cycles. Thus $\nu^+_{s,r}(\sigma)\le 2$ and lemma holds.
\end{proof}

\section{Finite symmetric groups, part II: decomposition properties}\label{S3}

In this section we include results on decompositions into products of multiple disjoint $s$-cycles, with $s\in\mathbb N^{\ast}$. We begin with the small values $s=2$ and $s=3$ and then we continue with arbitrary $s\ge 4$. 

In what follows we use without any extra comment that for each $x\in\mathbb Z$ we have identities $x=\lfloor x\rfloor=\lceil x\rceil$ and for each $x\in\mathbb R\setminus\mathbb Z$ we have strict inequalities $\lfloor x\rfloor<x<\lceil x\rceil=\lfloor x\rfloor+1<x+1$. So, $x\le\lceil x\rceil\le \lfloor x\rfloor +1\le x+1$ for each $x\in\mathbb R$.

\begin{proposition}\label{PR1}
Let $(l,r)\in \mathbb N^{\ast}\times(\mathbb N^{\ast}\setminus\{1\})^2$ and $\sigma\in S_l$. Then the following properties hold.

\medskip
{\bf (1)} If $\sigma\in A_l$, then we have 
$$\nu_{2,r}(\sigma)\le 2\Bigl\lceil\frac{\nu_{2,1}(\sigma)}{2r}\Bigr\rceil\le 2\Bigl\lceil\frac{l-1}{2r}\Bigr\rceil.$$

{\bf (2)} If $\sigma\notin A_l$, then we have 
$$\nu_{2,r}(\sigma)\le\Bigl\lceil\frac{\nu_{2,1}(\sigma)+1}{2r}\Bigr\rceil+\Bigl\lceil\frac{\nu_{2,1}(\sigma)-1}{2r}\Bigr\rceil\le\Bigl\lceil\frac{l}{2r}\Bigr\rceil+\Bigl\lceil\frac{l-2}{2r}\Bigr\rceil.$$
\end{proposition}

\begin{proof}
We write $\sigma=\prod_{i=1}^{\c_{\textup{odd}}(\sigma)} \theta_i\prod_{j=\c_{\textup{odd}}(\sigma)+1}^{\c(\sigma)} \vartheta_j$ into a product of disjoint cycles, each $\theta_i$ being a cycle of odd length and each $\vartheta_j$ being a cycle of even length. Based on Proposition \ref{P2}(1) applied to each such cycle, for every $i\in \llbracket1,\c_{\textup{odd}}(\sigma)\rrbracket$ we write $\theta_i=\sigma_{1,i}\sigma_{2,i}$ with $o(\sigma_{1,i})=o(\sigma_{2,i})=2$ and $\nu_{2,1}(\sigma_{1,i})=\nu_{2,1}(\sigma_{2,i})=\frac{\nu_{2,1}(\theta_i)}{2}$ and for every $j\in \llbracket\c_{\textup{odd}}(\sigma)+1,\c(\sigma)\rrbracket$ we write $\vartheta_j=\sigma_{1,j}\sigma_{2,j}$ with $o(\sigma_{1,j})=o(\sigma_{2,j})=2$ if $\n(\vartheta_j)\ge 4$ and with $\{o(\sigma_{1,j}),o(\sigma_{2,j})\}=\{1,2\}$ if $\n(\vartheta_j)=2$ and such that for $j$ odd we have $\nu_{2,1}(\sigma_{1,j})=\nu_{2,1}(\sigma_{2,j})+1=\frac{\nu_{2,1}(\vartheta_j)+1}{2}$ and for $j$ even we have $\nu_{2,1}(\sigma_{1,j})+1=\nu_{2,1}(\sigma_{2,j})=\frac{\nu_{2,1}(\vartheta_j)+1}{2}$. Thus $\sigma=\varsigma_1\varsigma_2$, where 
$\varsigma_1:=\prod_{i=1}^{\c(\sigma)} \sigma_{1,i}$ and $\varsigma_2:=\prod_{i=1}^{\c(\sigma)} \sigma_{2,i}$ have order at most $2$ and we have $\nu_{2,1}(\varsigma_1)=\nu_{2,1}(\varsigma_2)=\frac{\nu_{2,1}(\sigma)}{2}$ if $\c_{\textup{even}}(\sigma)$ is even and $\nu_{2,1}(\varsigma_1)=\nu_{2,1}(\varsigma_2)+1=\frac{\nu_{2,1}(\sigma)+1}{2}$ if $\c_{\textup{even}}(\sigma)$ is odd. As 
$$ \nu_{2,r}(\varsigma_1)+\nu_{2,r}(\varsigma_2)=\Bigl\lceil\frac{\nu_{2,1}(\varsigma_2)}{r}\Bigr\rceil+\Bigl\lceil\frac{\nu_{2,1}(\varsigma_2)}{r}\Bigr\rceil=\left\{
\begin{array}{l}
 2\bigl\lceil\frac{\nu_{2,1}(\sigma)}{2r}\bigr\rceil\;\quad\quad\quad\quad\quad\quad\;\,\; \textup{if}\;\sigma\in A_l\\
 \bigl\lceil\frac{\nu_{2,1}(\sigma)+1}{2r}\bigr\rceil+ \bigl\lceil\frac{\nu_{2,1}(\sigma)-1}{2r}\bigr\rceil\;\;\textup{if} \;\sigma\notin A_l,
\end{array}\right.$$
parts (1) and (2) follow from the inequality $\nu_{2,r}(\sigma)\le \nu_{2,r}(\varsigma_1)+\nu_{2,r}(\varsigma_2)$.
\end{proof}

\begin{example}\normalfont\label{EX1.5}
Suppose that $n\in\mathbb N^{\ast}\setminus\{1,2,3\}$, $l=2^n$, and $\sigma\in S_l$. 

\medskip
{\bf (1)} As $\lceil\frac{l}{2^{n-2}}\rceil\le 4$, we have $\nu_{2,2^{n-3}}(\sigma)\le 8$ by Proposition \ref{PR1}(1) and (2). 

\smallskip
{\bf (2)} Suppose that $\nu_{2,1}(\sigma)\ge l-2^{n-3}+1$. If $n_{2,2^{n-3}}(\sigma)\le 7$, then $\sigma$ can be written as a product of at most $7\cdot2^{n-3}$ transpositions, hence $\nu_{2,1}(\sigma)\le 7\cdot2^{n-3}=l-2^{n-3}$, a contradiction. From this and part (1) we get that $\nu_{2,2^{n-3}}(\sigma)=8$.
\end{example}

We have the following analog with a similar proof of Proposition \ref{PR1} for $3$-cycles.

\begin{proposition}\label{PR2}
For $(l,r)\in \mathbb N^{\ast}\times (\mathbb N^{\ast}\setminus\{1\})$ and $\sigma\in A_l$ the following properties hold.

\medskip
{\bf (1)} If $\nu_{3,1}(\sigma)$ is even, then we have 
$$\nu_{3,r}(\sigma)\le 2\Bigl\lceil\frac{\nu_{3,1}(\sigma)}{2r}\Bigr\rceil=2\Bigl\lceil\frac{\n(\sigma)-\c_{\textup{odd}}(\sigma)}{4r}\Bigr\rceil\le 2\Bigl\lceil\frac{4\lfloor\frac{l}{4}\rfloor}{4r}\Bigr\rceil=2\Bigl\lceil\frac{\lfloor\frac{l}{4}\rfloor}{r}\Bigr\rceil.$$

{\bf (2)} If $\nu_{3,1}(\sigma)$ is odd, then $\nu_{3,r}(\sigma)$ is less than or equal to
$$\Bigl\lceil\frac{\nu_{3,1}(\sigma)+1}{2r}\Bigr\rceil+\Bigl\lceil\frac{\nu_{3,1}(\sigma)-1}{2r}\Bigr\rceil=\Bigl\lceil\frac{\n(\sigma)-\c_{\textup{odd}}(\sigma)+2}{4r}\Bigr\rceil+\Bigl\lceil\frac{\n(\sigma)-\c_{\textup{odd}}(\sigma)-2}{4r}\Bigr\rceil.$$ 
In particular, we have $\nu_{3,r}(\sigma)\le\bigl\lceil\frac{\lfloor\frac{l+2}{4}\rfloor}{r}\bigr\rceil+\bigl\lceil\frac{\lfloor\frac{l-2}{4}\rfloor}{r}\bigr\rceil$.

\smallskip
{\bf (3)} If $l\ge 3r$, then $\nu^+_{3,r}(\sigma)-\nu_{3,r}(\sigma)\in\{0,1,2\}$.
\end{proposition}

\begin{proof}
For $i\in\{1,2,3,4\}$, we define $\c_i\in\mathbb N$ as follows. Let $\c_1:=\sum_{t=1}^{\lfloor\frac{l-1}{4}\rfloor} \c_{4t+i}(\sigma)$ and for $i\in\{2,3,4\}$ let $\c_i:=\sum_{t=0}^{\lfloor\frac{l-i}{4}\rfloor} \c_{4t+i}(\sigma)$. The sum $\c_2+\c_4=\c_{\textup{even}}(\sigma)$ is even. For $i\in\{2,4\}$, we write $\c_i=2s_i+\epsilon_i$ with $s_i\in\mathbb N$ and $\epsilon_i\in\{0,1\}$; thus $\epsilon_2=\epsilon_4$. Let $N:=\c_1+s_4+s_2+\c_3+\epsilon_2\in\mathbb N$. From the first identity of Lemma \ref{P3}(1) applied to cycles we get that $\c_1$, $s_2$, and $s_4$ contribute by even numbers to $\nu_{3,1}$, that each cycle that contributes to $\c_3$ contributes by an odd number to $\nu_{3,1}(\sigma)$, and that, if $\epsilon_2$ is $2$, then it contributes by an odd number to $\nu_{3,1}(\sigma)$. Hence $\nu_{3,1}(\sigma)$ is even iff $\c_3+\epsilon_2$ is even.

As parts (1) and (3) are clear if $\sigma$ is the identity permutation, we can assume that $N>0$.

We write $\sigma=\prod_{i=1}^N \theta_i$ as a product of disjoint permutations that are either cycles or are products of two disjoint cycles of even length. Based on Theorem \ref{P8}(2) and (3) and Lemma \ref{P3}(1), we can assume that the following properties hold.

\medskip\noindent
{\bf $\bullet$} For each $i\in \llbracket1,\c_1\rrbracket$, $\theta_i$ is a cycle of length in $1+4\mathbb N^{\ast}$ and therefore we can write $\theta_i=\vartheta_{1,i}\vartheta_{2,i}$ with $o(\vartheta_{1,i})=o(\vartheta_{2,i})=3$ and $\nu_{3,1}(\vartheta_{1,i})=\nu_{3,1}(\vartheta_{2,i})=\frac{\nu_{3,1}(\theta_i)}{2}$ also equal to $\frac{\n(\theta_i)-1}{4}$. 

\smallskip\noindent
{\bf $\bullet$} For each $i\in \llbracket\c_1+1,\c_1+s_4\rrbracket$, $\theta_i$ is a product of two disjoint cycles of length in $4\mathbb N^{\ast}$ and therefore we can write $\theta_i=\vartheta_{1,i}\vartheta_{2,i}$ with $o(\vartheta_{1,i})=o(\vartheta_{2,i})=3$ and $\nu_{3,1}(\vartheta_{1,i})=\nu_{3,1}(\vartheta_{2,i})=\frac{\nu_{3,1}(\theta_i)}{2}=\frac{\n(\theta_i)}{4}$. 

\smallskip\noindent
{\bf $\bullet$} For each $i\in \llbracket\c_1+s_4+1,\c_1+s_4+s_2\rrbracket$, $\theta_i$ is a product of two disjoint cycles of length in $2+4\mathbb N$ and therefore we can write $\theta_i=\vartheta_{1,i}\vartheta_{2,i}$ with $o(\vartheta_{1,i})=o(\vartheta_{2,i})=3$ and $\nu_{3,1}(\vartheta_{1,i})=\nu_{3,1}(\vartheta_{2,i})=\frac{\nu_{3,1}(\theta_i)}{2}=\frac{\n(\theta_i)}{4}$. 

\smallskip\noindent
{\bf $\bullet$} For each $i\in \llbracket\c_1+s_4+s_2+1,\c_1+s_4+s_2+\c_3\rrbracket$, $\theta_i$ is a cycle of length in $3+4\mathbb N$ and therefore we can write $\theta_i=\vartheta_{1,i}\vartheta_{2,i}$ with $\{o(\vartheta_{1,i}),o(\vartheta_{2,i}\}\subset\{1,3\}$ and the pair $\bigl(\nu_{3,1}(\vartheta_{1,i}),\nu_{3,1}(\vartheta_{2,i})\bigr)$ is $\bigl(\frac{\nu_{3,1}(\theta_i)+1}{2},\frac{\nu_{3,1}(\theta_i)-1}{2}\bigr)=\bigl(\frac{\n(\theta_i)+1}{4},\frac{\n(\theta_i)-3}{4}\bigr)$ if the expression $i-\c_1-s_4-s_2-\c_3$ is odd and is $\bigl(\frac{\nu_{3,1}(\theta_i)-1}{2},\frac{\nu_{3,1}(\theta_i)+1}{2}\bigr)=\bigl(\frac{\n(\theta_i)-3}{4},\frac{\n(\theta_i)+1}{4}\bigr)$ if $i-\c_1-s_4-s_2-\c_3$ is even. 

\smallskip\noindent
{\bf $\bullet$} If $\epsilon_2=1$, then $\theta_N$ is a product of two disjoint cycles of even lengths not congruent modulo $4$ and therefore, as we have $\n(\theta_N)\in 2+4\mathbb N^{\ast}$, we have a product decomposition $\theta_N=\vartheta_{1,N}\vartheta_{2,N}$ with $o(\vartheta_{1,N})=o(\vartheta_{2,N})=3$ and the pair $\bigl(\nu_{3,1}(\vartheta_{1,N}),\nu_{3,1}(\vartheta_{2,N})\bigr)$ is equal to $\bigl(\frac{\nu_{3,1}(\theta_N)+1}{2},\frac{\nu_{3,1}(\theta_N)-1}{2}\bigr)=\bigl(\frac{\n(\theta_i)+2}{4},\frac{\n(\theta_i)-2}{4}\bigr)$ if $\c_3$ is even and to $\bigl(\frac{\nu_{3,1}(\theta_N)-1}{2},\frac{\nu_{3,1}(\theta_N)+1}{2}\bigr)=\bigl(\frac{\n(\theta_i)-2}{4},\frac{\n(\theta_i)+2}{4}\bigr)$ if $\c_3$ is odd. 

\medskip
Thus $\sigma=\varsigma_1\varsigma_2$, where 
$\varsigma_1:=\prod_{i=1}^N \vartheta_{1,i}$ and $\varsigma_2:=\prod_{i=1}^N \vartheta_{2,i}$ have order $3$. Clearly, $\nu_{3,r}(\sigma)\le \nu_{3,r}(\varsigma_1)+\nu_{3,r}(\varsigma_2)$. As $\nu_{3,1}(\sigma)\equiv \c_3+\epsilon_2\pmod{2}$, we have identities $\nu_{3,1}(\varsigma_1)=\nu_{3,1}(\varsigma_2)=\frac{\nu_{3,1}(\sigma)}{2}$ if $\nu_{3,1}(\sigma)$ is even and $\nu_{3,1}(\varsigma_1)=\nu_{3,1}(\varsigma_2)+1=\frac{\nu_{3,1}(\sigma)+1}{2}$ if $\nu_{3,1}(\sigma)$ is odd. Hence the sum $\nu_{3,r}(\varsigma_1)+\nu_{3,r}(\varsigma_2)$ is equal to
$$\Bigl\lceil\frac{\nu_{3,1}(\varsigma_2)}{r}\Bigr\rceil+\Bigl\lceil\frac{\nu_{3,1}(\varsigma_2)}{r}\Bigr\rceil=\left\{
\begin{array}{l}
 2\lceil\frac{\nu_{3,1}(\sigma)}{2r}\rceil\quad\quad\quad\quad\quad\quad\quad\;\,\textup{if} \;\nu_{3,1}(\sigma)\;\;\textup{is even}\\
 \lceil\frac{\nu_{3,1}(\sigma)+1}{2r}\rceil+ \lceil\frac{\nu_{3,1}(\sigma)-1}{2r}\rceil\quad\,\textup{if} \;\nu_{3,1}(\sigma)\;\;\textup{is odd}.
\end{array}\right.$$
As $\nu_{3,1}(\sigma)\le \lfloor\frac{l}{2}\rfloor$ by Lemma \ref{P3}(2), if $\nu_{3,1}(\sigma)$ is even we have $\nu_{3,1}(\sigma)\le 2 \lfloor\frac{l}{4}\rfloor$ and if $\nu_{3,1}(\sigma)$ is odd we have $\nu_{3,1}(\sigma)+1\le 2 \lfloor\frac{l+2}{4}\rfloor$ and $\nu_{3,1}(\sigma)-1\le 2 \lfloor\frac{l-2}{4}\rfloor$. 

From the last paragraph we get that parts (1) and (2) hold.

For part (3), clearly we have $\nu_{3,r}(\sigma)\le\nu^+_{3,r}(\sigma)$. We already know that $\sigma$ is a product of at most $\nu_{3,r}(\sigma)$ permutations of order $3$ in such a way that: (i) there exists $j\in\{0,1,2\}$ such that precisely $\nu_{3,r}(\sigma)-j$ factors have support of cardinality $3r$, (ii) if $j=1$, one factor has support of cardinality $3<3r$, and (iii) if $j=2$, there exists $t\in\llbracket1,\ldots,r-1\rrbracket$ with the property that either (iii.a) two factors have support of cardinality $3t$ or (iii.b) we have $t\in\llbracket2,\ldots,r-1\rrbracket$ and one factor has support of cardinality $3t$ and the other one has support of cardinality $3t-3$. If $j=0$, then $\nu^+_{3,r}(\sigma)=\nu_{3,r}(\sigma)$. If $j=1$, then $\nu^+_{3,r}(\sigma)\le\nu_{3,r}(\sigma)-1+2=\nu_{3,r}(\sigma)+1$ by Lemma \ref{F1}. If $j=2$, then $\nu^+_{3,r}(\sigma)\le\nu_{3,r}(\sigma)-2+2+2=\nu_{3,r}(\sigma)+2$ by Lemma \ref{F1}.\footnote{If $6r\ge 3$ and $j=2$, then in the case (iii.a) we have $\nu^+_{3,r}(\sigma)\le\nu_{3,r}(\sigma)-2+2$ by Lemma \ref{F2} and hence $\nu^+_{3,r}(\sigma)=\nu_{3,r}(\sigma)$.} So part (3) holds.\end{proof}

\begin{definition}\label{D3}
Let $(l,r)\in \mathbb N^{\ast}\times (\mathbb N^{\ast}\setminus\{1\})$.\index{3-deviation}

\medskip
{\bf (1)} By the 3-deviation factor of $(r,l)$\index{3-deviation!3-deviation factor of $(r,l)$} we mean 
$$\digamma_{r,l}:=l-2r\Bigl[\max\Bigl(2\Bigl\lceil\frac{\lfloor\frac{l}{4}\rfloor}{r}\Bigr\rceil,\Bigl\lceil\frac{\lfloor\frac{l+2}{4}\rfloor}{r}\Bigr\rceil+\Bigl\lceil\frac{\lfloor\frac{l-2}{4}\rfloor}{r}\Bigr\rceil\Bigr)-2\Bigr].$$

{\bf (2)} By the 3-deviation factor of $\sigma\in A_l$\index{3-deviation!3-deviation factor of $\sigma$} we mean $\digamma_{r,\sigma}:=\n(\sigma)-2r\bigl(\nu_{3,r}(\sigma)-2\bigr)$.

\smallskip
{\bf (3)} By the 3-deviation exponent of $(r,l)$\index{3-deviation!3-deviation exponent of $(r,l)$} we mean 
$$\wp_{r,l}:=\Bigl\lfloor\frac{l-\digamma_{r,l}}{2r}\Bigr\rfloor+2.$$
\end{definition}

Thus $\digamma_{r,l}\in\mathbb Z$ is the largest such that we have an identity 
$$\max\Bigl(2\Bigl\lceil\frac{\lfloor\frac{l}{4}\rfloor}{r}\Bigr\rceil,\Bigl\lceil\frac{\lfloor\frac{l+2}{4}\rfloor}{r}\Bigr\rceil+\Bigl\lceil\frac{\lfloor\frac{l-2}{4}\rfloor}{r}\Bigr\rceil\Bigr)=\Bigl\lfloor\frac{l-\digamma_{r,l}}{2r}\Bigr\rfloor+2.$$
Similarly, $\digamma_{r,\sigma}\in\mathbb Z$ is the largest such that $\nu_{3,r}(\sigma)=\lfloor\frac{\n(\sigma)-\digamma_{r,\sigma}}{2r}\rfloor+2$.

\begin{proposition}\label{PR3}
Let $(l,r)\in \mathbb N^{\ast}\times (\mathbb N^{\ast}\setminus\{1\})$. We write first $l=4s+i$ with $(s,i)\in\mathbb N\times \{0,1,2,3\}$ and second $s=rt+j$ with $(t,j)\in\mathbb N\times \llbracket0,r-1\rrbracket$. Then the following properties hold. 

\medskip
{\bf (1)} We have 
$$2\Bigl\lceil\frac{\lfloor\frac{l}{4}\rfloor}{r}\Bigr\rceil-\Bigl\lceil\frac{\lfloor\frac{l+2}{4}\rfloor}{r}\Bigr\rceil-\Bigl\lceil\frac{\lfloor\frac{l-2}{4}\rfloor}{r}\Bigr\rceil=\begin{cases} \bigl\lceil\frac{\lfloor\frac{l}{4}\rfloor}{r}\bigr\rceil-\bigl\lceil\frac{\lfloor\frac{l}{4}\rfloor-1}{r}\bigr\rceil\in\{0,1\} \quad\;\, {\rm if}\; l\equiv 0\; \textup{or}\;1\pmod{4}\\
\bigl\lceil\frac{\lfloor\frac{l}{4}\rfloor}{r}\bigr\rceil-\bigl\lceil\frac{\lfloor\frac{l}{4}\rfloor+1}{r}\bigr\rceil\in\{-1,0\} \;\;\, {\rm if}\; l\equiv 2\; \textup{or}\;3\pmod{4}.
\end{cases}$$

{\bf (2)} We have
$$\digamma_{r,l}=\begin{cases}
2r+i\quad\quad {\rm if}\; (i,j)\in\{2,3\}\times\{0\}\\
4r+i+4j\quad\quad {\rm if}\; (i,j)\in\{0,1\}\times\{0\}\;\textup{or}\; j\in\llbracket1,r-1\rrbracket.\end{cases}$$

{\bf (3)} We have $\digamma_{r,l}\ge 2r+2\ge6$. Moreover, $\digamma_{r,l}=6$ iff $(r,i,j)=(2,2,0)$, $\digamma_{r,l}=7$ iff $(r,i,j)=(2,3,0)$, and $\digamma_{r,l}=8$ iff $(r,i,j)\in\{(3,2,0)\}\cup \bigl(\{(2,0)\}\times\llbracket1,r-1\rrbracket\Bigr)$. In particular, for $r\ge 3$ we have $\digamma_{r,l}\ge 2r+2\ge 8$.

\smallskip
{\bf (4)} If $\sigma\in A_l$, then $\digamma_{r,\n(\sigma)}\le \digamma_{r,\sigma}$ and $\nu_{3,r}(\sigma)\le\wp_{r,\n(\sigma)}\le\wp_{r,l}$. In particular, $\nu_{3,r}(\sigma)\le \wp_{r,\n(\sigma)}\le\lfloor\frac{\n(\sigma)-2r-2}{2r}\rfloor+2= \lfloor\frac{\n(\sigma)-2}{2r}\rfloor+1$.
\end{proposition}

\begin{proof}
Part (1) follows from the identity
\begin{equation}\label{EQ4}
2\Bigl\lceil\frac{\lfloor\frac{l}{4}\rfloor}{r}\Bigr\rceil-\Bigl\lceil\frac{\lfloor\frac{l+2}{4}\rfloor}{r}\Bigr\rceil-\Bigl\lceil\frac{\lfloor\frac{l-2}{4}\rfloor}{r}\Bigr\rceil=\begin{cases} \lceil\frac{j}{r}\rceil-\lceil\frac{j-1}{r}\rceil \quad\quad {\rm if}\; i\in\{0,1\}\\
\lceil\frac{j}{r}\rceil-\lceil\frac{j+1}{r}\rceil \quad\quad {\rm if}\; i\in\{2,3\}.\end{cases}
\end{equation}

For part (2), let $M:=\max\bigl(2\lceil\frac{\lfloor\frac{l}{4}\rfloor}{r}\rceil,\lceil\frac{\lfloor\frac{l+2}{4}\rfloor}{r}\rceil+\lceil\frac{\lfloor\frac{l-2}{4}\rfloor}{r}\rceil\bigr)$. We have $l=4rt+4j+i$. We compute $\digamma_{r,l}=4rt+4j+i-2r(M-2)$ using two disjoint cases as follows.

{\bf Case 1: $j=0$ and $i\in\{2,3\}$.} We have $M=\Bigl\lceil\frac{\lfloor\frac{l+2}{4}\rfloor}{r}\Bigr\rceil+\Bigl\lceil\frac{\lfloor\frac{l-2}{4}\rfloor}{r}\Bigr\rceil=1+2\lceil\frac{\lfloor\frac{l}{4}\rfloor}{r}\rceil$ by Equation (\ref{EQ4}). Thus $M-2=2t+1-2=2t-1$ and hence $\digamma_{r,l}=2r+i$.

{\bf Case 2: $i\in\{0,1\}$ or $j\in\llbracket1,r-1\rrbracket$.} We have $M=2\lceil\frac{\lfloor\frac{l}{4}\rfloor}{r}\rceil=2t$ by Equation (\ref{EQ4}). Thus $M-2=2t-2$ and hence $\digamma_{r,l}=4r+i+4j$.

From these three cases we get that part (2) holds.

Part (3) follows from part (2). 

For part (4), the first two inequalities follow from Proposition \ref{PR2} applied to $l=\n(\sigma)$ and the third inequality follows from definitions. If $r\ge 3$ (resp.\ $r=2$), then $\digamma_{\n(\sigma),l}\ge 8$ (resp.\ $\digamma_{\n(\sigma),l}\ge 6$) by part (3). Hence $\wp_{r,\n(\sigma)}\le \lfloor\frac{\n(\sigma)-8}{2r}\rfloor+2$ (resp.\ $\nu_{3,2}(\sigma)\le \lfloor\frac{\n(\sigma)-6}{4}\rfloor+2=\lfloor\frac{\n(\sigma)+2}{4}\rfloor$). From this and the second inequality we get that the fourth (resp.\ fifth) inequality also holds. So part (4) holds.\end{proof}

The following theorem draws conclusions based on the results and the methods of proofs of several parts of Section \ref{S2}. 

\begin{theorem}\label{T2}
Let $s\in\mathbb N^{\ast}\setminus\{1\}$. Let $l\in\mathbb N^{\ast}$ with $l\ge s$. Then the following properties hold for $\sigma\in A_l$ and $\varsigma\in S_l$.

 \medskip
{\bf (1)} If $\lfloor\frac{3l}{4}\rfloor\le s$, then for each $r\in\mathbb N^{\ast}$ we have $\nu_{s,r}(\sigma)\le 2$.

\smallskip
{\bf (2)} If $\lfloor\frac{l}{2}\rfloor\le s$ and $s$ is odd, then for each $r\in\mathbb N^{\ast}$ we have $\nu_{s,r}(\sigma)\le 3$.

\smallskip
{\bf (3)} If $\lfloor\frac{3l}{8}\rfloor\le s$, then for each $r\in\mathbb N^{\ast}$ we have $\nu_{s,r}(\sigma)\le 4$.

\smallskip
{\bf (4)} If $s=2$, then for each $r\in\mathbb N^{\ast}$ we have $\nu_{2,r}(\varsigma)\le 2\bigl\lceil\frac{\lfloor\frac{l}{2}\rfloor}{r}\bigr\rceil$.

\smallskip
{\bf (5)} If $s=2$, $\sigma\in A_l$ is conjugate to its inverse in $A_l$, and $r\in 2\llbracket1,\lfloor\frac{l}{4}\rfloor\rrbracket$, then $\nu_{2,r}^{\textup{even}}(\sigma)\le 2\bigl\lceil\frac{\lfloor\frac{l}{4}\bigr\rfloor}{r}\rceil$, i.e., $\sigma$ is a product of at most $2\bigl\lceil\frac{\lfloor\frac{l}{4}\rfloor}{r}\bigr\rceil$ permutations in $A_l$ of order at most $2$ and of support of cardinality at most $2r$.

\smallskip
{\bf (6)} Suppose that $l\ge s=4$. Let $t:=\lfloor\frac{l}{4}\rfloor\in\mathbb N^{\ast}$. Then for each $r\in\mathbb N^{\ast}$ we have $\nu_{4,r}(\sigma)\le 2\lceil\frac{t}{r}\rceil$ if 4$\mid l$ and $\nu_{4,r}(\sigma)\le 2\lceil\frac{t}{r}\rceil+2$ if $4\nmid l$. In particular, we have $\nu_{4,r}(\sigma)\le 2\lfloor\frac{\n(\sigma)}{4r}\rfloor+2$ if $4\mid\n(\sigma)$ and $\nu_{4,r}(\sigma)\le 2\lfloor\frac{\n(\sigma)}{4r}\rfloor+4$ if $4\nmid\n(\sigma)$.

\smallskip
{\bf (7)} Suppose that $l\ge s=5$. Let $t:=\lfloor\frac{l}{5}\rfloor\in\mathbb N^{\ast}$. Then for each $r\in\mathbb N^{\ast}$ we have $\nu_{5,r}(\sigma)\le 2\lceil\frac{t}{r}\rceil$ if $l\notin\{8,9,14\}$ and $\nu_{5,r}(\sigma)\le 2\lceil\frac{t}{r}\rceil+2$ if $l\in\{8,9,14\}$. In particular, we have $\nu_{5,r}(\sigma)\le 2\lfloor\frac{\n(\sigma)}{5r}\rfloor+2$ if $\n(\sigma)\ge 5$ and $\n(\sigma)\notin\{8,9,14\}$ and $\nu_{5,r}(\sigma)\le 2\lfloor\frac{\n(\sigma)}{5r}\rfloor+4$ if $\n(\sigma)\in\{8,9,14\}$.

\smallskip
{\bf (8)} If $s\ge 6$ and $l=st+j$ with $(t,j)\in\mathbb N^{\ast}\times\llbracket0,s-1\rrbracket$, then for each $r\in\mathbb N^{\ast}$ we have $\nu_{s,r}(\sigma)\le 2\lceil\frac{t}{r}\rceil+\jmath\epsilon_j$ where $\epsilon_0:=0$ and for $j\ge 1$ the product $\jmath\epsilon_j\in\{1,2\}$ is as in Theorem \ref{P11}(2). In particular, we have $\nu_{s,r}(\sigma)\le 2\lfloor\frac{\n(\sigma)}{sr}\rfloor+2+\jmath\epsilon_j\le 2\lfloor\frac{\n(\sigma)}{sr}\rfloor+4$.

 \smallskip
{\bf (9)} Suppose that $s\ge 4$. Then for each $r\in\mathbb N^{\ast}$ we have 
$$\nu_{s,r}(\sigma)\le 2\frac{\lfloor\frac{3\n(\sigma)}{4}\rfloor}{(s-1)r}+2+\frac{2(r-1)}{r}\Bigl\lceil\log_s\Bigl(\max\bigl(1,\bigl\lfloor\frac{3\n(\sigma)}{4}\bigr\rfloor-s+2\bigr)\Bigr)\Bigr\rceil.$$ 
In particular,
$$\nu_{s,1}(\sigma)\le 2\Bigl\lfloor\frac{\lfloor\frac{3\n(\sigma)}{4}\rfloor}{s-1}\Bigr\rfloor+2$$ 
and $\nu_{s,r}(\sigma)< \frac{3\n(\sigma)}{2(s-1)r}+2+2\bigl\lceil\log_s\bigl(\max\bigl(1,\frac{3\n(\sigma)}{4}-s+2\bigr)\bigr)\bigr\rceil$ for $r\ge 2$.

\smallskip
{\bf (10)} Suppose that $s\ge 5$ is odd and $l\ge (2r-1)s$. Then 
$$\nu^+_{s,r}(\sigma)\le 2\frac{\lfloor\frac{3\n(\sigma)}{4}\rfloor}{(s-1)r}+2+\frac{2(r-1)}{r}\Bigl\lceil\log_s\Bigl(\max\bigl(1,\bigl\lfloor\frac{3\n(\sigma)}{4}\bigr\rfloor-s+2\bigr)\Bigr)\Bigr\rceil.$$

\smallskip
{\bf (11)} Suppose that $s\ge 5$ is odd and $\sigma$ is a cycle of odd length (resp.\ is the disjoint product of a transposition and a cycle of even length). Then 
$$\nu_{s,r}(\sigma)\le \frac{\n(\sigma)}{(s-1)r}+1+\frac{(r-1)}{r}\Bigl\lceil\log_s\Bigl(\max\bigl(1,\n(\sigma)-s+2\bigr)\Bigr)\Bigr\rceil$$ 
(resp. $\nu_{s,r}(\sigma)\le \frac{\n(\sigma)}{(s-1)r}+2+\frac{(r-1)}{r}\bigl\lceil\log_s\bigl(\max\bigl(1,\n(\sigma)-s+2\bigr)\bigr)\bigr\rceil)$). If moreover we have $l\ge (2r-1)s$, then 
$$\nu^+_{s,r}(\sigma)\le \frac{\n(\sigma)}{(s-1)r}+2+\frac{2(r-1)}{r}\Bigl\lceil\log_s\Bigl(\max\bigl(1,\n(\sigma)-s+2\bigr)\Bigr)\Bigr\rceil$$ 
(resp. $\nu^+_{s,r}(\sigma)\le \frac{\n(\sigma)}{(s-1)r}+4+\frac{2(r-1)}{r}\bigl\lceil\log_s\bigl(\max\bigl(1,\n(\sigma)-s+2\bigr)\bigr)\bigr\rceil$).
\end{theorem}

\begin{proof}
Parts (1) to (3) follow from Theorem \ref{P5}(1) to (3) (respectively). 

Parts (4) and (5) follows from Proposition \ref{P2}(2) and \ref{P2}(3) (respectively). 

Parts (6) and (8) follow from Theorem \ref{P11}.

Part (7) follows from Theorem \ref{P9} if $l\notin\{8,9,14\}$ and from Theorem \ref{P11}(2) if $l\in\{8,9,14\}$. 

To prove parts (9) and (10) we follow the method of \cite{HR2}, Thm.\ 1. So we rely on Property \ref{P6} but only after applying Theorem \ref{P5}(1) so that the notation is simplified and we can next apply Lemma \ref{L2}. 

We can assume that $l=\n(\sigma)$ and we first define recursively two sequences $(l_i)_{i=0}^t$ and $(m_i)_{i=1}^t$ in $\mathbb N$ with $t\in\mathbb N$ as follows. Let $l_0:=\lfloor\frac{3l}{4}\rfloor$. Assuming that $(l_i)_{i=0}^j$ and $(m_i)_{i=1}^j$ have been defined with $j\in\mathbb N$, we consider two cases. If $l_j<s$, then we take $t=j$. If $l_j\ge s$, we write $l_j=sm_{j+1}+r_{j+1}$ with $(m_{j+1},r_{j+1})\in\mathbb N^{\ast}\times\llbracket0,s-1\rrbracket$ and define $l_{j+1}:=m_{j+1}+r_{j+1}$. 

From Theorem \ref{P5}(1) we get that we can write $\sigma=\theta_0\theta_2$ as a product of two $l_0$-cycles. If $t\ge 1$, then from Property \ref{P6} we get that for $\iota\in\{0,2\}$ we can write $\theta_{\iota}^{\iota-1}=\vartheta^+_{\iota,2}\vartheta_{\iota,1}$, where $\vartheta_{\iota,1}$ is a product of $m_1$ disjoint $s$-cycles and $\vartheta^+_{\iota,2}$ is an $l_1$-cycle. By induction on $j\in\{1,t\}$ we get that we can write 
$$\theta_{\iota}^{\iota-1}=\vartheta^+_{\iota,j+1}\prod_{i=j}^1\vartheta_{\iota,i}$$
with each $\vartheta_{\iota,i}$ as a product of $m_i$ disjoint $s$-cycles and with $\vartheta^+_{\iota,j+1}$ an $l_j$-cycle. 

Let $\vartheta^+:=(\vartheta^+_{0,t+1})^{-1}\vartheta^+_{2,t+1}$. We write $\vartheta^+=\vartheta_1\vartheta_2$ as a product of two $s$-cycles by Lemma \ref{L2} applied to $(s,t)$ equal to $(l_t,s)$. As 
$$\sigma=\Bigl(\prod_{i=1}^t\vartheta_{0,i}^{-1}\Bigr)\vartheta_1\vartheta_2\Bigl(\prod_{i=t}^1\vartheta_{2,i}\Bigr),$$
it follows that we have inequalities $\nu_{s,r}(\sigma)\le 2+2\sum_{i=1}^t \lceil\frac{m_i}{r}\rceil$ and, if $l\ge (2r-1)s$, $\nu^+_{s,r}(\sigma)\le 2+2\sum_{i=1}^t \lceil\frac{m_i}{r}\rceil$ by Lemma \ref{F2} applied to $j=1+\sum_{i=1}^t \lceil\frac{m_i}{r}\rceil$. 

By induction on $j\in\{1,\ldots,t\}$ one checks that we have 
$l_0=\sum_{i=1}^j (s-1)m_i+l_j$. Taking $j=t$ we get inequalities 
$$(s-1)\sum_{i=1}^t m_i\le l_0=l_t+(s-1)\sum_{i=1}^t m_i\le (s-1)\Bigl(1+\sum_{i=1}^t m_i\Bigr).$$ 
Thus
$\sum_{i=1}^t m_i\le \frac{l_0}{s-1}\le 1+\sum_{i=1}^t m_i$. Hence $\sum_{i=1}^t \frac{m_i}{r}\le \frac{l_0}{(s-1)r}\le \frac{1}{r}+\sum_{i=1}^t \frac{m_i}{r}$. As for $i\in\llbracket1,t\rrbracket$ we have $\lceil\frac{m_i}{r}\rceil\le \frac{r-1}{r}+\frac{m_i}{r}$, we get that $\sum_{i=1}^t\lceil\frac{m_i}{r}\rceil\le \frac{t(r-1)}{r}+\sum_{i=1}^t \frac{m_i}{r}$. It follows that 
\begin{equation}\label{EQ5}
\nu_{s,r}(\sigma)\le 2+2\frac{t(r-1)}{r}+2\sum_{i=1}^t \frac{m_i}{r}\le 2+2\frac{t(r-1)}{r}+2\frac{l_0}{(s-1)r}.
\end{equation}

Based on this, to end the proof of part (9) it suffices to show that we have an inequality $t\le t_0$, where $t_0:=\bigl\lceil\log_s\bigl(\max\bigl(1,l_0-s+2\bigr)\bigr)\bigr\rceil$. To check this we can assume $l_0\ge s$; so $t_0\ge 1$ and we have $l_0\le s^{t_0}+s-2$. 

By induction on $i\in\llbracket1,\min(t_0,t)\rrbracket$ we show that $l_i\le s^{t_0-i}+s-2$. The base of the induction holds. For the inductive step we assume that $i\in\llbracket1,\min(t_0,t)-1\rrbracket$ and $l_{i-1}\le s^{t_0-i+1}+s-2$. If $l_{i-1}\in\llbracket s^{t_0-i+1},s^{t_0-i+1}+s-2\rrbracket$, then $m_{i}=s^{t_0-i}$ and $r_i\le s-2$, hence $l_{i+1}=m_i+r_i\le s^{t_0-i}+s-2$. If 
$l_{i-1}<s^{t_0-i+1}-1$, then $m_{i}\le s^{t_0-i}-1$ and $r_i\le s-1$, hence $l_{i+1}=m_i+r_i\le s^{t_0-i}+s-2$. This ends the inductive step and thus the induction. 

As $l_{\min(t,t_0)}\le s^-+s-1=s-1<s$, we have $t=\min(t,t_0)$. So $t\le t_0$ and part (9) holds.

Part (10) is proved in the same way as part (9) as Inequalities (\ref{EQ5}) hold with $\nu_{s,r}(\sigma)$ replaced by $\nu^+_{s,r}(\sigma)$.

Part (11) is proved in the same way as parts (9) and (10), the main difference being: for a cycle $\sigma$ we do not have to perform the first step with two cycles of length at most $\lfloor\frac{3\n(\sigma)}{4}\rfloor$ and hence the term $2\frac{\lfloor\frac{3\n(\sigma)}{4}\rfloor}{(s-1)r}$ gets replaced by $\frac{\n(\sigma)}{(s-1)r}$. For $\nu^+$ we have to keep the factors $2$ as per Lemma \ref{F1}.
\end{proof}

\begin{corollary}\label{C4}
Let $(s,r)\in (\mathbb N^{\ast}\setminus\{1\})^2$ with $s$ odd. For $l\in\mathbb N$ with $l\ge (2r-1)s$ let 
$$\nu^+_{s,r}(l):=\max\bigl(\nu^+_{s,r}(\sigma)|\sigma\in A_l\bigr)$$
and let $\nu^+_{s,r}(l-\cycle)$ be
$$\max\bigl(\nu^+_{s,r}(\sigma)|\sigma\in A_l\,\textup{is a disjoint product of a transposition and a cycle}\bigr).$$
Let $N\in\mathbb N^{\ast}$. Then the following properties hold.

\medskip
{\bf (1)} If $s\ge 5$, then we have $\limsup \limits_{rN\to\infty} \frac{\nu^+_{s,r}(srN)}{srN}\le\frac{3s}{2(s-1)}$.

\smallskip
{\bf (2)} We have $\limsup \limits_{rN\to\infty} \frac{\nu^+_{3,r}(3rN)}{3rN}\le\frac{3}{2}$.

\smallskip
{\bf (3)} If $s\ge 5$, then we have $\limsup \limits_{rN\to\infty} \frac{\nu^+_{s,r}(srN-\cycle)}{srN}\le\frac{s}{s-1}$.
\end{corollary}

\begin{proof}
Part (1) follows from Theorem \ref{T2}(10). 

Part (2) follows from Proposition \ref{PR2}(1) to (3).

Part (3) follows from Theorem \ref{T2}(11). 
\end{proof}

The next lemma is used in the study of the $\pi_{2,|K|^2}(K)$s with $|K|\ge 5$ odd.

\begin{lemma}\label{L3}
Let $k\in 1+2\mathbb N^{\ast}$. Then the following properties hold for 
$$\nu_{k,1}(k^2):=\max\bigl(\nu_{k,1}(\sigma)\bigl|\sigma\in A_{k^2}\bigr).$$ 

{\bf (1)} We have $\nu_{3,1}(9)=4$.

\smallskip
{\bf (2)} If $k\ge 5$, then $\nu_{k,1}(k^2)\in\{\frac{3k+1}{2},\frac{3k+3}{2}\}$.

\smallskip
{\bf (3)} If $k\ge 9$ and $3\mid k$, then $\nu_{k,1}(k^2)=\frac{3k+1}{2}$.
\smallskip

\end{lemma}

\begin{proof}
Part (1) follows directly from Lemma \ref{P3}(1) and (2).

For part (2), let $\sigma\in A_{k^2}$. As $k$ is odd, from Theorem \ref{P5}(1) we get that we can write $\sigma=\theta_1\theta_2\theta_3$ as a product of three $\frac{k^2+1}{2}$-cycles. We write $\frac{k^2+1}{2}=k\frac{k-1}{2}+\frac{k+1}{2}$. As $\frac{k-1}{2}+\frac{k+1}{2}=k$, from Property \ref{P6} we get that for $i\in\{1,2,3\}$ each $\theta_i$ is a product of $\frac{k-1}{2}+1=\frac{k+1}{2}$ permutations that are $k$-cycles, the first $\frac{k-1}{2}$ permutations actually disjoint. Thus $\nu_{k,1}(\sigma)\le 3\frac{k+1}{2}=\frac{3k+3}{2}$. Hence $\nu_{k,1}(k^2)\le \frac{3k+3}{2}$.

We show that the assumption that $\nu_{k,1}(k^2)\le\frac{3k-1}{2}$ leads to a contradiction. From this assumption and \cite{HKL2}, Thm.\ 3.3 we get the relations 
$$k^2\le\Bigl\lfloor\frac{2k}{3}\frac{3k-1}{2}\Bigr\rfloor+1=\Bigl\lfloor\frac{3k^2-k}{3}\Bigr\rfloor+1=k^2+1+\Bigl\lfloor\frac{-k}{3}\Bigr\rfloor\le k^2+1-2=k^2-1,$$ a contradiction. Thus $\frac{3k+1}{2}\le\nu_{k,1}(k^2)$. 

Part (2) follows from the concluding inequalities of the prior two paragraphs.

Part (3) is a particular case of \cite{HKL2}, Thm.\ 3.4.\end{proof}

Lemma \ref{L3} refines a bit Theorem \ref{P4}(1) and (2) or \cite{HR2}, Thms.\ 1 and 2 in the particular case it considers.

For $s\in \llbracket4,l\rrbracket$ one can obtain upper bounds for the $\nu_{s,1}(\sigma)$s that are similar to the inequalities of Lemmas \ref{P1}(2) and \ref{P3}(2). To elaborate on this, we first introduce extra notation as follows.

\begin{notation}\normalfont\label{N2}
Let $s\in\mathbb N^{\ast}\setminus\{1,2,3\}$ and $r\in\mathbb N^{\ast}$. 

\medskip
{\bf (1)} Let $\overline{\nu}_{s,r}(\sigma):=\Bigl\lceil\frac{\sum_{i=s}^l c_i(\sigma)\lfloor\frac{i}{s}\rfloor}{r}\Bigr\rceil$.

\smallskip
{\bf (2)} Let $\underline{\c}=(\c_2,\ldots,c_{s-1})\in\mathbb N^{s-2}$. If $s$ is odd, then we assume that the sum $\underline{\c}_{
\textup{even}}:=\sum_{i\in\llbracket 1,\lfloor\frac{s-1}{2}\rfloor\rrbracket} \c_{2i}$ is even. For $l\in\mathbb N^{\ast}$ with $l\ge\sum_{i=2}^{s-1} i\c_i$, let $\nu^{(l)}_{s,r}(\underline{\c})$ be the smallest number such that we can write a permutation in $S_l$ which is a product of disjoint cycles of length in $\llbracket2,s-1\rrbracket$, the number of $i$-cycles being equal to $\c_i$ for each $i\in\llbracket2,s-1\rrbracket$, as a product of $\nu^{(l)}_{s,r}(\underline{\c})$ permutations that are products of at most $r$ disjoint $s$-cycles.

\smallskip
{\bf (3)} Let $l\in\mathbb N^{\ast}$ with $l\ge s$. Let $\sigma\in S_l$ if $s$ is even and $\sigma\in A_l$ if $s$ is odd. For $i\in\llbracket2,s-1\rrbracket$, let 
$\c_{\equiv i\pmod{s-1}}(\sigma):=\sum_{j=0}^{\lfloor\frac{l-i}{s-1}\rfloor}\c_{i+j(s-1)}(\sigma)$.

\smallskip
{\bf (4)} For $t\in\mathbb N$, let $\sigma_s^{(t)}$ in $S_l$ be a permutation, uniquely determined up to conjugation, defined recursively as follows. Let $\sigma_s^{(0)}:=\sigma$. We define $\sigma_s^{(1)}$ by the fact that for each $i\in\llbracket2,l\rrbracket$ we have
$$\c_i(\sigma_s^{(1)}):=\sum_{j\in\llbracket2,l\rrbracket;j-(s-1)\lfloor\frac{j}{s}\rfloor=i} \c_j(\sigma).$$
If $t\ge 2$, then $\sigma_s^{(t)}:=(\sigma_s^{(t-1)})_s^{(1)}$.
\end{notation}

As for $j\in\llbracket2,l\rrbracket$ we have $j-(s-1)\lfloor\frac{j}{s}\rfloor\equiv j\pmod{s-1}$, by induction on $t\in\mathbb N$ we get that for each $i\in\llbracket2,s-1\rrbracket$ we have
\begin{equation}\label{EQ6}
\c_{\equiv i\pmod{s-1}}(\sigma_s^{(t)})=\c_{\equiv i\pmod{s-1}}(\sigma).
\end{equation}

\begin{proposition}\label{PR4}
Let $s\in\mathbb N^{\ast}\setminus\{1,2,3\}$, $r\in \mathbb N^{\ast}\setminus\{1\}$, and $t\in\mathbb N^{\ast}$ with $l\ge s$. If $s$ is even let $\sigma\in S_l$ and if $s$ is odd let $\sigma\in A_l$. Then the following properties hold.

\medskip
{\bf (1)} We have an inequality 
$$\nu_{s,1}(\sigma)\le \nu^{(l)}_{s,1}\bigl(\c_{\equiv 2\pmod{s-1}}(\sigma),\ldots,\c_{\equiv s-1\pmod{s-1}}(\sigma)\bigr)+\sum_{i=s}^l \c_i(\sigma)\Bigl\lceil\frac{i-s+1}{s-1}\Bigr\rceil.$$

{\bf (2)} We have an inequality $\nu_{s,r}(\sigma)\le\nu_{s,r}\bigl(\sigma_s^{(1)}\bigr)+\overline{\nu}_{s,r}(\sigma)$.

\smallskip
{\bf (3)} There exists a smallest $N=N(\sigma,s)\in\mathbb N$ such that for each $t\in\mathbb N$ with $t\ge N$ and every $i\in\llbracket2,l\rrbracket$ we have $\c_i\bigl(\sigma_s^{(t)}\bigr)=0$ if $i\ge s$ and $\c_i\bigl(\sigma_s^{(t)}\bigr)=\c_{\equiv i\pmod{s-1}}(\sigma)$ if $i\in\llbracket2,s-1\rrbracket$. In particular, $\sigma_s^{(t)}$ and $\sigma_s^{(N)}$ are conjugate in $S_l$.

\smallskip
{\bf (4)} With $N$ as in part (3), we have 
$$\nu_{s,r}(\sigma)\le\nu^{(l)}_{s,r}\bigl(\c_{\equiv 2\pmod{s-1}}(\sigma),\ldots,\c_{\equiv s-1\pmod{s-1}}(\sigma)\bigr)+\sum_{i=1}^{N-1} \overline{\nu}_{s,r}\bigl(\sigma_s^{(i)}\bigr).$$
\end{proposition}

\begin{proof}
Part (1) follows from the fact that for each $i\in\mathbb N^{\ast}$ with $i\ge s$, by writing $i=t(s-1)+r$ with $(t,r)\in\mathbb N\times\llbracket 1,s-1\rrbracket$, so $t=\lceil\frac{i-s+1}{s-1}\rceil$, every $i$-cycle is a product of $t$ disjoint $s$-cycles and an $r$-cycle by the proof of Theorem \ref{T2}(9).

For part (2), for $i\in\mathbb N^{\ast}$ with $i\ge s$ we write $i=sn+m$ with $(n,m)\in\mathbb N\times\llbracket0,s-1\rrbracket$; every $i$-cycle is a product of $n$ disjoint $s$-cycles and an $(n+m)$-cycle by Property \ref{P6}. Performing this for each $i$-cycle with $i\ge s$, we get a product decomposition $\sigma=\theta\vartheta$ with $\theta$ a product of $\sum_{i=s}^l c_i(\sigma)\lfloor\frac{i}{s}\rfloor$ disjoint $s$-cycles and with $\vartheta$ in the conjugacy class of $\sigma_s^{(1)}$. From this part (2) follows.

Part (3) follows from the fact that for each pair $(i,j)\in\mathbb N^{\ast}\times\mathbb N$ with $i\ge s$, every $i$-cycle in the cyclic decomposition of $\sigma_s^{(j)}$ either does not contribute to the cyclic decomposition of $\sigma_s^{(j+1)}$ or it contributes with a cycle of smaller length, and this ensures that we have $N\in\llbracket0,l-s+1\rrbracket$.

For part (4), by adding the inequalities we get that from part (2) applied to $\sigma_s^{(t)}$ for each $t\in\llbracket0,N-1\rrbracket$, it follows that $\nu_{s,r}(\sigma)\le \nu_{s,r}\bigl(\sigma_s^{(N)}\bigr)+\sum_{i=1}^{N-1} \overline{\nu}_{s,r}\bigl(\sigma_s^{(i)}\bigr)$. But from part (3) and Equation (\ref{EQ6}) we get that we have the following identity $\nu_{s,r}\bigl(\sigma_s^{(N)}\bigr)=\nu^{(l)}_{s,r}\bigl(\c_{\equiv 2\pmod{s-1}}(\sigma),\ldots,\c_{\equiv s-1\pmod{s-1}}(\sigma)\bigr)$. So part (4) holds.
\end{proof}

\begin{example}\normalfont\label{EX2} {\bf (1)} Suppose that $s=4$. Let $(\c_2,\c_3,l)\in\mathbb N^2\times\mathbb N^{\ast}$ with $l\ge 2\c_2+3\c_3$. Then $\nu_{4,1}(1,0)=3$ and for $(\c_2,\c_3)\neq (1,0)$ with $\c_3$ even (resp.\ odd) we have $\nu_{4,1}(\c_2,c_3)=\c_2+\c_3$ (resp.\ $\nu_{4,1}(\c_2,c_3)=\c_2+\c_3+1$) by \cite{HR1}, Cor.\ 3.4(ii). 

\smallskip
{\bf (2)} Suppose that $s\ge 5$ is an even (resp.\ odd) integer. Let the $(s-1)$-tuple $(\c_2,\c_3,\ldots,\c_{s-1},l)\in\mathbb N^{s-2}\times\mathbb N^{\ast}$ with $l\ge \c:=\sum_{i=2}^{s-1}\c_i$ (resp.\ with $l\ge \c:=\sum_{i=2}^{s-1}\c_i$ and an even sum $\sum_{i=1}^{\lfloor\frac{s-1}{2}\rfloor}\c_{2i}$). If $\c=1$ and the only $i\in\llbracket2,s-1\rrbracket$ with $\c_i=1$ is even (resp.\ odd), then $\nu_{s,1}(\c_2,\c_3,\ldots,\c_{s-1})=3$ (resp.\ $\nu_{s,1}(\c_2,\c_3,\ldots,\c_{s-1})=2$) and in all other cases we have $\nu_{s,1}(\c_2,\c_3,\ldots,\c_{s-1})\le \c+1$ (resp.\ $\nu_{s,1}(\c_2,\c_3,\ldots,\c_{s-1})\le \c$) by \cite{HR1}, Thm.\ 2.2(iii) (resp.\ Thm.\ 2.2(ii)). 

\smallskip
{\bf (3)} Let $(\c_2,\c_3,\c_4,l)\in\mathbb N^3\times\mathbb N^{\ast}$ with $l\ge 2\c_2+3\c_3+4\c_4$ and $\c_2+\c_4$ even. Let $(t,r)\in\mathbb N\times\llbracket0,4\rrbracket$ be such that $2\c_2+3\c_3+4\c_4=5t+r$. Let $\sigma\in S_l$ be such that $\c_i(\sigma)=\c_i$ for each $i\in\{2,3,4\}$ and $\n(\sigma)=2\c_2+3\c_3+4\c_4$. If $r=0$, then $\nu_{5,1}(\c_2,\c_3,\c_4)\le 2t$ by Theorem \ref{P11}(1). If $r\in\{1,2\}$, then $\nu_{5,1}(\c_2,\c_3,\c_4)\le 2t+1$ by Theorem \ref{P11}(2) applied to $\sigma$ and $5$-cycles. If $t\in\{3,4\}$, then $\nu_{5,1}(\c_2,\c_3,\c_4)\le 2t+2$ by Theorem \ref{P11}(2). 
\end{example}

We have the following application of Proposition \ref{PR4}.

\begin{corollary}
Let $s\in\mathbb N^{\ast}\setminus\{1,2,3\}$ and $(n,l)\in (\mathbb N^{\ast})^2$ be such that $l\ge s^n$. Let $\sigma$ be an $s^n$-cycle. Then the following properties hold.

\medskip
{\bf (1)} We have $\nu_{s,1}(\sigma)=\frac{s^n-1}{s-1}$.

\smallskip
{\bf (2)} For $r\in\mathbb N^{\ast}\setminus\{1\}$ we have
$$\nu_{s,r}(\sigma)-\Bigl\lceil\frac{s^n-1}{r(s-1)}\Bigr\rceil\in\Bigl\llbracket0,\sum_{i=0}^{n-1}\Bigl\lceil\frac{s^i}{r}\Bigr\rceil-\Bigl\lceil\frac{s^n-1}{r(s-1)}\Bigr\rceil\Bigr\rrbracket\subset \llbracket0,n-1\rrbracket.$$
\end{corollary}

\begin{proof}
By induction on $t\in\llbracket1,n-1\rrbracket$ we get that $\sigma_s^{(t)}$ is an $s^{n-t}$-cycle. So $\sigma_s^{(n)}$ is the identity permutation. From Proposition \ref{PR4}(1) we get that $\nu_{s,1}(\sigma)\le\frac{s^n-1}{s-1}$. As an $s$-cycle is a product of $s-1$ transpositions, we get that $\nu_{s,1}(\sigma)\ge\frac{\nu_{2,1}(\sigma)}{s-1}=\frac{s^n-1}{s-1}$. Thus $\nu_{s,1}(\sigma)=\frac{s^n-1}{s-1}$. So part (1) holds.

Therefore $\lceil\frac{s^n-1}{r(s-1)}\rceil\le\nu_{s,r}(\sigma)$. From Proposition \ref{PR4}(4) we get the inequality $\nu_{s,r}(\sigma)\le\sum_{i=0}^{n-1}\lceil\frac{s^i}{r}\rceil$. So the belonging relation of part (2) holds. As the difference $\sum_{i=0}^{n-1}\Bigl\lceil\frac{s^i}{r}\Bigr\rceil-\Bigl\lceil\frac{s^n-1}{r(s-1)}\Bigr\rceil$ is a natural number at most equal to $\frac{n(r-1)}{r}$, it is an element of $\llbracket0,n-1\rrbracket$. So part (2) holds.
\end{proof}

\section{On representations attached to actions over finite fields}\label{S4}

For a set $A$, let $\End(A)$ be the monoid of functions from $A$ to $A$. For a field extension $K\rightarrow L$, let $\star_L$ be the pullback to $\Spec L$ of (morphisms of) schemes over $\Spec K$. Let $\char(K)$ be the characteristic of $K$.

\phantomsection{Let $X$ be a reduced scheme of finite type over $\Spec K$. If $M$ is a smooth monoid scheme over $\Spec K$ equipped with an action $\Theta: M\times_{\Spec K} X\rightarrow X$ on $X$, then for $h\in M(L)$ let $\Theta_h\in\End(X_L)$ be such that $\Theta_h(x)=\Theta(h,x)$. We obtain a monoid homomorphism}\label{PH15-}
$$\varrho_{\Theta_L}:M(L)\rightarrow\End\bigl(X(L)\bigr)$$
defined by the rule $\varrho_{\Theta_L}(h)\mapsto\Theta_h$. If $M$ is in fact a smooth group scheme over $\Spec K$, then we denote also by $\varrho_{\Theta_L}$ the factorization $\varrho_{\Theta_L}:M(L)\rightarrow\perm\bigl(X(L)\bigr).$

\phantomsection{In this and the next two paragraphs we assume that $K$ is perfect. Let $G$ be a smooth linear algebraic group over $\Spec K$. Let $G^0$ be the connected component of the identity element of $G(K)$. Let $R_u(G^0)$ be the unipotent radical of $G^0$.}\label{PH15} 

\phantomsection{Let $H:=G^0/R_u(G^0)$; it is a reductive group over $\Spec K$. Let $H^{\sc}$ be the simply connected semisimple group cover of the derived group $H^{\der}$ of $H$. Let $Z_H$ be the kernel of the central isogeny $H^{\sc}\rightarrow H^{\der}$; it is a finite flat group scheme over $\Spec K$ of multiplicative type. Let $o(H)$ be the order of $Z_H$.}\label{PH15j}

As $K$ is perfect, each connected smooth unipotent group over $\Spec K$, such as $R_u(G)$, is split (see \cite{Bor}, Ch.\ V, Prop.\ 15.5 (ii)). 

\phantomsection{In this paragraph we assume that $K$ is finite. So $H$ (or $H^{\der}$ or $H^{\sc}$) is quasi-split (see \cite{Bor}, Ch.\ V, Prop.\ 16.6); let $B^+$ be a Borel subgroup of $H$ and let $T$ be a maximal torus of it. Let $B^-$ be the Borel subgroup of $H$ that contains $T$ and is opposite to $B^+$; so $B^+\cap B^-=T$. For $\star\in\{+,-\}$, let $U^{\star}$ be the unipotent radical of $B^{\star}$; so $B^{\star}=U^{\star}\rtimes T$. The Whitehead group of $H^{\sc}$ is trivial, i.e., $H^{\sc}(K)$ is generated by subgroups $U(K)$ with $U$ a subgroup of $H^{\sc}$ isomorphic to $\mathbb G_{\a,K}$ by the Chevalley--Steinberg classification of simply connected semisimple groups over $\Spec K$ (e.g., see \cite{St}, Ch.\ 11, Lem.\ 64 and \cite{T}, Subsubsect.\ 1.1.2).}\label{PH16} 

\begin{theorem}\label{T3}
Suppose that $K$ is finite and the order of the finite \'etale group scheme $G/G^0$ over $\Spec K$ is odd. We also assume that one of the following conditions holds.

\medskip
{\bf (1)} Both $|K|$ and $o(H)$ are odd and $H$ is semisimple. 

\smallskip
{\bf (2)} We have $4\mid |K|$.

\medskip
Let $\Lambda:G\times_{\Spec K} X\rightarrow X$ be an action of $G$ on $X$. Then for each $h\in G(K)$ we have $\Lambda_h\in\Alt\bigl(X(K)\bigr)$, i.e., the homomorphism $\varrho_{\Lambda}:G(K)\rightarrow \perm\bigl(X(K)\bigr)$ factors through $\Alt\bigl(X(K)\bigr)$.
\end{theorem}

\begin{proof}
We can assume that $|X(K)|\ge 2$. Let $h\in G(K)$. To show that the permutation $\Lambda_h\in\perm\bigl(X(K)\bigr)$ is even we can perform the following operation.

\medskip
{\bf (O)} We replace $h$ by $h^s$ with $s\in \mathbb N$ odd.

\medskip
As $|G(K)/G^0(K)|$ is odd, by performing (O) we can assume that $h\in G^0(K)$. By replacing $G$ with $G^0$ and $\Lambda$ with its restriction to $G^0\times_{\Spec K} X$, we can assume that $G=G^0$ is connected; let $\varpi:G\rightarrow H$ be the quotient homomorphism. 

Assume condition (1) holds. As $R_u(G)$ is connected, we have a short exact sequence $1\rightarrow R_u(G)(K)\rightarrow G(K)\rightarrow H(K)\rightarrow 1$ by Lang's theorem. As the order of the group $H^1\bigl(\Gal(\overline{K}/K),Z_H(\overline{K})\bigr)$ is odd and we have an exact complex $1\rightarrow Z_H(K)\rightarrow H^{\sc}(K)\rightarrow H(K)\rightarrow H^1\bigl(\Gal(\overline{K}/K),Z_H(\overline{K})\bigr)$, by performing (O) we can assume that there exists $h^{\sc}\in H^{\sc}(K)$ such that the images of $h^{\sc}$ and $h$ in $H(K)$ coincide. So $h\in G(K)^{-}:=\varpi(K)^{-1}\bigl(\Im(H^{\sc}(K)\rightarrow H(K))\bigr)$. Note that $1\rightarrow R_u(G)(K)\rightarrow G(K)^{-}\rightarrow \Im\bigl(H^{\sc}(K)\rightarrow H(K)\bigr)\rightarrow 1$ is a short exact sequence. The orders of the groups $R_u(G)(K)$ and $\mathbb G_{\a,K}(K)$ are odd and the group $H^{\sc}(K)$, hence also $\Im\bigl(H^{\sc}(K)\rightarrow H(K)\bigr)$, being generated by subgroups isomorphic to the additive group $K$, is generated by elements of odd orders. From the last two sentences it follows that $G(K)^{-}$ is generated by elements of odd orders. Thus the image of $G(K)^{-}$ in $\perm\bigl(X(K)\bigr)/\Alt\bigl(X(K)\bigr)$ is trivial. Hence $\Lambda_h\in\Alt\bigl(X(K)\bigr)$.

Assume condition (2) holds. We have $H(K)=U^+(K)U^-(K)T(K)U(K)$ by \cite{BT}, Cor.\ 6.26. As $T$ lifts isomorphically to $G$, it suffices to treat the case when $G$ is either a torus or a connected unipotent group. 

If $G$ is a torus, then $G(K)$ has odd order as $\char(K)=2$, hence $\Lambda_h$ is even. 

Assume now that $G=G^0=R_u(G^0)$. It suffices to work with a scheme-theoretic $G$-orbit of a $K$-valued point of $X$; it is a reduced locally closed subscheme of $X$ of the form $G/U$ with $U$ a subgroup scheme of $G$.\footnote{The scheme-theoretic $G$-orbit is isomorphic to $G/U$ and hence, as $G$ is $K$-split, it is (isomorphic to) an affine space over $\Spec K$ by \cite{Ro}, Thm.\ 5 (see also \cite{Co}, Rmk.\ 3.10). Moreover, if $X$ is affine, then it is closed in $X$ (see \cite{Sp}, Prop.\ 2.4.14 or \cite{DG}, Ch.\ IV, Sect.\ 2, Subsect.\ 2, Cor.\ 2.7).} We can assume that $U\neq G$. The scheme-theoretic intersection $\cap_{x\in G(\overline{K})} xU_{\overline{K}}x^{-1}$ is a closed subgroup scheme of $G_{\overline{K}}$ which, by Galois descent, is the extension to $\overline{K}$ of a normal subgroup scheme $U_0$ of $G$. By replacing $G$ with $G/U_0$ we can assume that $U_0$ is the trivial subgroup scheme of $G$. Thus $G$ has a central subgroup $G_0\cong\mathbb G_{\a,K}$ that meets trivially the conjugates of $U$. Hence $G_0x\cong\mathbb G_{\a,K}$ for each $x\in X(K)$.

We show that $h$ defines an even permutation on each $\langle h\rangle G_0(K)$-orbit of a given $x\in X(K)$, where $\langle h\rangle$ is the subgroup of $G(K)$ generated by $h$. The elements of $\langle h\rangle$ that map $x$ to $G_0x\cong\mathbb G_{\a,K}$ are those mapping $G_0x$ to itself, so they form a monoid and hence a subgroup $\Gamma$ of $\langle h\rangle$. Elements of $\Gamma$ define automorphisms of $G_0x\cong\mathbb G_{\a,K}$ that commute with the translation action of $G_0\cong\mathbb G_{\a,K}$, so they are translations. This gives a cyclic group of translations, hence of order $1$ or $2$. Thus the action of $h$ on the $\langle h\rangle G_0(K)$-orbit of $x$ consists of $\langle h\rangle$-orbits permuted simply transitively by the group $G_0(K)/\Gamma$, which has order $|K|$ or $\frac{|K|}{2}$; as $4\mid |K|$, this order is even. So $\Lambda_h$ is a product of an even number of disjoint cycles of the same length, thus $\Lambda_h\in\Alt\bigl(X(K)\bigr)$.
\end{proof}

We have the following direct consequence of Theorem \ref{T3}.

\begin{corollary}\label{C5}
Suppose that $K$ is a finite field of odd characteristic and $X$ is affine and integral. Then $X$ is a flexible variety over $\Spec K$ iff $|X(K)|\le 1$.
\end{corollary}

\begin{proposition}\label{PR5}
Suppose that $K$ is finite. Let $M$ be a smooth monoid scheme over $\Spec K$ equipped with a faithful action $\Theta: M\times_{\Spec K} X\rightarrow X$. Let $M_1$ and $M_2$ be two smooth closed subschemes of $M$. If for each finite field $L$ that contains $K$ we have an inclusion $\varrho_{\Theta_L}\bigl(M_1(L)\bigr)\leqslant\varrho_{\Theta_L}\bigl(M_2(L)\bigr)$, then $M_1\leqslant M_2$.
\end{proposition}

\begin{proof}
Let $\overline{K}$ be an algebraic closure of $K$. As the action of $M$ on $X$ is faithful and $K$ is perfect (being finite), there exists a finite set $\{x_1,\ldots,x_j\}$ of $\overline{K}$-valued points of $X$ such that $\cap_{i=1}^j \Stab_{M_{\overline{K}}}(x_i)$ is the trivial submonoid scheme of $M_{\overline{K}}$\footnote{This only holds over perfect fields (see Example \ref{EX3}(1) and (2)). For imperfect fields, one needs to assume that $X$ is geometrically reduced.}, where, for $i\in \llbracket1,j\rrbracket$, $\Stab_{M_{\overline{K}}}(x_i)$ is the submonoid scheme of $M_{\overline{K}}$ that fixes $x_i$. Let $L_0$ be a finite field extension of $K$ such that $x_i\in X(L_0)$ for each $i\in \llbracket1,j\rrbracket$. This implies that $\varrho_{\Theta_{L}}$ is injective for each finite field $L$ that contains $L_0$, and hence we have an inclusion $M_1(L)\leqslant M_2(L)$ by hypothesis. Therefore $M_1(\overline{K})\leqslant M_2(\overline{K})$. Taking Zariski closures we get that $M_{1,\overline{K}}\leqslant M_{2,\overline{K}}$, hence $M_1\leqslant M_2$.\end{proof}

\begin{example}\normalfont\label{EX3}
Suppose that $K$ is an imperfect field of characteristic $p$; so there exists $\alpha\in K\setminus\{x^p|x\in K\}$. Let $L:=K(\sqrt[p]{\alpha})$. We also view $X_2:=\mathbb A^2_L$ and $X_1:=\mathbb A^1_L$ as schemes over $\Spec K$. 

\medskip
{\bf (1)} The rule $\bigl((\beta_1,\beta_2),(\gamma_1,\gamma_2)\bigr)\mapsto (\gamma_1+\beta_1+\sqrt[p]{\alpha}\beta_2,\gamma_2+\beta_2^p)$ defines an action $\Lambda_{X_2}:\mathbb G_{\a,K}^2\times_{\Spec K} X_2\rightarrow X_2$ over $\Spec K$. Note that this action comes from an action $\Lambda_{X_2}^{\prime}:\mathbb G_{\a,L}^2\times_{\Spec L} X_2\rightarrow X_2$ over $\Spec L$. The kernel $\Ker$ of $\Lambda_{X_2}^{\prime}$ (i.e., of the morphism $\mathbb G_{\a,L}^2\rightarrow\underline{Aut}_{\Spec L}(X_2)$ of functors from the category of schemes over $\Spec L$ to the category of groups that defines $\Lambda_{X_2}^{\prime}$) is a closed subgroup scheme of $\mathbb G_{\a,L}^2$ by \cite{SGA3-1}, Exp.\ VI${}_B$, Exs.\ 6.2.4(c). We have $\Ker\cong\pmb{\alpha}_{p,L}$, but $\Ker$ is not the pullback of a subgroup scheme of $\mathbb G_{\a,K}^2$. Hence $\Lambda$ is faithful but $\Ker_{\overline{K}}$ is equal to $\Stab_{G_{a,\overline{K}}}(x)$ for each $x\in X_2(\overline{K})$. 

\smallskip
{\bf (2)} Similarly, the rule $\bigl((\beta_1,\beta_2),\gamma_1\bigr)\mapsto \gamma_1+\beta_1+\sqrt[p]{\alpha}\beta_2$ defines a faithful action $\Lambda_{X_1}:\mathbb G_{\a,K}^2\times_{\Spec K} X_1\rightarrow X_1$ over $\Spec K$ with stabilizers of $\overline{K}$-valued points of $X_1$ isomorphic to $\mathbb G_{\a,\overline{K}}$.
\end{example}

\section{Affine type of invariants of finite subsets}\label{S5}

We begin by recalling affine spans and by introducing several invariants of finite subsets of $K^n$. For $Y\subset K^n$, let $\Span(Y)$ be the $K$-linear span of $Y$.

\begin{definition}\label{D4}
Let $K$ be a field and $(n,m)\in (\mathbb N^{\ast})^2$ be such that $m\le |K|^n$. Let $Y\subset K^n$ be such that $|Y|=m$. We write $Y=\{P_1,\ldots,P_m\}$.

\medskip
{\bf (1)} Let 
$$\langle Y\rangle_{\aff}:=\Bigl\{\sum_{i=1}^m \alpha_iP_i\Bigl|(\alpha_1,\ldots,\alpha_m)\in K^m,\sum_{i=1}^m\alpha_i=1\Bigr\}$$ 
be the affine span of $Y$ over $K$. By the affine dimension\index{affine dimension} of $Y$ we mean 
$$\d_Y:=\dim(\langle Y\rangle_{\aff})\in \llbracket0,\min(m-1,n)\rrbracket.$$

{\bf (2)} We say that the set $Y$ is or the points $P_1,\ldots,P_m$ are {\it affinely independent}\index{affinely independent} if one of the following equivalent conditions holds. 

\medskip\noindent
{\bf (2.a)} We have $\d_Y=m-1$. 

\smallskip\noindent
{\bf (2.b)} The $m$ vectors $(1,P_1),\ldots,(1,P_m)$ are linearly independent in $K^{n+1}$.

\smallskip\noindent
{\bf (2.c)} The $m-1$ vectors $P_2-P_1,\ldots,P_m-P_1$ are linearly independent in $K^n$.

\smallskip
{\bf (3)} If $Z$ is a subset of $Y$ which does not contain the zero vector of $K^n$, then by the big linearly independent partition number\index{linearly independent partition number!big linearly independent partition number} of $Z$ we mean the smallest $\L_Z\in\mathbb N$ such that $Z$ can be partitioned into $\L_Z$ subsets that are linearly independent. 

\smallskip
{\bf (4)} By the linearly independent partition number\index{linearly independent partition number} of $Y$ we mean
$$\l_Y:=\max\bigl(\L_{-P+(Y\setminus\{P\})}|P\in Y\bigr)\in\mathbb N.$$

{\bf (5)} By the collinear number $\c_Y\in\mathbb N^{\ast}$ of $Y$\index{collinear number} we mean the largest $|Z|$ with $Z$ a subset of $Y$ of collinear points.
\end{definition}

\begin{remark}\normalfont\label{R11.5}
Let $Y$ be a finite non-empty subset of $K^n$ which does not contain the zero vector. In the more general context of matroids, the invariant $\L_Y$ has been first studied in \cite{E}. For instance, \cite{E}, Thm.\ 1 proves the formula
\begin{equation}\label{EQ6.9}
\L_Y=\max\Bigl(\Bigl\lceil\frac{|Z|}{\rank(Z)}\Bigl\rceil\Bigl|Z\subset Y,Z\neq\emptyset\Bigr),
\end{equation}
where $\rank(Z)\in\{d_Z,\d_Z+1\}$ is the rank of $Z$, i.e., is the maximum number of linearly independent points in $Z$.
\end{remark}

Clearly, for each non-empty subset $Y$ of $K^n$ we have inequalities
\begin{equation}\label{EQ7}
\c_Y-1\le\l_Y\le |Y|-\d_Y.
\end{equation}

\begin{definition}\label{D5}
Let the field $K$ and the pair $(n,m)\in (\mathbb N^{\ast})^2$ be such that we have $\sqrt[n]{m}\le |K|$. Let $Y\subset K^n$ be such that $|Y|=m$.

\medskip
{\bf (1)} 
For each $a\in\AGL_n(K)$ and every $i\in \llbracket1,n\rrbracket$, let $m_{i,Y,a}\in \llbracket1,\min(m,|K|)\rrbracket$ be the cardinality of the image $Y_i$ of $a(Y)$ under the projection $\pi_i:\mathbb A^n_K\rightarrow\mathbb A^1_K$ on the $i$-th coordinate, and let 
$$\s_a(Y):=\sum_{i=1}^n (m_{i,Y,a}-1)\in \bigl\llbracket0,n\bigl(\min(m,|K|)-1\bigr)\bigr\rrbracket.$$ 
By the capacity\index{capacity} of Y we mean
$$\s_Y:=\min\bigl(\s_a(Y)|a\in\AGL_n(K)\bigr)\in\bigl\llbracket0,n\bigl(\min(m,|K|)-1\bigr)\bigr\rrbracket.$$

{\bf (2)} Let $P=(\alpha_1,\ldots,\alpha_n)\in Y$. For each $a\in\AGL_n(K)$ and every $i\in \llbracket1,n\rrbracket$, let 
$$\overline{m}_{i,Y,P,a}\in \llbracket0,m_{i,Y,a}\rrbracket\subset \llbracket0,\min(m,|K|)-1\rrbracket$$ be the cardinality of the subset
$$Z_i:=\{\pi_i\bigl(a(Q)\bigr)|Q\in Y,\pi_j\bigl(a(Q)\bigr)=\pi_j\bigl(a(P)\bigr)\;\forall j\in \llbracket i+1,n\rrbracket\}$$
of $Y_i$. Let 
$$\overline{\s}^P_a(Y):=\sum_{i=1}^n (\overline{m}_{i,Y,a}-1)\in \llbracket0,\s_a(Y)\rrbracket\subset \bigl\llbracket0,n\bigl(\min(m,|K|)-1\bigr)\bigr\rrbracket.$$ 
By the strict capacity\index{capacity} of Y at $P$\index{capacity!strict capacity at a point} we mean
$$\overline{\s}^P_Y:=\min\bigl(\overline{\s}^P_a(Y)|a\in\AGL_n(K)\bigr)\in \llbracket0,\s_Y\rrbracket\subset \bigl\llbracket0,n\bigl(\min(m,|K|)-1\bigr)\bigr\rrbracket.$$

{\bf (3)} By the strict capacity\index{capacity!strict capacity} of $Y$ we mean
$$\overline{\s}_Y:=\max(\overline{\s}^P_Y|P\in Y)\in \llbracket0,\s_Y\rrbracket\subset\bigl\llbracket0,n\bigl(\min(m,|K|)-1\bigr)\bigr\rrbracket.$$

{\bf (4)} Let $P\in Y$. By the secant number of $Y$ at $P$\index{secant number!secant number at a point} we mean the smallest $\q_Y^P\in\mathbb N$ such that there exists a line that passes through $P$ and intersects $Y\setminus\{P\}$ in $\q_Y^P$ points.

\smallskip
{\bf (5)} By the secant number of $Y$\index{secant number} we mean 
$$\q_Y:=\max(\q^P_Y|P\in Y)\in\bigl\llbracket0,\bigl(\min(m,|K|)-1\bigr)\bigr\rrbracket.$$

\smallskip
{\bf (6)} Let $P\in Y$. By the set of minimal $n$-tuple secants of $Y$ at $P$\index{set of minimal $n$-tuple secants} we mean the subset $\mathcal Q_Y^P\subset \llbracket0,|K|-1\rrbracket^n$ defined recursively as follows. If $n=1$, then $\mathcal Q_Y^P:=\{\q_Y^P\}$. If $n\ge 2$, then
$$\mathcal Q_Y^P:=\{\q_Y^P\}\times\cup_{\pi\in \Pi_Y^P} \mathcal Q_{\pi(Y)}^{\pi(P)}$$
where $\Pi_Y^P$ is the set of all linear projections $\pi:\mathbb A^n_K\rightarrow\mathbb A^{n-1}_K$ for which we have an identity $\bigl|[\pi^{-1}\bigl(\pi(P)\bigr)]\cap (Y\setminus\{P\})\bigr|=\q_Y^P$.

\smallskip
{\bf (7)} By the multisecant number of $Y$ at $P$\index{secant number!multisecant number at a point} we mean 
$$\mq^P_Y:=\min\Bigl(\sum_{i=1}^n q_i\Bigl|(q_1,\ldots,q_n)\in\mathcal Q_Y^P\Bigr)\in\bigl\llbracket0,n\bigl(\min(m,|K|)-1\bigr)\bigr\rrbracket.$$

{\bf (8)} By the multisecant number of $Y$\index{secant number!multisecant number} we mean 
$$\mq_Y:=\max(\mq^P_Y|P\in Y)\in\bigl\llbracket0,n\bigl(\min(m,|K|)-1\bigr)\bigr\rrbracket.$$
\end{definition}

In Definition \ref{D5}(1) we have $\s_a(Y)=0$ iff $m=1$ and iff $\s(Y)=0$.

For each $P\in Y$ we have $\overline{m}_{n,Y,P,a}=m_{n,Y,a}$. In particular, for $n=1$ we have $\s_a(Y)=\s_Y=\overline{\s}_Y^P=\overline{\s}_Y=\q_Y^P=\q_Y=|Y|-1$ for each $P\in Y$. 

For a later usage, we single out the relations
\begin{equation}\label{EQ8}
\max(\overline{\s}^P_Y|P\in Y)=\overline{\s}_Y\le\s_Y.
\end{equation}

The next example shows that in general the inequality $\overline{\s}_Y\le\s_Y$ is strict.

\begin{example}\normalfont\label{EX4}
Suppose that $K$ is a finite field, $n\ge 2$, and for $Y\subset K^n$ we have $|Y|\ge |K|^n-|K|^{n-1}$. If $|Y|=|K|^n-|K|^{n-1}$, then we assume that $K^n\setminus Y$ is not an affine set (equivalently, we have $\d_{K^n\setminus Y}=n$). Then for each linear projection $\pi:\mathbb A^n_K\rightarrow\mathbb A^1_K$, we have $\pi(Y)=\mathbb A^1_K(K)=K$. Thus for each $a\in\AGL_n(K)$ and $i\in \llbracket1,n\rrbracket$ we have $m_{i,Y,a}=|K|$. This implies that $\s_Y=n(|K|-1)$. If $Y\neq K^n$, then for each $P\in Y$, by considering $Q\in K^n\setminus Y$ and $a\in\AGL_n(K)$ such that the line generated by $a(P)$ and $a(Q)$ is the zero locus $x_2=\cdots=x_n=0$, we get that $\max(\overline{m}_{1,Y,P,a},\q_Y^P)\le |K|-2$. So $\overline{\s}_Y^P\le \overline{\s}_a^P(Y)\le n(|K|-1)-1$, and thus $\overline{\s}_Y<\s_Y$.
\end{example}

\begin{lemma}\label{F3}
Let $K$ be a field and $n\in\mathbb N^{\ast}$. Let $Y$ and $Z$ be finite subsets of $K^n$ with $\emptyset\neq Y\subset Z$. Then we have $\d_Y\le \d_Z$, $\s_Y\le\s_Z$, $\overline{\s}_Y\le\overline{\s}_Z$, and $\mq_Y\le\mq_Z$. 
\end{lemma}

\begin{proof}
The inequalities follow directly from the very definitions.
\end{proof}

\begin{proposition}\label{PR6}
Let $(n,m)\in (\mathbb N^{\ast})^2$. Let $K$ be a finite field with $|K|\ge\sqrt[n]{m}$. Let $Y\subset K^n$ be a subset with $m$ elements. Then the following properties hold.

\medskip
{\bf (1)} We have inequalities 
$$\max(\d_Y,\lceil \d_Y\sqrt[\d_Y]{m}\rceil -\d_Y)\le\s_Y\le \d_Y\min(m-\d_Y,|K|-1).$$ 

{\bf (2)} Suppose that $|K|=2$. Then $\overline{\s}_Y\le\s_Y=\d_Y$ and $m\in \llbracket\d_Y+1,2^{\d_Y}\rrbracket$. 

\smallskip
{\bf (3)} If $\d_Y=2$, then $\s_Y\in \llbracket2,2|K|-2\rrbracket$ and $m\in \llbracket\lceil\frac{\s_Y}{2}\rceil+2,\lfloor\frac{(\s_Y+2)^2}{4}\rfloor\rrbracket$.\end{proposition}

\begin{proof}
Up to an affine automorphism we can assume that $\langle Y\rangle_{\aff}$ is the zero locus $x_{\d_Y+1}=\cdots=x_n=0$ and $\s_Y=\s_{1_{\mathbb A^n_K}}(Y)$. 

Let $m_i:=m_{i,Y,1_{\mathbb A^n_K}}$ for $i\in \llbracket1,\d_Y\rrbracket$; so $\s_Y=\sum_{i=1}^{\d_Y} (m_i-1)$. As $m_i\le |K|$ for each $i\in \llbracket1,\d_Y\rrbracket$, we get that $\s_Y\le \d_Y(|K|-1)$.

If there exists $i\in \llbracket1,\d_Y\rrbracket$ such that $m_i=1$, then $\d_Y\le\d_Y-1$, a contradiction. Thus $m_i\in \llbracket2,|K|\rrbracket$ for each $i\in \llbracket1,\d_Y\rrbracket$. Hence $\d_Y\le\s_Y$. 

Clearly, $m\le\prod_{i=1}^{\d_Y} m_i\le |K|^{\d_Y}$. This and the arithmetic and geometric means inequality for $m_1,\ldots, m_{\d_Y}$ give $\s_Y\ge -\d_Y+\d_Y\sqrt[\d_Y]{\prod_{i=1}^{\d_Y} m_i}\ge -\d_Y+\d_Y\sqrt[\d_Y]{m}$. So $\lceil \d_Y\sqrt[\d_Y]{m}\rceil -\d_Y\le\s_Y$. 

If $a\in\AGL_n(K)$ is such that $a(P_1)=(0,\ldots,0)$ and for each $i\in \llbracket2,\d_Y+1\rrbracket$, $a(P_i)$ is the $i$-th point of the standard basis of $\langle Y\rangle_{\aff}$ over $K$, then we have an inequality $\s_a(Y)\le \sum_{i=1}^{\d_Y} (1+m-\d_Y-1)=\d_Y(m-\d_Y)$. As $\s_Y\le\s_a(Y)$ by the very definition of capacities, we get that $\s_Y\le \d_Y(m-\d_Y)$. Thus part (1) holds.

If $|K|=2$, then $\d_Y\le\s_Y\le \d_Y$ by part (1) and $m_1=\cdots=m_{\d_Y}=2$. Thus $\d_Y=s_Y$ and $m\le \prod_{i=1}^{\d_Y} m_i=2^{\d_Y}$. Hence, as $m\ge \d_Y+1$ and $\overline{\s}_Y\le\s_Y$ by Equation (\ref{EQ8}), part (2) holds.

Part (3) follows directly from part (1).\end{proof}

\begin{proposition}\label{PR7}
Let $K$ be a field and $n\in\mathbb N^{\ast}$. Let $Y\subset K^n$ be a non-empty subset and let $P\in Y$. Then the following properties hold.

\medskip
{\bf (1)} If $K$ is finite, then we have an inequality $$\q_Y\le \Bigl\lfloor\frac{|Y|-1}{\sum_{i=0}^{n-1} |K|^i}\Bigr\rfloor.$$

{\bf (2)} If $K$ is infinite and $n\ge 2$, then $\q_Y=0$.

\smallskip
{\bf (3)} Suppose that $n\ge 2$. For each linear projection $\pi:\mathbb A^n_K\rightarrow\mathbb A^{n-1}_K$ we have an inequality
$$\overline{\s}_Y^P\le \bigl|\pi^{-1}\bigl(\pi(P)\bigr)\cap Y\bigr|+\overline{\s}_{\pi(Y)}^{\pi(P)}.$$
In particular, by choosing $\pi$ such that $\bigl|\pi^{-1}\bigl(\pi(P)\bigr)\cap Y\bigr|=\q_Y^P$, we get that 
$$\overline{\s}_Y^P\le\q_Y^P+\overline{\s}_{\pi(Y)}^{\pi(P)}\le \q_Y^P+(n-1)(|K|-1).$$

{\bf (4)} For each $P\in Y$ we have an inequality $\overline{\s}_Y^P\le\mq_Y^P$. In particular, 
$$\overline{\s}_Y\le\mq_Y.$$

{\bf (5)} Suppose that $n\ge 2$ and $K$ is finite. Then
$$\overline{\s}_Y\le\q_Y+(n-1)(|K|-1)\le\Bigl\lfloor\frac{|Y|-1}{\sum_{i=0}^{n-1} |K|^i}\Bigr\rfloor+(n-1)(|K|-1).$$

{\bf (6)} If $K$ is infinite, then we have an inclusion $\mathcal Q_Y^P\subset \{0\}^{n-1}\times \llbracket0,|Y|-1\rrbracket$. In particular, $\mq_Y\le |Y|-1$. 

\smallskip
{\bf (7)} If $K$ is infinite, then we have an inequality $\overline{\s}_Y\le |Y|-\max(1,\d_Y)$.
\end{proposition}

\begin{proof}
For part (1), we can assume that $n\ge 2$. Let $q:=\bigl\lfloor\frac{|Y|-1}{\sum_{i=0}^{n-1} |K|^i}\bigr\rfloor$. As the number of lines that pass through $P$ is $\frac{|K|^n-1}{|K|-1}=\sum_{i=0}^{n-1} |K|^i$, there exits such a line that intersects $Y\setminus\{P\}$ in at most $q$ points; so $\q_Y^P\le q$ and part (1) holds. 

Part (2) holds as there exist infinitely many lines that pass through $P$.

Part (3) follows from the very definitions.

Part (4) holds for $n=1$ as $\overline{\s}_Y^P=\q_Y=\mq_Y^P$. The general case follows by induction on $n$ based on part (3).

For part (5), its first inequality follows from parts (3) and (4) and the fact that $\overline{\s}_{\pi(Y)}^{\pi(P)}\le (n-1)(|K|-1)$ and its second inequality follows from the first one and part (1). So part (5) holds. 

Part (6) follows from part (2). 

\phantomsection{For part (7), it suffices to show that we have $\overline{\s}_Y^P\le |Y|-\max(1,\d_Y)$ for each $P\in Y$. For this we can assume that $P=(0,\ldots,0)$ and $n=\d_Y\ge 2$. Let $\{Q_i|i\in \llbracket1,\ldots,\d_Y\rrbracket\}\subset Y$ be a subset such that $(Q_1,\ldots,Q_{\d_Y})$ is a $K$-basis of $\Span(Y)=K^n$. Let $W$ be the hyperplane of dimension $\d_Y-1$ that contains this basis and let $V$ be the $K$-vector subspace of $K^n$ of dimension $\d_Y-1$ that does not intersect $V$; so $P\in V$.}\label{PH13j} 

\phantomsection{In this paragraph we assume that $\d_Y\ge 3$. Then, as $K$ is infinite, there exists a line $\lambda$ in $V$ that passes through $P$ and such that $V\cap Y=\{P\}$. We consider the image $Y_1$ of $Y$ via the linear projection $\pi_1:\mathbb A^n_K\rightarrow\mathbb A^{n-1}_K$ which at the level of $K$-valued points is the $K$-linear quotient map $K^n\rightarrow K^n/\lambda$. Clearly, we have $\d_{Y_1}=\d_Y-1$ and $|Y_1|\le |Y|$. Moreover, we can choose $\lambda$ such that the subset $\{\pi_1(Q_i)|i\in \llbracket1,\ldots,\d_Y\rrbracket\}$ of $W_1:=\pi_1(W)$ (and thus also of $Y_1$) has $\d_Y$ points and contains a $K$-basis of $\Span(Y_1)$.}\label{PH14j} 

Repeating the process $\d_Y-2$ times, we end up with a subset $Y_{\d_Y-2}$ in $K^2$ with $\d_{\d_Y-2}=2$ and a line $W_{\d_Y-2}$ in $K^2$ with $|W\cap Y_{\d_Y-2}|\ge d_Y$ and $P\notin W_{\d_Y-2}$. Let $V_{\d_Y-2}$ be the line in $K^2$ parallel to $W_{\d_Y-2}$ that contains $P$ and we denote $s:=|V_{\d_Y-2}\cap Y_{\d_Y-2}|\in\mathbb N^{\ast}$. Considering the linear projection $\pi_{\d_Y-1}:\mathbb A^2_K\rightarrow\mathbb A^1_K$ which at the level of $K$-valued points is the $K$-linear quotient map $K^2\rightarrow K^2/V_{\d_Y-2}$, we have $|\pi_{\d_Y-1}(Y_{\d_Y-2})|\le |Y|-(s-1)-(\d_Y-1)$ and therefore we estimate $\overline{\s}_Y^P\le (s-1)+[|Y|-(s-1)-(\d_Y-1)-1]=|Y|-\d_Y$. So part (7) holds.\end{proof}

\section{Arithmetic partition functions on linear independence}\label{S6}

Next we introduce arithmetic partition functions\index{partitions function} that help in identifying polynomials of small degrees that define concrete functions from $Y$ to $K$.

\begin{definition}\label{D6}
{\bf (1)} Let $k\in\mathbb N^{\ast}\setminus\{1\}$. By the modulo $k$ non-selective linearly independent partitions function\index{partitions function!modulo $k$ non-selective linearly independent} we mean the arithmetic function 
$$\L_k:\mathbb N\rightarrow\mathbb N$$ 
defined recursively by the following rules.

\medskip\noindent
{\bf (1.a)} We have $\L_k(0):=0$.

\smallskip\noindent
{\bf (1.b)} If $s\in \llbracket k^r,k^{r+1}-1\rrbracket$ with $r\in\mathbb N$, then $\L_k(s):=1+\L_k(s-r-1)$.

\medskip
{\bf (2)} Let $K$ be a finite (resp.\ infinite) field and $d\in \mathbb N^{\ast}$. By the selective linearly independent in dimension $\le d$ over $K$ partitions function\index{partitions function!selective linearly independent in dimension $\le d$ over $K$} we mean the function 
$$\l^{[\le d]}_K:\llbracket0,|K|^d-1\rrbracket\rightarrow\mathbb N\;\;\;(\textup{resp.}\; \l^{[\le d]}_K:\mathbb N\rightarrow\mathbb N)$$ defined by the rule: for $s\in \llbracket0,|K|^d-1\rrbracket$ (resp.\ $s\in\mathbb N$), 
$$\l^{[\le d]}_K(s):=\max\bigl(\L_Y\bigl|Y\subset K^d\setminus\{(0,\ldots,0)\},\; |Y|=s\bigr).$$

{\bf (3)} Let $K$ be a finite (resp.\ infinite) field and $d\in \mathbb N^{\ast}$. By the selective linearly independent in dimension $d$ over $K$ partitions function\index{partitions function!selective linearly independent in dimension $d$ over $K$} we mean the function 
$$\l^{[d]}_K:\llbracket d+1,|K|^d-1\rrbracket\rightarrow\mathbb N\;\;\;(\textup{resp.}\; \l^{[d]}_K:\mathbb N\rightarrow\mathbb N)$$ 
defined by the rule: for $s\in \llbracket d+1,|K|^d-1\rrbracket$ (resp.\ $s\in\mathbb N$), 
$$\l^{[d]}_K(s):=\max\bigl(\L_Y\bigl|Y\subset K^d\setminus\{(0,\ldots,0)\},\; |Y|=s,\; \d_Y=d\bigr).$$

{\bf (4)} Let $K$ be a field and $d\in \mathbb N^{\ast}$. By the selective linearly independent over $K$ partitions function\index{partitions function!selective linearly independent over $K$} we mean the function 
$$\l_K:\mathbb N\rightarrow\mathbb N$$ 
defined by the rule: for $s\in\mathbb N$, $\l_K(s)\in\mathbb N$ is the smallest such that each subset of a $K$-vector space formed by $s$ non-zero vectors can be partitioned into at most $\l_K(s)$ linearly independent parts.
\end{definition}

If $K$ is infinite, then by using finite sets of collinear non-zero vectors of arbitrary cardinality in $K$-vector spaces we get that $\l_K:\mathbb N\rightarrow\mathbb N$ is the identity function.

Clearly, the functions $\L_k$ and $\l_K$ are non-decreasing and surjective and the sequence $\bigl(\l^{[\le d]}_K(s)\bigr)_{d\ge\lceil\log_{|K|}(s+1)\rceil}$ if $K$ is finite and the sequence $\bigl(\l^{[\le d]}_K(s)\bigr)_{d\ge 1}$ if $K$ is infinite are non-decreasing.

Though it is not used in what follows, we note that for each $(k,l)\in\bigl(\mathbb N\setminus\{1\}\bigr)\times\mathbb N^{\ast}$, there exists a smallest $N_{k,l}\in\mathbb N^{\ast}$ such that for all integers $s\ge N_{k,l}$ we have $\L_k(s)\le \frac{s}{l}$. 

The next lemma justifies the terminology of Definition \ref{D6} and the fact that for a finite field $K$ we have an identity
\begin{equation}\label{EQ9}
\l_K(s)=\max\bigl(\l^{[\le d]}_K(s)\bigl|d\in\mathbb N^{\ast}, s+1\le |K|^d\bigr).
\end{equation}

\begin{lemma}\label{L4}
Let $K$ be a finite field and $d\in\mathbb N^{\ast}$. Then for each $s\in\mathbb N$ we have 
$$\l_K(s)\le\L_{|K|}(s)\le s.$$
\end{lemma}

\begin{proof}
We prove this by induction on $s\in\mathbb N$. As $\l_K(0)=\L_{|K|}(0)=0$, the base of the induction holds. Assuming that $s\in \llbracket|K|^r,|K|^{r+1}-1\rrbracket$ with $r\in\mathbb N$ and that $\l_K(s-r-1)\le\L_{|K|}(s-r-1)\le s-r-1$, to prove that $\l_K(s)\le\L_{|K|}(s)$ we consider $s$ non-zero vectors $P_1,\ldots,P_s$ in some $K$-vector space $V$. Let $Y:=\{P_1,\ldots,P_s\}$. As $s\ge |K|^r$, $\Span(Y)$ has dimension at least $r+1$ and hence there exists a linearly independent subset $Z_0$ of $Y$ with $|Z_0|=r+1$. As $Y\setminus Z_0$ can be partitioned in linearly independent parts $Z_1,\ldots,Z_j$ for some $j\in \llbracket0,\l_K(s-r-1)\rrbracket$, $Y=\sqcup_{i=0}^j Z_j$ is a partition in at most $1+j\le 1+\l_K(s-r-1)$ linearly independent parts. Based on this, our inductive assumption, and Definition \ref{D6}(4) and (1.b) we estimate
$$\l_K(s)\le 1+\l_K(s-r-1)\le 1+\L_{|K|}(s-r-1)=\L_{|K|}(s)\le s-r\le s.$$ 
This ends the induction and hence the proof.\end{proof}

We group together some basic monotone properties and inequalities as follows.

\begin{lemma}\label{F4}
Let $K$ be a field and $n\in\mathbb N^{\ast}$. Let $Y$ and $Z$ be finite subsets of $K^n$ with $\emptyset\neq Y\subset Z$. Then the following properties hold.

\medskip
{\bf (1)} We have $\l_Y\le\l_Z$. If $Z$ does not contain the zero vector of $K^n$, then $\L_Y\le\L_Z$.

\smallskip
{\bf (2)} If $Y$ is collinear, then $\l_Y=\s_Y=\overline{\s}_Y=\mq_Y=|Y|-1$.

\smallskip
{\bf (3)} Let $d\in\mathbb N^{\ast}$. Then for each $s\in \llbracket d+1,|K|^d-1\rrbracket$ we have an inequality 
$$\l_K^{[d]}(s)\le 1+\l_K^{[\le d]}(s-d).$$ 

{\bf (4)} Let $d\in\mathbb N^{\ast}$. If $Y$ is non-empty, then
$$\l_Y\le \l_K^{[\d_Y]}(|Y|-1)\le 1+\l_K^{[\le\d_Y]}(|Y|-\d_Y-1)\le |Y|-\d_Y.$$

{\bf (5)} Let $(d,s)\in\mathbb N^{\ast}\times\mathbb N$ with $s\ge d+1$. If $|K|\ge s-d+2$, then $\l_K^{[d]}(s)=s-d+1$. If moreover $|K|\ge s+1$, then $\l_K^{[\le d]}(s)=\l_K(s)=s$.

{\bf (6)} The functions $\l_K^{[d]}$ and $\l_K^{[\le d]}$ with $d\in\mathbb N^{\ast}$ and $\l_K$ depend only on $|K|$.
\end{lemma}

\begin{proof}
Part (1) follows directly from the very definitions. 

For part (2), as $Y$ is collinear, for each $P\in Y$, every non-empty linearly independent subset of $-P+(Y\setminus\{P\})$ has cardinality $1$, hence $\l_Y=|Y|-1$. Clearly, $\s_Y=\overline{\s}_Y=\mq_Y=|Y|-1$. So part (2) holds. 

To prove part (3) we can assume that $n=d=\d_Y$ and $Y\subset K^d\setminus\{(0,\ldots,0)\}$ with $|Y|=s$. Let $W\subset Y$ be linearly independent with $|W|=d$. Then $Y\setminus W$ can be partitioned in at most $\l_K^{[\le d]}(s-d)$ subsets that are linearly independent, hence $Y$ can be partitioned in at most $1+\l_K^{[\le d]}(s-d)$ subsets that are linearly independent. From this and Definition \ref{D6}(3) we get that part (3) holds.

To prove part (4) we can assume that $n=\d_Y$. Let $P\in Y$. We consider $W:=-P+(Y\setminus\{P\})\subset K^{\d_Y}$. As $\L_W\le \l_K^{[\d_Y]}(|W|)$ by Definition \ref{D6}(3) and $|W|=|Y|-1$, part (4) follows from part (3) and Lemma \ref{L4}.

For part (5), the inequality $\l_K^{[d]}(s)\le s-d+1$ follows from part (4) applied to a subset $Y$ of $K^d$ with $|Y|=s+1$ and $\d_Y=d$. 

As $|K|\ge s-d+2$, we can take $Y\subset K^d\setminus\{(0,\ldots,0)\}$ with $\d_Y=d$ and $|Y|=s$ such that it contains $s-d+1$ points on a line and additional $d-1$ points so that $\d_Y=d$. If $Z$ is the subset of $Y$ formed by the collinear $s-d+1$ points, then $\L_Y\ge \L_Z$ by part (1) and $\L_Z=s-d+1$. Hence $\L_Y\ge s-d+1$ and therefore $\l_K^{[d]}(s)\ge s-d+1$. Thus $\l_K^{[d]}(s)=s-d+1$. If moreover $|K|\ge s+1$, then $\l_K^{[i]}(s)=s-i+1$ for each $i\in \llbracket1,d\rrbracket$ and this implies that $\l_K^{[\le d]}(s)=\l_K(s)=s$. So part (5) holds.

Part (6) is clear if $K$ is finite and it follows from part (5) if $K$ is infinite.\end{proof}

The next elementary abstract lemma is used to compute the first values of the functions of Definition \ref{D6}. Before stating it, we need some definitions.

\begin{definition}\label{D7}
 Let $(r,s) \in (\mathbb N^{\ast}\setminus\{1\})\times \mathbb N$ and $A$ a set with $|A|\ge s$. Let $\mathcal F_{r,s}(A)$ be the set whose elements are sets $\mathcal I=\{I_1,\ldots,I_r\}$ with $I_1,\ldots,I_r$ pairwise disjoint finite subsets of $A$ with $\sum_{i=1}^r |I_i|=s$.\footnote{We have $|\mathcal I|\in \llbracket1,r\rrbracket$ as some of the subsets could be empty.} 
 
\medskip
{\bf (1)} If $\min(r,s)\ge 2$, by a $2$-removal operation\index{removal operation} on $\mathcal I\in\mathcal F_{r,s}(A)$, to be called shortly removal operation, we mean the replacement of $\mathcal I$ by $\mathcal I'=\{I_1',\ldots,I_r'\}\in \mathcal F_{r,s-2}(A)$ for which there exists $(i_1,i_2)\in \llbracket1,r\rrbracket^2$ with $i_1<i_2$, $P_1\in I_{i_1}$, and $P_2\in I_{i_2}$ such that $I_{i_1}':=I_{i_1}\setminus\{P_1\}$, $I_{i_2}':=I_{i_1}\setminus\{P_2\}$, and $I_j':=I_j$ for each $j\in \llbracket1,r\rrbracket\setminus\{i_1,i_2\}$. If $\min(|I_{i_1}|,|I_{i_2}|)\ge\max(|I_j|j\in \llbracket1,r\rrbracket\setminus\{i_1,i_2\})$ we call it a simple removal operation\index{removal operation!simple removal operation}.
 
 \smallskip
{\bf (2)} The partial order $\leq$ on $\mathcal F_{r,s}(A)$ is defined by the rule: if $\mathcal I=\{I_1,\ldots,I_r\}$ and $\mathcal J=\{J_1,\ldots,j_r\}$ are in $\mathcal F_{r,s}(A)$, then $\mathcal I\leq\mathcal J$ iff there exists $(\sigma,\varsigma)\in S_r^2$ such that $|I_{\sigma(1)}|\ge\cdots\ge|I_{\sigma(r)}|$, $|J_{\varsigma(1)}|\ge\cdots \ge|I_{\varsigma(r)}|$, and for each $i\in \llbracket1,r\rrbracket$ we have $\sum_{j=1}^i |I_{\sigma(j)}|\le\sum_{j=1}^i |J_{\varsigma(j)}|$.
\end{definition}

A removal operation on $\mathcal I\in\mathcal F_{r,s}(A)$ exists iff $\mathcal I$ has at least two non-empty sets, and in such a case we have $\min(r,s)\ge 2$. We have $|\mathcal F_{r,0}(A)|=1$. If $A$ is a finite set, then $|\mathcal F_{r,1}(A)|=|A|$, $|\mathcal F_{1,2}(A)|=\binom{|A|}{2}$, and $|\mathcal F_{r,2}(A)|=|A|(|A|-1)$ for $r\ge 2$. 

We view the next lemma as the abstract essence of the construction in Section \ref{S7} of functions $f:Y\rightarrow K$ with $Y\subset K^2$ that are products of linear polynomials.

\begin{lemma}\label{L5} Let $(r,s) \in (\mathbb N^{\ast}\setminus\{1\})\times \mathbb N$ and $A$ a set with $|A|\ge s$. Let
$$\l=\l_{r,s}^A:\mathcal F_{r,s}(A)\rightarrow \llbracket0,s\rrbracket$$
be the function defined as follows. For $\mathcal I=\{I_1,\ldots,I_r\}\in\mathcal F_{r,s}(A)$ we consider all possible (finite) sequences of removal operations that start at $\mathcal I$ and terminate at a set $\mathcal I'=\{I_1',\ldots,I_r'\}\in\mathcal F_{r,s'}(A)$ with at most one element which is not an empty set; let $\l(\mathcal I)$ be the smallest possible value for $s'$. Then the following properties hold.

\medskip
{\bf (1)} With respect to the partial order of Definition \ref{D7}(2), the function $\l$ is non-decreasing.

\smallskip
{\bf (2)} In a sequence of only simple removal operations, the sequence starting at $\mathcal I$ terminates at a set in $\mathcal F_{r,\l(\mathcal I)}(A)$.

\smallskip
{\bf (3)} We have $\l(\mathcal I)=0$ iff $s$ is even and $2|I_i|\le s$ for each $i\in \llbracket1,r\rrbracket$.

\smallskip
{\bf (4)} We have $\l(\mathcal I)=1$ iff $s$ is odd and $2|I_i|\le s+1$ for each $i\in \llbracket1,r\rrbracket$.

\smallskip
{\bf (5)} We have $\l(\mathcal I)\ge 2$ iff there exists $i\in \llbracket1,r\rrbracket$ such that
$|I_i|>\lceil\frac{s}{2}\rceil$, and in such a case $\l(\mathcal I)=2|I_i|-s.$
\end{lemma}

\begin{proof} Note that if we only care about the cardinalities of the sets that are elements of $\mathcal I$, then a simple removal operation is essentially unique.

Parts (1) and (2) together are equivalent to proving that if we have an inequality $\mathcal I=\{I_1,\ldots,I_r\} \leq \mathcal J=\{J_1,\ldots,J_r\}$ and if terminating sequences of removal operations for $\mathcal I$ and $\mathcal J$ produce $\mathcal I'\in\mathcal F_{r,s'}(A)$ and $\mathcal J'\in\mathcal F_{r,s''}(A)$ (respectively), then the inequality $s'\le s''$ holds provided all operations in the sequence that produces $\mathcal I'$ are simple removal operations.

We prove the above statement by induction on $s$. The base of the induction, for $s=0$ and $s=1,$ is obvious. For $s\ge 2$, assuming that the statement holds for $s-2$, the passage from $s-2$ to $s$ goes as follows. For $\mathcal I\leq\mathcal J$ we can assume, by reindexing the sets, that $|I_1|\geq |I_2| \geq \cdots \geq |I_r|$ and $|J_1|\geq |J_2| \geq \cdots \geq |J_r|$. Then for every $i\in\llbracket1,r\rrbracket$ we have $\sum_{j=1}^i |I_j|\leq \sum_{j=1}^i |J_j|$ by Definition \ref{D7}(2). Suppose the first removal operation in the sequence for $\mathcal J$ removes elements in $J_{i_1}$ and $J_{i_2}$ with $(i_1,i_2)\in \llbracket1,r\rrbracket^2$, $i_1<i_2$. Then $|J_{i_2}|>0,$ therefore $\sum_{j=1}^{i_2-1} |I_j|\leq \sum_{j=1}^{i_2-1} |J_j|<s,$ which implies that $|I_{i_2}|>0.$ So we can perform a removal operation for $\mathcal I$ with the same indices $i_1$ and $i_2$. We also have $|I_{i_1}|\geq |I_{i_2}|>0$. Denote by $\mathcal J^{i_1,i_2}$ the result of the removal operation that removes elements $J_{i_1}$ and $J_{i_2}$, and denote similarly $\mathcal I^{i_1,i_2}$ and $\mathcal I^{1,2}$. Then for these three elements of $\mathcal F_{r,s-2}(A)$, we have $ \mathcal I^{1,2}\leq\mathcal I^{i_1,i_2} \leq\mathcal J^{i_1,i_2}$, and we can apply the induction hypothesis to $\mathcal I^{1,2}$ and $\mathcal J^{i_1,i_2}$ to obtain the result. This completes the induction.

\phantomsection{To prove parts (3), (4), and (5), let 
$$N_{\mathcal I}:=\max(|I_i||i\in \llbracket1,r\rrbracket)\;\;\;\textup{and}\;\;\;D_{\mathcal I}:=N_{\mathcal I}-\lceil \frac{s}{2}\rceil$$ 
and we denote by $\mathcal I^{\O}$ an arbitrary sequence in $\mathcal F_{r,s-2}$ obtained from $\mathcal I$ after a simple removal operation.}\label{PH16n}

We first prove that a simple removal operation does not change the value of truth of the statement $\mathfrak S_{\mathcal I}$ that there exists $i\in \llbracket1,r\rrbracket$ with $|I_i|> \lceil\frac{s}{2}\rceil$.

If such an $i$ exists, then we check that the positive value of $D_{\mathcal I}$ does not change after a simple removal operation. We have $N_{\mathcal I}=|I_i|$ and $I_i$ is the unique set in $\mathcal I$ of cardinality $N_{\mathcal I}$. After a simple removal operation, $I_i$ has one fewer element and thus $N_{\mathcal I^{\O}}=N_{\mathcal I}-1$ and hence $D_{\mathcal I^{\O}}=N_{\mathcal I}-1-\lceil\frac{s-2}{2}\rceil=D_{\mathcal I}$. 

Suppose that $D_{\mathcal I}\leq 0$. If $N_{\mathcal I^{\O}}\le N_{\mathcal I}-1$, then $D_{\mathcal I^{\O}}\le N_{\mathcal I}-1-\lceil\frac{s-2}{2}\rceil=D_{\mathcal I}\leq 0$. If $N_{\mathcal I^{\O}}=N_{\mathcal I}$, then we must have had at least $3$ sets of cardinality $N_{\mathcal I}\geq 1$ in $\mathcal I$. So $s\geq 3N_{\mathcal I}$, and hence $\lceil\frac{s-2}{2}\rceil \geq \lceil\frac{3N_{\mathcal I}-2}{2}\rceil \geq N_{\mathcal I}=N_{\mathcal I^{\O}}$, i.e., $D_{\mathcal I^{\O}}\le 0$.

From the last two paragraphs we get that indeed the value of truth of the statement $\mathfrak S_{\mathcal I}$ does not change after simple removal operations.

Now, for $\mathcal I$ we consider three cases, that correspond to parts (3), (4), and (5): (i) $D_{\mathcal I}\leq 0$ and $s$ is even, (ii) $D_{\mathcal I}\leq 0$ and $s$ is odd, and (iii) $D_{\mathcal I}>0$. Simple removal operations preserve this division into cases by the last paragraph, and in the last case the largest set $I_i$ remains the largest and the value of $D_{\mathcal I}$ does not change. Based on this and part (2), to prove parts (3) to (5) we can assume that $s=\l(\mathcal I)$. So $\mathcal I$ has at most one non-empty subset and thus $\l(\mathcal I)$ is $0$ for part (3), is $1$ for part (4), and is $s=2|I_i|-s$ for part (5). So parts (3) to (5) hold.
\end{proof}

The first application of Lemma \ref{L5} is as follows. 

\begin{proposition}\label{PR8} {\bf (1)} Let $k\in\mathbb N^{\ast}\setminus\{1\}$. Then we have identities $\L_k(s)=s$ for each $s\in \llbracket0,k-1\rrbracket$ and $\L_k(s)=k-1+\lfloor\frac{s+1-k}{2}\rfloor$ for each $s\in \llbracket k,k^2-1\rrbracket$. 

\smallskip
{\bf (2)} Let $K$ be a finite field and $d\in\mathbb N^{\ast}$. Then the following properties hold.

\medskip\noindent
{\bf (2.a)} We have $\l_K(s)=\l^{[\le d]}_K(s)=s$ for each $s\in \llbracket0,|K|-1\rrbracket$. 

\smallskip\noindent
{\bf (2.b)} If $d\ge 2$, then $\l_K(s)=\l^{[\le d]}_K(s)=|K|-1$ for each $s\in \llbracket|K|,2|K|-2\rrbracket$.

\smallskip\noindent
{\bf (2.c)} If $d\ge 2$, then $\l_K(s)=\l^{[\le d]}_K(s)=\lfloor\frac{s+1}{2}\rfloor$ for each $s\in \llbracket2|K|-1,|K|^2-1\rrbracket$.
\end{proposition}

\begin{proof}
As $\L_k(0)=0$ and $\L_k(s)=1+\L_k(s-1)$ for $s\in \llbracket1,k-1\rrbracket$, an induction on $s\in \llbracket0,k-1\rrbracket$ gives that $\L_k(s)=s$. Similarly, as we have $\L_k(k-1)=k-1$, $\L_k(k)=1+\L_k(k-2)=k-1$, and $\L_k(s)=2+\L_k(s-2)$ for $s\in \llbracket k,k^2-1\rrbracket$, an induction on $s\in \llbracket k,k^2-1\rrbracket$ gives that $\L_k(s)=k-1+\lfloor\frac{s+1-k}{2}\rfloor$. So part (1) holds.

Clearly, $\l^{[\le 1]}_K(s)=s$ for each $s\in \llbracket0,|K|-1\rrbracket$. From this, Equation (\ref{EQ9}), the inequality $\l_K(s)\le s$ (see Lemma \ref{L4}), and the fact that the sequence $\bigl(\l^{[\le d]}_K(s)\bigr)_{d\ge 1}$ is non-decreasing for $s\in \llbracket0,|K|-1\rrbracket$, we get that $\l_K(s)=\l^{[\le d]}_K(s)=s$ for each $(s,d)\in \llbracket0,|K|-1\rrbracket\times\mathbb N^{\ast}$. So part (2.a) holds.

One can use Equation (\ref{EQ6.9}) to prove parts (2.b) and (2.c). For readers convenience we include self-contained proves of them.

We now assume that $s\in \llbracket|K|,|K|^2-1\rrbracket$ and consider a subset $Y=\{P_1,\ldots,P_s\}$ of non-zero points in a $K$-vector space $V$ of dimension $d\ge 2$. Let $\lambda_1,\ldots,\lambda_r$ be all the $1$-dimensional subspaces of $V$ that have non-trivial intersection with $Y$. Let $\mathcal I:=\{\lambda_i\cap Y|i\in \llbracket1,r\rrbracket\}$ and $N:=\max(\lambda_i\cap Y|i\in \llbracket1,r\rrbracket)$. As $Y=\sqcup_{i=1}^r \lambda_i\cap Y$, we have $\mathcal I\in\mathcal F_{r,s}(V)$ and $N\in \llbracket0,|K|-1\rrbracket$. 

Let $l:=\l_{r,s}^V(\mathcal I)\in \llbracket0,s\rrbracket$ be as in Lemma \ref{L5}; so $s-l=2j$ where $j\in\mathbb N$ is the largest number of simple removal operations possible. We consider a partition $Y=\sqcup_{i=1}^{j+l} Y_i$ into linearly independent subsets of $V$ such that $|Y_1|=\cdots=|Y_j|=2$ and $|Y_{j+1}|=\cdots=|Y_{j+l}|=1$ by the very definition of $l$.

Assume that $s\in \llbracket|K|,2|K|-2\rrbracket$. If $l=0$, then $j+l=\frac{s}{2}\le |K|-1$. If $l=1$, then $s\le 2|K|-3$ and hence $j+l=\frac{s+1}{2}\le |K|-1$. If $l\ge 2$, then $l=2N-s$ by Lemma \ref{L5}(5) and $j+l=\frac{s+l}{2}=N\le |K|-1$. Thus $\l_K(s)\le |K|-1$. If $r=2$ and $\lambda_1\setminus\{0\}\subset Y$, then $Y$ cannot be partitioned into at most $|K|-2$ linearly independent subsets and hence $|K|-1\le \l_K^{[\le 2]}(s)\le \l_K^{[\le d]}(s)\le\l_K(s)$. From the last two sentences we get that part (2.b) holds.

Assume that $s\in \llbracket2|K|-1,|K|^2-1\rrbracket$. So $2N\le 2(|K|-1)\le s-1$. From this and Lemma \ref{L5}(3) and (4) we get that $l\in\{0,1\}$. As $\l_K(s)\le j+l=\frac{s+l}{2}\le\lfloor\frac{s+1}{2}\rfloor$, we get that $\l_K(s)\le\lfloor\frac{s+1}{2}\rfloor$. Clearly, $\lfloor\frac{s+1}{2}\rfloor\le \l_K^{[\le 2]}(s)\le \l_K^{[\le d]}(s)\le \l_K(s)$. From the last two sentences we get that part (2.c) holds.\end{proof}

\section{Polynomial functions on finite subsets}\label{S7}

Recall that the Kronecker delta $\delta_{ij}\in\{0,1\}$ defined for two elements $i$ and $j$ of some set is $1$ iff $i=j$. For a field $K$, $n\in\mathbb N^{\ast}$, and $P\in Y\subset K^n$, let 
$$f^P_Y:Y\rightarrow K$$ 
be the function defined by $f_Y^P(Q)=\delta_{PQ}$ for each $Q\in Y$. 

Next example illustrates the relevance of the invariants of Definition \ref{D5}.

\begin{example}\normalfont\label{EX5}
Let $K$ be a field, $n\in\mathbb N^{\ast}$, and $Y\subset K^n$ a non-empty finite subset. Let $P=(\alpha_1,\ldots,\alpha_n)\in Y(K)$. For $i\in \llbracket1,n\rrbracket$ let $Z_i\subset Y_i\subset K$ be as Definition \ref{D5}(1) and (2) for the automorphism $a=1_{\mathbb A^n_K}\in\AGL_n(K)$; so $Y_i$ is the set formed by the $i$-th coordinates of points in $Y$. We have $Z_n=Y_n$, $\sum_{i=1}^n (|Y_i|-1)=\s_{1_{\mathbb A^n_K}}(Y)$, and $\sum_{i=1}^n (|Z_i|-1)=\overline{\s}^P_{1_{\mathbb A^n_K}}(Y)$. 

\medskip
{\bf (1)} Let
\begin{equation*}
\chi_Y^P:=\prod_{i=1}^n\prod_{\beta_i\in Y_i\setminus\{\alpha_i\}} \frac{x_i-\beta_i}{\alpha_i-\beta_i}\in K[x_1,\ldots,x_n].
\end{equation*}
So $\deg(\chi_Y^P)=\s_{1_{\mathbb A^n_K}}(Y)$, $\chi^P_Y$ is a product of linear polynomials, and $\chi_Y^P$ represents the function $f_Y^P$, i.e., we have $\chi_Y^P(Q)=\delta_{PQ} $ for each $Q\in Y\setminus\{P\}$. In particular, if $\s_Y=\s_{1_{\mathbb A^n_K}}(Y)$ then $\deg(\chi_Y^P)=\s_Y$. 

\smallskip
{\bf (2)} Suppose that $n\ge 2$. Let 
\begin{equation*}
\overline{\chi}_Y^P:=\prod_{i=1}^n\prod_{\beta_i\in Z_i\setminus\{\alpha_i\}} \frac{x_i-\beta_i}{\alpha_i-\beta_i}\in K[x_1,\ldots,x_n].
\end{equation*}
So $\deg(\overline{\chi}_Y^P)=\overline{\s}_{1_{\mathbb A^n_K}}^P(Y)$, $\overline{\chi}_Y^P$ is a product of linear polynomials, $\overline{\chi}_Y^P$ also represents the function $f^P_Y$, and we have $\deg(\overline{\chi}_Y^P)\le\deg(\chi_Y^P)$. In particular, if $\overline{\s}^P_Y=\overline{\s}^P_{1_{\mathbb A^n_K}}(Y)$, then $\deg(\overline{\chi}_Y^P)=\overline{\s}_Y^P$. 

\smallskip
{\bf (3)} If $K$ is finite, then for
\begin{equation}\label{EQ10}
\chi^P:=\chi^P_{K^n}=\overline{\chi}^P_{K^n}=\prod_{i=1}^n\prod_{\beta_i\in K\setminus\{\alpha_i\}} \frac{x_i-\beta_i}{\alpha_i-\beta_i}\in K[x_1,\ldots,x_n]
\end{equation}
we have $\deg(\chi^P)=n(|K|-1)$ and $\chi^P(Q)=\delta_{PQ}$ for each $Q\in K^n\setminus\{P\}$.
\end{example}

\begin{definition}\label{D8}
Let the field $K$ and $(n,m)\in (\mathbb N^{\ast})^2$ be such that $\sqrt[n]{m}\le |K|$. Let $Y\subset K^n$ be such that $|Y|=m$. 

\medskip
{\bf (1)} By the representation degree of $Y$\index{representation degree} we mean the smallest $\r_Y\in\mathbb N$ such that each function $Y\rightarrow K$ can be represented by a polynomial in $K[x_1,\ldots,x_n]$ of degree at most $\r_Y$.

\smallskip
{\bf (2)} Let $Y\in P$. By the representation degree of $Y$ at $P$\index{representation degree!representation degree at a point} we mean the smallest $\r_Y^P\in\mathbb N$ such that the function $f^P_Y$ is represented by a polynomial in $K[x_1,\ldots,x_n]$ of degree $\r_Y^P$.
\end{definition}

\begin{proposition}\label{PR9}
Let $K$ be a field and $(n,m)\in (\mathbb N^{\ast})^2$. Let $Y\subset K^n$ be such that $|Y|=m$. Then the following properties hold.

\medskip
{\bf (1)} We have $\r_Y=\max(\r_Y^P|P\in Y)$.

\smallskip
{\bf (2)} For each $P\in Y$ we have $\r_Y^P\le\overline{\s}_Y^P$.

\smallskip
{\bf (3)} We have $\r_Y\le\overline{\s}_Y$.

\smallskip
{\bf (4)} Let $Z$ be a non-empty subset of $Y$. For each $P\in Z$ we have $\r_Z^P\le\r_Y^P$. In particular, $\r_Z\le\r_Y$.
\end{proposition}

\begin{proof}
By the very definitions we have $\r_Y^P\le\r_Y$. Thus $\r_Y\ge\max(\r_Y^P|P\in Y)$. As for each function $h:Y\rightarrow K$ we have 
an identity $h=\sum_{P\in Y} h(P)f^P_Y$, it follows that $\r_Y\le\max(\r_Y^P|P\in Y)$. So part (1) holds.

Part (2) follows from Example \ref{EX5}(2) as up to affine automorphisms we can assume that $\overline{\s}_Y^P=\overline{\s}^P_{1_{\mathbb A^n_K}}(Y)$.

Part (3) follows from parts (1) and (2).

Part (4) follows from the definitions.\end{proof}

Directly from Propositions \ref{PR9}(3) and \ref{PR7}(4) and Equation (\ref{EQ8}) we get the following consequence.

\begin{corollary}\label{C6}
Let $K$ be a field, $m\in\mathbb N^{\ast}$, and $Y$ a non-empty subset of $K^n$. Then we have inequalities
$$\r_Y\le\overline{\s}_Y\le\min(\s_Y,\mq_Y).$$
\end{corollary}

Proposition \ref{PR7}(1) and (3) suggests the following abstract notion.

\begin{definition}\label{D9}
{\bf (1)} Let $k\in\mathbb N^{\ast}\setminus\{1\}$. By the modulo $k$ strict capacity function\index{capacity!modulo $k$ strict capacity function} we mean the function 
$$\overline{\s}_k:\cup_{n\in\mathbb N^{\ast}} \{n\}\times \llbracket1,k^n\rrbracket \rightarrow\mathbb N$$ 
defined recursively by the following rules.

\medskip\noindent
{\bf (1.a)} We have $\overline{\s}_k(1,m):=m-1$ for each $m\in \llbracket1,k\rrbracket$.

\smallskip\noindent
{\bf (1.b)} For each $n\in\mathbb N^{\ast}\setminus\{1\}$ and every $m\in \llbracket1,k^n\rrbracket$ we have 
$$\overline{\s}_k(n,m):=\Bigl\lfloor\frac{m-1}{\sum_{i=0}^{n-1} k^i}\Bigr\rfloor+\overline{\s}_k\biggl(n-1,\min\Bigl(k^{n-1},m-\Bigl\lfloor\frac{m-1}{\sum_{i=0}^{n-1} k^i}\Bigr\rfloor\Bigr)\biggr).$$

{\bf (2)} By the modulo $\infty$ strict capacity function\index{capacity!modulo $\infty$ strict capacity function} we mean the function 
$$\overline{\s}_{\infty}: (\mathbb N^{\ast})^2\rightarrow\mathbb N$$ 
defined by the rule: $\overline{\s}_{\infty}(n,m)=m-1$ for each $(n,m)\in (\mathbb N^{\ast})^2$.
\end{definition}

\begin{lemma}\label{L6}
Let $(k,n)\in (\mathbb N^{\ast}\setminus\{1\})\times\mathbb N^{\ast}$. Let $m\in \llbracket1,k^n\rrbracket$. Then the following properties hold.

\medskip
{\bf (1)} For each $s\in \llbracket1,m\rrbracket$ we have inequalities 
$$\overline{\s}_k(n,m)+s-m\le\overline{\s}_k(n,s)\le\overline{\s}_k(n,m)\le n(k-1).$$

{\bf (2)} We have $\overline{\s}_k(n,m)=n(k-1)$ iff $m=k^n$.

\smallskip
{\bf (3)} Let $l\in \llbracket1,n\rrbracket$ be the smallest with $m\le \sum_{i=0}^l k^i$ and let $s:=\min(k^l,m)$. Then the following properties hold.

\medskip\noindent
{\bf (3.a)} We have an identity $\overline{\s}_k(n,m)=\overline{\s}_k(l,s)$. In particular, if $s=k^l$, then $\overline{\s}_k(n,m)=l(k-1)$.

\smallskip\noindent
{\bf (3.b)} Suppose that $s=m\in \llbracket1+\sum_{i=0}^{l-1} k^i,k^l-1\rrbracket$ and $l\ge 2$. Then 
$$\overline{\s}_k(n,m)=\Bigl\lfloor\frac{m-1}{\sum_{i=0}^{l-1} k^i}\Bigr\rfloor+(l-1)(k-1).$$

\noindent
{\bf (3.c)} Suppose that $l=1$, i.e., $m\in \llbracket1,k+1\rrbracket$. Then $\overline{\s}_k(n,m)=\min(k,m)-1$.
\end{lemma}

\begin{proof}
We have $\overline{\s}_k(n,1)=0$ for each $n\in\mathbb N^{\ast}$. Thus for parts (1) and (2) we can assume that $m\ge 2$ and $s=m-1$. We prove both parts by induction on $n\in\mathbb N^{\ast}$. The case $n=1$ follows from Definition \ref{D9}(1.a). 

For $n\ge 2$, for the passage from $n-1$ to $n$ let $t:=\bigl\lfloor\frac{m-1}{\sum_{i=0}^{n-1} k^i}\bigr\rfloor\in \llbracket0,k-1\rrbracket$ and $r:=\bigl\lfloor\frac{m-2}{\sum_{i=0}^{n-1} k^i}\bigr\rfloor\in \llbracket0,k-2\rrbracket$. We have $t\in\{r,r+1\}$.

Assume that $t=r$. Hence $t\le k-2$. Based on Definition \ref{D9}(1.b) and the inductive hypothesis we get that
$$\overline{\s}_k(n,m)-1=t-1+\overline{\s}_k\bigl(n-1,\min(k^{n-1},m-t)\bigr)\le t+\overline{\s}_k\bigl(n-1,\min(k^{n-1},m-1-t)\bigr)$$
$$=\overline{\s}_k(n,m-1)\le t+\overline{\s}_k\bigl(n-1,\min(k^{n-1},m-t)\bigr)=\overline{\s}_k(n,m)\le t+(n-1)(k-1)$$
$$\le k-2+(n-1)(k-1)=n(k-1)-1.$$

If $t=r+1\in\mathbb N^{\ast}$, then $m-1=t\sum_{i=0}^{n-1} k^i$ and hence 
$$m-t=1-t+t\sum_{i=0}^{n-1} k^i\ge 1+tk^{n-1}>k^{n-1}.$$ 
Again based on Definition \ref{D9}(1.b) and the inductive hypothesis we get
$$\overline{\s}_k(n,m)-1=t-1+\overline{\s}_k(n-1,k^{n-1})=r+\overline{\s}_k(n-1,k^{n-1})=\overline{\s}_k(n,m-1)$$
$$<t+\overline{\s}_k(n-1,k^{n-1})=\overline{\s}_k(n,m)\le t+(n-1)(k-1)\le n(k-1).$$
If $\overline{\s}_k(n,m)=n(k-1)$, then $t=k-1$, hence $m=k^n$. 

This ends the inductive step and the induction. So parts (1) and (2) hold.

For part (3) let $d:=n-l\in \llbracket0,n-1\rrbracket$ and we use the same $t\in \llbracket0,k-1\rrbracket$. We prove part (3.a) by induction on $d$. If $d=0$, then $n=l$, $s=m=k^l$, and therefore $\overline{\s}_k(n,m)=\overline{\s}_k(l,s)=l(k-1)$; so the base of the induction holds. For $d\ge 1$, the passage from $d-1$ to $d$ goes as follows. As $l=n-d\le n-1$ and $m\le \sum_{i=0}^l k^i$, we have $t=0$. Hence $\overline{\s}_k(n,m)=\overline{\s}_k\bigl(n-1,\min(k^{n-1},m)\bigr)$. As $s=\min\bigl(k^l,\min(k^{n-1},m)\bigr)$ and $n-1-l=d-1$, from the inductive hypothesis we get the following identity $\overline{\s}_k\bigl(n-1,\min(k^{n-1},m)\bigr)=\overline{\s}_k(l,s)$. Thus $\overline{\s}_k(n,m)=\overline{\s}_k(l,s)$. This ends the inductive step and the induction. So part (3.a) holds.

For part (3.b), let $q:=\bigl\lfloor\frac{m-1}{\sum_{i=0}^{l-1} k^i}\bigr\rfloor\in \llbracket1,k-1\rrbracket$. We write $m-1=q(\sum_{i=0}^{l-1} k^i)+r_1$ with $r_1\in \llbracket0,\sum_{i=1}^{l-1} k^i\rrbracket$. Then $m-1-q=q(\sum_{i=1}^{l-1} k^i)+r_1$ is greater than or equal to $k^{l-1}$ as $l\ge 2$. Thus $\overline{\s}_k(l,s)=q+\overline{\s}_k(l-1,k^{l-1})=q+(l-1)(k-1)$ by Definition \ref{D8}(2) and part (2). From this and part (3.a) we get that part (3.b) holds.

For part (3.c), we have $s=\min(k,m)$. As $\overline{\s}_k(1,s)=s-1$ by Definition \ref{D9}(1.a), part (3.c) follows from part (3.a).\end{proof}

\begin{theorem}\label{T4}
Let $n\in\mathbb N^{\ast}$ and $K$ a finite field. Let $Y$ be a non-empty subset of $K^n$. Let $l\in \llbracket1,n\rrbracket$ be the smallest with $|Y|\le \sum_{i=0}^l |K|^i$. Let $s:=\min(|K|^l,|Y|)$. Then the following properties hold.

\medskip
{\bf (1)} We have inequalities 
$$\r_Y\le\overline{\s}_Y\le\min\bigl(\overline{\s}_{|K|}(l,s),\mq_Y\bigr)\le\max\bigl(\overline{\s}_{|K|}(l,s),\mq_Y\bigr)\le l(|K|-1).$$ 

{\bf (2)} If moreover $|Y|< |K|^l$, then $\max\bigl(\overline{\s}_{|K|}(l,s),\mq_Y\bigr)\le l(|K|-1)-1$.

\smallskip
{\bf (3)} For each $(l',s')\in (\mathbb N^{\ast})^2$ with $l'\ge l$ and $s'\ge s$ we have $\r_Y\le\overline{\s}_Y\le\overline{\s}_{|K|}(l',s')$.
\end{theorem}

\begin{proof}
Let $q:=\bigl\lfloor\frac{|Y|-1}{\sum_{i=0}^{n-1} |K|^i}\bigr\rfloor$. 

For part (1), we prove the inequality $\overline{\s}_Y\le\overline{\s}_{|K|}(n,|Y|)$ by induction on $n\in\mathbb N^{\ast}$. For $n=1$ we have $\overline{\s}_Y=\overline{\s}_{|K|}(1,|Y|)=|Y|-1$ by Lemma \ref{F4}(2) and Definition \ref{D9}(1.a), so the base of the induction holds. For $n\ge 2$, for the passage from $n-1$ to $n$ let $P\in Y$. Proposition \ref{PR7}(1) and (3) gives that $\q_Y^P\le q$ and that for each linear projection $\pi:\mathbb A^n_K\rightarrow\mathbb A^{n-1}_K$ with $|\pi^{-1}\bigl(\pi(P)\bigr)\cap Y|=\q_Y^P$ we have $\overline{\s}_Y^P\le \q_Y^P+\overline{\s}_{\pi(Y)}^{\pi(P)}$. Clearly, $|\pi(Y)|\le \min(|Y|-\q_Y^P,|K|^{n-1})$. From the last two sentences, the inductive assumption, Lemma \ref{L6}(1), and Definition \ref{D9}(1.b) we get that
$$\overline{\s}_Y^P\le \q_Y^P+\overline{\s}_{\pi(Y)}^{\pi(P)}\le\q_Y^P+\overline{\s}_{|K|}\bigl(n-1,|\pi(Y)|\bigr)\le\q_Y^P+\overline{\s}_{|K|}\bigl(n-1,\min(|K|^{n-1},|Y|-\q_Y^P)\bigr)$$
$$\le q+\overline{\s}_{|K|}\bigl(n-1,\min(|K|^{n-1},|Y|-\q)\bigr)=\overline{\s}_{|K|}(n,|Y|).$$
Therefore $\overline{\s}_Y\le\overline{\s}_{|K|}(n,|Y|)$. This ends the inductive step and the induction.

Based on the last paragraph, Corollary \ref{C6}, and Lemma \ref{L6}(1) and (3.a), for part (1) it suffices to prove that $\mq_Y^P\le l(|K|-1)$.

To prove $\mq_Y^P\le l(|K|-1)$ we can assume that $n>l$. Hence $q=0$ and thus $\q_Y^P=0$ and $\mq_Y^P\le (n-1)(|K|-1)$. If $l\le n-2$, then the argument can be repeated to give that $\q_{\pi(Y)}^{\pi(P)}=0$. By performing the argument $n-l$ times it follows that for each $n$-tuple $(q_1,\ldots,q_n)\in\mathcal Q_Y^P$ we have $q_1=\cdots=q_{n-l}=0$ and thus 
$$\sum_{i=0}^n q_i=\sum_{i=n-l+1}^n q_i\le \sum_{i=n-l+1}^n (|K|-1)=l(|K|-1).$$ 
Hence the inequality $\mq_Y^P\le l(|K|-1)$ holds. So part (1) holds.

For part (2), the image of $Y$ under each linear projection $\mathbb A^n_K\rightarrow\mathbb A^l_K$ has cardinality less than $|K|^l$ and this implies that $q_{n-l+1}\le |K|-2$ (cf.\ Example \ref{EX4}; hence $\sum_{i=0}^n q_i\le l(|K|-1)-1$ and $\mq_Y^P\le l(|K|-1)-1$. So part (2) holds.

As we have $\overline{\s}_{|K|}(l,s)=\overline{\s}_{|K|}(l',s)\le \overline{\s}_{|K|}(l',s')$ by Lemma \ref{L6}(3.a) and (1), part (3) follows from part (1).\end{proof}

\begin{example}\normalfont\label{EX6}
Suppose that $Y\subset\mathbb F_3^2$ with $|Y|=6$. Let $Z:=\mathbb F_3^2\setminus Y$.

\medskip
{\bf (1)} Suppose that $Z$ is collinear. To check that $\r_Y=3$, up to affine automorphisms we can assume that $Z$ is the zero locus $x_2=2$. Let $P:=(0,0)\in Y$. We show that the assumption that $f^P_Y$ is represented by a polynomial $g\in K[x_1,\ldots,x_n]$ of degree at most $2$ leads to a contradiction. This assumption implies that there exists a quintuple $(\alpha,\beta,\gamma,\delta_1,\delta_2)\in\mathbb F_3^5$ such that the polynomial $$g:=\alpha x_1^2+\beta x_1x_2+\gamma x_2^2+\delta_1x_1+\delta_2x_2+1$$ 
represents $f^P_Y$. The identities $g(1,0)=g(2,0)=0$ imply that $\delta_1=0$ and $\alpha=2$. So $g(x_1,1)$ is a polynomial of degree $2$ in $x_1$ that has three roots in $\mathbb F_3$, a contradiction. Thus $\r_Y\ge 3$. As $\r_Y\le 3$ by Theorem \ref{T4}(1), we have $\r_Y=3$. It is easy to see that $\s_Y=3\ge\mq_Y$. As $\r_Y=\s_Y=3$, from the last two sentences and Corollary \ref{C6} we get that $\r_Y=\overline{\s}_Y=\mq_Y=\s_Y=3$.

\smallskip
{\bf (2)} Suppose that $Z$ is non-collinear. We have $\s_Y=4$ by Example \ref{EX4}. To compute the invariants $\r_Y$, $\overline{\s}_Y$, and $\mq_Y$, up to affine automorphisms we can assume that $Z=\{(1,2),(2,1),(2,2)\}$. Let $P\in Y$. Up to the linear automorphism $\e(x_2,x_1)$ we can assume that $P$ is $(0,0)$, $(0,1)$, $(0,2)$ or $(1,1)$ and hence $f^P_Y$ is represented by $2x^2+xy+2y^2+1$, $2xy+2y^2+2y$, $2y^2+y$ or $xy$ (respectively). Thus $\r^P_Y\le 2$. Clearly, $\r_Y^P>1$. So $\r_Y^P=2$. If $Q\in\{(0,0),(2,0),(0,2)\}$, then $\mathcal Q^Q_Y=\{(0,2)\}$, hence $\mq_Y^Q=2\ge\overline{\s}_Y^Q$. If $Q\in\{(1,1),(1,0),(0,1)\}$, then $\mathcal Q^Q_Y=\{(1,2)\}$ and thus $\mq_Y^Q=3$; moreover it is easy to see that $\overline{\s}_Y^Q=3$. Hence $\mq_Y=\overline{\s}_Y=3$. Therefore we have
$$\r_Y=2<\overline{\s}_Y=\mq_Y=3<\s_Y=4.$$
\end{example}

The following basic result refines \cite{Poo}, Lem.\ 2.1 from multiple points of view and can be viewed as an applications of Lemma \ref{L5} and Theorem \ref{T4}.

\begin{theorem}\label{T5}
Let $(n,m)\in (\mathbb N^{\ast})^2$. Let $K$ be a field with $|K|\ge\sqrt[n]{m}$. Let $Y\subset K^n$ be a subset with $m$ elements. Then the following properties hold.

\medskip
{\bf (1)} We have inequalities 
$$\r_Y\le\min(\l_Y,\overline{\s}_Y)\le\min\bigl(\l^{[\d_Y]}_K(m-1),\overline{\s}_Y\bigr)\le\min\bigl(1+\l^{[\le \d_Y]}_K(m-\d_Y-1),\overline{\s}_Y\bigr).$$
In particular, $\r_Y\le\min(\l_Y,\overline{\s}_Y)\le\min\bigl(m-\d_Y,\d_Y(|K|-1)\bigr)$.

\smallskip
{\bf (2)} Suppose that $|K|=2$. Then 
$$\r_Y\le\min(\l_Y,\overline{\s}_Y)\le\min\bigl(\l_K^{[\d_Y]}(m-1),\d_Y\bigr)$$
$$\le\min\bigl(1+\l_K^{[\le\d_Y]}(m-\d_Y-1),\d_Y\bigr)\le \min\Bigl(\bigl\lfloor\frac{m-\d_Y+2}{2}\bigr\rfloor,\d_Y\Bigr)\le\min\Bigl(\bigl\lfloor\frac{m}{2}\bigr\rfloor,n\Bigr).$$\end{theorem}

\begin{proof}
As $\r_Y\le\overline{\s}_Y$ by Proposition \ref{PR9}(3) and based on Lemma \ref{F4}(4) and Inequalities (\ref{EQ7}), for part (1) it suffices to show that for each $P\in Y$, the function $f_Y^P$ is represented by a polynomial $h$ with $\deg(h)\le\l_Y$. 

\phantomsection{Let $l:=\L_{-P+(Y\setminus\{P\})}\in \llbracket1,\l_Y\rrbracket$. Let $-P+(Y\setminus\{P\})=\sqcup_{j=1}^lZ_j$ be a partition into linearly independent sets. As for $j\in \llbracket1,l\rrbracket$, $(P+Z_j)\cup\{P\}\subset Y$ is affinely independent, there exists a linear polynomial $g_j\in R$ with $g_j(P+Z_j)=\{0\}$ and $g_j(P)=1$. As $h:=\prod_{j=1}^l g_j\in R$ represents $f_Y^P$ and $\deg(h)=l\le\l_Y$, part (1) holds.}\label{PH16z}

For part (2), let $q:=m-\d_Y-1\in\mathbb N$. Based on part (1) we can assume that $\d_Y\ge 2$. We have $\l_K^{[\le \d_Y]}(q)\le\L_2(q)$ by Lemma \ref{L4}. As $\L_2(0)=0$, $\L_2(i)=1$ for $i\in\{1,2\}$, $\L_2(3)=2$, and $\L_2(i)\le 1+\L_2(i-3)$ if $i\ge 4$, it follows that for $i\in\mathbb N\setminus\{1,3\}$ we have $\L_2(i)\le \lfloor\frac{i}{2}\rfloor$. From this and Proposition \ref{PR8}(1) we get that 
$$1+\l_K^{[\le \d_Y]}(q)\le 1+\L_2(q)\le 1+\Bigl\lfloor\frac{q}{2}\Bigr\rfloor=\Bigl\lfloor\frac{m-\d_Y+1}{2}\Bigr\rfloor$$ 
if $q\notin\{1,3\}$. So part (2) follows from part (1) and Proposition \ref{PR6}(2) if $q\notin\{1,3\}$. 

If $q=1$, then $m=\d_Y+2$ and $1+\L_2(q)=2$ and hence, as $\d_Y\ge 2$, we have $m\ge 4$ and $\min(2,\d_Y)=\min\bigl(\lfloor\frac{m-\d_Y+2}{2}\rfloor,\d_Y\bigr)=2\le\min(\lfloor\frac{m}{2}\rfloor,n)$. Similarly, if $q=3$, then $m=\d_Y+4$ and $1+\L_2(q)$=3, and hence, as $\d_Y\ge 2$, we have $m\ge 6$ and $\min(3,\d_Y)=\min(\lfloor\frac{m-\d_Y+2}{2}\rfloor,\d_Y)\le\min(\lfloor\frac{m}{2}\rfloor,n)$, part (2) holds.\end{proof}

\begin{corollary}\label{C7} 
Let $K$ be a field. Let $Y$ be a non-empty subset of $K^n$. Then the following properties hold. 

\medskip
{\bf (1)} If $Y$ is collinear, then for each $P\in Y$ we have $\r_Y^P=\r_Y=|Y|-1$.

\smallskip
{\bf (2)} If $\c_Y\ge\bigl\lceil\frac{|Y|+1}{2}\bigr\rceil$, then $\r_Y=\c_Y-1$.

\smallskip
{\bf (3)} If $\c_Y<\bigl\lceil\frac{|Y|+1}{2}\bigr\rceil$, then $\c_Y-1\le\r_Y\le\bigl\lfloor\frac{|Y|+1}{2}\bigr\rfloor$.
\end{corollary}

\begin{proof}
For part (1), we have $\r_Y^P\le\r_Y\le |Y|-1$ by Lemma \ref{F4}(2) and Theorem \ref{T5}(1). We show that the assumption that $f_Y^P$ can be represented by a polynomial $g\in K[x_1,\ldots,x_n]$ with $\deg(g)<|Y|-1$ leads to a contradiction. Up to affine automorphisms we can assume that $P=(0,\ldots,0)$ and that $Y$ is contained in the zero locus $x_2=\cdots=x_n=0$. So the polynomial $h(x_1):=g(x_1,0,\ldots,0)\in K[x_1]$ is such that $h(0)=1$ and it has $|Y|-1$ distinct zeros. Hence $\deg(h)\ge |Y|-1>\deg(g)$, a contradiction. Thus $\r_Y^Y\ge |Y|-1$. We conclude that part (1) holds.

For parts (2) and (3), let $Z$ be a collinear subset of $Y$ with $|Z|=\c_Y$. We have $\r_Z=|Z|-1=\c_Y-1$ by Lemma \ref{F4}(2). From this and Proposition \ref{PR9}(4) we get that $\c_Y-1\le\r_Y$. For $P\in Y$, to study $\r_Y^P$, up to affine automorphisms we can assume that $P=(0,\ldots,0)$.

For part (2), we write $|Y|=2s+r$ with $(s,r)\in\mathbb N\times\{0,1\}$. Let $t\in \llbracket0,s+r-1\rrbracket$ be such that $|Z|=s+1+t$; so $\c_Y-1=s+t$. We consider two cases as follows.

{\bf Case 1: $P\notin Z$.} Writing $Y\setminus (Z\cup\{P\})=\{Q_1,\ldots,Q_{s+r-t-2}\}$, for every $i\in \llbracket1,s+r-t-2\rrbracket$, let $g_i\in R$ be a linear polynomial such that $g_i(P)=1$ and $g_i(Q_i)=0$. If $g_0\in R$ is a linear polynomial such that $g_0(P)=1$ and $g_0(Z)=\{0\}$, then the product $\prod_{i=0}^{s+r-t-2} g_i$ has degree $s+r-t-1$ and represents $f_Y^P$. Thus $\r_Y^P\le s+r-t-1\le s+t=\c_Y-1$.

{\bf Case 2: $P\in Z$.} We partition $Y\setminus\{P\}=\sqcup_{i=1}^{s+t} W_i$ into linearly independent subsets such that for each $i\in\llbracket1,\min(s+t,s-t+r)\rrbracket$ we have $|W_i|=2$ and $|W_i\cap Z|=1$ and for every $i\in \llbracket \min(s+t,s-t+r)+1,\max(s+t,s-t+r)$ we have $|W_i|=1$\footnote{We have $s+t<\min(s+t,s-t+r)$ iff $t=0$ and $r=1$, in which case $s-t+r=s+t+1$.}; therefore we have $\L_{-P+(Y\setminus \{P\})}\le s+t$. From the proof of Theorem \ref{T5}(1) we get that $\r_Y^P\le\L_{-P+(Y\setminus \{P\})}$. Thus $\r_Y^P\le s+t=\c_Y-1$. 

From the two cases we get that part (2) holds.

For part (3), as $\c_Y\le \lceil\frac{|Y|-1}{2}\rceil$, from Lemma \ref{L5}(3) and (4) applied to the set formed by the intersections of lines passing through the origin $P$ with $Y\setminus\{P\}$ we get that we can partition $Y\setminus\{P\}$ into $\bigl\lfloor\frac{|Y|+1}{2}\bigr\rfloor$ linearly independent subsets, among which the first $\bigl\lfloor\frac{|Y|-1}{2}\bigr\rfloor$ subsets have $2$ elements and correspond to simple removal operations. So as in Case 2 we get that $\r_Y^P\le \bigl\lfloor\frac{|Y|+1}{2}\bigr\rfloor$ from which part (3) follows. 
\end{proof}

The following application of Theorem \ref{T5}(1) is used in what follows.

\begin{corollary}\label{C8}
Let $(n,m)\in (\mathbb N^\ast\setminus\{1\})\times (\mathbb N^\ast\setminus\{1,2\})$. Let $K$ be a field such that $|K|\ge m-1$. Let $\underline{P}=(P_1,\ldots,P_m)\in\mathbb D_{n,m}(K)$ be such that for the set $Y:=\{P_1,\ldots,P_m\}\subset K^n$ there exists a linear projection $\pi:\mathbb A^n_{K}\rightarrow\mathbb A^1_{K}$ with $|\pi(Y)|=m-1$. Then the following properties hold.

\medskip
{\bf (1)} If $\underline{Q}=(Q_1,\ldots,Q_m)\in\mathbb D_{n,m}(K)$ is such that $\pi(P_i)=\pi(Q_i)$ for every $i\in \llbracket1,m\rrbracket$, then there exists $a\in\TGA_n(K)[m-2]$ such that $a(\underline{P})=\underline{Q}$.

\smallskip
{\bf (2)} Suppose that $|K|\ge m$. Let $j=j(m,K)\in \llbracket2,m-1\rrbracket$ be the largest such that the inequality $m\le\frac{2|K|+j^2-j-2}{2(j-1)}$ holds. Then there exists $b\in\STGA_n(K)[m-j]$ such that the first coordinates of the $b(P_i)$s with $i\in\llbracket1,m\rrbracket$ are distinct.
\end{corollary}

\begin{proof}
Up to linear automorphisms we can assume that $\pi$ is the projection on the first coordinate and that by writing $P_i=(\alpha_i,\beta_i,O_i)\in K^n$ with $O_i\in K^{n-2}$ for $i\in\llbracket1,m\rrbracket$, we have $(\alpha_1,\ldots,\alpha_{m-1})\in\mathbb D_{1,m-1}(K)$, $\alpha_m=\alpha_{m-1}$, and $\beta_{m-1}\neq \beta_m$. 

To prove part (1), for each $(j,l)\in\{m-1,m\}\times \llbracket2,n\rrbracket$, let $\gamma_{j,l}\in K$ be the $l$-th coordinate of $Q_j$; hence $Q_j=(\alpha_j,\gamma_{j,2},\ldots,\gamma_{j,n})$. Up to a linear automorphism that does not change the first coordinates we can also assume that $\gamma_{m-1,l}\neq\gamma_{m,l}$ for each $l\in \llbracket2,n\rrbracket$. 

\phantomsection{As $\beta_{m-1}\neq\beta_m$ and $\gamma_{m-1,l}\neq\gamma_{m,l}$, the system of linear equations}\label{PH16y}
$$x\beta_{m-1}+y-\gamma_{m-1,l}=x\beta_m+y-\gamma_{m,l}=0$$ 
in the indeterminates $x$ and $y$ has a unique solution $(\delta_l,\delta_l^{\prime})\in K^{\ast}\times K$. Let 
$$a:=\e\bigl(x_1,\delta_2 x_2+f_2(x_1),\ldots,\delta_nx_n+f_n(x_1)\bigr)\in\TGA_n(K)$$ 
with $(f_2,\ldots,f_n)\in K[x_1]^{n-1}$ be such that $f_l(\alpha_{m-1})=\delta_l^{\prime}$ and $f_l(\alpha_1),\ldots,f_l(\alpha_{m-2})$ are uniquely determined by the identity $a(P_i)=Q_i$ for each $i\in \llbracket1,m\rrbracket$; we can assume that for each $l\in \llbracket2,n\rrbracket$ we have $\deg(f_l)\le m-2$ by Theorem \ref{T5}(1) applied to $n=1$; as $a\in\TGA_n(K)[m-2]$, part (1) holds.

To prove part (2), we can assume that $n=2$. We take $b$ to be a composite $b_1b_2$ with $b_1=\e(x_1+\alpha x_2,x_2)\in\GL_2(K)$ for some $\alpha\in K^{\ast}$ and $b_2=\e\bigl(x_1,x_2+g(x_1)\bigr)$ with $g\in K[x_1]$ of degree at most $m-j$; clearly $b\in\SGA_2(K)[m-j]$. We choose $g$ such that the second coordinate of $b_2(P_l)$ is $0$ for each $l\in \llbracket1,m-j+1\rrbracket$ by Lagrange interpolation. Up to interchanging the indices $m-1$ and $m$ we can assume that the second coordinate of $b_2(P_m)$ is non-zero. 

Let $I$ be the set of those $i\in \llbracket m-j+2,m\rrbracket$ such that the second coordinate of $b_2(P_m)$ is non-zero; so $m\in I$. Denoting $l:=|I|\in \llbracket1,j-1\rrbracket$, we write $I=\{i_1,\ldots,i_l\}$ with $i_1<i_2<\cdots<i_l=m$. For each $\alpha\in K^{\ast}$, the first coordinates of the $b(P_i)$s with $i\in \llbracket1,m\rrbracket\setminus I$ are distinct. By induction on $s\in \llbracket1,l-1\rrbracket$ we get that there exists a subset $I_s$ of $K^{\ast}$ of cardinality $m-l+s-1$ such that if $\alpha\in K^{\ast}\setminus\cup_{t=1}^s I_t$, then the first coordinates of the $b(P_i)$s with $i\in \llbracket1,m\rrbracket\setminus\{i_{s+1},\ldots,i_l\}$ are distinct. Similarly, as the first coordinates of $b(P_{m-1})$ and $b(P_m)$ are distinct for each $\alpha\in K^{\ast}$, there exists a set $I_l$ with $m-2$ elements such that if $\alpha\in K^{\ast}\setminus\cup_{t=1}^l I_t$, the first coordinates of $b(P_1)$ to $b(P_m)$ are distinct. We estimate 
$$\sum_{t=1}^l |I_t|\le (m-j+1)+(m-j+2)+\cdots+(m-3)+(m-2)+(m-2)$$
$$=\frac{m(m-1)}{2}-\frac{(m-j)(m-j+1)}{2}-1=(j-1)m-\frac{j(j-1)}{2}-1\le |K|-2,$$ where the last inequality holds as it is equivalent to the inequality $m\le\frac{2|K|+j^2-j-2}{2(j-1)}$. Thus the set $K^{\ast}\setminus\cup_{t=1}^s I_t$ is non-empty and for an element $\alpha$ of it we get that $b$ has the desired properties. So part (2) holds. 
\end{proof}

We have the following consequence of Proposition \ref{PR8}(2), Lemma \ref{L6}(3.b), and Theorems \ref{T4}(1)  and \ref{T5}(1) that involves bounds independent of $Y$ for $n=2$.

\begin{corollary}\label{C9}
Let $m\in\mathbb N^{\ast}\setminus\{1\}$. Then for each finite field $K$ with $|K|\ge\sqrt{m}$, each subset $Y\subset K^2$ with $|Y|=m$, and each function $f:Y\rightarrow K$, there exists a polynomial $g\in K[x_1,x_2]$ of degree at most $\varphi(m,|K|)$ that represents $f$, with $\varphi(m,|K|)$ independent of $Y$ and defined as follows.

\medskip
{\bf (1)} If $m\in \llbracket2,|K|\rrbracket$, then $\varphi(m,|K|):=m-1$.

\smallskip
{\bf (2)} If $m\in \llbracket|K|+1,2|K|-1\rrbracket$, then $\varphi(m,|K|):=|K|-1$.

\smallskip
{\bf (3)} If $m\in \llbracket2|K|,2|K|+2\rrbracket$ with $|K|\ge 3$, then $\varphi(m,|K|):=|K|$.

\smallskip
{\bf (4)} If $m\in \llbracket2|K|+3,|K|^2-1\rrbracket$ with $|K|\ge 4$, then $\varphi(m,|K|):=|K|+\lfloor\frac{m-1}{|K|+1}\rfloor-1$.

\smallskip
{\bf (5)} If $m=|K|^2$, then $\varphi(m,|K|):=2|K|-2$.
\end{corollary}

\begin{proof}
Part (1) (resp.\ (5)) follows from Theorem \ref{T5}(1) applied to $n=2$ and $\d_Y=1$ (resp.\ $\d_Y=2$). 

For part (2), as $m-1\in \llbracket|K|-1,2|K|-2\rrbracket$ we have $\l_K(m-1)=|K|-1$ by Proposition \ref{PR8}(2.a) and (2.b). From this and Theorem \ref{T5}(1) we get that part (2) holds.

Parts (3) and (4) follow from Theorem \ref{T4}(1) and Lemma \ref{L6}(3.b) applied to $n=l=2$ and $k=|K|$.\end{proof}

\begin{remark}\normalfont\label{R1}
{\bf (1)} In Corollary \ref{C9}(2) to (5) we have $\d_Y=2$.

\smallskip
{\bf (2)} If in Corollary \ref{C9}(1) (resp.\ \ref{C9}(2)), the set $Y$ is collinear (resp.\ contains the set $Z$ of $K$-valued points of a line in $\mathbb A^2_K$), then there exists a function $f:Y\rightarrow K$ which cannot be represented by a polynomial $g\in K[x_1,x_2]$ of degree less than $\varphi(m,|K|)$ by Corollary \ref{C7}(1) applied to $Y$ (resp.\ $Z$). So Corollaries \ref{C9}(1) and (2) are optimal for each fixed value of $m$. Similarly, for $|K|=3$, Corollary \ref{C9}(3) is optimal by Example \ref{EX6}(1). It is well-known that Corollary \ref{C9}(5) is optimal.
\end{remark}

\section{Selective shifts}\label{S8}

In this section we introduce certain `nice' permutations $a(K)\in\perm(K^n)$ with $a\in\TGA_n(K)$ for finite fields $K$. To begin with, the following definition is a natural extrapolation of the $K$-linear combinations of the proof of Theorem \ref{T5}(1).

\begin{definition}\label{D10}
Let $n\in\mathbb N^{\ast}$. We consider an (inner) direct sum decomposition $K^n=V\oplus W$ of $K$-vector spaces. Let $v_0\in V\setminus\{0\}$ and a non-empty subset $W_0$ of $W$. Let $\chi_{W_0}:W\rightarrow K$ be the characteristic function (sending $W_0$ to $1$ and $W\setminus W_0$ to $0$). We call the permutation $\shift^{v_0+W_0}_{V\oplus W}\in\perm(K^n)$ defined by the rule $v+w\mapsto v+\chi_{W_0}(w)v_0+w$ for $(v,w)\in V\times W$ as the selective shift\index{selective shift} of the selection of $W_0\subset W$ and the shift $v_0\in V$, or simply as a selective shift. If $W_0=\{w_0\}$ has only one element, then we denote $\shift^{v_0+W_0}_{V\oplus W}$ simply by $\shift^{v_0+w_0}_{V\oplus W}$
\end{definition}

Note that if $V=K^n$, so $W$ is the zero subspace of $K^n$, then $W_0=\{0\}$ and the permutation $\shift^{v_0+0}_{V\oplus W}$ is the translation by $v_0$.

\begin{lemma}\label{F5}
Suppose that $K$ is a finite field of characteristic $p$. Let $n\in\mathbb N^{\ast}\setminus\{1\}$. Let $d\in \llbracket0,n-1\rrbracket$. We consider an (inner) direct sum decomposition $K^n=V\oplus W$ of $K$-vector spaces with $\dim_K(W)=d$. Then the following properties hold for a selective shift $\shift^{v_0+W_0}_{V\oplus W}\in\perm(K^n)$ with $v_0\in V\setminus\{0\}$ and $W_0\subset W$.

\medskip
{\bf (1)} The permutation $\shift^{v_0+W_0}_{V\oplus W}$ is a product of $\frac{|W_0||K|^{n-d}}{p}$ disjoint $p$-cycles and its support is $V_0+W_0$.

\smallskip
{\bf (2)} There exists $a^{v_0+W_0}_{V\oplus W}\in\STGA_n(K)$ such that $a^{v_0+W_0}_{V\oplus W}(K)=\shift^{v_0+W_0}_{V\oplus W}$ and $\ell(a^{v_0+W_0}_{V\oplus W})\le \max\bigl(1,d(|K|-1)\bigr)$. Moreover, if $W_0=\{w_0\}$ has only one element, then $\ell(a^{v_0+w_0}_{V\oplus W})=\max\bigl(1,d(|K|-1)\bigr)$.

\smallskip
{\bf (3)} 
The permutation $\shift^{v_0+W_0}_{V\oplus W}$ is odd iff $(|K|,d)=(2,n-1)$ and $|W_0|$ is odd and iff it is a product of $|W_0|$ disjoint transpositions with $|W_0|$ odd.
\end{lemma}

\begin{proof}
For each $\overline{v}\in V/\mathbb F_pv_0$, let $P_{\overline{v}}\in V$ be such that $P_{\overline{v}}+\mathbb F_pv_0=\overline{v}$. As $V\cong K^{n-d}$ and we have a product decomposition into disjoint $p$-cycles
\begin{equation}\label{EQ11}
\shift^{v_0+W_0}_{V\oplus W}=\prod_{w\in W_0} \prod_{\overline{v}\in V/\mathbb F_pv_0} \bigl((P_{\overline{v}},w)\; (P_{\overline{v}}+v_0,w)\;\cdots\; (P_{\overline{v}}+(p-1)v_0,w)\bigr),
\end{equation}
parts (1) and (3) hold.
 
For part (2), up to special affine automorphisms we can assume that $V$ and $W$ are linear subspaces defined by the linear equations $x_{n-d+1}=\cdots=x_n=0$ and $x_1=\cdots =x_{n-d}=0$ (respectively) and that $v_0=(1,0,\ldots,0)$. 

Let $h\in K[x_{n-d+1},\ldots,x_n]$ be a polynomial of degree at most $d(|K|-1)$ that represents the function $\chi_{W_0}$ of Definition \ref{D10} by Theorem \ref{T4}(1). If $W_0=\{w_0\}$, then $\deg(h)=d(|K|-1)$ by Example \ref{EX5}(3). As we can take $a^{v_0+W_0}_{V\oplus W}$ to be the automorphism $\e(x_1+h,x_2,\ldots,x_n)$, part (2) holds.
\end{proof}

Due to their relevance to multiple places in the monograph, for each triple $(l,x,y)\in \mathbb N\times [1,\infty)^2$ we introduce the real numbers
\begin{equation}\label{EQ12}
E_{l,x,y}:=x^2(x-y)^2+lx(x^2-xy+1)
\end{equation}
and
\begin{equation}\label{EQ13}
E_{l,x}:=E_{l,x,0}=x^4+lx(x^2+1).
\end{equation}

\begin{proposition}\label{PR10}
Let $n\in\mathbb N^{\ast}\setminus\{1\}$ and $K$ a finite field. Let $p:=\char(K)$ and $k:=|K|-1$. Then the following properties hold.

\medskip
{\bf (1)} For each $s\in \llbracket1,p-1\rrbracket$ there exists 
$$a_s\in\STGA_n(K)[E_{n-2,k,s-1}]=\STGA_n(K)[k^2(k+1-s)^2+k(k^2+k-ks+1)(n-2)]$$ 
such that $a(K)$ is a $(2s+1)$-cycle in $\Alt(K^n)$ with $\d_{\supp\bigl(a(K)\bigr)}=2$.

\smallskip
{\bf (2)} If $|K|=3$ (equivalently, $k=2$), then for each $3$-cycle $\theta\in\Alt(K^n)$ of collinear (resp.\ non-collinear) points there exists $b\in\STGA_n(K)[2n-2]$ (resp.\ $b\in\STGA_n(K)[4n]$) such that $b(K)=\theta$.

\smallskip
{\bf (3)} If $|K|\ge 4$ (equivalently, $k\ge 3$), then for each $3$-cycle $\theta\in\Alt(K^n)$ of collinear (resp.\ non-collinear) points there exists $c\in\STGA_n(K)[4E_{n-2,k}]$ (resp.\ $c\in\STGA_n(K)[E_{n-2,k}]$) such that $c(K)=\theta$.
\end{proposition}

\begin{proof}
Let $q\in\mathbb N^{\ast}$ be such that $p^q=|K|=k+1$. The prime field of $K$ is $\mathbb F_p$ and we identify $\mathbb F_p=\{0,\ldots,p-1\}$. For $\alpha\in K$, let $P_{\alpha}:=(\alpha,0,\ldots,0)\in K^n$; so $P_p=P_0$. 

Let $V$ and $W$ be the $K$-vector subspaces of $K^n$ defined by the linear equations $x_2=\cdots=x_n=0$ and $x_1=0$ (respectively). Therefore $K^n=V\oplus W$ and $\dim_K(W)=n-1$. Let $h:=\prod_{i=3}^n (1-x_i^k)\in\mathbb F_p [x_3,\ldots,x_n]$ and
$$a:=\e\bigl(x_1+(1-x_2^k)h,x_2,\ldots,x_n\bigr)\in\STGA_n(K).$$ So $a(K)=\shift^{(1,0,\ldots,0)}_{V\oplus W}$ is a product $\prod_{j=1}^{p^{q-1}} \theta_j$ of $p^{q-1}$ disjoint permutations that are $p$-cycles by Lemma \ref{F5}(1) and (2) and $Y:=\supp\bigl(a(K)\bigr)=\{P_{\alpha}|\alpha\in K\}$ by Equation (\ref{EQ11}). We can assume that $Y_1:=\supp(\theta_1)=\{P_0,\ldots,P_{p-1}\}$; so we have $\theta_1=(P_0\; P_1\;\cdots\; P_{p-1})$. 

Let $f_s(x_1):=\prod_{i=1}^s (x_1-i)\in\mathbb F_p[x_1]$, $g_s(x_1):=\frac{x_1^{p^q}-x_1}{f_s(x_1)}\in\mathbb F_p[x_1]$, and 
$$b_s:=\e\bigl(x_1,x_2+g_s(x_1),x_3,\ldots,x_n)\bigr)\in\STGA_n(K)[k+1-s].$$ 
It follows that $b_s(K)$ fixes each element of the set $Y\setminus\{P_1,\ldots,P_s\}$ and we have $b_s(P_{\alpha})=\bigl(\alpha,g_s(\alpha),0,\ldots,0\bigr)\neq P_{\alpha}$ for each $\alpha\in \llbracket1,s\rrbracket$. Let 
$$c_s:=b_sa^{-1}b_s^{-1}\in\STGA_n(K);$$ 
so $c_s(K)$ is a product $\prod_{j=2}^{p^{q-1}+1} \theta_j^{-1}$ of $p^{q-1}$ disjoint permutations that are $p$-cycles, with $Y_2:=\supp(\theta_{p^{q-1}+1})=b_s(Y_1)$. Let 
$$a_s:=ac_s=ab_sa^{-1}b_s^{-1}\in\STGA_n(K).$$
As $|Y_1\cap Y_2|=p-s$ and 
\begin{equation}\label{EQ13.1}
a_s(K)=\theta_1\theta_{p^{q-1}+1}^{-1},
\end{equation}
$a_s(K)$ is a $(2s+1)$-cycle by Property \ref{P10} applied to $(t,s)=(s,p-s)$. We have $\d_{\supp\bigl(a(K)\bigr)}=\d_{Y_1\cup Y_2}=2$.

\phantomsection{One computes that $c_s=\e(x_1-F_{s,1},x_2-F_{s,2},x_3,\ldots,x_n)$, where}\label{PH16m}
$$F_{s,1}:=\bigl[1-\bigl(x_2-g_s(x_1)\bigr)^k\bigr]h,\;\;\;F_{s,2}:=g_s(x_1)-g_s\bigl(x_1-[1-(x_2-g_s(x_1)^k)h]\bigr).$$
As $\deg(F_{s,1})=k\deg(g_s)+\deg(h)\le\deg(F_{s,2})=\deg(g_s)[k\deg(g_s)+\deg(h)]$, $\deg(g_s)=k+1-s$, and $\deg(h)=k(n-2)$, we get that 
$$\deg(F_{s,2})=k(k+1-s)(k-s+n-1)$$ 
and hence $\ell(c_s)\le k(k+1-s)^2+k(k+1-s)(n-2)$.
A similar computation gives that $a_s=\e(x_1-F_{s,3},x_2-F_{s,2},x_3,\ldots,x_n)$, where
$$F_{s,3}:=F_{s,1}-[1-(x_2-F_{s,2})^k]h.$$
As $\deg(F_{s,1})\le\deg(F_{s,2})<\deg(F_{s,3})=k\deg(F_{s,2})+\deg(h)$, we get the inequality $\ell(a_s)\le E_{n-2,k,s-1}$. So part (1) holds.\footnote{The upper bounds for $\ell(c_s)$ and $\ell(a_s)$ are written without factors in order to emphasize the linear growth in $n$ when $k$ and $s$ are fixed.}

If $p=3$ and $q=s=1$, then $K\cong\mathbb F_3$, $a(K)\in\STGA_n(\mathbb F_3)[2n-2]$ is a $3$-cycle of collinear points and $c_1(\mathbb F_3)\in\STGA_n(\mathbb F_3)[4n]$ is a $3$-cycle of non-collinear points; thus, as $\AGL_n(K)$ acts transitively on the set of $3$-cycles of collinear (resp.\ non-collinear) points in $K^n$, we get that part (2) holds.

For part (3), let $Q_1:=b_1(P_1)=\bigl(1,g_1(1),0,\ldots,0\bigr)\in K^n\setminus\{P_{\alpha}|\alpha\in K\}$. We have $a_1\in\STGA_n(K)[E_{n-2,k}]$
with $a_1(K)=(P_1\; P_2\; Q_1)$ a $3$-cycle of non-collinear points. Up to affine automorphisms we can assume that $\theta$ is either $a_1(K)$ or $(P_1\; P_2\; P_{\alpha})$ with $\alpha\in K\setminus\{1,2\}$. If $\theta=a_1(K)$ then we can take $c:=a_1$. From now we assume that $\theta=(P_1\; P_2\; P_{\alpha})$ with $\alpha\in K\setminus\{1,2\}$. 

Let $h_0(x_1):=\frac{g_1(1)}{(\alpha-1)(\alpha-2)}(x_1-1)(x_1-2)\in K[x_1]$. Clearly, $h_0(1)=h_0(2)=0$, $h_0(\alpha)=g_1(1)\neq 0$, and $\deg(h_0)=2$. With $\beta:=\frac{\alpha-1}{g_1(1)}\in K^{\ast}$, let
$$c_0:=\e\bigl(x_1-\beta x_2-\beta h_0(x_1),x_2+h_0(x_1),x_3,\ldots,x_n\bigr)\in\STGA_n(K)[2].$$ 
As $c_0$ fixes $P_1$ and $P_2$ and maps $P_{\alpha}$ to $Q_1$, for $c:=c_0^{-1}a_1c_0\in\STGA_n(K)$ we have $c(K)=\theta$. We compute that
$c_0^{-1}=\e\bigl(x_1+\beta x_2-\beta h_0(x_1),x_2-h_0(x_1),x_3,\ldots,x_n\bigr)$, 
$a_1c_0=\e(x_1-F_{1,4},x_2-F_{1,5},x_3,\ldots,x_n)$, and $c:=\e(x_1-F_{1,6},x_2-F_{1,7},x_3,\ldots,x_n)$ with $F_{1,4}:=F_{1,3}+\beta\bigl(x_2-F_{1,2}+h_0(x_1-F_{1,3})\bigr)$, $F_{1,5}:=F_{1,2}-h_0(x_1-F_{1,3})$, 
$$F_{1,6}:=-\beta x_2+\beta h_0(x_1)+F_{1,4}\bigl(x_1+\beta x_2-\beta h_0(x_1),x_2-h_0(x_1)\bigr),$$
$$F_{1,7}:=h_0+F_{1,5}\bigl(x_1+\beta x_2-\beta h_0(x_1),x_2-h_0(x_1)\bigr).$$

By computing the degrees of $F_{1,6}$ and $F_{1,7}$ or drectly from the relations $\ell(c_0)=2$, $c=c_0^{-1}a_1c_0$, and $\ell(a_1)\le E_{n-2,k}$ and Inequality (\ref{EQ3}), we get the following inequalities $\ell(c)\le 4\ell(a_1)\le 4E_{n-2,k}$. So part (3) holds.
\end{proof}

\begin{example}\normalfont\label{EX7}
Suppose that $K=\mathbb F_2$ and $n\ge 3$. In $K^n$, let $O:=(0,\ldots,0)$, $O_1:=(1,0,\ldots,0)$, $O_2:=(0,1,0,\ldots,0)$, and $O_3:=(0,0,1,0,\ldots,0)$. For $i\in \llbracket1,n\rrbracket$, let $y_i:=x_i+1$. Let $a:=\e\bigl(x_1+\prod_{i=2}^n y_i,x_2,\ldots,x_n\bigr)\in\TGA_n(\mathbb F_2)[n-1]$ and 
$$b:=\e\Bigl(x_1,x_2,x_3+y_1y_2\prod_{i=4}^n y_i,x_4,\ldots,x_n\Bigr).$$ Then $a(\mathbb F_2)=(O\; O_1)$ and $b(\mathbb F_2)=(O\; O_3)$. Let $a_1:=ba\in\TGA_n(\mathbb F_2)$; we have $a_1(\mathbb F_2)=(O\; O_1\; O_3)$. As $n-1< 2n-3$ and 
$$a_1=\e\Bigl(x_1+\prod_{i=2}^n y_i,x_2,x_3+y_1y_2\prod_{i=4}^n y_i+y_2^2y_3\prod_{i=4}^n y_i^2,x_4,\ldots,x_n\Bigr),$$
we have $a_1\in\TGA_n(\mathbb F_2)[2n-3]$.\footnote{Note that in Proposition \ref{PR10}(1) applied to $p=2$, $n\ge 2$, and $s=1$ we also have an automorphism $a_1\in\TGA_n(\mathbb F_2)[2n-3]$ with $a_1(\mathbb F_2)$ a $3$-cycle.} Let $$c:=\e\Bigl(x_1,x_2,x_3+y_1x_2\prod_{i=4}^n y_i,x_4,\ldots,x_n\Bigr)\in\TGA_n(\mathbb F_2)[n-1].$$ Then $c(\mathbb F_2)=(O_2\; O_2+O_3)$. Let $a_2:=ca\in\TGA_n(\mathbb F_2)$; so $a_2(\mathbb F_2)$ is a product of two disjoint transpositions, $\supp\bigl(a_2(\mathbb F_2)\bigr)$ is affinely independent and we have $a_2\in\TGA_n(\mathbb F_2)[2n-3]$ as 
$$a_2=\e\Bigl(x_1+\prod_{i=2}^n y_i,x_2,x_3+y_1x_2\prod_{i=4}^n y_i+y_2x_2y_3\prod_{i=4}^n y_i^2,x_4,\ldots,x_n\Bigr).$$
Similarly, for $a_3:=\e(x_1+\prod_{i=3}^n y_i,x_2,\ldots,x_n)\in\TGA_n(\mathbb F_2)[n-2]$ we get that $a_3(\mathbb F_2)=(0\; O_1)(O_2\;O_1+O_2)$ is a product of two disjoint transpositions with $\supp\bigl(a_3(\mathbb F_2)\bigr)$ affinely dependent.
\end{example}

\begin{example}\normalfont\label{EX8}
Let $n\in \mathbb N^{\ast}\setminus\{1\}$ and $K$ a field with $|K|\ge 3$. Let $P_1$, $P_2$, and $P_3$ be three distinct and collinear points in $K^n$. 

\medskip
{\bf (1)} We check that there exists an automorphism $a\in\STGA_n(K)[2]$ such that $Q_1:=a(P_1)$, $Q_2:=a(P_2)$, and $Q_3:=a(P_3)$ are non-collinear. Up to linear automorphisms, we can assume that we have $P_1=(0,\ldots,0)$, $P_2=(1,0,\ldots,0)$, and $P_3=(\alpha,0,\ldots,0)$ with $\alpha\in K\setminus\{0,1\}$, so we can take $a:=\e(x_1,x_2+x_1^2,x_3,\ldots,x_n)$. 

Assume that $K$ is finite. If $l\in\mathbb N^{\ast}$ is such that the $3$-cycle $(Q_1\;Q_2\;Q_3)$ in $\Alt(K^n)$ is the image of an element $b\in\TGA_n(K)[l]$ (resp.\ $b\in\STGA_n(K)[l]$), then $(P_1\;P_2\;P_3)=a(K)^{-1}b(K)a(K)$ is the image of an element in $\TGA_n(K)[4l]$ (resp.\ in $\STGA_n(K)[4l]$). Similarly, if $s\in\mathbb N^{\ast}$ is such that $(P_1\;P_2\;P_3)\in\Alt(K^n)$ is the image of an element $c\in\TGA_n(K)[s]$ (resp.\ $c\in\STGA_n(K)[s]$), then $(Q_1\;Q_2\;Q_3)=a(K)c(K)a(K)^{-1}$ is the image of an element in $\TGA_n(K)[4s]$ (resp.\ in $\STGA_n(K)[4s]$). 

Assume that $K=\mathbb F_3$; so $\alpha=2$, $Q_1=(0,0)$, $Q_2=(1,1)$, and $Q_3=(2,1)$. Let $a_1:=\e(x_1,x_2+2x_1^2+2,x_3,\ldots,x_n)\in\STGA_n(\mathbb F_3)[2]$. We have $a_1(P_i)=Q_i$ for $i\in\{2,3\}$ and $Q_4:=a_3(P_1)=(0,2)$. Clearly, $Q_1+Q_4=Q_2+Q_3$. By Proposition \ref{PR10}(2) we can take $c\in\STGA_n(\mathbb F_3)[2n-2]$. Let $\sigma:=(Q_1\; Q_2)(Q_3\;Q_4)\in\Alt(\mathbb F_3^n)$. Let $a_2:=a^{-1}a_1=\e(x_1,x_2+x_1^2+2,x_3,\ldots,x_n)\in\STGA_n(\mathbb F_3)[2]$ and 
$$\sigma=(Q_1\;Q_2\;Q_3)(Q_2\;Q_3\;Q_4)=(aca^{-1}a_1ca_1^{-1})(\mathbb F_3)=(aca_2ca_1^{-1})(\mathbb F_3),$$
it follows that $\ell_{n,\mathbb F_3}(\sigma)\le 8(2n-2)^2$.

\smallskip
{\bf (2)} Suppose that $|K|\ge 4$ and that $P_4:=P_2+P_3-P_1\in K^n\setminus\{P_1,P_2,P_3\}$. Let $(Q_1,Q_2,Q_3,Q_4)\in\mathbb D_{n,4}(K)$ be non-collinear and such that $Q_1+Q_4=Q_2+Q_3$. We check that there exists $a\in\STGA_n(K)[3]$ such that $a(P_i)=Q_i$ for each $i\in \llbracket1,4\rrbracket$. Up to linear automorphisms we can assume that $n=2$, that $P_1$ to $P_3$ are as in part (1), and that $Q_1=P_1$, $Q_2=(1,1)$, and $Q_3=(\alpha,\alpha^2)$. So $P_4=(\alpha+1,0)$, $\alpha\neq -1$, and $Q_4=(\alpha+1,\alpha^2+1)$. Then we can take 
$$a:=\e\Bigl(x_1,x_2+x_1^2-\frac{2}{\alpha+1}x_1(x_1-1)(x_1-\alpha)\Bigr).$$ 

Assume that $K$ is finite. If $l\in\mathbb N^{\ast}$ is such that $(Q_1\;Q_2)(Q_3\;Q_4)\in\Alt(K^n)$ is the image of an element $b\in\TGA_n(K)[l]$ (resp.\ $b\in\STGA_n(K)[l]$), then the permutation $(P_1\;P_2)(P_3\;P_4)=a(K)^{-1}b(K)a(K)$ is the image of an element in $\TGA_n(K)[4l]$ (resp.\ in $\STGA_n(K)[4l]$). Similarly, if $s\in\mathbb N^{\ast}$ is such that $(P_1\;P_2)(P_3\;P_4)$ is the image of an element $c\in\TGA_n(K)[s]$ (resp.\ $c\in\STGA_n(K)[s]$), then $(Q_1\;Q_2)(Q_3\;Q_4)$ is the image of an element in $\TGA_n(K)[4s]$ (resp.\ in $\STGA_n(K)[4s]$). If $K=\mathbb F_4$, then we can take $s=3(n-1)$ by Lemma \ref{F5}(1) and (2). 
\end{example}

\section{Transitivity properties related to the actions $\mathbb T_{n,m}(K)$}\label{S9}

We have $\AGL_n(K)=\TGA_n(K)[1]$ and $\ASL_n(K)=\STGA_n(K)[1]$. 

\begin{notation}\normalfont\label{N3}
Let $(n,m)\in (\mathbb N^{\ast})^2$ and $K$ a field. For $i\in\llbracket1,\min(n,m-1)\rrbracket$ let
$$\mathbb D_{n,m}^{\d=i}(K):=\bigl\{(P_1,\ldots,P_m)\in\mathbb D_{n,m}(K)|\d_{\{P_1,\ldots,P_m\}}=i\bigr\},$$ 
for $i\in\llbracket1,\min(n,m-1)-1\rrbracket$ let
$$\mathbb D_{n,m}^{\d\ge i}(K):=\bigl\{(P_1,\ldots,P_m)\in\mathbb D_{n,m}(K)|\d_{\{P_1,\ldots,P_m\}}\ge i\bigr\},$$ 
and for $i\in\llbracket2,\min(n,m-1)\rrbracket$ let
$$\mathbb D_{n,m}^{\d\le i}(K):=\bigl\{(P_1,\ldots,P_m)\in\mathbb D_{n,m}(K)|\d_{\{P_1,\ldots,P_m\}}\le i\bigr\}.$$
\end{notation}

Recall $\varrho_{n,K}:\GA_n(K)\rightarrow\perm(K^n)$ is the representation that defines $\mathbb T_{n,m}(K)$; so for $a\in\GA_n(K)$, $\varrho_{n,K}(a)$ is the bijection $a(K):\mathbb A^n_K(K)\rightarrow \mathbb A^n_K(K)$, which, under the identification $\mathbb A^n_K(K)=K^n$, we view as a permutation $a(K):K^n\rightarrow K^n$. Hence the subgroups of $\GA_n(K)$ act via conjugation on (conjugacy classes of) $\perm(K)$.

We first recall the following two simple properties which are mostly well-known.

\begin{lemma}\label{F6}
Suppose that $K$ is a finite field. Let $a\in\AGL_n(K)\setminus\{1_{\mathbb A^n_K}\}$. Then the following properties hold.

\medskip
{\bf (1)} We have $\n\bigl(a(K)\bigr)\in\{|K|^n\}\cup\{|K|^n-|K|^l|l\in \llbracket0,n-1\rrbracket\}$. 

\smallskip
{\bf (2)} If there exists $l\in \llbracket0,n-1\rrbracket$ such that $\n\bigl(a(K)\bigr)=|K|^n-|K|^l$, then the complement $K^n\setminus\supp\bigl(a(K)\bigr)$ is an affine set of dimension $l$.

\smallskip
{\bf (3)} We have $|K|^n-|K|^{n-1}\le \n\bigl(a(K)\bigr)$.
\end{lemma}

\begin{proof}
For parts (1) and (2), we can assume that there exists $P\in K^n$ fixed by $a(K)$. By replacing $a$ by $b^{-1}ab$, where $b\in\ASL_n(K)$ is the translation by $P$, we can assume that $P=0$. So $a\in\GL_n(K)$ and thus $K^n\setminus\supp\bigl(a(K)\bigr)$ is a $K$-linear subspace $V$ of $K^n$. Let $l:=\dim_K(V)$; as $a\neq 1_{\mathbb A^n_K}$, we have $l\in\llbracket0,n-1\rrbracket$. Then $\n\bigl(a(K)\bigr)=|K|^n-|K|^l$. So parts (1) and (2) hold.

Part (3) follows from part (1).
\end{proof}

\begin{lemma}\label{F7} Let $(n,m)\in (\mathbb N^{\ast}\setminus\{1\})\times\mathbb N^{\ast}$ and $K$ a field be such that $\sqrt[n]{m}\le |K|$. Then the following properties hold.
 
\medskip
{\bf (1)} If $m\le n+1$ (resp.\ $m\le n$), then the group $\AGL_n(K)$ (resp.\ $\ASL_n(K)$) acts transitively on $\mathbb D_{n,m}^{\d=m-1}(K)$. In particular, $\AGL_n(K)$ acts transitively on $\mathbb D^{\d=2}_{n,3}(K)$ and $\ASL_n(K)$ acts $2$-transitively on $K^n$.

\smallskip
{\bf (2)} If $K$ is infinite (resp.\ finite with $K|\ge 3$), then the action of $\AGL_n(K)$ (or $\ASL_n(K)$) on $\mathbb D_{n,3}^{\d=1}(K)$ has infinitely many (resp.\ has precisely $|K|-2$) orbits. In particular, if $|K|=3$, then $\ASL_n(K)$ acts transitively on $\mathbb D_{n,3}^{\d=1}(K)$.

\smallskip
{\bf (3)} Suppose that $|K|=2$. Then the following properties hold.

\medskip\noindent
{\bf (3.a)} We have an identity $\AGL_n(K)=\ASL_n(K)$.

\smallskip\noindent
{\bf (3.b)} The group $\AGL_n(\mathbb F_2)$ acts transitively on $\mathbb D_{n,3}(\mathbb F_2)$.

\smallskip\noindent
{\bf (3.c)} The group $\AGL_n(\mathbb F_2)$ acts transitively on $\mathbb D^{\d=2}_{n,4}(\mathbb F_2)$. 

\smallskip\noindent
{\bf (3.d)} If $n\ge 3$, then $\AGL_n(\mathbb F_2)$ acts transitively on linearly dependent subsets of $\mathbb F_2^n$ with $5$ elements.
\end{lemma}

\begin{proof}
Clearly, part (1) holds.

Part (2) follows from the fact that the set $K\setminus\{0,1\}$ is in bijection to the set of orbits via the rule $\alpha\mapsto (P_0,P_1,P_{\alpha})$, where for $\beta\in K$, $P_\beta:=(\beta,0,\ldots,0)\in K^n$.

Part (3.a) follows from the fact that the group $K^{\ast}$ is trivial.

Part (3.b) follows from part (1) and the identity $\mathbb D_{n,3}(\mathbb F_2)=\mathbb D^{\d=2}_{n,3}(\mathbb F_2)$.

For part (3.c) we can assume that $n=2$, and in this case part (3.c) holds as $\GL_2(\mathbb F_2)\cong\perm\bigl(\mathbb F_2^2\setminus\{(0,0)\}\bigr)$.

For part (3.d), let $Y$ be a linearly dependent subset of $\mathbb F_2^n$ with $|Y|=5$. We have $\d_Y=3$, so up to affine automorphisms we can assume that $Y$ is contained in the zero locus $x_4=\cdots=x_n=0$. Let $Z\subset\mathbb F_2^3\times\{(0,\ldots,0)\}$ be a fixed subset with $|Z|=5$. Part (3.b) gives that there exists $a\in\AGL_3(\mathbb F_2)\times\{1_{\mathbb A^{n-3}_{\mathbb F_2}}\}$ with $a\bigl([\mathbb F_2^3\times\{(0,\ldots,0)\}]\setminus Y\bigr)=[\mathbb F_2^3\times\{(0,\ldots,0)\}]\setminus Z$. Thus $a(Y)=Z$. So part (3.d) holds.\end{proof}

Similar to the left action $\mathbb T_{n,m}$, we have a left action
$$\perm(K^n)\times\mathbb D_{n,m}(K)\rightarrow \mathbb D_{n,m}(K)$$
defined by the rule $(\sigma,\underline{P})\mapsto \sigma(\underline{P}):=\bigl(\sigma(P_1),\ldots,\sigma(P_m)\bigr)$.

\begin{definition}\label{D10.2}
Let $K$ be a finite field and $n\in\mathbb N^{\ast}$. We call 
$$\Perm(K^n):=\varrho_{n,K}\bigl(\TGA_n(K)\bigr)$$ the group of tame permutations\index{group of tame permutations} of $K^n$.
\end{definition}

\begin{theorem}\label{T6}
Suppose that $K$ is finite and $n\in\mathbb N^{\ast}\setminus\{1\}$. Then the following properties hold.

\medskip
{\bf (1)} If $|K|$ is odd (resp.\ $K\cong\mathbb F_2$), then there exists $a\in\AGL_n(K)\setminus\STGA_n(K)$ (resp.\ $a\in \STGA_n(K)[n-1]$) such that the permutation $a(K)\in\perm(K^n)$ is odd with $\n\bigl(a(K)\bigr)=|K|^n-|K|^{n-1}$ (resp.\ is a transposition).

\smallskip
{\bf (2)} If either $|K|$ is odd or $K\cong\mathbb F_2$, then $\Perm(K^n)=\perm(K^n)$.

\smallskip
{\bf (3)} If $4\mid |K|$, then $\Perm(K^n)=\Alt(K^n)$.

\smallskip
{\bf (4)} If $|K|\ge 3$, then we have $\varrho_{n,K}\bigl(\STGA_n(K)\bigr)=\varrho_{n,K}\bigl(\SGA_n(K)\bigr)=\Alt(K^n)$. 
\end{theorem}

\begin{proof}
If $|K|$ is odd, let $\gamma$ be a generator of the multiplicative cyclic group $K^{\ast}$ and $a:=\e(\gamma x_1,x_2,\ldots,x_n)$; then $a(K)$ is a product of $|K|^{n-1}$ disjoint $(|K|-1)$-cycles and hence $a(K)$ is an odd permutation with $\bigl|\n\bigl(a(K)\bigr)\bigr|=|K|^n-|K|^{n-1}$. 

If $|K|=2$, then let $a$ be $\e(x_1+x_2,x_1)$ if $n=2$ and be as in Example \ref{EX7} if $n\ge 3$; so $a(K)$ is a transposition. Thus part (1) holds.

For part (2), we first recall that $\varrho_{n,K}\bigl(\STGA_n(K)\bigr)$ and hence also $\Perm(K^n)$ contains a $3$-cycle of non-collinear points by Proposition \ref{PR10}(1).

Based on Lemma \ref{F7}(1) and $\AGL_n(K)\leqslant\TGA_n(K)$, the action of $\Perm(K^n)$ on $K^n$ is primitive. So, as $\Perm(K^n)$ contains a $3$-cycle of non-collinear points and $|K|^n\ge 6$ except when $|K|=n=2$, Jordan's Theorem (e.g., see either the references for \cite{Jo}, Thm.\ 1.1 or \cite{Jo}, Thm.\ 1.2) gives that $\Alt(K^n)\leqslant\Perm(K^n)$ except when $|K|=n=2$. If $|K|=n=2$, then $\AGL_2(K)=\perm(K^2)$, so $\Perm(K^2)=\perm(K^2)$. Thus $\Alt(K^n)\lhd\Perm(K^n)$. 

We include a second proof that $\Alt(K^n)\lhd\Perm(K^n)$. It suffices to show that $\Alt(K^n)\lhd\varrho_{n,K}\bigl(\STGA_n(K)\bigr)$. Let $P$, $Q$, and $O$ be distinct points in $K^n$. Based on Lemma \ref{F7}(1) and the fact that $\varrho_{n,K}\bigl(\STGA_n(K)\bigr)$ contains a $3$-cycle of non-collinear points, we have $(P\; Q\; O)\in\varrho_{n,K}\bigl(\STGA_n(K)\bigr)$ if $P$, $Q$, and $O$ are non-collinear. If $P$, $Q$, and $O$ are collinear, then we have $(P\; Q\; O)\in\varrho_{n,K}\bigl(\STGA_n(K)\bigr)$ by Example \ref{EX8}(1). As $\Alt(K^n)$ is generated by $3$-cycles, it follows that $\Alt(K^n)$ is a subgroup of $\varrho_{n,K}\bigl(\STGA_n(K)\bigr)$.

Part (2) follows from part (1) and $\Alt(K^n)\lhd\Perm(K^n)$. 

If $4\mid |K|$, then, as $\TGA_n(K)$ is generated by subgroups of $K$-valued points of connected linear groups, more precise by $\AGL_n(K)$ and the $\mathbb G_{\a,K}(K)$ subgroups that contain $\e(x_1,\ldots,x_{n-1},x_n+f)$ with $f\in K[x_1,\ldots,x_{n-1}]$, we have $\Perm(K^n)\leqslant\Alt(K^n)$ by Theorem \ref{T3}. So part (3) holds.

For part (4), as $\Alt(K^n)\lhd\varrho_{n,K}\bigl(\STGA_n(K)\bigr)$ and $\STGA_n(K)\leqslant\SGA_n(K)$, it suffices to prove that $\varrho_{n,K}\bigl(\SGA_n(K)\bigr)\leqslant\Alt(K^n)$. As $\SGA_n(K)$ is generated by $\mathbb G_{\a,K}(K)=K$ subgroups defined by faithful actions $\mathbb G_{\a,K}\times_{\Spec K}\mathbb A^n_K\rightarrow\mathbb A^n_K$, it suffices to show that the image of each such subgroup under $\varrho_{n,K}$ is contained in $\Alt(K^n)$ but this is a particular case of Theorem \ref{T3}. Thus part (4) holds.\end{proof}

\subsection{Proof of Theorem \ref{T1}}\label{S9.1}
Let $m\in\mathbb N^{\ast}$. As \cite{Sr}, Thm.\ 1.2 implies that infinite fields are $T_{n,m}$ fields and hence infinitely transitive and as the proof of \cite{KZ}, Lem.\ 5.5 can be adapted to prove that they are also $ST_{n,m}$ fields, to prove Theorem \ref{T1} we can assume that $|K|\in\mathbb N^{\ast}$. 

Let $(\underline{P},\underline{Q})\in\mathbb D_{n,m}(K)^2$. Let $\sigma\in\perm(K^n)$ be such that $\sigma(\underline{P})=\underline{Q}$. 

If $m\le |K|^n-2$, then we can choose $\sigma\in\Alt(K^n)$ and therefore, as we have $\Alt(K^n)\lhd\varrho_{n,K}\bigl(\STGA_n(K)\bigr)$ by Theorem \ref{T6}(4) if $|K|\ge 3$ and by Theorem \ref{T6}(2) and $\STGA_n(\mathbb F_2)=\TGA_n(\mathbb F_2)$ if $|K|=2$, $K$ is an $ST_{n,m}$ field and thus a $T_{n,m}$ field. 

If $m\in\{|K|^n-1,|K|^n\}$, then $\sigma$ is uniquely determined and is an arbitrary permutation in $\perm(K^n)$; from this and Theorem \ref{T6}(2) to (4) we get first that $K$ is not an $ST_{n,m}$ field iff $|K|\ge 3$ and second that $K$ is not a $T_{n,m}$ field iff $4\mid |K|$. 

The last two paragraphs imply that Theorem \ref{T1}(1) and (2) holds. 

Theorem \ref{T1}(3) follows from Theorem \ref{T1}(2) as a finite field $K$ is infinitely transitive iff it is a $T_{n,m}$ field for each $m\in \llbracket1,|K|^n\rrbracket$. So Theorem \ref{T1} holds.

\subsection{Proof of Corollary \ref{C1}}\label{S9.2}
Suppose that $|K|\in\mathbb N^{\ast}$. Let $\sigma\in\perm(K^n)$ be such that $\sigma(Y)=Z$. As $|K|^n\ge 4$, we have $|Z|\ge 2$ or $|K|^n\setminus Z|\ge 2$. Thus by replacing $\sigma$ by $\theta\sigma$ with $\theta\in\perm(K^n)$ such that $\theta(Z)=Z$, we can assume that $\sigma\in\Alt(K^n)$. So the Corollary \ref{C1}(1) follows from Theorem \ref{T6}(2) and (4) if $|K|\in\mathbb N^{\ast}$. 

If $|K|=\infty$, then $K$ is an $ST_{n,m}$ field and thus there exists $a\in\STGA_n(K)$ such that $a(Y)=Z$. So Corollary \ref{C1}(1) holds. 

Corollary \ref{C1}(2) follows from Corollary \ref{C1}(1). 

Corollary \ref{C1}(3) follows from Corollary \ref{C1}(2). Thus Corollary \ref{C1} holds. 

\section{Regular length functions}\label{S10}

We first recall some basic language on length functions in the multiplicative context with codomain $\mathbb N^{\ast}$.

\begin{definition}\label{D11}
Let $\Gamma$ be an abstract group. Let $\mathcal L:\Gamma\rightarrow\mathbb N^{\ast}$ be a function.

\medskip
{\bf (1)} We say that $\mathcal L$ is a length function\index{length function} (on $\Gamma$) if the following three axioms hold. 

\medskip
{\bf A1.} The identity element of $\Gamma$ is mapped by $\mathcal L$ to $1$.

\medskip
{\bf A2.} For each $x\in\Gamma$ we have $\mathcal L(x)=\mathcal L(x^{-1})$.

\smallskip
{\bf A3.} For each pair $(x,y)\in\Gamma^2$ we have $\mathcal L(xy)\le\mathcal L(x)\mathcal L(y)$.

\medskip
{\bf (2)} If $\mathcal L$ is a length function, then we say that it is L-regular\index{length function!L-regular} if for each triple $(x,y,z)\in\Gamma^3$, the strict inequality $\mathcal L(y)\mathcal L(x^{-1}z)<\mathcal L(z)\mathcal L(x^{-1}y)$ implies the identity $\mathcal L(z)\mathcal L(x^{-1}y)=\mathcal L(x)\mathcal L(y^{-1}z)$.

\smallskip
{\bf (3)} If $\mathcal L$ is a length function, then we say that it is regular\index{length function!regular} if for each triple $(x,y,z)\in\Gamma^3$, the strict inequalities $\mathcal L(xy)<\mathcal L(x)\mathcal L(y)$ and $\mathcal L(y^{-1}z)<\mathcal L(y)\mathcal L(z)$ imply the inequality $\mathcal L(xz)<\mathcal L(x)\mathcal L(z)$.
\end{definition}

\begin{remark}\normalfont\label{R2}
Let $\mathcal L:\Gamma\rightarrow\mathbb N^{\ast}$ be a length function. 

\medskip
{\bf (1)} The kernel $\Ker(\mathcal L):=\{x\in\Gamma|\mathcal L(x)=1\}$
of $\mathcal L$ is a subgroup of $\Gamma$ which in general is not normal. If $(x,y,z)\in\Ker(\mathcal L)^2\times\Gamma$, then we use the following identity 
$$\mathcal L(xzy)=\mathcal L(z)$$ 
(e.g., by \cite{Pro}, Lem.\ 2.1(c)) without any extra comment. 

\smallskip
{\bf (2)} Note that $\ln\mathcal L:\Gamma\rightarrow [0,\infty)$ is a semigauge in the terminology of \cite{Pro}, Sect.\ 2 and following \cite{L}, Sect.\ 2 one considers the function $\mathcal H:\Gamma^2\rightarrow [0,\infty)$ that measures `the deviation from equality' in Axiom A3 defined by
$$\mathcal H(x,y):=\frac{1}{2}\Bigl[\ln\mathcal L(x)+\ln\mathcal L(y)-\ln\mathcal L(x^{-1}y)\Bigr].$$ 

{\bf (3)} The length function $\mathcal L$ is regular (resp.\ L-regular) iff $\ln\mathcal L:\Gamma\rightarrow [0,\infty)$ is regular in the sense of \cite{Pro}, Sect.\ 4, Def.\ (resp.\ is such that Lyndon's Axiom A4 of \cite{L}, Sect.\ 2 holds for it, i.e., $[\mathcal H(x,y)<\mathcal H(x,z)]\Rightarrow [\mathcal H(y,z)=\mathcal H(x,y)]$). 

{\bf (4)} If $\mathcal L$ is $L$-regular, then it is regular by \cite{Pro}, Ex.\ 4.2(2) applied to the semigauge $\ln\mathcal L$. 

{\bf (5)} Basic properties of regular length function are summarized in \cite{Pro}, Lems.\ 4.3, 4.5, 4.6, and 4.7. In particular, the uniqueness of reduced product decompositions into atoms up to a suitable equivalence relation holds in $\Gamma$ (see \cite{Pro}, Lem.\ 4.3(e); cf.\ the terminology of \cite{Pro}, Sect.\ 2).
\end{remark}

In order to study the length function $\ell_{\GA_2(K)}$ we first reformulate \cite{vdK}, Thms.\ 1 and 2 via the following definition. 

\begin{definition}\label{D12}
For $a\in \GA_2(K)\setminus\AGL_2(K)$ let 
$$\NAT_a:=\bigl\{(f,c)\in K[x]\times \GL_2(K)|\deg(f)\ge 2\;\;\textup{and}\;\;a=c\e\bigl(x_1,x_2+f(x_1)\bigr)c^{-1}\bigr\}$$
and
$$\RNAT_a:=\bigl\{(f,c)\in x^2K[x]\times \GL_2(K)|\deg(f)\ge 2\;\;\textup{and}\;\;a=c\e\bigl(x_1,x_2+f(x_1)\bigr)c^{-1}\bigr\}$$

{\bf (1)} We say that $a$ is non-affine triangular\index{automorphism!non-affine triangular} if it belongs to the set 
$$\NAT(K):=\{a\in\GA_2(K)\setminus\AGL_2(K)|\NAT_a\neq\emptyset\}$$
and we say that $a$ is reduced non-affine triangular\index{automorphism!reduced non-affine triangular} if it belongs to the set 
$$\RNAT(K):=\{a\in\GA_2(K)\setminus\AGL_2(K)|\RNAT_a\neq\emptyset.\}$$

{\bf (2)} If $(f,c)\in\NAT_a$, then we say that $f$ or $(f,c)$ defines $a$.

\smallskip
{\bf (3)} If $a$ is non-affine triangular, then by its direction\index{direction} we mean the $K$-valued point $\dir(a)\in\mathbb P^1(K)$ that defines the one-dimensional vector space over an algebraic closure $\overline{K}$ of K spanned by all $a\bigl((\alpha,\beta)\bigr)-(\alpha,\beta)$ with $(\alpha,\beta)\in \overline{K}^2$.
\end{definition}

\begin{remark}\normalfont\label{R3} {\bf (1)} For $a\in\GA_2(K)$, $a$ is non-affine triangular (resp.\ reduced non-affine triangular) iff $a^{-1}$ is so. Moreover, if $a=c\e\bigl(x_1,x_2+f(x_1)\bigr)c^{-1}$ is non-affine triangular, then $\dir(a)=\dir(a^{-1})$ is the $K$-valued point of $\mathbb P^1_K$ defined by $c\bigl((0,1)\bigr)$.

\smallskip
{\bf (2)} Let $B_2(K)$ be the subgroup of $\GL_2(K)$ of invertible upper triangular matrices. For $a\in\NAT(K)$, the set $\NAT_a$ is a trivial right torsor under $B_2(K)$ via the right action $\NAT_a\times\B_2(K)\rightarrow\NAT_a$ defined by the rule 
$$\left((f,c),\begin{bmatrix} 
\beta & \gamma\\
0 & \delta\\
\end{bmatrix}\right)\mapsto \bigl(\delta^{-1} f(\beta x_1),c\e(\beta x_1,\gamma x_1+\delta x_2)\bigr);$$
here $(\beta,\gamma,\delta)\in K^{\ast}\times K\times K^{\ast}$. This is so as, with $b:=\e\bigl(x_1,x_2+f(x_1)\bigr)$, the identity $cbc^{-1}=\tilde{c}[(\tilde{c}^{-1}c)b(\tilde{c}^{-1}c)^{-1}]\tilde{c}^{-1}$ for $\tilde c:=c\e(\beta x_1,\delta x_2+\gamma x_1)$ gives 
$$(\tilde{c}^{-1}c)b(\tilde{c}^{-1}c)^{-1}=\e\bigl(x_1,x_2+\delta^{-1} f(\beta x_1)\bigr).$$
In particular, if $\RNAT_a\neq\emptyset$, then $\RNAT_a=\NAT_a$. 

{\bf (3)} For each $(\alpha,\beta)\in K^2$ and $(b,c)\in\AGL_2(K)\times\GL_2(K)$, we have an identity 
$$c\e\bigl(x_1,x_2+f(x_1)\bigr)c^{-1}c\e\bigl(x_1,x_2+\alpha x_1+\beta\bigr)c^{-1}b=c\e\bigl(x_1,x_2+f(x_1)+\alpha x_1+\beta\bigr)c^{-1}b$$ 
and thus for each $d\in\mathbb N^{\ast}$ we have an identity of subsets
$$\bigl\{c\e\bigl(x_1,x_2+f(x_1)\bigr)|\deg(f)\le d\bigr\}\AGL_2(K)=\bigl\{c\e\bigl(x_1,x_2+f(x_1)\bigr)c^{-1}|\deg(f)\le d\bigr\}\Gamma_c$$
of $\GA_2(K)$, where $\Gamma_c$ is a subset of $\AGL_2(K)$ of representatives of right cosets of the subgroup $c\bigl\{\e\bigl(x_1,x_2+\alpha x_1+\beta\bigr)|(\alpha,\beta)^2\in K^2\bigr\}c^{-1}$
of $\AGL_2(K)$.
\end{remark}

\begin{definition}\label{D13}
Let $K$ be a field and $a\in\GA_2(K)=\TGA_2(K)$. 

\medskip
{\bf (1)} By the triangular index set\index{triangular index set} of $a$ we mean the non-empty subset $\TIS_a$ of $\mathbb N$ formed by integers $j\in\mathbb N$ for which there exists $(c_1,c_2,\ldots,c_{j+1})\in \AGL_2(K)^{j+1}$ and $(f_1,\ldots,f_j)\in K[x]^j$ with $\deg(f_i)\ge 2$ for each $i\in \llbracket1,j\rrbracket$ such that by denoting $b_i:=\e\bigl(x_1,x_2+f_i(x_1)\bigr)$ for each $i\in \llbracket1,j\rrbracket$ we have $a=c_1b_1c_2b_2\cdots c_jb_jc_{j+1}$. 

\smallskip
{\bf (2)} By the Furter's length\index{Furter's length} $\j_a$ of $a$ we mean $\min(\TIS_a)$.\footnote{See \cite{F1}, Sect.\ 4, Thm.\ and Def. and \cite{F2}, Sect.\ 1 c.; strictly speaking, \cite{F1}, Sect.\ 3 onwards assumes the field is $\mathbb C$ but \cite{F1}, Sect.\ 4, Thm.\ holds over each field cf.\ also \cite{F2}, Sect.\ 1 c.}

\smallskip
{\bf (3)} We use the notation of part (1) for $j=\j_a$. Replacing $c_i$ by their linear parts and conjugating $b_i$ by suitable translations, we can assume that $c_i\in \GL_2(K)$ for all $i\in\llbracket1,j\rrbracket,$ and only $c_{j+1}\in \AGL_2(K)$. By denoting $c:=\prod_{i=1}^{j+1} c_i$ and $a_l:=(\prod_{i=1}^l c_i)b_l(\prod_{i=1}^l c_i)^{-1}\in\SGA_2(K)$ for each $l\in \llbracket1,j\rrbracket$, we get a product decomposition 
$$a=a_1\cdots a_jc$$
with $a_i$ as a non-affine triangular automorphism for every $i\in \llbracket1,j\rrbracket$ and $\break c\in\AGL_2(K)$, which we call a standard product decomposition\index{standard product decomposition} of $a$. 

Moreover, if after making $c_1$ to be in $\GL_2(K)$ we consider a sum decomposition $f_1=g_1+h_1$ with $h_1\in x^2K[x]$ and $g_1\in K[x]$ of degree at most $1$ and replace $(b_1,c_1)$ by $\bigl(\e(x_1,h_1(x_1)),\e(x_1,x_2+g_1(x_1))c_1\bigr)$, we can assume that $c_1\in\GL_2(K)$ and $f_1\in x^2K[x]$. Repeating the process for $i\in\llbracket2,j\rrbracket$, we get a product decomposition 
$$a=a_1\cdots a_jc$$
with $a_i$ as a reduced non-affine triangular automorphism for every $i\in \llbracket1,j\rrbracket$ and $c\in\AGL_2(K)$, which we call a reduced standard product decomposition\index{standard product decomposition!reduced} of $a$. 
\end{definition}

\phantomsection{We consider the following subset $\mathfrak C:=\{\e(\alpha x_1+x_2,x_1)|\alpha\in K\}$ of $\GA_2(K)$.}\label{PH86}

\begin{proposition}\label{PR11}
Let $K$ be a field and $a\in\GA_2(K)$. Let $a=a_1\cdots a_jc$ be a standard product decomposition of $a$ with $j\in\mathbb N$, each $a_i$ with $i\in \llbracket1,j\rrbracket$ as a non-affine triangular automorphism, and $c\in\AGL_2(K)$. For each $i\in \llbracket1,j\rrbracket$ we consider a pair $(f_i,c_i)\in K[X_1]\times\GL_2(K)$ that defines $a_i$. Then the following properties hold.

\medskip
{\bf (1)} If $a=a_1\cdots a_jc$ is a standard product decomposition of $a$, so $j=\j_a$, then we have $\dir(a_{i-1})\neq \dir(a_i)$ for each $i\in\llbracket2,j\rrbracket$.

\smallskip
{\bf (2)} For each $i\in \llbracket1,j\rrbracket$ we have $\ell(a_i)=\pi(a_i)=\deg(f_i)\ge 2$ and we can assume that $c_i\in\{1_{\mathbb A^2_K}\}\cup\mathfrak C$.

\smallskip
{\bf (3)} If $\dir(a_{i-1})\neq \dir(a_i)$ for each $i\in\llbracket2,j\rrbracket$ and $c_1\notin B_2(K)$, then we have a product decomposition $$\prod_{i=1}^j a_i=\Bigl[\prod_{i=1}^j \tilde c_i\e\bigl((x_1,x_2+g_i(x_1)\bigr)\Bigr]\tilde c_{j+1}$$
with $\tilde c_{j+1}\in\GL_2(K)$ and for each $i\in \llbracket1,j\rrbracket$ with $\deg(g_i)=\deg(f_i)$ and $\tilde c_i\in\mathfrak C$.

\smallskip
{\bf (4)} Suppose that $\dir(a_{i-1})\neq \dir(a_i)$ for each $i\in\llbracket2,j\rrbracket$. Then the following properties hold.

\medskip\noindent
{\bf (4.a)} We have an identity
\begin{equation}\label{EQ14}
\ell(a)=\ell(a_1\cdots a_j)=\pi(a_1\cdots a_j)=\prod_{i=1}^j \deg(f_i)=\prod_{i=1}^j \ell(a_i).
\end{equation}

\noindent
{\bf (4.b)} If $\dir(a_s)\neq\dir(ca_1c^{-1})$, then $\ell(a^l)=\ell(a)^l$ for each $l\in\mathbb N^{\ast}$.
\end{proposition}

\begin{proof}
For part (1), we have $c_{i-1}\bigl((0,1)\bigr)\neq c_i\bigl((0,1)\bigr)$ as otherwise we can replace $j$ by $j-1$ and $j$ would not be $\j_a$. From this and Remark \ref{R3}(1) we get that part (1) holds.

The identities of part (2) follow from very definitions. As 
$$\GL_2(K)\setminus B_2(K)=\bigsqcup_{\alpha\in K} \begin{bmatrix} 
\alpha & 1\\
1 & 0\\
\end{bmatrix}B_2(K),$$ either $c_i\in B_2(K)$, or there exists $\alpha_i\in K$ such that $c_i\in\begin{bmatrix} 
\alpha & 1\\
1 & 0\\
\end{bmatrix}B_2(K)$. Based on this and Remark \ref{R3}(2) we can assume that, as automorphisms in $\GA_2(K)$, either $c_i=1_{\mathbb A^2_K}$ or $c_i\in\mathfrak C$. So part (2) holds.

We prove part (3) by induction on $j$. The case $j=0$ is clear and the case $j=1$ follows from part (2) applied to $i=1$. For $j\ge 2$, for the inductive step from $j-1$ to $j$, we can assume that $c_1\in\mathfrak C$. As $\dir(a_1)\neq \dir(a_2)$, we have $c_1^{-1}c_2\notin B_2(K)$. Based on this and the inductive assumption, we can rewrite the product $c_1^{-1}\prod_{i=2}^j a_i=\bigl[\prod_{i=2}^j [c_1^{-1}c_i]\e\bigl(x_1,x_2+f_i(x_1)\bigr)[c_1^{-1}c_i]^{-1}\bigr]c_1^{-1}$ as a product $[\prod_{i=2}^j \tilde c_i\e\bigl(x_1,x_2+g_i(x_1)\bigr)]\tilde c_{j+1}$
with $\tilde c_{j+1}\in\GL_2(K)$ and for each $i\in \llbracket2,j\rrbracket$ with $\deg(g_i)=\deg(f_i)$ and $\tilde c_i\in\mathfrak C$. Therefore
$$\prod_{i=1}^j a_i=c_1\e\bigl(x_1,x_2+f_1(x_1)\bigr)c_1^{-1}\prod_{i=2}^j a_i=\Bigl[\prod_{i=1}^j \tilde c_i\e\bigl(x_1,x_2+g_i(x_1)\bigr)\Bigr]\tilde c_{j+1}$$
with $\tilde c_1:=c_1$ and $g_1:=f_1$. This ends the induction and thus the proof of part (3).

The first two identities of part (4.a) are clear as $n=2$ and the last identity of part (4.a) follows from part (2). 

The third identity of part (4.a) follows from \cite{F3}, Prop.\ 2. We include a second proof of it based on the earlier result \cite{W}, Prop.\ 1.9. To prove the identity $\pi(a_1\cdots a_j)=\prod_{i=1}^j \deg(f_i)$, when $c_1=1_{\mathbb A^2_K}$ (resp.\ when $c_1\notin B_2(K)$), we rewrite $a_1\cdots a_j$ as a product $[\prod_{i=1}^j \tilde c_i\e\bigl(x_1,x_2+g_i(x_1)\bigr)]\tilde c_{j+1}$ with $\tilde c_{j+1}\in\GL_2(K)$, with $\deg(g_i)=\deg(f_i)$ and $\tilde c_i\in\mathfrak C$ for each $i\in \llbracket2,j\rrbracket$ , and moreover $(g_1,\tilde c_1)=(f_1,1_{\mathbb A^2_K})$ (resp.\ and $\deg(g_1)=\deg(f_1)$ and $\tilde c_1\in\mathfrak C$) by part (3) applied to $a_2\cdots a_j$ (resp.\ to $a_1\cdots a_j$). From this and \cite{W}, Prop.\ 1.9 we get directly the following identities $\pi(a_1\cdots a_j)=\prod_{i=1}^j \deg(g_i)=\prod_{i=1}^j \deg(f_i)$. So part (4.a) holds.

For part (4.b), for $i\in\{1,\ldots,j\}$ and $t\in\llbracket1,l-1\rrbracket$ let $a_{i+jt}:=c^ta_ic^{-t}$. We have $a^l=(\prod_{i=1}^{jl} a_i)c^l$ and $\dir(a_i)\neq\dir(a_{i+1})$ for each $i\in\llbracket1,jl-1\rrbracket$. Thus $\ell(a^l)$ is equal to $\prod_{i=1}^{jl}\deg(a_i)=[\prod_{i=1}^j\deg(a_i)]^l=\ell(a)^l$ by part (4.a). So part (4.b) holds.
\end{proof}

One can check that \cite{vdK}, Thms.\ 1 and 2 are equivalent to the existence of (reduced) standard product decompositions for which Equation (\ref{EQ14}) holds and hence the latter form the mentioned reformulation. Note that Proposition \ref{PR11}(4.b) is equivalent to \cite{F3}, Prop.\ 3.

\begin{proposition}\label{PR12}
For each field $K$, the length function $\ell_{\GA_2(K)}$ is regular. 
\end{proposition}

\begin{proof}
Let $(a_x,a_y,a_z)\in\GA_2(K)^3$ be a triple such that $\ell(a_xa_y)<\ell(a_x)\ell(a_y)$ and $\ell(a_y^{-1}a_z)<\ell(a_y)\ell(a_z)$; so $(a_x,a_y,a_x)\in [\GA_2(K)\setminus\AGL_2(K)]^3$. We consider standard product decomposition $a_x^{-1}=(\prod_{i=1}^{j_x} a_{i,x})b_x$, $a_y=(\prod_{i=1}^{j_y} a_{i,y})b_y$, and $a_z=(\prod_{i=1}^{j_z} a_{i,z})b_z$; so $(j_x,j_y,j_z):=(\j_{a_x},\j_{a_y},\j_{a_z})$, $(b_x,b_y,b_z)\in\AGL_2(K)^3$, and for each $\star\in\{x,y,z\}$ and every $i\in \llbracket1,j_{\star}\rrbracket$, the automorphism $a_{i,\star}\in\GA_2(K)\setminus\AGL_2(K)$ is non-affine triangular with $\ell(a_{i,\star})\ge 2$. 

As $\ell(a_xa_y)<\ell(a_x)\ell(a_y)$ and $a_xa_y=b_x^{-1}(\prod_{i=j_x}^{1} a_{i,x}^{-1})(\prod_{i=1}^{j_y} a_{i,y})b_y$ we get that $\ell\bigl((\prod_{i=j_x}^{1} a_{i,x}^{-1})(\prod_{i=1}^{j_y} a_{i,y})\bigr)<\ell(\prod_{i=j_x}^{1} a_{i,x}^{-1})\ell(\prod_{i=1}^{j_y} a_{i,y})$. If $\dir(a_{1,x}^{-1})\neq\dir(a_{1,y})$, then from Proposition \ref{PR11}(1) and (4.a) we get that $\ell\bigl((\prod_{i=j_x}^{1} a_{i,x}^{-1})(\prod_{i=1}^{j_y} a_{i,y})\bigr)$ is equal to $\ell(\prod_{i=j_x}^{1} a_{i,x}^{-1})\ell(\prod_{i=1}^{j_y} a_{i,y})$, a contradiction. So $\dir(a_{1,x}^{-1})=\dir(a_{1,y})$. We similarly argue, based on $a_y^{-1}a_z=b_y^{-1} (\prod_{i=j_y}^{1} a_{i,y}^{-1})(\prod_{i=1}^{j_z} a_{i,z})b_z$ and the inequality $\ell(a_y^{-1}a_z)<\ell(a_y)\ell(a_z)$, that we have $\dir(a_{1,y}^{-1})=\dir(a_{1,z})$. As $\dir(a_{1,y})=\dir(a_{1,y}^{-1})$ by Remark \ref{R3}(1), by transitivity we get that $\dir(a_{1,x}^{-1})=\dir(a_{1,z})$. 

As $\dir(a_{1,x}^{-1})=\dir(a_{1,z})$, $\ell(a_{1,x}^{-1})\ge 2$, and $\ell(a_{1,z})\ge 2$, we have inequalities $\ell(a_{1,x}^{-1}a_{1,z})\le\max\bigl(\ell(a_{1,x}^{-1}),\ell(a_{1,z})\bigr)<\ell(a_{1,x}^{-1})\ell(a_{1,z})$. So, as we have a product decomposition $a_xa_z=b_x^{-1}(\prod_{i=j_x}^2 a_{i,x}^{-1})(a_{1,x}^{-1}a_{1,z})(\prod_{i=2}^{j_z} a_{i,z})b_z,$ based on the identity $\ell(a)=\ell(a^{-1})$ with $a\in\GA_2(K)$ we estimate
$$\ell(a_xa_z)\le \Bigl[\prod^2_{i=j_x} \ell(a_{i,x})\Bigr]\ell(a_{1,x}^{-1}a_{1,z})\Bigl[\prod_{i=2}^{j_z} \ell(a_{i,z})\Bigr]$$
$$<\Bigl[\prod^2_{i=j_x} \ell(a_{i,x})\Bigr]\ell(a_{1,x})\ell(a_{1,z})\Bigl[\prod_{i=2}^{j_z} \ell(a_{i,z})\Bigr]=\Bigl[\prod^1_{i=j_x} \ell(a_{i,x})][\prod_{i=1}^{j_z} \ell(a_{i,z})\Bigr]=\ell(a_x)\ell(a_z),$$
where the last equality follows from Equation (\ref{EQ14}) applied to $a_x^{-1}$ and $a_z$. As $\ell(a_xa_z)<\ell(a_x)\ell(a_z)$, by the very definition we get that $\ell_{\GA_2(K)}$ is regular.\end{proof}

\begin{lemma}\label{F8}
Let $K$ be a field. If $n\in\mathbb N^{\ast}\setminus\{1,2\}$, then the length function $\ell_{\GA_n(K)}$ is not regular. 
\end{lemma}

\begin{proof}
It suffices to prove this for $n=3$. For a triple $(f,g,h)\in K[x]^3$ such that $\min\bigl(\deg(f),\deg(g),\deg(h)\bigr)\ge 2$, let $(a,b,c)\in\STGA_n(K)^3$ be defined by $\break a:=\e\bigl(x_1,x_2,x_3+f(x_2)\bigr)$, $b:=\e\bigl(x_1,x_2,x_3+g(x_1)\bigr)$, and $c:=\e\bigl(x_1,x_2+h(x_1),x_3)\bigr)$. Hence $ab=\e\bigl(x_1,x_2,x_3+g(x_1)+f(x_2)\bigr)$, $b^{-1}c=\e\bigl(x_1,x_2+h(x_1),x_3-g(x_1)\bigr)$, $ac=\e\bigl(x_1,x_2+h(x_1),x_3+f(x_2+h(x_1))\bigr)$, and $c^{-1}a^{-1}=\e\bigl(x_1,x_2-h(x_1),x_3-f(x_2)\bigr)$. So $\ell(a)=\deg(f)$, $\ell(b)=\deg(g)$, $\ell(c)=\deg(h)$, $\ell(ab)=\max\bigl(\deg(f),\deg(g)\bigr)$, $\ell(b^{-1}c)=\max\bigl(\deg(g),\deg(h)\bigr)$, and $\ell(ac)=\deg(f)\deg(h)$. Hence, as we have $\ell(ab)<\ell(a)\ell(b)$ and $\ell(b^{-1}c)<\ell(b)\ell(c)$ while $\ell(ac)=\ell(a)\ell(c)$, the lemma holds.
\end{proof}

\section{Length functions on permutation groups}\label{S11}

In this section and the next two sections we assume that $K$ is a finite field. We have the following first variant of Definition \ref{D2}(2) and (3).

\begin{definition}\label{D14}
Let $(n,m)\in (\mathbb N^{\ast})^2$ and the finite field $K$ be such that $|K|\ge\sqrt[n]{m}$. 

\medskip
{\bf (1)} For $(\underline{P},\underline{Q})\in\mathbb D_{n,m}(K)^2\in\mathbb N^{\ast}\cup\{\infty\}$, let $\pi^{\E}_{\underline{P},\underline{Q}}$ be the infimum of the set $\{\ell(a)|a\in\TGA_n(K), a(\underline{P})=\underline{Q}, a(K)\in\Alt(K^n)\}$.

\smallskip
{\bf (2)} Let $\pi^{\E}_{n,m}(K):=\sup\bigl(\{\pi^{\E}_{\underline{P},\underline{Q}}|(\underline{P},\underline{Q})\in\mathbb D_{n,m}(K)^2\}\bigr)\in\mathbb N^{\ast}\cup\{\infty\}$\index{$\pi^{\E}_{n,m}(K)$ main invariant}.
\end{definition}

Clearly, we have $\pi_{n,m}(K)\le\pi_{n,m}^{\E}(K)\le\pi_{n,m}^{\S}(K)$, with the last inequality by Theorem \ref{T6}(4) and the identity $\TGA_{n,m}(\mathbb F_2)=\STGA_{n,m}(\mathbb F_2)$.

We use Notation \ref{N1} in the context of $\perm(K^n)\cong S_{|K|^n}$ and its subgroup $\Perm(K^n)=\varrho_{n,K}\bigl(\TGA_n(K)\bigr)$ introduced in Definition \ref{D10.2}. 

\phantomsection{The function $$\ell_{n,K}:\Perm(K^n)\rightarrow\mathbb N^{\ast}$$\index{length function!$\ell_{n,K}$ degree length function} that maps $\sigma\in\Perm(K^n)$ to}\label{PH16d}
$$\ell_{n,K}(\sigma):=\min\bigl(\ell(a)\bigl|a\in\TGA_n(K),\; a(K)=\sigma\bigr)$$
is a length function as it satisfies the following axioms.

\medskip\noindent
{$\pmb{A1_{n,K}}$.} For $\sigma\in\Perm(K^n)$, we have $\ell_{n,K}(\sigma)=1$ iff $\sigma\in\AGL_n(K)$.

\smallskip\noindent
{$\pmb{A2_{n,K}}$.} For each $\sigma\in\Perm(K^n)$, we have $\ell_{n,K}(\sigma)=\ell_{n,K}(\sigma^{-1})$.

\smallskip\noindent
{$\pmb{A3_{n,K}}$.} For each pair $(\sigma,\theta)\in\Perm(K^n)^2$, we have $\ell_{n,K}(\sigma\theta)\le\ell_{n,K}(\sigma)\ell_{n,K}(\theta)$ by Equation (\ref{EQ3}).

\medskip
We have a variant of the restriction of $\ell_{n,K}$ to $\Alt(K^n)$ as follows.

The function $$\ell^{\S}_{n,K}:\Alt(K^n)\rightarrow\mathbb N^{\ast}$$\index{length function!$\ell^{\S}_{n,K}$ degree length function} 
that maps $\sigma\in\Alt(K^n)$ to
$$\ell^{\S}_{n,K}(\sigma):=\min\bigl(\ell(a)\bigl|a\in\STGA_n(K),\; a(K)=\sigma\bigr)$$
is also a length function as it satisfies the following axioms.

\medskip\noindent
{$\pmb{A1^{\S}_{n,K}}$.} For $\sigma\in\Alt(K^n)$, we have $\ell^{\S}_{n,K}(\sigma)=1$ iff $\sigma\in\ASL_n(K)$.

\smallskip\noindent
{$\pmb{A2^{\S}_{n,K}}$.} For each $\sigma\in\Alt(K^n)$, we have $\ell^{\S}_{n,K}(\sigma)=\ell^{\S}_{n,K}(\sigma^{-1})$.

\smallskip\noindent
{$\pmb{A3^{\S}_{n,K}}$.} For each pair $(\sigma,\theta)\in\Alt(K^n)^2$, we have $\ell^{\S}_{n,K}(\sigma\theta)\le\ell^{\S}_{n,K}(\sigma)\ell^{\S}_{n,K}(\theta)$ by Equation (\ref{EQ3}).

\medskip
Clearly, for each $\sigma\in\Alt(K^n)$ we have $\ell_{n,K}(\sigma)\le \ell^{\S}_{n,K}(\sigma)$.

\phantomsection{The difficulties in the study of the length function $\ell_{n,K}$ stem from the facts that (i) in general $\AGL_n(K)$ is not a normal subgroup of $\GA_n(K)$ or $\TGA_n(K)$ or $\perm(K^n)$ and (ii) $\Ker(\varrho_{n,K})$ neither is contained in $\AGL_n(K)$ nor contains $\AGL_n(K)$.\footnote{In other words, $\Ker(\varrho_{n,K})$ is not a convex subgroup of $\GA_n(K)$ for the semigauge $\ln\ell_{\GA_n(K)}$ in the sense of \cite{Pro}, Sect.\ 3, Def.\ before Lem.\ 3.1. So \cite{Pro}, Thm.\ 3.4 does not apply in the context of $\varrho_{n,K}$ and hence we do not know if the two epimorphisms $\GA_n(K)\rightarrow\Im(\varrho_{n,K})$ and $\TGA_n(K)\rightarrow\Perm(K^n)$ induced by $\varrho_{n,K}$ are semigauge maps in the sense of \cite{Pro}, Sect.\ 3, Def.\ after Lem.\ 3.1 with respect to the semigauges $\ln\ell_{\GA_n(K)}$ and $\ln\ell_{n,K}$.} These difficulties explain why, while upper bounds for $\ell_{n,K}(\sigma)$ with $\sigma\in\Perm(K^n)\setminus\AGL_n(K)$ are omnipresent in this monograph, the rare cases when $\ell_{n,K}(\sigma)$ is bounded from below rely on the next proposition that uses the ideal}\label{PH17} 
$$\mathfrak I:=\bigl(x_i^{|K|}-x_i|i\in \llbracket1,n\rrbracket\bigr)$$ 
of $R=K[x_1,\ldots,x_n]$.

\begin{proposition}\label{PR13}
Let $\sigma\in\Perm(K^n)$. Let $a=\e(f_1,\ldots,f_n)\in\TGA_n(K)$ be such that $a(K)=\sigma$ and we write $a^{-1}=\e(g_1,\ldots,g_n)$ (so we have a $2n$-tuple $(f_1,\ldots,f_n,g_1,\ldots,g_n)\in R^{2n}$). Then the following properties hold.

\medskip
{\bf (1)} For each $i\in \llbracket1,n\rrbracket$, the image of $f_i$ (or $g_i$) in the quotient ring $\bar R:=R/\mathfrak I$ depends only on $\sigma$ and not on $a$.

\smallskip
{\bf (2)} Suppose that for each $(i,j)\in \llbracket1,n\rrbracket^2$, the partial degree of $f_i$ with respect to $x_j$ is at most $|K|-1$. Then either $\ell_{n,K}(\sigma)=\ell(a)$ or $\max\bigl(|K|,\pi(a)\bigr)\le\ell_{n,K}(\sigma)<\ell(a)$. In particular, if $\pi(a)<|K|\le\pi(a^{-1})$, then $|K|\le\ell_{n,K}(\sigma)\le\pi(a^{-1})=\ell(a)$. 

\smallskip
{\bf (3)} Suppose that for each $(i,j)\in \llbracket1,n\rrbracket^2$, the partial degrees of both $f_i$ and $g_i$ with respect to $x_j$ is at most $|K|-1$. Then $\ell_{n,K}(\sigma)=\ell(a)$. In particular, if $\ell(a)<|K|$, then $\ell_{n,K}(\sigma)=\ell(a)$.

\smallskip
{\bf (4)} Suppose that $n=2$, $\j_a=2$, and $\ell(a)=|K|$. Then $\ell_{2,K}(\sigma)=|K|$.

\smallskip
{\bf (5)} If $\sigma\in\Alt(K^n)$ and $a\in\STGA_n(K)$, then parts (1) to (4) hold with $\ell_{n,K}(\sigma)$ replaced by $\ell^{\S}_{n,K}(\sigma)$.

\smallskip
{\bf (6)} For a selective shift $\shift^{v_0+w_0}_{V\oplus W}\in\Perm(K^n)$ as in Definition \ref{D10} and an automorphism $a^{v_0+w_0}_{V\oplus W}$ as in Lemma \ref{F5}(2), we have identities 
$$\ell_{n,K}(\shift^{v_0+w_0}_{V\oplus W})=\ell^{\S}_{n,K}(\shift^{v_0+w_0}_{V\oplus W})=\ell(a^{v_0+w_0}_{V\oplus W})=\dim_K(W)(|K|-1).$$
\end{proposition}

\begin{proof}
Let $b=\e(h_1,\ldots,h_n)\in\TGA_n(K)$ be such that $b(K)=\sigma$. Let the $n$-tuple $(y_1,\ldots,y_n)\in R^n$ be such that for $c:=ba^{-1}\in\TGA_n(K)$ we have the identity $c=\e(y_1,\ldots,y_n)$. As $c(K)$ is the identity permutation of $K^n$, we have $x_i+\mathfrak I=y_i+\mathfrak I$ for each $i\in \llbracket1,n\rrbracket$. From this and $b=ca$ we get that $h_i+\mathfrak I=f_i+\mathfrak I$ for each $i\in \llbracket1,n\rrbracket$. 

From the last paragraph and its analog for $(a^{-1},\sigma^{-1})$ instead of $(a,\sigma)$ we get that part (1) holds.

For part (2) it suffices to prove that if $\ell_{n,K}(\sigma)<\ell(a)$ then the inequality $\max\bigl(|K|,\pi(a)\bigr)\le\ell_{n,K}(\sigma)$ holds. We can assume that $\ell(b)=\ell_{n,K}(\sigma)$. As we have $h_i+\mathfrak I=f_i+\mathfrak I$ for each $i\in \llbracket1,n\rrbracket$ by part (1), from the hypothesis on partial degrees we get that $\pi(b)\ge\pi(a)$. As $\ell(b)\ge\pi(b)$, we get that $\pi(a)\le\ell_{n,K}(\sigma)$. So to complete the proof of part (2) it suffices to show that the assumption $|K|>\ell_{n,K}(\sigma)$ leads to a contradiction. This assumption implies that $\pi(b)<|K|$. So for each $(i,j)\in \llbracket1,n\rrbracket^2$, the partial degree of $h_i$ with respect to $x_j$ is at most $|K|-1$. From this and the identity $h_i+\mathfrak I=f_i+\mathfrak I$ for each $i\in \llbracket1,n\rrbracket$ it follows that $b=a$, hence $\ell(b)=\ell_{n,K}(\sigma)$ which contradicts the inequality $\ell_{n,K}(\sigma)<\ell(a)$. So part (2) holds. 

Part (3) follows from part (2) applied to $a$ and $a^{-1}$: as it is impossible to have $\max\bigl(|K|,\pi(a),\pi(a^{-1})\bigr)\le\ell_{n,K}(\sigma)<\ell(a)$, we get that $\ell_{n,K}(\sigma)=\ell(a)$.

For part (4), we consider a standard product decomposition $a=a_1a_2c$ with $(a_1,a_2)\in\NAT(K)^2$ with $\dir(a_1)\neq\dir(a_2)$ and $c\in\AGL_2(K)$. By replacing $(\sigma,a)$ by $\bigl(\sigma c(K)^{-1},ac^{-1}\bigr)$ we can assume that $c=1_{\mathbb A^2_K}$. So $a=a_1a_2$. As $\GL_2(K)$ acts $2$-transitively on lines passing through the origin, we can assume that we have $a=\e\bigl(x_1+g(x_2),x_2\bigr)\e\bigl(x_1,x_2+f(x_1)\bigr)$ with $(f,g)\in K[x]^2$ such that $\min\bigl(\deg(f),\deg(g)\bigr)\ge 2$. As $a=\e\bigl(x_1+g(x_2+f(x_1)),x_2+f(x_1)\bigr)$, we have $\ell(a)=\deg(f)\deg(g)$ by Proposition \ref{PR11}(4.a); thus $\deg(f)\deg(g)=|K|$. We show that the assumption that $\ell_{2,K}(\sigma)<|K|$ leads to a contradiction. This assumption implies that there exists a quadruple $(g_1,g_2,f_1,f_2)\in R^4$ such that $$\e\Bigl(x_1+g(x_2+f(x_1))+\sum_{i=1}^2 (x_i^{|K|}-x_i)g_i(x_1,x_2),x_2+f(x_1)+\sum_{i=1}^2 (x_i^{|K|}-x_i)f_i(x_1,x_2)\Bigr)$$
is an automorphism $b\in\GA_2(K)[|K|-1]$. As $x_2+f(x_1)+\sum_{i=1}^2 (x_i^{|K|}-x_i)f_i(x_1,x_2)$ has degree at most $|K|-1$, we get that we can assume that $f_1=f_2=0$. As $b$ is an automorphism, for each $(\alpha,\beta)\in L^2$, with $L$ an algebraic closure of $K$, the system of two equations
$$x_1+g(x_2+f(x_1))+\sum_{i=1}^2 (x_i^{|K|}-x_i)g_i(x_1,x_2)-\alpha=x_2+f(x_1)-\beta=0$$
in the indeterminates $x_1$ and $x_2$ defines a scheme over $\Spec L$ isomorphic to $\Spec L$. Denoting $\gamma:=\alpha-g(\beta)$, this system is equivalent to the single equation
$$x_1+(x_1^{|K|}-x_1)g_1\bigl(x_1,\beta-f(x_1)\bigr)+\bigl(\beta^{|K|}-f(x_1)^{|K|}-\beta+f(x_1)\bigr)g_2\bigl(x_1,\beta-f(x_1)\bigr)=\gamma$$
in the indeterminate $x_1$. It follows that for each $\beta\in L$ the polynomial
\begin{equation}\label{EQ14.5}
(x_1^{|K|}-x_1)g_1\bigl(x_1,\beta-f(x_1)\bigr)+\bigl(\beta^{|K|}-f(x_1)^{|K|}-\beta+f(x_1)\bigr)g_2\bigl(x_1,\beta-f(x_1)\bigr)\in L[x_1]
\end{equation}
has degree $1$. So there exists a unique pair $(h,\delta)\in (L[x])\times L^{\ast}$ such that 
$$(x_1^{|K|}-x_1)g_1(x_1,x_2)+(x_2^{|K|}-x_2)g_2(x_1,x_2)=h\bigl(x_2+f(x_1)\bigr)+x_1\delta;$$
from the uniqueness part and Galois theory we get that $(h,\delta)\in (K[x])\times K^{\ast}$.
From this and the inequality $\ell(b)<|K|$ we get that $\deg(h)=\deg(g)$. Moreover, as the polynomial $h\bigl(x_2+f(x_1)\bigr)+\delta x_1$ represents the zero function on $K^2$ by Equation (\ref{EQ14.5}), we get that $\delta=0$, a contradiction. 

Part (5) follows from parts (1) to (4) and definitions.

The first (resp.\ second) identity of part (6) follows from part (2) applied to $(\sigma,a)=(\shift^{v_0+w_0}_{V\oplus W},a^{v_0+w_0})$ (resp.\ from Lemma \ref{F5}(2)).
\end{proof}

We have the following direct consequence of Lemma \ref{F8} and Proposition \ref{PR13}.

\begin{corollary}\label{C10}
Suppose that $n\in\mathbb N^{\ast}\setminus\{1,2\}$. If the finite field $K$ has at least $3$ elements, then the length functions $\ell_{n,K}$ and $\ell^{\S}_{n,K}$ are not regular.
\end{corollary}

\begin{proof}
To prove this we can assume that $n=3$. In the proof of Lemma \ref{F8} we take $f(x)=g(x)=h(x)=x^2$. So $ac=\e(x_1,x_2+x_1^2,x_3+x_1^4+2x_1^2x_2+x_2^2)$. Clearly, the partial degrees of the polynomials in $K[x_1,x_2,x_3]$ that define $a$, $b$, $c$, $ab$, $b^{-1}c$, their inverses, and $c^{-1}a^{-1}$ are at most $|K|-1$ and $\pi(ac)=4$. Thus $\ell_{3,K}\bigl(\star(K)\bigr)=\ell(\star)$ for each $\star\in\{a,b,c,ab,b^{-1}c\}$ and, if $|K|\ge 5$, for $\star=ac$ by Proposition \ref{PR13}(3). If $|K|=4$, then $\ell_{3,K}\bigl(a(K)c(K)\bigr)=\ell(ac)$ by Proposition \ref{PR13}(2) applied to $c^{-1}a^{-1}$. 

In this paragraph we assume that $|K|=3$. We show that the assumption that $\ell_{3,K}\bigl(a(K)c(K)\bigr)\neq\ell(ac)$ leads to a contradiction. From this assumption, the fact that $\ell(ac)=4$, and Proposition \ref{PR13}(2) applied to $c^{-1}a^{-1}$ we get that $\ell_{3,K}\bigl(a(K)c(K)\bigr)=3$. This implies that there exist constants $\alpha_{ij}\in K$ indexed by pairs $(i,j)\in\llbracket1,6\rrbracket\times\{1,2,3\}$ such that by denoting $F_i:=\sum_{j=1}^3 \alpha_{ij} (x_j^3-x_j)$ for $i\in \llbracket1,6\rrbracket$, 
$b_1:=\e(x_1+F_1,x_2+x_1^2+F_2,x_3+x_1^2+2x_1^2x_2+x_2^2+F_3)\in\GA_3(K)[3]$ has its inverse equal to $b_2:=\e(x_1+F_4,x_2-x_1^2+F_5,x_3-x_2^2+F_6)\in\GA_3(K)[3]$. As the Jacobian determinant of $b_1$ is in $K^{\ast}$, we get that $\alpha_{12}=\alpha_{13}=\alpha_{23}=0$. Writing $b_1b_2=\e(F_7,F_8,F_9)$, we have $F_7=x_1=x_1+F_4+\alpha_{11}(x_1^3+F_4^3-x_1-F_4)$. Thus $F_4=0$ and $\alpha_{11}=0$. As the Jacobian determinant of $b_2$ is in $K^{\ast}$ and $F_4=0$, we get that $\alpha_{53}=0$. As $b_1$ and $b_2$ are automorphisms we get that $\alpha_{22}=\alpha_{52}=0$ and by computing $F_8$ we get that $\alpha_{51}=-\alpha_{21}$; so $F_5=-F_2$. Writing $b_2b_1=\e(F_{10},F_{11},F_{12})$, it follows that $F_{12}=x_3$ is the sum of the polynomials $x_3-x_2^2+[\sum_{j=1}^3 \alpha_{6j}(x_j^3-x_j)]+x_1^2+2x_1^2(x_2-x_1^2-F_2)+(x_2-x_1^2-F_2)^2$ and
$$\alpha_{31}(x_1^3-x_1)+\alpha_{32}(x_2^3-x_1^6-F_2^3-x_2+x_1^2+F_2)+\alpha_{33}(x_3^3-x_2^6+F_6^3-x_3+x_2^2-F_6).$$
From this and the fact that $F_2=\alpha_{21}(x_1^3-x_1)$, by considering the coefficients of $x_1^4$ and $x_1x_2$ we get that $\alpha_{21}^2=1$ and $\alpha_{21}=0$, a contradiction.

Thus we have $\ell_{3,K}\bigl(a(K)c(K)\bigr)=\ell(ac)=4$ in all cases.

From this and the proof of Lemma \ref{F8} we get that the identities
$$2=\ell_{3,K}\bigl(a(K)b(K)\bigr)<\ell_{3,K}\bigl(a(K)\bigr)\ell_{3,K}\bigl(b(K)\bigr)=4,$$ 
$$2=\ell_{3,K}\bigl(b(K)^{-1}c(K)\bigr)<\ell_{3,K}\bigl(b(K)\bigr)\ell_{3,K}\bigl(c(K)\bigr)=4,$$ 
and $\ell_{3,K}\bigl(a(K)c(K)\bigr)=\ell_{3,K}\bigl(a(K)\bigr)\ell_{3,K}\bigl(c(K)\bigr)=4$ hold. 

As the permutations $a(K)$, $b(K)$, and $c(K)$ are even, above we can replace $\ell_{3,K}$ by $\ell_{3,K}^{\S}$ by Proposition \ref{PR13}(5). So the corollary holds.
\end{proof}

\begin{example}\normalfont\label{EX9}
For a pair $(n,r)\in (\mathbb N^{\ast}\setminus\{1\})^2$ we consider the automorphism $a:=\e(x_1,x_2-x_1^r+1,x_3-x_2^r+1,\ldots,x_n-x_{n-1}^r+1)\in\STGA_n(K)$ and the permutation $\sigma:=a(K)\in\Alt(K^n)$. We have $a^{-1}=\e(f_1,\ldots,f_n)\in\STGA_n(K)$, where $(f_1,\ldots,f_n)\in R^n$ is defined recursively by the following rules $f_1:=x_1$ and $f_{i+1}:=x_{i+1}+f_i^r-1$ for $i\in \llbracket2,n\rrbracket$. By induction on $i\in \llbracket1,n\rrbracket$ we get that $\deg(f_i)=r^{i-1}$. Thus Equation (\ref{EQ1}) gives
$$\ell(a)=\ell(a^{-1})=\pi(a^{-1})=r^{n-1}=\pi(a)^{n-1}.$$
Proposition \ref{PR13}(2) and (5) gives that $\ell_{n,K}(\sigma)=\ell^{\S}_{n,K}(\sigma)=r^{n-1}$ if $r^{n-1}\le |K|$ and $|K|\le\ell_{n,K}(\sigma)=\ell^{\S}_{n,K}(\sigma)\le r^{n-1}$ if $r<|K|<r^{n-1}$. If we have $r=|K|-1\ge 2$, $n\ge 3$, and $p:=\char(K)$, then $\sigma$ is a product of disjoint permutations that are $p^j$-cycles with $j\in \llbracket1,n-1\rrbracket$ arbitrary; for instance, for $n=3$, 
$$\supp(\sigma)=\{(\alpha,\beta,\gamma)\in K^3|\alpha\beta=0\},$$ 
$\n(\sigma)=|K|(2|K|-1)$, and in the cyclic decomposition of $\sigma$ the number of $p^2$-cycles is $\frac{|K|}{p}$ and hence the umber of $p$-cycles is $\frac{2|K|^2}{p}-\frac{|K|(p+1)}{p}$. 
\end{example}

\begin{example}\normalfont\label{EX10}
Let $(r,s)\in (\mathbb N^{\ast}\setminus\{1,2\})^2$ and $p:=\char(K)$. Let 
$$a:=\e(x_1,x_2-x_1^r,x_3-x_2^s)\in\STGA_3(K)\;\;\textup{and}\;\;\sigma:=a(K)\in\Alt(K^3).$$ 
We have $\pi(a)=\max(r,s)$,
$$a^{-1}=\e\bigl(x_1,x_2+x_1^r,x_3+(x_2+x_1^r)^s\bigr)\in\STGA_3(K)[rs],$$
and $\pi(a^{-1})=rs$. We assume that $p|s-1$ and $r(s-1)\ge |K|$. Let
$$b:=\e\Bigl(x_1,x_2+x_1^r,x_3+x_1^{rs-|K|+1}+sx_2x_1^{r(s-1)-|K|+1}+\sum_{i=2}^s\binom{s}{i}x_2^ix_1^{r(s-i)}\Bigr)$$
in $\STGA_3(K)$. Clearly, $b(K)=a^{-1}(K)=\sigma^{-1}$. Thus $b^{-1}(K)=\sigma$. Defining
$$f(x_1,x_2):=-x_1^{rs-|K|+1}-s(x_2-x_1^r)x_1^{r(s-1)-|K|+1}-\sum_{i=2}^s\binom{s}{i}(x_2-x_1^r)^ix_1^{r(s-i)}$$
in $K[x_1,x_2]$, we have $b^{-1}=\e\bigl(x_1,x_2-x_1^r,x_3+f(x_1,x_2)\bigr)$. The coefficient of $x_1^{rs}$ in $f$ is $1-s$ and hence it is $0$. This implies that 
$$\ell_{3,K}(\sigma)\le\ell(b^{-1})=\ell(b)\le\max\bigl(rs-|K|+1,r(s-1)+1\bigr)<\ell(a)=\pi(a^{-1})=rs.$$ Thus, if $\max(s,r)<|K|$ (hence $p^2\mid |K|$), then 
$$\pi(a)< |K|\le\ell_{3,K}(\sigma)\le \ell^{\S}_{3,K}(\sigma)<rs=\pi(a^{-1})$$ by Proposition \ref{PR13}(2) and (5).
\end{example}

\begin{remark}\normalfont\label{R4} 
{\bf (1)} Let $n\ge 2$ and $l\ge |K|$ be integers, $a\in\STGA_n(K)$, and $\break\sigma:=a(K)\in\Alt(K^n)$. If $\pi(a)=\pi(a^{-1})=l$, then the inequality $\ell^{\S}_{K,n}(\sigma)\le\ell(a)$ is strict in general. Like for $a:=\e(x_1,x_2+x_1^l-x_1^{l-|K|+1},x_3,\ldots,x_n)\in\STGA_n[l]$, we have $\ell^{\S}_{K,n}(\sigma)=1<\ell(a)=l=\pi(a)=\pi(a^{-1})$. Similarly, writing $l=s|K|+r$ with $(s,r)\in\mathbb N^{\ast}\times \llbracket0,|K|-1\rrbracket$, for $n\ge 3$ and $a:=\e(x_1,x_2,x_3+x_1^rx_2^{s|K|},x_4,\ldots,x_n)$, we have $\ell^{\S}_{K,n}(\sigma)\le r+s<\ell(a)=l=\pi(a)=\pi(a^{-1})$; moreover, if $r+s<|K|$, then we have $\ell_{K,n}(\sigma)=\ell^{\S}_{K,n}(\sigma)=r+s$ by Proposition \ref{PR13}(3) and (5) applied to the automorphism $b:=\e(x_1,x_2,x_3+x_1^rx_2^s,x_4,\ldots,x_n)\in\STGA_n(K)$ that satisfies $b(K)=\sigma$. Hence Proposition \ref{PR13}(3) is optimal.

\smallskip
{\bf (2)} Example \ref{EX10} also gives that no one of the relations of the conclusions of Proposition \ref{PR13}(2) can be strengthened.
\end{remark}

\section{On Furter's lengths when $n=2$}\label{S12}

The following definition is suggested by Definition \ref{D13}.

\begin{definition}\label{D15}
Let $\sigma\in\Perm(K^2)$. 

\medskip
{\bf (1)} By the Furter's length of $\sigma$\index{Furter's length!Furter's length of a permutation} we mean the smallest $\j_{\sigma}\in\mathbb N$ for which there exists $a\in\GA_2(K)$ such that $a(K)=\sigma$ and $\j_{\sigma}$ is the Furter's length $\j_a$ of $a$.

\smallskip
{\bf (2)} An automorphism $a\in\GA_2(K)$ is called a minimal representation\index{minimal representation} of $\sigma$ if the following properties hold for it.

\medskip\noindent
{\bf (2.a)} We have $\sigma=a(K)$.

\smallskip\noindent
{\bf (2.b)} It has a reduced standard product decomposition $a=a_1\cdots a_{\j_a} c$ such that for each $i\in \llbracket1,\j_a\rrbracket$ we have $\ell(a_i)\le |K|-1$.

\medskip
{\bf (3)} An automorphism $b\in\GA_2(K)$ is called a strict minimal representation\index{minimal representation!strict} of $\sigma$ if it is a minimal representation of $\sigma$ with $\j_a=\j_\sigma$.

\smallskip
{\bf (4)} By the Furter's length of $K$\index{Furter's length!Furter's length of a field} we mean
$$\j_{K}:=\max\bigl(\j_{\sigma}|\sigma\in\Perm(K^2)\bigr)\in\mathbb N.$$
\end{definition}

Basic properties of minimal representations are grouped together as follows.

\begin{proposition}\label{PR14}
Let $\sigma\in\Perm(K^2)$ and let $j:=\j_{\sigma}$. Let $a\in\GA_2(K)$ be such that $a(K)=\sigma$. Then the following properties hold.

We assume that one of the following two conditions holds.

\medskip
{\bf (1)} The automorphism $a$ is a minimal representation of $\sigma$ provided one of the following two conditions holds.

\medskip\noindent
{\bf (1.a)} 
We have $\j_a=j$ and for each $b\in\GA_2(K)$ such that $b(K)=\sigma$ and $\j_b=j$, the inequality $\ell(a)\le\ell(b)$ holds.

\smallskip\noindent
{\bf (1.b)} We have $\ell(a)=\ell_{2,K}(\sigma)$ and for each $b\in\GA_2(K)$ such that $b(K)=\sigma$ and $\ell(b)=\ell_{2,K}(\sigma)$, the inequality $\j_b\ge\j_a$ holds.

\medskip
{\bf (2)} The automorphism $a$ is a strict minimal representation of $\sigma$ and $\ell_{2,K}(\sigma)$ is a product of $j$ primes provided one of the following two conditions holds.

\medskip\noindent
{\bf (2.a)} Condition (1.a) holds and $\ell(a)$ is a prime.

\smallskip\noindent
{\bf (2.b)} Condition (1.b) holds and $\ell_{2,K}(\sigma)$ is a product of at most $j$ primes.
\end{proposition}

\begin{proof}
If $j=0$, then $\ell_{2,K}(\sigma)=1$ and it follows that $a\in\AGL_2(K)$, therefore $\ell(a)=1=\ell_{2,K}(\sigma)$, $\j_a=0$, and clearly $a$ is a strict minimal representation of $\sigma$. Thus to prove both parts we can assume that $j\in\mathbb N^{\ast}$. Condition (1.a) or (1.b) holds for both parts.

We consider a reduced standard product decomposition $a=a_1\cdots a_{\j_a}c$. By Proposition \ref{PR11}(3) we have $\ell(a)=\prod_{l\in \llbracket1,\j_a\rrbracket} \ell(a_l)$.

We show that the assumption that there exists $i\in \llbracket1,\j_a\rrbracket$ such that $\ell(a_i)\ge |K|$ leads to a contradiction. Let $(f_i,c_i)\in x^2K[x]\times \GL_2(K)$ be a pair that defines $a_i$. We have $\deg(f_i)=\ell(a_i)\ge |K|$. Let $g_i\in K[x]$ be a polynomial of minimal degree $d_i$ such that $g_i(\alpha)=f_i(\alpha)$ for each $\alpha\in K$. We have $d_i\le |K|-1$ by Lagrange interpolation. Let $a_i':=c_i\e\bigl(x_1,x_2+g_i(x_1)\bigr)c_i^{-1}\in\GA_2(K)$ and
$$a':=a_1\cdots a_{i-1}a_i'a_{i+1}\cdots a_{\j_a}c\in\GA_2(K).$$ As $a_i'(K)=a_i(K)$, it follows that $a '(K)=a(K)=\sigma$. 

We consider two disjoint cases as follows.

{\bf Case 1: $d_i\in \llbracket2,|K|-1\rrbracket$.} We have $\ell(a_i')=d_i<|K|\le \deg(f_i)=\ell(a_i)$. As $\dir(a_i')=\dir(a_i)$, Proposition \ref{PR11}(3) also gives $\ell(a')=\ell(a_i')\prod_{l\in \llbracket1,\j_a\rrbracket\setminus\{i\}} \ell(a_l)$. Therefore $\ell(a')<\prod_{l=1}^{\j_a} \ell(a_l)=\ell(a)$. Hence $\ell(a)>\ell_{2,K}(a)$ and, if condition (1.b) holds, we reached a contradiction. The Furter's length of $a'$ belongs to $\llbracket j,\j_a\rrbracket$. If condition (1.a) holds, then $\j_a=j$ and hence the Furter's length of $a'$ is $j$; thus we have $\ell(a)\le\ell(a')$, a contradiction.

{\bf Case 2: $d_i\in\{0,1,-\infty\}$.} Then $a_i'\in\AGL_2(K)$. Based on this and Definition \ref{D13}(2) and (3) we get that $\j_{a'}<\j_a$, hence $j<\j_a$ and thus, if condition (1.a) holds, we reached a contradiction. From Inequality (\ref{EQ3}) we get that $\ell(a')\le\prod_{l\in \llbracket1,\j_a\rrbracket\setminus\{i\}} \ell(a_l)<\ell(a)$, hence, if condition (1.b) holds, we reached a contradiction.

We conclude that $\ell(a_i)\le |K|-1$ for each $i\in \llbracket1,\j_a\rrbracket$ and hence $a$ is a minimal representation of $\sigma$. So part (1) holds.

Assume first that condition (2.b) holds. As $\ell_{2,K}(\sigma)=\ell(a)=\prod_{l\in \llbracket1,\j_a\rrbracket} \ell(a_i)$ and each $\ell(a_i)\ge 2$ is a product of one or more primes for each $i\in \llbracket1,\j_a\rrbracket$, it follows that $\ell_{2,K}(\sigma)$ is divisible by a product of $\j_ a$ primes. As $j\le \j_a$, we get that $j=\j_a$ and that each $\ell(a_i)$ is a prime. So $a$ is a strict minimal representation of $\sigma$ and $\ell_{2,K}(\sigma)$ is a product of $j$ primes.

Assume now that condition (2.a) holds. Then $\j_a=1$ and hence $j=1$. From $\ell(a_1)<|K|$ and Proposition \ref{PR13}(3) we get that $\ell(a)=\ell(a_1)=\ell_{2,K}\bigl(a_1(K)\bigr)$. From this and the identity $a_1(K)=\sigma c(K)^{-1}$ which implies that $\ell_{2,K}(\sigma)=\ell_{2,K}\bigl(a_1(K)\bigr)$, we get that condition (1.b) holds and hence $a$ is a strict minimal representation of $\sigma$ and $\ell(a)$ is a prime and thus a product of $j$ primes.
\end{proof}

For $n\in\mathbb N^{\ast}$ let $\mathbb M_n$ be the multiplicative submonoid of $\mathbb N^{\ast}$ generated by prime numbers less of equal to $n$.

\begin{corollary}\label{C11}
The image of the length function $\ell_{2,K}:\Perm(K^2)\rightarrow\mathbb N^{\ast}$ is contained in $\mathbb M_{|K|-1}$ and contains all primes less than or equal to $|K|-1$. 
\end{corollary}

\begin{proof}
Let $\sigma\in\Perm(K^2)$. Let $a\in\GA_2(K)$ be a minimal representation of $\sigma$ with $\ell(a)=\ell_{2,K}(\sigma)$ by Proposition \ref{PR14}(1.b). We consider a standard product decomposition $a=a_1\cdots a_{\j_a} c$ such that for each $i\in \llbracket1,\j_a\rrbracket$ we have $\ell(a_i)\le |K|-1$ by Definition \ref{D15}(2.b). From this and Proposition \ref{PR11}(4.a) we get the identity $\ell(a)=\prod_{i=1}^{\j_a} \ell(a_i)$. From the last two sentences we get that $\ell_{2,K}(\sigma)\in\mathbb M_{|K|-1}$.

The fact that the image contains all the primes less than or equal to $|K|-1$ follows from Proposition \ref{PR13}(3). Thus the corollary holds.
\end{proof}

If $n=2$ we have two additional length functions as follows.

The function $$\ell^{-}_{2,K}:\Perm(K^2)\rightarrow\mathbb N^{\ast}$$\index{length function!$\ell^{-}_{2,K}$ degree length function} 
that maps $\sigma\in\Perm(K^2)$ to
$$\ell^{-}_{2,K}(\sigma):=2^{\j_{\sigma}}$$
is also a length function as we have the following axioms. 

\medskip\noindent
{$\pmb{A1^-_{2,K}}$.} For $\sigma\in\Perm(K^2)$, we have $\ell^-_{2,K}(\sigma)=1$ iff $\sigma\in\AGL_2(K)$.

\smallskip\noindent
{$\pmb{A2^-_{2,K}}$.} For each $\sigma\in\Perm(K^2)$, we have $\ell^-_{2,K}(\sigma)=\ell^-_{2,K}(\sigma^{-1})$.

\smallskip\noindent
{$\pmb{A3^-_{2,K}}$.} For each pair $(\sigma,\theta)\in\Perm(K^2)^2$, we have $\ell^-_{2,K}(\sigma\theta)\le\ell^-_{2,K}(\sigma)\ell^-_{2,K}(\theta)$ as $\j_{\sigma\theta}\le \j_{\sigma}+\j_{\theta}$.

\medskip
The function $$\ell^{+}_{2,K}:\Perm(K^2)\rightarrow\mathbb N^{\ast}$$\index{length function!$\ell^{+}_{2,K}$ degree length function} 
that maps $\sigma\in\Perm(K^2)$ to
$$\ell^{+}_{2,K}(\sigma):=(|K|-1)^{\j_{\sigma}}$$
is also a length function as the axioms $A1^+_{2,K}$ to $A3^+_{2,K}$ analogous to $A1^-_{2,K}$ to $A3^-_{2,K}$ hold.

\begin{definition}\label{D16}
Let $K$ a finite field and $n\in\mathbb N^{\ast}\setminus\{1\}$. 

\medskip
{\bf (1)} We call $\ell_{n,K}$ as the degree length function on $\Perm(K^n)$ and we call $\ell^{\S}_{n,K}$ as the special degree length function\index{length function!special degree length} on $\Alt(K^n)$.

\smallskip
{\bf (2)} Suppose that $n=2$. We call $\ell^{-}_{2,K}$ (resp.\ $\ell^{\S,-}_{2,K}$) as the minimal\index{length function!minimal degree length} (resp.\ minimal special\index{length function!special minimal degree length}) degree length function on $\Perm(K^2)$.
\end{definition}

Terminology of Definition \ref{D16}(2) is justified by the following corollary.

\begin{corollary}\label{C12}
For $\sigma\in\Perm(K^2)$ (resp.\ $\sigma\in\Alt(K^2)$) we have inequalities
$$\ell^-_{2,K}(\sigma)\le\ell_{2,K}(\sigma)\le\ell^+_{2,K}(\sigma)\le (|K|-1)^{\j_K}.$$
\end{corollary}

\begin{proof}
Let $(a,b)\in\GA_2(K)^2$ be such that $a$ and $b$ are minimal representations of $\sigma$ with $\j_a=\j_{\sigma}$ by Proposition \ref{PR14}(1.a) and $\ell(b)=\ell_{2,K}(\sigma)$ by Proposition \ref{PR14}(1.b). We have the following inequalities $\ell_{2,K}(\sigma)\le\ell(a)\le (|K|-1)^{\j_{\sigma}}\le (|K|-1)^{\j_K}$ and $\ell^-_{2,K}(\sigma)=2^{\j_{\sigma}}\le 2^{\j_b}\le\ell(b)=\ell_{2,K}(\sigma)$, where the last inequality follows from Proposition \ref{PR11}(4.a). So the corollary holds.\end{proof}

\begin{corollary}\label{C13}
Suppose that $|K|=3$. Then $\ell_{2,K}=\ell^-_{2,K}=\ell^+_{2,K}$ and each permutation $\sigma\in\perm(K^2)$ has a strict minimal representation.
\end{corollary}

\begin{proof}
The first part follows from Corollary \ref{C12}. The second part follows from the first part and Proposition \ref{PR14}(2.b).
\end{proof}

In order to get examples of and applications to permutations $\sigma\in\Perm(K^2)$ with $\j_{\sigma}\ge 2$, we first prove the following lemma which supplements Lemma \ref{F6}.

\begin{lemma}\label{L7}
Let $\sigma\in\Perm(K^2)$ with $\j_{\sigma}\in\{1,2\}$ (thus $|K|\ge 3$). Then the following properties hold.

\medskip
{\bf (1)} Suppose that $j_{\sigma}=1$. Then the following properties hold.

\medskip\noindent
{\bf (1.a)}
We have $\n(\sigma)\in\{|K|i\bigl|i\in\llbracket1,|K|-1\rrbracket\}\sqcup\llbracket|K|^2-|K|+1,|K|^2\rrbracket$.

\smallskip\noindent
{\bf (1.b)} If $|K|\mid\n(\sigma)$ and $\n(\sigma)\notin\{|K|^2-|K|,|K|^2\}$, then $\supp(\sigma)$ is a union of disjoint lines (i.e., of lines that have the same uniquely determined direction).

\smallskip\noindent
{\bf (1.c)} If $\n(\sigma)=|K|^2-|K|$, then $K^2\setminus\supp(\sigma)$ is either a line or up to affine automorphisms is the graph of a polynomial function from $K$ to $K$ of degree in $\llbracket2,|K|-1\rrbracket$. 

\smallskip\noindent
{\bf (1.d)} Let $p:=\char(K)$. If $\n(\sigma)=|K|^2-|K|$ and $K^2\setminus\supp(\sigma)$ is the graph of a polynomial function from $K$ to $K$ of degree in $\llbracket2,|K|-1\rrbracket$, then $\sigma$ is a product of $\frac{|K|^2-|K|}{p}$ disjoint $p$-cycles whose supports are collinear with the same uniquely determined direction.

\medskip
{\bf (2)} Suppose that $j_{\sigma}=2$. Then $\n(\sigma)\in\llbracket2|K|-1,|K|^2\rrbracket$.
\end{lemma}

\begin{proof} 
For part (1) we write $\sigma=c(K)a(K)c(K)^{-1}b(K)$, where $a=\e\bigl(x_1,x_2+f(x_1)\bigr)$ for some $f(x)\in K[x]$ with $\deg(f)\in\llbracket2,|K|-1\rrbracket$ and $(b,c)\in\AGL_2(K)\times\GL_2(K)$. By replacing $(\sigma,b)$ by $\bigl(c(K)^{-1}\sigma c(K),c^{-1}bc\bigr)$ we can assume that $c=1_{\mathbb A^2_K}$. We write $b=\e(\gamma_{11}x_1+\gamma_{12}x_2+\delta_1,\gamma_{21}x_1+\gamma_{22}x_2+\delta_2)$ with $(\gamma_{11},\gamma_{12},\delta_1,\gamma_{21},\gamma_{22},\delta_2)\in K^6$ such that $\gamma_{11}\gamma_{22}-\gamma_{12}\gamma_{21}\neq 0$. 

\phantomsection{For $P\in K^2$ we have $\sigma(P)=P$ iff $b(P)=a^{-1}(P)$. Thus $\n(\sigma)$ is $|K|^2$ minus the cardinality of the solution set $\mathcal S$ in $K^2$ of the system of two equations}\label{EXT1}
$$x_1-\gamma_{11}x_1-\gamma_{12}x_2-\delta_1=0=x_2-f(x_1)-\gamma_{21}x_1-\gamma_{22}x_2-\delta_2$$
in the indeterminates $x_1$ and $x_2$. So for part (1.a) it suffices to show that we have $|\mathcal S|\in \{|K|i\bigl|i\in\llbracket1,|K|-1\rrbracket\}\sqcup\llbracket0,|K|-1\rrbracket$. We consider three disjoint cases as follows. 

{\bf Case 1: $\gamma_{12}\neq 0$.} We have 
$$\mathcal S=\Bigl\{\Bigl(\alpha,\frac{\alpha-\gamma_{11}\alpha-\delta_1}{\gamma_{12}}\Bigr)\Bigl|\alpha\in K,\,f(\alpha)+\gamma_{21}\alpha+(\gamma_{22}-1)\frac{\alpha-\gamma_{11}\alpha-\delta_1}{\gamma_{12}}+\delta_2=0\Bigr\}.$$

{\bf Case 2: $(\gamma_{11},\gamma_{12})=(1,0)$.} We have $\gamma_{22}\neq 0$. If $\delta_1\neq 0$, then $\mathcal S=\emptyset$. If $\delta_1=0$ and $\gamma_{22}\neq 1$, then $\mathcal S=\bigl\{\bigl(\alpha,\frac{f(\alpha)+\gamma_{21}\alpha+\delta_2}{1-\gamma_{22}}\bigr)\bigl|\alpha\in K\bigr\}$. If $(\delta_1,\gamma_{22})=(0,1)$, then 
$\mathcal S=\{\alpha\in K|f(\alpha)+\gamma_{21}\alpha+\delta_2=0\}\times K$.

{\bf Case 3: $\gamma_{11}\neq 1$ and $\gamma_{12}=0$.} We have $x_1=\frac{\delta_1}{1-\gamma_{11}}$ and $x_2$ is a solution of a linear equation in $x_2$ and hence $|\mathcal S|\in\{0,1,|K|\}$.

As in all cases we have $|\mathcal S|\in\{|K|i\bigl|i\in\llbracket1,|K|-1\rrbracket\}\sqcup\llbracket0,|K|-1\rrbracket$, part (1.a) holds.

For part (1.b), we have $|K|\mid\n(\sigma)$ and $\n(\sigma)\notin\{|K|^2-|K|,|K|^2\}$ iff $\frac{|\mathcal S|}{|K|}\in\mathbb N^{\ast}\setminus\{1\}$ and this can happen only in Case 2 with $(\delta_1,\gamma_{22})=(0,1)$ in which case $\mathcal S$ is a union of lines parallel to the $x_2$-axis. So part (1.b) holds.

For part (1.c), we have $\n(\sigma)=|K|^2-|K|$ iff $|\mathcal S|=|K|$ and this can happen only in Case 2 with $(\delta_1,\gamma_{22})=(0,1)$ and $\mathcal S$ a line, in Case 3 with $\mathcal S$ a line, and in Case 2 with $\delta_1=0$, $\gamma_{22}\neq 1$, and $\mathcal S$ as the graph of the polynomial function $K\rightarrow K$ defined by the rule $\alpha\mapsto\frac{f(\alpha)+\gamma_{21}\alpha+\delta_2}{1-\gamma_{22}}$. So part (1.c) holds.

For part (1.d), we are in Case 2 with $(\gamma_{11},\gamma_{12},\delta_1)=(1,0,0)$ and therefore $ab=\e\bigl(x_1,x_2+g(x_1)\bigr)$ with $g\in K[x]$ of the same degree as $f$ and having only one root in $K$; therefore $\sigma$ is a product of $\frac{|K|^2-|K|}{p}$ disjoint $p$-cycles whose supports are collinear in the direction $(0:1)$. So part (1.d) holds.

For part (2) we write $\sigma=c_1(K)a_1(K)c_1(K)^{-1}c_2(K)a_2(K)c_2(K)^{-1}b(K)$, where $a_1=\e\bigl(x_1,x_2+f(x_1)\bigr)$ and  $a_2=\e\bigl(x_1+g(x_2),x_2\bigr)$ for some pair $(f,g)\in K[x]^2$ with $(\deg(f),\deg(g))\in \llbracket2,|K|-1\rrbracket^2$ and $(b,c_1,c_2)\in\AGL_2(K)\times\GL_2(K)^2$. As $\GL_2(K)$ acts $2$-transitively on the set of directions in $K^2$, by replacing $(\sigma,b)$ with $\bigl(c(K)\sigma c(K)^{-1},cbc^{-1}\bigr)$ where $c:=c_3c_1^{-1}$ for some $c_3\in\GL_2(K)$ such that $c$ fixes $(1,0)$ and $c_3c_1^{-1}c_2$ fixes $e_2$, we can assume that $c_1=c_2=1_{\mathbb A^2_K}$. 

Let $a:=a_1a_2\in\GA_2(K)$; we have $\ell(a)=\deg(f)\deg(g)$. As 
$$a^{-1}=a_2^{-1}a_1^{-1}=\e\bigl(a_1-g(x_2-f(x_1)),x_2-f(x_1)\bigr),$$ 
for $P\in K^2$ we have $\sigma(P)=P$ iff $b(P)=a^{-1}(P)$. We write $b$ as above. Thus $\n(\sigma)$ is $|K|^2$ minus the cardinality of the solution set $\mathcal T$ in $K^2$ of the system of two equations 
$$x_1-g\bigl(x_2-f(x_1)\bigr)-\gamma_{11}x_1-\gamma_{12}x_2-\delta_1=0=x_2-f(x_1)-\gamma_{21}x_1-\gamma_{22}x_2-\delta_2$$
in the indeterminates $x_1$ and $x_2$. So for part (2) it suffices to show that we have $|\mathcal T|\in\llbracket0,\deg(f)\deg(g)\rrbracket\subset\llbracket0,(|K|-1)^2\rrbracket$. 

Suppose that $\gamma_{22}\neq 1$. Then from the second equation we get directly that $x_2=h(x_1)$, where $h(x):=\frac{f(x)+\gamma_{21}x+\delta_2}{1-\gamma_{22}}\in K[x]$. So $\deg(h)=\deg(g)$ and 
$$\mathcal T=\bigl\{\bigl(\alpha,h(\alpha)\bigr)|\alpha\in K,\,\alpha(1-\gamma_{11})-g\bigl(\gamma_{22}h(\alpha)+\gamma_{21}\alpha+\delta_2\bigr)-\gamma_{12}h(\alpha)-\delta_1\bigr\}.$$ 
Hence $|\mathcal T|\in\llbracket0,\deg(f)\deg(g)\rrbracket$.

If $\gamma_{22}=1$, then $\mathcal T$ is formed by pairs $(\alpha,\beta)\in K^2$ with $\alpha$ a solution of the equation $f(x_1)+\gamma_{21}x_1+\delta_2=0$ and with $\beta$ a solution of the resulting equation  $\alpha-g(\gamma_{21}\alpha+x_2+\delta_2)-\gamma_{11}\alpha-\gamma_{12}x_2-\delta_1
=0$; hence $|\mathcal T|\in\llbracket0,\deg(f)\deg(g)\rrbracket$.

From the last two paragraphs we get that part (2) holds.\end{proof}

Directly from Lemma \ref{L7}(1.a) and (2) and Corollary \ref{C12} we get the following consequence.

\begin{corollary}\label{C13.5}
Suppose that $|K|\ge 3$. If $\sigma\in\Perm(K^2)$ is such that we have $\n(\sigma)\in\llbracket2,|K|-1\rrbracket\cup\llbracket|K|+1,2|K|-2\rrbracket$, then $\j_{\sigma}\ge 3$ and $\ell_{2,K}(\sigma)\ge \ell^{-}_{2,K}(\sigma)\ge 8$.
\end{corollary}

\begin{example}\normalfont\label{EX11}
Suppose that $|K|=3$. 

\medskip
{\bf (1)} If $a:=\e(2x_1,x_2+x_1^2+1)\in\GA_2(K)[2]$, then $\sigma:=a(K)$ is a disjoint product of a $6$-cycle and a $3$-cycle, the support of the $3$-cycle being collinear; we have $\n(\sigma)=9$, $o(\sigma)=6$, $\j_{\sigma}=1$, and $\ell_{2,K}(\sigma)=2$. 

\smallskip
{\bf (2)} If $b:=\e(2x_1,2x_2+x_1^2)\in\SGA_2(K)[2]$, then $\varsigma:=b(K)$ is a product of four disjoint transpositions whose supports generate four distinct directions; we have $\n(\varsigma)=8$, $o(\varsigma)=2$, $\j_{\varsigma}=1$, $\ell_{2,K}(\varsigma)=\ell^{\S}_{2,K}(\varsigma)=2$, and the only fixed point of $\sigma$ is $(0,0)$. 

\smallskip
{\bf (3)} If $c:=\e\bigl(x_1+x_2,x_2+2x_1+1+(x_1+x_2)^2)\bigr)\in\GA_2(K)[2]$, then $\varsigma_1:=c(K)$ is an $8$-cycle; we have $\n(\varsigma_1)=8=o(\varsigma_1)$, $\j_{\varsigma_1}=1$, $\ell_{2,K}(\varsigma_1)=2$, and the only fixed point of $\sigma$ is $(0,2)$. 

\smallskip
{\bf (4)} If $a_1:=\e\bigl(2x_1+x_2,x_2+(2x_1+x_2)^2-1\bigr)\in\GA_2(K)[2]$, then $\sigma_1:=a_1(K)$ is a disjoint product of a $5$-cycle and a transposition, the support of the $5$-cycle being the union of two lines; we have $\n(\sigma_1)=7$, $o(\sigma_1)=10$, $\j_{\sigma_1}=1$, $\ell_{2,K}(\sigma_1)=2$, and the fixed points of $\sigma$ are $(1,-1)$ and $(-1,1)$. 

\smallskip
{\bf (5)} If $b_1:=\e(x_1,2x_2-x_1^2)\in\GA_2(K)[2]$, then $\varsigma_2:=b_1(K)$ is a product of three disjoint transposition whose supports have the same direction; we have $\n(\varsigma_2)=6$, $o(\varsigma_2)=2$, $\j_{\varsigma_2}=1$, $\ell_{2,K}(\varsigma_2)=2$, and $\sigma$ fixes $(0,0)$, $(1,1)$, and $(2,1)$. 

\smallskip
{\bf (6)} If $\n(\sigma)=4$, then we have $\j_{\sigma}\ge 3$ and $\ell_{2,K}(\sigma)\ge 8$ by Corollary \ref{C13.5}. 
\end{example}

\begin{remark}\normalfont\label{R5}
If $|K|\ge 3$, equivalently $\AGL_2(K)\neq\Perm(K^2)$, then we have $\j_{K}\in\mathbb N^{\ast}\setminus\{1\}$ and $\j_{K}$ is the smallest $j\in\mathbb N^{\ast}$ such that we have an identity
$$\Perm(K^2)=\{\sigma\in\Perm(K)| \j_{\sigma}\in\{0,1\}\}^j.$$
There exists $\sigma\in\Perm(K)$ with $\n(\sigma)=|K|$ and $\j_{\sigma}=1$ by Proposition \ref{PR13}(6).
\end{remark}

\begin{example}\normalfont\label{EX12}
Let $p:=\char(K)$. Let $q\in\mathbb N^{\ast}$ be such that $|K|=p^q$. We assume that $q\ge 2$. Let $(d_1,d_2)\in\{1,\ldots,q-1\}^2$. For $\iota\in\{1,2\}$, let $V_{\iota}$ be an $\mathbb F_p$-vector subspace of $K$ of dimension $d_{\iota}$. Let $(v_1,v_2)\in (V_1\setminus\{0\})\times (V_2\setminus\{0\})\subset K^2$. Let the polynomial $f_{V_{\iota}}(x)\in K[x]$ of degree at most $|K|-1$ be such that $f_{V_{\iota}}(V_{\iota})=\{1\}$ and $f_{V_{\iota}}(K\setminus V_{\iota})=\{0\}$. We have $f_{V_{\iota}}(x)\in K^{\ast}\prod_{\gamma\in K\setminus V_{\iota}} (x-\gamma)$. Thus $\deg(f_{V_{\iota}})=|K|-|V_{\iota}|$. For instance, if $d_1=1$ and $V_1=\mathbb F_p\subset K$, then $f_{V_1}=\frac{x^{|K|}-x}{x^p-x}$ has degree $|K|-p$. Let $$b_1:=\e\bigl(x_1,x_2+f_{V_1}(x_1)v_2\bigr)\in\GA_2(K)\;\;\textup{and}\;\; b_2:=\e\bigl(x_1+f_{V_2}(x_2)v_1,x_2\bigr)\in\GA_2(K).\footnote{Note that $b_2=\shift^{v_1+V_2}_{K\oplus K}$ and $b_1=\e(x_2,x_1)\shift^{v_2+V_1}_{K\oplus K}\e(x_2,x_1)$.}$$
For $Y\in\bigl\{V_1\times V_2,V_1\times (K\setminus V_2),(K\setminus V_1)\times V_2,(K\setminus V_1)\times (K\setminus V_2)\bigr\}$ we have $b_{\iota}(Y)=Y$. Based on this, for 
$$(c_1,c_2):=(b_1b_2,b_2b_1)\in\GA_2(K)^2$$ 
we get that $c_1(K)=c_2(K)$ is a product of $p^{q+d_1-1}+p^{q+d_2-1}-p^{d_1+d_2-1}$ disjoint $p$-cycles, each $p$-cycle having collinear support in one of the three directions $(1:0)$, $(0:1)$, and $(v_1:v_2)$. For $(a_1,a_2,a_3,a_4):=(b_1,b_2,b_1^{-1},b_2^{-1})$ and 
$$a:=a_1a_2a_3a_4\in\GA_2(K),$$ 
$a(K)$ is the identity permutation of $K^2$, we have $\dir(a_i)\neq\dir(a_{i+1})$ for $i\in\{1,2,3\}$, and hence Proposition \ref{PR11}(4.a) gives $\ell(c_1)=\ell(c_2)=(|K|-|V_1|)(|K|-|V_2|)$ and
$$\ell(a)=[(|K|-|V_1|)(|K|-|V_2|)]^2=(p^{2q}-p^{q+d_1}-p^{q+d_2}+p^{d_1+d_2})^2.$$ 
Also, we have $\j_{c_1}=\j_{c_2}=2$: the ``$\le$'' inequality is clear while the ``$\ge$'' inequality follows from the fact that $\ell(c_{\iota})\ge p^q$. 
If $d_1=d_2=q-1$ (resp.\ $d_1=d_2=1$), then the number of disjoint $p$-cycles is $\frac{|K|^2(2p-1)}{p^3}$ (resp.\ is $2|K|-p$).

We have 
$$\n\bigl(c_{\iota}(K)\bigr)=p^{q+d_1}+p^{q+d_2}-p^{d_1+d_2}\le p^{q+q-1}+p^{q+q-1}-p^{2q-2}.$$ Hence $\n\bigl(c_{\iota}(K)\bigr)\le\frac{2|K|^2}{p}-\frac{|K|^2}{p^2}\le |K|^2-|K|$. Thus we have $\n\bigl(c_{\iota}(K)\bigr)=|K|^2-|K|$ iff $|K|=4$. From the last sentence, Lemma \ref{L7}(1.c), and the description of $c_{\iota}(K)$ we get that $\j_{c_{\iota}(K)}\ge 2$. Clearly, $\j_{c_{\iota}(K)}\le 2$. Thus $\j_{c_{\iota}(K)}=2$. Also, $|K|\mid\n\bigl(c_{\iota}(K)\bigr)$ iff $d_1+d_2\ge q$. 
\end{example}

\begin{example}\normalfont\label{EX13}
We assume that $|K|\ge 3$. Let $p:=\char(K)$. Let $l$ be the smallest prime divisor of $|K|-1$ and let $\alpha\in K$ be such that for the solution set $\mathcal S_0$ in $K$ of the equation $x^l-\alpha=0$ we have $|\mathcal S_0|=l$. Let $c:=\e(x_1,x_2+x_1^{|K|-1}-1)$. For $\epsilon\in\{-1,1\}$ let $b_{\epsilon}:=\e\bigl(x_1+\epsilon(x_2^l-\alpha),x_2\bigr)$ and 
$$a_{\epsilon}:=b_{\epsilon}c=\e\bigl(x_1+\epsilon(x_2+x_1^{|K|-1}-1)^l-\epsilon\alpha,x_2+x_1^{|K|-1}-1\bigr).$$
Let $\sigma_{\epsilon}:=a_{\epsilon}(K)$. Clearly, $\j_{\sigma_{\epsilon}}\le 2$. We have $\n(\sigma_{\epsilon})=|K|^2-|\mathcal S|$, where $\mathcal S$ is the solution set in $K^2$ of the system of two equations 
$$x_1+\epsilon(x_2+x_1^{|K|-1}-1)^l-\epsilon\alpha-x_1=0=x_2+x_1^{|K|-1}-1-x_2$$
in the indeterminates $x_1$ and $x_2$. As $\mathcal S=K^{\ast}\times\mathcal S_0$, we have $|S|=l(|K|-1)$. So $\n(\sigma)=|K|^2-l(|K|-1)<|K|^2-|K|$ and $|K|\nmid\n(\sigma)$. From this and Lemma \ref{L7}(1.a) we get that $\j_{\sigma_{\epsilon}}\neq 1$. As $\j_{\sigma_{\epsilon}}\neq 0$ by Lemma \ref{F6}(2), we get that $\j_{\sigma_{\epsilon}}=2$.
\end{example}

\begin{example}\normalfont\label{EX14}
{\bf (1)} From the proof of Proposition \ref{PR10}(1) we get that for each $3$-cycle $\sigma\in\Alt(K^2)$ of non-collinear points there exists a quadruple $(b_1,b_2,b_3,b_4)$ in $(\GA_2(K)[|K|-1])^4$ for which we have identities $\sigma=b_1b_2b_1^{-1}b_2^{-1}=b_3b_4b_3^{-1}b_4^{-1}$ and $\j_{b_1}=\j_{b_2}=\j_{b_3}=\j_{b_4}=1$ and moreover $\dir(b_1)=\dir(b_3)$, $\dir(b_2)$, and $\dir(b_4)$ are the directions given by the three distinct lines of the triangle formed by the three non-collinear points. So for $(a_1,a_2,\ldots,a_8):=(b_1,b_2,b_1^{-1},b_2^{-1},b_4,b_3,b_4^{-1},b_3^{-1})$, we have $\j_{a_1}=\cdots=\j_{a_8}=1$ and the product $\prod_{i=1}^8 a_i(K)$ is the identity permutation.

\smallskip
{\bf (2)} The proof of Proposition \ref{PR10}(1) also gives an analog of the previous paragraph with $8$ replaced by $6$, $3$-cycles replaced by $2p-1$-cycles, and non-collinear supports replaced by supports of the form $\{O,P_1,\ldots,P_{p-1},Q_1,\ldots,Q_{p-1}\}$ that have the property that there exists a pair $(v,w)\in K^2$ of linearly independent vectors such that we have $P_i-O=iv$ and $Q_i-O=iw$ for each $i\in\llbracket1,p-1\rrbracket$.
\end{example}

\begin{proposition}\label{PR15}
Suppose that $|K|$ is not a prime. Then the length functions $\ell^{-}_{2,K}$ and $\ell^{+}_{2,K}$ are not regular.
\end{proposition}

\begin{proof} It suffices to consider the case of $\ell^{-}_{2,K}$.

We write $|K|=p^q$ with $p$ a prime and $q\in\mathbb N^{\ast}\setminus\{1\}$. Let the quadruple $(b_1,b_2,c_1)\in\GA_2(K)^3$ be as in Example \ref{EX12} with $d_1=d_2=1$. Let $b:=b_1^{-1}b_2$. We have $\ell(b_1)=\ell(b_2)=p^q-p$. Moreover, both $b(K)=b_1(K)^{-1}b_2(K)=b_2(K)b_1(K)^{-1}$ and $c_1(K)$ are products of $2|K|-p$ disjoint $p$-cycles whose supports have three distinct directions. If $|K|=4$ we have $\n\bigl(b(K)\bigr)=\n\bigl(c_1(K)\bigr)=|K|^2-|K|=12$ and if $|K|\neq 3$ we have $\n\bigl(b(K)\bigr)=\n\bigl(c_1(K)\bigr)<|K|^2-|K|=12$ by Example \ref{EX12}.

Let $\sigma_a:=b_1(K)$, $\sigma_b:=b(K)$, and $\sigma_c:=b_2(K)$. We check that the three identities $\ell^{-}_{2,K}(\sigma_a)=2$, $\ell^{-}_{2,K}(\sigma_b)=4$, and $\ell^{-}_{2,K}(\sigma_c)=2$ hold. The ``$\le$'' inequalities are clear. By the symmetry between $b_1(K)$ and $b_2(K)$, to prove the ``$\ge$'' inequalities it suffices to show that $\ell_{2,K}^{-}\bigl(c_1(K)\bigr)=$ and hence that $\j_{c_1(K)}=2$. Based on the description of $c_1(K)$ in Example \ref{EX12}, we have $\j_{c_1(K)}\neq 1$ by Lemma \ref{L7}(1) and $\j_{c_1(K)}\neq 0$ by Lemma \ref{F6}(2) and (3). Thus $\j_{c_1(K)}\ge 2$. As $c_1=b_1b_2$, we have $\j_{c_1(K)}\le 2$. So $\j_{c_1(K)}=2$ and the three identities hold.

Based on the following identities $\sigma_a\sigma_b=b_2(K)=\sigma_c$, $\sigma_b^{-1}\sigma_c=b_1(K)=\sigma_a$, and $\sigma_a\sigma_c=b_1(K)b_2(K)=c_1(K)$, we have $\ell^{-}_{2,K}(\sigma_a\sigma_b)=2<8=\ell^{-}_{2,K}(\sigma_a)\ell^{-}_{2,K}(\sigma_b)$, $\ell^{-}_{2,K}(\sigma_b^{-1}\sigma_c)=2<8\le\ell^{-}_{2,K}(\sigma_b^{-1})\ell^{-}_{2,K}(\sigma_c)$, and $\ell^{-}_{2,K}(\sigma_a\sigma_c)=4=\ell^{-}_{2,K}(\sigma_a)\ell^{-}_{2,K}(\sigma_c)$. Hence $\ell^{-}_{2,K}$ is not regular.\end{proof}

\begin{corollary}\label{C14}
Suppose that $|K|=4$. Then $\ell_{2,K}$ is not regular.
\end{corollary}

\begin{proof}
We use the proof of Proposition \ref{PR15}. We have $\ell(b_1)=\ell(b_2)=|K|-2=2$. 

The corollary follows from the proof of Proposition \ref{PR15} once we check that we have $\ell_{2,K}(\star)=\ell_{2,K}^-(\star)$ for each $\star\in\{\sigma_a,\sigma_b,\sigma_c,\sigma_a\sigma_b,\sigma_b^{-1}\sigma_c,\sigma_a\sigma_c\}$. Proposition \ref{PR13}(3) gives that this holds if $\star\in\{\sigma_a,\sigma_c,\sigma_a\sigma_b,\sigma_b^{-1}\sigma_c\}$. If $\star\in\{\sigma_b,\sigma_a\sigma_c\}$, then $\j_{\star}=2$ by the proof of Proposition \ref{PR15} and $\ell_{2,K}(\star)\le 4=\ell(b)=\ell(c_1)$. Hence $\ell_{2,K}(\star)$ is a product of at most $j_{\star}$ primes. From this and Proposition \ref{PR14}(2.b) we get that $\ell_{2,K}(\star)$ is a product of $2$ primes and hence we must have $\ell_{2,K}(\star)=4$.
\end{proof}

For refinement purposes in what follows for $n=2$ we use the notation of the following second variant of Definition \ref{D2}(2) and (3).

\begin{definition}\label{D17}
Let $K$ be a finite field. Let $m\in\llbracket1,|K|^2\rrbracket$ if $4\nmid |K|$ and let $m\in\llbracket1,|K|^2-2\rrbracket$ if $4\nmid |K|$. Let $j\in\llbracket0,\j_K\rrbracket$. 

\medskip
{\bf (1)} For $(\underline{P},\underline{Q})\in\mathbb D_{2,m}(K)^2$, let $\pi^{\le j}_{\underline{P},\underline{Q}}$ (resp.\ $\pi^{\S,\le j}_{\underline{P},\underline{Q}}$) in $\mathbb N^{\ast}\cup\{\infty\}$ be the infimum of the set 
$$\{\ell(a)|a\in\GA_2(K), a(\underline{P})=\underline{Q},\j_a\le j\}$$ 
(resp.\ $\{\ell^{\S}(a)|a\in\SGA_2(K), a(\underline{P})=\underline{Q},\j_a\le j\}$).

\smallskip
{\bf (3)} Let $\pi^{\le j}_{2,m}(K)\in\mathbb N^{\ast}\cup\{\infty\}$\index{$\pi^{\le j}_{2,m}(K)$ main invariant} (resp.\ $\pi^{\S,\le j}_{2,m}(K)\in\mathbb N^{\ast}\cup\{\infty\}$\index{$\pi^{\S,\le j}_{2,m}(K)$ main invariant}) be the supremum of the set $\{\pi^{\le j}_{\underline{P},\underline{Q}}|(\underline{P},\underline{Q})\in\mathbb D_{2,m}(K)^2\}$ (resp.\ $\{\pi^{\S,\le j}_{\underline{P},\underline{Q}}|(\underline{P},\underline{Q})\in\mathbb D_{2,m}(K)^2\}$).
\end{definition}

\section{Estimating length functions via conjugacy classes}\label{S13}

In this section we study the variation of the values of $\ell_{n,K}$ and $\ell_{n,K}^{\S}$ on unions of conjugacy classes in order to obtain upper bounds for these values. 

\begin{notation}\normalfont\label{N4}
For a non-empty subset $\nabla$ of $\Perm(K^n)$, let 
$$\omega(\nabla):=\min\left(\ell_{n,K}(\sigma)|\sigma\in \nabla\right)\in\mathbb N^{\ast}\;\;\;\textup{and}\;\;\;\Omega(\nabla):=\max\bigl(\ell_{n,K}(\sigma)|\sigma\in\nabla\bigr)\in\mathbb N^{\ast}.$$
If $\nabla\subset\Alt(K^n)$, we also define
$$\omega^{\S}(\nabla):=\min\left(\ell^{\S}_{n,K}(\sigma)|\sigma\in \nabla\right)\in\mathbb N^{\ast}\;\;\;\textup{and}\;\;\;\Omega^{\S}(\nabla):=\max\bigl(\ell^{\S}_{n,K}(\sigma)|\sigma\in\nabla\bigr)\in\mathbb N^{\ast}.$$
For $s\in \llbracket2,|K|^n\rrbracket$ and $t\in\bigl\llbracket1,\lfloor\frac{|K|^n}{s}\rfloor\bigr\rrbracket$, let $t\mathcal Y_s=t\mathcal Y_{s,n,K}$ be the conjugacy class of $\perm(K^n)$ formed by $t$ disjoint $s$-cycles. Let $\mathcal Y_s:=1\mathcal Y_s$. If $3\le s\le n+1$ (resp.\ $2\le s\le |K|$), let $\mathcal Y_s^{\d=s-1}$ (resp.\ $\mathcal Y_s^{\d=1}$) be the subset of $\mathcal Y_s$ formed by $s$-cycles whose support is affinely independent (resp.\ collinear). Similarly, for $4\le s\le n+1$ (resp.\ $3\le s\le |K|$) let $\mathcal Y_s^{\d\le s-2}$ (resp.\ $\mathcal Y_s^{\d\ge 2}$) be the subset of $\mathcal Y_s$ formed by $s$-cycles whose support is affinely dependent (resp.\ non-collinear). If $s$ is odd, let $\Pi^{\S}_{n,s\textup{-cycle}}(K):=\{\ell^{\S}_{n,K}(\sigma)|\sigma\in\mathcal Y_s\}$.
\end{notation}

If $\nabla\subset\Alt(K^n)$, then we have inequalities $\omega(\nabla)\le\omega^{\S}(\nabla)$ and $\Omega(\nabla)\le\Omega^{\S}(\nabla)$. 

If $|K|=2$, then $\STGA_n(K)=\TGA_n(K)$ and hence $\ell_{n,K}=\ell_{n,K}^{S}$; in particular, for $\nabla\subset\Alt(K^n)$ we have $\omega(\nabla)=\omega^{\S}(\nabla)$ and $\Omega(\nabla)=\Omega^{\S}(\nabla)$.

\begin{remark}\normalfont\label{R6}
Referring to (QP2) of Section \ref{S1} with $s$ odd if $4\mid |K|$, we have an inclusion $\Pi_{n,s\textup{-cycle}}(K)\subset \llbracket\omega(\mathcal Y_s),\Omega(\mathcal Y_s)\rrbracket$ and identities $\omega(\mathcal Y_s)=\min\bigl(\Pi_{n,s\textup{-cycle}}(K)\bigr)$ and $\Omega(\mathcal Y_s)=\max\bigl(\Pi_{n,s\textup{-cycle}}(K)\bigr)$. Similarly, if $s$ is odd, we have an inclusion $\Pi^{\S}_{n,s\textup{-cycle}}(K)\subset \llbracket\omega^{\S}(\mathcal Y_s),\Omega^{\S}(\mathcal Y_s)\rrbracket$ and identities $\omega^{\S}(\mathcal Y_s)=\min\bigl(\Pi^{\S}_{n,s\textup{-cycle}}(K)\bigr)$ and $\Omega^{\S}(\mathcal Y_s)=\max\bigl(\Pi^{\S}_{n,s\textup{-cycle}}(K)\bigr)$.
\end{remark}

\begin{example}\normalfont\label{EX15}
{\bf (1)} We have $\omega(\mathcal Y_{|K|^n-1})=1$ as there exists $a\in\GL_n(K)$ such that $a(K)$ is a $(|K|^n-1)$-cycle. For instance, if $L$ is a finite field extension of $K$ with $[L:K]=n$, by identifying $L=K^n$ as $K$-vector spaces, we can take $a$ such that $a(K)\in\Perm(K^n)$ is the multiplication by a generator of the multiplicative cyclic group $L^{\ast}$.

\smallskip
{\bf (2)} Let $p:=\char(K)$. Let $b:=\e(x_1+1,x_2-x_1^{p-1}+1,\ldots,x_n-x_{n-1}^{p-1}+1)$ in $\STGA_n(K)$ with $n\in\mathbb N^{\ast}\setminus\{1\}$. Then $b(K)\in\Alt(K^n)$ is a product of $\frac{|K|^n}{p^n}$ disjoint $p^n$-cycles and indexed by elements of the quotient group $(K/\mathbb F_p)^n$ and therefore $\omega\bigl(\frac{|K|^n}{p^n}\mathcal Y_{p^n}\bigr)\le\omega^{\S}\bigl(\frac{|K|^n}{p^n}\mathcal Y_{p^n}\bigr)\le \ell(a)=(p-1)^{n-1}$, where the equality holds by Example \ref{EX9} applied to $r=p-1$ and $a:=\e(x_1-1,x_2,\ldots,x_n)b$.
\end{example}

\begin{definition}\label{D18}
Let $\nabla$ be a non-empty subset of $\Perm(K^n)$ (resp.\ of $\Alt(K^n)$). We call the real number $\eth(\nabla):=\sqrt{\frac{\Omega(\nabla)}{\omega(\nabla)}}\in [1,\infty)$ (resp.\ $\eth^{\S}(\nabla):=\sqrt{\frac{\Omega^{\S}(\nabla)}{\omega^{\S}(\nabla)}}\in [1,\infty)$) the augmentation factor\index{augmentation factor} of $\nabla$.
\end{definition}

\begin{lemma}\label{L8} 
Let $n\in\mathbb N^{\ast}\setminus\{1\}$. Then the following properties hold.

\medskip
{\bf (1)} Let $\mathcal C$ be a non-trivial conjugacy class of $\Perm(K^n)$ or $\Alt(K^n)$. Then we have an inequality $\eth(\mathcal C)\le\pi_{n,\n(\mathcal C)}(K)$.

\smallskip
{\bf (2)} Let $\mathcal C$ be a non-trivial conjugacy class of $\Alt(K^n)$. Then we have inequalities $\eth(\mathcal C)\le\pi^{\E}_{n,\n(\mathcal C)}(K)\le\pi^{\S}_{n,\n(\mathcal C)}(K)$.
\end{lemma}

\begin{proof}
For part (1), let the pair $(\sigma,\theta)\in\mathcal C^2$ be such that $\ell_{n,K}(\sigma)=\omega(\mathcal C)$ and $\ell_{n,K}(\theta)=\Omega(\mathcal C)$. Let $\varsigma\in\Perm(K^n)$ be such that $\theta=\varsigma\sigma\varsigma^{-1}$. By listing the elements of $\supp(\sigma)=\{P_i|i\in \llbracket1,\n(\sigma)\rrbracket\}$, we have $\supp(\theta)=\{\varsigma(P_i)|i\in \llbracket1,\n(\sigma)\rrbracket\}$. 

Let $(a,b)\in\TGA_n(K)[\omega(\mathcal C)]\times\TGA_n(K)[\pi_{n,\n(\mathcal C)}(K)]$ be such that $a(K)=\sigma$ and $b(P_i)=\varsigma(P_i)$ for each $i\in \llbracket1,\n(\sigma)\rrbracket$. For 
$$c:=bab^{-1}\in \TGA_n(K)[\omega(\mathcal C)\pi_{n,\n(\mathcal C)}(K)^2]$$ 
we have $c(K)=\theta$, thus $\ell_{n,K}(\theta)\le \omega(\mathcal C)\pi_{n,\n(\mathcal C)}(K)^2$. Hence $\frac{\Omega(\mathcal C)}{\omega(\mathcal C)}\le \pi_{n,\n(\mathcal C)}(K)^2$ from which part (1) follows.

The proof of part (2) is the same as of part (1), with $\pi_{n,\n(\mathcal C)}(K)$ replaced by $\pi^{\E}_{n,\n(\mathcal C)}(K)$ or $\pi^{\S}_{n,\n(\mathcal C)}(K)$ and with the quadruple $\bigl(\TGA_n(K),\Omega(\mathcal C),\omega(\mathcal C),\ell_{n,K}\bigr)$ replaced by $\bigl(\STGA_n(K),\Omega^{\S}(\mathcal C),\omega^{\S}(\mathcal C),\ell^{\S}_{n,K}\bigr)$.
\end{proof}

\begin{lemma}\label{F9} Let $n\in\mathbb N^{\ast}\setminus\{1\}$. Then the following properties hold.

\medskip
{\bf (1)} Let $m\in \llbracket2,n+1\rrbracket$. If $4\mid |K|$, then we also assume that $m$ is odd. Then there exists $\kappa^{\d=m-1}_{n,K;m}\in\mathbb N^{\ast}$ (resp. $\kappa^{\S,\d=m-1}_{n,K;m}\in\mathbb N^{\ast}$) such that we have an identity $\ell_{n,K}(\mathcal Y_m^{\d=m-1})=\{\kappa^{\d=m-1}_{n,K;m}\}$ (resp.\ $\ell^{\S}_{n,K}(\mathcal Y_m^{\d=m-1})=\bigl\{\kappa^{\S,\d=m-1}_{n,K;m}\bigr\}$). 

\smallskip
{\bf (2)} Suppose that $|K|=3$. Then there exists $\kappa^{\d=1}_{n,K;3}\in\mathbb N^{\ast}$ (resp.\ $\kappa^{\S,\d=1}_{n,K;3}\in\mathbb N^{\ast}$) such that we have an identity $\ell_{n,K}(\mathcal Y_3^{\d=1})=\{\kappa^{\d=1}_{n,K;3}\}$ (resp.\ $\ell^{\S}_{n,K}(\mathcal Y_3^{\d=1})=\{\kappa^{\S,\d=1}_{n,K;3}\}$).

\smallskip
{\bf (3)} Suppose that $|K|=2$. Then there exist positive integers $\kappa_{n,K;3}$ and $\kappa^{\d=2}_{n,K;4}$ such that we have identities $\ell_{n,K}(\mathcal Y_3)=\{\kappa_{n,K;3}\}$ and $\ell_{n,K}(\mathcal Y_4^{\d=2})=\{\kappa^{\d=2}_{n,K;4}\}$.
\end{lemma}

\begin{proof}
For $(a,b)\in\AGL_n(K)\times \STGA_n(K)$ we have $aba^{-1}\in\STGA_n(K)$ and $\ell(b)=\ell(aba^{-1})$. Based on this, parts (1) and (2) follow from Lemma \ref{F7}(1) and (2) (respectively).

Part (3) follows directly from Lemma \ref{F7}(3.b) and (3.c).
\end{proof}

\begin{proposition}\label{PR16}
Suppose that $4\nmid |K|$ and $n\in\mathbb N^{\ast}\setminus\{1\}$.. Then the following properties hold.

\medskip
{\bf (1)} There exists $\kappa_{n,K}\in\mathbb N^{\ast}$ such that $\ell_{n,K}(\mathcal Y_2)=\{\kappa_{n,K}\}$.

\smallskip
{\bf (2)} Let $\mathcal C$ be a non-trivial conjugacy class of $\perm(K^n)$. Then $\Omega(\mathcal C)\le \kappa_{n,K}^{\nu_{2,1}(\mathcal C)}$.

\smallskip
{\bf (3)} For each $s\in \llbracket3,|K|^n\rrbracket$ we have an inequality $\Omega(\mathcal Y_s)\le \kappa_{n,K}^{s-1}$.

\smallskip
{\bf (4)} We have $\kappa_{n,\mathbb F_2}=n-1$ and $\kappa_{n,\mathbb F_2;3}\le 2n-3$.

\smallskip
{\bf (5)} Suppose that $K=\mathbb F_2$. The for each $s\in \llbracket3,2^n\rrbracket$ we have inequalities 
$$\Omega(\mathcal Y_s)\le\Omega(\mathcal Y_2^{s-1})\le (2n-3)^{\lfloor\frac{s-1}{2}\rfloor}(n-1)^{s-1-2\lfloor\frac{s-1}{2}\rfloor}.$$

{\bf (6)} Let $p$ be an odd prime and $q\in\mathbb N^{\ast}$. Then $\kappa_{n,\mathbb F_{p^q}}^{p^{q-1}(p-1)}\ge (n-1)(|K|-1)$ and hence $\lim_{n\rightarrow\infty} \kappa_{n,\mathbb F_{p^q}}=\infty$.

\smallskip
{\bf (7)} We have $\kappa_{n,\mathbb F_3;3}^{\d=1}=\kappa_{n,\mathbb F_3;3}^{\S,\d=1}=2(n-1)$ and $\kappa_{n,\mathbb F_3;3}^{\S,\d=2}\le 4n$.
\end{proposition}

\begin{proof}
As $\mathcal Y_2=\mathcal Y_2^{\d=1}$, by taking $\kappa_{n,K}:=\kappa_{n,K;2}^{\d=1}$ (see Lemma \ref{F9}(1)), we get that part (1) holds. 

For part (2), let $\sigma\in\mathcal C$ be such that $\ell_{n,K}(\sigma)=\Omega(\mathcal C)$. Writing $\sigma=\prod_{i=1}^{\nu_{2,1}(\mathcal C)} \tau_i$ as a product of transpositions, we have $\ell_{n,K}(\sigma)\le\prod_{i=1}^{\nu_{2,1}(\mathcal C)}\ell_{n,K}(\tau_i)=\kappa_{n,K}^{\nu_{2,1}(\mathcal C)}$ by Axiom $A3_{n,K}$. So part (2) holds.

Part (3) follows from part (2) applied to $\mathcal C=\mathcal Y_s$ and $\nu_{2,1}(\mathcal Y_s)=s-1$.

For part (4), for $a=\e\bigl(x_1+\prod_{i=2}^n (1+x_i),x_2,\ldots,x_n\bigr)\in\TGA_n(\mathbb F_2)$ we have $a(\mathbb F_2)=(O\; O_1)$ where $O:=(0,\ldots,0)$ and $O_1:=(1,0,\ldots,0)$ (cf.\ Example \ref{EX7}). Clearly, $\ell(a)=n-1$. As $\ell_{n,K}(O\; O_1)=\ell(a)$ by Proposition \ref{PR13}(3), we get that $\kappa_{n,\mathbb F_2}=n-1$. The inequality $\kappa_{n,\mathbb F_2;3}\le 2n-3$ follows from Example \ref{EX7} if $n\ge 3$ and from $\kappa_{2,\mathbb F_2}=1$ if $n=2$. So part (4) holds. 

For $n=2$, as $\kappa_{2,\mathbb F_2}=1$ by part (4), part (5) follows from part (2). So to prove part (5) we can assume that $n\ge 3$. 

If $\sigma\in 2\mathcal Y_2$ has affinely dependent (resp.\ affinely independent) support, then we have $\ell_{n,\mathbb F_2}(\sigma)=n-2$ (resp.\ $\ell_{n,\mathbb F_2}(\sigma)\le 2n-3$). To check this, up to conjugation with affine automorphisms by Lemma \ref{F7}(3.c) (resp.\ \ref{F7}(1)), we can assume that $\sigma$ is $a_2(\mathbb F_2)$ (resp.\ $a_3(\mathbb F_2)$) of Example \ref{EX7} and the equality (resp.\ inequality) follows from Proposition \ref{PR13}(6) and Example \ref{EX7} (resp.\ from Example \ref{EX7}). 

As $\mathcal Y_2^2=2\mathcal Y_2\sqcup\mathcal Y_3$, from the prior paragraph and $\kappa_{n,\mathbb F_2;3}\le 2n-3$ (see part (4)) we get inequalities $\omega(\mathcal Y_2^2)\le n-2$ and $\Omega(\mathcal Y_2^2)\le 2n-3$. As $\kappa_{n,\mathbb F_2}=n-1$ by part (4) and $\Omega(\mathcal Y_2^2)\le 2n-3$, by induction on $s\in \llbracket3,2^n\rrbracket$ we get the second inequality $\Omega(\mathcal Y_2^{s-1})\le (2n-3)^{\lfloor\frac{s-1}{2}\rfloor}(n-1)^{-1-2\lfloor\frac{s-1}{2}\rfloor}$. As $\mathcal Y_s\subset\mathcal Y_2^{s-1}$, the first inequality of part (5) also holds. So part (5) holds.

For part (6), we take $\mathcal C:=p^{q-1}\mathcal Y_p$. For $a:=\e\bigl(x_1+\prod_{i=2}^n(1-x_i^{|K|-1}),x_2,\ldots,x_n\bigr)$ in $\STGA_n(\mathbb F_{p^q})$ we have $\sigma=a(\mathbb F_{p^q})\in\mathcal C$ by Lemma \ref{F5}(1) and (2) and identities $\ell_{n,K}(\sigma)=\ell(a)=(n-1)(|K|-1)$ by Proposition \ref{PR13}(6). Thus, as we have $\nu_{2,1}(\mathcal C)=p^{q-1}(p-1)$, part (6) follows from part (2).

For part (7), its identities follow from the fact that there exists $\sigma\in\mathcal Y_3^{\d=1}$ with $\ell^{\S}_{n,K}(\sigma)=\ell_{n,K}(\sigma)=n-1$ by Proposition \ref{PR13}(6) applied to $d=n-1$ and its inequality follows from Proposition \ref{PR10}(2).
\end{proof}

\begin{lemma}\label{L8.5} If $|K|$ is odd, then we have inequalities 
$$\max\bigl(\omega(\mathcal Y_{|K|^n-2}),\omega(\mathcal Y_{|K|^n})\bigr)\le\kappa_{n,K}\le\min\bigl(\Omega(\mathcal Y_{|K|^n-2}),\Omega(\mathcal Y_{|K|^n})\bigr)\le\pi_{n,|K|^n-2}.$$
\end{lemma}

\begin{proof} Let $a\in\GL_n(K)$ be such that $a(K)\in\mathcal Y_{|K|^n-1}$ by Example \ref{EX15}(1). For $\iota\in\{0,2\}$ we consider a transposition $\tau_{\iota}\in\Perm(K^n)$ such that for $\theta_{\iota}:=a(K)\tau_{\iota}$ we have $\theta_{\iota}\in\mathcal Y_{|K|^n-2+\iota}\subset\Alt(K^n)$. We get that 
$$\kappa_{n,K}=\ell_{n,K}(\tau_{\iota})=\ell_{n,K}(\theta_{\iota})\in\Pi_{n,|K|^n-2+\iota\textup{-cycle}}(K)$$ and from this the lemma follows.
\end{proof}

We exemplify how to get smaller upper bounds than Lemma \ref{L8.5} for $\kappa_{n,K}$ if $|K|\equiv 3\; (\textup{mod}\; 4)$.

\begin{example}\normalfont\label{EX16}
Let $n\in\mathbb N^{\ast}\setminus\{1\}$, $k:=|K|-1$, and the linear automorphism $b:=\e(-x_1,x_2,\ldots,x_n)\in\GL_n(K)$. If $|K|$ is odd, then $b(K)=\prod_{i=1}^s\tau_i$ is a product of $s:=\frac{|K|^{n-1}(|K|-1)}{2}=\frac{k(k+1)^{n-1}}{2}$ disjoint transpositions. 

\medskip
{\bf (1)} Assume that $|K|\equiv 3\; (\textup{mod}\; 4)$. Then $\frac{k}{2}\equiv 1\pmod{2}$ and thus $s$ is odd. Hence $\theta:=\prod_{i=2}^{s} \tau_i$ is even. We write $\theta=\prod_{i=2}^{s} \theta_i$ as a product of $s-1$ permutations that are $3$-cycles by Lemma \ref{P3}(3); we can assume that all these $3$-cycles have non-collinear supports. As $\tau_1=\theta b(K)^{-1}$, Proposition \ref{PR10}(2) and (3) gives for $|K|\ge 7$ that
$$\kappa_{n,K}=\ell_{n,K}(\tau_1)\le (E_{n-2,k})^{\frac{k(k+1)^{n-1}}{2}-1}$$ 
and for $|K|=3$, i.e., $k=2$, that $\kappa_{n,\mathbb F_3}=\ell_{n,K}(\tau_1)\le (4n)^{3^{n-1}-1}$.

\smallskip
{\bf (2)} The last inequality can be improved for $n\in\{2,3\}$ as follows.
Assume that $|K|=3$, i.e., $k=2$. Let $l:=\frac{s-1}{2}=\frac{3^{n-1}-1}{2}\in\mathbb N^{\ast}$; we can assume that the numbering of the $\tau_i$s is such that $\theta=\prod_{i=1}^l (\tau_{2i}\tau_{2i+1})$ is a product of $l$ permutations that are products of two disjoint transpositions of the form $(Q_1\;Q_2)(Q_3\;Q_4)$ with $Q_1+Q_4=Q_2+Q_3$. From Example \ref{EX8}(1) we get that $\ell_{2,K}(\tau_{2i}\tau_{2i+1})\le 8(2n-2)^2$. Based on this and Axiom ($A3_{n,K}$) we get that
$$\kappa_{n,\mathbb F_3}=\ell_{2,K}(\tau_1)\le [2\sqrt{2}(2n-2)]^{3^{n-1}-1}.$$
For $n\in\{2,3\}$ we get the better estimates $\kappa_{2,\mathbb F_3}\le 32=2^5$ and $\kappa_{3,\mathbb F_3}\le 128^4=2^{28}$.\footnote{First author and Inna Sysoeva checked using a Fortran code that $\kappa_{2,\mathbb F_3}=32$; thus for $\tau\in\mathcal Y_2$ we have $\j_{\tau}=5$. Note that Corollary \ref{C13.5} and Example \ref{EX16}(2) only give that $\j_{\tau}\in\{3,4,5\}$.\label{foot13}}\end{example}

The next general lemma is the essence of our estimations of values of length functions via unions of conjugacy classes.

\begin{lemma}\label{L9} 
Let $(n,r)\in (\mathbb N^{\ast}\setminus\{1\})\times\mathbb N^{\ast}$. If $|K|=2$, then we assume that $n\ge 3$. Let $(\mathcal C_i)_{i\in \llbracket1,r\rrbracket}$ be a sequence of distinct non-trivial conjugacy classes of $\Alt(K^n)$. For $\mathcal U:=\cup_{i=1}^r \mathcal C_i$ let $\n(\mathcal U)$ and $\mu_t(\mathcal U)$ with $t\in\llbracket\max\bigl(5-\n(\mathcal U),0\bigr),|K|^n-\n(\mathcal U)\rrbracket$ be as in Theorem \ref{P12}(7). If $4\mid |K|$ we assume that $\n(\mathcal U)\le |K|^n-2$. Let 
$$\omega_0(\mathcal U):=\max\bigl(\omega(\mathcal C_i)|i\in\llbracket1,r\rrbracket\bigr),\;\;\;\omega^{\S}_0(\mathcal U):=\max\bigl(\omega^{\S}(\mathcal C_i)|i\in\llbracket1,r\rrbracket\bigr),$$
$$\eth_0(\mathcal U):=\max\bigl(\eth(\mathcal C_i)|i\in\llbracket1,r\rrbracket\bigr),\;\;\textup{and}\;\;\eth^{\S}_0(\mathcal U):=\max\bigl(\eth^{\S}(\mathcal C_i)|i\in\llbracket1,r\rrbracket\bigr).$$
Then the following properties hold.

\medskip
{\bf (1)} Let $\varkappa\in\Perm(K^n)$. Let $t(\varkappa)\in \llbracket0,|K|^n-\n(\varkappa)\rrbracket$ be the smallest such that we have $\n(\varkappa)+t(\varkappa)\ge\max\bigl(5,\n(\mathcal U)\bigr)$. Then we have the following properties.

\medskip\noindent
{\bf (1.a)} Suppose that $\varkappa\in\Alt(K^n)$. Let $\nu_{\mathcal U}(\varkappa)\in \llbracket1,\mu_{\n(\varkappa)+t(\varkappa)-\n(\mathcal U)}(\mathcal U)\rrbracket$ be the smallest such that $\varkappa=\prod_{i=1}^{\nu_{\mathcal U}(\varkappa)} \theta_i$ with each $\theta_i\in\mathcal U$. Then we have inequalities
\begin{equation}\label{EQ15}
\ell_{n,K}(\varkappa)\le [\omega_0(\mathcal U)\eth_0(\mathcal U)^2]^{\nu_{\mathcal U}(\varkappa)}\le [\omega_0(\mathcal U)\eth_0(\mathcal U)^2]^{\mu_{\n(\varkappa)+t(\varkappa)-\n(\mathcal U)}(\mathcal U)}
\end{equation}
and
\begin{equation}\label{EQ16}
\ell^{\S}_{n,K}(\varkappa)\le [\omega_0^{\S}(\mathcal U)\eth_0^{\S}(\mathcal U)^2]^{\nu_{\mathcal U}(\varkappa)}\le [\omega_0^{\S}(\mathcal U)\eth_0^{\S}(\mathcal U)^2]^{\mu_{\n(\varkappa)+t(\varkappa)-\n(\mathcal U)}(\mathcal U)}.
\end{equation}

\noindent
{\bf (1.b)} Suppose that $|K|$ is odd and $\varkappa$ is odd. Then we have an inequality
\begin{equation*}
\ell_{n,K}(\varkappa)\le [\omega_0(\mathcal U)\eth_0(\mathcal U)^2]^{\mu_{\max\bigl(\n(\mathcal U),\min(|K|^n,|K|^n-|K|^{n-1}+\n(\varkappa)-2\bigr)-\n(\mathcal U)}(\mathcal U)}.
\end{equation*}

\medskip\noindent
{\bf (1.c)} Suppose that $|K|$ is odd, $\varkappa$ is odd, and $2\n(\mathcal U)\le |K|^n$. Then we have inequalities
\begin{equation*}
\ell_{n,K}(\varkappa)\le \pi^{\E}_{n,\n(\mathcal U)}(K)[\omega_0(\mathcal U)\eth_0(\mathcal U)^2]^{\mu_{|K|^n-2\n(\mathcal U)}(\mathcal U)}.
\end{equation*}

\smallskip\noindent
{\bf (1.d)} Suppose that $\varkappa\in\Alt(K^n)$ and $2\n(\mathcal U)\le |K|^n$. If $(n,|K|)=(3,2)$ then we also assume that $\mathcal U=\mathcal Y_3$ (so $\n(\mathcal U)=3$). Then we have an inequality
\begin{equation*}
\ell_{n,K}(\varkappa)\le \pi^{\E}_{n,\n(\mathcal U)}(K)[\omega_0(\mathcal U)\eth_0(\mathcal U)^2]^{\mu_{|K|^n-2\n(\mathcal U)}(\mathcal U)}
\end{equation*}
and
\begin{equation*}
\ell^{\S}_{n,K}(\varkappa)\le \pi^{\S}_{n,\n(\mathcal U)}(K)[\omega_0^{\S}(\mathcal U)\eth_0^{\S}(\mathcal U)^2]^{\mu_{|K|^n-2\n(\mathcal U)}(\mathcal U)}.
\end{equation*}

\noindent
{\bf (1.e)} Suppose that $|K|=2$ and $\varkappa$ is odd. If $\c_2(\varkappa)=0$ (resp.\ $\c_2(\varkappa)\ge 1$) then we have an inequality
\begin{equation*}
\ell_{n,K}(\varkappa)\le (n-1)[\omega_0(\mathcal U)\eth_0(\mathcal U)^2]^{\mu_{\n(\varkappa)+t(\varkappa)-\epsilon-\n(\mathcal U)}(\mathcal U)},
\end{equation*}
where $\epsilon\in\{0,1\}$ (resp.\ $\epsilon\in\{0,1,2\}$) is $1$ iff $t(\varkappa)=0$ and $\n(\varkappa)-1\ge \max\bigl(5,\n(\mathcal U)\bigr)$ (resp.\ is positive iff f $t(\varkappa)=0$ and $\n(\varkappa)-1\ge \max\bigl(5,\n(\mathcal U)\bigr)$, in which case it is the largest such that $\n(\varkappa)-\epsilon\ge \max\bigl(5,\n(\mathcal U)\bigr)$). 

\smallskip
{\bf (2)} Let $m\in \llbracket2,|K|^n\rrbracket$. Let $\mu_{m,\mathcal U}\in\mathbb N^{\ast}$ be the smallest such that 
$$\{\varkappa\in\Alt(K^n)|\n(\varkappa)\le m\}\subset\mathcal U^{\mu_m,\mathcal U}.$$ Then for each $\varkappa\in\Alt(K^n)$ with $\n(\varkappa)\le m$, 
we have inequalities
\begin{equation}\label{EQ17}
\ell_{n,K}(\varkappa)\le [\omega_0(\mathcal U)\eth_0(\mathcal U)^2]^{\mu_{m,\mathcal U}}
\end{equation}
and
\begin{equation}\label{EQ18}
\ell^{\S}_{n,K}(\varkappa)\le [\omega^{\S}_0(\mathcal U)\eth_0^{\S}(\mathcal U)^2]^{\mu_{m,\mathcal U}}.
\end{equation}
Moreover, for each $l\in\llbracket1,|K|^n-2\rrbracket$ we have inequalities
\begin{equation}\label{EQ19}
\ell_{n,K}(\varkappa)\le \pi^{\E}_{n,l}(K)[\omega_0(\mathcal U)\eth_0(\mathcal U)^2]^{\mu_{|K|^n-l,\mathcal U}},
\end{equation}
and
\begin{equation}\label{EQ20}
\ell^{\S}_{n,K}(\varkappa)\le \pi^{\S}_{n,l}(K)[\omega^{\S}_0(\mathcal U)\eth_0^{\S}(\mathcal U)^2]^{\mu_{K|^n-l,\mathcal U}}.
\end{equation}
\end{lemma}

\begin{proof}
For part (1.a), for each $i\in \llbracket1,\nu_{\mathcal U}(\varsigma)\rrbracket$ let $j_i\in\llbracket1,r\rrbracket$ be such that $\theta_i\in\mathcal C_{j_i}$. 

We first note that $\ell_{n,K}(\theta_i)\le\Omega(\mathcal C_{j_i})=\omega(\mathcal C_{j_i})\eth(\mathcal C_{j_i})^2\le\omega_0(\mathcal U)\eth_0(\mathcal U)^2$ and that $\ell^{\S}_{n,K}(\theta_i)\le\Omega^{\S}(\mathcal C_{j_i})=\omega(\mathcal C_{j_i})\eth(\mathcal C_{j_i})^2\le\omega^{\S}_0(\mathcal U)\eth^{\S}_0(\mathcal U)^2$. Based on these, Inequalities (\ref{EQ15}) and (\ref{EQ16}) follow from Axioms $A3_{n,K}$ and $A3^{\S}_{n,K}$. So part (1.a) holds.

For parts (1.b), (1.c), and (1.e), $\varkappa$ is odd and we consider $a\in\TGA_n(K)$ such that $a(K)$ is odd and as in Theorem \ref{T6}(1); thus $a(K)\varkappa$ is even. 

For part (1.b) we have $a\in\AGL_n(K)$ with $\n\bigl(a(K)\bigr)=|K|^n-|K|^{n-1}$ and hence we can choose $a$ such that 
$\n\bigl(a(K)\varkappa\bigr)\le\min\bigl(|K|^n,|K|^n-|K|^{n-1}+\n(\varkappa)-2\bigr)$.
As $\n(\varkappa)\ge 2$ and $|K|\ge 3$, we have 
$$\min\bigl(|K|^n,|K|^n-|K|^{n-1}+\n(\varkappa)-1\bigr)\ge |K|^n-|K|^{n-1}\ge 6>5.$$
This implies that 
$$\n\bigl(\a(K)\varkappa\bigr)+t\bigl(\a(K)\varkappa\bigr)\le\max\bigl(\n(\mathcal U),\min(|K|^n,|K|^n-|K|^{n-1}+\n(\varkappa)-2)\bigr).$$ 
Thus from Equation (\ref{EQ15}) applied to $a(K)\varkappa$ we get that part (1.b) holds.

For part (1.c), we have $a\in\AGL_n(K)$ and let $b\in\TGA_n(K)[\pi^{\E}_{n,\n(\mathcal U)}(K)]$ with $\n\bigl(b(K)a(K)\varkappa\bigr)\le |K|^n-\n(\mathcal U)$ and $b(K)$ even; we have $\ell_{n,K}\bigl(b(K)\bigr)\le\pi^{\E}_{n,\n(\mathcal U)}(K)$. As $|K|^n\ge\max\bigl(9,2\n(\mathcal U)\bigr)$ by hypotheses, we have $|K|^n-\n(\mathcal U)\ge\max\bigl(5,\n(\mathcal U)\bigr)$. From this and the relations $\ell_{n,K}(\varkappa)=\ell_{n,K}\bigl(a(K)\varkappa\bigr)\le\ell_{n,K}\bigl(b(K)\bigr)\ell_{n,K}\bigl(b(K)a(K)\varkappa\bigr)$ by Axioms $A1_{n,K}$ and $A3_{n,K}$, we get that part (1.c) follows from Equation (\ref{EQ15}) applied to $b(K)a(K)\varkappa$.

The argument that part (1.d) follows from part (1.a) is the same as the one of the previous paragraph applied to $a=1_{\mathbb A^n_K}$.

If $|K|=2$, then from Proposition \ref{PR16}(1) and (4) we get that we can choose $a\in\TGA_n(K)[n-1]$ such that $a(K)$ is a transposition and when $\c_2(\varsigma)=0$ (resp.\ $\c_2(\varsigma)\ge 1$) we have $\n\bigl(a(K)\varkappa\bigr)=\n(\varkappa)-1$ (resp. $\n\bigl(a(K)\varkappa\bigr)=\n(\varkappa)-2$) and hence from part (1.a) applied to $a(K)\varkappa$ we get that part (1.e) holds.

For part (2), we use that for $\sigma\in\mathcal U$ we have an inequality $\ell_{n,K}(\sigma)\le \omega_0(\mathcal U)\eth_0(\mathcal U)^2$. 

We consider a product decomposition $\varkappa=\prod_{i=1}^{\mu_m(\mathcal U)} \sigma_i$ by the definition of $\mu_m(\mathcal U)$, where $\sigma_i\in\mathcal U$ for each $i\in \llbracket1,\mu_m(\mathcal U)\rrbracket$. By Axiom $A3_{n,K}$ we have an inequality $\ell_{n,K}(\varkappa)\le\prod_{i=1}^{\mu_m(\mathcal U)} \ell_{n,K}(\sigma_i)$ from which Inequality (\ref{EQ17}) follows.

Let $b\in\TGA_n(K)[\pi^{\E}_{n,l}(K)]$ be such that $\n\bigl(b(K)\varkappa\bigr)\le |K|^n-l$ and $b(K)$ is even; we have $\ell_{n,K}\bigl(b(K)\bigr)\le\pi^{\E}_{n,l}(K)$. We write $b(K)\varkappa=\prod_{i=1}^{\mu_{|K|^n-l,\mathcal U}} \vartheta_i$ by the definition of $\mu_{|K|^n-l,\mathcal U}$, with $\vartheta_i\in\mathcal U$ for each $i\in \llbracket1,\mu_{|K|^n-l,\mathcal U}\rrbracket$. By Axiom $A3_{n,K}$ we have a similar inequality $\ell_{n,K}(\varkappa)\le \ell_{n,K}\bigl(b(K)\bigr)\prod_{i=1}^{\mu_{|K|^n-l,\mathcal U}} \ell_{n,K}(\vartheta_i)$ from which Inequality (\ref{EQ19}) follows.

Inequalities (\ref{EQ18}) and (\ref{EQ20}) are proved in the same way as Inequalities (\ref{EQ17}) and (\ref{EQ19}) (respectively) using the inequality $\ell^{\S}_{n,K}(\sigma)\le \omega^{\S}_0(\mathcal U)\eth^{\S}_0(\mathcal U)^2$, an automorphism $b\in\STGA_n(K)[\pi^{\S}_{n,l}(K)]$, and Axiom $A3^{\S}_{n,K}$. 
\end{proof}

\begin{example}\normalfont\label{EX17}
Let $k:=|K|-1$. Let $s\in \llbracket2,|K|^n\rrbracket$. If $4\mid |K|$ we assume that $s$ is odd and $s\le |K|^n-2$. We have $\n(\mathcal Y_s)=s$. We consider three cases as follows.

{\bf Case 1: $s>2$.} If $o\in\bigl\llbracket0,\lfloor\frac{|K|^n}{s}\rfloor-1\bigr\rrbracket$, then for $s\ge 4$ we have $\mu_{so}(\mathcal Y_s)\le 2(o+1)$ by Theorem \ref{P11}(1). If $s\in \llbracket3,2p-1\rrbracket$ is odd, then from Proposition \ref{PR10}(1) we get the inequality $\omega^{\S}(\mathcal Y_s)\le E_{n-2,k,\frac{s-3}{2}}$.

{\bf Case 2: $s=3$.} Then $\mu_m(\mathcal Y_3)\le\lfloor\frac{m+3}{2}\rfloor$ for each $m\in \llbracket0,|K|^n-3\rrbracket$ by Lemma \ref{P3}(2). If $|K|=2$, then $\omega(\mathcal Y_3)=\Omega(\mathcal Y_3)\le 2n-3$ by the proof of Proposition \ref{PR16}(5). If $|K|=3$, then $\omega^{\S}(\mathcal Y_3)\le 2n-2$ and $\Omega^{\S}(\mathcal Y_3)\le 4n$ by Proposition \ref{PR16}(7). If $|K|\ge 4$, then $\omega^{\S}(\mathcal Y_3)\le E_{n-2,k}$ and $\Omega^{\S}(\mathcal Y_3)\le 4E_{n-2,k}$ by Proposition \ref{PR10}(3) and Example \ref{EX8}.

\phantomsection{{\bf Case 3: $s=5$ and $|K|^n\ge 15$.} Then $\mu_m(\mathcal Y_5)\le 2\lfloor\frac{m+5}{5}\rfloor$ by Theorem \ref{T2}(7) for each $m\in\{0,1,2,5,6,7,8\}\cup\llbracket10,|K|^n-5\rrbracket$. If $\char(K)>2$, then we have
$$\omega^{\S}(\mathcal Y_5)\le E_{n-2,k,1}=(k-1)^2k^2+k(k^2-k+1)(n-2)$$ 
by Proposition \ref{PR10}(1) applied to $5$-cycles. If $|K|=2$ (resp.\ $4\mid |K|$), then we have $\Omega(\mathcal Y_5)\le (2n-3)^2$ (resp.\ $\Omega^{\S}(\mathcal Y_5)\le (4E_{n-2,k})^2$) based on Case 2 and the fact that each $5$-cycle is a product of two $3$-cycles. Hence from Inequalities (\ref{EQ15}) and (\ref{EQ16}) applied to $\mathcal C=\mathcal Y_5$ and the relations $\eth(\mathcal Y_5)\le\pi_{n,5}(K)$ by Lemma \ref{L8}(1), we get that for each $\varkappa\in\Alt(K^n)$ we have 
$\ell_{n,K}(\varkappa)\le [E\pi_{n,5}(K)^2]^{\nu_{5,1}(\varkappa)}$,
with $E$ equal to $E_{n-2,k,1}$ if $\char(K)>2$, to $(2n-3)^2=E_{n-2,1}$ if $|K|=2$, and to $(4E_{n-2,k})^2$ if $4\mid |K|$. Here we used that $\nu_{\mathcal Y_5}(\varkappa)=\nu_{5,1}(\varkappa).$ 
We similarly argue that $\ell^{\S}_{n,K}(\varkappa)\le [E\pi^{\S}_{n,5}(K)^2]^{\nu_{5,1}(\varkappa)}$. Recall that $\nu_{5,1}(\varkappa)$ is bounded in Theorem \ref{P4}(3) and Theorem \ref{T2}(7) and (9).}\label{PH93}
\end{example}

\section{On Nagata's automorphism}\label{S14}

For a field $K$, let $\mathcal N=\mathcal N_K\in\GA_3(K)$ be Nagata's automorphism\index{Nagata's automorphism} of \cite{N}, Sect.\ 2.1 given by 
$$\mathcal N:=\e(x_1-2x_1x_2x_3-2x_2^3-x_1^2x_3^3-2x_1x_2^2x_3^2-x_2^4x_3,x_2+x_1x_3^2+x_2^2x_3,x_3).$$ 
We have $\mathcal N^{-1}=\e(x_1+2x_1x_2x_3+2x_2^3-x_1^2x_3^3-2x_1x_2^2x_3^2-x_2^4x_3,x_2-x_1x_3^2-x_2^2x_3,x_3)$. If $\char(K)=2$, then we have $\mathcal N=\e(x_1+x_1^2x_3^3+x_2^4x_3,x_2+x_1x_3^2+x_2^2x_3,x_3)$ and $\mathcal N^{-1}=\e(x_1+x_1^2x_2^3+x_2^4x_3,x_2+x_1x_3^2+x_2^2x_3,x_3)$.

\phantomsection{Let $\Delta:=x_1x_3+x_2^2\in K[x_1,x_2,x_3]$. The $K$-algebra automorphism of $K[x_1,x_2,x_2]$ that defines $\mathcal N$ fixes both $x_3$ and $\Delta$. Based on this, it is easy to see that for $l\in\mathbb Z$ we have the following simpler iteration identity}\label{EXT4}
$$\mathcal N^l=\e(x_1-2l\Delta x_2-l^2\Delta^2x_3,x_2+l\Delta x_3,x_3).$$

Clearly, for each non-zero integer $l$ not divisible by $\char(K)$ we have $\ell(\mathcal N^l)=5$. 

Also, $\mathcal N$ has infinite order iff $\char(K)=0$. If $\char(K)$ is a prime $p$, then $\mathcal N$ has order $p$ and therefore $\mathcal N(K)\in\perm(K^3)$ is a product of disjoint $p$-cycles. 

If $\char(K)=2$, then
$$\{P\in K^3|\mathcal N(P)=P\}=\bigl(K^2\times\{0\}\bigr)\cup\Bigl\{\Bigl(\frac{\beta^2}{\gamma},\beta,\gamma\Bigr)\Bigl|(\beta,\gamma)\in K\times K^{\ast}\Bigr\}$$
is the zero locus $x_3(x_1x_3+x_2^2)=0$ in $K^3$. If $\char(K)\neq 2$, then 
\begin{equation}\label{EQ20.1}
\{P\in K^3|\mathcal N(P)=P\}=\bigl(K\times\{(0,0)\}\bigr)\cup\Bigl\{\Bigl(\frac{-\beta^2}{\gamma},\beta,\gamma\Bigr)\Bigl|(\beta,\gamma)\in K\times K^{\ast}\Bigr\}.
\end{equation}
is the zero locus $x_1x_3+x_2^2=0$ in $K^3$. 

If $\char(K)=0$, then both $x_1-2x_1x_2x_3-2x_2^3-x_1^2x_3^3-2x_1x_2^2x_2^2-x_2^2x_3$ and $x_2+x_1x_3^2+x_2^2x_3$ are wild coordinates of $R=K[x_1,x_2,x_3]$, i.e., for each tame automorphism $\e(f_1,f_2,f_3)\in\TGA_3(K)$ we have 
$$\{x_1-2x_1x_2x_3-2x_2^3-x_1^2x_3^3-2x_1x_2^2x_2^2-x_2^2x_3,x_2+x_1x_3^2+x_2^2x_3\}\cap\{f_1,f_2,f_3\}=\emptyset$$ by \cite{UY}, Thm.\ 4 or its generalization \cite{UY}, Thm.\ 5. Thus $\mathcal N\in\GA_3(K)\setminus\TGA_3(K)$ if $\char(K)=0$, which was first proved in \cite{SU}, Cor.\ 9. 

If $\char(K)>0$, then Nagata's conjecture that $\mathcal N\in\GA_3(K)\setminus\TGA_3(K)$ (see \cite{N}, Conj.\ 3.1) is still open.

From now on we assume that $K$ is a finite field. Let $p:=\char(K)$.

If $p\neq 2$, then $\n\bigl(\mathcal N(K)\bigr)=|K|^3-|K|^2$ and hence for each $l\in\mathbb N^{\ast}\setminus p\mathbb N^{\ast}$ we have $\mathcal N^l(K)\in\frac{|K|^3-|K|^2}{p}\mathcal Y_p$. If $p=2$, then $\n\bigl(\mathcal N(K)\bigr)=|K|^3-2|K|^2+|K|$ and hence $\mathcal N(K)=\mathcal N^l(K)\in\frac{|K|^3-2|K|^2+|K|}{2}\mathcal Y_2$ for each odd $l\in\mathbb N^{\ast}$.

For $\gamma\in K$, let 
$$\mathcal N_{\gamma}:=\e(x_1-2x_1x_2\gamma-2x_2^3-x_1^2\gamma^3-2x_1x_2^2\gamma^2-x_2^4\gamma,x_2+x_1\gamma^2+x_2^2\gamma)\in\GA_2(K).$$

We have the following consequence of various parts of Sections \ref{S9}, \ref{S11}, and \ref{S12}.

\begin{corollary}\label{F16}
Let $\gamma\in K^{\ast}$ and $l\in\llbracket1,p-1\rrbracket$. If $|K|\ge 3$, then $\ell_{2,K}(\mathcal N^l_{\gamma})=4$ and $\j_{\mathcal N^l_{\gamma}}=2$.
\end{corollary}

\begin{proof} As $\mathcal N_{\gamma}^{t}:=\e(x_1-2tx_1x_2\gamma-2tx_2^3-t^2x_1^2\gamma^3-2t^2x_1x_2^2\gamma^2-t^2x_2^4\gamma,x_2+tx_1\gamma^2+tx_2^2\gamma)$ for each $t\in\mathbb Z$, we have $\ell_{2,K}(\mathcal N^l_{\gamma})\le 4$. For $|K|\ge 5$ we have $\ell_{2,K}(\mathcal N^l_{\gamma})=4$ by Proposition \ref{PR13}(3).

We have $\n(\mathcal N_{\gamma})=|K|^2-|K|$ and $\mathcal N_{\gamma}(K)$ is a product of $\frac{|K|^2-|K|}{p}$ disjoint $p$-cycles. From this and the fact that $K^2\setminus\supp(\mathcal N_{\gamma})$ is non-collinear by Equation (\ref{EQ20.1}), we get that $\j_{\mathcal N_{\gamma}}\ge 1$ by Lemma \ref{F6}(2). 

For $P=(\alpha,\beta)\in K^2$ we have $P\in\supp\bigl(\mathcal N_{\gamma}(K)\bigr)$ iff $\alpha\beta+\gamma^2\neq 0$. If $p\ge 3$, then for $\alpha\beta+\gamma^2\neq 0$ the vectors 
$\mathcal N_{\gamma}(P)-P=\bigl(-2\beta(\alpha\gamma+\beta^2)-\gamma(\alpha\gamma+\beta^2)^2,\gamma(\alpha\gamma+\beta^2)\bigr)$
and $\mathcal N^2_{\gamma}(P)-P=\bigl(-4\beta(\alpha\gamma+\beta^2)-4\gamma(\alpha\gamma+\beta^2)^2,2\gamma(\alpha\gamma+\beta^2)\bigr)$ are not proportional; thus each $p$-cycle of $\mathcal N_{\gamma}$ is non-collinear. If $p=2$ and $4\mid |K|$, then for $\alpha\beta+\gamma^2\neq 0$ the line passing through $P$ and $\mathcal N_{\gamma}(P)$ has slope $\frac{1}{\alpha\gamma+\beta^2}$ which clearly depends on $P$. Thus for $|K|\ge 3$ we have $\j_{\mathcal N_{\gamma^l}}\ge 2$ by Lemma \ref{L7}(1.c) and (1.d). 

As $\ell_{2,K}(\mathcal N^l_{\gamma})\le 4$, $\ell_{2,K}\bigl(\mathcal N^l_{\gamma}(K)\bigr)$ is a product of at most $2$ primes. From this and Proposition \ref{PR14}(2.b) we get that $\ell_{2,K}\bigl(\mathcal N^l_{\gamma}(K)\bigr)=4$ and $\j_{\mathcal N^l_{\gamma}}=2$.\footnote{The inequality $\j_{\mathcal N_{\gamma^l}}\le 2$ also follows directly from the proof of \cite{N}, Part 2, Thm.\ 1.4.}\end{proof} 

\begin{remark}\normalfont\label{R12}
Lemma \ref{L7}(1.d) does not hold for $\j_{\sigma}=2$ as one can see by considering $\sigma:=\mathcal N_{\gamma}(K)$ with $\gamma\in K^{\ast}$. 
\end{remark}

We use the $\mathcal N_{\gamma}$s to exemplify that Furter's lengths can jump in families.

\begin{example}\normalfont\label{EX19}
If $p=2$, then $\mathcal N_0=1_{\mathbb A^2_K}$ and thus $\j_{\mathcal N_0}=0$. If $p>2$, then $\mathcal N_0=\e(x_1-2x_2^3,x_2)$ and thus $\j_{\mathcal N_0}=1$ and $\mathcal N_0(K)$ is a product of $\frac{|K|^2-|K|}{p}$ disjoint $p$-cycles. If $\gamma\in K^{\ast}$, then $\j_{\mathcal N_{\gamma}}=2$ by Corollary \ref{F16}.
\end{example}

We use a variation of $\mathcal N$ and the $\mathcal N_{\gamma}$s to exemplify situations in which one can get smaller degree lengths than what one gets via the usual usage of selective shifts.

\begin{example}\normalfont\label{EX19.1}
For $t\in\mathbb N^{\ast}\setminus\{1\}$ with $t<|K|$ and $(\alpha_1,\ldots,\alpha_{t-1})\in\mathbb D_{1,t-1}(K^{\ast})$, let $f(x):=x\prod_{i=1}^{t-1} (x-\alpha_i)\in K[x]$. We consider the automorphism 
$$\mathcal N_f=\mathcal N_{f,K}:=\e\bigl(x_1-2f(\Delta)x_2-f(\Delta)^2x_3,x_2+f(\Delta)x_3,x_3\bigr)\in\GA_3(K).$$
As above we argue that for each $l\in\mathbb Z$ not divisible by $\char(K)$ we have an identity $\ell(\mathcal N_f)=4t+1$. If $\char(K)$ is a prime $p$, then $\mathcal N_f^p=1_{\mathbb A^3_K}$ and $\supp(\mathcal N_f)$ is the complement in $K^3$ of the zero locus $\Delta\prod_{i=1}^{t-1} (\Delta-\alpha_i)=0$ in $K^3$ if $\char (K)\neq 2$ and of the zero locus $x_3\Delta\prod_{i=1}^{t-1} (\Delta-\alpha_i)=0$ in $K^3$ if $\char (K)=2$. Thus, for $\gamma\in K^{\ast}$, 
$$\mathcal N_{f,\gamma}:=\e\bigl(x_1-2f(x_1\gamma+x_2^2)x_2-f(x_1\gamma+x_2^2)^2\gamma,x_2+f(x_1\gamma+x_2^2)\gamma\bigr)\in\GA_2(K)$$
is such that $\mathcal N_{f,\gamma}(K)$ is a product of $\frac{|K|^2-t|K|}{p}$ disjoint $p$-cycles and $\ell(\mathcal N_f)=4t$, which for $4t+1<|K|$ is a smaller degree length than what one usually gets via selective shifts for the same number of disjoint $p$-cycles.
\end{example}

\begin{remark}\normalfont\label{R12.1}
Let $(f,g,h,l)\in K[x]^3\times\mathbb N^{\ast}$ be such that $l\ge\deg(g)+1$ and $fgh\neq 0$. We check that
$$\e\left(x_1+\frac{f(x_2)}{x_3}-\frac{f\bigl(x_2+g(x_1+\frac{f(x_2)}{x_3})x_3^l h(x_3)\bigr)}{x_3},x_2+g\Bigl(x_1+\frac{f(x_2)}{x_3}\Bigr)x_3^l h(x_3),x_3\right)$$
is an automorphism $\mathcal W_{f,g,h,l}\in\GA_3(K)$. The second coordinate is a polynomial in $K[x_1,x_2,x_3]$ as $l\geq \deg g.$ The first coordinate is a polynomial $K[x_1,x_2,x_3]$ as the difference $f(x_2)-f\bigl(x_2+g(x_1+\frac{f(x_2)}{x_3})x_3^l h(x_3)\bigr)$ is divisible by $g\bigl(x_1+\frac{f(x_2)}{x_3})x_3^l h(x_3)$ and hence, as $l\geq \deg(g) +1$, by $z$. As endomorphisms of $K[x_1,x_2,x_3]\bigl[\frac{1}{x_3}\bigr]$, we have an identity $\mathcal W_{f,g,h,l}=\mathcal W_1^{-1}\mathcal W_2\mathcal W_1$, where $\mathcal W_1:=\e\bigl(1+\frac{f(x_2)}{x_3},x_2,x_3\bigr)$, $\mathcal W_2:=\e\bigl(x_1,x_2+g(x_1)x_3^lh(x_3),x_3\bigr)$, and $\mathcal W_1^{-1}=\e\bigl(1-\frac{f(x_2)}{x_3},x_2,x_3\bigr)$ have Jacobian determinants equal to $1$. Thus $\mathcal W_{f,g,h,l}\in\End_3(K)$ has Jacobian determinant $1$ and hence $\mathcal W_{f,g,h,l}\in\GA_3(K)$. We have $\mathcal W_{f,g,h,l}^{-1}=\mathcal W_1\mathcal W_2^{-1}\mathcal W_1^{-1}=\mathcal W_{-f,-g,h,l}$.

We have $\mathcal N=\mathcal N_{x}=\mathcal W_{x^2,x,1,2}$ and $\mathcal N_{x^l}=\mathcal W_{x^2,x^l,1,l+1}$ for each $l\in\mathbb N^{\ast}\setminus\{1\}$.
\end{remark}

\section{Projective and tame type of invariants of finite subsets}\label{S15} 

We first study linear projections and their relations to directions parametrized by projective spaces and then we introduce the tame type and the tame number of finite subsets of $K^n$. 

\begin{notation}\normalfont\label{N5}
Let $(n,m)\in (\mathbb N^\ast\setminus\{1\})^2$ and $K$ a field. For a non-zero vector $P\in K^n$, let $(P)\in\mathbb P^{n-1}_K(K)$ be the point defined by $P$. For $\lambda\in\mathbb P_K^{n-1}(K)$, let $V_{\lambda}$ be the $1$-dimensional subspace of $K^n$ that defines it and let $\pi_{\lambda}:K^n\rightarrow K^n/V_{\lambda}$ be the $K$-linear projection. Let $J_m:=\{(i,j)\in \llbracket1,m\rrbracket^2|i< j\}$. 
\end{notation}

\begin{definition}\label{D19}
Using Notation \ref{N5}, we consider a subset $Y=\{P_1,\ldots,P_m\}$ of $K^n$ with $|Y|=m$.

\medskip
{\bf (1)} By the directional function of $Y$\index{direction!directional function} mean the function
$$\psi_Y:J_m\to\mathbb P_K^{n-1}(K)$$ 
defined by the rule $(i,j)\mapsto (P_i-P_j)$.

\smallskip
{\bf (2)} By the minimal directional number of $Y$\index{direction!minimal directional number} we mean 
$$\dir_Y:=\min\bigl(|\psi_Y^{-1}(\lambda)|\lambda\in\mathbb P_K^{n-1}(K)\bigr)\in\mathbb N.$$

{\bf (3)} By the weak affine dimension of $Y$\index{affine dimension!weak affine dimension} we mean the smallest $\w_Y\in\llbracket1,n\rrbracket$ such that there exists a linear projection $\pi:\mathbb A^n_K\rightarrow\mathbb A_K^{\w_Y}$ with $|\pi(Y)|=m$.
\end{definition}

\begin{lemma}\label{L10}
We assume that the field $K$ is finite and we use Notation \ref{N5} and Definition \ref{D19}. Let $N:=\Bigl\lfloor\frac{m(m-1)}{2(\sum_{i=0}^{n-1} |K|^i)}\Bigr\rfloor\in\mathbb N$. Let $l\in\mathbb N^{\ast}$ be the smallest such that $\sum_{i=0}^l |K|^i>\frac{m(m-1)}{2}$. Then the following properties hold.

\medskip
{\bf (1)} For each $\lambda\in\mathbb P_K^{n-1}(K)$ we have an identity 
$$\psi_Y^{-1}(\lambda)=\sqcup_{O\in K^n/V_{\lambda}} \{(i,j)\in J_m|\pi_{\lambda}(P_i)=\pi_{\lambda}(P_j)=O\}$$ and hence 
$$|\psi_Y^{-1}(\lambda)|=\sum_{O\in K^n/V_{\lambda}} \frac{|\pi_{\lambda}^{-1}(O)\cap Y|\bigl(|\pi_{\lambda}^{-1}(O)\cap Y|-1\bigr)}{2}.$$

{\bf (2)} We have $\dir_Y\le N$. 

\smallskip
{\bf (3)} If $l\le n-1$, then $\psi_Y$ is not surjective (i.e., $\dir_Y=0$).

\smallskip
{\bf (4)} If $\lambda\in\mathbb P_K^{n-1}(K)\setminus\Im(\psi_Y)$ then we have $|\pi_{\lambda}(Y)|=m$ and if $\lambda\in\Im(\psi_Y)$ then we have $|\pi_{\lambda}(Y)|\le m-1$. In particular, we have $\w_Y=n$ iff $\psi_Y$ is surjective.

\smallskip
{\bf (5)} If $l\le n-1$, then $\w_Y\le l$, i.e., there exits an $n-l$ dimensional subspace $V$ of $K^n$ such that for the $K$-linear map $\pi:K^n\rightarrow K^n/V$ we have $|\pi(Y)|=m$.

\smallskip
{\bf (6)} Suppose that $\psi_Y$ is surjective (so $l=n$) and $m(m-1)\le 4\sum_{i=1}^{n-1} |K|^i$. Then $\dir_Y=1$ and there exist at least two points $\lambda\in\mathbb P_K^{n-1}(K)$ such that $|\psi_Y^{-1}(\lambda)|=1$.

\smallskip
{\bf (7)} We have $\bigl\lfloor\frac{|K|(|K|-1)}{2(|K|+1)}\bigr\rfloor=\bigl\lfloor\frac{|K|-2}{2}\bigr\rfloor$. In particular, if $n=2$ and $m=|K|$, then $N= \bigl\lfloor\frac{|K|-2}{2}\bigr\rfloor$.

\smallskip
{\bf (8)} If $n=2$ and $m=|K|$, then $\dir_Y\in\bigl\llbracket0,\bigl\lfloor\frac{|K|-2}{2}\bigr\rfloor\bigr\rrbracket$.
\end{lemma}

\begin{proof}
Parts (1) and (4) are clear. 

Part (2) follows from the identities $|\mathbb P_K^{n-1}(K)|=\sum_{i=0}^{n-1} |K|^i$ and $|J_m|=\frac{m(m-1)}{2}$ that involve the cardinalities of target and the source of $\psi_Y$. 

Part (3) follows from part (2) and the identity $N=0$.

Part (5) follows from an iterated application of parts (3) and (4). 

For part (6), if there exists at most one $\lambda\in\mathbb P_K^{n-1}(K)$ such that $|\psi_Y^{-1}(\lambda)|=1$, then $\frac{m(m-1)}{2}=|J_m|\ge 1+2(|\mathbb P_K^{n-1}(K)|-1)=1+2\sum_{i=1}^{n-1} |K|^i$, a contradiction to $\frac{m(m-1)}{2}\le 2\sum_{i=1}^{n-1} |K|^i$. So part (6) holds.

Part (7) holds as we have strict inequalities $\frac{|K|-2}{2}<\frac{|K|(|K|-1)}{2(|K|+1)}<\frac{|K|-1}{2}$.

Part (8) holds as $N=\lfloor\frac{|K|-2}{2}\rfloor$ by part (7) and $\dir_Y\le N$ by part (2).\end{proof}

\begin{definition}\label{D20}
Let $(n,m)\in (\mathbb N^{\ast}\setminus\{1\})^2$. Let $K$ be a field with $|K|\ge\sqrt[n]{m}$. Let $(d,s)\in \llbracket1,\min(n,m-1)\rrbracket\times\mathbb N^{\ast}\subset (\mathbb N^{\ast})^2$. Let $Y\subset K^n$ be a subset with $|Y|=m$. 

\medskip
{\bf (1)} We say that $Y$ has tame type\index{tame type} $(d,s)$ if there exists $a\in\STGA_n(K)[s]$ (equivalently $a\in\TGA_n(K)[s]$) and a linear projection $\pi:\mathbb A^n_K\rightarrow\mathbb A^d_K$ such that we have $\bigl|\pi\bigl(a(Y)\bigr)\bigr|=m$ and $\d_{\pi\bigl(a(Y)\bigr)}=d$.

\smallskip
{\bf (2)} Suppose that $n=2$. For $j\in\mathbb N$, we say that $Y$ has extended tame type\index{tame type!extended tame type} $(1,s,j)$ if there exists $a\in\SGA_2(K)[s]$ (equivalently $a\in\GA_2(K)[s]$) and a linear projection $\pi:\mathbb A^2_K\rightarrow\mathbb A^1_K$ such that $\j_a\le j$ and $\bigl|\pi\bigl(a(Y)\bigr)\bigl|=m$.

\smallskip
{\bf (3)} If $K$ is finite, then by the virtual affine dimension of $Y$ or of the pair $(m,K)$\index{affine dimension!virtual affine dimension} we mean the smallest $\v_Y=\v_{m,|K|}\in \llbracket1,\d_Y\rrbracket$ such that $|Y|\le |K|^{\v_Y}$. If $K$ is infinite, by the virtual dimension of $Y$ or of the pair $(m,K)$ we mean $\v_Y=\v_{m,|K|}:=1$.

\smallskip
{\bf (4)} By the tame number\index{tame number} of $Y$ we mean the smallest $\t_Y\in\mathbb N^{\ast}$ such that $Y$ has tame type $(\v_Y,\t_Y)$.
\end{definition}

Clearly, $Y$ has tame type $(\w_Y,1)$ but, if $\w_Y\ge 2$, it does not have tame type $(\w_Y-1,1)$. Moreover, 
$$\v_Y\le\w_Y\le\d_Y.$$ 
Also, we have $\t_Y=1$ iff $\v_Y=\w_Y$.

The next proposition provides the first computations and estimates of the tame numbers in terms of collinear and minimal directional numbers.

\begin{proposition}\label{PR16.5}
Let $n\in\mathbb N^{\ast}\setminus\{1\}$. Let $K$ be a finite field with $|K|\ge 3$. Let $Y$ be a finite subset of $K^n$ with $3\le|Y|\le |K|$. Then the following properties hold.

\medskip
{\bf (1)} If $\c_Y\ge |Y|-1$, then $\t_Y=1$.

\smallskip
{\bf (2)} Suppose that $n=2$ and $\dir_Y<\frac{|K|+4}{4}$ (equivalently, and $\dir_Y<\bigl\lceil\frac{|K|+4}{4}\bigr\rceil$). Then $\t_Y\le |Y|-2$ and $Y$ has extended tame type $(1,|Y|-2,1)$.

\smallskip
{\bf (3)} Suppose that $n=2$ and $\dir_Y\in\bigl\llbracket\bigl\lceil\frac{|K|+4}{4}\bigr\rceil,\bigl\lfloor\frac{|K|-2}{2}\bigr\rfloor\bigr\rrbracket$. Let 
$$r_{0;Y}\in\bigl\{\lceil\sqrt{2\dir_Y}\rceil,\lceil\sqrt{2\dir_Y}\rceil+1\bigr\}$$ 
be the smallest such that $2\dir_Y\le r_{0;Y}^2-r_{0;Y}$. Then $\t_Y\le \max(1,|Y|-r_{0;Y})(|Y|-2)$ and $Y$ has extended tame type $\bigl(1,\max(1,|Y|-r_{0;Y})(|Y|-2),2\bigr)$.

\smallskip
{\bf (4)} If $n=2$ and $Y$ is contained in the union of two distinct lines in $K^n$, then $\t_Y\le |Y|-2$ and $Y$ has extended tame type $(1,|Y|-2,1)$.

\smallskip
{\bf (5)} If $|Y|=3$, then $\t_Y=1$.

\smallskip
{\bf (6)} If $n=2$ and $|Y|=4$, then $\t_Y\in\{1,2\}$ if $|K|\in\{4,5\}$ and $\t_Y=1$ if $|K|\ge 7$.
\end{proposition}

\begin{proof} 
Let $m:=|Y|$. Based on Lemma \ref{L10}(5) we can assume that $n=2$ even for parts (1) and (5). We write $Y=\{P_1,\ldots,P_m\}$. The case $\c_Y=m$ is clear so we can assume that $\c_Y\le m-1$.

For part (1), we have $\c_Y=m-1$. Up to affine automorphisms we can assume that there exists a subset $W=\{\alpha_1,\ldots,\alpha_{m-1}\}\subset K$ with $|W|=m-1$ and $(\alpha_m,\beta)\in K\times K^{\ast}$ such that $Y=(W\times\{0\})\cup\{(\alpha_m,\beta)\}$. 

We can assume that $\alpha_m\in W$. As $(W\setminus\{\alpha_m\}\cap K^{\ast})\le m-2\le |K|-2$, there exists $\gamma\in K^{\ast}$ such that $\alpha_m+\beta\gamma\notin W\setminus\{\alpha_m\}$. So for $a:=\e(x_1+\gamma x_2,x_2)\in\SL_2(K)$, the first coordinates of the $m$ points in $a(Y)$ are distinct; so $\t_Y=1$ and part (1) holds.

For parts (2) and (3), let $(1:0)\in\mathbb P^1_K(K)$ and $N:=|\psi_Y^{-1}(1:0)|$. 

By writing $P_{i}=(\alpha_i,\beta_i)\in K^2$ for $i\in \llbracket1,m\rrbracket$, let $Z:=\{\beta_i|i\in \llbracket1,m\rrbracket\}\subset K$. Let $s:=|Z|\in \llbracket1,m\rrbracket$. We list $Z=\{\gamma_1,\ldots,\gamma_s\}$ such that by defining integers $q_i:=|\{j\in \llbracket1,m\rrbracket|\beta_j=\gamma_i\}|$ for $i\in \llbracket1,s\rrbracket$, we have $q_1\ge q_2\ge\cdots\ge q_s\ge 1$. We also have $\sum_{i=1}^s q_i=m$. As 
\begin{equation}\label{EQ22}
m-s=\sum_{i=1}^s (q_i-1)\le\sum_{i=1}^s \binom{q_i}{2}=N,
\end{equation} 
we get that $s\ge m-N$.

To prove parts (2) and (3), up to affine automorphisms we can assume that $N=\dir_Y$. Hence $s\ge 2$.

For part (2), we have $4N-4<|K|$ and we show that there exists a function $h:Z\rightarrow K$ such that we have an injective function $\llbracket1,m\rrbracket\rightarrow K$ defined by the rule $i\mapsto\alpha_i+h(\beta_i)$. Taking $h(\gamma_1)\in K$ arbitrarily, it is possible to define $h(\gamma_j)$ inductively for $j\in \llbracket2,s\rrbracket$ provided for every $j\in \llbracket2,s\rrbracket$, by defining $l_j:=q_j(\sum_{i=1}^{j-1} q_i)$ we have an inequality $l_j<|K|$. If $q_j=1$, then $l_j\le m-1<|K|$. If $q_j\ge 2$, then 
\begin{equation}\label{EQ22.9}
l_j\le \sum_{i=1}^{j-1} q_i^2\le \sum_{i=1}^{j-1} 4\binom{q_i}{2}\le 4\Bigl[N-\binom{q_j}{2}\Bigr]\le 4N-4.
\end{equation} 
As $4N-4<|K|$, we get that $l_j<|K|$. So $h$ exists.\footnote{If $4\mid |K|$ and we have $4N-4=|K|$, then the resulting inequality $l_j\le |K|$ is strict in all cases except when $q_1=\cdots=q_j=2$, $q_{j+1}=\cdots=q_s=1$, $j=N$, and $m=j+s$.} 

If $s\le m-1$, let $f\in K[x]$ be the Lagrange interpolating polynomial of degree at most $s-1\le m-2$ such that for every $l\in \llbracket1,s\rrbracket$ we have $f(\gamma_l)=h(\gamma_l)$; therefore $a:=\e(x_1+f(x_2),x_2)\in\SGA_2(K)[s-1]\subset \SGA_2(K)[m-2]$ and the first coordinates of $a(P_1),\ldots,a(P_m)$ are distinct. If $s=m$, then the same holds for $a:=\e(x_2,-x_1,x_3,\ldots,x_n)\in\SL_2(K)$. Thus $\t_{Y}\le m-2$ and $Y$ has extended tame type $(1,m-2,1)$. So part (2) holds.

For part (3), we have $4N-4\ge |K|$ and $N\le \lfloor\frac{|K|-2}{2}\rfloor$; hence $2N\le |K|-2$. 

We show that there exists a function $h:Z\rightarrow K$ such that the fibers of the function $\hbar:\llbracket1,m\rrbracket\rightarrow K$ defined by the rule $i\mapsto\alpha_i+h(\beta_i)$ have cardinality at most $2$ and moreover the number of fibers of cardinality $2$ is at most $\lfloor\frac{N}{2}\rfloor$. For $j\in\llbracket2,s\rrbracket$, let $m_j:=\sum_{i=1}^{j-1} q_j$ and let $\hbar_j:=\hbar|\llbracket1,m_j\rrbracket$ be the restriction of $\hbar$ to $\llbracket1,m_j\rrbracket$. For $\iota\in\{1,2\}$, let $l_{j,\iota}$ be the number of fibers of $\hbar_j$ of cardinality $\iota$; we have $l_{j,1}+2l_{j,2}=m_j=\sum_{i=1}^{j-1} q_i$. Taking $h(\gamma_1)\in K$ arbitrarily, it is possible to define $h(\gamma_j)$ inductively for $j\in \llbracket2,s\rrbracket$ provided for every $j\in \llbracket2,s\rrbracket$ with $q_j\notin\{1,3\}$ we have $l_{j,2}\lceil\frac{q_j}{2}\rceil+(l_{j,1}+l_{j,2})\lfloor\frac{q_j}{2}\rfloor<|K|$, for $q_j=3$ we have $l_{j,2}+2(l_{j,1}+l_{j,2})<|K|$, and for $q_j=1$ we have $(l_{j,1}+l_{j,2})q_j<|K|$. For  instance, if $j\in \llbracket2,s\rrbracket$ with $q_j\notin\{1,3\}$, this is so as the contribution of $q_j$ to $\lfloor\frac{N}{2}\rfloor$ is at least $\lfloor\frac{q_j(q_j-1)}{2}\rfloor$ which is greater than or equal to $\lceil\frac{q_j}{2}\rceil$.

If $q_j=1$, then $l_{j,1}+l_{2,j}\le m-1<|K|$. 

If $q_j\ge 2$ but $q_j\neq 3$, then from the inequality $\bigl\lceil\frac{q_j}{2}\bigr\rceil+\bigl\lfloor\frac{q_j}{2}\bigr\rfloor\le q_j$ we get that
$$l_{j,2}\bigl\lceil\frac{q_j}{2}\bigr\rceil+(l_{j,1}+l_{j,2})\bigl\lfloor\frac{q_j}{2}\bigr\rfloor\le (l_{j,1}+2l_{j,2})\frac{q_j}{2}=\frac{q_j}{2}\Bigl(\sum_{i=1}^{j-1} q_i\Bigr).$$ 
Based on this and the inequality $q_j(\sum_{i=1}^{j-1} q_i)<4N-4$ (see Equation (\ref{EQ22.9})) we get that 
$l_{j,2}\lceil\frac{q_j}{2}\rceil+(l_{j,1}+l_{j,2})\lfloor\frac{q_j}{2}\rfloor\le 2N-2\le |K|-4<|K|$. 

If $q_j=3$, then $\sum_{i=1}^{j-1} q_i\le \sum_{i=1}^{j-1} \frac{q_i(q_i-1)}{2}\le N$ and thus we estimate
$$l_{j,2}+2(l_{j,1}+l_{j,2})=2l_{j,1}+3l_{j,2}=\frac{3}{2}(l_{j,1}+2l_{j,2})+\frac{1}{2}l_{j,1}\le \Bigl\lfloor\frac{3}{2}(l_{j,1}+2l_{j,2})+\frac{N}{2}\Bigr\rfloor$$ 
$$=\Bigl\lfloor\frac{N+q_jm_j}{2}\Bigr\rfloor=\Bigl\lfloor\frac{N}{2}+\frac{q_j}{2}\Bigl(\sum_{i=1}^{j-1} q_i\Bigr)\Bigr\rfloor\le \Bigl\lfloor\frac{N}{2}+\frac{1}{2}\Bigl(\sum_{i=1}^{j-1} q_i^2\Bigr)\Bigr\rfloor\le\Bigl\lfloor\frac{N}{2}+\frac{1}{2}\Bigl[\sum_{i=1}^{j-1} 3\binom{q_i}{2}\Bigr]\Bigr\rfloor$$
$$\le \Bigl\lfloor\frac{N}{2}+\frac{3}{2}\Bigl[N-\binom{q_j}{2}\Bigr]\Bigr\rfloor=\Bigl\lfloor2N-\frac{9}{2}\Bigr\rfloor=2N-5\le |K|-7<|K|.$$

We conclude that $h$ exists.

As 
$$2\dir_Y=2N=\sum_{i=1}^s q_i(q_i-1)=-m+\sum_{i=1}^s q_i^2$$
$$\le -m+(q_i-s-1)^2+\sum_{i=1}^{s-1} 1^2=-m+s-1+(m-s+1)^2,$$ it follows that $m-s+1\ge r_{0;Y}$ and therefore $m-r_{0;Y}\ge s-1\ge 1$.

Let $f\in K[x]$ of degree at most $s-1\le m-r_{0:Y}$ be such that for each $l\in \llbracket1,s\rrbracket$ we have $f(\gamma_l)=h(\gamma_l)$. Therefore 
$$a:=\e(x_1+f(x_2),x_2)\in\SGA_2(K)[s-1]\subset \SGA_2(K)[m-r_{0:Y}]$$ 
and we have 
$$\bigl|\psi^{-1}_{a(Y)}\bigr|\le\Bigl\lfloor\frac{N}{2}\Bigr\rfloor\le\Bigl\lfloor\frac{\lfloor\frac{|K|-2}{2}\rfloor}{2}\Bigr\rfloor\le\Bigl\lfloor\frac{|K|-2}{4}\Bigr\rfloor<\frac{|K|+4}{4}.$$
Based on this and part (2) we get that $a(Y)$ has extended tame type $\bigl(1,|Y|-2,1\bigr)$. Thus $Y$ has extended tame type $\bigl(1,(|Y|-r_{0;Y})(|Y|-2),2\bigr)$. So part (3) holds.

For part (4), let $\lambda_1$ and $\lambda_2$ be two distinct lines in $K^n$ with $Y\subset\lambda_1\cup\lambda_2$. For $\iota\in\{1,2\}$ let $m_{\iota}:=|Y\cap\lambda_{\iota}|$. We have $m_1+m_2\in\{m,m+1\}$ and the equality $m_1+m_2=m+1$ can happen only when the intersection $\lambda_1\cap\lambda_2$ is a point of $Y$. Based on part (1) we can assume that $\c_Y\le m-2$; thus $\min(m_1,m_2)\ge 2$ and $m\ge 4$. So $|K|\ge 4$. We prove that $\dir_Y<\frac{|K|+4}{4}$ by considering two disjoint cases.

{\bf Case 1: $|K|\ge 7$.} Let $(\lambda_3,\lambda_4)\in\mathbb P^1_K(K)^2$ be such that $\lambda_{\iota}$ has the direction given by $\lambda_{\iota+2}$. If $\lambda_3=\lambda_4$, then $|\psi_Y^{-1}(\lambda_3)|\ge \frac{m_1(m_1-1)}{2}+\frac{m_2(m_2-1)}{2}$ and thus $\dir_Y\le \bigl\lfloor\frac{m(m-1)-m_1(m_1-1)-m_2(m_2-1)}{2|K|}\bigr\rfloor$. If $\lambda_3\neq\lambda_4$, then $|\psi_Y^{-1}(\lambda_{\iota+2})|\ge \frac{m_{\iota}(m_{\iota}-1)}{2}$ and hence we have $\dir_Y\le \bigl\lfloor\frac{m(m-1)-m_1(m_1-1)-m_2(m_2-1)}{2(|K|-1)}\bigr\rfloor$. Thus, as $m\le m_1+m_2$, to prove that $\dir_Y<\frac{|K|+4}{4}$ it suffices to show that $4\frac{m_1m_2}{(|K|-1)}-4<|K|$. As we have $m_1+m_2\le m+1\le |K|+1$, we get that $4m_1m_2\le (|K|+1)^2$. So it suffices to show that $|(K|+1)^2<(|K|-1)(|K|+4)$ which holds as $|K|\ge 7$. 

{\bf Case 2: $|K|\in\{4,5\}$.} As $m\le |K|$, we have $\dir_Y\le\bigl\lfloor\frac{|K|(|K|-1)}{2(|K|+1)}\bigr\rfloor=1$ by Lemma \ref{L10}(2) and hence $\dir_Y\le 1<\frac{|K|+4}{4}$.

As $\dir_Y<\frac{|K|+4}{4}$ by the two cases above, part (4) follows from part (2). 

Part (5) follows from part (1).

For part (6), based on part (1) we can assume that $\c_Y=2$. Based on this and Case 2 above we can assume that $|K|\ge 7$. Up to affine automorphisms we can assume that $P_1=(0,0)$, $P_2=(1,0)$, $P_3=(0,1)$, and $P_4=(\alpha,\beta)$ with $\alpha\beta(\alpha+\beta-1)\in K^{\ast}$. If $\alpha=\beta=1$, then for $\gamma\in K\setminus\{-1,0,1\}$ and $a:=\e(x_1+\gamma x_2,x_2)$ the first coordinates of the $4$ points $a(Y)$ are distinct and thus $\t_Y=1$. If $\alpha\neq 1$ or $\beta\neq 1$, then by the symmetry under the involution $\e(x_2,x_1)$ of $\mathbb A^2_K$ we can assume that $\beta\neq 1$ and thus for $\gamma\in K^{\ast}\setminus\bigl\{1,\frac{-\alpha}{\beta},\frac{1-\alpha}{\beta},\frac{\alpha}{1-\beta}\bigr\}$ and $a:=\e(x_1+\gamma x_2,x_2)$ the first coordinates of the $4$ points $a(Y)$ are distinct. So part (6) holds.\footnote{If $|K|=4=|Y|$ and $\c_Y=2=n$, then we have $\t_Y=1$ iff we can write $Y=\{P_1,P_2,P_3,P_4\}$ with $P_1+P_3=P_2+P_4$.}\end{proof}

\begin{remark}\normalfont\label{R12.2}
Suppose that $n=2$ and $Y$ is a finite subset of $K^2$ such that $3\le|Y|\le |K|$. Lemma \ref{L10}(8) and Proposition \ref{PR16.5}(3) can be adapted to the situation when only the fibers of $\psi_Y$ above points of $\mathbb A^1(K)\subset\mathbb P^1_K(K)$ are considered, i.e., when the directional function $\psi_Y$ is replaced by the {\it affine directional function}\index{direction!affine directional function}
$$\psi_Y^0:\psi_Y^{-1}\bigl(\mathbb A^1(K)\bigr)\rightarrow\mathbb A^1_K(K)=\mathbb P^1_K(K)\setminus\{(1:0)\},$$ the role of $\bigl\lfloor\frac{|K|-2}{2}\bigr\rfloor$ being replaced by $\bigl\lfloor\frac{|K|(|K|-1)}{2|K|}\bigr\rfloor=\bigl\lfloor\frac{|K|-1}{2}\bigr\rfloor$. This is so as the proof of Proposition \ref{PR16.5}(3) still works if we only assume that $2N\le |K|-1$.
\end{remark}

We have the following consequence of the inequality of Display (\ref{EQ22}).

\begin{corollary}\label{C16}
Let $K$ be a finite field with at least $3$ elements. Let $m\in \llbracket3,|K|\rrbracket$. Let $Y$ be a non-empty subset of $K^2$ with $|Y|=m$. Let $N_1$ and $N_2$ be the largest two values of the cardinalities of the fibers of the directional function $\psi_Y:J_m\rightarrow\mathbb P_K^1(K)$. Then we have an inequality $\s_Y\ge 2m-2-N_1-N_2$. 
\end{corollary}

\begin{proof}
Up to affine automorphisms we can assume that $\s_Y=\s_{1_{\mathbb A^2_K}}(Y)$. We define $n_1:=\psi_Y^{-1}\bigl((0:1)\bigr)$ and $n_2:=\psi_Y^{-1}\bigl((1:0)\bigr)$. Display (\ref{EQ22}) implies the inequalities $m_{1,Y,1_{\mathbb A^2_K}}\ge m-n_1$ and $m_{2,Y,1_{\mathbb A^2_K}}\ge m-n_2$. Hence the corollary follows from $\s_Y=m_{1,Y,1_{\mathbb A^2_K}}-1+m_{2,Y,1_{\mathbb A^2_K}}-1\ge 2m-2-n_1-n_2\ge 2m-2-N_1-N_2$.
\end{proof}

We have the following general lemma.

\begin{lemma}\label{F11.5}
Let $(n,m)\in (\mathbb N^\ast\setminus\{1\})^2$ and a field $K$ be such that $|K|\ge m$. For $(\underline{P},\underline{Q})=\bigl((P_1,\ldots,P_m),(Q_1,\ldots,Q_m)\bigr)\in\mathbb D^{\d=1}_{n,m}(K)^2$ the following properties hold. 

\medskip
{\bf (1)} Suppose there exists $(\alpha_1,\beta_1,\ldots,\alpha_m,\beta_m)\in K^{2m}$ such that $P_i=(\alpha_i,0,\ldots,0)$ and $Q_i=(0,\beta_i,0,\ldots,0)$ for each $i\in\llbracket1,m\rrbracket$. Let $(f,g)\in K[x]^2$ be the pair of Lagrange interpolating polynomials of degree at most $m-1$ such that $f(\alpha_i)=\beta_i$ and $g(\beta_i)=\alpha_i$  for each $i\in\llbracket1,m\rrbracket$.  If $n=2$, then 
$$\pi^{\S,\le 2}_{\underline{P},\underline{Q}}=\deg(f)\deg(g)\le (m-1)^2.$$ 
In particular, if $n=2$, then either $\pi^{\S,\le 2}_{\underline{P},\underline{Q}}=1$ or $4\le\pi^{\S,\le 2}_{\underline{P},\underline{Q}}\le (m-1)^2$.

\smallskip
{\bf (2)}  We have $\pi^{\S}_{\underline{P},\underline{Q}}\le (m-1)^2$.

\smallskip
{\bf (3)}  If $3\le m=|K|$, then $\pi^{\S}_{\underline{P},\underline{Q}}\le (|K|-2)^2$.

\smallskip
{\bf (4)} If $n=2$ and $3\le m=|K|$, then $\pi^{\S,\le 2}_{\underline{P},\underline{Q}}\le (|K|-2)^2$.
\end{lemma}

\begin{proof}
Up to affine automorphisms in $\ASL_n(K)$ we can assume that the extra hypothesis of part (1) holds for parts (2) to (4) as well. So to prove parts (2) to (4) we can assume that $n=2$  and we use the pair $(f,g)$ of part (1). 

As $m\ge 2$, we have $\deg(f)\le 1$ iff $\deg(f)=1$. Moreover, $\deg(f)=1$ iff $\e\bigl(f(x_1)\bigr)\in\GA_1(K)$ and hence iff $\e\bigl(g(x_1)\bigr)=\e\bigl(f(x_1)\bigr)^{-1}\in\GA_1(K)$. Thus $\deg(f)=1$ iff $\deg(g)=1$. Therefore we have either an identity $\deg(f)\deg(g)=1$ or an inequality $\min\bigl(\deg(f),\deg(g)\bigr)\ge 2$ that implies $\deg(f)\deg(g)\ge 4$.

For $a:=\e\bigl(x_1-g(x_2),x_2\bigr)\e\bigl(x_1,x_2+f(x_1)\bigr)\in\SGA_2(K)$ we have $a(\underline{P})=\underline{Q}$, $\ell^{\S}(a)\le\deg(f)\deg(g)$, and $\j_a\le 2$. Thus $\pi^{\S,\le 2}_{\underline{P},\underline{Q}}\le\deg(f)\deg(g)$. From this and the inequalities $\deg(f)\deg(g)\le (m-1)^2$ and $\pi^{\S}_{\underline{P},\underline{Q}}\le\pi^{\S,\le 2}_{\underline{P},\underline{Q}}$ we get that part (2) holds.

If $3\le m=|K|$, then $f$ and $g$ are permutation polynomials in $K[x]$ and hence they have degrees at most $|K|-2$ by Hermite's criterion of \cite{Di}, Sect.\ 1, Subsect.\ 11; thus parts (3) and (4) hold as we have $\deg(f)\deg(g)\le (|K|-2)^2$.

To show that $\pi^{\S,\le 2}_{\underline{P},\underline{Q}}\ge \deg(f)\deg(g)$ we can assume that $\deg(f)\deg(g)\neq 1$. Thus $\deg(f)\deg(g)\ge 4$ and hence $m\ge 3$. Let $c\in\GA_2(K)$ be such that $\j_c\le 2$, $\ell(c)=\pi^{\S,\le 2}_{\underline{P},\underline{Q}}$, and $c(\underline{P})=\underline{Q}$. If $\j_c=1$, then it is easy to see that $\deg(g)=1$, a contradiction. Thus $\j_c=2$. 

We write $c=c_1a_1c_1^{-1}c_2a_2c_2^{-2}b\in\GA_2(K)$ with $(b,c_1,c_2)\in\AGL_2(K)\times\GL_2(K)^2$, $a_1=\e\bigl(x_1,x_2+f_1(x_1)\bigr)$, and $a_2=\e(x_1+g_1(x_2),x_2)$ where $(f_1,g_1)\in K[x]^2$ satisfies $\bigl(\deg(f_1),\deg(g_1)\bigr)\in \llbracket1,|K|-1\rrbracket^2$ and $\deg(f_1)\deg(g_1)<\deg(f)\deg(g)$. As $\GL_2(K)$ acts $2$-transitively on the set of directions in $K^2$, by replacing $(c,b)$ with $\bigl(c_4cc_4^{-1},c_4bc_4^{-1}\bigr)$ where $c_4:=c_3c_1^{-1}\in\GL_2(K)$ for some $c_3\in\GL_2(K)$ such that $c_4$ fixes $(1,0)$ and $c_3c_1^{-1}c_2$ fixes $e_2$, we can assume that $c_1=c_2=1_{\mathbb A^2_K}$ and $c\bigl(c_4(\underline{P})\bigr)=c_4(\underline{Q})$. So the automorphism $a_3:=\e\bigl(x_1-g_1(x_2+f_1(x_1)),x_2+f_1(x_1)\bigr)$ maps $b\bigl(c_4(\underline{P})\bigr)$ to $c_4(\underline{Q})$. It is easy to see that this implies that $f$ is a linear combination of $1$, $x$, and $f_1(\alpha_1x+\beta_1)$ with coefficients in $K$ for a suitable pair $(\alpha_1,\beta_1)\in K^{\ast}\times K$ and that $g$ is a linear combination of $1$, $x$, and $g_1(\gamma_1x+\delta_1)$ with coefficients in $K$ for a suitable pair $(\gamma_1,\delta_1)\in K^{\ast}\times K$. Hence $\deg(f_1)\ge \deg(f)$ and $\deg(g_1)\ge\deg(g)$. Thus $\deg(f_1)\deg(g_1)\ge \deg(f)\deg(g)$. As $\ell(a)=\deg(f_1)\deg(g_1)$ by Proposition \ref{PR11}(4.a), it follows that $\pi^{\S,\le 2}_{\underline{P},\underline{Q}}\ge\deg(f)\deg(g)$.

Hence $\pi^{\S,\le 2}_{\underline{P},\underline{Q}}=\deg(f)\deg(g)$ and part (1) holds.\end{proof}

We have the following generalization of Lemma \ref{F11.5}(1) and (2). 

\begin{proposition}\label{PR17}
Let $(n,m)\in (\mathbb N^\ast\setminus\{1\})^2$ and $(d_1,d_2,s_1,s_2)\in (\mathbb N^{\ast})^4$. Let $K$ be a field. Let $(\underline{P},\underline{Q})=\bigl((P_1,\ldots,P_m),(Q_1,\ldots,Q_m)\bigr)\in\mathbb D_{n,m}(K)^2$ be such that $Y_1:=\{P_1,\ldots,P_m\}$ has tame type $(d_1,s_1)$ and $Y_2:=\{Q_1,\ldots,Q_m\}$ has tame type $(d_2,s_2)$. Let $d:=\min(d_1,d_2)$. If $n\ge d_1+d_2$, then the following properties hold.\footnote{For instance, if $K$ is finite and $m\le \lfloor\sqrt{|K|^n}\rfloor$, then we can take $d_1=d_2=\v_{m,|K|}$ and $s_i=\t_{Y_i}$ for $i\in\{1,2\}$.}

\medskip
{\bf (1)} We have $\pi^{\S}_{\underline{P},\underline{Q}}\le \prod_{\iota=1}^2 s_{\iota}\min\bigl(\l^{[d_{\iota}]}_K(m-1),\overline{\s}_{|K|}(d,m)\bigr)$. In particular, we have $\pi^{\S}_{\underline{P},\underline{Q}}\le \prod_{\iota=1}^2 s_{\iota}\min\bigl(m-d_{\iota},d(|K|-1)\bigr)$ and, if $n=2$ and $Y_{\iota}$ has extended tame type $(1,s_{\iota},j_{\iota})$ for $\iota\in\{1,2\}$ (so $d=d_1=d_2$ and $m\le |K|$), we have the inequality $$\pi^{\S,\le 2+j_1+j_2}_{\underline{P},\underline{Q}}\le s_1s_2(m-1)^2.$$ 

{\bf (2)} Suppose that $s_1=s_2=1$. Then we have inequalities $d_1\le\d_{Y_1}$, $d_2\le\d_{Y_2}$, and $\pi^{\S}_{\underline{P},\underline{Q}}\le \prod_{\iota=1}^2 \min\bigl(m-d_{\iota},\overline{\s}_{|K|}(d,m)\bigr)$. In particular, we have the following inequality $\pi^{\S}_{\underline{P},\underline{Q}}\le \prod_{\iota=1}^2 \min\bigl(m-d_{\iota},d(|K|-1)\bigr)$ and, if $n=2$ and $d_1=d_2=1$, we have the inequality $\pi^{\S,\le 2}_{\underline{P},\underline{Q}}\le (m-1)^2$.

\smallskip
{\bf (3)} Suppose that $s_1=s_2=1$ and $(d_1,d_2)=(\d_{Y_1},\d_{Y_2})$. Then we have an inequality $\pi^{\S}_{\underline{P},\underline{Q}}\le \prod_{\iota=1}^2 \min(\l_{Y_{\iota}},\overline{\s}_{Y_{\iota}})$ and, if $n=2$ and $d_1=d_2=1$, we also have $\pi^{\S,\le 2}_{\underline{P},\underline{Q}}\le \prod_{\iota=1}^2 \min(\l_{Y_{\iota}},\overline{\s}_{Y_{\iota}})$.
\end{proposition}

\begin{proof}
Up to special linear automorphisms we can assume that for the two projections $\pi_1:\mathbb A^n_K\rightarrow\mathbb A^{d_1}_K$ on the first $d_1$ coordinates and $\pi_2:\mathbb A^n_{K}\rightarrow\mathbb A^{d_2}_{K}$ on coordinates $d_1+1$ to $d_1+d_2$ there exists $(a_1,a_2)\in\STGA_n(K)[s_1]\times \STGA_n(K)[s_2]$ such that for $\iota\in\{1,2\}$ we have $\left|\pi_{\iota}\bigl(a_{\iota}(Y_{\iota})\bigr)\right|=m$ and $\d_{\pi_{\iota}\bigl(a_{\iota}(Y_{\iota})\bigr)}=d_{\iota}$. Based on Theorem \ref{T5}(1) applied to $\pi_1\bigl(a_1(Y_1)\bigr)\subset K^{d_1}$, let $(f_1,\ldots,f_{n-d_1})\in K [x_1,\ldots,x_{d_1}]^{n-d_1}$ be such that for each $j\in \llbracket1,n-d_1\rrbracket$ we have $\deg(f_j)\le\min\bigl(\l_K^{[d_1]}(m-1),\overline{\s}_{\pi_1(a_1(Y_1))}\bigr)$ and for every $i\in \llbracket1,m\rrbracket$ we have $b_1\bigl(a_1(P_i)\bigr)=\bigl(\pi_1(a_1(P_i)),\pi_2(a_2(Q_i)),0,\ldots,0\bigr)$ with $b_1:=\e(x_1,\ldots,x_{d_1},x_{d_1+1}+f_1,\ldots,x_n+f_{n-d_1})\in\STGA_n(K)$. Based on Theorem \ref{T5}(1) applied to $\pi_2\bigl(a_2(Y_2)\bigr)\subset K^{d_2}$, let 
$(g_1,\ldots,g_{n-d_2})\in K [x_{d_1+1},\ldots,x_{d_1+d_2}]^{n-d_2}$ 
be such for each $j\in \llbracket1,n-d_2\rrbracket$ we have $\deg(g_j)\le\min\bigl(\l_K^{[d_2]}(m-1),\overline{\s}_{\pi_2(a_2(Y_2))}\bigr)$ and for every $i\in \llbracket1,m\rrbracket$ we have $b_2\bigl(a_2(Q_i)\bigr)=\bigl(\pi_1(a_1(P_i)),\pi_2(a_2(Q_i)),0,\ldots,0\bigr)$, with
$$b_2:=\e(x_1+g_1,\ldots,x_{d_1}+g_{d_1},x_{d_1+1},\ldots,x_{d_1+d_2},x_{d_1+d_2+1}+g_{d_1+1},\ldots,x_n+g_{n-d_2})$$ 
in $\STGA_n(K)$. 

For $a:=a_2^{-1}b_2^{-1}b_1a_1\in\STGA_n(K)$ we have $a(\underline{P})=\underline{Q}$. If $n=2$, then we have $\j_a\le 2+\j_{a_1}+\j_{a_2}$.

We have $\overline{\s}_{\pi_{\iota}\bigl(a_{\iota}(Y_{\iota})\bigr)}\le\overline{\s}_{|K|}(d_{\iota},m)$ by Theorem \ref{T4}(1) and $\l_K^{[d_{\iota}]}(m-1)\le m-d_{\iota}$ by Lemma \ref{F4}(4). 

We check that we have $\overline{\s}_{|K|}(d_{\iota},m)=\overline{\s}_{|K|}(d,m)\le d(|K|-1)$. If $K$ is infinite, then this is clear as $\overline{\s}_{|K|}(d_{\iota},m)=\overline{\s}_{|K|}(d,m)=m-1$. If $K$ is finite, then $m\le |K|^{d_{\iota}}$ for each $\iota\in\{1,2\}$ and hence $m\le |K|^d$; thus $\overline{\s}_{|K|}(d_{\iota},m)=\overline{\s}_{|K|}(d,m)\le d(|K|-1)$ by Lemma \ref{L6}(1) and (3.a). 

Based on the above inequalities and Inequality (\ref{EQ3}) we conclude that we have an inequality $\ell(a)\le \prod_{\iota=1}^2 s_{\iota}\min\bigl(\l^{[d_{\iota}]}_K(m-1),\overline{\s}_{|K|}(d,m)\bigr)$. So part (1) holds.

For part (2), based on part (1) it suffices to prove the inequality $d_{\iota}\le \d_{Y_{\iota}}$ but this follows from the fact that $\pi_{\iota}$ is a linear projection. 

For part (3), we can take $a_1=a_2=1_{\mathbb A^n_K}$. For $\iota\in\{1,2\}$ we can replace the upper bounds $\l_K^{[d_2]}(m-1)$ and $\overline{\s}_{|K|}(d_{\iota},m)$ of the invariants $\l_{\pi_{\iota}\bigl(a_{\iota}(Y_{\iota})\bigr)}$ and $\overline{\s}_{\pi_{\iota}\bigl(a_{\iota}(Y_{\iota})\bigr)}$ of the unknown set $\pi_{\iota}\bigl(a_{\iota}(Y_{\iota})\bigr)$ with $\l_{Y_{\iota}}$ and $\overline{\s}_{Y_{\iota}}$ by Theorem \ref{T5}(1). So part (3) holds.
\end{proof}

We have the following application of Proposition \ref{PR17}(1) that supplements Lemma \ref{P7} when $n=2$ and the conjugacy classes have at most $|K|$ elements in a way that leads to substantial improvements in the next sections.

\begin{corollary}\label{C15}
Let $\mathcal C$ be a non-trivial conjugacy class of $\Perm(K^2)$ or $\Alt(K^2)$ with $3\le\n(\mathcal C)\le |K|$. For $i\in\{1,2\}$ let $\sigma_i\in\mathcal C$ and let $Y_i:=\supp(\sigma_i)$. Then the following properties hold.

\medskip
{\bf (1)} We have the following inequalities $\ell_{2,K}(\sigma_2)\le \ell_{2,K}(\sigma_1)\t_{Y_1}^2\t_{Y_2}^2(m-1)^4$ 
and $\ell^{\S}_{2,K}(\sigma_2)\le \ell^{\S}_{2,K}(\sigma_1)\t_{Y_1}^2\t_{Y_2}^2(m-1)^4$.

\smallskip
{\bf (2)} Suppose that for $i\in\{1,2\}$, $Y_i$ has extended tame type $(1,\t_{Y_i},\j_i)$. Then $\j_{\sigma_2}\le \j_{\sigma_1}+4+2\j_1+2\j_2$.
\end{corollary}

\begin{proof}
Let $m:=\n(\mathcal C)$. We write $Y_1=\{P_1,\ldots,P_m\}$ and $Y_2=\{Q_1,\ldots,Q_m\}$ in such a way that for each $\theta\in\Perm(K^n)$ with $\theta(P_i)=Q_i$ for every $i\in\llbracket1,m\rrbracket$ we have $\sigma_2=\theta\sigma_1\theta^{-1}$. Let $a\in\SGA_2(K)[\t_{Y_1}\t_{Y_2}(m-1)^2]$ be such that for $\vartheta:=a(K)$ we have $\vartheta(P_i)=Q_i$ for each $i\in\llbracket1,m\rrbracket$ by Proposition \ref{PR17}(1) applied to the tame types $(1,\t_{Y_1})$ and $(1,\t_{Y_2})$. As $\sigma_2=\vartheta\sigma_1\vartheta^{-1}$, from Axioms $A3_{2,K}$ and $A3^{\S}_{2,K}$ we get that $\ell_{2,K}(\sigma_2)\le \ell_{2,K}(\sigma_1)\ell_{2,K}\bigl(a(K)\bigr)^2$ and $\ell^{\S}_{2,K}(\sigma_2)\le \ell^{\S}_{2,K}(\sigma_2)\ell^{\S}_{2,K}\bigl(a(K)\bigr)^2$. From this and the inequalities $\ell_{2,K}\bigl(a(K)\bigr)\le \ell^{\S}_{2,K}\bigl(a(K)\bigr)\le \t_{Y_1}\t_{Y_2}(m-1)^2$ we get that part (1) holds.

Part (2) follows from the fact that we can take $a$ such that $\j_a\le 2+\j_1+\j_2$ by the combination of Proposition \ref{PR17}(1) and (3).
\end{proof}

Before starting with upper bounds we first mention basic identities. Clearly, 
\begin{equation}\label{EQ21}
\pi_{n,p^{nq}-1}(\mathbb F_{p^q})=\pi_{n,p^{nq}}(\mathbb F_{p^q}).
\end{equation}

\begin{lemma}\label{L11}
Suppose that $(n,m)\in (\mathbb N^\ast\setminus\{1\})\times (\mathbb N^\ast\setminus\{1,2\})$. Let $K$ be a field with $|K|\ge\sqrt[n]{m}$. Then the following three statements are equivalent.

\medskip
{\bf (1)} We have $\pi_{n,m}(K)=1$.

\smallskip
{\bf (2)} We have $K\cong\mathbb F_2$ and either $m=3$ or $(n,m)=(2,4)$.

\smallskip
{\bf (3)} We have $\pi^{\S}_{n,m}(K)=1$.
\end{lemma}

\begin{proof}
We prove that $(1)\Rightarrow (2)$. As $n\ge 2$, $K^n$ has three non-collinear points. As each $a\in\AGL_n(K)$ maps non-collinear points to non-collinear points, it follows that $K^n$ does not have three collinear points, hence $K\cong\mathbb F_2$. If $n\ge 3$, there exist sets of four distinct points in $\mathbb F_2^n$ that are affinely independent and no $a\in\AGL_n(K)$ can map them to $\mathbb F_2^2\times\{(0,\ldots,0)\}$; so $\pi_{n,4}(K)>1$ and hence $\pi_{n,4+i}(K)>1$ if $i\in\mathbb N$. Thus for $n\ge 3$ we have $m=3$. For $n=2$ we have $m\in\{3,4\}$. Hence $(1)\Rightarrow (2)$.

We prove that $(2)\Rightarrow (3)$. As $|K|=2$, this is the same as proving $(2)\Rightarrow (1)$. We have $\pi_{n,3}(\mathbb F_2)=1$ by Lemma \ref{F7}(3.b). From this and Equation (\ref{EQ21}) applied to $n=2$ or Lemma \ref{F7}(3.c) applied to $n=2$ we get that $\pi_{2,4}(\mathbb F_2)=1$. So $(2)\Rightarrow (3)$.

Clearly, $(3)\Rightarrow (1)$.
\end{proof}

\section{Automorphisms representing multiple disjoint cycles}\label{S16}

For a finite field $K$ of characteristic $p$ and $n\in\mathbb N^{\ast}\setminus\{1\}$, in this section we identify triples $(s,r,l)\in [\mathbb N^{\ast}\setminus\{1\}]\times (\mathbb N^{\ast})^2$ for which there exists $a\in\STGA_n(K)[l]$ such that $a(K)\in r\mathcal Y_s$, with a special emphasis on the cases when either $s\in\llbracket3,2p-1\rrbracket$ is odd or $p=s=2$ and $r$ is as large as possible (see Propositions \ref{PR18} and \ref{PR19}). Applications to the case $(K,s)=(\mathbb F_2,2)$ are presented in Proposition \ref{PR20}.

\begin{proposition}\label{PR18}
Let $(n,s,q,l)\in (\mathbb N^{\ast}\setminus\{1\})^2\times (\mathbb N^{\ast})^2$, $r\in \llbracket1,\lfloor\frac{|K|^n}{s}\rfloor\rrbracket$, and $K=\mathbb F_{p^q}$. Then the following properties hold.

\medskip
{\bf (1)} There exists $a\in\STGA_n(K)[l]$ such that $a(K)\in r\mathcal Y_s$, i.e., $a(K)$ is a product of $r$ disjoint $s$-cycles, provided one of the following disjoint conditions holds.

\medskip\noindent
{\bf (1.a)} We have $p>2$, $s\in \llbracket5,2p-1\rrbracket$ is odd, $r\le\frac{(p-1)|K|^{n-1}}{2p}$, and 
$l=E_{n-2,|K|-1,\frac{s-3}{2}}$ (see Equation (\ref{EQ12})), i.e., 
$$l=\Bigl(|K|-\frac{s-1}{2}\Bigr)^2(|K|-1)^2+(n-2)(|K|-1)\Bigl[(|K|-1)^2+|K|-(|K|-1)\frac{s-1}{2}\Bigr].$$
If $n=2$, then we also have $\j_a\le 4$.

\smallskip\noindent
{\bf (1.b)} We have $s=p$, $\frac{|K|^{n-1}}{p}\mid r$, $r<\frac{|K|^n}{p}$, and $l=\max\bigl(1,|K|-\frac{pr}{|K|^{n-1}}\bigr)$. If $n=2$, then for $r=\frac{|K|^{n-1}(|K|-1)}{p}$ we also have $\j_a=0$ and for $r=\frac{i|K|^{n-1}}{p}$ with $i\in\llbracket1,|K|-2\rrbracket$ we also have $\j_a=1$.

\medskip
{\bf (2)} Suppose that $p>2$, $s\in \llbracket5,p\rrbracket$ is odd, $r\le\lfloor\frac{p}{s}\rfloor\frac{|K|^{n-1}}{p}$, and 
$l=E_{n-2,|K|-1,\frac{s-3}{2}}$. Then there exists $b\in\STGA_n(K)[l]$ such that $\w_{\supp\bigl(b(K)\bigr)}\le n-1$ and $b(K)\in r\mathcal Y_s$. In particular, if $n=2$, then $\t_{\supp\bigl(b(K)\bigr)}=1$ and $\j_b\le 4$.
\end{proposition}

\begin{proof}
Let $P_1:=(1,0,\ldots,0)\in K^n$ and $P_2:=(0,1,0,\ldots,0)\in K^n$. 

For parts (1.a) and (2), let $K=\mathbb F_p^q$ be an identification of $\mathbb F_p$-vector spaces under which $1\in K$ is identified with the $q$-tuple $(1,0,\ldots,0)$. Identifying also $\mathbb F_p=\{0,\ldots,p-1\}$, let 
$$Z:=\{2,4,\ldots,p-1\}\;\;\;\textup{and}\;\;\;Y_0:=Z\times \mathbb F_p^{q-1}\subset K.$$ 
As $|Y_0|\times |K|^{n-2}=\frac{(p-1)|K|^{n-1}}{2p}\ge r$, there exists $Y_1\subset Y_0\times K^{n-2}\subset K^{n-1}$ such that $|Y_1|=r$. Let $Y:=\{0\}\times Y_1\subset K^n$; so $|Y|=r$. The sets $Y$, $iP_1+Y$ with $i\in\mathbb F_p^{\ast}$ and $P_2+Y$ are pairwise disjoint. Based on Equation (\ref{EQ10}) applied to $n-1$, let 
$$c_1:=\e\bigl(x_1+\sum_{P\in Y_1} \chi^P(x_2,\ldots,x_n),x_2,\ldots,x_n\bigr)\in\STGA_n(K)[(n-1)(|K|-1)].$$ 
Let $h_1:K\rightarrow\{0,1\}$ be the function such that $h_1^{-1}(1)=\llbracket1,\frac{s-1}{2}\rrbracket$. Let $g_1\in K[x_1]$ of degree at most $|K|-1$ be such that it represents $h_1$ by Lagrange interpolation. Let 
$$c_2:=\e\bigl(x_1,x_2+g_1(x_1),x_3,\ldots,x_n\bigr)\in\STGA_n(K)[|K|-1].$$ 

For each element $(\alpha_1,\ldots,\alpha_n)\in K^n$ we have identities
$$c_1(K)(\alpha_1,\alpha_2,\ldots,\alpha_n)=\left\{
\begin{array}{l}
(\alpha_1+1,\alpha_2,\ldots,\alpha_n) \;\;\;\;\;\;\;\; \textup{if} \;\;\; (\alpha_2,\ldots,\alpha_n)\in Y_1\\
(\alpha_1,\alpha_2,\ldots,\alpha_n) \quad\quad\;\;\;\;\;\;\; \textup{if} \;\;\; (\alpha_2,\ldots,\alpha_n)\notin Y_1,
\end{array}\right.$$
$$c_2(K)(\alpha_1,\alpha_2,\ldots,\alpha_n)=\left\{
\begin{array}{l}
(\alpha_1,\alpha_2+1,\alpha_3,\ldots,\alpha_n) \;\;\;\; \textup{if} \;\;\;\; \alpha_1\in \llbracket1,\frac{s-1}{2}\rrbracket\\
(\alpha_1,\alpha_2,\alpha_3,\ldots,\alpha_n) \quad\quad\;\;\; \textup{if} \;\;\; \alpha_1\in K\setminus \llbracket1,\frac{s-1}{2}\rrbracket.
\end{array}\right.$$
So $c_1(K)=\prod_{Q\in Y} \theta_Q$ with each $\theta_Q$ a $p$-cycle with $\supp(\theta_Q)=Q+\mathbb F_pP_1$ and $(c_2c_1c_2^{-1})(K)=\prod_{Q\in Y} \vartheta_Q$ with each $\vartheta_Q$ a $p$-cycle with $\supp(\vartheta_Q)=c_2\bigl(\supp(\theta_Q)\bigr)$. 

Let $Q\in Y$. We have $\supp(\theta_Q)\cap\supp(\vartheta_Q)=\{Q+iP_1|i\in \mathbb F_p\setminus \llbracket1,\frac{s-1}{2}\rrbracket\}$ and hence $|\supp(\theta_Q)\cap\supp(\vartheta_Q)|=p-\frac{s-1}{2}$. If $Q_1\in Y\setminus\{Q\}$, then we check that $\supp(\theta_Q)\cap\supp(\vartheta_{Q_1})=\emptyset$. This is clear if there exists $i\in \llbracket3,n\rrbracket$ such that the $i$-th coordinates of $Q$ and $Q_1$ are distinct. If the $i$-th coordinates of $Q$ and $Q_1$ are the same for each $i\in \llbracket3,n\rrbracket$, then the second coordinates of $Q$ and $Q_1$ are distinct in $Z$, hence $\supp(\theta_Q)\cap\supp(\vartheta_{Q_1})=\emptyset$ follows from the fact that $Z\cap (Z+1)=\emptyset$.

For the commutator $a:=c_2c_1c_2^{-1}c_1^{-1}\in\STGA_n(K)$ we have the product decomposition $a(K)=\prod_{Q\in Y} \vartheta_Q\theta_Q^{-1}$ by the last paragraph. From this, either via a direct computation based, with $O:=Q+P_2$, on the identities 
$$\vartheta_Q=\Bigl(Q\;O+P_1\;O+2P_1\;\cdots\;O+\frac{s-1}{2}P_1\;Q+\frac{s+1}{2}P_1\;Q+\frac{s+3}{2}P_1\;\cdots\;Q+(p-1)P_1\Bigr)$$ and $\theta_Q^{-1}=\bigl(Q+(p-1)P_1\;Q+(p-2)P_1\;\cdots\;Q+P_1\;Q\bigr)$ or from the identity $|\supp(\theta_Q)\cap\supp(\vartheta_Q)|=p-\frac{s-1}{2}$ and Property \ref{P10} we get that $a(K)$ is the product
\begin{equation}\label{EQ21.5}
\prod_{Q\in Y} \Bigl(O+P_1\;O+2P_1\;\cdots\; O+\frac{s-1}{2}P_1\; Q+\frac{s+1}{2}P_1\; Q+\frac{s-1}{2}P_1\;\cdots\;Q+P_1\Bigr)
\end{equation} 
of $r$ disjoint $s$-cycles. Thus, as $a\in\STGA_n(K)[l]$ (cf.\ proof of Proposition \ref{PR10}(1) applied to $s$ equal to $\frac{s-1}{2}$) and, if $n=2$, we have $\j_a\le 4$, part (1.a) holds. 

For part (2), let $Z_0:=\{0,s,2s,\ldots,s\lfloor\frac{p}{s}\rfloor-s\}$. We define
$$Z:=-\Bigl(\frac{s-1}{2}\Bigr)^{-1}Z_0\;\;\;\textup{and}\;\;\;Y_0:=Z\times\mathbb F_p^{q-1}\subset K.$$ 
As $Z_0\cap (Z_0-\frac{s-1}{2})=\emptyset$, we also have $Z\cap (Z+1)=\emptyset$.

As $|Y_0|\times |K|^{n-2}=\lfloor\frac{p}{s}\rfloor\frac{|K|^{n-1}}{p}\ge r$, there exists $Y_1\subset Y_0\times K^{n-2}\subset K^{n-1}$ such that $|Y_1|=r$. Let $Y$ and $a$ be defined as above. To check that $\w_{\supp\bigl(a(K)\bigr)}=n-1$ we can assume that $r=\lfloor\frac{p}{s}\rfloor\frac{|K|^{n-1}}{p}$ and $Y_1=Y_0\times K^{n-2}$. We have 
$$\supp\bigl(a(K)\bigr)=\Bigl[\Bigl\llbracket1,\frac{s-1}{2}\Bigr\rrbracket\times (1+Y_0)\bigsqcup \Bigl\llbracket1,\frac{s+1}{2}\Bigr\rrbracket\times Y_0\Bigr]\times K^{n-2}$$
by Equation (\ref{EQ21.5})
Hence for 
$$b:=\e\Bigl(x_1-\frac{s-1}{2}x_2,x_2,x_3,\ldots,x_n\Bigl)a\e\Bigl(x_1+\frac{s-1}{2}x_2,x_2,x_3,\ldots,x_n\Bigl)\in\SGA_2(K)[l],$$ $\supp\bigl(b(K)\bigr)$ is
$$\Bigl[\Bigl(\Bigl\llbracket\frac{3-s}{2},0\Bigr\rrbracket+Z_0\Bigr)\times\mathbb F_p^{q-1}\times (1+Y_0)\bigsqcup \Bigl(\Bigl\llbracket1,\frac{s+1}{2}\Bigr\rrbracket+Z_0\Bigr)\times\mathbb F_p^{q-1}\times Y_0\Bigr]\times K^{n-2}.$$
As $s\le p$, we have $\frac{s+1}{2}<p+\frac{3-s}{2}$ and thus the restriction of the linear projection $\mathbb A^n_K\rightarrow\mathbb A^{n-1}_K$ that omits the second coordinate to $\supp\bigl(b(K)\bigr)$ is injective. Thus $\w_{\supp\bigl(b(K)\bigr)}\le n-1$. If $n=2$, then $\j_a\le 4$. So part (2) holds.

For part (1.b), we consider a subset $Z\subset K$ with $|Z|=|K|-\frac{pr}{|K|^{n-1}}\le |K|-1$. We take $a:=\e\bigl(x_1,x_2+g,x_3,\ldots,x_n\bigr)\in\STGA_n(K)[s]$ with $g:=\prod_{\alpha\in Z} (x_1-\alpha)$. We have $a(K)=\prod_{Q\in(K\setminus Z)\times Z_0\times K^{n-2}} \bigl(Q\; Q+g(Q)P_2\;\cdots\;Q+(p-1)g(Q)P_2\bigr)$, with $Z_0\subset K$ a subset that maps bijectively onto $K/\mathbb F_p$. 

Moreover, $\ell_{n,K}\bigl(a(K)\bigr)=\ell(a)=\deg(g)=|Z|$ by Proposition \ref{PR13}(3) and $\j_a=0$ if $|Z|=1$ and $\j_a=1$ if $|Z|\ge 2$. As $|Z|=|K|-\frac{pr}{|K|^{n-1}}\ge 2$ except when $r=\frac{|K|^{n-1}(|K|-1)}{p}$, we get that part (1.b) holds.
\end{proof}

\begin{proposition}\label{PR19}
Let $(n,q,l)\in (\mathbb N^{\ast}\setminus\{1\})\times (\mathbb N^{\ast})^2$, $K=\mathbb F_{p^q}$, and $r\in \bigl\llbracket1,\lfloor\frac{|K|^n}{3}\rfloor\bigr\rrbracket$. Then the following properties hold.

\medskip
{\bf (1)} There exists $a\in\STGA_n(K)[l]$ such that $a(K)\in r\mathcal Y_3$, i.e., $a(K)$ is a product of $r$ disjoint $3$-cycles, provided one of the following disjoint conditions holds.

\medskip\noindent
{\bf (1.a)} We have $r\le\frac{|K|^{n-1}\lfloor\frac{p}{2}\rfloor}{p}$ and 
$$l=E_{n-2,|K|-1}=(|K|-1)^4+(n-2)(|K|-1)[(|K|-1)^2+1].$$
Moreover, if $n=2$, then we can also assume that $\j_a\le 4$.

\smallskip\noindent
{\bf (1.b)} We have $p\neq 3$, $r=\frac{|K|^n-|K|^{n-2}}{3}$, and $l=1$.

\smallskip\noindent
{\bf (1.c)} We have $3\mid |K|-2$, $r=\frac{|K|^n-|K|^{n-2t}}{3}$ with $t\in \llbracket1,\lfloor\frac{n-1}{2}\rfloor\rrbracket$, and $l=1$.

\medskip
{\bf (2)} If $l=E_{n-2,|K|-1}$ and either $r\le\lfloor\frac{p}{3}\rfloor\frac{|K|^{n-1}}{p}$ with $p\ge 3$ or $r\le\frac{|K|^{n-1}}{4}$ with $p=2$ and $|K|\ge 4$, then there exists $b\in\STGA_n(K)[l]$ such that $b(K)\in r\mathcal Y_3$ and we have an inequality $\w_{\supp\bigl(b(K)\bigr)}\le n-1$. Moreover, if $n=2$, then $\t_{\supp\bigl(b(K)\bigr)}=1$ and $\j_b\le 4$.
\end{proposition}

\begin{proof}
Part (1.a) is proved similarly to Proposition \ref{PR18}(1.a) with the same $Y_0$ if $p>2$ and with $Y_0:=\{0\}$ if $p=2$, with $Y$ constructed from $Y_0$ in the same manner, with $g_1(x_1):=1-x_1^{|K|-1}$, and with $a(K)=\prod_{Q\in Y} (Q+P_1+P_2\; Q+P_2\; Q+P_1)$ by Equation (\ref{EQ21.5}). 

For parts (1.b) and (1.c), let $c\in\SL_2(K)$ be such that it has order $3$ but $1$ is not an eigenvalue of it; we recall that we identify it with the linear automorphism in $\SGA_2(K)$ it defines. For instance, as $p\neq 3$, we can take $c=\begin{bmatrix} 
0 & -1\\
1 & -1\\
\end{bmatrix}$. Part (1.b) holds by taking $a:=c\times 1_{\mathbb A^{n-2}_{K}}$.

Part (1.c) holds by taking $a:=c\times\cdots\times c\times 1_{\mathbb A_{K}^{n-2l}}$.

Part (2) is proved similarly to Proposition \ref{PR18}(2) with the same $Y_0$ if $p>2$. If $p=2$, then we consider $\alpha\in K\setminus\mathbb F_2$. Let $Z_0$ be an $\mathbb F_2$-vector subspace of $K$ which is a direct supplement of $\{0,1,\alpha,\alpha+1\}$. Let $Z:=\alpha^{-1}Z_0\subset K$. As $Z_0\cap (Z_0+\alpha)=\emptyset$, we also have $Z\cap (Z+1)=\emptyset$. By taking $Y_0=Z\times K^{n-1}$ and
$$b:=\e(x_1+\alpha x_2,x_2,x_3,\ldots,x_n)a\e(x_1-\alpha x_2,x_2,x_3,\ldots,x_n),$$
the remaining part of the proof is the same as of Proposition \ref{PR18}(2).\end{proof}

\begin{remark}\normalfont\label{R7}
If $n\in\mathbb N^{\ast}\setminus\{1\}$, $3\mid |K|-1$, and $\alpha\in K^{\ast}$ has order $3$, then for $a:=\e(\alpha x_1,\ldots,\alpha x_{n-t},x_{n-t+1},\ldots,x_n)\in\GL_n(K)$, we have $\a(K)\in\frac{|K|^n-|K|^t}{3}\mathcal Y_3$.
\end{remark}

The variant of Lemma \ref{L9}(2) for transpositions requires the following lemma.

\begin{lemma}\label{F11} Let $n\in\mathbb N^{\ast}\setminus\{1\}$. Let $r\in \llbracket1,2^{n-1}\rrbracket$. We write $r=2^ls$ with the pair $(l,s)\in\mathbb N\times (2\mathbb N^{\ast}-1)$ such that $s$ is odd. Then there exists $$a\in\TGA_n(\mathbb F_2)[\max(1,n-l-1)]\subset \TGA_n(\mathbb F_2)[n-1]$$ 
such that $a(\mathbb F_2)$ is a product of $r$ disjoint transpositions.
\end{lemma}

\begin{proof}
We take $a:=\e\bigl(x_1,\ldots,x_l,x_{l+1}+f(x_{l+2},\ldots,x_n),x_{l+2},\ldots,x_n\bigr)$, where the polynomial $f\in\mathbb F_2[x_{l+2},\ldots,x_n]$ represents a function $h:\mathbb F_2^{n-l-1}\rightarrow\mathbb F_2$ with $|h^{-1}(1)|=s$. As $\deg(f)\le n-l-1$ by Theorem \ref{T5}(1), the fact holds.
\end{proof}

\begin{proposition}\label{PR20}
Let $n\in\mathbb N^{\ast}\setminus\{1,2\}$, $r\in \llbracket1,2^{n-2}\rrbracket$, and $\sigma\in\perm(\mathbb F_2^n)$. We write $r=2^ls$ with $(l,s)\in\mathbb N\times (2\mathbb N^{\ast}-1)$. Then the following properties hold.

\medskip
{\bf (1)} Suppose that $\sigma$ is even. Then we have inequalities
\begin{equation}\label{EQ21.1}
\ell_{n,\mathbb F_2}(\sigma)\le (n-l-1)^{2\lceil\frac{\nu_{2,1}(\sigma)}{2r}\rceil} \pi_{n,{2r}}(\mathbb F_2)^{4\lceil\frac{\nu_{2,1}(\sigma)}{2r}\rceil}
\end{equation}
and
\begin{equation}\label{EQ21.2}
\ell_{n,\mathbb F_2}(\sigma)\le (n-l-1)^{2\lceil\frac{2^n-2r-1}{2r}\rceil} \pi_{n,{2r}}(\mathbb F_2)^{1+4\lceil\frac{2^n-2r-1}{2r}\rceil}.
\end{equation}

{\bf (2)} Suppose that $\sigma$ is odd. Then we have inequalities
\begin{equation}\label{EQ21.3}
\ell_{n,\mathbb F_2}(\sigma)\le \frac{(n-1)}{(n-l-1)}[(n-l-1)\pi_{n,{2r}(\mathbb F_2)}^2]^{\lceil\frac{\nu_{2,1}(\sigma)+1}{2r}\rceil+\lceil\frac{\nu_{2,1}(\sigma)-1}{2r}\rceil}
\end{equation}
and
\begin{equation}\label{EQ21.4}
\ell_{n,\mathbb F_2}(\sigma)\le \pi_{n,{2r}}(\mathbb F_2)\frac{(n-1)}{(n-l-1)}[(n-l-1)\pi_{n,{2r}}(\mathbb F_2)^2]^{\lceil\frac{2^n-2r}{2r}\rceil+\lceil\frac{2^n-2r-2}{2r}\rceil}.
\end{equation}

{\bf (3)} Suppose that $n\ge 4$, $r=2^{n-3}$, and $\sigma$ is even. Then we have inequalities 
$$\ell_{n,\mathbb F_2}(\sigma)\le 2^{2\lceil\frac{\nu_{2,1}(\sigma)}{2^{n-2}}\rceil} \pi_{n,{2^{n-2}}}(\mathbb F_2)^{4\lceil\frac{\nu_{2,1}(\sigma)}{2^{n-2}}\rceil}$$ 
and $\ell_{n,\mathbb F_2}(\sigma)\le 64\pi_{n,{2^{n-2}}}(\mathbb F_2)^{13}$.

\smallskip
{\bf (4)} Suppose that $n\ge 4$, $r=2^{n-3}$, and $\sigma$ is odd. Then we have the inequalities
$$\ell_{n,\mathbb F_2}(\sigma)\le 2^{\lceil\frac{\nu_{2,1}(\sigma)+1}{2^{n-2}}\rceil+\lceil\frac{\nu_{2,1}(\sigma)-1}{2^{n-2}}\rceil-1} (n-1)\pi_{n,{2^{n-2}}}(\mathbb F_2)^{2\lceil\frac{\nu_{2,1}(\sigma)+1}{2^{n-2}}\rceil+2\lceil\frac{\nu_{2,1}(\sigma)-1}{2^{n-2}}\rceil}$$
and $\ell_{n,\mathbb F_2}(\sigma)\le 32(n-1)\pi_{n,{2^{n-2}}}(\mathbb F_2)^{13}$.
\end{proposition}

\begin{proof}
For parts (1) and (2), as in the proof of Proposition \ref{PR1}(1) and (2) we write $\sigma=\varsigma_1\varsigma_2$, where 
$\varsigma_1$ and $\varsigma_2$ have order $2$ and we have $\nu_{2,1}(\varsigma_1)=\nu_{2,1}(\varsigma_2)=\frac{\nu_{2,1}(\sigma)}{2}$ if $\sigma$ is even and $\nu_{2,1}(\varsigma_1)=\nu_{2,1}(\varsigma_2)+1=\frac{\nu_{2,1}(\sigma)+1}{2}$ if $\sigma$ is odd. For $i\in\{1,2\}$, let $q_i:=\lceil\frac{\nu_{2,1}(\varsigma_i)}{r}\rceil$. We write $\varsigma_i=\prod_{j=1}^{q_i} \theta_{i,j}$ where $o(\theta_{i,j})=2$ and $\n(\theta_{i,j})=2r$ for each $i\in\{1,2\}$ and every $j\in \llbracket1,q_i-1\rrbracket$, $o(\theta_{1,q_1})=2$, $\n(\theta_{1,q_1})\in 2\llbracket1,r\rrbracket$, $o(\theta_{2,q_2})\in\{1,2\}$, and $\n(\theta_{2,q_2})=\n(\theta_{1,q_1})$ if $\sigma$ is even and $\n(\theta_{2,q_2})=\n(\theta_{1,q_1})-2$ if $\sigma$ is odd. 

If $\sigma$ is even, then $q_1=q_2$ and based on Lemma \ref{F2} applied to $\sigma=\varsigma_1\varsigma_2$ and $(l,s,r)=(2^n,2,r)$, we can assume that $n(\theta_{1,q_1})=n(\theta_{2,q_2})=2r$.\footnote{Here is where the hypotheses $r\in \llbracket1,2^{n-2}\rrbracket$ is used. One can get an analog version for $r\in \llbracket2^{n-2},2^{n-1}-1\rrbracket$.}

In this paragraph we assume that $\sigma$ is odd. Thus we have $q_2\in\{q_1-1,q_1\}$. If $q_2=q_1-1$, then $\n(\theta_{2,q_2})=2r$ and $\theta_{1,q_1}$ is a transposition; so $\n(\theta_{1,q_1})=2$. If $q_2=q_1$ (so $r\ge 2$), then $o(\theta_{2,q_2})=2$ and based on Lemma \ref{F2} applied to $(l,s,r)=(2^n,2,r)$ and to $\sigma\tau=\varsigma_1(\varsigma_2\tau_2)$ with $\tau_2$ a transposition whose support is disjoint from $\supp(\theta_{2,q_2})$, we can assume that $n(\theta_{1,q_1})=n(\theta_{2,q_2})+2=2r$. 

For $i\in \llbracket1,r\rrbracket$ let $\mathcal C_i:=i\mathcal Y_2$. We have inequalities $\omega(\mathcal C_i)\le n-1$ if $i\in \llbracket1,r-1\rrbracket$, $\omega(\mathcal C_i)\le n-2$ if $i\in \llbracket1,r-1\rrbracket\cap 2\mathbb N^{\ast}$, and $\omega(\mathcal C_r)\le n-l-1$ by Lemma \ref{F11}. From this and Lemma \ref{L8}(1) applied to $\mathcal C_r$ we get the inequality 
$$\ell_{n,\mathbb F_2}(\theta_{i,j})\le (n-l-1)\pi_{n,2r}(\mathbb F_2)^2$$ 
if one of the following conditions holds: (i) we have $j\in \llbracket1,q_1-1\rrbracket$, or (ii) $\sigma$ is even, or (iii) $\sigma$ is odd, $q_2=q_1$, and $(i,j)=(1,q_1)$. Similarly, if $\sigma$ is odd, then we have $\ell_{n,\mathbb F_2}(\theta_{1,q_1})\le (n-1)\pi_{n,2r}(\mathbb F_2)^2$ and $\ell_{n,\mathbb F_2}(\theta_{2,q_2})\le (n-l-1)\pi_{n,2r}(\mathbb F_2)^2$ if $q_2=q_1+1$ and we have $\ell_{n,\mathbb F_2}(\theta_{2,q_2})\le (n-1)\pi_{n,2r}(\mathbb F_2)^2$ if $q_2=q_1$.

From the last two sentences and Axiom $A3_{n,\mathbb F_2}$, if $\sigma$ is even we get that
$$\ell_{n,\mathbb F_2}(\sigma)\le (n-l-1)^{q_1+q_2}\pi_{n,2r}(\mathbb F_2)^{2q_1+2q_2},$$
and if $\sigma$ is odd we get that
$$\ell_{n,\mathbb F_2}(\sigma)\le (n-l-1)^{q_1+q_2-1}(n-1)\pi_{n,2r}(\mathbb F_2)^{2q_1+2q_2}.$$ So Inequalities (\ref{EQ21.1}) and (\ref{EQ21.3}) hold.

For Inequalities (\ref{EQ21.2}) and (\ref{EQ21.4}), we consider $b\in\TGA_n(\mathbb F_2)[\pi_{n,2r}(\mathbb F_2)]$ such that $\n\bigl(b(\mathbb F_2)\sigma\bigr)\le 2^n-2r$; so we have $\ell_{n,\mathbb F_2}\bigl(b(\mathbb F_2)\bigr)\le\pi_{n,2r}(\mathbb F_2)$ and (by Theorem \ref{T6}(4)) $b(\mathbb F_2)$ is even. From Inequalities (\ref{EQ21.1}) and (\ref{EQ21.3}) applied to the permutation $b(\mathbb F_2)\sigma$ with $\nu_{2,1}\bigl(b(\mathbb F_2)\sigma\bigr)\le\n\bigl(b(\mathbb F_2)\sigma\bigr)-1\le 2^n-2r-1$ and (by Axioms $A2_{n,\mathbb F_2}$ and $A3_{n,\mathbb F_2}$) the inequality $\ell_{n,\mathbb F_2}(\sigma)\le \ell_{n,\mathbb F_2}\bigl(b(\mathbb F_2)\bigr)\ell_{n,\mathbb F_2}\bigl(b(\mathbb F_2)\sigma\bigr)$, we get that Inequalities (\ref{EQ21.2}) and (\ref{EQ21.4}) hold. So parts (1) and (2) hold.

For parts (3) and (4) we have $l=n-3$ and thus $n-l-1=2$. So parts (3) and (4) follow from parts (1) and (2) (respectively) and the identity $\lceil\frac{2^n-2^{n-2}}{2^{n-2}}\rceil=3$.\end{proof}

\section{Upper bounds for small $m$}\label{S17} 

Our upper bounds for the $\pi_{n,m}(K)$s when $|K|\ge m$ or $n$ is greater than suitable functions on $K$ and $m$ are grouped in the next two theorems.

We begin with the following general lemma.

\begin{lemma}\label{L11.5}
Suppose that the pair $(n,m)\in (\mathbb N^{\ast}\setminus\{1\})^2$ is such that $m\in\llbracket 3,n+1\rrbracket$. Let $(d_1,d_2)\in\llbracket1,n\rrbracket^2$. Let $(\underline{P},\underline{Q})=\bigl((P_1,\ldots,P_m),(Q_1,\ldots,Q_m)\bigr)\in\mathbb D_{n,m}(K)^2$ be such that for $Y_1:=\{P_i|i\in\llbracket 1,m\rrbracket\}$ and $Y_2:=\{Q_i|i\in\llbracket 1,m\rrbracket\}$ we have $\d_{Y_1}=d_1$ and $\d_{Y_2}=d_2$. Then the following properties hold.

\medskip
{\bf (1)} If $(d_1,d_2,n)\neq (m-1,m-1,m-1)$ or $|K|=2$, then 
$$\pi_{\underline{P},\underline{Q}}\le \pi^{\S}_{\underline{P},\underline{Q}}\le \r_{Y_1}\r_{Y_2}\le\min(\l_{Y_1},\overline{\s}_{Y_1})\min(\l_{Y_2},\overline{\s}_{Y_2})\le (m-d_1)(m-d_2).$$ 
In particular, we have 
$\pi_{\underline{P},\underline{Q}}\le \min\bigl(m-\v_{m,|K|},\overline{\s}_{|K|}(\v_{m,|K|},m))\bigr)^2$.

\smallskip
{\bf (2)} If $d_1=d_2=n=m-1$ and $|K|\neq 2$, then $\pi_{\underline{P},\underline{Q}}=1$ and $\pi^{\S}_{\underline{P},\underline{Q}}\in\{1,2\}$.

\smallskip
{\bf (3)} Suppose that $m\ge 6$, $\{d_1,d_2\}=\{1,2\}$, and $\min(\c_{Y_1},\c_{Y_2})\le m-2$. Then $\pi^{\S}_{\underline{P},\underline{Q}}\le (m-1)(m-3)$.
\end{lemma}

\begin{proof}
For $\iota\in\{1,2\}$, up to a special linear automorphism we can assume that $\langle Y_{\iota}\rangle_{\aff}$ is the zero locus $x_{\d_{Y_{\iota}}+1}=\cdots=x_n=0$ and let $X_{\iota}$ be the zero locus $x_1=\cdots=x_{\d_{Y_{\iota}}}=0$ in $\mathbb A^n_K$. For a suitable function $f:Y_{\iota}\rightarrow X_{\iota}(K)$, the graph consists of $m$ points that are affinely independent. Such a function $f$ is represented by an $(n-\d_{Y_{\iota}})$-tuple $(f_{\d_{Y_{\iota}}+1},\ldots,f_n)\in K[x_1,\ldots,x_{\d_{Y_{\iota}}}]^{n-\d_{Y_{\iota}}}$ with 
$$\deg(f_i)\le\r_{Y_{\iota}}\le \min(\l_{Y_{\iota}},\overline{\s}_{Y_{\iota}})\le\min\bigl(m-\d_{Y_{\iota}},\overline{\s}(\v_{Y_{\iota}},m)\bigr)$$ 
for each $i\in \llbracket\d_{Y_{\iota}}+1,n\rrbracket$ by Theorem \ref{T5}(1). Thus for 
$$c_{\iota}:=\e(x_1,\ldots,x_{\d_{Y_{\iota}}},x_{\d_{Y_{\iota}+1}}+f_{\d_{Y_{\iota}+1}},\ldots,x_n+f_n)\in\STGA_n(K)[\r_{Y_{\iota}}],$$ the sets $\{c_1(P_i)|i\in\llbracket1,m\rrbracket\}$ and $\{c_2(Q_i)|i\in\llbracket1,m\rrbracket\}$ are affinely independent. 

Let $c\in\AGL_n(K)$ be such that $c\bigl(c_1(\underline{P})\bigr)=c_2(\underline{Q})$. For $a:=c_2^{-1}cc_1\in\TGA_n(K)$, we have $a(\underline{P})=\underline{Q}$ and $\ell(a)\le \r_{Y_1}\r_{Y_2}$. We can choose $c_1$ and $c_2$ such that $c$ has Jacobian determinant $1$, and thus we have $a\in\STGA_n(K)[\r_{Y_1}\r_{Y_2}]$, in all cases except when $\d_{Y_1}=\d_{Y_2}=n=m-1$ and $|K|\neq 2$. So part (1) holds.

For part (2), as $\d_{Y_1}=\d_{Y_2}=n$, up to special linear automorphisms we can assume that $P_1=(0,\ldots,0)$, that $(P_2,\ldots,P_{n+1})$ is the standard $K$-basis of $K^n$, that $Q_i=P_i$ for each $i\in \llbracket0,n\rrbracket$, and that $Q_{n+1}=\alpha P_n$ for some $\alpha\in K\setminus\{0,1\}$. We can also assume that $n=2$. 

We define automorphisms $b_1:=\e(x_1+\alpha x_2,x_2)$, $b_2:=\e\bigl(x_1,x_2+\alpha^{-1}x_1(x_1-1)\bigr)$, and $b_3:=\e(x_1-x_2,x_2)$. For $b:=b_3b_2b_1\in\SGA_2(K)[2]$ and $c:=\e(x_1,\alpha x_2)\in\GL_2(K)$, we have  $b(\underline{P})=c(\underline{P})=\underline{Q}$. So part (2) holds.

For part (3), as $\pi^{\S}_{\underline{P},\underline{Q}}=\pi^{\S}_{\underline{Q},\underline{P}}$, to fix the notation we can assume that $d_1=1$ and $d_2=2$. So $\c_{Y_1}=m$, $\r_{Y_1}=m-1$ by Corollary \ref{C7}(1), and $\c_{Y_2}\le m-2$. As $m\ge 6$, from Corollary \ref{C7}(2) and (3) we get that $\r_{Y_2}\le m-3$. Thus $\r_{Y_1}\r_{Y_2}\le (m-1)(m-3)$ and hence $a\in \STGA_n(K)[(m-1)(m-3)]$ by part (1). So part (3) holds.\end{proof}

\begin{theorem}\label{T7}
Let $(n,m)\in (\mathbb N^\ast\setminus\{1\})^2$ and a field $K$ be such that $|K|\ge \sqrt[n]{m}$. Let $N:=\min\bigl(m-\v_{m,|K|},\overline{\s}_{|K|}(\v_{m,|K|},m)\bigr)$. Then the following properties hold.

\medskip
{\bf (1)} Let $l\in\mathbb N^{\ast}$ be the smallest such that $\sum_{i=0}^l |K|^i>\frac{m(m-1)}{2}$. Suppose that $n\ge 2l$. Then we have the inequality
$$\pi^{\S}_{n,m}(K)\le N^2.$$ 
In particular, we have $\pi^{\S}_{n,m}(K)\le\min\bigl(m-\v_{m,|K|},\v_{m,|K|}(|K|-1)\bigr)^2\le (m-1)^2$ and if moreover $m\neq |K|^{\v_{m,|K|}}$ (equivalently, and $m<|K|^{\v_{m,|K|}}$) then we have inequalities $\pi^{\S}_{n,m}(K)\le\min\bigl(m-\v_{m,|K|},\v_{m,|K|}(|K|-1)-1\bigr)^2\le (m-1)^2$.\footnote{If $\v_{m,|K|}\ge 3$, then we have $\min\bigl(m-\v_{m,|K|},\v_{m,|K|}(|K|-1)-1\bigr)=\v_{m,|K|}(|K|-1)-1$ for $m\neq |K|^{\v_{m,|K|}}$ and we have $\min\bigl(m-\v_{m,|K|},\v_{m,|K|}(|K|-1)\bigr)=\v_{m,|K|}(|K|-1)$ for $m=|K|^{\v_{m,|K|}}$.}

\smallskip
{\bf (2)} Suppose that $3\le m$. Then $\bigl(\pi_{m-1+i,m}(K)\bigr)_{i\in\mathbb N}$ and $\bigl(\pi^{\S}_{m-1+i,m}(K)\bigr)_{i\in\mathbb N}$ are non-increasing sequences in $\llbracket1,N^2\rrbracket$.\footnote{Note that the implicit inequalities of part (2) do not follow directly from part (1) when $\frac{m(m-1)}{2}\ge\sum_{i=0}^{\lfloor\frac{n}{2}\rfloor} |K|^i$, i.e., in the following cases: (i) $m\in\{3,5,7,9\}$, $n\in\{m-1,m\}$, and $K\cong\mathbb F_2$; (ii) $m\in\{8,10,12\}$, $n=m-1$, and $K\cong\mathbb F_2$; (iii) $n=3$, $m=4$, and $|K|\in\{2,3,4,5\}$; (iv) $n=5$, $m=6$, and $|K|\in\{2,3\}$; (v) $n\in\{6,7\}$, $m=6$, and $K\cong\mathbb F_2$.}

\smallskip
{\bf (3)} Suppose that $3\le m=|K|$. Then $\bigl(\pi_{m-1+i,m}(K)\bigr)_{i\in\mathbb N}$ and $\bigl(\pi^{\S}_{m-1+i,m}(K)\bigr)_{i\in\mathbb N}$ are non-increasing sequences in $\llbracket1,(|K|-2)(|K|-1)\rrbracket$.
\end{theorem}

\begin{proof}
Let $(\underline{P},\underline{Q})=\bigl((P_1,\ldots,P_m),(Q_1,\ldots,Q_m)\bigr)\in\mathbb D_{n,m}(K)^2$. We consider the sets $Y_1:=\{P_i|i\in \llbracket1,m\rrbracket\}$ and $Y_2:=\{Q_i|i\in \llbracket1,m\rrbracket\}$.

For part (1), given $\iota\in\{1,2\}$, let $\pi_{\iota}:\mathbb A^n_{K}\rightarrow\mathbb A^l_{K}$ be a linear projection such that $|\pi_{\iota}(Y_{\iota})|=m$ by Lemma \ref{L10}(5). Let $\v:=\v_{m,|K|}$. 

We have inequalities $\v\le\d_{\pi_{\iota}(Y_{\iota})}\le\min(l,m-1,\d_{Y_{\iota}})$ by $m>|K|^{\v-1}$ and Proposition \ref{PR17}(2). Also, $Y_{\iota}$ has tame type $(\d_{\pi_{\iota}(Y_{\iota})},1)$. As $\sqrt[\v]{m}\le |K|$, from Theorem \ref{T4}(1) we get that $\overline{\s}_{|K|}(\d_{\pi_{\iota}(Y_{\iota})},m)$ is at most $\v(|K|-1)$ and for $\sqrt[\v]{m}< |K|$ is at most $\v(|K|-1)-1$. Therefore there exists $a\in\STGA_n(K)$ such that $a(\underline{P})=\underline{Q}$ and the inequality $\ell(a)\le\prod_{\iota=1}^2\min\bigl(\l_{\pi_{\iota(Y_{\iota})}},\overline{\s}_{\pi_{\iota}(Y_{\iota})}\bigr)\le\prod_{\iota=1}^2\min\bigl(m-\d_{\pi_{\iota}(Y_{\iota})},\overline{\s}_{|K|}(\d_{\pi(Y_{\iota})},m)\bigr)$ holds by Proposition \ref{PR17}(1), from which part (1) follows.

For parts (2) and (3), if $n\ge m-1$, then up to special linear automorphisms we can assume that the union $Y_1\cup Y_2$ is contained in the zero locus $x_m=\cdots=x_n=0$ of $\mathbb A^n_K$ and hence the sequences $\bigl(\pi_{m-1+i,m}(K)\bigr)_{i\in\mathbb N}$ and $\bigl(\pi^{\S}_{m-1+i,m}(K)\bigr)_{i\in\mathbb N}$ are non-increasing. Thus we can assume that $n=m-1$. 

For part (2), we check that either $N\ge 2$ or $N=\pi^{\S}_{m-1,m}(K)=1$. As $m\ge 3$, this is clear if $\v=1$. If $\v\ge 2$, then the field $K$ is finite. As $m\ge 3$, for $|K|\ge 3$ we estimate that $\overline{\s}_{|K|}(\v,m)\ge\overline{\s}_{|K|}(\v,3)=2$ by Lemma \ref{L6}(1), (3.a), and (3.c). If $m\ge 4$ and $|K|=2$, then we similarly estimate that $\overline{\s}_{2}(\v,m)\ge\overline{\s}_{|K|}(2,4)=2$. If $m-\v=1$ and $|K|\ge 3$ (resp.\ $m-\v=1$, $m\ge 4$, and $|K|=2$), then we have $m=\v+1>|K|^{\v-1}\ge 3^{\v-1}$ and hence $\v>3^{\v-1}-1$ (resp.\ $m=\v+1>2^{\v-1}$ and $\v\ge 3$ and hence $\v>2^{\v-1}-1$), a contradiction. Thus $N\ge 2$ except when $|K|=\v=2$ and $m=3$. If $|K|=\v=2$ and $m=3$, then from Lemma \ref{L11} we get that $N=\pi^{\S}_{m-1,m}(K)=1$. Based on the last two sentences, Lemma \ref{L11.5}(1) and (2) gives that $\pi^{\S}_{m-1,m}(K)\le N^2$. So part (2) holds.

For part (3), let $(\underline{P},\underline{Q})=\bigl((P_1,\ldots,P_m),(Q_1,\ldots,Q_m)\bigr)\in\mathbb D_{n,m}(K)^2$ and let $Y_1:=\{P_i|i\in\llbracket 1,m\rrbracket\}$ and $Y_2:=\{Q_i|i\in\llbracket 1,m\rrbracket\}$. As $m\le |K|$, we have $\v=1$. If $d_{Y_1}=\d_{Y_2}=1$, then $\c_{Y_1}=\c_{Y_2}=m$ and $\pi^{\S}_{\underline{P},\underline{Q}}\le (|K|-2)^2$ by Lemma \ref{F11.5}(2) and (4). If $\max(\d_{Y_1},\d_{Y_2})\ge 2$, then $\pi^{\S}_{\underline{P},\underline{Q}}\le (|K|-1)(|K|-2)$ by Lemma \ref{L11.5}(1) and (2). From the last two sentences we get that $\pi^{\S}_{m-1,m}(K)\le (|K|-1)(|K|-2)$. So part (3) holds.
\end{proof}

\begin{theorem}\label{T8}
Let $n\in \mathbb N^\ast\setminus\{1\}$. Then the following properties hold.

\medskip
{\bf (1)} Let $K$ be a finite field such that $|K|\ge 3$. Let $r=r(|K|)\in \llbracket3,|K|\rrbracket$ be the largest such that $4\Bigl\lfloor\frac{\frac{r(r-1)}{2}}{|K|+1}\Bigr\rfloor-4<|K|$. Let $r_0=r_0(|K|)\in\Bigl\{\Bigl\lceil\sqrt{\frac{|K|}{2}}\Bigr\rceil,\Bigl\lceil\sqrt{\frac{|K|}{2}}\Bigr\rceil+1\Bigr\}$ be the smallest such that $2\lceil\frac{|K|+4}{4}\rceil\le r_0^2-r_0$. Then the following properties hold.

\medskip\noindent
{\bf (1.a)} We have $\min\bigl(|K|,\lceil\frac{|K|}{\sqrt{2}}\rceil+1\bigr)\le r\le \min\bigl(|K|,\lceil\frac{|K|}{\sqrt{2}}\rceil+2\bigr)$.
 
\smallskip\noindent
{\bf (1.b)} If $m\in\llbracket3,r\rrbracket$, then for each subset $Y$ of $K^n$ with $m$ elements we have an inequality $\t_Y\le m-2$ and thus $\pi^{\S}_{n,m}(K)\le (m-2)^2(m-1)^2$. Moreover, if $n=2$, then $Y$ has extended tame type $(1,m-2,1)$ and $\pi^{\S,\le 4}_{n,m}(K)\le (m-2)^2(m-1)^2$.

\smallskip\noindent
{\bf (1.c)} If $m\in \llbracket r+1,|K|\rrbracket$ (hence $|K|\ge 4$), then each subset $Y$ of $K^2$ with $m$ elements has extended tame type $\bigl(1,(m-r_0)(m-2),2\bigr)$ and thus we have inequalities $\t_Y\le (m-r_0)(m-2)$ and $\pi^{\S,\le 6}_{2,m}(K)\le (m-r_0)^2(m-2)^2(m-1)^2$.

\smallskip\noindent
{\bf (1.d)} If $m\in \llbracket r+1,|K|\rrbracket$, then there exists a smallest $r_1=r_1(m,|K|)\in\llbracket 3,r-1\rrbracket$ such that we have $4\Bigl\lfloor\frac{\frac{m(m-1)}{2}-\frac{r_1(r_1-1)}{2}}{|K|}\Bigr\rfloor-4<|K|$.\footnote{If $\epsilon\in(\frac{\sqrt{2}}{2},1]$ and $m\approx\epsilon|K|$, then $r_1\approx\sqrt{\epsilon^2-\frac{1}{2}}|K|$.} 

\smallskip\noindent
{\bf (1.e)} If $m\in \llbracket r+1,|K|\rrbracket$, then each subset $Y$ of $K^2$ with $m$ elements has extended tame type $\bigl(1,(r_1-2)(r_1-1)(m-2),3\bigr)$ and thus $\t_Y\le (r_1-2)(r_1-1)(m-2)$ and 
$$\pi^{\S,\le 8}_{2,m}(K)\le (r_1-2)^2(r_1-1)^2(m-2)^2(m-1)^2\le (r-3)^2(r-2)^2(m-2)^2(m-1)^2.$$

{\bf (2)} If $n>2$ and $|K|\ge m\ge 4$, then $\pi^{\S}_{n,m}(K)\le (m-1)(m-2)^2$.

\smallskip
{\bf (3)} If $|K|+1>\max\bigl(m,\frac{m^2-m}{4}\bigr)$ and $m\ge 4$, then $\pi^{\le 3}_{2,m}(K)\le (m-1)(m-2)^2$.

\smallskip
{\bf (4)} Let $\jmath:=\min(m,|K|)-1\in\mathbb N^{\ast}$. If $|K|\ge\frac{m^2-m}{4}$, then $\pi_{n,m}(K)\le\jmath^2$.\footnote{Part (4) improves on part (3) except when $2|K|+1=\frac{m^2-m}{2}$, i.e., when $(m-2)(m+1)=4|K|$. By elementary arguments the only solutions are $(m,|K|)\in \{(6,7),(11,27)\}$.}
\end{theorem}

\begin{proof}
As $3\sqrt{2}<5$, we have
$$\frac{2(\frac{|K|}{\sqrt{2}}+1)(\frac{|K|}{\sqrt{2}}+2)}{|K|+1}-4=\frac{|K|^2+(3\sqrt{2}-4)|K|}{|K|+1}<|K|,$$
hence the inequality $\min(|K|,\lceil\frac{|K|}{\sqrt{2}}\rceil+1)\le r$ holds. A similar argument shows that $\frac{2(\frac{|K|}{\sqrt{2}}+2)(\frac{|K|}{\sqrt{2}}+3)}{|K|+1}-4>|K|$. Thus $r<\lceil\frac{|K|}{\sqrt{2}}\rceil+3$, i.e., $r\le \lceil\frac{|K|}{\sqrt{2}}\rceil+2$. So part (1.a) holds.

For parts (1.b) to (4), let $(\underline{P},\underline{Q})=\bigl((P_1,\ldots,P_m),(Q_1,\ldots,Q_m)\bigr)\in\mathbb D_{n,m}(K)^2$. We consider the sets $Y_1:=\{P_i|i\in \llbracket1,m\rrbracket\}$ and $Y_2:=\{Q_i|i\in \llbracket1,m\rrbracket\}$.

For part (1.b), let $N:=\bigl\lfloor\frac{r(r-1)}{2(|K|+1)}\bigr\rfloor$; we have $4N-4<|K|$, equivalently, we have $N<\frac{|K|+4}{4}$.

As $|K|^2\ge r^2>\frac{r(r-1)}{2}\ge \frac{m(m-1)}{2}$, we have $\sum_{i=0}^2 |K|^i>\frac{m(m-1)}{2}$. Thus $\w_Y\le 2$ by Lemma \ref{L10}(5). We can assume that for the linear projection $\pi_1:\mathbb A^n_{K}\rightarrow\mathbb A^2_{K}$ on the first two coordinates we have $|\pi_1(Y_1)|=m$. As $4 N-4<|K|$, from Proposition \ref{PR16.5}(2) applied to $\pi_1(Y_1)$ we get that $\t_{\pi_1(Y_1)}\le m-2$ and $\pi_1(Y_1)$ has extended tame type $(1,m-2,1)$; hence $Y_1$ has tame type $(1,m-2)$ and, if $n=2$, extended tame type $(1,m-2,1)$. Similarly, $\t_{Y_2}\le m-2$ and $Y_2$ has tame type $(1,m-2)$ and, if $n=2$, extended tame type $(1,m-2,1)$. By Proposition \ref{PR17}(1) there exists $a\in\STGA_n(K)[(m-2)^2(m-1)^2]$ such that $a(\underline{P})=\underline{Q}$ and, if $n=2$, we have $\j_a\le 2+1+1=4$; so part (1.b) holds.

For parts (1.c) and (1.d), as $r+1\le |K|$ we have $r\ge \lceil\frac{|K|}{\sqrt{2}}\rceil+1$ by part (1.a). Thus $m-r_0\ge r+1-r_0\ge\bigl\lceil\frac{|K|}{\sqrt{2}}\bigr\rceil-\bigl\lceil\sqrt{\frac{|K|}{2}}\bigr\rceil\ge 1$ as $|K|\ge 4$.

For part (1.c) we take $n=2$. If $4\dir_Y-4<|K|$, then $Y$ has extended tame type $(1,m-2,1)$ by Proposition \ref{PR16.5}(2) and therefore it also has extended tame type $\bigl(1,(m-r_0)(m-2),2\bigr)$. Thus we can assume that $\dir_Y\ge\lceil\frac{|K|+4}{4}\rceil$. Based on this and Lemma \ref{L10}(8) we get that $\dir_Y\in\bigl\llbracket\lceil\frac{|K|+4}{4}\rceil,\lfloor\frac{|K|-2}{2}\rfloor\bigr\rrbracket$. 

Proposition \ref{PR16.5}(3) gives that $Y$ has extended tame type $\bigl(1,(m-r_{0;Y})(m-2),2\bigr)$, where $r_{0;Y}\in\{\lceil\sqrt{2\dir_Y}\rceil,\lceil\sqrt{2\dir_Y}\rceil+1\}$ is the smallest such that $2\dir_Y\le r_{0;Y}^2-r_{0;Y}$. As $2\dir_Y\ge 2\lceil\frac{|K|+4}{4}\rceil$, we have $r_0\le r_{0;Y}$. Hence $Y$ also has extended tame type $\bigl(1,(m-r_0)(m-2),2\bigr)$. Based on this and Proposition \ref{PR17}(1) we get that part (1.c) holds.

For part (1.d), we first note that $|K|\ge 4$ implies that $r\ge 4$. Let 
$$N_1:=\frac{\frac{m(m-1)}{2}- \frac{(r-2)(r-1)}{2}}{|K|}\in\mathbb Q.$$
As $r\ge \lceil\frac{|K|}{\sqrt{2}}\rceil+1>\frac{|K|}{\sqrt{2}}+1$, we estimate 
$$N_1\le \frac{\frac{|K|(|K|-1)}{2}- \frac{(r-2)(r-1)}{2}}{|K|}\le \frac{|K|-1}{2}-\frac{\frac{|K|}{\sqrt{2}}(\frac{|K|}{\sqrt{2}}-1)}{2|K|}=\frac{|K|}{4}-\frac{\sqrt{2}-1}{2\sqrt{2}}<\frac{|K|}{4}$$
and therefore $4N_1-4<|K|$. Hence $r_1\in\llbracket3,r-1\rrbracket$ exists. So part (1.d) holds.

For part (1.e) we take $n=2$. By part (1.b), the set $\{P_1,\ldots,P_{r_1}\}$ has extended tame type $(1,r_1-2,1)$. So there exists an automorphism $c_{r_1}\in\SGA_2(K)[r_1-2]$ such that the first coordinates of $c_{r_1}(P_1),\ldots,c_{r_1}(P_{r_1})$ are distinct and $\j_{c_{r_1}}\le 1$. 

Let $b_{r_1}=\e\bigl(x_1,x_2+g(x_1)\bigr)\in\SGA_2(K)[r_1-1]$ be such that for $a_{r_1}:=b_{r_1}c_{r_1}$, the second coordinates of $a_{r_1}(P_1),\ldots,a_{r_1}(P_{r_1})$ are $0$; we have $\j_{a_{r_1}}\le 2$. This implies that there exists $\lambda\in\mathbb P_K^1(K)\setminus\{(1:0)\}$ with 
$$|\psi_{a_{r-1}(Y_1)}^{-1}(\lambda)|\le\Bigl\lfloor\frac{\frac{m(m-1)}{2}-\frac{r_1(r_1-1)}{2}}{|K|}\Bigr\rfloor<\frac{|K|+4}{4}.$$ So, as in the proof of part (1.b), we get that there exists $a_m\in\SGA_2(K)[m-2]$ such that for $b:= a_ma_{r_1}$, the first coordinates of $b(P_1),\ldots,b(P_m)$ are distinct and $\j_b\le 3$. Thus $Y_1$ has extended tame type $\bigl(1,(r_1-2)(r_1-1)(m-2),3\bigr)$ and we have $\t_{Y_1}\le (r_1-2)(r_1-1)(m-2)$. A similar argument gives that $Y_2$ has extended tame type $\bigl(1,(r_1-2)(r_1-1)(m-2),3\bigr)$. 

We have $\pi_{\underline{P},\underline{Q}}^{\S,\le 8}\le (r_1-2)^2(r_1-1)^2(m-2)^2(m-1)^2]$ by Proposition \ref{PR17}(1); so part (1.e) holds.

For part (2), for $i\in\{1,2\}$ we have $\w_{Y_i}\le n-1$ by Lemma \ref{L10}(5). Up to linear automorphisms we can assume that for $i\in\{1,2\}$ the image $Z_i$ of $Y_i$ under the projection to the first $n-1$ coordinates has $m$ elements. To fix the ideas, we can assume that if only one of the sets $Z_1$ and $Z_2$ is collinear, then that set is $Z_1$. We consider three disjoint cases as follows.

{\bf Case 1: $Z_1$ and $Z_2$ are collinear.} Both $Y_1$ and $Y_2$ have tame type $(1,1)$. Thus we have an inequality $\pi^{\S}_{\underline{P},\underline{Q}}\le (m-1)^2$ by Proposition \ref{PR17}(2).

{\bf Case 2: $Z_1$ is collinear and $Z_2$ is non-collinear.} Up to a linear automorphism we can assume that the $n$-th coordinates of the $P_i$s are distinct. Based on Corollary \ref{C7}(2) and (3), let $f\in K[x_1,\ldots,x_{n-1}]$ of degree at most $m-2$ be such that for the automorphism $c:=\e(x_1,\ldots,x_{n-1},f+x_n)\in\STGA_n(K)[m-2]$, the $n$-th coordinates of $c(Q_i)$ and $P_i$ are equal for each $i\in \llbracket1,m\rrbracket$. There exists $b=\e\bigl(x_1+h_1(x_n),\ldots,x_{n-1}+h_{n-1}(x_n),x_n\bigr)$ with $h_1,\ldots,h_{n-1}$ Lagrange interpolation polynomials of degrees at most $m-1$ such that we have $b(P_i)=c(Q_i)$ for all $i\in \llbracket1,m\rrbracket$ by Theorem \ref{T5}(1). For $a:=c^{-1}b\in\STGA_n(K)[(m-2)(m-1)]$ we have $a(\underline{P})=\underline{Q}$. So $\pi^{\S}_{\underline{P},\underline{Q}}\le (m-2)(m-1)$.

{\bf Case 3: $Z_1$ and $Z_2$ are non-collinear.} Let $(\delta_1,\ldots,\delta_m)\in\mathbb D_{1,m}(K)$. As in Case 2 we argue that there exists $(c_1,c_2)\in\STGA_n(K)[m-2]^2$ such that the $n$-th coordinates of $c_1(P_i)$ and $c_2(Q_i)$ are $\delta_i$ for every $i\in \llbracket1,m\rrbracket$. Let $(h_1,\ldots,h_{n-1})$ be the tuple of Lagrange interpolation polynomials in $K[x]$ of degrees at most $m-1$ such that for $b:=\e\bigl(x_1+h_1(x_n),\ldots,x_{n-1}+h_{n-1}(x_n),x_n\bigr)\in\STGA_n(K)[m-1]$ we have $b\bigl(c_1(\underline{P})\bigr)=c_2(\underline{Q})$. For $a:=c_2^{-1}bc_1\in\STGA_n(K)[(m-1)(m-2)^2]$ we have $a(\underline{P})=\underline{Q}$. So $\pi^{\S}_{\underline{P},\underline{Q}}\le (m-1)(m-2)^2$.

Part (2) follows from the three cases above as for $m\ge 4$ we have inequalities $(m-1)(m-2)<(m-1)^2<(m-1)(m-2)^2$.

For part (3) and the case $n=2$ of part (4), up to affine automorphisms we can assume that the projection $\pi_2:\mathbb A^2_{K}\rightarrow\mathbb A^1_{K}$ on the second coordinate is such that for $\iota\in\{1,2\}$ the sum $\Sigma_{\iota}:=\sum_{\beta\in K} \frac{|\pi_2(K)^{-1}(\beta)\cap Y_{\iota}|\bigl(|\pi_2(K)^{-1}(\beta)\cap Y_{\iota}|-1\bigr)}{2}$ is at most $N:=\bigl\lfloor\frac{m(m-1)}{2(1+|K|)}\bigr\rfloor$ by Lemma \ref{L10}(2). As $\frac{m(m-1)}{2(1+|K|)}<2$ by our hypotheses, we have $N\in\{0,1\}$. Writing $P_i=(\alpha_{i1},\alpha_{i2})\in K^2$ and $Q_i=(\beta_{i1},\beta_{i2})\in K^2$, for the two sets $Y_3:=\{\alpha_{i2}|i\in \llbracket1,m\rrbracket\}\subset K$ and $Y_4:=\{\beta_{i2}|i\in \llbracket1,m\rrbracket\}\subset K$ we have $\{|Y_3|,|Y_4|\}\subset\{m-1,m\}$ as $\{\Sigma_1,\Sigma_2\}\subset\{0,1\}$. We prove part (3), so $|K|\ge m$, by considering three disjoint cases as follows.

{\bf Case I: $|Y_3|=|Y_4|=m$.} As $Y_1$ and $Y_2$ have tame type $(1,1)$, we have an inequality $\pi^{\S,\le 2}_{\underline{P},\underline{Q}}\le (m-1)^2$ by Proposition \ref{PR17}(2).

{\bf Case II: $|Y_3|=|Y_4|=m-1$.} We consider $(\delta_1,\ldots,\delta_m)\in\mathbb D_{1,m}(K)$. For $j\in\{1,2\}$ there exists a pair $(\delta_j,f_j)\in K^{\ast}\times K[x_2]$ with $\deg(f_j)\le m-2$ such that for $c_j:=\e\bigl(\delta_jx_1+f_j(x_2),x_2\bigr)\in\GA_2(K)[m-2]$ we have $c_1(P_i)=(\delta_j,\alpha_{i2})$ and $c_2(Q_i)=(\delta_j,\beta_{i2})$ for each $i\in \llbracket1,m\rrbracket$ by Corollary \ref{C8}(1) applied to $n=2$. If $b:=\e(x_1,x_2+h_1(x_1))\in\SGA_2(K)$ with $h_1\in K[x_1]$ the Lagrange interpolation polynomial of degree at most $m-1$ such that $b\bigl(c_1(\underline{P})\bigr)=c_2(\underline{Q})$, then for the automorphism $a:=c_2^{-1}bc_1\in\GA_2(K)[(m-1)(m-2)^2]$ we have $a(\underline{P})=\underline{Q}$ and $\j_a\le 3$; so $\pi^{\le 3}_{\underline{P},\underline{Q}}\le (m-1)(m-2)^2$.

{\bf Case III: $|Y_3|\neq |Y_4|$.} To fix the ideas we assume that $|Y_3|=m$. By taking $\delta_i=\alpha_{i2}$ for each $i\in \llbracket1,m\rrbracket$ and $c_1:=\e(x_2,x_1)$, Case II applies and gives an automorphism $a\in\GA_2(K)[(m-2)(m-1)]$ such that $a(\underline{P})=\underline{Q}$ and $\j_a\le 2$. So $\pi^{\le 2}_{\underline{P},\underline{Q}}\le (m-2)(m-1)$.

Similar to part (2), part (3) follows from Cases I to III.

We prove part (4) for $n=2$. As $|K|\ge\frac{m^2-m}{4}$, for $j\in\{1,2\}$ there exists a linear projection $\pi_{0j}:\mathbb A^2_{K}\rightarrow\mathbb A^1_{K}$ with $|\pi_{0j}(Y_j)|=m$ or there exist two linear projections $\pi_{1j},\pi_{2j}:\mathbb A^2_{K}\rightarrow\mathbb A^1_{K}$ with $|\pi_{1j}(Y_j)|=|\pi_{2j}(Y_j)|=m-1$ by Lemma \ref{L10}(4) and (6). In particular, we have $|K|\ge m-1$. 

If $\pi_{01}$ and $\pi_{02}$ exist, then $\pi^{\S,\le 2}_{\underline{P},\underline{Q}}\le (m-1)^2$ by Proposition \ref{PR17}(2) applied to the tame types $(1,1)$. 

\phantomsection{If $\pi_{11},\pi_{21},\pi_{12}$, and $\pi_{22}$ exist (e.g., this holds if $|K|=m-1\le 3$), let $\iota\in\{1,2\}$ be such that the unique pair $(i_1,j_1)\in \llbracket1,m\rrbracket^2$ with $i_1<j_1$ and $\pi_{11}(P_{i_1})=\pi_{11}(P_{j_1})$ is distinct from the unique pair $(i_2,j_2)\in \llbracket1,m\rrbracket^2$ with $i_2<j_2$ and $\pi_{\iota2}(Q_{i_2})=\pi_{\iota2}(Q_{j_2})$; to fix the ideas we assume that $\iota=1$, so}\label{PH92} 
$$\underline{O}:=\bigl((\pi_{11}(P_1),\pi_{12}(Q_1)),\ldots,(\pi_{11}(P_m),\pi_{12}(Q_m)\bigr)\in\mathbb D_{2,m}(K).$$ 
There exists $a_j\in\GA_2(K)[m-2]$ such that $a_j(\underline{P})=\underline{O}$ by Corollary \ref{C8}(1) applied to $n=2$; we can assume that $\j_{a_j}\le 1$ by the proof of Corollary \ref{C8}(1). Thus for $a:=a_2^{-1}a_1\in\GA_2(K)[(m-2)^2]$ we have $a(\underline{P})=\underline{Q}$ and $\j_a\le 2$; so $\pi^{\le 2}_{\underline{P},\underline{Q}}\le (m-2)^2$.

If only one of $\pi_{01}$ and $\pi_{02}$ exists, then the last two paragraphs can be combined to give the existence of $a:=a_2^{-1}a_2\in\GA_2(K)[(m-2)(m-1)]$ such that $a(\underline{P})=\underline{Q}$ and $\j_a\le 2$; so $\pi^{\le 2}_{\underline{P},\underline{Q}}\le (m-2)(m-1)$. 

From the last three paragraphs we get that $\pi_{2,m}(K)\le (m-1)^2$ if $|K|\ge m$ and $\pi_{2,m}(K)\le (m-2)^2$ if $|K|=m-1$. Thus part (4) holds if $n=2$.

As $\frac{m(m-1)}{2}\le 2|K|<1+|K|+|K|^2$, to prove part (4) for $n\ge 3$ we can assume that for the projection $\pi_1:\mathbb A^n_{K}\rightarrow\mathbb A^2_{K}$ on the first two coordinates we have $|\pi_1(Y_1)|=|\pi_1(Y_2)|=m$ by Lemma \ref{L10}(5) applied to $l=2$. The remaining part of the proof is as in the case $n=2$ but applied to $\bigl(\pi_1(Y_1),\pi_1(Y_2)\bigr)$ instead of to $(Y_1,Y_2)$ using the case $n\ge 3$ of Proposition \ref{PR17}(2) and Corollary \ref{C8}(1).\end{proof}

For the next applications we need the following lemma. 

\begin{lemma}\label{L12}
We write $|K|=p^q$ with $p$ a prime and $q\in\mathbb N^{\ast}$. Then the following properties hold.

\medskip
{\bf (1)} For $n\in\mathbb N^{\ast}\setminus\{1\}$ there exists a selective shift $\sigma_0\in p^{q-1}\mathcal Y_p$ and hence we have inequalities 
$$\omega(p^{q-1}\mathcal Y_p)\le \omega^{\S}(p^{q-1}\mathcal Y_p)\le (n-1)(p^q-1)\le\Omega(p^{q-1}\mathcal Y_p)\le \Omega^{\S}(p^{q-1}\mathcal Y_p).$$ 

{\bf (2)} If $q=n=2$ and $p$ is odd, then 
$$\omega(\mathcal Y_{p^2})\le \omega^{\S}(\mathcal Y_{p^2})\le E_{n-2,p^2-1,\frac{p-3}{2}}(p-1)^2(|K|-1).$$
\end{lemma}

\begin{proof}
Part (1) follows from Lemma \ref{F5}(2) and Proposition \ref{PR13}(6) applied to $\bigl(\dim_K(W),K\bigr)$ equal to $(n-1,p)$.

For part (2), we define $\theta:=\bigl((0,0)\,(1,0)\,\ldots\,(p-1,0)\bigr)\in\mathcal Y_p$, $\sigma_0:=\shift_{K\times K}^{(\alpha,0)}$ with $\alpha\in K\setminus\mathbb F_p$, and $\sigma:=\theta\sigma_0$. We have $\supp(\sigma)=\supp(\sigma_0)=K\times\{0\}$ and $\sigma\in\mathcal Y^{\d=1}_{p^2}$. Let $\vartheta\in\mathcal Y_p$ be such that we have $\ell^{\S}_{2,K}(\vartheta)\le E_{n-2,k,\frac{p-3}{2}}$ by Proposition \ref{PR10}(1) applied to $s=\frac{p-1}{2}$; in fact we can assume that we have $\vartheta=a(K)\theta(K)^{-1}a(K)^{-1}$ with $\break a=\e\bigl(x_1,x_2+f(x_1)\bigr)\in\SGA_2(K)[p-1]$ by either Equation (\ref{EQ13.1}) or Equation (\ref{EQ21.5}) applied to $r=1$. So $\ell_{2,K}^{\S}(\theta)\le E_{n-2,k,\frac{p-3}{2}}(p-1)^2$. From this, the inequality $\ell_{2,K}^{\S}(\sigma)\le \ell_{2,K}^{\S}(\theta)\ell_{2,K}^{\S}(\sigma_0)$, and $\ell_{2,K}^{\S}(\sigma_0)=|K|-1$ by Lemma \ref{PR13}(6), we get that $\ell_{2,K}^{\S}(\sigma)\le E_{n-2,p^2-1,\frac{p-3}{2}}(p-1)^2(|K|-1).$ from which part (2) follows.\footnote{From Lemma \ref{L10}(5) applied to $(n,l)=(2,1)$ and the hypotheses $q=2$ we get that $\w_{\supp(\theta)}=\w_{\supp(\vartheta)}=1$ and thus $\supp(\theta)$ and $\supp(\vartheta)$ have tame type $(1,1)$. From this and Corollary \ref{C15}(1) we only get that $\ell^{\S}_{2,K}(\theta)\le E_{n-2,k,\frac{p-3}{2}}(p-1)^4$.}
\end{proof}

We have the following direct consequence of Theorem \ref{T8}(1.b) and prior results.

\begin{corollary}\label{C17}
Let $K$ be a finite field with $|K|\ge 5$. Let $\epsilon\in\{1,2\}$ be such that the largest $t_0\in\llbracket3,|K|\rrbracket$ with $|K|>4\Bigl\lfloor\frac{\frac{t_0(t_0-1)}{2}}{|K|+1}\Bigr\rfloor-4$ is $\bigl\lceil\frac{|K|}{\sqrt{2}}\bigr\rceil+\epsilon$ (cf.\ Theorem \ref{T8}(1.a)). Let $\varepsilon\in\llbracket0,\epsilon\rrbracket=\{0,1,\epsilon\}$ and
$$N=N_{|K|,\varepsilon}:=\Bigl\lceil\frac{|K|}{\sqrt{2}}\Bigr\rceil+\varepsilon-1\le|K|+\varepsilon-\epsilon-1.$$ 
Let $p:=\char(K)$ and let $s\in \llbracket3,2p-1\rrbracket$ be odd and such that for $r:=\Bigl\lfloor \frac{\lceil\frac{|K|}{\sqrt{2}}\rceil+\varepsilon}{s}\Bigr\rfloor$ we have $r\in\mathbb N^{\ast}$. We consider the non-negative rational number 
$$t:=\begin{cases} \lfloor\frac{p}{s}\rfloor\frac{|K|}{pr}\quad\quad\quad\quad\quad {\rm if}\; p\ge 3\\
\frac{|K|}{4r} \quad\quad\quad\quad\quad\quad\;\; {\rm if}\; p=2.\end{cases}$$
Let 
$$N_1=N_{1,|K|,\varepsilon}:=\begin{cases}[N(N-1)]^2[(|K|-1)N\sqrt{N-1}]^{4\wp_{r,|K|^2-\lceil\frac{|K|}{\sqrt{2}}\rceil-\varepsilon}}\quad\; {\rm if}\; t\ge 1\\
[N(N-1)]^2[(|K|-1)N(N-1)]^{4\wp_{r,|K|^2-\lceil\frac{|K|}{\sqrt{2}}\rceil-\varepsilon}}\quad\; {\rm if}\; t<1.\end{cases}$$ 

Let $\sigma\in\Alt(K^2)$ and 
$$N_{\sigma}:=\begin{cases} [(|K|-1)N\sqrt{N-1}]^{4\nu_{3,r}(\sigma)}\quad\quad\quad\quad\quad\quad\quad\; {\rm if}\; s=3\;\textup{and}\; t\ge 1\\
\left[(|K|-1)\bigl(|K|-\frac{s-1}{2}\bigr)N^2(N-1)\right]^{2\nu_{s,r}(\sigma)}\quad\; {\rm if}\; s\ge 5\;\textup{and}\;t\ge 1\;(\textup{so}\; p\ge 5)\\
[(|K|-1)N(N-1)]^{4\nu_{3,r}(\sigma)}\quad\quad\quad\quad\quad\quad\quad\; {\rm if}\; s=3\;\textup{and}\;t<1\\
\left[(|K|-1)\bigl(|K|-\frac{s-1}{2}\bigr)N^2(N-1)^2\right]^{2\nu_{s,r}(\sigma)}\;\;\; {\rm if}\; s\ge 5\;\textup{and}\;t<1.\end{cases}$$
Then the following properties hold.

\medskip
{\bf (1)} We have an inequality $\ell_{2,K}(\sigma)\le N_{\sigma}$. Moreover, for $s=3$ and $t\ge 1$ we have $\ell_{2,K}(\sigma)\le [(|K|-1)N\sqrt{N-1}]^{4\wp_{r,\n(\sigma)}}$ and $\j_{\sigma}\le 10\wp_{r,\n(\sigma)}$.

\smallskip
{\bf (2)} Suppose that $s=3$. Then we also have inequalities $\ell_{2,K}(\sigma)\le N_1$ with $\j_{\sigma}\le 10\wp_{r,|K|^2-\lceil\frac{|K|}{\sqrt{2}}\rceil-\varepsilon}+4$ if $t\ge 1$ and $\j_{\sigma}\le 12\wp_{r,|K|^2-\lceil\frac{|K|}{\sqrt{2}}\rceil-\varepsilon}+4$ if $t<1$.
\end{corollary}

\begin{proof}
We have $3\le s\le sr\le \lceil\frac{|K|}{\sqrt{2}}\rceil+\varepsilon$. Thus $sr-1\le N\le |K|+1$.

For $i\in \llbracket1,r\rrbracket$, let $\mathcal C_i:=i\mathcal Y_s$. We note that for each prime $p\ge 3$ we have inequalities $\frac{p-1}{2p}\ge \frac{1}{3}\ge\frac{1}{s}$ and for all primes $p$ we have inequalities $\frac{\lfloor\frac{p}{2}\rfloor}{p}\ge \frac{1}{3}\ge\frac{1}{s}$. Using this and the inequality $sr\le |K|$, we consider a permutation $\vartheta_j\in\mathcal C_{i_j}$ such that we have $\t_{\supp(\vartheta_j)}=1$, $\ell^{\S}_{2,K}(\vartheta_j)\le \bigl(|K|-\frac{s-1}{2}\bigr)^2(|K|-1)^2$, and $\j_{\vartheta_j}\le 4$ by Proposition \ref{PR18}(1.a) if $s\ge 5$ and Proposition \ref{PR19}(1.a) if $s=3$. 

We write $\sigma=\prod_{j=1}^{\nu_{s,r}(\sigma)} \theta_j$ with each $\theta_j\in\mathcal C_{i_j}$ for some $i_j\in \llbracket1,r\rrbracket$. 

If $t\ge 1$ and $p\ge s\ge 5$ (resp.\ and $s=3$), then as $i_j\le r$ we can assume that for each $j\in\llbracket1,\nu_{s,r}(\sigma)\rrbracket$ we also have $\t_{\supp(\vartheta_j)}=1$ by Proposition \ref{PR18}(2) (resp.\ \ref{PR19}(2)). From Corollary \ref{C15}(1) we get the following inequalities $\j_{\theta_j}\le 10$ and $\ell^{\S}_{2,K}(\theta_j)\le\ell^{\S}_{2,K}(\vartheta_j)\t_{\supp(\theta_j)}^2(si_j-1)^4$. Also, as $s\i_j\le sr\le\lceil\frac{|K|}{\sqrt{2}}\rceil+\varepsilon$, we have $\t_{\supp(\theta_j)}\le si_j-2$ by Theorem \ref{T8}(1.b). Thus, as $si_j\le sr\le N+1$, we estimate
$$\ell^{\S}_{2,K}(\theta_j)\le\ell^{\S}_{2,K}(\vartheta_j)(sr-2)^2(sr-1)^4\le \Bigl(|K|-\frac{s-1}{2}\Bigr)^2(|K|-1)^2(N-1)^2N^4.$$ 

If $t<1$, then based on Lemma \ref{L8}(2) we estimate
$$\ell^{\S}_{2,K}(\theta_j)\le \Omega^{\S}(i_j\mathcal Y_s)\le\omega^{\S}(\mathcal C_{i_j})\pi_{2,si_j}(K)^2\le \Bigl(|K|-\frac{s-1}{2}\Bigr)^2(|K|-1)^2\pi_{2,si_j}(K)^2.$$
So $\ell^{\S}_{2,K}(\theta_j)\le \left(|K|-\frac{s-1}{2}\right)^2(|K|-1)^2\pi_{2,sr}(K)^2$. As $\pi_{2,sr}(K)\le (sr-1)^2(sr-2)^2$ by Theorem \ref{T8}(1.b), we get that $\pi_{2,sr}(K)\le N^2(N-1)^2$. It follows that we have $\ell^{\S}_{2,K}(\theta_j)\le\left(|K|-\frac{s-1}{2}\right)^2(|K|-1)^2N^4(N-1)^4$ and a similar argument gives that $\j_{\theta_j}\le 12$ for each $j\in \llbracket1,\nu_{s,r}(\sigma)\rrbracket$. 

From the last two paragraphs and Axiom $A3_{2,K}$ we get that $\ell_{2,K}(\sigma)\le N_\sigma$ and $j_{\sigma}\le 10\nu_{s,r}(\sigma)$ if $t\ge 1$ and $j_{\sigma}\le 12\nu_{s,r}(\sigma)$ if $t<1$. As $\nu_{3,r}(\sigma)\le\wp_{r,\n(\sigma)}$ by Proposition \ref{PR3}(4), part (1) holds.

For part (2), let $a\in\SGA_2\bigl[\pi_{2,\lceil\frac{|K|}{\sqrt{2}}\rceil+\varepsilon}(K)\bigr]\subset \SGA_2(K)[N^2(N-1)^2]$ be such that $\supp\bigl(\sigma a(K)^{-1}\bigr)\le |K|^2-\bigl\lceil\frac{|K|}{\sqrt{2}}\bigr\rceil-\varepsilon$ and $\j_a\le 4$ by Theorem \ref{T8}(1.b). As we have $\nu_{3,r}\bigl(\sigma a(K)^{-1}\bigr)\le \wp_{r,|K|^2-\lceil\frac{|K|}{\sqrt{2}}\rceil-\varepsilon}$ by Proposition \ref{PR3}(4), there exists a product decomposition $\sigma a(K)^{-1}=\prod_{j=1}^{\wp}\sigma_j$ with $\wp\in\bigl\llbracket1,\wp_{r,|K|^2-\lceil\frac{|K|}{\sqrt{2}}\rceil-\varepsilon}\bigr\rrbracket$ and with each $\sigma_j\in\cup_{i=1}^r \mathcal C_i$. We estimate
$$\j_{\sigma}\le \j_{a(K)}+\j_{\sigma a(K)^{-1}}\le \j_a+\sum_{j=1}^{\wp} \j_{\sigma_j}\le 4+\wp\max(\j_{\sigma_j}|j\in\llbracket1,\wp\rrbracket)\le\begin{cases}  10\wp+4\quad\,\textup{if}\;\; t\ge 1\\
12\wp+4\quad\,\textup{if}\;\; t<1.\end{cases}$$ 
Therefore, as above we argue that $\ell_{2,K}(\sigma)\le N_1$ and that $\j_{\sigma}\le 10\wp_{r,|K|^2-\lceil\frac{|K|}{\sqrt{2}}\rceil-\varepsilon}+4$ if $t\ge 1$ and $\j_{\sigma}\le 12\wp_{r,|K|^2-\lceil\frac{|K|}{\sqrt{2}}\rceil-\varepsilon}+4$ if $t<1$. So part (2) holds.
\end{proof}

We have the following direct consequence of Theorem \ref{T8}(1.c) and Lemma \ref{L3}.

\begin{corollary}\label{C18}
Let $p$ be a prime and $q\in\mathbb N^{\ast}$. Let $r=r(p^q)\in \llbracket3,p^q\rrbracket$ be the largest such that $p^q>4\bigl\lfloor\frac{\frac{r(r-1)}{2}}{p^q+1}\bigr\rfloor-4$. Let $r_0=r_0(p^q)\in\bigl\{\lceil\sqrt{\frac{p^q}{2}}\rceil,\lceil\sqrt{\frac{p^q}{2}}\rceil+1\bigr\}$ be the smallest such that $2\bigl\lceil\frac{p^q+4}{4}\bigr\rceil\le r_0^2-r_0$. If $q=1$ and $p\ge 5$ let $\epsilon:=5$, and if $q\ge 2$ and $p\ge 3$ let $\epsilon:=10$. Let $\nu^+_{p,p^{q-1}}(|K|^2)$ be as in Corollary \ref{C4}. Then the following properties hold.

\medskip
{\bf (1)} If $(p,q)\neq (3,1)$ and $p\ge 3$, then we have inequalities 
\begin{equation}\label{EQ22.0}
\j_{\mathbb F_{p^q}}\le\begin{cases} 18p^q\quad\quad\;\;\;\,\textup{if}\;\; \epsilon=10\\
\frac{27p^q+27}{2}\quad\,\;\;\textup{if}\;\; \epsilon=5,\end{cases}
\end{equation}
\begin{equation}\label{EQ22.1}
\Omega^{\S}(\mathcal Y_{p^q})\le\begin{cases} (p^q-r_0)^2(p^q-2)^2(p^q-1)^{\epsilon}\quad\quad\quad\quad\quad\;\textup{if}\;\; q\ge 1\\
E_{n-2,p^2-1,\frac{p-3}{2}}(p-1)^2(|K|-1)\,\quad\quad\quad\,\;\;\textup{if}\;\; q=2,\end{cases}
\end{equation}
\begin{equation}\label{EQ22.2}
\Omega^{\S}(p^{q-1}\mathcal Y_p)\le (p^q-r_0)^2(p^q-2)^2(p^q-1)^5,
\end{equation}
\begin{equation}\label{EQ22.3}
\ln \bigl(\pi_{2,p^{2q}}(\mathbb F_{p^q})\bigr)\le \frac{3p^q+3}{2}\ln\bigl[(p^q-r_0)^2(p^q-2)^2(p^q-1)^{\epsilon}\bigr],
\end{equation}
and
\begin{equation}\label{EQ22.4}
\begin{split}
&\ln \bigl(\pi_{2,p^{2q}}(\mathbb F_{p^q})\bigr)\le\nu^+_{p,p^{q-1}}(|K|^2)\ln\bigl[(p^q-r_0)^2(p^q-2)^2(p^q-1)^5\bigr]\\
&\le \begin{cases}2p^q\ln\bigl[(p^q-r_0)^2(p^q-2)^2(p^q-1)^5\bigr]< 18p^q\ln(p^q)\quad\quad\quad\quad\quad\quad\quad\quad\;\textup{if}\;\; p\ge 5\\
\Bigl(2\lceil\frac{3|K|^2-1}{4|K|}\rceil+2\Bigr)\ln\bigl[(p^q-r_0)^2(p^q-2)^2(p^q-1)^5\bigr]< \frac{27p^q+72}{2}\ln(p^q)\,\,\;\;\textup{if}\;\; p=3.\end{cases}
\end{split}
\end{equation}
Moreover, if $p=3$ and $q\ge 3$, then $3p^q+3$ can be replaced by $3p^q+1=3^{q+1}+1$.

\smallskip
{\bf (2)} If $p=2$ and $q\ge 2$, then we have inequalities 
\begin{equation*}\label{EQ22.6}
\j_{\mathbb F_{2^q}}\le 18\cdot 2^q,
\end{equation*}
\begin{equation*}\label{EQ22.7}
\Omega^{\S}(2^{q-1}\mathcal Y_2)\le (2^q-r_0)^2(2^q-2)^2(2^q-1)^5,
\end{equation*}
\begin{equation}\label{EQ22.8}
\ln \bigl(\pi_{2,2^{2q}-2}(\mathbb F_{2^q})\bigr)\le 2^{q+1}\ln\bigl[(2^q-r_0)^2(2^q-2)^2(2^q-1)^5\bigr]< 18\cdot 2^q \ln(2^q).
\end{equation}

{\bf (3)} We have $\pi_{2,9}(\mathbb F_3)\le 2^{12}$ and $\j_{\mathbb F_3}\le 12$.
\end{corollary}

\begin{proof}
We have $r\in\lceil\frac{p^q}{\sqrt{2}}\rceil+2$ by Theorem \ref{T8}(1.a). If $p\ge 5$, then we have $\nu^+_{p,p^{q-1}}(|K|^2)\le 2p^q$ by Theorem \ref{P11}(1). If $p=3$, then we have
$$\nu^+_{p,p^{q-1}}(|K|^2)=\nu^+_{3,3^{q-1}}(3^{2q})\le 2\Bigl\lceil\frac{\lfloor\frac{3^{2q}-1}{4}\rfloor}{3^{q-1}}\Bigr\rceil+2<\frac{3(3^{2q}-1)}{2\cdot3^q}+4.$$ 
based on Proposition \ref{PR2}(1) to (3) applied to $l=3^{2q}\equiv 1 \pmod{4}$. Therefore for Inequalities (\ref{EQ22.3}), (\ref{EQ22.4}), and (\ref{EQ22.8}) it suffices to prove that if $(p,q)\neq (3,1)$ and $p\ge 3$, then for each $\sigma\in\perm(\mathbb F_{p^q}^2)$ we have two inequalities 
\begin{equation}\label{EQ23.1}
\ln\bigl(\ell_{2,\mathbb F_{p^q}}(\sigma)\bigr)\le \frac{3p^q+3}{2}\ln\bigl[(p^q-r_0)^2(p^q-2)^2(p^q-1)^{\epsilon}\bigr]
\end{equation}
and 
\begin{equation}\label{EQ23.2}
\ln\bigl(\ell_{2,\mathbb F_{p^q}}(\sigma)\bigr)\le \nu^+_{p,p^{q-1}}(|K|^2)\ln\bigl[(p^q-r_0)^2(p^q-2)^2(p^q-1)^5\bigr],
\end{equation}
and if $p=2$, then for each $\sigma\in\Alt(\mathbb F_{2^q}^2)$ we have the inequality
\begin{equation}\label{EQ23.3}
\ln\bigl(\ell_{2,\mathbb F_{2^q}}(\sigma)\bigr)\le 2^{q+1}\ln\bigl[(2^q-r_0)^2(2^q-2)^2(2^q-1)^5\bigr].
\end{equation}

\medskip
To prove Inequalities (\ref{EQ23.1}) to (\ref{EQ23.3}) we can assume that $\sigma$ is even by Theorem \ref{T6}(1) and (3). 

We prove that Inequality (\ref{EQ23.1}), Inequality (\ref{EQ22.0}) for $\epsilon=5$, and Inequality (\ref{EQ22.1}) hold. From Lemmas \ref{L3}(2) and \ref{L12}(2) we get that it suffices to show that we have an inequality $\Omega^{\S}(\mathcal Y_{p^q})\le (p^q-r_0)^2(p^q-2)^2(p^q-1)^{\epsilon}$ and, if $\epsilon=5$, that the Furter's length of each permutation in $\mathcal Y_{p^q}$ is at most $9$. Let $\theta\in\mathcal Y_{p^q}$ and let $Y:=\supp(\theta)$. Let $\vartheta\in\mathcal Y_{p^q}$ be such that $Z:=\supp(\vartheta)=\mathbb F_{p^q}\times\{0\}$. We have $|Y|=|Z|=p^q$.

The set $Y$ has extended tame type $\bigl(1,(p^q-r_0)(p^q-2),2\bigr)$ by Theorem \ref{T8}(1.c) and we have $\t_Z=1$. These and Corollary \ref{C15} imply that $\j_{\theta}\le\j_{\vartheta}+8$ and 
$$\ell^{\S}_{2,p^q}(\theta)\le [(p^q-r_0)^2(p^q-2)^2(p^q-1)^4]\ell^{\S}_{2,\mathbb F_{p^q}}(\vartheta).$$
Thus it suffices to prove that we can choose $\vartheta$ such that $\ell^{\S}_{2,\mathbb F_{p^q}}(\vartheta)\le (p^q-1)^{\epsilon-4}$, and if $\epsilon=5$, that $\j_{\vartheta}\le 1$

If $\epsilon=5$, then $q=1$ and we can take $\vartheta$ such that $\ell^{\S}_{2,\mathbb F_{p^q}}(\vartheta)=p-1=(p^q-1)^{\epsilon-4}$ and $\j_{\vartheta}=1$ by Proposition \ref{PR13}(6) applied to $(n,d)=(2,1)$ and $K=\mathbb F_p$; hence $\j_{\theta}\le 9$. So Inequality (\ref{EQ22.0}) for $\epsilon=5$ holds. 

So we can assume that $\epsilon=10$, i.e., $q\ge 2$. We choose the permutation $\vartheta$ such that, based on Theorem \ref{P11}(1), once we have a product decomposition $\vartheta=\vartheta_1\vartheta_2$ with $(\vartheta_1,\vartheta_2)\in (p^{q-1}\mathcal Y_p)^2$ satisfying $\supp(\vartheta_1)\cup\supp(\vartheta_2)\subset Z$, $\vartheta_1$ is a selective shift and we have $\ell^{\S}_{2,K}(\vartheta_1)=p^q-1$ and $\j_{\vartheta_1}\le 1$ by Proposition \ref{PR13}(6) applied to $(n,d)=(2,1)$ and $K=\mathbb F_{p^q}$. As $\t_{\supp(\vartheta_1)}=\t_{\supp(\vartheta_2)}=1$, from Corollary \ref{C15} we get that $\ell^{\S}_{2,p^q}(\vartheta_2)\le (p^q-1)^5$ and $\j_{\vartheta_2}\le 5$. Thus $\ell^{\S}_{2,p^q}(\vartheta)\le (p^q-1)^6$ by Axiom $A3^{\S}_{2,\mathbb F_{p^q}}$ and $\j_{\vartheta}\le 6$. So Inequalities (\ref{EQ23.1}) and (\ref{EQ22.1}) hold.

We prove that Inequality (\ref{EQ23.2}), Inequality (\ref{EQ22.0}) for $\epsilon=10$, and Inequality (\ref{EQ22.2}) hold. It suffices to show that for each permutation $\theta\in p^{q-1}\mathcal Y_p$ we have inequalities $\ell^{\S}_{2,\mathbb F_{p^q}}(\theta)\le (p^q-r_0)^2(p^q-2)^2(p^q-1)^5$ and $\j_{\theta}\le 9$. This follows from the fact that we can assume that $\ell^{\S}_{2,\mathbb F_{p^q}}(\vartheta)\le p^q-1$ and $\j_{\vartheta}=1$ by Proposition \ref{PR13}(6). Thus Inequality (\ref{EQ23.2}), Inequality (\ref{EQ22.0}) for $\epsilon=10$, and Inequality (\ref{EQ22.2}) hold. So part (1) holds.

We prove that Inequality (\ref{EQ23.3}) and part (2) hold. As $\lceil\frac{\nu_{2,1}(\sigma)+1}{2^q}\rceil\le\lceil\frac{2^{2q}}{2^q}\rceil=2^q$ and $\nu_{2,1}(\sigma)$, as in the proof of Proposition \ref{PR20}(1) applied to $(n,r)=(2q,2^{q-1})$ we argue that $\sigma$ is a product of at most $2^{q+1}$ permutations in $2^{q-1}\mathcal Y_2$. So it suffices to show that for each $\theta\in 2^{q-1}\mathcal Y_2$ we have $\ell^{\S}_{2,\mathbb F_{2^q}}(\theta)\le (2^q-r_0)^2(2^q-2)^2(2^q-1)^5$ and $\j_{\theta}\le 9$. This is proved in the same we proved Inequality (\ref{EQ23.1}), by taking $\vartheta\in 2^{q-1}\mathcal Y_2$ to have collinear support and satisfy $\ell^{\S}_{2,\mathbb F_{2^q}}(\vartheta)=2^q-1$ and $\j_{\vartheta}\le 1$ by Proposition \ref{PR13}(6) applied to $(n,d)=(2,1)$ and $K=\mathbb F_{2^q}$. So Inequality (\ref{EQ23.3}) and part (2) hold.

For part (3), we have $\Alt(\mathbb F_3^2)=\mathcal Y_3^4$ by Lemma \ref{L3}(1). From this and Proposition \ref{PR16}(7) we get that part (3) holds.\end{proof}

We exemplify some of the results of this section for finite fields with $3\le |K|\le 16$.

\begin{example}\normalfont\label{EX17.5}  Let $K$, $r$, and $r_0$ be as in Theorem \ref{T8}. Let $Y$ be a subset of $K^2$ with $|Y|=|K|$. 

\medskip
{\bf (1)} Suppose that $|K|\in\{3,4,5,7,9\}$. Then $r=|K|$. By Theorem \ref{T8}(1.b) we have $\pi^{\S,\le 4}_{2,|K|}(K)\le(|K|-2)^2(|K|-1)^2$ and $Y$ has extended tame type $(1,|K|-2,1)$. For each $\sigma\in\frac{|K|}{p}\mathcal Y_p$ we have $\j_{\sigma}\le 7=2+1+4$ by Proposition \ref{PR17}(1) applied as in the proof of Corollary \ref{C18}. As in the proof of Corollary \ref{C18}(1) we argue that for $|K|=5$ we have $\j_K\le 10.5(|K|+1)=63$, that for $|K|=7$ we have $\j_K\le 10.5(|K|+1)=84$, and that for $|K|=9$ we have $\j_K\le 14|K|=126$. 

\smallskip
{\bf (2)} Suppose that $|K|=8$. Then $r=|K|-1=7$ and $r_0=3$. By Theorem \ref{T8}(1.c) we have  $\pi^{\S,\le 6}_{2,8}(K)\le(8-3)^2(8-2)^2(8-1)^2=210^2$ and $Y$ has extended tame type $(1,30,2)$.

\smallskip
{\bf (3)} Suppose that $|K|\in\{11,13\}$. Then $r=|K|-1$ and $r_0=4$. By Theorem \ref{T8}(1.c) we have  $\pi^{\S,\le 6}_{2,|K|}(K)\le(|K|-4)^2(|K|-2)^2(|K|-1)^2$ and $Y$ has extended tame type $\bigl(1,(|K|-4)(|K|-2),2\bigr)$.

\smallskip
{\bf (4)} If $|K|\ge 16$, then $\frac{(|K|-2)(|K|-1)}{2(|K|+1)}>\frac{|K|+8}{2}$ and hence $r\le |K|-2$. 

\smallskip
{\bf (5)} Suppose that $|K|=16$. Then $r=13$ and $r_0=4$. Hence $Y$ has extended tame type $(1,168,2)$ by Theorem \ref{T8}(1.c). 
\end{example}

\section{Polynomial upper bounds for medium $m$}\label{S18}

In this section we use Weil restriction of scalars in order to obtain upper bounds for the $\pi_{n,m}^{\S}(K)$s that are polynomial in $m$ with $m\in \llbracket1,\|K|^{\frac{n}{2}}\rrbracket$ if $n$ is even (see Proposition \ref{PR21}) and with $m\in\bigl\llbracket1,\bigl\lfloor\frac{|K|^\frac{n}{2}}{2\sqrt{2}}\bigr\rfloor\bigr\rrbracket$ if $n$ is odd (see Proposition \ref{PR23}).

Let $l\in\mathbb N^{\ast}$. Recall that for a finite field extension $K\rightarrow L$ of degree $l$ we have a Weil restriction of scalars functor between categories of quasi-projective schemes
$$\R_{L/K}: (\textup{quasi-proj\;} L\textup{-schemes})\rightarrow (\textup{quasi-proj\;} K\textup{-schemes})$$
characterized by functorial identities $\R_{L/K}(\mathcal X)(S)=\mathcal X\bigl(\Spec(S\otimes_K L)\bigr)$, where $S$ is a commutative $K$-algebra and $\mathcal X$ is a quasi-projective scheme over $\Spec L$. 

The functor sends affine schemes to affine schemes and 
$\R_{L/K}(\mathbb A^n_L)\cong\mathbb A_K^{ln}$; one can specify such an isomorphism using a $K$-basis $\mathcal{B}=(\alpha_1,\ldots,\alpha_l)$ of $L$ over $K$ as for each $K$-algebra $S$ we have a functorial $S$-linear map $\Psi_S:(S\otimes_K L)^n\cong S^{ln}$ such that for each $(i,j)\in \llbracket1,n\rrbracket\times \llbracket1,l\rrbracket$, the $(i-1)l+j$-th coordinate of $\Psi_S(\underline{\beta_1},\ldots,\underline{\beta_n})$ is the $j$-th coordinate of $\underline{\beta_i}\in S\otimes_K L$ under the $S$-linear isomorphism $S\otimes_K L\rightarrow S^l$ defined by $\mathcal B$. With respect to the group scheme structures on $\mathbb A^n_L$ and $\mathbb A^{ln}_K$ defined by $\mathbb G^n_L$ and $\mathbb G^{ln}_K$ (respectively), the monoid homomorphism $\R_{L/K}:\End_n(L)\rightarrow\End_{ln}(K)$ is in fact a ring homomorphism. 

If $K\rightarrow K_1$ is a field extension linearly disjoint from $K\rightarrow L$, so $L_1:=L\otimes_K K_1$ is a field, then for each morphism $\mathcal Q:\mathcal X\rightarrow Y$ of quasi-projective $L$-schemes, we have a functorial identification
\begin{equation}\label{EQ23}
\bigl(\R_{L/K}(\mathcal Q)\bigr)_{K_1}\cong \R_{L_1/K_1}(\mathcal Q_{L_1}).
\end{equation}
For basic properties and notation on Weil restriction see \cite{CGP}, Part IV, App.\ A.5 (other common notation in the literature is $\Res_{L/K}$ or $\prod_{\Spec L/\Spec K}$).

\begin{lemma}\label{L13} With the notation of this section, the following properties hold. 

\medskip
{\bf (1)} For each $d\in\mathbb N$, the ring homomorphism $\R_{L/K}$ sends homogeneous endomorphisms of degree $d$ to homogeneous endomorphisms of degree $d$.

\smallskip
{\bf (2)} The ring homomorphism $\R_{L/K}:\End_n(L)\rightarrow\End_{ln}(K)$ is injective.

\smallskip
{\bf (3)} For each $a\in \GA_n(L)$ we have $\pi\bigl(\R_{L/K}(a)\bigr)=\pi(a)$ and therefore also $\ell\bigl(\R_{L/K}(a)\bigr)=\ell(a)$.

\smallskip
{\bf (4)} We have $\R_{L/K}\bigl(\SGA_n(L)\bigr)\subset \SGA_{nl}(K)$ and $\R_{L/K}\bigl(\TGA_n(L)\bigr)\subset \TGA_{nl}(K)$.
\end{lemma}

\begin{proof}
Part (1) is straightforward from the construction.

For part (2), based on Equation(\ref{EQ23}) we can replace $K$ by $K(x)$ and $L$ by $L(x)$, so we can assume that $K$ and $L$ are infinite, in which case the endomorphisms are determined by their action on rational points and the assertion is clear.

Part (3) follows from parts (1) and (2). 

Part (4) follows directly from the definitions. Note that for the automorphism $b=\e\bigl(x_1,\ldots,x_{n-1},x_n+f(x_1,\ldots,x_{n-1})\bigr)\in \SGA_n(L)$, $\R_{L/K}(b)$ is an automorphism $\e(x_1,\ldots,x_{(n-1)l},x_{(n-1)l+1}+g_1,\ldots,x_{nl}+g_l)$ with $(g_1,\ldots,g_l)\in K[x_1,\ldots,x_{(n-1)l}]^l$ which is the composite of $l$ commuting triangular automorphisms in $\SGA_{nl}(K)$.\end{proof}

We consider two cases, depending on whether $n$ is even or odd and we begin by computing the required values of the strict capacity functions. 

\begin{lemma}\label{L14}
Let $(k,n)\in (\mathbb N^{\ast}\setminus\{1\})^2$. Let $m:=\lceil\frac{k^{\frac{n}{2}}}{\sqrt{2}}\rceil\in\mathbb N^{\ast}$. Let $d\in\mathbb N^{\ast}$ be such that $m\le k^d$. Then the following properties hold.

\medskip
{\bf (1)} Suppose that $n=2l$ is even. Then the following properties hold.

\medskip\noindent
{\bf (1.a)} If $k\in\{2,3\}$ and $l=1$, then $m=k$ and $\overline{\s}_k(d,m)=\overline{\s}_k(d,m+1)=m-1$.

\smallskip\noindent
{\bf (1.b)} If $k\ge 4$ or $l\ge 2$, then $k^{l-1}<m<k^l$ and hence $d\ge l$. Moreover, if $l\ge 2$ then for $k\ge 3$ we have $\sum_{i=0}^{l-1} k^i<m<k^l-1$ and for $k\ge 4$ we have $2\sum_{i=0}^{l-1} k^i<m<k^l-2$.

\smallskip\noindent
{\bf (1.c)} If $k\ge 4$ and $l=1$, then $\overline{\s}_k(d,m+1)=m$ and for each $n_1\in\mathbb N$ we have
$\overline{\s}_k(d+n_1,mk^{n_1})=m-1+n_1(k-1)$. 

\smallskip\noindent
{\bf (1.d)} Suppose that $l\ge 2$. Then there exists $\varepsilon\in\{-1,0\}$ such that 
$$\overline{\s}_k(d,m) =(l-1)(k-1)+\Bigl\lfloor\frac{k}{\sqrt{2}}\Bigr\rfloor+\varepsilon\le \overline{\s}_k(d,m+1) \le (l-1)(k-1)+\Bigl\lfloor\frac{k}{\sqrt{2}}\Bigr\rfloor.$$ 

\noindent
{\bf (1.e)} If $l\ge 2$ and $k\ge 3$, then for each $n_1\in\mathbb N$ and every $j\in\{0,1,2\}$ we have 
$$\overline{\s}_k\bigl(d+n_1,(m+j)k^{n_1}\bigr)\le (l+n_1-1)(k-1)+\Bigl\lfloor\frac{k}{\sqrt{2}}\Bigr\rfloor.$$

\noindent
{\bf (1.f)}
Suppose that $l\ge 2$ and $k=2$. Then for $l\in\{2,3,4\}$ we have $m=2^l-2^{l-2}$ and for $l\ge 5$ we have a strict inequality $m< 2^l-2 ^{l-2}$. Moreover, we have identities $\overline{\s}_k(d+n_1,m2^{n_1})=\overline{\s}_k\bigl(d+n_1,2^{n_1}(m+2^{l-2}-1)\bigr)=l+n_1-1$ for each $n_1\in\mathbb N$. In particular, for $l\ge 3$, with $\varepsilon$ as in part (1.d), we have $\varepsilon=-\bigl\lfloor\frac{2}{\sqrt{2}}\bigr\rfloor=-1$.

\medskip
{\bf (2)} Suppose that $n=2l+1$ is odd. Then the following properties hold.

\medskip\noindent
{\bf (2.a)} If $k=2$, then $m=2^l$ and $\overline{\s}_2(d,m)=l$.

\smallskip\noindent
{\bf (2.b)} Suppose that $k>2$. Then $k^l<m<k^{l+1}$. Moreover, there exists $\varepsilon\in\{-1,0\}$ such that $\overline{\s}_k(d,m) =l(k-1)+\bigl\lfloor\sqrt{\frac{k}{2}}\bigr\rfloor+\varepsilon\le l(k-1)+\bigl\lfloor\sqrt{\frac{k}{2}}\bigr\rfloor$.\end{lemma}

\begin{proof}
We write $m=\frac{k^\frac{n}{2}}{\sqrt{2}}+\epsilon$ with $\epsilon\in [0,1)$; if $n$ is even, then $\epsilon\notin\mathbb Q$.

Parts (1.a) and (2.a) are clear based on Lemma \ref{L6}(2) and (3.c).

For part (1.b), as $k\ge 2$ we have $m>\frac{k^l}{\sqrt{2}}>k^{l-1}$. We show that the assumption that $m> k^l-1$ leads to a contradiction. This assumption gives $k^l+\epsilon\sqrt{2}>\sqrt{2}k^l-\sqrt{2}$ and thus $k^l<\frac{\sqrt{2}(1+\epsilon)}{\sqrt{2}-1}<\frac{2\sqrt{2}}{\sqrt{2}-1}<6.83$, hence either $l=1$ and $k\in\{4,5,6\}$ or $k=l=2$. But in all these four cases one checks directly that $m=k^l-1$, a contradiction. 

If $k\ge 3$, then $\sum_{i=0}^{l-1} k^i=\frac{k^l-1}{k-1}<\frac{k^l}{k-1}\le \frac{k^l}{2}<\frac{k^l}{\sqrt{2}}<m$. If $k\ge 4$, then
$2\sum_{i=0}^{l-1} k^i=2\frac{k^l-1}{k-1}<2\frac{k^l}{k-1}\le 2\frac{k^l}{3}<\frac{k^l}{\sqrt{2}}<m$.

If $k\ge 3$ and $l\ge 2$, then we show that the assumption that $m> k^l-2$ leads to a contradiction. This assumption gives $k^l+\epsilon\sqrt{2}>\sqrt{2}k^l-2\sqrt{2}$ and hence $k^l<\frac{\sqrt{2}(2+\epsilon)}{\sqrt{2}-1}<\frac{3\sqrt{2}}{\sqrt{2}-1}<10.25<4^2$. Thus $(k,l)=(3,2)$ and we have $m=7=3^2-2$, a contradiction. So $m\le k^l-2$ and $m<k^l-1$. As $\frac{4\sqrt{2}}{\sqrt{2}-1}<13.66<4^2$, a similar argument gives that if $k\ge 4$ and $l\ge 2$, then $m<k^l-2$. So part (1.b) holds. 

For part (1.c), as $m<k$ by part (1.b), we have identities $\overline{\s}_k(d,m)=m-1$ and $\overline{\s}_k(d,m+1)=m$ by Lemma \ref{L6}(3.c). Next we assume that $n_1>0$. As $1<m<k$ by part (1.b), we have strict inequalities $(m-1)(k^{n_1}-1)>0$, i.e., $mk^{n_1}-m+1> k^{n_1}$, and $(m-1)(\sum_{i=0}^{n_1} k^i)<mk^{n_1}-1<m(\sum_{i=0}^{n_1} k^i)$. Therefore $\bigl\lfloor\frac{mk^{n_1}-1}{\sum_{i=0}^{n_1} k^i}\bigr\rfloor=m-1$ and $\min(mk^{n_1}-m+1,k^{n_1})=k^{n_1}$; hence $\overline{\s}_k(d+n_1,m)=m-1+n_1(k-1)$ by Definition \ref{D9}(1.b) and Lemma \ref{L6}(2). So part (1.c) holds.

For parts (1.d) to (1.f), for $j\in\{0,1,2,3\}$ we consider 
$$q_j:=\frac{m-1+j}{\sum_{i=0}^{l-1} k^i}.$$ 
We estimate
\begin{equation*}
\begin{aligned}
\frac{\sqrt{2}q_j}{k}&=\frac{(k-1)[k^l+\sqrt{2}(\epsilon+j-1)]}{k^{l+1}-k}<\frac{k^{l+1}-k^l+\sqrt{2}j(k-1)}{k^{l+1}-k}\\
&=1+\frac{(k-1)[\sqrt{2}j-k(\sum_{i=0}^{l-2} k^i)]}{k^{l+1}-k}.\end{aligned}
\end{equation*}
So we get that $\frac{\sqrt{2}q_j}{k}<1$ except when $l=2$ and $(j,k)\in\{(2,2),(3,3),(3,4)\}$. If $l=2$ and $k=2$ (resp.\ $k=4$), then $m=3$ (resp.\ $m=12$) and hence $\frac{\sqrt{2}q_2}{k}=\frac{4\sqrt{2}}{6}<1$ (resp.\ $\frac{\sqrt{2}q_3}{k}=\frac{14\sqrt{2}}{20}<1$). So the inequality $\frac{\sqrt{2}q_j}{k}<1$ holds if $(l,j,k)\neq (2,3,3)$.

Similarly, we have
$$\frac{\sqrt{2}q_0}{k-\sqrt{2}}=\frac{k^{l+1}-k^l-k(\sqrt{2}-\sqrt{2}\epsilon)-\sqrt{2}\epsilon+\sqrt{2}}{k^{l+1}-\sqrt{2}k^l-k+\sqrt{2}}>1$$
as the last inequality holds due to the fact that 
$$k^l(\sqrt{2}-1)>k(\sqrt{2}-1)>k(\sqrt{2}-1)+\sqrt{2}\epsilon(1-k)=k(\sqrt{2}-\sqrt{2}\epsilon-1)+\sqrt{2}\epsilon.$$

So for part (1.d) we have inequalities $\frac{k}{\sqrt{2}}-1<q_0<q_1<\frac{k}{\sqrt{2}}$; hence there exists $\varepsilon\in\{-1,0\}$ such that $\lfloor q_0\rfloor=\lfloor\frac{k}{\sqrt{2}}\rfloor-\varepsilon$. Thus 
$$m+1-q_1> m-q_0>\frac{k^l-k}{\sqrt{2}}>k^{l-1},$$ where the last inequality holds as it is equivalent to the inequality $k>\sqrt{2}+\frac{1}{k^{l-1}}$. From the inequalities $m+1-q_1> m-q_0>k^{l-1}$, Definition \ref{D9}(1.b), and Lemma \ref{L6}(2) we get that we have the following identities $\overline{\s}_k(d,m) =(l-1)(k-1)+\lfloor q_0\rfloor$ and $\overline{\s}_k(d,m+1)=(l-1)(k-1)+\lfloor q_1\rfloor$. From this, the inequality $\lfloor q_1\rfloor\le\lfloor\frac{k}{\sqrt{2}}\rfloor$, and the monotony part of Lemma \ref{L6}(1) we get that part (1.d) holds.

For part (1.e), we have $k^{l+n_1-1}<mk^{n_1}<(m+1)k^{n_1}<k^{l+n_1}$ and $q_1>1$ by part (1.b). So, as $l\ge 2$ and $k\ge 3$, if $(l,k)\neq (2,3)$ then we have $q_3<\frac{k}{\sqrt{2}}$ and thus we estimate
$$\frac{(m+2)k^{n_1}-1}{\sum_{i=0}^{l+n_1-1} k^i}< \frac{(m+2)k^{n_1}}{\sum_{i=n_1}^{l+n_1-1} k^i}=\frac{m+2}{\sum_{i=0}^{l-1} k^i}=q_3<\frac{k}{\sqrt{2}}.$$

If $(l,k)=(2,3)$, then $m=7$, $m+2=9$, and for each $n_1\in\mathbb N$ we compute
$$\overline{\s}_k\bigl(d+n_1,(m+2)k^{n_1}\bigr)=\overline{\s}_2(n_1+2,3^{n_1+2})=2(n_1+2)=(l+n_1-1)(k-1)+\Bigl\lfloor\frac{k}{\sqrt{2}}\Bigr\rfloor,$$
where the first two identities hold by Lemma \ref{L6}(3.a).

From these estimates and identities, Definition \ref{D6}(1.b), and Lemma \ref{L6}(1) we get that part (1.e) holds.

For part (1.f), if $l=2$ then $m=2^l-2^{l-2}=3$, if $l=3$ then $m=2^l-2^{l-2}=6$, and if $l=4$ then $m=2^l-2^{l-2}=12$. If $l\ge 5$, then $\frac{\epsilon}{2^{l-1}}<\frac{1}{2^{l-1}}\le\frac{1}{16}<\frac{3}{2}-\sqrt{2}$ and thus $\frac{\epsilon}{2^{l-1}}+\sqrt{2}<\frac{3}{2}$. By multiplying the last inequality by $2^{l-1}$ we get that $m=2^{l-1}\sqrt{2}+\epsilon<3\cdot 2^{l-2}=2^l-2^{l-2}$. So the first statement of part (1.f) holds. From $m\le 2^l-2^{l-2}$ and part (1.b) we get that $2^{l-1}<m\le m+2^{l-2}-1< 2^l$. So $2^{l+n_1-1}<2^{n_1}m\le 2^{n_1}(m+2^{l-2}-1)<2^{l+n_1}$. From this, Lemma \ref{L6}(3.a), and the monotony part of Lemma \ref{L6}(1) we get that the second statement of part (1.f) also holds. So part (1.f) holds.

For part (2.b), we have $m>k^l\sqrt{\frac{k}{2}}> k^l$. As $k-\sqrt{\frac{k}{2}}> 1> k^{-l}$, we get that $k^l\sqrt{\frac{k}{2}}< k^{l+1}-1$; hence $m<k^{l+1}$ and the first sentence of part (2.b) holds. 

Let $q:=\frac{m-1}{\sum_{i=0}^{l} k^i}$. Thus $q<\frac{k^l\sqrt{\frac{k}{2}}}{1+k^l}<\sqrt{\frac{k}{2}}$. As 
$$\sum_{i=1}^l (k^i-k^{i-\frac{1}{2}})\ge 3-\sqrt{3}>\sqrt{2}-1>\sqrt{2}(1-\epsilon)-1,$$ 
it follows that $q>\frac{\sqrt{k}-1}{\sqrt{2}}$. So as in the proof of part (1.d) we get that there exists $\varepsilon\in\{-1,0\}$ such that $\lfloor q\rfloor=\bigl\lfloor\sqrt{\frac{k}{2}}\bigr\rfloor+\varepsilon$. 

We check that the inequality $(k^l-1)\bigl(\sqrt{\frac{k}{2}}-1\bigr)>1+\epsilon$ holds. For this we can assume that $l=1$ and $k=3$ and in this case it holds as $m=4$, $\lfloor q\rfloor=0$, and $\epsilon=-1$. The last inequality implies that $m-\lfloor q\rfloor\ge k^l\sqrt{\frac{k}{2}}-\sqrt{\frac{k}{2}}-\epsilon>k^l$. Hence from Definition \ref{D6}(1.b) and Lemma \ref{L6}(3.a) we get that $\overline{\s}_k(d,m) =l(k-1)+\bigl\lfloor\sqrt{\frac{k}{2}}\bigr\rfloor+\varepsilon$. So the second sentence of part (2.b) also holds. So part (2.b) holds.\end{proof}

The even case is easier, so we start with it. Essentially, we combine a particular case of Lemma \ref{L13} with a minor modification of an argument in the proof of Theorem \ref{T8}(1.b).

\begin{proposition}\label{PR21} Let a finite field $K$ and $l\in\mathbb N^{\ast}$ be such that $|K|^l\ge 3$. Let $n:=2l$ and let $r\in\llbracket3,|K|^l\rrbracket$ be the largest such that $|K|^l>4\Bigl\lfloor\frac{\frac{r(r-1)}{2}}{|K|^l+1}\Bigr\rfloor-4$. Then the following properties hold.

\medskip
{\bf (1)} Let $m\in\llbracket3,r\rrbracket$. Then the following properties hold.

\medskip\noindent
{\bf (1.a)} Suppose that $|K|\ge 3$. Then $\pi^{\S}_{n,m}(K)\le N_0^2N_1^2$, where 
$$N_i:=\min\Bigl(\l^{[\le \min(m-i-1,l)]}_K(m-i-1),\overline{\s}_{|K|}(l,m-i)\Bigr)$$ 
for $i\in\{0,1\}$. In particular, we have inequalities
\begin{equation*}
\begin{split}
\pi^{\S}_{n,m}(K)&\le \min\bigl((m-2)^2(m-1)^2,\overline{\s}_{|K|}(l,m)^4\bigr)\\
&\le\begin{cases} (m-2)^2(m-1)^2\;\quad\quad\quad\quad\quad\quad\quad\quad\quad\quad\quad\quad\quad\quad\quad\;\quad\, {\rm if}\; l=1\\
\min\bigl((m-2)(m-1), \bigl[(l-1)(|K|-1)+\bigl\lfloor\frac{|K|}{\sqrt{2}}\bigr\rfloor\bigr]^2\bigr)^2\;\quad\quad\;\, {\rm if}\; l\ge 2.\end{cases}
\end{split}
\end{equation*}

\noindent
{\bf (1.b)} Suppose that $|K|=2$ and $l\ge 2$. If $l=2$, then we also assume that $m=3$. Then we have inequalities
$$\pi^{\S}_{n,m}(K)\le \min\Bigl(\bigl\lfloor\frac{m-1}{2}\bigr\rfloor,l-1\Bigr)^2\min\Bigl(\bigl\lfloor\frac{m}{2}\bigr\rfloor,l-1\Bigr)^2\le (l-1)^4.$$

{\bf (2)} Suppose that $m\in\llbracket r+1,|K|^l\rrbracket$. Let $r_0\in\Bigl\{\bigl\lceil\sqrt{\frac{|K|^l}{2}}\bigr\rceil,\bigl\lceil\sqrt{\frac{|K|^l}{2}}\bigr\rceil+1\Bigr\}$ be the smallest such that $2\bigl\lceil\frac{|K|^l+4}{4}\bigr\rceil\le r_0^2-r_0$. Then 
$$\pi^{\S}_{n,m}(K)\le (m-r_0)^2(m-2)^2(m-1)^2\le (m-3)^2(m-2)^2(m-1)^2.$$
\end{proposition}

\begin{proof} Let $(\underline{P},\underline{Q})=\bigl((P_1,\ldots,P_m),(Q_1,\ldots,Q_m)\bigr)\in\mathbb D_{n,m}(K)^2$. We consider the subsets $Y_1:=\{P_1,\ldots,P_m\}$ and $Y_2:=\{Q_1,\ldots,Q_m\}$ of $K^n$.

Let $K\rightarrow L$ be a field extension of degree $l$. Let $\mathcal B$ be a $K$-basis of $L$; it allows us to identify $K^l=L$ and hence $K^n=K^l\times K^l=L\times L=L^2$. To emphasize that elements in $L=K^l$ are $l$-tuples of elements in $K$ we underline them. 

For part (1.a) it suffices to show that there exists $b\in\STGA_n(K)[N_0^2N_1^2]$ such that $b(\underline{P})=\underline{Q}$. 

Up to an $a_1\in\SL_2(L)$, based on Lemma \ref{L10}(1) and (2) applied to the triple $(n,l,K)=(2,1,L)$ we can assume that 
$$4\Bigl(\sum_{j=1}^{|L|} \frac{q_i(q_i-1)}{2}\Bigr)-4\le 4\frac{m(m-1)}{2|L|+2}-4\le 4\frac{r(r-1)}{2|L|+2}-4<|L|,$$ 
where $q_i:=|Z_i|\in \llbracket0,m\rrbracket$ for each $i\in \llbracket1,|L|\rrbracket$, with $Z_1,\ldots,Z_{|L|}$ as the fibers of the restriction to $Y_1$ of the projection $\pi_2:\mathbb A^2_L\rightarrow\mathbb A^1_L$ on the second coordinate. Let $Y_3:=\pi_2(Y_1)\subset L=K^l$.

If $|Y_3|=m$, then we define $Y_5:=Y_3$ and we note that $Y_1$ has tame type $(\d_{Y_5},1)$ with $\d_{Y_5}\in \llbracket1,\min(m-1,l)\rrbracket$ and hence also tame type $(\d_{Y_5},N_1)$. 

Next we assume that $|Y_3|\le m-1$ and we find a tame type of $Y_1$ as follows. 

Reindexing we can assume that $q_1\ge\cdots\ge q_{|L|}$. The proof of Theorem \ref{T8}(1.b) applied to the triple $(n,K,r)=(2,L,r)$ shows that there exists a rule of the form $(\underline{\beta}_1,\underline{\beta}_2)\mapsto \bigl(\underline{\beta}_1+f(\underline{\beta}_2), \underline{\beta}_2\bigr)$ that defines $a\in\SGA_2(L)$ such that $\{a(P_1),\ldots,a(P_m)\}$ is a subset of $L^2$ with distinct first coordinates in $L$. We can assume that we have $a\in\SGA_2(L)[m-2]$ by Theorem \ref{T5}(1) applied to $(n,K)=(2,L)$; hence we have $\ell\bigl(\R_{L/K}(a)\bigr)=\ell(a)\le m-2$ by Lemma \ref{L13}(3). 

Based on Theorem \ref{T5}(1) we can replace the components of $\R_{L/K}(a)$ by other polynomials in the variables $x_{l+1},\ldots,x_{n}$, giving us an automorphism in 
$$\STGA_n(K)\Bigl[\min\Bigl(\l^{[\d_{Y_3}]}_K(|Y_3|-1),\overline{\s}_{Y_3}\Bigr)\Bigr]\subset \STGA_n(K)[N_1]$$ 
that acts the same on $K^n$ as $\R_{L/K}(a)$. Thus $Y_1\subset K^{2n}$ has tame type $(\d_{Y_5},N_1)$, where $Y_5$ is the image of $a(Y)$ under the projection $\mathbb A^n_K\rightarrow\mathbb A^l_K$ on the first $l$ coordinates. We have $|Y_5|=m$, $\d_{Y_5}\in \llbracket1,\min(m-1,l)\rrbracket$, $\d_{Y_3}\in \llbracket1,\min(m-2,l)\rrbracket$, and $\overline{\s}_{Y_3}\le \overline{\s}_{|K|}(l,m-1)$ by Theorem \ref{T4}(3).

Regardless of what $|Y_3|$ is, $Y_1$ has tame type $(\d_{Y_5},N_1)$ and $|Y_5|\le m\le r\le |K|^l$. Thus $\overline{\s}_{Y_5}\le \overline{\s}_{|K|}(l,m)$ by Theorem \ref{T4}(3). From this and definitions we get that for $l=1$ we have $\overline{\s}_{Y_5}\le m-1$. As $\min\bigl(|K|^l,\bigl\lceil\frac{|K|^l}{\sqrt{2}}\bigr\rceil+1\bigr)\le r\le \min\bigl(|K|^l,\bigl\lceil\frac{|K|^l}{\sqrt{2}}\bigr\rceil+2\bigr)$ by Theorem \ref{T8}(1.a), for $(|K|,l)\notin\{(2,1),(3,1)\}$ we have the belonging relation $r\in\bigl\{\bigl\lceil\frac{|K|^l}{\sqrt{2}}\bigr\rceil+1,\bigl\lceil\frac{|K|^l}{\sqrt{2}}\bigr\rceil+2\bigr\}$ by Lemma \ref{L14} (1.b). Thus for $l\ge 2$, we get the inequality 
$\max(\overline{\s}_{Y_3},\overline{\s}_{Y_5})\le\overline{\s}_{|K|}(l,r)\le (l-1)(|K|-1)+\bigl\lfloor\frac{|K|}{\sqrt{2}}\bigr\rfloor$ by Lemma \ref{L14}(1.e) applied to $n_1=0$.

Based on the last two paragraphs and their analogs for $\underline{Q}$, the existence of $b$ for part (1.a) follows from Proposition \ref{PR17}(1). So part (1.a) holds.

Part (1.b) is proved entirely in the same way as part (1.a) except that, when we have $|Y_3|\le m-1$, based on Theorem \ref{T5}(2) we can replace the components of $\R_{L/K}(a)$ by other polynomials in the variables $x_{l+1},\ldots,x_{n}$, giving us an automorphism in 
$$\STGA_n(K)\left[\min\left(\l^{[\d_{Y_3}]}_K(|Y_3|-1),\overline{\s}_{Y_3}\right)\right]\subset \STGA_n(K)\Bigl[\min\Bigl(\Bigl\lfloor\frac{m-1}{2}\bigr\rfloor,l,\overline{\s}_{Y_3}\Bigr)\Bigr].$$
We only have to add that, as $|K|=2$, $l\ge 2$, and $(l,m)\neq (2,4)$, we have inequalities $\overline{\s}_{Y_3}\le l-1$ and $\overline{\s}_{Y_5}\le l-1$ by Lemmas \ref{L14}(1.f) and \ref{L6}(1); in particular, $\min(l,\overline{\s}_{Y_3})\le l-1$.

For part (2), we have $|L|=8$ or $|L|\ge 11$ by Example \ref{EX17.5}(1); thus $r_0\ge 3$. Let $b\in\SGA_2(L)[(m-r_0)^2(m-2)^2(m-1)^2]$ be such that $b(\underline{P})=\underline{Q}$ by Theorem \ref{T8}(1.c). For the automorphism $a:=\R_{L/K}(b)\in\SGA_n(K)$ we have $a(\underline{P})=\underline{Q}$ and $\ell(a)=\ell(b)\le (m-r_0)^2(m-2)^2(m-1)^2$ by Lemma \ref{L13}(3). So part (2) holds.\end{proof}

The way we proved Proposition \ref{PR21} allows us to generalize a weaker form of it to the following `relative' setting that is used in the case when $n$ is odd.

\begin{proposition}\label{PR22} Let $(l,n_1)\in\mathbb N^{\ast}\times\mathbb N$, $n:=2l$, and $K$ a finite field. We consider the largest $r\in\llbracket3,|K|^l\rrbracket$ such that $|K|^l>4\Bigl\lfloor\frac{\frac{r(r-1)}{2}}{|K|^l}\Bigr\rfloor-4$. Suppose that for each $O\in K^{n_1}$ we have an integer $m_O\in \llbracket1,|K|^l\rrbracket$ and two $m_O$-tuples $(P_{1,O},\ldots P_{m_O,O})$ and $(Q_{1,O},\ldots Q_{m_O,O})$ of distinct points in the subset $K^n\times \{O\}$ of $K^{n+n_1}$. Let $m:=\max(m_O|O\in K^{n_1})$ and $k:=|K|-1$. Then the following properties hold.

\medskip
{\bf (1)} Suppose that $m\in\llbracket3,r\rrbracket$. Then there exists $(f_1,\ldots,f_n)\in K[x_1,\ldots,x_{n+n_1}]^n$ such that the endomorphism $a:=\e(f_1,\ldots,f_n,x_{n+1},\ldots,x_{n+n_1})$ is in $\STGA_{n+n_1}(K)$ 
and the following properties hold.

\medskip\noindent
{\bf (1.a)} We have $a(P_{i,O})=Q_{i,O}$ for each $O\in K^{n_1}$ and every $i\in \llbracket1,m_O\rrbracket$.

\smallskip\noindent
{\bf (1.b)} If $l=1$, then $\ell(a)\le (n_1k+1)^2(m-2)^2(m-1)^2$.

\smallskip\noindent
{\bf (1.c)} If $l\ge 2$ and $|K|\ge 3$, then $\sqrt{\ell(a)}$ is bounded from above by
$$(n_1k+1)\min\Bigl(m-2,(l+n_1-1)k+\Bigl\lfloor\frac{|K|}{\sqrt{2}}\Bigr\rfloor\Bigr)\min\Bigl(m-1,(l+n_1-1)k+\Bigl\lfloor\frac{|K|}{\sqrt{2}}\Bigr\rfloor\Bigr).$$

\noindent
{\bf (1.d)} Suppose that $l\ge 2=|K|$. If $l=2$ we assume that $m=3$. Then 
$$\ell(a)\le (n_1+1)^2\Bigl[\min\bigl(\bigl\lfloor\frac{m-1}{2}\bigr\rfloor,l+n_1-1\bigr)\Bigr]^2\Bigl[\min\bigl(\bigl\lfloor\frac{m}{2}\bigr\rfloor,l+n_1-1\bigr)\Bigr]^2.$$

\medskip
{\bf (2)} Suppose that $m\in\llbracket r+1,|K|^l\rrbracket$. Let $r_0\in\Bigl\{\bigl\lceil\sqrt{\frac{|K|^l}{2}}\bigr\rceil,\bigl\lceil\sqrt{\frac{|K|^l}{2}}\bigr\rceil+1\Bigr\}$ be the smallest such that $2\bigl\lceil\frac{|K|^l+4}{4}\bigr\rceil\le r_0^2-r_0$. Then there exists $(f_1,\ldots,f_n)$ in $K[x_1,\ldots,x_{n+n_1}]^n$ such that the endomorphism $a:=\e(f_1,\ldots,f_n,x_{n+1},\ldots,x_{n+n_1})$ is in $\STGA_{n+n_1}(K)$,
we have an identity $a(P_{i,O})=Q_{i,O}$ for each $O\in K^{n_1}$ and every $i\in \llbracket1,m_O\rrbracket$, and we have inequalities
$$\ell(a)\le (n_1k+1)^2[(m-r_0)(m-2)(m-1)]^2\le (n_1k+1)^2[(m-3)(m-2)(m-1)]^2.$$
\end{proposition}

\begin{proof} Essentially, we only repeat the proof of Proposition \ref{PR21} fiberwise, with respect to the projection on the last $n_1$ `parameter' coordinates. 

All constants in our choice of the automorphisms are now polynomials in the parameter indeterminates of partial degree at most $|K|-1$ in every indeterminate $x_i$ with $i\in \llbracket n+1,n+n_1\rrbracket$. Then all the arguments go through, except for the upper bounds on the degree length of the automorphisms which require two adjustments; the fact that this is so even for part (2) follows from Remark \ref{R12.2}. 

First, for both parts, instead of $a_1$s in $\SL_2(L)$ in the beginning, we have automorphisms $a_1=\e(g_1,\ldots,g_n,x_{n+1},\ldots,x_{n+n_1})$ in $\STGA_{n+n_1}(K)$ with $\ell(a_1)\le n_1k+1$. 

Second, only for part (1), for the other automorphisms, if $l\ge 2$ then their degree lengths are bounded by $\min(m-1-i,\overline{s}_{|K|}(l+n_1,r|K|^{n_1}-i)$ with $i\in\{0,1\}$ and the strict capacities $\overline{s}_{|K|}(l+n_1,r|K|^{n_1})$ are computed in Lemma \ref{L14}(1.e) if $|K|\ge 3$ and Lemma \ref{L14}(1.f) if $|K|=2$ as $r\le\min(|K|^l,\lceil \frac{|K|^l}{\sqrt{2}}\rceil+2)$ by Theorem \ref{T8}(1.a).

With the field extension $K\rightarrow L$ as in the proof of Proposition \ref{PR21}, identifying $K^{n+n_1}=L^2\times K^{n_1}$ and using $\mathbb A^2_L=\Spec(L[\underline{x_1},\underline{x_2}])$, we detail on the construction of the two $a_1$s in the context of the $P_{i,O}$s points for part (1) only.

For $O\in K^{n_1}$ let $\underline{\gamma_O}\in L$ be such that the cardinality $N_O=N_O(\underline{\gamma_O})$ of the set 
$$\bigl\{(i,j)\in \llbracket1,m_O\rrbracket^2|i< j,\underline{x_1}+\gamma_{O}\underline{x_2}\; \textup{takes same value in}\; L\; \textup{at}\;P_{i,O}\; \textup{and}\;P_{j,O}\bigr\}$$ 
is as small as possible. Thus $4N_O-4\le 4\frac{m_O(m_O-1)}{2|K|^l}-4\le 4\frac{r(r-1)}{2|K|^l}-4<|L|$.\footnote{This is the place where we use that $r$ is defined using $\frac{\frac{r(r-1)}{2}}{|K|^l}$ and not the `usual' $\frac{\frac{r(r-1)}{2}}{|K|^l+1}$.}

We take an automorphism
$$a_1:=\e\bigl(\underline{x_1}+\underline{f_1,\ldots,f_l}\cdot\underline{x_2},\underline{x_2},x_{n+1},\ldots,x_{n+n_1}\bigr)\in\STGA_{n+n_1}(K)$$
with $(f_1,\ldots,f_l)\in K[x_{n+1},\ldots,x_{n+n_1}]^l$, $\underline{\bigl(f_1(O),\ldots,f_l(O)\bigr)}=\underline{\gamma_{O}}\in L=K^l$ for every $O\in K^{n_1}$, and $\deg(f_i)\le n_1k$ for each $i\in \llbracket1,n-1\rrbracket$ by Theorem \ref{T5}(1); therefore we have an inequality $\ell(a_1)\le n_1k+1$.
\end{proof}

Now we consider the case when $n$ is odd.

\begin{proposition}\label{PR23} Let $(l,m)\in (\mathbb N^{\ast})^2$ with $m\ge 3$ and let $n:=2l+1$. Let $K$ be a finite field and $k:=|K|-1$. Let $r\in\llbracket3,|K|^l\rrbracket$ be the largest such that $4\Bigl\lfloor\frac{\frac{r(r-1)}{2}}{|K|^l}\Bigr\rfloor-4<|K|^l$. If $l=1$ we assume that $|K|\ge 4$. Then the following properties hold.

\medskip
{\bf (1)} If $4\left(\frac{m(m-1)}{2|K|}-1\right)<|K|^{2l}$ (equivalently, $2m^2-2m<|K|^n+4|K|$) and $m\le |K|r$, then for $|K|\ge 3$ we have
\begin{equation}\label{EQ24}
\pi^{\S}_{n,m}(K)\le \Bigl[k|K|\bigl[kl+\bigl\lfloor\sqrt{|K|}\bigr\rfloor\bigr]\min\Bigl(m-2,kl+\Bigl\lfloor\frac{|K|}{\sqrt{2}}\Bigr\rfloor\Bigr)\min\Bigl(m-1,kl+\Bigl\lfloor\frac{|K|}{\sqrt{2}}\Bigr\rfloor\Bigr)\Bigr]^2
\end{equation}
and for $|K|=2$ we have $\pi_{n,m}(K)\le 4l^2\bigl[\min(\lfloor\frac{m-1}{2}\rfloor,l)\bigr]^2\bigl[\min(\lfloor\frac{m}{2}\rfloor,l)\bigr]^2\le 4l^6$.

\smallskip
{\bf (2)} If $m\le \lfloor \frac{|K|^{\frac{n}{2}}}{\sqrt{2}}\rfloor$ and $|K|\ge 3$, then we have
\begin{equation}\label{EQ25}
\pi^{\S}_{n,m}(K)\le \Bigl[k|K|\Bigl(kl+\Bigl\lfloor\sqrt{\frac{|K|}{2}}\Bigr\rfloor\Bigr)\min\Bigl(m-2,kl+\Bigl\lfloor\frac{|K|}{\sqrt{2}}\Bigr\rfloor\Bigr)\min\Bigl(m-1,kl+\Bigl\lfloor\frac{|K|}{\sqrt{2}}\Bigr\rfloor\Bigr)\Bigr]^2.
\end{equation}

{\bf (3)} Let $r_0\in\Bigl\{\bigl\lceil\sqrt{\frac{|K|^l}{2}}\bigr\rceil,\bigl\lceil\sqrt{\frac{|K|^l}{2}}\bigr\rceil+1\Bigr\}$ be the smallest such that we have an inequality $2\lceil\frac{|K|^l+4}{4}\rceil\le r_0^2-r_0$. If $4\left(\frac{m(m-1)}{2|K|}-1\right)<|K|^{2l}$ and $m\le |K|^{l+1}$, then we have
\begin{equation*}
\pi^{\S}_{n,m}(K)\le \bigl[k|K|(kl+\lfloor\sqrt{|K|}\rfloor)(m-r_0)(m-2)(m-1)\bigr]^2.
\end{equation*}
\end{proposition}

\begin{proof} If $m\le \lfloor \frac{|K|^{\frac{n}{2}}}{\sqrt{2}}\rfloor$, then we have $2m^2-2m<2m^2\le |K|^n<|K|^n+4|K|$ and $m\le \frac{|K|^{\frac{n}{2}}}{\sqrt{2}} =\frac{|K|^l}{\sqrt{2}}\cdot \sqrt{|K|} <\lceil\frac{|K|^l}{\sqrt{2}} \rceil |K|\le |K|r$ by Theorem \ref{T8}(1.a). So we prove parts (1) and (2) together.

We have $(2m-1)^2<2(|K|^n+4|K|)+1$. Thus $m<\frac{\sqrt{2|K|^n+8|K|+1}+1}{2}$. The inequality $\frac{\sqrt{2|K|^n+8|K|+1}+1}{2}\le |K|^{\frac{n}{2}}$ is equivalent to $8|K|\le 2|K|^n-4|K|^{\frac{n}{2}}$ and hence to $4\le |K|^{\frac{n}{2}-1}(|K|^{\frac{n}{2}}-2)$, which holds as $l\ge 2$ or $|K|\ge 4$. We conclude that $m<|K|^{\frac{n}{2}}<|K|^{l+1}$ and $\frac{m-1}{\sum_{i=0}^l |K|^i}<\sqrt{|K|}$. Therefore for $|K|\ge 3$ we have $\overline{\s}_{|K|}(2l,m)\le kl+\lfloor\sqrt{|K|}\rfloor$ by Definition \ref{D9}(1.b) and Lemma \ref{L6}(1). If $m\le \lfloor \frac{|K|^{\frac{n}{2}}}{\sqrt{2}}\rfloor$ and $|K|\ge 3$, then from Lemmas \ref{L14}(2.b) and \ref{L6}(1) we get that $\break\overline{\s}_{|K|}(2l,m)\le kl+\bigl\lfloor\sqrt{\frac{|K|}{2}}\bigr\rfloor$. If $|K|=2$, then $2m^2-2m<2^n+8$ is equivalent to $m^2-m<2^{2l}+4$ and (as $l\ge 2$) to $m\le 2^l$; so $\overline{\s}_{|K|}(2l,m)\le l$ by Lemma \ref{L6}(1).

Let $(\underline{P},\underline{Q})=\bigl((P_1,\ldots,P_m),(Q_1,\ldots,Q_m)\bigr)\in\mathbb D_{n,m}(K)^2$; we get two subsets $Y_1:=\{P_1,\ldots,P_m\}$ and $Y_2:=\{Q_1,\ldots,Q_m\}$ of $K^n$. It suffices to show that there exists $a\in\STGA_n(K)$ such that $a(\underline{P})=\underline{Q}$, $\ell(a)$ for $|K|\ge 3$ is less than or equal to the right-hand side of Inequality (\ref{EQ24}), and, in case $m\le \lfloor \frac{|K|^{\frac{n}{2}}}{\sqrt{2}}\rfloor$, even of Inequality (\ref{EQ25}), and $\ell(a)$ for $|K|=2$ is bounded by $4l^2\bigl[\min(\lfloor\frac{m-1}{2}\rfloor,l)\bigr]^2\bigl[\min(\lfloor\frac{m}{2}\rfloor,l)\bigr]^2$.

For $\underline{\alpha}=(\alpha_1,\ldots,\alpha_{2l})\in K^{2l}$ let $a_{\underline{\alpha}}:=\e\bigl(x_1,\ldots,x_{2l},x_{n}+\sum_{i=1}^{2l}\alpha_ix_i\bigr)\in\SL_n(K)$; for $s\in\{1,2\}$ we choose $\underline{\beta}_s\in K^{2l}$ that produces the smallest total number $N_s$ of subsets of cardinality $2$ of $a_{\underline{\beta}_s}(Y_s)$ whose images under the projection $\pi_1:\mathbb A^n_K\rightarrow\mathbb A^1_K$ on the last coordinate have cardinality $1$. As for distinct points $O_1$ and $O_2$ in $K^n$, the last coordinates of $a_{\underline{\alpha}}(O_1)$ and $a_{\underline{\alpha}}(O_2)$ are made equal for at most $|K|^{2l-1}$ choices of $\underline{\alpha}\in K^{2l}$, we have an inequality
$$N_s\le \Bigl\lfloor\frac{m(m-1)}{2|K|}\Bigr\rfloor.$$ Up to special linear automorphisms, to show the existence of $a$ we can assume that $\underline{\beta}_1=\underline{\beta}_2=(0,\ldots,0)$.

Similarly, for each $\underline{f}:=(f_1,\ldots,f_{2l})\in K[x_n]^{2l}$ with $\pi(\underline{f})\le k$ we consider the automorphism $a_{\underline{f}}:=\e\bigl(x_1+f_1(x_n),\ldots,x_{2l}+f_{2l}(x_n),x_n\bigr)\in\STGA_n(K)[k]$. We show that there exists $\underline{f}_s\in K[x_n]^{2l}$ such that $\bigl|\pi_{2l}\bigl(a_{\underline{f}_s}(Y_s)\bigr)\bigr|=m$, where $\pi_{2l}:\mathbb A^n_K\rightarrow\mathbb A^{2l}_K$ is the projection on the first $2l$ coordinates. We order the cardinalities of the intersection of $Y_s$ with the fibers of $\pi_1$ in decreasing order: $q_{1,s}\geq q_{2,s}\geq\cdots \geq q_{|K|,s}$. Thus we have identities $N_s=\sum_{i=1}^{|K|} \frac{q_{i,s}(q_{i,s}-1)}{2}$ and $\sum_{i=1}^{|K|} q_{i,s}=m$. We are shifting the fibers independently, starting from the larger ones, and it suffices to show for every $j\in \llbracket1,|K|\rrbracket$ with $q_{j,s}\ge 1$ we have an inequality $(\sum_{i=1}^{j-1}q_{i,s})q_{j,s} <|K|^{2l}$. If $q_j=1$, this follows from the inequalities $\sum_{i=1}^{j-1} q_{i,s}\le m-1<|K|^{\frac{n}{2}}<|K|^{2l}$. If $q_{j,s}\geq 2$, then
$$\Bigl(\sum_{i=1}^{j-1}q_{i,s}\Bigr)q_{j,s}\le 4\Bigl(\sum_{i=1}^{|K|}\frac{q_{i,s}(q_{i,s}-1)}{2}-1\Bigr)\le 4\Bigl(\frac{m(m-1)}{2|K|}-1\Bigr)<|K|^{2l},$$
as the strict inequality holds by our hypotheses. 

For $h\in K[x_1,\ldots,x_{2l}]$ of degree at most $\overline{\s}_{|K|}(2l,m)$, let
$$b_h:=\e\bigl(x_1,\ldots,x_{2l},x_n+h\bigr)\in\STGA_n(K)[\overline{\s}_{|K|}(2l,m)].$$ 
As $m<|K|^{l+1}< \sum_{i=0}^{l+1} |K|^i$ and $m\le |K|r$, from Theorem \ref{T4}(1) applied to $l+1$ and sets with $m$ points we get that there exists a pair $(h_1,h_2)\in K[x_1,\ldots,x_{2l}]^2$ for which $b_{h_1}$ and $b_{h_2}$ are defined and the following two properties hold: (i) for each $i\in \llbracket1,m\rrbracket$ the last coordinate of $b_{h_1}\bigl(a_{\underline{f}_1}(P_i)\bigr)$ coincides with the last coordinate of $b_{h_2}\bigl(a_{\underline{f}_2}(Q_i)\bigr)$, and (ii) for each $\alpha\in K$ there exist at most $r$ points $b_{h_1}\bigl(a_{\underline{f}_1}(P_i)\bigr)$ with the last coordinate $\alpha$.

If $l=1$, then $m-2<m-1\le r-1\le k\le kl+\bigl\lfloor\frac{|K|}{\sqrt{2}}\bigr\rfloor$. Based on this, for $|K|\ge 3$ and $l\ge 1$ (resp.\ $|K|=2$ and $l\ge 2$), let 
$$c\in\STGA_n(K)\Bigl[|K|^2\min\Bigl(m-2,kl+\Bigl\lfloor\frac{|K|}{\sqrt{2}}\Bigr\rfloor\Bigr)^2\min\Bigl(m-1,kl+\Bigl\lfloor\frac{|K|}{\sqrt{2}}\Bigr\rfloor\Bigr)^2\Bigr]$$ 
(resp.\ $c\in\STGA_n(K)[4[\min(\lfloor\frac{m-1}{2}\rfloor,l)]^2[\min(\lfloor\frac{m}{2}\rfloor,l)]^2]$) be such that we have an identity $(cb_{h_1}a_{\underline{f}_1})(\underline{P})=(b_{h_2}a_{\underline{f}_2})(\underline{Q})$ by the case $n_1=1$ in Proposition \ref{PR22}(1.b) if $l=1$ and in Proposition \ref{PR22}(1.c) (resp.\ Proposition \ref{PR22}(1.d)) if $l\ge 2$. 

For $a:=(a_{\underline{f}_2})^{-1}(b_{h_2})^{-1}cb_{h_1}a_{\underline{f}_1}\in\STGA_n(K)$ 
we have $a(\underline{P})=\underline{Q}$. As we have $\ell(a)\le \ell(c)\ell(a_{\underline{f}_1})\ell(a_{\underline{f}_1})\ell(b_{h_1})\ell(b_{h_2})$ by Inequality (\ref{EQ3}), the fact that $\ell(a)$ is bounded from above by the required values follows from the above upper bounds of $\ell(c)$, $\ell(a_{\underline{f}})$, and $\overline{\s}_{|K|}(2l,m)$. So parts (1) and (2) hold.

Part (3) is proved similarly to part (1), the only differences being that in (ii) we have to replace $r$ by $|K|^l$ and that for the existence of $c$ we have to quote the case $n_1=1$ of Proposition \ref{PR22}(2).
\end{proof}

We have the following consequence of Propositions \ref{PR21}(1.a) and \ref{PR23} and Theorems \ref{T7}(2) and \ref{T8}(1.b) and (2).

\begin{corollary}\label{C19}
Let $n\in\mathbb N^{\ast}\setminus\{1\}$, $K$ a finite field, and
$m\in\bigl\llbracket3,\lfloor \frac{|K|^{\frac{n}{2}}}{\sqrt{2}}\rfloor\bigr\rrbracket$. Then the following properties hold.

\medskip
{\bf (1)} If $|K|\geq 4$, then $\pi^{\S}_{n,m}(K)\le\frac{182}{49}(m-2)^6(m-1)^4$.

\smallskip
{\bf (2)} If $|K|=3$, then 
$$\pi^{\S}_{n,m}(K)\le \min\bigl(36(m-2)^4(m-1)^2,36(5m+1)(m-2)^2(m-1)^2\bigr).$$

{\bf (3)} If $|K|=2$, then 
$$\pi_{n,m}(K)\le\max\Bigl(\frac{1}{16}(m-3)^{6},(m-2)(m-1)^2,\frac{3^6m^3}{2^4}\Bigr)\le\frac{3^6}{2^4}m^3.$$
\end{corollary}

\begin{proof}
If $n=2$, then $3\le m\le\frac{|K|}{\sqrt{2}}$ implies that $|K|\ge 5$; therefore we have $\pi^{\S}_{n,m}(K)\le (m-2)^2(m-1)^2$ by Theorem \ref{T8}(1.b) and the corollary holds. Thus we can assume that $n\ge 3$.

If $m\le |K|$, then $\pi^{\S}_{n,m}(K)\le (m-2)(m-1)^2$ by Theorem \ref{T8}(2) for $m\ge 4$ and by Theorem \ref{T8}(1.b) for $m=3$ and the corollary holds. So we can assume that $|K|+1\le m\le \frac{|K|^{\frac{n}{2}}}{\sqrt{2}}$. 

If $\frac{m(m-1)}{2}<\sum_{i=0}^l |K|^i$, then $\pi^{\S}_{n,m}(K)\le (m-1)^2$ by Theorem \ref{T7}(1) and the corollary holds. So we can assume that $\sum_{i=0}^l |K|^i\le \frac{m(m-1)}{2}$.

If $n\ge m-1$, then $\pi^{\S}_{n,m}(K)\le (m-1)^2$ by Theorem \ref{T7}(2) and the corollary holds. So we can assume that $n\le m-2$; thus $m\ge n+2\ge 5$. 

As $n\ge 3$ and $m\le \frac{|K|^{\frac{n}{2}}}{\sqrt{2}}$, Propositions \ref{PR21}(1) and \ref{PR23}(1) and (2) apply. As the upper bounds for $\pi^{\S}_{n,m}(K)$ in Proposition \ref{PR23}(1) are weaker than the ones in Proposition \ref{PR21}(1), we can assume that $n=2l+1$ is odd and we get that Inequality (\ref{EQ25}) holds for $|K|\ge 3$ and that $\pi^{\S}_{n,m}(K)\le 4l^6$ for $|K|=2$. 

As the inequality $\sum_{i=0}^l |K|^i\le \frac{m(m-1)}{2}$ involves $l$ and so leads to inequalities that depend on $|K|$, we consider three disjoint cases as follows.

{\bf Case 1: $|K|\ge 4$.} So we have $5\le |K|+1\le m$, $n=2l+1\le m-2$, and $\sum_{i=0}^l |K|^i\le \frac{m(m-1)}{2}$. If $l=1$, then 
$$l(|K|-1)+\Bigl\lfloor\sqrt{\frac{|K|}{2}}\Bigr\rfloor=|K|-1+\Bigl\lfloor\sqrt{\frac{|K|}{2}}\Bigr\rfloor\le m-2+\Bigl\lfloor\sqrt{\frac{m-1}{2}}\Bigr\rfloor$$
$$< m-2+\frac{m-2}{2}=\frac{3(m-2)}{2}$$
and thus Inequality (\ref{EQ25}) gives that $\pi^{\S}_{n,m}(K)\le\frac{9}{4}(m-2)^6(m-1)^4$. 

Assume now that $l\ge 2$. As $8(\sum_{i=0}^l |K|^i)+1\le 4m(m-1)+1=(2m-1)^2$, for $l=2$ we get that $2\sqrt{2}|K|+1<2m-1$ and for $l\ge 3$ we get first that $2\sqrt{2}|K|^{\frac{l}{2}}+3<2m-1$ and second that $|K|^{\frac{l}{2}}<\frac{m-2}{\sqrt{2}}$; for the case $l\ge 3$ we used the following inequality $12\sqrt{2}|K|^{\frac{l}{2}}< 8|K|^{l-1}+8|K|^{l-2}$ which holds as for $|K|\ge 4$ we have $\frac{12\sqrt{2}}{8}<2.13<2.5\le\sqrt{|K|}+\frac{1}{\sqrt{|K|}}$.

For $l=2$ we have $|K|< \frac{m-1}{\sqrt{2}}$ and $|K|-1<\frac{m-2}{\sqrt{2}}$, hence
$$(|K|-1)|K|\Bigl[2(|K|-1)+\Bigl\lfloor\sqrt{\frac{|K|}{2}}\Bigr\rfloor\Bigr]\le  \frac{(m-2)(m-1)}{2}\Bigl[\sqrt{2}(m-2)+\frac{\sqrt{m-1}}{\sqrt[4]{8}}\Bigr]$$
$$<\frac{(m-2)(m-1)}{2}\Bigl[\sqrt{2}(m-2)+\frac{2\sqrt{2}}{7}(m-2)\Bigr]=\frac{9\sqrt{2}}{7}(m-2)^2(m-1)^.$$
Thus Inequality (\ref{EQ25}) gives that $\pi^{\S}_{n,m}(K)\le\frac{182}{49}(m-2)^6(m-1)^4$. 

For $l=3$ we have $|K|^{\frac{3}{2}}< \frac{m-2}{\sqrt{2}}$ and thus
$$(|K|-1)|K|\Bigl[3(|K|-1)+\Bigl\lfloor\sqrt{\frac{|K|}{2}}\Bigr\rfloor\Bigr]\le \bigl(3+\frac{\sqrt{2}}{2}\bigr)|K|^3<\frac{3+\frac{\sqrt{2}}{2}}{2}(m-2)^2.$$
Thus Inequality (\ref{EQ25}) gives that $\pi^{\S}_{n,m}(K)< 3.436 (m-2)^6(m-1)^4$. 

If $l\ge 4$, then $l|K|^3\le |K|^l$ for $|K|\ge 4$ and $(l+1)|K|^3\le |K|^l$ for $|K|\ge 5$. Hence 
$$(|K|-1)|K|\Bigl[l(|K|-1)+\Bigl\lfloor\sqrt{\frac{|K|}{2}}\Bigr\rfloor\Bigr]<|K|^l=\bigl(|K|^{\frac{l}{2}}\bigr)^2<\frac{(m-2)^2}{2}.$$
Thus Inequality (\ref{EQ25}) gives that $\pi^{\S}_{n,m}(K)\le\frac{1}{4}(m-2)^6(m-1)^4$.

As $\max\bigl(\frac{9}{4},\frac{182}{49},3.46,\frac{1}{4}\bigr)=\frac{182}{49}$, we get that part (1) holds.
 
{\bf Case 2: $|K|\ge 3$.} Then we have $l(3-1)+\bigl\lfloor\sqrt{\frac{3}{2}}\bigr\rfloor=2l+1=n\le m-2$. As $\sum_{i=0}^l 3^i\le \frac{m(m-1)}{2}$, we get that $4(3^{l+1}-1)+1\le (2m-1)^2$, equivalently that $\sqrt{4\cdot 3^{l+1}-3}\le 2m-1$. So $2\cdot 3^{\frac{l+1}{2}}-1<2m-1$ and thus $3^{\frac{l+1}{2}}<m$. 
Hence 
\begin{equation}\label{EQ25.3}
l\le \min\Bigl(\frac{m-3}{2},\lfloor2\log_3(m)\rfloor-1\Bigr).
\end{equation}
We check that $2l+1\le \sqrt{5m+1}$. As $m\ge 2l+3$, to check this we can assume that $l\ge 3$. If $l=3$, then $3^{\frac{3+1}{2}}=9<m$ and thus $2l+1=7<\sqrt{50}\le \sqrt{5m}$. If $l=4$, then $15.58<3^{\frac{3+1}{2}}<m$ and thus $2l+1=9=\sqrt{81}\le \sqrt{5m+1}$. If $l\ge 5$, then $4l^2+4l< 5\cdot 3^{\frac{l+1}{2}}<5m$ and therefore $2l+1=\sqrt{4l^2+4l+1}<\sqrt{5m}$. 

Thus $2l+1\le \min(m-2),\sqrt{5m+1})$. As the Inequality (\ref{EQ25}) is equivalent to $\pi^{\S}_{n,m}(K)\le 36(2l+1)^2(m-2)^2(m-1)^2$, it follows that part (2) holds.
 
{\bf Case 3: $|K|=2$.} Then $2l=n-1\le m-3$, so $l\le\frac{m-3}{2}$. As $\sum_{i=0}^l 2^i\le \frac{m(m-1)}{2}$, we get that $8(2^{l+1}-1)+1\le (2m-1)^2$, i.e., $\sqrt{2^{l+4}-7}\le 2m-1$. So $2^{\frac{l}{2}+2}-1<2m-1$ and thus $2^{\frac{l}{2}+1}<m$. Therefore 
\begin{equation}\label{EQ25.3}
l\le \min\Bigl(\frac{m-3}{2},\lfloor2\log_2(m)\rfloor-2\Bigr).
\end{equation}
We check that $l\le\frac{3\sqrt{m}}{2}$. As $m\ge 5$ and $\frac{3\sqrt{5}}{2}>3$, based on Equation (\ref{EQ25.3}) to check this we can assume that $\frac{m-3}{2}>3$, i.e., $m\ge 10$. For $m\ge 10$, as $\frac{3\sqrt{5}}{2}>4.5$, we can also similarly assume that $\frac{m-3}{2}>4.5$, i.e., $m\ge 13$, and thus that $l\ge 6$. But for $l\ge 6$ we have the inequalities $l\le 3\cdot 2^{\frac{l-2}{4}}=\frac{3}{2} 2^{\frac{l+2}{4}}<\frac{3}{2}\sqrt{m}$. From this and Equation (\ref{EQ25.3}) we get that the inequality $l\le\frac{3\sqrt{m}}{2}$ holds.

Therefore $\pi_{n,m}(K)\le 4l^6\le\min\bigl(\frac{1}{16}(m-3)^6,\frac{3^6}{2^4}m^3\bigr)$. So part (3) holds.\end{proof}

\begin{remark}\normalfont\label{R8}
While the bounds for the $\pi^{\S}_{n,m}(K)$s in Proposition \ref{PR23}(1) and (2) are weaker than in Proposition \ref{PR21}(1), they are still polynomial in $m$ for $m$ up to a suitable multiple of $|K|^{\frac{n}{2}}$ (as we can assume that $|K|< m$ by Theorem \ref{T8}(2)).
\end{remark}

\section{Non-polynomial upper bounds for medium $m$: principles}\label{S19}

In this section we obtain upper bounds for the $\pi_{n,m}(K)$s that pertain to numbers $m\in\bigl\llbracket2,2\lceil \frac{|K|^{n-1}}{4}\rceil\bigr\rrbracket$ (see Corollary \ref{C20}). 

Proposition \ref{PR21} is probably close to optimal for $n=2$. However, for $n>2$ we have $\frac{n}{2}<n-1$ and one can extend the range of $m$ for which we can get a polynomial, or close to polynomial, upper bound up to $2\bigl\lceil\frac{|K|^{n-1}}{4}\bigl\rceil$. This is done using the same ideas, but the construction is considerably more technical. To make the induction on $n$ work, we work in a `relative' setting, cf.\ Proposition \ref{PR22}.

Again we begin by computing the required values of the strict capacity functions. For $(k,n)\in (\mathbb N^{\ast}\setminus\{1\})^2$ we define
\begin{equation}\label{EQ26}
q_k(n):= \begin{cases} \frac{k+1}{2} \quad {\rm if}\; n=2\; \; \textup{and}\; k\equiv 1 \pmod{4}\\
\frac{k-2}{2} \quad \textup{if either}\; k\in 2+4\mathbb N\; \textup{with}\; n\ge 3 \; \textup{or}\; 4\mid k\\
\frac{k}{2} \quad \quad{\rm if}\; n=2\; \; \textup{and}\; k\equiv 2 \pmod{4}\\
\frac{k-1}{2} \quad \textup{in all other cases}.
 \end{cases}
\end{equation}
Note that we have identities
\begin{equation}\label{EQ27}
q_k(n)=q_k\bigl(\min(3,n)\bigr)
\end{equation}
and
\begin{equation}\label{EQ28}
q_k(n)=q_k(2)=\frac{k-1}{2} \quad \textup{if}\; k\equiv 3 \pmod{4}.
\end{equation}

\begin{lemma}\label{L15}
Let $(k,n,n_1)\in (\mathbb N^{\ast}\setminus\{1\})^2\times\mathbb N$. Let $m:=2k^{n_1}\lceil \frac{k^{n-1}}{4} \rceil\in\mathbb N^{\ast}$. If $k\ge 3$, let $q:=\bigl\lfloor\frac{m-1}{\sum_{i=0}^{n+n_1-2} k^i}\bigr\rfloor\in\mathbb N$. Let $l\in\mathbb N^{\ast}$ be such that $m\le k^l$. Then 
$$\overline{\s}_k(l,m)=(k-1)(n+n_1-2)+q_k(n).$$ More precisely, the following properties hold.

\medskip
{\bf (1)} If $k=2$ and $n=2$, then $m=2^{n_1+1}$ and $\overline{\s}_2(l,m)=n_1+1$.\footnote{If $k=2(2l+1)$ with $l\in\mathbb N^{\ast}$ and $n=2$, then $k^{n_1}<m=(\frac{k}{2}+1)k^{n_1}<k^{n_1+1}$, we have strict inequalities $\frac{k}{2}<\bigl\lfloor\frac{m-1}{\sum_{i=0}^{n_1} k^i}\bigr\rfloor<\frac{k}{2}+1$, and thus $\overline{\s}_k(l,m)=n_1(k-1)+\frac{k}{2}=n_1(k-1)+q_k(2)$. Such a $k$ is not the cardinal of a finite field.}

\smallskip
{\bf (2)} If $k=2$ and $n\ge 3$, then $m=2^{n+n_1-2}$ and $\overline{\s}_2(l,m)=n+n_1-2$. 

\smallskip
{\bf (3)} If $k\ge 3$, then we have inequalities $k^{n+n_1-2}<\frac{k^{n+n_1-1}}{2}\le m< k^{n+n_1-1}$ and $m-k^{n+n_1-2}\ge\lceil\frac{k}{2}\rceil$. 

\smallskip
{\bf (4)} If $k\ge 4$ and $4\mid k^{n-1}$ (i.e., if either $k\in 2+4\mathbb N^{\ast}$ and $n\ge 3$ or $4\mid k$), then $m=\frac{k^{n+n_1-1}}{2}$, $q=\lfloor\frac{(m-1)(k-1)}{2m-1}\rfloor=\frac{k-2}{2}$, and $\overline{\s}_k(l,m)=(n+n_1-2)(k-1)+\frac{k-2}{2}$. 

\smallskip
{\bf (5)} If $k=4s+r+1$ with $(s,r)\in\mathbb N\times\{0,2\}$, then for $n=2$ we have the identity $\overline{\s}_k(l,m)= n_1(k-1)+\frac{k+1-r}{2}$.

\smallskip
{\bf (6)} If $k\ge 3$ is odd and $n\ge 3$, then for each integer $l\ge n+n_1-1$ we have $\overline{\s}_k(l,m)=(n+n_1-2)(k-1)+\frac{k-1}{2}$.
\end{lemma}

\begin{proof}
For parts (1) and (2), the values of $m$ are clear and the identities for $\overline{\s}_2(l,m)$ follow from Lemma \ref{L6}(1). So parts (1) and (2) hold.

To prove part (3) we can assume that $n_1=0$. Thus $m\ge \frac{k^{n-1}}{2}>k^{n-2}$. Also, $m\le 2\frac{k^{n-1}+3}{4}=\frac{k^{n-1}+3}{2}\le k^{n-1}$. We have $\frac{k^{n-1}+3}{2}=k^{n-1}$ iff $(k,n)=(3,2)$. But for $(k,n)=(3,2)$ we have $m=2<k$. Therefore $m<k^{n-1}$. The inequalities $m-k^{n-2}\ge\frac{k^{n-2}(k-2)}{2}\ge\frac{k}{2}$ give that $m-k^{n-2}\ge\lceil\frac{k}{2}\rceil$. So part (3) holds.

For part (4), the first two identities are clear. The inequality $\frac{k-2}{2}\le \frac{(m-1)(k-1)}{2m-1}$ is equivalent to $(k-2)(2m-1)\le 2(m-1)(k-1)$ and thus to $k\le 2m$, which holds as $2m=k^{n+n_1-1}\ge k$. We have $\frac{k-2}{2}\in\mathbb N^{\ast}$ and $\frac{(m-1)(k-1)}{2m-1}<\frac{k-1}{2}$. From the last two sentences we get that $q=\frac{k-2}{2}$. Note that $m-q=\frac{k^{n+n_1-1}-k+2}{2}\ge k^{n+n_1-2}$ as $(k-2)k^{n+n_1-2}\ge k-2$. From the last two sentences, Definition \ref{D9}(1.b), and Lemma \ref{L6}(2), it follows that $\overline{\s}_k(l,m)=(n+n_1-2)(k-1)+\frac{k-2}{2}$. So part (4) holds.

For parts (5) and (6), we write $k^{n-1}-1=4s+r$ with $s\in\mathbb N$ and $r\in\{0,2\}$. We have $s=0$ iff $n=2$ and $k=3$. We compute $m-1=k^{n_1}(2s+2)-1$. To prove parts (5) and (6) we can assume that $l=n+n_1-1$ by Lemma \ref{L6}(3.a) and part (3). As $k\ge 3$ is odd, from part (3) we get that 
$$m\ge \frac{k^{n+n_1-1}+1}{2}=1+\frac{k^{n+n_1-1}-1}{2}\ge 1+\frac{k^{n+n_1-1}-1}{k-1}=1+\sum_{i=0}^{l-1} k^l.$$ Hence $\overline{\s}_k(n+n_1-1,m)=q+(n+n_1-2)(k-1)$ by Lemma \ref{L6}(3.b). So it suffices to show that $q=q_k(n)$.

For part (5), we have $\frac{k+3-r}{2}=2s+2$ and hence $q=\bigl\lfloor\frac{k^{n_1}(k+3-r)-2}{2\sum_{i=0}^{n_1} k^i}\bigr\rfloor$. 

Assume that $n_1=0$. Then $q=\frac{k+1-r}{2}=q_k(2)$.

Assume that $n_1>0$ and $r=0$. Then $k^{n_1}(k+3)\le (k+2)(\sum_{i=0}^{n_1} k^i)$ gives that $q\le \lfloor\frac{k^{n_1}(k+3)}{2\sum_{i=0}^{n_1} k^i}\rfloor\le \lfloor\frac{k+2}{2}\rfloor$. The inequality $k^{n_1}(k+3)-2\ge (k+1)(\sum_{i=0}^{n_1} k^i)$ is equivalent to $k^{n_1+1}+3k^{n_1}-2\ge k^{n_1+1}+2(\sum_{i=0}^{n_1} k^i)+1$ and therefore also to $k^{n_1}-1\ge 2(\sum_{i=0}^{n_1-1} k^i)=\frac{2(k^{n_1}-1)}{k-1}$, which holds as $k\ge 3$. From the last two sentences and $\frac{k+1}{2}\in\mathbb N$ we get that $q=\bigl\lfloor\frac{k^{n_1}(k+3)-2}{2\sum_{i=0}^{n_1} k^i}\bigr\rfloor=\frac{k+1}{2}=q_k(2)$.

Assume that $n_1>0$ and $r=2$. Clearly, we have $k^{n_1}(k+1)-2<k(\sum_{i=0}^{n_1} k^i)$. The inequality $k^{n_1}(k+1)-2> (k-1)(\sum_{i=0}^{n_1} k^i)$ is equivalent to the inequality $k^{n_1+1}+k^{n_1}-2> k^{n_1+1}-1$, which holds as $k^{n_1}\ge k>1$. From the last two sentences and $\frac{k-1}{2}\in\mathbb N$ we get that $q=\bigl\lfloor\frac{k^{n_1}(k+1)-2}{2\sum_{i=0}^{n_1} k^i}\bigr\rfloor=\frac{k-1}{2}=q_k(2)$.

We have $r=0$ iff $k\equiv 1 \pmod{4}$ and we have $r=2$ if $k\equiv 3 \pmod{4}$. From this and the last three paragraphs we get that part (5) holds.

For part (6), we have $\frac{k^{n-1}+3-r}{2}=2s+2$ and hence $q=\bigl\lfloor\frac{k^{n_1}(k^{n-1}+3-r)-2}{2\sum_{i=0}^{n+n_1-2} k^i}\bigr\rfloor$. The fact that $q=q_k(n)=\lfloor\frac{k-1}{2}\rfloor$ follows from the chain of inequalities 
$$k^{n_1}(k^{n-1}+3-r)-2< k^{n+n_1-1}+3k^{n_1}\le k^{n+n_1-1}+k^{n_1+1}< (k+1)\Bigl(\sum_{i=0}^{n+n_1-2} k^i\Bigr),$$
where $k\ge 3$ is used in the second (non-strict) inequality, and the inequality $\frac{k^{n_1}(k^{n-1}+3-r)-2}{2\sum_{i=0}^{n+n_1-2} k^i}>\frac{k-1}{2}$. So part (6) holds.
\end{proof}

We apply Lemma \ref{L15} to $k=|K|$ with $K$ a finite field.

\begin{proposition}\label{PR24} Let $(n,n_1)\in (\mathbb N^{\ast}\setminus\{1\})\times\mathbb N$. Suppose that for each point $O\in K^{n_1}$ we have $m_O$ distinct points $\{P_{1,O},\ldots, P_{m_O,O} \}$ in $K^n\times \{O\} \subset K^{n+n_1}, $ with $m_O\in \llbracket0,2\lceil \frac{|K|^{n-1}}{4} \rceil\rrbracket$. For each $i\in \llbracket0,n-2\rrbracket$ let $\varepsilon_i\in\{0,1\}$ be $1$ iff $n_1+i>0$ and there exists $O\in K^{n_1}$ such that $m_O\le |K|^i$. Let $k:=|K|-1$. Then there exists an $n$-tuple $(f_1,\ldots,f_n)\in K[x_1,\ldots,x_{n+n_1}]^n$ such that $a:=\e(f_1,\ldots,f_n,x_{n+1},\ldots,x_{n+n_1})$ is an automorphism in $\STGA_{n+n_1}(K)$ with the following properties.

\medskip
{\bf (1)} We have the inequality
$$\ell(a)\le [k(n+n_1-2)+q_{|K|}(2)][k(n+n_1-2)+q_{|K|}(3)]^{n-2}\prod_{j=0}^{n-2}[k(n_1+j)+1-\varepsilon_j].$$

\smallskip
{\bf (2)} For every $O\in K^{n_1}$ the projection function $\{a(P_{i,O})|i\in \llbracket1,m_O\rrbracket\}\rightarrow K^{n-1}$ to the second to the $n$-th coordinates is injective.
\end{proposition}

\begin{proof}
Let $m_n:=\min(m_O|O\in K^{n_1})\in\mathbb N$ and $M_n:=\max(m_O|O\in K^{n_1})\in\mathbb N$; we have inequalities $0\le m_n\le M_n\le 2\bigl\lceil \frac{|K|^{n-1}}{4}\bigr\rceil$. Let $M_{n,n_1}:=M_n|K|^{n_1}$. We can assume that $M_n\ge 1$. We have 
$$M_{n,n_1}\le 2\Bigl\lceil \frac{|K|^{n-1}}{4}\Bigr\rceil|K|^{n_1}\le |K|^{n+n_1-1}$$ 
by Lemma \ref{L15}(1) to (3). With the $q_{|K|}(n)$s as in Display (\ref{EQ26}), for each $l\in\mathbb N^{\ast}$ with $|K|^l\ge 2|K|^{n_1}\lceil \frac{|K|^{n-1}}{4}\rceil$ we have 
\begin{equation}\label{EQ29}
\overline{\s}_{|K|}(l,M_{n,n_1})\le \overline{\s}_{|K|}\Bigl(l,2|K|^{n_1}\Bigl\lceil \frac{|K|^{n-1}}{4}\Bigr\rceil\Bigr)=k(n+n_1-2)+q_{|K|}(n)
\end{equation} 
by the monotony part in Lemma \ref{L6}(1) and Lemma \ref{L15} applied to $k=|K|$.

We use induction on $n\ge 2$. For $n=2$, we have $M_2\le |K|$ and the base of the induction is proved by a slight modification of the proof of Proposition \ref{PR22}(1). 

First, for $O\in K^{n_1}$ let $\gamma_O\in K$ be such that the cardinality $N_O=N_O(\gamma_O)$ of the set 
$$\{(i,j)\in \llbracket1,m_O\rrbracket^2|i< j,x_1+\gamma_Ox_2\; \textup{takes same value at}\;P_{i,O}\; \textup{and}\;P_{j,O}\}$$
is as small as possible. Thus $N_O\le \frac{m_O(m_O-1)}{2|K|}$. Let 
$$a_1:=e\bigl(x_1+g_1(x_3,\ldots, x_{2+n_1})x_2,x_2,\ldots,x_{2+n_1}\bigr)\in \STGA_{2+n_1}(K)$$ 
be such that $g_1(O)=\gamma_O$ for each $O\in K^{n_1}$ with $m_O\ge 2$; therefore we have $\ell(a_1)\le kn_1+1-\varepsilon_0$ by Theorem \ref{T4}(1) and (2). 

Second, for each $O\in K^{n_1}$ we consider the function $\{a_1(P_{i,O})|i\in \llbracket1,m_O\rrbracket\}\rightarrow K$ induced by the projection on the first coordinate and order its fibers in decreasing order of cardinality: $q_{1,O}\geq q_{2,O}\geq \cdots \geq q_{|K|,O}$. We have $\sum_{i=1}^{|K|} q_{i,O}=m_O\le |K|$. We are looking for an automorphism of the form $$a_2:=\e\bigl(x_1,x_2+g_2(x_1,x_3\ldots,x_{2+n_1}),x_3,\ldots,x_{2+n_1}\bigr)$$ 
to get distinct second coordinates for each set $\bigl\{a_2\bigl(a_1(P_{i,O})\bigr)|i\in \llbracket1,m_O\rrbracket\bigr\}$. It exists if we can show that for each $O\in K^{n_1}$, for every $s\in \llbracket1,|K|\rrbracket$ we have an inequality $(\sum_{i=1}^{s-1}q_{i,O})q_{s,O}<|K|$. 

If $q_{s,O}=0$, then the inequality clearly holds. If $q_{s,O}=1$, then the inequality holds as $\sum_{i=1}^{s-1}q_{i,O}\le m_O-1\le M_2-1<|K|$. If $q_{s,O}\geq 2,$ then 
$$\Bigl(\sum_{i=1}^{s-1}q_{i,O}\Bigr)q_{s,O}\le 4\sum_{i=1}^{s-1}\frac{q_{i,O}(q_{i,O}-1)}{2} \le 4\Bigl(\frac{m_O(m_O-1)}{2|K|}-1\Bigr)$$
and $|K|\geq 3$; so $m_O\le \frac{|K|}{2}+2$ and hence $4\bigl(\frac{m_O(m_O-1)}{2|K|}-1\bigr)$ is at most 
$$4\left(\frac{(\frac{|K|}{2}+2)(\frac{|K|}{2}+1)}{2|K|}-1\right)=\frac{|K|}{2}-1+\frac{4}{|K|}<|K|.$$

Therefore $g_2$ exists. We can assume that $\deg(g_2)\le\overline{\s}(n_1+1,M_{2,n_1})$ by Theorem \ref{T4}(1) and (3). Thus we can assume that $\ell(a_2)\le kn_1+q_k(2)$ by Display (\ref{EQ29}).

As $a:=a_2a_1$ satisfies all requirements, the base of the induction holds.

For $n\geq 3$, assuming that the statement holds for elements in $\llbracket2,n-1\rrbracket$, we prove it for $n$. As above, we start by considering for each $O\in K^{n_1}$ an $(n-1)$-tuple $\underline{\gamma_O}=(\gamma_{1,O},\ldots,\gamma_{n-1,O})\in K^{n-1}$ so that the cardinality $N_O=N_O(\underline{\gamma_O})$ of the set 
$$\bigl\{(i,j)\in \llbracket1,m_O\rrbracket^2|i< j,x_r+\gamma_{r,O}x_n\; \textup{has same value at}\;P_{i,O}\; \textup{and}\;P_{j,O}\;\forall r\in \llbracket1,n-1\rrbracket\bigr\}$$ 
is as small as possible. Thus $N_O\le\frac{m_O(m_O-1)}{2|K|^{n-1}}$.
We take an automorphism
$$a_1=\e(x_1+f_1x_n,\ldots, x_{n-1}+f_{n-1}x_n, x_n,x_{n+1},\ldots,x_{n+n_1})\in\STGA_{n+n_1}(K)$$
such that for each $i\in\llbracket1,n-1\rrbracket$ we have $f_i\in K[x_{n+1},\ldots,x_{n+n_1}]$ with $f_i(O)=\gamma_{i,O}$ for every $O\in K^{n_1}$ for which we have $m_O\ge 2$; as above we argue that we have $\ell(a_1)\le kn_1+1-\varepsilon_0$.

We search for an automorphism
$$a_2=\e\bigl(x_1,\ldots,x_{n-1},x_n+g_2(x_1,\ldots,x_{n-1},x_{n+1},\ldots,x_{n+n_1}),x_{n+1},\ldots,x_{n+n_1}\bigr)$$ 
in $\STGA_{n+n_1}(K)$ such that for each $O\in K^{n_1}$ and every $\alpha\in K$, for the set $Y_{\alpha,O}$ of elements $i\in \llbracket1,m_O\rrbracket$ such that the last $n_1+1$ coordinates of $a_2\bigl(a_1(P_{i,O})\bigr)$ form the point $(\alpha,O)\in K^{1+n_1}$ we have $|Y_{\alpha,O}|\le 2\lceil \frac{|K|^{n-2}}{4}\rceil$. This will allow us to apply the induction hypothesis in order to get the result for $n$.

For this to work, similar to the base case, we arrange the fibers of the projection $\{a_1(P_{i,O})|i\in \llbracket1,m_O\rrbracket\}\rightarrow K^{n-1}$ to the first $n-1$ coordinates in decreasing order of cardinality: $q_{1,O}\geq q_{2,O}\geq \cdots \geq q_{|K|^{n-1},O}$. We have $\sum_{i=1}^{|K|^{n-1}} q_{i,O}=m_O$. As above, $g_2$ exists if we have $\frac{(\sum_{i=1}^{s-1}q_{i,O})q_{s,O}}{|K|}<2\lceil\frac{|K|^{n-2}}{4}\rceil$ for every $s\in \llbracket1,|K|^{n-1}\rrbracket$.

If $q_{s,O}=0$ the inequality clearly holds and if $q_{s,O}=1$ the inequality follows from the inequalities $\sum_{i=1}^{s-1}q_{i,O}\le m_O-1< 2\lceil \frac{|K|^{n-1}}{4}\rceil \le |K|\cdot 2\lceil \frac{|K|^{n-2}}{4}\rceil$. If $q_{s,O}\geq 2$, it suffices to prove the inequality $(\sum_{i=1}^{s-1}q_{i,O})q_{s,O}<\frac{|K|^{n-1}}{2}$. As in the proof of the base of the induction we argue that
$$\Bigl(\sum_{i=1}^{s-1}q_{i,O}\Bigr)q_{s,O}\le 4\left(\frac{m_O(m_O-1)}{2|K|^{n-1}}-1\right)=\frac{2m_O(m_O-1)}{|K|^{n-1}}-4$$
and thus it suffices to prove the inequality 
\begin{equation}\label{EQ30}
\frac{2m_O(m_O-1)}{|K|^{n-1}}-4<\frac{|K|^{n-1}}{2}.
\end{equation}
If $|K|$ is even, then $2m_O\le|K|^{n-1}$ and thus $\frac{2m_O(m_O-1)}{|K|^{n-1}}-4\le\frac{|K|^{n-1}}{2}-5<\frac{|K|^{n-1}}{2}$. If $|K|$ is odd, then $m_O\le \frac{|K|^{n-1}}{2}+2,$ so
$$\frac{2m_O(m_O-1)}{|K|^{n-1}}-4\le \frac{2(\frac{1}{2}|K|^{n-1}+2)(\frac{1}{2}|K|^{n-1}+1)}{|K|^{n-1}}-4=\frac{1}{2}|K|^{n-1}-1+\frac{4}{|K|^{n-1}}.$$
As $|K|\geq 3$ and $n\geq 3$, this is less than $\frac{1}{2}|K|^{n-1}$. Hence Inequality (\ref{EQ30}) holds.

Therefore $g_2$ exists. We can assume that $\deg(g_2)\le\overline{\s}(n+n_1-1,M_{n,n_1})$ by Theorem \ref{T4}(1) and (3). Thus we can assume that $\ell(a_2)\le k(n+n_1-2)+q_k(n)$ by Display (\ref{EQ29}); so $\ell(a_2)\le k(n+n_1-2)+q_k(3)$ by Equation (\ref{EQ28}).

So $a:=a_2a_1$ satisfies $\ell(a)\le [kn_1+1-\varepsilon_0][k(n+n_1-2)+q_k(3)]$. Moreover, $a$ puts the points in the position to apply the induction hypothesis for a new $(n,n_1)$ being $(n-1,n_1+1)$, where for each pair $(\alpha,O)\in K\times K^{n_1}=K^{n_1+1}$ we have $m_{\alpha,O}:=|Y_{\alpha,O}|\le 2\bigl\lceil \frac{|K|^{n-2}}{4}\bigr\rceil$ and the set of $m_{\alpha,O}$ points in $K^n$ being the subset of $\{a_2\bigl(a_1(P_{i,O})\bigr)|i\in \llbracket1,m_O\rrbracket\}$ formed by points whose last $n_1+1$ coordinates are $(\alpha,O)$. Thus we define $m_{n-1}:=\min\bigl(m_{\alpha,O}|(\alpha,O)\in K\times K^{n_1}\bigr)\in\bigl\llbracket0,2\bigl\lceil \frac{|K|^{n-2}}{4}\bigr\rceil\bigr\rrbracket$ and $M_{n-1}:=\max\bigl(m_{\alpha,O}|(\alpha,O)\in K\times K^{n_1}\bigr)\in\bigl\llbracket0,2\bigl\lceil\frac{|K|^{n-2}}{4}\bigr\rceil\bigr\rrbracket$. Moreover,
$$M_{n-1,n_1+1}:=M_{n-1}|K|^{n_1+1}\le 2|K|^{n_1+1}\Bigl\lceil \frac{|K|^{n-2}}{4}\Bigr\rceil.$$ If there exists $l\in \llbracket1,n-2\rrbracket$ such that $m_n\le |K|^l$, then for each $O\in K^{n_1}$ with $m_O\le |K|^l$ we can find $\alpha\in K$ such that $m_{\alpha,O}\le |K|^{l-1}$ and hence $m_{n-1}\le |K|^{l-1}$; this allows us to apply the inductive hypotheses even for the $\varepsilon$s. In particular, if such an $l$ exists, we have $\varepsilon_l=\cdots=\varepsilon_{n-2}=1$.\footnote{If we have $\frac{(\sum_{i=1}^{s-1}q_{i,O})q_{s,O}}{k}<2\lceil\frac{|K|^{n-2}}{4}\rceil$ for every $s\in \llbracket1,|K|^{n-1}\rrbracket$, then we can assume that there exists $\alpha\in K$ with $m_{\alpha,O}=0$ and thus we can assume that $\epsilon_1=\cdots=\epsilon_{n-2}=1$. Each improvement (via an assumption or argument) of the inequality $M_{n-j,n_1+j}\le 2|K|^{n_1+j}\lceil \frac{|K|^{n-1-j}}{4} \rceil$ could be used to improve the factor $k(n+n_1-2)+q_{|K|}(n-j)$ in the upper bound for $\ell(a)$.}

This completes the induction and the proof of the proposition.
\end{proof}

Proposition \ref{PR24} has several applications, starting with the following one.

\begin{theorem}\label{T9} Suppose that $(n,n_1)\in (\mathbb N^{\ast}\setminus\{1\})\times \mathbb N$, $K$ is a finite field, and for each $O\in K^{n_1}$ we have two $m_O$-tuples $(P_{1,O},\ldots,P_{m,O})$ and $(Q_{1,O},\ldots,Q_{m_O,O})$ of distinct points in $K^n\times \{O\} \subset K^{n+n_1}$, with $m_O\in \llbracket1,2\lceil \frac{|K|^{n-1}}{4}\rceil\rrbracket$. Let $k:=|K|-1$. For every $i\in \llbracket0,n_1-2\rrbracket$ let $\varepsilon_i\in\{0,1\}$ be $1$ iff $n_1+i>0$ and there exists $O\in K^{n_1}$ such that $m_O\le |K|^i$. Then there exists an $n$-tuple $(f_1,\ldots,f_n)\in K[x_1,\ldots,x_{n+n_1}]^n$ such that $a:=\e(f_1,\ldots,f_n,x_{n+1},\ldots,x_{n+n_1})$ is an automorphism in $\STGA_{n+n_1}(K)$ for which the following properties hold.

\medskip
{\bf (1)} The degree length $\ell(a)$ is less than or equal to
$$[k(n+n_1-2)+q_{|K|}(2)]^{2n+1}[k(n+n_1-2)+q_{|K|}(3)]^{n^2-n-2}\prod_{j=0}^{n-2}[k(n_1+j)+1-\varepsilon_j]^{2(j+1)}.$$

\smallskip
{\bf (2)} For each $O\in K^{n_1}$ and every $i\in \llbracket1,m_O\rrbracket$ we have $a(P_{i,O})=Q_{i,O}$. 
\end{theorem}

\begin{proof}
Let $\star\in\{P,Q\}$. We use induction on $n\ge 2$. 

For the base of the induction we have $n=2$ and in this case Proposition \ref{PR24} gives an automorphism $a_{1,\star}\in\STGA_{2+n_1}(K)\bigl[\bigl(kn_1+1-\varepsilon_0\bigr)\bigl(n_1k+q_{|K|}(2)\bigr)\bigr]$ with the property that for each $O\in K^{n_1}$, the second coordinates of the $m_O$ points in the set $\{a_{1,\star}(\star_{i,O})|i\in \llbracket1,m_O\rrbracket\}$ are distinct. 

Then we construct 
$$a_{2,\star}=\e\bigl(x_1+g_{1,\star}(x_2,\ldots,x_{2+n_1}),x_2,\ldots,x_{2+n_1}\bigr)\in\STGA_{2+n_1}(K)[n_1k+q_{|K|}(2)]$$ 
such that for every $O\in K^{n_1}$ the following properties hold: (i) for each $i\in \llbracket1,m_O\rrbracket$, the first coordinate of $a_{2,Q}\bigl(a_{1,Q}(Q_{i,O})\bigr)$ equals that of $a_{2,P}\bigl(a_{1,P}(P_{i,O})\bigr)$, and (ii) the  $m_O$ points in the set $\bigl\{a_{2,P}\bigl(a_{1,P}(P_{i,O})\bigr)|i\in \llbracket1,m_O\rrbracket\bigr\}$ have distinct first coordinates. 

\medskip\noindent
This is possible as by Theorem \ref{T4}(3) and Lemma \ref{L15}(1), (4), and (5) we can assume that 
$$\deg(g_{1,\star})\le\overline{\s}_{|K|}\bigl(n_1+1,\sum_{O\in K^{n_1}} m_O\bigr)\le \overline{\s}_{|K|}\Bigl(n_1+1,2|K|^{n_1}\Bigl\lceil \frac{|K|}{4}\Bigr\rceil\Bigr)=n_1k+q_{|K|}(2)$$
by the same argument as in the proof of Proposition \ref{PR24}.
 
Next we construct $a_3=\e\bigl(x_1,x_2+g_2(x_1,x_3,\ldots,x_{2+n_1}),x_3,\ldots,x_{2+n_1}\bigr)$ also in $\STGA_n(K)[n_1k+q_{|K|}(2)]$ such that for each point $O\in K^{n_1}$, for every $i\in \llbracket1,m_O\rrbracket$ the second coordinates of $(a_3a_{2,P}a_{1,P})(P_{i,O})$ and $(a_{2,Q}a_{1,Q})(Q_{i,O})$ are equal. 

As the automorphism $a:=a_{1,Q}^{-1}a_{2,Q}^{-1}a_3a_{2,P}a_{1,P}$ has the required properties, the base of the induction holds.

With $n\ge 3$, we want to prove the existence of $a$ assuming the statement holds for elements in the set $\llbracket2,n-1\rrbracket$. First, Proposition \ref{PR24} gives for each $\star\in\{P,Q\}$ an automorphism $a_{1,\star}\in\STGA_n(K)$ of degree length less than or equal to
$$[k(n+n_1-2)+q_{|K|}(2)][k(n+n_1-2)+q_{|K|}(3)]^{n-2}\prod_{j=0}^{n-2}[k(n_1+j)+1-\varepsilon_j]$$
such that for each $O\in K^{n_1}$, the projection 
$$\{a_{1,\star}(\star_{i,O})|i\in \llbracket1,m_O\rrbracket\}\rightarrow K^{n-1}$$
on the second to the $n$-th coordinates is injective. 

Then we find an automorphism $$a_{2,\star}=\e\bigl(x_1+h_{1,\star}(x_2,\ldots,x_{n+n_1}),x_2,\ldots,x_{n+n_1}\bigr)\in\STGA_n(K)$$ 
such that $\ell(a_{2,\star})\le k(n+n_1-2)+q_{|K|}(n)=k(n+n_1-2)+q_{|K|}(3)$, for each $O\in K^{n_1}$ the first coordinate of $a_{2,P}\bigl(a_{1,P}(P_{i,O})\bigr)$ equals that of $a_{2,Q}\bigl(a_{1,Q}(Q_{i,O})\bigr)$ for every $i\in \llbracket1,m_O\rrbracket$ and moreover, for every $\alpha\in K$ the number of $i\in \llbracket1,m_O\rrbracket$ such that $a_{2,P}\bigl(a_{1,P}(P_{i,O})\bigr)$ has the first coordinate $\alpha$ is at most $\lceil \frac{|K|^{n-2}}{4}\rceil$. This is possible as, again by Theorem \ref{T4}(3), we can assume that 
$$\deg(h_{1,\star})\le\overline{\s}_{|K|}\bigl(n+n_1-1,\sum_{O\in K^{n_1}} m_O\bigr)\le \overline{\s}_{|K|}\Bigl(n+n_1-1,2|K|^{n_1}\Bigl\lceil \frac{|K|^{n-1}}{4}\Bigr\rceil\Bigr),$$
hence $\deg(h_1)\le k(n+n_2-2)+q_{|K|}(n)=k(n+n_2-2)+q_{|K|}(3)$ by Lemma \ref{L15} and Equation (\ref{EQ27}).

Then we switch mentally (without actually introducing any extra automorphisms) the roles of $x_1$ and $x_n$ and apply the induction hypothesis, for the new $(n,n_1)$ being $(n-1,n_1+1)$, which gives us an automorphism of the form
$$a_3=\e(x_1,h_2,\ldots,h_n,x_{n+1},\ldots,x_{n+n_1})\in\STGA_{n+n_1}(K)$$ 
with $(h_2,\ldots,h_n)\in K[x_1,\ldots,x_{n+n_1}]^{n-1}$ and such that for each $O\in K^{n_1}$ and every $i\in\llbracket1,m_O\rrbracket$ we have $a_3\bigl(a_{2,P}(a_{1,P}(P_{i,O}))\bigr)=a_{2,Q}\bigl(a_{1,Q}(Q_{i,O})\bigr)$.
By taking $a:=a_{1,Q}^{-1}a_{2,Q}^{-1}a_3a_{2,P}a_{1,P}$,we get that part (2) holds.

Inequality (\ref{EQ3}) gives $\ell(a)\le \ell(a_{1,P})\ell(a_{2,P})\ell(a_{1,Q})\ell(a_{2,Q})\ell(a_3)$. Hence we have
$$\frac{\ell(a)}{\ell(a_3)}\le [k(n+n_2-2)+q_{|K|}(2)]^2[k(n+n_1-2)+q_{|K|}(3)]^{2n-2}\prod_{j=0}^{n-2}[k(n_1+j)+1-\varepsilon_j]^2.$$
As $[k(n+n_1-2)+q_{|K|}(2)]^{2n-1}[k(n+n_1-2)+q_{|K|}(3)]^{n^2-3n}$ times the product
$\prod_{j=0}^{n-3} [k(n_1+1+j)+1-\varepsilon_{j+1}]^{2(j+1)}=\prod_{j=1}^{n-2} [k(n_1+j)+1-\varepsilon_j]^{2j}$
are greater than or equal to $\ell(a_3)$ by the inductive assumption, part (1) also holds. So the inductive step holds. This completes the induction and the proof of the theorem.\end{proof}

We can easily deduce from the above theorem an estimate for the $\pi^{\S}_{n,m}(K)$s with $m\in \llbracket1,2\lceil \frac{|K|^{n-1}}{4}\rceil\rrbracket$ by setting $n_1=0$ and squaring the bound. But we can do somewhat better, by using Proposition \ref{PR24} more directly, as follows.

\begin{theorem}\label{T10} Suppose that $(n,n_1)\in (\mathbb N^{\ast}\setminus\{1\})\times \mathbb N$, $K$ is a finite field, and for each $O\in K^{n_1}$ we have two $m_O$-tuples $(P_{1,O},\ldots,P_{m,O})$ and $(Q_{1,O},\ldots,Q_{m_O,O})$ of distinct points in $K^n\times \{O\} \subset K^{n+n_1}$, with $m_O\in \llbracket1,2\lceil\frac{|K|^{n-1}}{4}\rceil\rrbracket$. Let $k:=|K|-1$. For each $i\in \llbracket0,n_1-2\rrbracket$ let $\varepsilon_i\in\{0,1\}$ be $1$ iff $n_1+i>0$ and there exists $O\in K^{n_1}$ such that $m_O\le |K|^i$. Then there exists an $n$-tuple $(f_1,\ldots,f_n)\in K[x_1,\ldots,x_{n+n_1}]^n$ such that $b:=\e(f_1,\ldots,f_n,x_{n+1},\ldots,x_{n+n_1})$ is in $\STGA_{n+n_1}(K)$ and the following properties hold.

\medskip
{\bf (1)} The degree length $\ell(b)$ is less than or equal to 
$$[k(n+n_1-2)+q_{|K|}(2)]^{2n+3}[k(n+n_1-2)+q_{|K|}(3)]^{n^2+n-6}\prod_{j=0}^{n-2}[k(n_1+j)+1-\varepsilon_j]^{2(j+2)}.$$

{\bf (2)} The identity $b(P_{i,O})=Q_{i,O}$ holds for each $O\in K^{n_1}$ and every $i\in \llbracket1,m_O\rrbracket$.\end{theorem}

\begin{proof}
By Proposition \ref{PR24} applied for $\star\in\{P,Q\}$ to the set $\{\star_{i,O}|i\in \llbracket1,m_O\rrbracket\}$ of $m_O$ distinct points indexed by $O\in K^{n_1}$, there exists $a_{1,\star}\in \STGA_{n+n_1}(K)$ with 
$$\ell(a_1)\le [k(n+n_1-2)+q_{|K|}(2)][k(n+n_1-2)+q_{|K|}(3)]^{n-2}\prod_{j=0}^{n-2}[k(n_1+j)+1-\varepsilon_j].$$
such that the image of the projection $\{a_{1,\star}(\star_{i,O})|i\in \llbracket1,m_O\rrbracket\}\rightarrow K^{n-1}$ on the second to the $n$-th coordinates has cardinality $m_O$ for each $O\in K^{n_1}$.

Let $a\in\STGA_{n+n_1}(K)$ be as in Theorem \ref{T9} and such that for every point $O\in K^{n_1}$ and $i\in \llbracket1,m_O\rrbracket$ it sends $a_{1,P}(P_{i,O})$ to $a_{1,Q}(Q_{i,O})$. 

For $b:=a_{1,Q}^{-1}aa_{1,P}\in \STGA_{n+n_1}(K)$ we have $b(P_{i,O})=Q_{i,O}$ for every $O\in K^{n_1}$ and each $i\in \llbracket1,m_O\rrbracket$. So part (2) holds.

Inequality (\ref{EQ3}) gives $\ell(b)\le\ell(a_{1,P})\ell(a_{1,Q}) \ell(a)$. Based on this and the upper bounds for $\ell(a_{1,\star})$ and $\ell(a)$, it follows that $\ell(b)$ satisfies the required inequality. So part (1) holds.
\end{proof}

\section{Non-polynomial upper bounds for medium $m$: estimates}\label{S20}

For applying Theorem \ref{T10} to $n_1=0$ we first prove the following inequalities.

\begin{lemma}\label{L16} For $(n,k)\in (\mathbb N^{\ast}\setminus\{1,2\})\times\mathbb N^{\ast}$ let 
$$C_{n,k}:=[k(n-2)+q_{k+1}(2)]^{2n+3}[k(n-2)+q_{k+1}(3)]^{n^2+n-6} \prod \limits_{i=1}^{n-2}(ik+1)^{2(i+2)}.$$
Let $$\varepsilon_{n,k,2}:=\frac{2\min\bigl(n(2-k)+4\ln(n-2)+\frac{k}{3},n(3-k)-\frac{5-k}{3}\bigr)}{k}$$
and $\overline{C}^+_{n,k}:=\overline{C}^+_{n,k,1}\overline{C}^+_{n,k,2}$, where
$$\overline{C}^+_{n,k,1}:=k^{2n^2+4n-9}(n-1)^{2n^2+n-2}e^{\frac{\varepsilon_{n,k,2}}{2}}$$
and
\begin{equation}\label{EQ31}
\overline{C}^+_{n,k,2}:=\begin{cases} e^{-\frac{n^2+n-6}{n-1}-\frac{n^2+n-6}{2(n-1)^2}} \,\;\;\quad\quad\quad\quad\quad\quad\quad\quad\quad\quad\quad\quad\quad\quad\;\;\,\;\;\, {\rm if}\;\; k=1\\
e^{-\frac{(k+1)(n^2+3n-3)}{2k(n-1)}-\frac{(k+1)^2(n^2+3n-3)}{8k^2(n-1)^2}} \;\;\;\quad\quad\quad\quad\quad\quad\quad\quad\quad\,\;\; {\rm if}\;\;4\mid k+1\\
e^{-\frac{(k-1)(2n+3)}{2k(n-1)}-\frac{(k+1)(n^2+n-6)}{2k(n-1)}-\frac{(k-1)^2(2n+3)}{8k^2(n-1)^2}-\frac{(k+1)^2(n^2+n-6)}{8k^2(n-1)^2}} \;\quad {\rm if}\;\; k\in 5+4\mathbb N\\
e^{-\frac{n^2+3n-3}{2(n-1)}-\frac{n^2+3n-3}{8(n-1)^2}}
\quad\quad\quad\quad\quad\quad\quad\quad\quad\quad\quad\quad\quad\quad\quad\;\, {\rm if}\;\; 4\mid k-2\\
e^{-\frac{(k-2)(2n+3)}{2k(n-1)}-\frac{n^2+n-6}{2(n-1)}-\frac{(k-2)^2(2n+3)}{8k^2(n-1)^2}-\frac{n^2+n-6}{8(n-1)^2}} \,\;\quad\quad\quad\quad\quad\;\;\;{\rm if}\;\; 4\mid k.
\end{cases}
\end{equation}
Let $C^+_{n,k}:=C^+_{n,k,1}C^+_{n,k,2}$, where
$$C^+_{n,k,1}:=k^{2n^2+4n-9}(n-2)^{n^2+3n-3}(n-1)^{n^2-2n+1}e^{\frac{\varepsilon_{n,k,2}}{2}}$$
and
\begin{equation*}
C^+_{n,k,2}:=\begin{cases} e^{\frac{2n+3}{n-2}} \;\;\quad\quad\quad\quad\quad\quad\quad\quad\;\;\; {\rm if}\;\;\: k=1\\
e^{\frac{(k-1)(n^2+3n-3)}{2k(n-2)}} \;\;\quad\quad\quad\quad\quad\; {\rm if}\;\;4\mid k+1\\
e^{\frac{(k+1)(2n+3)}{2k(n-2)}+\frac{(k-1)(n^2+n-6)}{2k(n-2)}} \;\quad\, {\rm if}\;\;k\in 5+4\mathbb N\\\
e^{\frac{n^2+3n-3}{2(n-2)}} \quad\quad\quad\quad\quad\quad\quad\;\;\; {\rm if}\;\; 4\mid k-2\\
e^{\frac{(k+2)(2n+3)}{2k(n-2)}+\frac{n^2+n-6}{2(n-2)}} \;\quad\quad\quad\;\;{\rm if}\;\; 4\mid k.
\end{cases}
\end{equation*}
Then the following properties hold.

\medskip
{\bf (1)} We have inequalities $C_{n,k}< \overline{C}^+_{n,k}$ and $C_{n,k}< C^+_{n,k}$.

\smallskip
{\bf (2)} If $k=1$, then we have an inequality 
$$C_{n,1}<(n-2)^{n^2+n-6}(n-1)^{n^2+4}e^{\min\bigl(n+4\ln(n-2)+\frac{1}{3},2n-\frac{4}{3}\bigr)}.$$

{\bf (3)} If $4\mid k-2$, then we have an inequality
$$C_{n,k}<k^{2n^2+4n-9}\Bigl(n-\frac{3}{2}\Bigr)^{n^2+3n-3}(n-1)^{n^2-2n+1}e^{\frac{\varepsilon_{n,k,2}}{2}}.$$

\smallskip
{\bf (4)} If $4\mid k$, then we have an inequality
$$C_{n,k}<k^{2n^2+4n-9}(n-2)^{2n+3}\Bigl(n-\frac{3}{2}\Bigr)^{n^2+n-6}(n-1)^{n^2-2n+1}e^{\frac{(k+2)(2n+3)}{2k(n-2)}+\frac{\varepsilon_{n,k,2}}{2}}.$$
\end{lemma}

\begin{proof}
For every $i\in \llbracket1,n-2\rrbracket$ we have $1+ik=ik (1+\frac{1}{ik})$. The inequality $(1+\frac{1}{x})^{xy}< e^y$ for $(x,y)\in (0,\infty)\times (0,\infty)$ applied to $(x,y)=\bigl(ik,\frac{2(i+2)}{ik}\bigr)$ gives 
$$\left[1+\frac{1}{ik}\right]^{2(i+2)}< e^{\frac{2(i+2)}{ik}}= e^{\frac{2+\frac{4}{i}}{k}}.$$
By induction on $n\in (\mathbb N^{\ast}\setminus\{1,2\})$ we get that $2n+4\sum_{i=2}^{n-2}\frac{1}{i}\le 3n-\frac{5}{3}$ and that the inequality is strict iff $n\neq 5$. Also, $2n+4\sum_{i=2}^{n-2}\frac{1}{i}\le 2n+4\int_1^{n-2} \frac{dt}{t}=2n+4\ln(n-2)$ and the inequality is strict if $n\ge 4$. We conclude that
$$ \prod \limits_{i=1}^{n-2} \left[1+\frac{1}{ik}\right]^{2(i+2)}< e^{\frac{2n+4\sum_{i=2}^{n-2}\frac{1}{i}}{k}}\le e^{\frac{\min\bigl(2n+4\ln(n-2),3n-\frac{5}{3}\bigr)}{k}}.$$
Moreover, $ \prod \limits_{i=1}^{n-2} i^{2(i+2)}$ equals $1$ if $n=3$ and equals $e^{\sum \limits_{i=2}^{n-2} 2(i+2)\ln i}$ if $n\geq 4$. We check that for each integer $n\ge 3$ we have 
\begin{equation}\label{EQ31.5}
\sum \limits_{i=2}^{n-2} 2(i+2)\ln i <(n-1)^2\ln (n-1)-n+\frac{1}{3}.
\end{equation}
We estimate
\begin{equation*} 
\begin{aligned}
&\sum \limits_{i=2}^{n-2} 2(i+2)\ln i \le \int \limits_{2}^{n-1} (2x+4)\ln(x)dx=\Bigl(x^2\ln(x)-\frac{x^2}{2}+4[x\ln(x)-x]\Bigr)\Bigl|_2^{n-1}\\
&=(n-1)^2\ln(n-1)-\frac{(n-1)^2}{2}+4(n-1)\ln(n-1)-4(n-1)-12\ln(2)+10\\
&=(n-1)^2\ln(n-1)+\frac{n-1}{2}[8\ln(n-1)-8-n+1]-12\ln(2)+10\\
&<(n-1)^2\ln (n-1)-n+\frac{1}{3}+\frac{n-1}{2}[8\ln(n-1)-6-n+1]+2.349.
\end{aligned}
\end{equation*}
As for $n\ge 18$ we have $\frac{n-1}{2}[8\ln(n-1)-6-n+1]+2.349<0$, Inequality (\ref{EQ31.5}) holds for $n\ge 18$. For $n\in\llbracket3,17\rrbracket$, one checks via direct computations that Inequality (\ref{EQ31.5}) holds. So  Inequality (\ref{EQ31.5}) holds for $n\ge 3$.

Thus, for all $n\ge 3$ we have an inequality
$$ \prod \limits_{i=1}^{n-2} i^{2(i+2)}< e^{(n-1)^2\ln (n-1)-(n-\frac{1}{3})}=(n-1)^{n^2-2n+1}e^{-n+\frac{1}{3}}.$$ 
As $\prod \limits_{i=1}^{n-2} k^{2(i+2)}=k^{n^2+n-6}$, we conclude that 
\begin{equation*}
\prod \limits_{i=1}^{n-2}(ik+1)^{2(i+2)}< k^{n^2+n-6}(n-1)^{n^2-2n+1} e^{\frac{\varepsilon_{n,k,2}}{2}}.
\end{equation*}
\phantomsection{Thus, for $\overline{C}_{n,k,1}:=[k(n-1)]^{n^2+3n-3}\prod \limits_{i=1}^{n-2}(ik+1)^{2(i+2)}$, we have
$\overline{C}_{n,k,1}<\overline{C}_{n,k,1}^+$. Similarly, for $C_{n,k,1}:=[k(n-2)]^{n^2+3n-3}\prod \limits_{i=1}^{n-2}(ik+1)^{2(i+2)}$, we have
$C_{n,k,1}<C_{n,k,1}^+$.}\label{PH94} 

Let
$$\overline{C}_{n,k,2}:=\Bigl[1-\frac{k-q_{k+1}(2)}{k(n-1)}\Bigr]^{2n+3}\Bigl[1-\frac{k-q_{k+1}(3)}{k(n-1)}\Bigr]^{n^2+n-6}.$$ 
Based on the formulas in Display (\ref{EQ26}) applied to $|K|=k+1$, we have
\begin{equation*}
\overline{C}_{n,k,2}:= \begin{cases} \bigl(1-\frac{1}{n-1}\bigr)^{n^2+n-6} \,\;\quad\quad\quad\quad\quad\quad\quad\quad\quad\quad {\rm if}\;\; k=1\\
\bigl(1-\frac{k+1}{2k(n-1)}\bigr)^{n^2+3n-3} \,\;\quad\quad\quad\quad\quad\quad\quad\quad\, {\rm if}\;\;4\mid k+1\\
\bigl(1-\frac{k-1}{2k(n-1)}\bigr)^{2n+3}\bigl(1-\frac{k+1}{2k(n-1)}\bigr)^{n^2+n-6} \,\quad {\rm if}\;\; k\ge 5 \; \;\textup{and}\;\;4\mid k+3\\
\bigl(1-\frac{1}{2(n-1)}\bigr)^{n^2+3n-3} \;\;\;\quad\quad\quad\quad \quad\quad\quad\quad\;{\rm if}\;\; 4\mid k-2\\
\bigl(1-\frac{k-2}{2k(n-1)}\bigr)^{2n+3}\bigl(1-\frac{1}{2(n-1)}\bigr)^{n^2+n-6} \,\;\;\quad {\rm if}\;\; 4\mid k.
 \end{cases}
\end{equation*}

For $x>1$, we have 
\begin{equation}\label{EQ32}
\Bigl(1-\frac{1}{x}\Bigr)^x=e^{x\ln(1-\frac{1}{x})}<e^{-1-\frac{1}{2x}-\frac{1}{3x^2}}<e^{-1-\frac{1}{2x}}
\end{equation}
based on the Taylor series $x\ln(1-\frac{1}{x})=-\sum_{i=1}^{\infty}\frac{1}{ix^{i-1}}$. Based on this we get that 
$\overline{C}_{n,k,2}<\overline{C}^+_{n,k,2}$. As $C_{n,k}=\overline{C}_{n,k,1}\overline{C}_{n,k,2}$, the first inequality of part (1) holds. 

Let
$$C_{n,k,1}:=\Bigl[1+\frac{q_{k+1}(2)}{k(n-2)}\Bigr]^{2n+3}\Bigl[1+\frac{q_{k+1}(3)}{k(n-2)}\Bigr]^{n^2+n-6}.$$ 
Based on the formulas in Display (\ref{EQ26}) applied to $|K|=k+1$, we have
\begin{equation*}
C_{n,k,2}:= \begin{cases} \bigl(1+\frac{1}{n-2}\bigr)^{2n+3}\,\;\;\quad\quad\quad\quad\quad\quad\quad\quad\quad\quad\quad {\rm if}\; k=1\\
\bigl(1+\frac{k-1}{2k(n-2)}\bigr)^{n^2+3n-3} \,\,\;\quad\quad\quad\quad\quad\quad\quad\quad\, {\rm if}\;\; 4\mid k+1\\
\bigl(1+\frac{k+1}{2k(n-2)}\bigr)^{2n+3}\bigl(1+\frac{k-1}{2k(n-2)}\bigr)^{n^2+n-6} \,\;\quad {\rm if}\;\; k\ge 5 \; \;\textup{and}\;\;4\mid k+3\\
\bigl(1+\frac{1}{2(n-2)}\bigr)^{n^2+3n-3} \,\,\quad\quad\quad\quad \quad\quad \quad\quad\quad {\rm if}\;\; 4\mid k-2\\
\bigl(1+\frac{k+2}{2k(n-2)}\bigr)^{2n+3}\bigl(1+\frac{1}{2(n-2)}\bigr)^{n^2+n-6} \,\;\;\;\quad{\rm if}\;\; 4\mid k.
 \end{cases}
\end{equation*}
The inequality $(1+\frac{1}{x})^{xy}< e^y$ for $(x,y)\in (0,\infty)\times (0,\infty)$ gives $C_{n,k,2}<C^+_{n,k,2}$. As $C_{n,k}=C_{n,k,1}C_{n,k,2}$, the second inequality of part (1) holds. So part (1) holds.

Parts (2) and (3) follow from the inequality $C_{n,k}<C_{n,k,1}^+C_{n,k,2}$. 

Part (4) follows from the inequality $C_{n,k,2}<e^{\frac{(k+2)(2n+3)}{2k(n-2)}}\bigl(1-\frac{1}{2(n-2)}\bigr)^{n^2+n-6}$.\end{proof}

Based on Lemma \ref{L16}(1) applied to $k=|K|-1$ and $n\ge 3$, we restate the case $n_1=0$ of Theorem \ref{T10} as follows. 

\begin{corollary}\label{C20}
Let $n\in\mathbb N^{\ast}\setminus\{1,2\}$ and $K$ a finite field. With $k:=|K|-1$, let $C_{n,k}$, $\overline{C}^+_{n,k}$, and $C^+_{n,k}$ be as in Lemma \ref{L16}. Let $\overline{D}^+_{n,k}:=\lfloor\overline{C}^+_{n,K}\rfloor\in\mathbb N^{\ast}$ and $D^+_{n,k}:=\lfloor C^+_{n,k}\rfloor\in\mathbb N^{\ast}$. Let $D^{-}_{n,k}:=\lfloor C_{n,k}\rfloor$ and $D_{n,k}:=\min(\overline{D}^+_{n,k},D^+_{n,k})$. If we have $m\in \llbracket2,2\lceil \frac{|K|^{n-1}}{4}\rceil\rrbracket$, then $\pi^{\S}_{n,m}(K)\le D^{-}_{n,k}\le D_{n,k}$.
\end{corollary}

\section{Non-polynomial upper bounds for all permutations}\label{S21}

In this section we apply Sections \ref{S13} and \ref{S20} in order to obtain, via transpositions for $K=\mathbb F_2$ and via $s$-cycles with $s$ odd in the set $\llbracket3,2p-1\rrbracket$ for all finite fields, upper bounds for all the values of $\ell_{n,K}$.

From Corollary \ref{C20} and Proposition \ref{PR20}(3) and (4) we get directly the following consequence.

\begin{corollary}\label{C21}
Suppose that $n\ge 4$. Let $\sigma\in\perm(\mathbb F_2^n)$. Then the following properties hold.

\medskip
{\bf (1)} Suppose that $\sigma$ is even. Then we have inequalities
$$\ell_{n,\mathbb F_2}(\sigma)\le 2^{2\lceil\frac{\nu_{2,1}(\sigma)}{2^{n-2}}\rceil} (D_{n,1}^-)^{4\lceil\frac{\nu_{2,1}(\sigma)}{2^{n-2}}\rceil}\le 2^{2\lceil\frac{\nu_{2,1}(\sigma)}{2^{n-2}}\rceil} D_{n,1}^{4\lceil\frac{\nu_{2,1}(\sigma)}{2^{n-2}}\rceil}$$
and $\ell_{n,\mathbb F_2}(\sigma)\le 64 (D_{n,1}^-)^{13}\le 64D_{n,1}^{13}$.

{\bf (2)} Suppose that $\sigma$ is odd. Then we have inequalities
$$\ell_{n,\mathbb F_2}(\sigma)\le 2^{\lceil\frac{\nu_{2,1}(\sigma)+1}{2^{n-2}}\rceil+\lceil\frac{\nu_{2,1}(\sigma)-1}{2^{n-2}}\rceil-1} (n-1)(D_{n,1}^-)^{2\lceil\frac{\nu_{2,1}(\sigma)+1}{2^{n-2}}\rceil+2\lceil\frac{\nu_{2,1}(\sigma)-1}{2^{n-2}}\rceil}$$
and $\ell_{n,\mathbb F_2}(\sigma)\le 32(n-1)(D_{n,1}^-)^{13}\le 32(n-1)D_{n,1}^{13}$.\end{corollary}

For applications, the combination of Corollary \ref{C20} and Proposition \ref{PR2} require the following lemma.

\begin{lemma}\label{L17}
Let $n\in\mathbb N^{\ast}\setminus\{1,2\}$ and $K$ a finite field be such that $|K|\ge 3$ or $n\ge 5$.\footnote{If one defines $\digamma_{1,l}$ and $\wp_{1,l}$ by the same formulas (see Definition \ref{D3}(1) and (2)), then for $(n,|K|)=(4,2)$ we have $l=4$, $r=1$, $|K|^n-l=12$, $\digamma_{2,24}=4r=4$, and $\wp_{r,|K|^n-l}=\wp_{1,12}=6$.} Let $l:=2\lceil \frac{|K|^{n-1}}{4}\rceil$ and $r:=\lfloor \frac{l}{3}\rfloor$. Then the following properties hold.

\medskip
{\bf (1)} If $|K|^{n-2}\ge 24$, then we have 
$$\wp_{r,|K|^n-l}=\Bigl\lfloor\frac{|K|^n-l-\digamma_{r,|K|^n-l}}{2r}\Bigr\rfloor+2\le\Bigl\lfloor\frac{|K|^n-l}{2r}\Bigr\rfloor+2\le 3|K|.$$

{\bf (2)} If $(n,|K|)\neq (3,4)$, then $\wp_{r,|K|^n-l}=\bigl\lfloor\frac{|K|^n-l-\digamma_{r,|K|^n-l}}{2r}\bigr\rfloor+2\le 3|K|$.

\smallskip
{\bf (3)} If $|K|$ is odd, then $\lfloor\frac{|K|^n-l-4}{2r}\rfloor+2\le 3|K|$.
\end{lemma}

\begin{proof}
Let $\iota\in\{0,1,2,3\}$ be such that $l=\frac{|K|^{n-1}+\iota}{2}$ and let $\varepsilon\in\{0,1,2\}$ be such that $r=\frac{l-\varepsilon}{3}$. Defining $\mho:=\iota-2\varepsilon\in \llbracket-4,3\rrbracket$, we have $r=\frac{|K|^{n-1}+\mho}{6}$. Thus for each $\epsilon\in\llbracket0,\digamma_{r,|K|^n-l}\rrbracket$ we have that $\Bigl\lfloor\frac{|K|^n-l-\epsilon}{2r}\Bigr\rfloor+2$ is equal to 
\begin{equation}\label{EQ33}
\Bigl\lfloor\frac{|K|^n-\frac{|K|^{n-1}+\iota+2\epsilon}{2}}{\frac{|K|^{n-1}+\mho}{3}}\Bigr\rfloor+2=\Bigl\lfloor\frac{3|K|^n+\frac{|K|^{n-1}-3\iota-6\epsilon+4\mho}{2}}{|K|^{n-1}+\mho}\Bigr\rfloor=3|K|+1+\Bigl\lfloor\frac{\frac{-|K|^{n-1}}{2}-\hslash}{|K|^{n-1}+\mho}\Bigr\rfloor
\end{equation}
with $\hslash:= \frac{3}{2}\iota+3\epsilon+(3|K|-1)\mho$. 

Based on Equation (\ref{EQ33}) and $|K|^{n-1}+\mho\ge 8+\mho>0$, to prove part (1) it suffices to show that $0< |K|^{n-1}+2\hslash$, equivalently, that $-\hslash<\frac{|K|^{n-1}}{2}$. To check this we can assume that $\mho<0$; hence $\varepsilon\in\{1,2\}$. We have $-\hslash\le -\mho(3|K|-1)$. Therefore, if $|K|^{n-2}\ge -6\mho$, then $\frac{|K|^{n-1}}{2}\ge -3\mho|K|>-\hslash$, and it follows that $\wp_{r,|K|^n-l}\le 3|K|$. Thus part (1) holds.

We have $r\ge 3$ iff $|K|^{n-1}\ge 17$. Thus, if $|K|^{n-2}\ge 24$, we have inequalities $|K|^{n-1}\ge 24|K|\ge 48>17$ and hence $r\ge 3$; from this and Proposition \ref{PR3}(3) we get that $\digamma_{r,|K|^n-l}\ge 8>4$. Based on this and the proof of part (1), to prove parts (2) and (3) we can assume that $|K|^{n-2}\le 24$, $\mho<0$, $\varepsilon\in\{1,2\}$, and $|K|^{n-1}+2\hslash\le 0$.

If $n=3$ and $|K|\in\llbracket5,23\rrbracket$ is a prime, then $|K|^2\equiv 1\pmod{12}$, hence $\iota=3$ and $\varepsilon=1$; so $\mho=1$ and this case is not possible. Similarly, the case $n=3=|K|$ is not possible as $\varepsilon=0$. 

Hence either $n=4$ and $|K|=3$ or $2\mid |K|$. If $n=4$ and $|K|=3$, then $l=14$, $\iota=0$, $r=4$, and $3^4-14-4=63$; so $\lfloor\frac{|K|^n-l-4}{2r}\rfloor+2=9$. So part (3) holds.

Assume now that $2\mid |K|$; so $\iota=0$. As $(n,|K|)\notin\{(3,2),(4,2),(3,4)\}$ and $|K|^{n-2}\le 24$, we have $(n,|K|)\in\{(3,16),(3,8),(4,4),(5,2),(6,2)\}$. 

If $(n,|K|)=(3,16)$, then $l=128$, $r=42$, $|K|^n-l=3968$, $\digamma_{42,3968}=4r=168$ by Proposition \ref{PR3}(2), and $\wp_{42,3968}=\lfloor\frac{3968}{84}\rfloor=46<48$. 

If $(n,|K|)=(3,8)$, then $l=32$, $r=10$, $|K|^n-l=480$, $\digamma_{10,480}=4r=40$ by Proposition \ref{PR3}(2), and $\wp_{10,480}=\lfloor\frac{480}{20}\rfloor=24$. 

If $(n,|K|)=(4,4)$, then $l=32$, $r=10$, $|K|^n-l=224$, $\digamma_{10,224}=4r=40$ by Proposition \ref{PR3}(2), and $\wp_{10,224}=\lfloor\frac{224}{20}\rfloor=11<12$. 

If  $(n,|K|)=(5,2)$, then $l=8$, $r=2$, $|K|^n-l=24$, $\digamma_{2,24}=4r=8$ by Proposition \ref{PR3}(2), and $\wp_{2,24}=\lfloor\frac{24}{4}\rfloor=6$. 

If  $(n,|K|)=(6,2)$, then $l=16$, $r=5$, $|K|^n-l=48$, $\digamma_{5,48}=4r=20$ by Proposition \ref{PR3}(2), and $\wp_{5,48}=\lfloor\frac{48}{10}\rfloor=4<6$. 

So part (2) holds.\end{proof}

Similar to Corollary \ref{C17}, we have the following direct consequence of Corollary \ref{C20} and prior results. 

\begin{corollary}\label{C22}
Let $n\in\mathbb N^{\ast}\setminus\{1,2\}$, $l:=\lfloor\frac{n}{2}\rfloor$, $K$ a finite field, $p:=\char(K)$, and $k:=|K|-1$. Let $s\in \llbracket3,2p-1\rrbracket$ be odd. Let $r\in\mathbb N^{\ast}$. If $s\ge 5$, let $\epsilon\in\{2,3\}$ be such that $\epsilon=3$ iff $s=5$. Let $\sigma\in\Alt(K^n)$. There exists $a\in\STGA_n(K)[\lfloor N_{\sigma}\rfloor]$ such that $\sigma=a(K)$, where $N_{\sigma}$ is defined on cases as follows with $D_{n,k}$ as in Corollary \ref{C20}, $E_{n-2,k,\frac{s-3}{2}}$ as in Equation (\ref{EQ12}), and $E_{n-2,k}$ as in Equation (\ref{EQ13}).

\medskip
{\bf (1)} Suppose that $s=3$ and $\bigl\lfloor \frac{2\lceil \frac{|K|^{n-1}}{4}\rceil}{3}\bigr\rfloor\ge 2$.\footnote{Part (1) can be adapted to each $s\ge 5$ and parts (1) to (3) can be adapted to $r=1$.} We take $r:=\bigl\lfloor \frac{2\lceil \frac{|K|^{n-1}}{4}\rceil}{3}\bigr\rfloor$. Then 
$$N_{\sigma}:=\min\Bigl((E_{n-2,k}D^2_{n,k})^{\nu_{3,r}(\sigma)},D_{n,k}(E_{n-2,k}D_{n,k}^2)^{\wp_{r,|K|^n-2\lceil\frac{|K|^{n-1}}{4}\rceil}}\Bigr).$$ In particular, we have $N_{\sigma}\le (E_{n-2,k}D_{n,k}^2)^{\wp_{r,\n(\sigma)}}$.

\smallskip
{\bf (2)} Suppose that $s=3$. Then we take $N_{\sigma}$ as follows.

\medskip
{\bf (2.a)} Suppose that $|K|\ge 3$. If $n=2l$ (resp.\ $n=2l+1$) and $2\le r\le \Bigl\lfloor \frac{\lceil \frac{|K|^l}{\sqrt{2}}\rceil}{3}\Bigr\rfloor$, then $N_{\sigma}:=\min\Bigl(\bigl(E_{n-2,k}k^8l^8\bigr)^{\nu_{3,r}(\sigma)},k^4l^4\bigl(E_{n-2,k}k^8l^8\bigr)^{\wp_{r,|K|^n-\lceil \frac{|K|^l}{\sqrt{2}}\rceil}}\Bigr)$
(resp.\ $N_{\sigma}$ is the minimum of the two numbers $[E_{n-2,k}k^4(k+1)^4(kl+\lfloor\sqrt{\frac{k+1}{2}}\rfloor)^4(kl+\lfloor\frac{k+1}{\sqrt{2}}\rfloor)^8]^{\nu_{3,r}(\sigma)}$ and $E_{n-2,k}^{-\frac{1}{2}}[E_{n-2,k}k^4(k+1)^4(kl+\lfloor\sqrt{\frac{k+1}{2}}\rfloor)^4(kl+\lfloor\frac{k+1}{\sqrt{2}}\rfloor)^8]^{\frac{1}{2}+\wp_{r,|K|^n-\lceil \frac{|K|^l}{\sqrt{2}}\rceil}}$).

\smallskip
{\bf (2.b)} Suppose that $|K|=2$. If $n=2l$ (resp.\ $n=2l+1$) and $2\le r\le \Bigl\lfloor \frac{\lceil \frac{|K|^l}{\sqrt{2}}\rceil}{3}\Bigr\rfloor$, then $N_{\sigma}:=\min\Bigl(\bigl(E_{n-2,k}(l-1)^4\bigr)^{\nu_{3,r}(\sigma)},k^4l^4\bigl(E_{n-2,k}(l-1)^4\bigr)^{\wp_{r,|K|^n-\lceil \frac{|K|^l}{\sqrt{2}}\rceil}}\Bigr)$
(resp.\ $N_{\sigma}:=\min\Bigl([4l^6E_{n-2,k}]^{\nu_{3,r}(\sigma)},E_{n-2,k}^{-\frac{1}{2}}[4l^6E_{n-2,k}]^{\frac{1}{2}+\wp_{r,|K|^n-\lceil \frac{|K|^l}{\sqrt{2}}\rceil}}\Bigr)$).

\medskip
{\bf (3)} Suppose that $s=3$. If $6\le 3r\le |K|$, then 
$$N_{\sigma}:=\Bigl[E_{n-2,k}(3r-1)^2(3r-2)^4\Bigr]^{\nu_{3,r}(\sigma)}\le \Bigl[E_{n-2,k}(3r-1)^2(3r-2)^4\Bigr]^{\wp_{r,\n(\sigma)}}.$$

\smallskip
{\bf (4)} If $s\ge 5$, then $N_{\sigma}:=\Bigl(E_{n-2,k,\frac{s-3}{2}}D_{n,k}^2\Bigr)^{\nu_{s,r}(\sigma)}.$

\smallskip
{\bf (5)} Suppose that $s\ge 5$. If we have $n=2l$ (resp.\ $n=2l+1$) and $r\le\bigl\lfloor \frac{\lceil \frac{|K|^l}{\sqrt{2}}\rceil}{s}\bigr\rfloor$, then $N_{\sigma}:=\bigl(E_{n-2,k,\frac{s-3}{2}}k^8l^8\bigr)^{\nu_{s,r}(\sigma)}$ (resp.\ 
$$N_{\sigma}:=\Bigl[E_{n-2,k,\frac{s-3}{2}}k^4(k+1)^4\Bigl(kl+\Bigl\lfloor\sqrt{\frac{k+1}{2}}\Bigr\rfloor\Bigr)^4\Bigl(kl+\Bigl\lfloor\frac{k+1}{\sqrt{2}}\Bigr\rfloor\Bigr)^8\Bigr]^{\nu_{s,r}(\sigma)}).$$

{\bf (6)} If $s\ge 5$ and $sr\le |K|$, then $N_{\sigma}:=\Bigl[E_{n-2,k,\frac{s-3}{2}}(sr-1)^2(sr-2)^4\Bigr]^{\nu_{s,r}(\sigma)}$.
\end{corollary}

\begin{proof}
For $i\in \llbracket1,r\rrbracket$, let $\mathcal C_i:=i\mathcal Y_s$. Following the notation of Lemma \ref{L9}(2), let $\mathcal U:=\cup_{i=1}^r\mathcal C_i$, $\n(\mathcal U):=\max\bigl(\n(\mathcal C_i)\bigr)$, and $\omega^{\S}(\mathcal U):=\max\bigl(\omega^{\S}(\mathcal C_i)\|i\in \llbracket1,r\rrbracket\bigr)$. Clearly, $\n(\mathcal U)=sr$. Moreover, we have $\omega^{\S}(\mathcal U)\le E_{n-2,k,\frac{s-3}{2}}$ by Proposition \ref{PR18}(1.a) if $s>3$ and Proposition \ref{PR19}(1.a) if $s=3$ (recall $E_{n-2,k}=E_{n-2,k,0}$).

With the notation of and by Lemma \ref{L9}(2) we have inequalities 
$$\ell^{\S}_{n,K}(\sigma)\le\left[\omega^{\S}(\mathcal U)\pi^{\S}_{n,sr}(K)^2\right]^{\mu_{\n(\sigma),\mathcal U}},$$
$$\ell^{\S}_{n,K}(\sigma)\le\pi^{\S}_{n,sr}(K)\left[\omega^{\S}(\mathcal U)\pi^{\S}_{n,sr}(K)^2\right]^{\mu_{|K|^n-sr,\mathcal U}},$$
$$\ell^{\S}_{n,K}(\sigma)\le\pi^{\S}_{n,2\lceil \frac{|K|^{n-1}}{4}\rceil}(K)\left[\omega^{\S}(\mathcal U)\pi^{\S}_{n,2\lceil \frac{|K|^{n-1}}{4}\rceil}(K)^2\right]^{\mu_{|K|^n-2\lceil \frac{|K|^{n-1}}{4}\rceil,\mathcal U}}$$
only for part (1), and
$$\ell^{\S}_{n,K}(\sigma)\le\pi^{\S}_{n,sr}(K)\left[\omega^{\S}(\mathcal U)\pi^{\S}_{n,sr}(K)^2\right]^{\mu_{|K|^n-sr,\mathcal U}}$$
only for parts (2) and (5).

If $s=3$, then, as $\r\ge 2$, we have $\nu_{3,r}(\sigma)\le\wp_{r,\n(\sigma)}$ and for each $m\in \llbracket1,|K|^n\rrbracket$ we have $\mu_{m,\mathcal U}\le \wp_{r,m}$ by Proposition \ref{PR3}(4).

For parts (1) and (4) (resp.\ (3) and (6)), as $\pi^{\S}_{n,sr}(K)\le D_{n,k}$ by Corollary \ref{C20}, (resp.\ $\pi_{n,sr}(K)\le (sr-1)(sr-2)^2$ by Theorem \ref{T8}(2)) we conclude that these parts hold.

For parts (2.a) and (5) and $n=2l$ (resp.\ $2l+1$), we have $\pi_{n,sr}(K)\le k^4l^4$ (even $\pi_{n,sr}(K)\le [(l-1)(k-1)+\lfloor\frac{k+1}{\sqrt{2}}\rfloor]^4$ if $l\ge 2$) by Proposition \ref{PR21}(1.a) (resp.\ $\pi_{n,sr}(K)\le k^2(k+1)^2[kl+\lfloor\sqrt{\frac{k+1}{2}}\rfloor]^2(kl+\lfloor\frac{k+1}{\sqrt{2}}\rfloor)^4$ by Proposition \ref{PR23}(2)), and we similarly conclude that these parts hold.

For part (2.b) and $n=2l$ (resp.\ $2l+1$), we have $\pi_{n,sr}(K)\le (l-1)^4$ (resp.\ $\pi_{n,sr}(K)\le 4l^6$ by Proposition \ref{PR21}(1.b) (resp.\ Proposition \ref{PR23}(1)), and we similarly conclude that these parts hold.\footnote{The bounds based on Propositions \ref{PR21}(1.a) and \ref{PR3}(4) are simplified. Moreover, parts (3) to (5) also have variants that either define $N_{\sigma}$ via minima or are based on Proposition \ref{PR21}(2) and Proposition \ref{PR23}(3).}\end{proof}

In what follows we require inequalities for the real numbers $E_{l,x}$ defined by Equation (\ref{EQ13}) and the real numbers associated to the $\varepsilon_{n,k,2}$s of Lemma \ref{L16}.

\begin{lemma}\label{L18}
For $(l,x)\in [\mathbb N^{\ast}\setminus\{1\}]\times [1,\infty)$ the following inequalities hold.

\medskip
{\bf (1)} Let $N\in\mathbb N^{\ast}\setminus\{1\}$. If $x\ge N$, then $E_{l-1,x}\le\frac{N^3+N^2+1}{2N^3}lx^4\le \frac{13}{16}lx^4$.

\smallskip
{\bf (2)} We have $E_{l-1,1}\le l^{1.6}$ (i.e., $2l-1\le l^{1.6}$) for all $l\ge 2$ and $E_{l-1,1}\le l^{1.3}$ (i.e., $2l-1\le l^{1.3}$) for all $l\ge 9$.
\end{lemma}

\begin{proof}
Part (1) follows from the estimates
$$\frac{E_{l-1,x}}{lx^4}=\frac{1}{l}+(x^{-1}+x^{-3})\frac{l-1}{l}\le \frac{1}{l}+\frac{(N^2+1)(l-1)}{N^3l}=\frac{N^2+1}{N^3}+\frac{N^3-N^2-1}{N^3l}$$
$$\le \frac{N^2+1}{N^3}+\frac{N^3-N^2-1}{2N^3}=\frac{N^3+N^2+1}{2N^3}\le\frac{13}{16}.$$

\phantomsection{For part (2) we consider the two real-valued functions $\mathfrak f_1:[2,\infty)\to\mathbb R$ and $\mathfrak f_2:[9,\infty)\to\mathbb R$ defined by the rules $\mathfrak f_1(y):=y^{1.6}-2y$ and $\mathfrak f_2(y):=y^{1.3}-2y+1$. As $\mathfrak f_1'(y)=1.6y^{0.6}-2$ and $\mathfrak f_2'(y)=1.3y^{0.3}-2$ are increasing functions with $\mathfrak f_1'(2)>0$ and $\mathfrak f_2'(9)>0$, we get that $\mathfrak f_1$ and $\mathfrak f_2$ are increasing functions. From this and the inequalities $\mathfrak f_1(2)=2^{1.6}-3>0.03$ and $\mathfrak f_2(9)=9^{1.3}-17>0.39>0$ we get that all values of $\mathfrak f_1$ and $\mathfrak f_2$ are positive and hence part (2) holds.}\label{PH102}
\end{proof}

\begin{lemma}\label{L19}
Let $(n,k)\in (\mathbb N^{\ast}\setminus\{1,2\})\times (\mathbb N^{\ast}\setminus\{1\})$. We define real numbers 
$$\varepsilon_{n,k,1}:=-\frac{(k-2)(2n+3)}{k(n-1)}-\frac{n^2+n-6}{n-1}-\frac{(k-2)^2(2n+3)}{4k^2(n-1)^2}-\frac{n^2+n-6}{4(n-1)^2}$$ 
and
$$\epsilon_{n,k}:=-\frac{4n^2+8n-14}{k+1}-\frac{2n^2+4n-7}{(k+1)^2}+\varepsilon_{n,k,1}+\varepsilon_{n,k,2},$$
where $\varepsilon_{n,k,2}$ is as in Lemma \ref{L16}. Then the following properties hold.

\medskip
{\bf (1)} We have $3n-\frac{5}{3}\le 2n+4\ln(n-2)$ iff $n\in \llbracket4,9\rrbracket$.

\smallskip
{\bf (2)} If $n\in \llbracket4,9\rrbracket$, then $2n+4\ln(n-2)-3n+\frac{5}{3}\in (0.439,1.212)$ attains the largest value at $n=6$.

\smallskip
{\bf (3)} We have $\varepsilon_{n,k,2}+2n-\frac{2}{3}< \frac{4n+8\ln(n-2)+2.424}{k}<\frac{4n+8\ln(n-1)-\frac{8}{n-1}-\frac{4}{(n-1)^2}+2.424}{k}$.

\smallskip
{\bf (4)} We have $\epsilon_{n,2}< 4\ln(n-1)+5.0732-\frac{14n^2+37n}{9}$.

\smallskip
{\bf (5)} We have $\epsilon_{n,3}< \frac{8}{3}\ln(n-1)+2.4955-\frac{27n^2+94n}{24}$.

\smallskip
{\bf (6)} For $k\ge 4$ we have $\epsilon_{n,k}< \frac{8}{k}\ln(n-1)-\frac{n^2}{k+1}-3n-\frac{37}{48}$.

\smallskip
{\bf (7)} We have $\epsilon_{n,k}<0$.
\end{lemma}

\begin{proof}
Parts (1) and (2), we consider the real-valued function $\mathfrak f:[3,\infty)\rightarrow\mathbb R$ defined by $\mathfrak f(x):=x-\frac{5}{3}-4\ln(x-2)$. As $\mathfrak f'(x)=1-\frac{4}{x-2}$, $\mathfrak f$ is strictly decreasing on the interval $[3,6]$ and strictly increasing on the interval $[6,\infty)$. So $\mathfrak f$ attains its minimum at $6$. We have $\mathfrak f(10)=0.0155...>0$, $\mathfrak f(9)=-0.4503...$, $\mathfrak f(6)=\frac{13}{3}-4\ln(4)=-1.2118...$, $\mathfrak f(4)=-0.4392...$, and $\mathfrak f(3)=\frac{4}{3}>0$. Hence parts (1) and (2) hold.

As $\varepsilon_{n,k,2}+2n-\frac{2}{3}=\frac{2\min\bigl(2n+4\ln(n-2),3n-\frac{5}{3}\bigr)}{k}$, part (3) follows directly from parts (1) and (2) and the relations $\ln(n-2)-\ln(n-1)=\ln(1-\frac{1}{n-1})<-\frac{1}{n-1}-\frac{1}{2(n-1)^2}$.

For part (4), we have identities
$$\varepsilon_{n,2,1}=-\frac{n^2+n-6}{n-1}-\frac{n^2+n-6}{4(n-1)^2}=-n-\frac{9}{4}+\frac{13}{4(n-1)}+\frac{1}{(n-1)^2}$$
and $\epsilon_{n,2}=-\frac{14n^2}{9}-\frac{28n}{9}+\frac{49}{9}+\varepsilon_{n,2,1}+\varepsilon_{n,2,2}$. Moreover, 
$$1.212+\frac{49}{9}-\frac{9}{4}-\frac{3}{4(n-1)}-\frac{1}{(n-1)^2}+\frac{2}{3}<1.212+\frac{49}{9}-\frac{9}{4}+\frac{2}{3}<5.0732.$$
Based on the last two sentences and part (3) we get that part (4) holds.

For part (5), we have an inequality $\varepsilon_{n,3,1}<\varepsilon_{n,2,1}-\frac{2}{3}-\frac{1}{18(n-1)}-\frac{5}{3(n-1)}$
and an identity $\epsilon_{n,3}=-\frac{9n^2}{8}-\frac{9n}{4}+\frac{63}{16}+\varepsilon_{n,3,1}+\varepsilon_{n,3,2}$. 
Moreover, 
$$\frac{2.424}{3}+\frac{63}{16}-\frac{9}{4}-\frac{41}{36(n-1)}-\frac{1}{3(n-1)^2}<\frac{2.424}{3}+\frac{27}{16}=2.4955.$$
Based on the last two sentences and part (3) we get that part (5) holds.

For part (6), we have inequalities 
$$\varepsilon_{n,k,1}<\varepsilon_{n,2,1}-\frac{2(k-2)}{k}-\frac{1}{8(n-1)}\le-n-1-\frac{9}{4}+\frac{25}{8(n-1)}+\frac{1}{(n-1)^2}\le -n-\frac{23}{16}$$
(as $n\ge 3$) and $\epsilon_{n,k}<-\frac{n^2}{k+1}-\frac{8n}{k+1}-\frac{13}{k+1}+\varepsilon_{n,k,1}+\varepsilon_{n,k,2}$ as for $n\ge 3$ we have $4n^2+8n-14\ge n^2+8n+13$. 
Moreover, 
$$n\Bigl(-\frac{8}{k+1}+\frac{4}{k}\Bigr)+\frac{2}{3}-\frac{23}{16}-\frac{13}{k+1}+\frac{2.424}{k}<\frac{2}{3}-\frac{23}{16}=-\frac{37}{48}.$$
Based on the last two sentences and part (3) we get that part (6) holds.

Part (7) follows from parts (4) to (6).\end{proof}

We have the following consequence of parts of Lemmas \ref{L15}, \ref{L18}, and \ref{L19}.

\begin{corollary}\label{C23} Let $(n,k,s)\in (\mathbb N^{\ast}\setminus\{1,2\})\times\mathbb N^{\ast}\times (2\mathbb N^{\ast}+1)$. Let $\epsilon_{n,k}$ and $\epsilon_{n,k,1}$ be as in Lemma \ref{L19} and let $\epsilon_{n,k,2}$ be as in Lemma \ref{L16}. For $N\in\mathbb N^{\ast}\setminus\{1\}$ let $C_N:=\frac{N^3+N^2+1}{2N^3}$. Then the following inequalities hold.

\medskip
{\bf (1)} If $k\ge N$, then we have inequalities
$$E_{n-2,k,\frac{s-3}{2}}D_{n,k}^2\le C_Nk^{4n^2+8n-14}(n-1)^{4n^2+2n-3}[\overline{C}^+_{n,k,2}]^2e^{\epsilon_{n,k,2}}$$
$$\le C_Nk^{4n^2+8n-14}(n-1)^{4n^2+2n-3}e^{\varepsilon_{n,k,1}+\varepsilon_{n,k,2}}$$
$$<C_N(k+1)^{4n^2+8n-14}(n-1)^{4n^2+2n-3}e^{\epsilon_{n,k}}.$$
In particular, if $k\ge 2$, then 
$$E_{n-2,k,\frac{s-3}{2}}D_{n,k}^2\le \frac{13}{16}(k+1)^{4n^2+8n-14}(n-1)^{4n^2+2n-3}.$$

{\bf (2)} For $k=1$ we have the inequalities
$$E_{n-2,1}D_{n,1}^2\le (n-2)^{2n^2+2n-4}(n-1)^{2n^2+9.6}e^{2\min\bigl(n+\frac{1}{3},2n-\frac{4}{3}-4\ln(n-2)\bigr)}$$
$$< (n-1)^{4n^2+2n+5.6}e^{[-1-\frac{1}{2(n-1)}]\frac{2n^2+2n-4}{n-1}}e^{2\min\bigl(n+\frac{1}{3},2n-\frac{4}{3}-4\ln(n-2)\bigr)}$$
$$=(n-1)^{4n^2+2n+5.6}e^{-5-\frac{3}{n-1}+2\min\bigl(\frac{1}{3},n-\frac{4}{3}-4\ln(n-2)\bigr)}$$
$$\le(n-1)^{4n^2+2n+5.6}e^{\frac{2}{3}-5-\frac{3}{n-1}}<0.0132(n-1)^{4n^2+2n+5.6}.$$

{\bf (3)} For $k=1$ and $n\ge 10$ we have the inequalities
$$E_{n-2,1}D_{n,1}^2\le (n-2)^{2n^2+2n-4}(n-1)^{2n^2+9.3}e^{2n+\frac{2}{3}}$$
$$<(n-1)^{4n^2+2n+5.3}e^{[-1-\frac{1}{2(n-1)}]\frac{2n^2+2n-4}{n-1}}e^{2n+\frac{2}{3}}=(n-1)^{4n^2+2n+5.3}e^{\frac{2}{3}-5-\frac{3}{n-1}}$$
$$<(n-1)^{4n^2+2n+5.3}e^{\frac{2}{3}-5}<0.0132(n-1)^{4n^2+2n+5.3}.$$
\end{corollary}

\begin{proof}
Note that $E_{n-2,k,\frac{s-3}{2}}$ is a decreasing function in odd integers $s\ge 3$. So the first inequality of part (1) follows by a direct additions of exponents from the inequalities $E_{n-2,k,\frac{s-3}{2}}\le E_{n-2,k,0}=E_{n-2,k}\le C_Nk^4(n-1)$ (see Lemma \ref{L18}(1)) and $D_{n,k}\le\overline{C}_{n,k}^+$ (see Corollary \ref{C20}). 

The first inequalities of parts (2) and (3) follow by a direct additions of exponents from the inequalities $2n+4\ln(n-2)<3n-\frac{5}{3}$ for $n\ge 10$ (see Lemma \ref{L19}(1)) and $D_{n,1}<C_{n,1}$ (see Corollary \ref{C20} and Lemma \ref{L16}(2)) and from Lemma \ref{L18}(2).

The second inequality of part (1) follows from the fact that for $n\ge 3$ and $k\ge 2$ we have 
$$-\varepsilon_{n,k,1}=\frac{(k-2)(2n+3)}{k(n-1)}+\frac{n^2+n-6}{n-1}+\frac{(k-2)^2(2n+3)}{4k^2(n-1)^2}+\frac{n^2+n-6}{4(n-1)^2}$$
$$<\begin{cases}
\frac{n^2+n-6}{n-1}+\frac{n^2+n-6}{2(n-1)^2} \,\;\;\quad\quad\quad\quad\quad\quad\quad\quad\quad\quad\quad\quad\quad\quad\quad\quad\quad\quad\;\;\,\;\;\, {\rm if}\;\; k=1\\
\frac{(k+1)(n^2+3n-3)}{2k(n-1)}+\frac{(k+1)^2(n^2+3n-3)}{8k^2(n-1)^2} \;\;\;\quad\quad\quad\quad\quad\quad\quad\quad\quad\quad\quad\quad\,\;\; {\rm if}\;\;4\mid k+1\\
\frac{(k-1)(2n+3)}{2k(n-1)}+\frac{(k+1)(n^2+n-6)}{2k(n-1)}+\frac{(k-1)^2(2n+3)}{8k^2(n-1)^2}+\frac{(k+1)^2(n^2+n-6)}{8k^2(n-1)^2} \;\;\;\quad {\rm if}\;\; k\in 5+4\mathbb N\\
\frac{n^2+3n-3}{2(n-1)}+\frac{n^2+3n-3}{8(n-1)^2}
\quad\quad\quad\quad\quad\quad\quad\quad\quad\quad\quad\quad\quad\quad\quad\quad\quad\quad\quad\;\, {\rm if}\;\; 4\mid k-2.\end{cases}$$
which gets translated into $[\overline{C}^+_{n,k,2}]^2\le e^{\varepsilon_{n,k,1}}$ by Display (\ref{EQ31}) and the last inequality is strict iff $4\nmid k$. 

The other inequalities of parts (1) to (3) follow directly from the following estimates $(n-2)^{n-1}=(n-1)^{n-1}\left(1-\frac{1}{n-1}\right)^{n-1}<(n-1)^{n-1}e^{-1-\frac{1}{2(n-1)}}$ and 
$$k^{4n^2+8n-14}=(k+1)^{4n^2+8n-14}\bigl(1-\frac{1}{k+1}\bigr)^{4n^2+8n-14}$$
$$<(k+1)^{4n^2+8n-14}\bigl(1-\frac{1}{k+1}\bigr)^{4n^2+8n-14}<(k+1)^{4n^2+8n-14}e^{\varepsilon_{n,k}-\varepsilon_{n,k,1}-\varepsilon_{n,k,2}}$$
(see Inequality (\ref{EQ32})) and, for part (1),  Lemma \ref{L19}(2) and (7).
\end{proof}

\section{Exponential upper bounds for all $m$ via affine permutations}\label{S22}

In this section we provide upper bounds for the $\pi_{n,m}(\mathbb F_{p^q})$s with $m\in \llbracket3,p^{nq}\rrbracket$ if $4\nmid |K|$ and with $m\in \llbracket3,p^{nq}-2\rrbracket$ if $4\mid |K|$ and for the $\pi^{\S}_{n,m}(\mathbb F_{p^q})$s with $\break m\in \llbracket3,p^{nq}-2\rrbracket$ via product decompositions of permutations $\sigma\in\perm(\mathbb F_{p^q}^n)$ of the form $\sigma=\theta b(\mathbb F_{p^q})a(\mathbb F_{p^q})$ for triples $(\theta,a,b)\in\Alt(\mathbb F_{p^q}^n)\times\AGL_n(\mathbb F_{p^q})\times\STGA_n(\mathbb F_{p^q})[l]$ with $l\in\{1,n-1\}$.

\begin{lemma}\label{L20}
For $n\in\mathbb N^{\ast}\setminus\{1\}$ and $K$ a finite field the following properties hold.

\medskip
{\bf (1)} Each permutation $\sigma\in\perm(K^n)$ is a product $\theta a(K)$ with $a\in\AGL_n(K)$ and $\theta\in\perm(K^n)$ such that $\n(\theta)\le |K|^n-n-1$.

\smallskip
{\bf (2)} If $|K|\ge 3$ and $\sigma\in\Perm(K^n)$, then in part (1) we can assume that $\theta\in\Alt(K^n)$.

\smallskip
{\bf (3)} If $|K|$ is odd, then in part (1) we can assume that $\theta\in\perm(K^n)\setminus\Alt(K^n)$.

\smallskip
{\bf (4)} If $|K|=2$ and $n\ge 3$, then each $\sigma\in\perm(K^n)$ is a product $\theta_0 b(K)a(K)$ with $(\theta_0,a,b)\in\Alt(K^n)\times\AGL_n(K)\times\STGA_n(K)[n-1]$ and $\n(\theta_0)\le 2^n-n-1$.

\smallskip
{\bf (5)} Each permutation $\sigma\in\perm(K^n)$ is a product $\theta a(K)$ with $a\in\ASL_n(K)$ and $\theta\in\perm(K^n)$ such that $\n(\theta)\le |K|^n-n$.
\end{lemma}

\begin{proof}
Let $\sigma\in\perm(K^n)$ and $P_0\in K^n$. By induction on $i\in \llbracket1,n\rrbracket$ we show that there exist distinct points $P_1,\ldots,P_i$ in $K^n\setminus\{P_0\}$ such that both the families $(P_j-P_0)_{j\in \llbracket1,i\rrbracket}$ and $\bigl(\sigma(P_j)-\sigma(P_0)\bigr)_{j\in \llbracket1,i\rrbracket}$ are linearly independent. 

For $i=1$ we take $P_1\in K^n\setminus\{P_0\}$ arbitrarily, so the base of the induction holds. If $i\in \llbracket2,n\rrbracket$, for the inductive passage from $i-1$ to $i$, the complement $Z_i$ of
$$\Span(\{P_j-P_0|j\in \llbracket1,i-1\rrbracket\})\cup\sigma^{-1}\bigl(\Span(\{\sigma(P_j)-\sigma(P_0)|j\in \llbracket1,i-1\rrbracket\})\bigr)$$
in $K^n$ has at least $|K|^n-2|K|^{i-1}+(i-1)>0$ elements and thus it is non-empty and we can take $P_i$ arbitrarily in the set $Y_i:=P_0+Z_i$. This completes the inductive step and the induction. 

If $a\in\AGL_n(K)$ (resp.\ $a\in\ASL_n(K)$) is such that we have $a(P_j)=\sigma(P_j)$ for each $j\in \llbracket0,n\rrbracket$ (resp.\ $j\in \llbracket0,n-1\rrbracket$), then for $\theta:=\sigma [a(K)]^{-1}\in\perm(K^n)$ we have an identity $\sigma=\theta a(K)$ and an inclusion $\supp(\theta)\subset K^n\setminus\{\sigma(P_i)|i\in \llbracket0,n\rrbracket\}$ (resp.\ $\supp(\theta)\subset K^n\setminus\{\sigma(P_i)|i\in \llbracket0,n-1\rrbracket\}$), hence $\n(\theta)\le |K|^n-n-1$ (resp. $\n(\theta)\le |K|^n-n$) and part (1) (resp.\ (5)) holds. 

To prove parts (2) and (3) we can assume that $P_0=(0,\ldots,0)$ is the zero vector and, based on Theorem \ref{T6}(3), that $|K|$ is odd. For the passage from $n-1$ to $n$, $P_n$ is an arbitrary point in a subset $Y_n\subset K^n$ with $|Y_n|\ge |K|^n-2|K|^{n-1}+(n-1)$; for $P_n\in Y_n$, let $a_{P_n}\in\GL_n(K)$ be such that $a_{P_n}(P_j)=\sigma(P_j)$ for each $j\in \llbracket1,n\rrbracket$. 

We fix a point $Q\in Y_n$. The rule $P_n\mapsto a_Q^{-1}a_{P_n}$ defines an injective map $Y_n\rightarrow \nabla$, where $\nabla$ is the subgroup of $\GL_n(K)$ that fixes $P_i$ for each $i\in \llbracket1,n-1\rrbracket$. Let $\delta\in K^{\ast}$ be a generator of the multiplicative cyclic group $K^{\ast}$. As $|K|$ is odd, for each $c\in \nabla$ that has eigenvalue $\delta$, the permutation $c(K)$ is odd. Similarly, for each $c\in \nabla$ that has all eigenvalues $1$, the permutation $c(K)$ is even by Theorem \ref{T3}(1) applied to a subgroup $G=\mathbb G_{\a,K}^{n-1}$ of $\SL_n(K)$ with the property that $G(K)\leqslant\nabla$. Thus for $\Gamma\in\{\Alt(K^n),\perm(K^n)\setminus\Alt(K^n)\}$ we have $|\varrho_{n,K}(\nabla)\cap\Gamma|\ge |K|^{n-1}$. From this and the identity $|\varrho_{n,K}(\nabla)|=|K|^n-|K|^{n-1}$ we get the folowing inequality $|\varrho_{n,K}(\nabla)\cap\Gamma|\le |K|^n-2|K|^{n-1}<|Y_n|$. 

So there exists a point $P_{n,\Gamma}\in Y_n$ such that $a_Q(K)^{-1}a_{P_{n,\Gamma}}(K)\in\Gamma$; hence $a_{P_{n,\Gamma}}(K)\in a_Q(K)\Gamma$ and parts (2) and (3) hold.

It suffices to prove part (4) only when $\theta$ is odd. We take $b$ such that we have $b(K)=a_1(K)\shift^{(1,0,\ldots,0)}_{V\oplus W}[a_1(K)]^{-1}$, with the selective shift $\shift^{(1,0,\ldots,0)}_{V\oplus W}\in\Perm(K^n)$ as in the proof of Proposition \ref{PR10}(2) for $\dim_K(W)=n-1$ and therefore a transposition by Lemma \ref{F5}(1), with $a_1\in\AGL_n(K)$, and with $\supp\bigl(b(K)\bigr)=a_1(V)$ contained in $K^n\setminus\{\sigma(P_i)|i\in \llbracket0,n\rrbracket\}$; hence $b\in\STGA_n(K)[n-1]$. This is possible, as up to linear automorphisms we can assume that all entries of the $\sigma(P_i)$s are $1$ except that for each $i\in \llbracket1,n\rrbracket$, $\sigma(P_i)$ has the $i$-th entry $0$, so we can take $a_1:=1_{\mathbb A^n_K}$. As $b(K)$ is odd, $\theta_0:=\theta b(K)^{-1}\in\Alt(K^n)$ as we are assuming $n\ge 3$. As $\supp(\theta_0)\subset K^n\setminus\{\sigma(P_i)|i\in \llbracket0,n\rrbracket\}$ and $\sigma=\theta a(K)=\theta_0b(K)a(K)$, part (4) holds.\end{proof}

\begin{proposition}\label{PR25}
Let $n\in\mathbb N^{\ast}\setminus\{1\}$ and $q\in \llbracket1,n\rrbracket$. Let $K$ be a finite field. Let $m\in \llbracket|K|^{q-1}+1,|K|^q\rrbracket$ and $(\underline{P},\underline{Q})=\bigl((P_1,\ldots,P_m),(Q_1,\ldots,Q_m)\bigr)\in\mathbb D_{n,m}(K)^2$. Then the following properties hold.

\medskip
{\bf (1)} Suppose that $K=\mathbb F_2$ and $q\ge 2$. Let $\epsilon\in\{0,1\}$ be such that we have $\epsilon=0$ iff $m\in \llbracket2^{q-1}+1,2^q-q\rrbracket$. Then there exists a pair $(a,\theta)\in\AGL_n(\mathbb F_2)\times\perm(\mathbb F_2^n)$ such that $\nu_{2,1}(\theta)\le\min(m-q-\epsilon,2^n-q-1-\epsilon)$, $\n(\theta)\le\min(2m-2q-2\epsilon,2^n-q-\epsilon)$, and for $\sigma:=\theta a(\mathbb F_2)$ we have $\sigma(\underline{P})=\underline{Q}$.

\smallskip
{\bf (2)} If $|K|\ge 3$, then there exists a pair $(a,\theta)\in\AGL_n(K)\times\perm(K^n)$ such that $\nu_{2,1}(\theta)\le m-q-1$, $\n(\theta)\le\min(2m-2q-2,|K|^n-q-1)$, and for $\sigma:=\theta a(K)$ we have $\sigma(\underline{P})=\underline{Q}$.

\smallskip
{\bf (3)} Suppose that $|K|\ge 3$. If $q\in \llbracket1,n-1\rrbracket$ (resp.\ $q=n$ and we have either $m\in\llbracket2,|K|^n-2\rrbracket$ or $m\in\{|K|^n-1,|K|^n\}$ and the only permutation $\sigma\in\perm(K^n)$ with $\sigma(\underline{P})=\underline{Q}$ is in $\Alt(K^n)$), then there exists a pair $(a,\theta)\in\AGL_n(K)\times\Alt(K^n)$ such that for the permutation $\sigma:=\theta a(K)\in\perm(K^n)$ we have $\sigma(\underline{P})=\underline{Q}$ and the following inequalities $\nu_{2,1}(\theta)\le 2\lfloor\frac{m-q}{2}\rfloor$, $\n(\theta)\le 2m-2q-1$, and $\nu_{3,1}(\theta)\le m-q-1$ (resp.\ $\nu_{2,1}(\theta)\le\min(2\lfloor\frac{m-q}{2}\rfloor,|K|^n-n)$, $\n(\theta)\le\min(2m-2q-1,|K|^n-n)$, and $\nu_{3,1}(\theta)\le\min(m-q-1,\lfloor\frac{|K|^n-n}{2}\rfloor)$).

\smallskip
{\bf (4)} Suppose that $|K|\ge 3$ and $m\ge q+3$. If $q\in \llbracket1,n-1\rrbracket$ (resp.\ $q=n$ and either $n+3\le m\le |K|^n-2$ or $m\in\{|K|^n-1,|K|^n\}$ and the only permutation $\sigma\in\perm(K^n)$ with $\sigma(\underline{P})=\underline{Q}$ is in $\Alt(K^n)$), then there exists a pair $\break (a,\theta)\in\AGL_n(K)\times\Alt(K^n)$ such that for $\sigma:=\theta a(K)$ we have $\sigma(\underline{P})=\underline{Q}$ and inequalities $\nu_{2,1}(\theta)\le m-q-1$, $\n(\theta)\le 2m-2q-2$, and $\nu_{3,1}(\theta)\le 2\lfloor\frac{m-q-1}{2}\rfloor$ (resp.\ $\nu_{2,1}(\theta)\le\min(m-n,|K|^n-n-1)$, $\n(\theta)\le\min(2m-2n-2,|K|^n-n-1)$, and $\nu_{3,1}(\theta)\le\min(2\lfloor\frac{m-n-1}{2}\rfloor,\lfloor\frac{|K|^n-n-1}{2}\rfloor)$).\end{proposition}

\begin{proof}
To prove part (1) we note that $m\ge 2^{q-1}+1\ge q+1\ge q+\epsilon$. By induction on $i\in \llbracket2+\epsilon,q+\epsilon\rrbracket$ we show that there exists a subset $I_i\subset \llbracket2,m\rrbracket$ such that $|I_i|=i-1$ and the families $(P_j-P_1)_{j\in I_i}$ and $(Q_j-Q_1)_{j\in I_i}$ are linearly independent. Taking $I_{2+\epsilon}:=\{2,2+\epsilon\}$ the base of the induction holds. With $i\in \llbracket2+\epsilon,q+\epsilon-1\rrbracket$, for the passage from $i$ to $i+1$ we can assume that $I_i=\llbracket2,i\rrbracket$. Let $l\in \llbracket i+1,m\rrbracket$ be such that $P_l\notin\Span(\{P_j-P_1|j\in I_i\})$ and $Q_l\notin\Span(\{Q_j-Q_1|j\in I_i\})$, which is possible because for $\epsilon=0$ we have $m-i\ge 2^{q-1}+1-i> 2[2^{i-1}-i]=2^i-2i$ as $i\le q-1$ and for $\epsilon=1$ we have $m-i\ge 2^q-q+1-i> 2[2^{i-1}-i]=2^i-2i$ as $i\le q$; so we can take $I_{i+1}:=I_i\cup\{l\}$. This ends the induction. 

Let $I:=\{1\}\cup I_{q+\epsilon}$. Let $a\in\AGL_n(\mathbb F_2)$ be such that $a(P_i)=Q_i$ for each $i\in I$. Let $\sigma_0\in\perm(\mathbb F_2^n)$ be such that $\sigma_0(\underline{P})=\underline{Q}$. So for $\sigma_1:=\sigma_0 a(\mathbb F_2)^{-1}\in\perm(\mathbb F_2^n)$ we have $\n(\sigma_1)\le 2^n-q-\epsilon$ as $\sigma_1(Q_i)=Q_i$ for each $i\in I$ and $|I|=q+\epsilon$. 

Let $J:=\llbracket1,m\rrbracket\setminus I$ and $s:=|J|=m-q-\epsilon$. 

Let $A:=\mathbb F_2^n\setminus\{Q_i|i\in I\}$ and $l:=|A|=2^n-q-\epsilon$. Let $\sigma_A\in\perm(A)$ be such that for each $P\in A$ we have $\sigma_ A(P)=\sigma_1(P)$; so for each $i\in J$ we have $\sigma_A\bigl(a(P_i)\bigr)=\sigma_1\bigl(a(P_i)\bigr)=\sigma_0(P_i)=Q_i\in A$. There exists $\theta_A\in\perm(A)$ such that $\nu_{2,1}(\theta_A)\le \min(s,l-1)=\min(m-q-\epsilon,2^n-q-1-\epsilon)$ and $\theta_A\bigl(a(P_i)\bigr)=Q_i$ for each $i\in J$ by Lemma \ref{L1}(1) applied to $\sigma_A$. Taking $\theta\in\perm(\mathbb F_2^n)$ such that $\theta(Q_i)=Q_i$ for each $i\in I$ and $\theta(P)=\theta_A(P)$ for each $P\in A$, we have
$$\n(\theta)=\n(\theta_A)\le\min\bigl(2\nu_{2,1}(\theta_A),2^n-q-\epsilon\bigr)\le\min(2m-2q-2\epsilon,2^n-q-\epsilon)$$ 
and we get that part (1) holds for $\sigma:=\theta a(\mathbb F_2)\in\perm(\mathbb F_2^n)$. 

Part (2) (resp.\ (3) or (4)) is proved in the same way as part (1) with $\epsilon=1$ based on Lemma \ref{L1}(1) (resp.\ Lemma \ref{L1}(2) or (3) applied to $(l,s)$ being equal to $(|K|^n-q-1,m-q-1)$ and Lemma \ref{P3}(2)).\end{proof}

\begin{example}\normalfont\label{EX18}
Let $n\in\mathbb N^{\ast}\setminus\{1\}$. Let $K$ be a finite field. Let $r\in \llbracket3,n+1\rrbracket$. We have $\pi_{n,r}(K)\le (r-1)^2$ by Theorem \ref{T7}(2). Let $m\in \llbracket r+1,|K|^n\rrbracket$. If $4\mid |K|$ we assume that $m\le |K|^n-2$ and $r$ is odd. If $|K|=2$, then we assume that $n\ge 3$.

\medskip
{\bf (1)} Let $s:=\lceil\frac{|K|^n-n-1}{r}\rceil\in\mathbb N^{\ast}$. We have $rs\in \llbracket|K|^n-n-1,|K|^n-1\rrbracket$ as $r\le n+1$. For $\sigma\in\Perm(K^n)$ we consider a product decomposition $\sigma=\theta a(K)b(K)$ with $(\theta,a,b)\in\Alt(K^n)\times\AGL_n(K)\times\STGA_n(K)[n_0]$ such that $\supp(\theta)\le |K|^n-n-1$ and $n_0=1$ if $|K|\ge 3$ (by Lemma \ref{L20}(1) and (2)) and $n_0=n-1$ if $|K|=2$ (by Lemma \ref{L20}(4)). Let $(l,t)\in (\mathbb N^{\ast})^2$ be such that $lrt\in\llbracket rs,|K|^n\rrbracket$. For the $\theta\in\Alt(Y)$ with $Y\subset K^n$ such that $|Y|=rs\le lrt$ we have $\nu^+_{r,l}(\theta)\le 2t$ by Theorem \ref{P11}(1) if $r\ge 4$ and $\nu^+_{r,l}(\theta)\le 2+\nu_{r,l}(\theta)\le 2+2\bigl(\lfloor\frac{3t+2}{4}\rfloor+1\bigr)$ by Proposition \ref{PR2}(3) and Theorem \ref{P8}(2) and (3) if $r=3$. Recall that $l\mathcal Y_r$ is the conjugacy class of $\perm(K^n)$ formed by $l$ disjoint $r$-cycles. We have an inequality $\eth(l\mathcal Y_r)\le\pi_{n,lr}(K)$ by Lemma \ref{L8}(1). Moreover, by writing $\theta$ as a product of $\nu^+_{r,l}(\theta)$ permutations in $l\mathcal Y_r$, from Inequality (\ref{EQ3}) we get that $\ell_{n,K}(\sigma)\le n_0\Omega(l\mathcal Y_r)^{\nu^+_{r,l}(\theta)}$. From this and Lemma \ref{L8}(1) we get that
 there exists $c\in\TGA_n(K)[n_0\eth(l\mathcal Y_r)^{2\nu^+_{r,l}(\theta)}\omega(l\mathcal Y_r)^{\nu^+_{r,l}(\theta)}]$ with $c(K)=\sigma$. So we have inequalities
\begin{equation}\label{EQ34}
\ell_{n,K}(\sigma)\le n_0\eth(l\mathcal Y_r)^{2\nu_{r,l}^+(\theta)}\omega(l\mathcal Y_r)^{\nu_{r,l}^+(\theta)}\le n_0\pi_{n,lr}(K)^{2\nu_{r,l}^+(\theta)}\omega(l\mathcal Y_r)^{\nu_{r,l}^+(\theta)}.
\end{equation}
Similarly, if $m\le |K|^n-2$, $r$ is odd, and $\sigma\in\Alt(K^n)$, then we have inequalities
\begin{equation}\label{EQ35}
\ell^{\S}_{n,K}(\sigma)\le \eth^{\S}(l\mathcal Y_r)^{2\nu_{r,l}^+(\theta)}\omega^{\S}(l\mathcal Y_r)^{\nu_{r,l}^+(\theta)}\le \pi^{\S}_{n,lr}(K)^{2\nu_{r,l}^+(\theta)}\omega^{\S}(l\mathcal Y_r)^{\nu_{r,l}^+(\theta)}.
\end{equation}

{\bf (2)} Suppose that $m\le |K|^n-2$. Let $q\in \llbracket1,n\rrbracket$ be the unique integer such that $m\in \llbracket|K|^{q-1}+1,|K|^q\rrbracket$. Let $\epsilon\in\{0,1\}$ be such that $\epsilon=0$ iff $K\cong\mathbb F_2$ and $m\in \llbracket2^{q-1}+1,2^q-q\rrbracket$. Let $$N=N(n,m,K):=\min(2m+1-2q-2\epsilon,|K|^n+1-q-\epsilon).$$ 
It is easy to see that $N\le |K|^n+1-n$ and moreover we have $N\le |K|^n-n$ if $q\le n-1$ or $\epsilon=1$. Let $(\underline{P},\underline{Q})=\bigl((P_1,\ldots,P_m),(Q_1,\ldots,Q_m)\bigr)\in\mathbb D_{n,m}(K)^2$. Let $(a,\theta)\in\AGL_n(K)\times\perm(K^n)$ be such that $\theta\bigl(a(K)(\underline{P})\bigr)=\underline{Q}$, we have inequalities $\nu_{2,1}(\theta)\le m-q-\epsilon$ and $\n(\theta)\le\min(2m-2q-2\epsilon,|K|^n-q-\epsilon)=N-1$, and $\theta$ fixes a subset $Z$ of $\{Q_i|i\in\llbracket1,m\rrbracket\}$ with $m-q-\epsilon$ elements by Proposition \ref{PR25}(1) and (2) and its proof. If $\theta\in\Alt(K^n)$ let $\vartheta:=\theta$. If $\theta\not\in\Alt(K^n)$, let $\vartheta:=\tau\theta$ with $\tau\in\perm(K^n)$ a transposition with $\supp(\tau)\cap\{Q_i|i\in \llbracket1,m\rrbracket\}=\emptyset$ and $\n(\vartheta)\le \max\bigl(\n(\theta)+1,m+2-q-\epsilon\bigr)$; more precisely, if $\n(\theta)>m$, we have $\n(\vartheta)\le \n(\theta)+1$ and if $\n(\theta)\le m$, we have $\n(\vartheta)\le m+2-q-\epsilon$ by the existence of the set $Z$. So $\sigma:=\vartheta a(K)$ is such that $\sigma(\underline{P})=\underline{Q}$, $\vartheta\in\Alt(K^n)$, and $\n(\vartheta)\le N$. Let $s:=\lceil\frac{N}{r}\rceil\in\mathbb N^{\ast}$. If $N=|K|^n-n+1$ (so $(q,\epsilon)=(n,0)$) we assume that either $r\le n$ (so $N\le |K|^n-r+1$ and $rs\le |K|^n$) or $r=n+1$ and $rs\le |K|^n$. If $N\le |K|^n-n$, then $N\le |K|^n-r+1$ as $r\le n+1$ and thus $rs\le|K|^n$. So in all cases we have $rs\le |K|^n$. Let $(l,t)\in\mathbb N^2$ be such that $rlt\in\llbracket rs,|K|^n\rrbracket$. Taking $n_0:=1$, as in part (1) we argue that Inequalities (\ref{EQ34}) and (\ref{EQ35}) hold.\end{example}

\begin{theorem}\label{T11}
Let $(n,m)\in [\mathbb N^{\ast}\setminus\{1\}]\times \mathbb N^{\ast}$. Let $K$ be a finite field such that $|K|^n-1\ge m$. Let $k:=|K|-1$. Let $q\in \llbracket1,n\rrbracket$ be such that $m\in \llbracket|K|^{q-1}+1,|K|^q\rrbracket$. If $4\mid |K|$ and $q=n$, then we assume that $m\le |K|^n-2$. If $|K|\ge 3$, we assume that $m\ge q+3$ (i.e., $(m,q)\notin\{(2,1),(3,1)\}$ and $(m,q,|K|)\neq (4,2,3)$). Then the following properties hold.

\medskip
{\bf (1)} Suppose that $|K|\ge 4$. Then $\pi^{\S}_{n,m}(K)\le (4E_{n-2,k})^{2\lfloor\frac{m-q-1}{2}\rfloor}$ if $q\le n-1$. Moreover, if $q=n$ we have
$$\pi_{n,m}(K)\le (4E_{n-2,k})^{\min\bigl(2\lfloor\frac{m-n-1}{2}\rfloor,\lfloor\frac{|K|^n-n-1}{2}\rfloor\bigr)}.$$ 

{\bf (2)} Suppose that $K=\mathbb F_3$. Then we have inequalities $\pi^{\S}_{n,m}(\mathbb F_3)\le (4n)^{2\lfloor\frac{m-q-1}{2}\rfloor}$ if $q\le n-1$ and $\pi_{n,m}(\mathbb F_3)\le (4n)^{\min\bigl(2\lfloor\frac{m-n-1}{2}\rfloor,\lfloor\frac{3^n-n-1}{2}\rfloor\bigr)}$ if $q=n$.

\smallskip
{\bf (3)} Suppose that $K=\mathbb F_2$ and $m\ge 3$ (thus $q\ge 2$). Then we have inequalities $\pi_{n,m}(\mathbb F_2)\le (2n-3)^{\lfloor\frac{m-q}{2}\rfloor}(n-1)^{m-q-2\lfloor\frac{m-q}{2}\rfloor}$ if $m\in \llbracket2^{q-1}+1,2^q-q\rrbracket$ and $\break\pi_{n,m}(\mathbb F_2)\le (2n-3)^{\lfloor\frac{m-q-1}{2}\rfloor}(n-1)^{m-q-1-2\lfloor\frac{m-q-1}{2}\rfloor}$ if $m\in \llbracket2^q-q+1,2^q\rrbracket$.\end{theorem}

\begin{proof}
Let $(\underline{P},\underline{Q})=\bigl((P_1,\ldots,P_m),(Q_1,\ldots,Q_m)\bigr)\in\mathbb D_{n,m}(K)^2$. Let $\sigma\in\perm(K^n)$ be such that $\sigma(\underline{P})=\underline{Q}$. If $|K|$ is odd, $q=n$, and $m=|K|^n-1$, then based on Theorem \ref{T6}(1) we can assume up to linear automorphisms that $\sigma$ is even. 

Therefore for both parts (1) and (2) we can assume that there exists a pair $(a,\theta)\in\AGL_n(K)\times\Alt(K^n)$ such that $\sigma=\theta a(K)$ and moreover for $l:=\nu_{3,1}(\theta)$ we have $l\le 2\lfloor\frac{m-q-1}{2}\rfloor$ if $q\le n-1$ and $l\le\min\bigl(2\lfloor\frac{m-n-1}{2}\rfloor,\lfloor\frac{|K|^n-n-1}{2}\rfloor\bigr)$ if $q=n$ by Proposition \ref{PR25}(3) and (4). 

We write $\theta=\prod_{i=1}^l \theta_i$ with each $\theta_i\in\perm(K^n)$ either a $3$-cycle or the identity permutation. Based on Lemma \ref{F7}(1), for every $i\in \llbracket1,l\rrbracket$ there exists $\break b_i\in\STGA_n(K)[4E_{n-2,k}]$ with $b_i(K)=\theta_i$ by Example \ref{EX8}(1) and Proposition \ref{PR10}(3). For the composite $c:=(\prod_{i=1}^l b_i)a\in\TGA_n(K)$ we have $\ell(c)\le [4E_{n-2,k}]^l$ and $c(K)=\sigma$, hence $c(\underline{P})=\underline{Q}$. If $q\le n-1$ then we can assume that $a\in\ASL_n(K)$, hence $c\in\STGA_n(K)$. So part (1) holds.

By replacing $4E_{n-2,k}$ with $4n$ in the last paragraph based on Proposition \ref{PR10}(2), we also get that part (2) holds.

For part (3), we consider $l\in\{m-q,m-q-1\}$ such that $l=m-q$ iff we have $m\in \llbracket2^{q-1}+1,2^q-q\rrbracket$. We choose $\sigma$ such that there exists an $(l+1)$-tuple $(a,\tau_1,\ldots,\tau_l)\in\AGL_n(K)\times\perm(K^n)^l$ with $\sigma=a(K)\prod_{i=1}^l \tau_i$ and for each $i\in \llbracket1,l\rrbracket$ the permutation $\tau_i$ is either a transposition or the identity element by Proposition \ref{PR25}(1). 

As $\ell_{n,K}(\sigma)\le (2n-3)^{\lfloor\frac{l}{2}\rfloor}(n-1)^{l-2\lfloor\frac{l}{2}\rfloor}$ by Proposition \ref{PR16}(5), there exists an automorphism $b\in\TGA_n(K)[(2n-3)^{\lfloor\frac{l}{2}\rfloor}(n-1)^{l-2\lfloor\frac{l}{2}\rfloor}]$ such that $b(K)=\sigma$, hence $b(\underline{P})=\underline{Q}$. So part (3) holds.\end{proof}

\section{Subexponential upper bounds for $n=2$ and all $m$}\label{S23}

In this section we combine the results on medium $m$ with Propositions \ref{PR2} and \ref{PR19}(1.a) to get upper bounds for the finite $\pi_{2,m}(K)$s when $K$ is finite and the integer $m\in \{|K|^2-2,|K|^2\}$ takes the largest possible value for $K$. The bounds are of the form $\pi_{2,m}(K)\le e^{(C_1\sqrt{m}+C_2)[(\ln(C_3m)]}$ with $C_1$, $C_2$, and $C_3$ universal positive constants. We begin with a lemma required in the proof. 

\begin{lemma}\label{L21}
Let $K$ be a finite field. Then the following properties hold.

\medskip
{\bf (1)} We have the inequality
$$\frac
{\lfloor\frac{|K|}{\sqrt{2}}\rfloor\bigl(\lfloor\frac{|K|}{\sqrt{2}}\rfloor-1\bigr)}{(|K|-1)^2}\le 0.515625$$ 
and the equality holds only for $|K|=17$.

\smallskip
{\bf (2)} We have the inequality $\frac
{\lfloor\frac{|K|}{\sqrt{2}}\rfloor\bigl(\sqrt{\lfloor\frac{|K|}{\sqrt{2}}\rfloor-1}\bigr)}{(|K|-1)^{\frac{3}{2}}}\le\frac{\sqrt{2}}{2}=0.7071...$ and the equality holds only for $|K|=3$.

\smallskip
{\bf (3)} For $|K|\ge 13$ we have the inequality $\frac
{\lfloor\frac{|K|}{\sqrt{2}}\rfloor\bigl(\sqrt{\lfloor\frac{|K|}{\sqrt{2}}\rfloor-1}\bigr)}{(|K|-1)^{\frac{3}{2}}}<\sqrt{0.4}<0.633$.
\end{lemma}

\begin{proof}
For $|K|\le 16$ one checks by direct computations that $\frac{\lfloor\frac{|K|}{\sqrt{2}}\rfloor\bigl(\lfloor\frac{|K|}{\sqrt{2}}\rfloor-1\bigr)}{(|K|-1)^2}\le 0.5$. 

For $|K|=17$ we have $\lfloor\frac{17}{\sqrt{2}}\rfloor=12$ and $\frac{12\cdot 11}{16^2}=0.515625$. For $|K|=19$ we have $\lfloor\frac{19}{\sqrt{2}}\rfloor=13$ and $\frac{13\cdot 12}{18^2}<0.482$. 

For $|K|\ge 23$, it suffices to show that $|K|(|K|-\sqrt{2})\le 1.03125|K|^2-2.0625|K|$ which holds as $\frac{2.0625-\sqrt{2}}{0.03125}<20.75<|K|$. So part (1) holds.

Direct computations give that $\frac
{\lfloor\frac{|K|}{\sqrt{2}}\rfloor\bigl(\sqrt{\lfloor\frac{|K|}{\sqrt{2}}\rfloor-1}\bigr)}{(|K|-1)^{\frac{3}{2}}}\le \frac{\sqrt{2}}{2}$ for $|K|\in\{2,3,4,5\}$ and that the equality holds if $|K|=3$. We show that $|K|^2(|K|-\sqrt{2})< \sqrt{2}(|K|-1)^3$ for $|K|\ge 7$; this is equivalent to $|K|^2[(\sqrt{2}-1)|K|-2\sqrt{2}]+\sqrt{2}(3|K|-1)>0$ and hence it holds as $6.828...=\frac{2\sqrt{2}}{\sqrt{2}-1}<|K|$. So part (2) holds.

For part (3), the case $|K|=13$ is checked by a direct computation. For $|K|\ge 16$, it suffices to show that $|K|^2(|K|-\sqrt{2})< \frac{4\sqrt{2}}{5}(|K|-1)^3$; this is equivalent to $|K|^2[(\frac{4\sqrt{2}}{5}-1)|K|-\frac{12\sqrt{2}}{5}+\sqrt{2}]+\frac{4\sqrt{2}}{5}(3|K|-1)>0$ and hence it holds as we have $15.071...=\frac{\frac{12\sqrt{2}}{5}-\sqrt{2}}{\frac{4\sqrt{2}}{5}-1}<|K|$.\end{proof}

\begin{theorem}\label{T12} Suppose $K$ is a finite field and $m\in\llbracket122, |K|^2\rrbracket$. If $4\mid |K|$, then we also assume that $m\le |K|^2-2$. Then the following properties hold.

\medskip
{\bf (1)} We have $\pi_{2,m}(K)\le (1.3865m)^{5(3\sqrt{m}+7)}=e^{(15\sqrt{m}+35)\ln (1.3865m)}$ and also, if $m\le |K|^2-2$, $\pi^{\S}_{2,m}(K)\le (1.3865m)^{5(3\sqrt{m}+7)}=e^{5(3\sqrt{m}+7)\ln (1.3865m)}$.

\smallskip
{\bf (2)} If $m\ge \frac{|K|^2}{2}$, then 
$$\pi_{2,m}(K)\le (1.3865m)^{5(3\sqrt{m}+32.8053)}=e^{5(3\sqrt{m}+32.8053)\ln (1.3865m)}$$
and, if $m\le |K|^2-2$, 
$$\pi^{\S}_{2,m}(K)\le (1.3865m)^{5(3\sqrt{m}+32.8053)}=e^{5(3\sqrt{m}+32.8053)\ln (1.3865m)}.$$

{\bf (3)} If $4\nmid |K|$ we have 
$$\ln \pi_{2,|K|^2}(K)<\Bigl(15\sqrt{2}|K|+64+\frac{15\sqrt{2}}{|K|-2\sqrt{2}}\Bigr)\ln0.8326(|K|-1)+\ln 0.5533$$ 
and if $4\mid |K|$ we have 
$$\ln \pi^{\S}_{2,|K|^2-2}(K)<\Bigl(15\sqrt{2}|K|+64+\frac{15\sqrt{2}}{|K|-2\sqrt{2}}\Bigr)\ln0.8326(|K|-1)+\ln 0.5533.$$
\end{theorem}

\begin{proof}
We can assume that $m>|K|$ by Proposition \ref{PR21}(2), that $|K|\ge 13$ by $|K|\ge\sqrt{122}>11$, and that either $m\le |K|^2-2$ or $m=|K|^2$ by Equation (\ref{EQ21}). 

Let $(\underline{P},\underline{Q})=\bigl((P_1,\ldots,P_m),(Q_1,\ldots,Q_m)\bigr)\in\mathbb D_{2,m}(K)^2$. Let $\sigma\in\perm(K^2)$ be such that $\sigma(P_i)=Q_i$ for all $i\in \llbracket1,m\rrbracket$. If $m\le |K|^2-2$ (resp.\ for each $m$) we can assume that $\sigma\in\Alt(K^2)$ with $\n(\sigma)\le \min(2m-6,|K|^2)$ (resp.\ $\n(\sigma)\le |K|^2-3$) by Proposition \ref{PR25}(4) applied to $n=q=2$ (resp.\ by Lemma \ref{L20}(1) and (2)). 

Let $r:=\bigl\lfloor \frac{\lceil\frac{|K|}{\sqrt{2}}\rceil}{3}\bigr\rfloor$. We have $r\geq \frac{\lceil\frac{|K|}{\sqrt{2}}\rceil}{3}-\frac{2}{3}> \frac{|K|-\sqrt{8}}{3\sqrt{2}}$. As $|K|\ge 13$ we have $\frac{|K|-\sqrt{8}}{3\sqrt{2}}>2$. Thus $r\ge 3$. 

Let $N_{\sigma}$ be the minimum of $0.51563^2(|K|-1)^4[(0.4)^2(|K|-1)^{10}]^{\wp_{r,|K|^2-\lceil\frac{|K|}{\sqrt{2}}\rceil}}$ and 
$[(0.4)^2(|K|-1)^{10}]^{\wp_{r,\n(\sigma)}}$.
There exists $a\in\SGA_2(K)[\lfloor N_{\sigma}\rfloor]$ such that $\sigma=a(K)$ by Corollary \ref{C17}(1) and (2) applied to $\varepsilon=0$, $s=3$, and $t\ge 1$, Lemma \ref{L21}(1) and (3) which for $N$ as in Corollary \ref{C17} gives $N(N-1)\le 0.515625(|K|-1)^2$ and $N^2(N-1)\le 0.4(|K|-1)^3$, and Proposition \ref{PR3}(4). We have 
$$N_{\sigma}< \min\Bigl([0.8326(|K|-1)]^{10\wp_{r,\n(\sigma)}},0.5533[0.8326(|K|-1)]^{10\bigl(\frac{2}{5}+\wp_{r,|K|^2-\lceil\frac{|K|}{\sqrt{2}}\rceil}\bigr)}\Bigr),$$
as $0.4<0.8326^5$ and $0.8326^{-4}\cdot 0.51563^2<0.5533$.

So part (1) holds if $N_{\sigma}\le e^{(15\sqrt{m}+35)\ln (1.3865m)}$. To prove the equivalent statement that $\ln N_{\sigma}\le (15\sqrt{m}+35)\ln (1.3865m)$ and the refinements of it that give parts (2) and (3), we consider two disjoint cases as follows.

{\bf Case 1: $m < \frac{|K|^2}{2}$.} Hence $|K|\in (\sqrt{2m},m)$ and $\n(\sigma)\le 2m-6$. We have $\wp_{r,\n(\sigma)}\le \lfloor \frac{\n(\sigma)-2}{2r}\rfloor+1$ by Proposition \ref{PR3}(4) and $\lfloor \frac{\n(\sigma)-2}{2r}\rfloor \le \lfloor \frac{m-4}{r}\rfloor\le \frac{m-4}{r}$. Thus $\wp_{r,\n(\sigma)}\le\frac{m-4}{r}+1<\frac{3\sqrt{2}(m-4)}{|K|-\sqrt{8}}+1$.

We consider the function $\mathfrak f:[\sqrt{2m},m)\rightarrow\mathbb R$ defined by the rule 
$$\mathfrak f(x):=30\sqrt{2}(m-4)\frac{\ln 0.8326(x-1)}{x-\sqrt{8}}+10\ln 0.8326(x-1).$$ 
We have $\ln N_{\sigma}\le \mathfrak f(|K|)$ as $\ln N_{\sigma}\le 10\left(\frac{m-4}{r}+1\right)\ln0.8326(|K|-1)$ implies
\begin{equation}\label{EQ36}
\ln N_{\sigma}<10\Bigl(\frac{3\sqrt{2}(m-4)}{|K|-\sqrt{8}}+1\Bigr)\ln 0.8326(|K|-1).
\end{equation}
As $m\ge 122$, we have $x\ge \sqrt{244}>15.6$.

We compute $\mathfrak f'(x)=30\sqrt{2}(m-4)\left(\frac{1}{(x-1)(x-\sqrt{8})}-\frac{\ln 0.8326(x-1)}{(x-\sqrt{8})^2}\right)+\frac{10}{x-1}$.
With the expression within big parentheses being clearly negative, by estimating $m$ as $|K|$ from below and by using $\frac{x-\sqrt{8}}{x-4}\le 1.101$ for $x>15.6$, we get 
$$\mathfrak f'(x)\le 30\sqrt{2}\left(x-4\right)\left[\frac{1}{(x-1)(x-\sqrt{8})}-\frac{\ln 0.8326(x-1)}{(x-\sqrt{8})^2}\right]+\frac{10}{x-1}$$
$$<30\sqrt{2}\left(x-4\right)\left[\frac{1+0.367}{(x-1)(x-\sqrt{8})}-\frac{\ln0.8326(x-1)}{(x-\sqrt{8})^2}\right]$$
$$=\frac{30\sqrt{2}\left(x-4\right)}{(x-\sqrt{8})^2}\left[1.367\frac{x-\sqrt{8}}{x-1}-\ln 0.8326(x-1)\right].$$
Hence $\mathfrak f'(x)<0$ if $x\geq 6$. So the absolute maximum of $\mathfrak f$ is $\mathfrak f(\sqrt{2m})$. 

Therefore 
\begin{equation*}
\frac{\ln N_{\sigma}}{\ln0.8326(\sqrt{2m}-1)}< 10\Bigl[\frac{3\sqrt{2}(m-4)}{\sqrt{2m}-\sqrt{8}}+1\Bigr]=5\Bigl(6\sqrt{m}+14\Bigr).
\end{equation*}
Thus $\ln N_{\sigma}<5(6\sqrt{m}+13.332)\ln (0.8326\sqrt{2m}).$
So, as $2\cdot 0.8326^2<1.3865$, we have $\ln N_{\sigma}<(15\sqrt{m}+35)\ln 1.3865m$. So part (1) holds in this case.

{\bf Case 2: $m\geq \frac{|K|^2}{2}$.} We use $N_{\sigma}\le 0.5533[0.8326(|K|-1)]^{10\bigl(\frac{2}{5}+\wp_{r,|K|^2-\lceil\frac{|K|}{\sqrt{2}}\rceil}\bigr)}$ as it gives better results in this case. We write $\lceil\frac{|K|}{\sqrt{2}}\rceil=3r+s$ with $s\in\{0,1,2\}$. We have $\digamma_{r,|K|^2-\lceil\frac{|K|}{\sqrt{2}}\rceil}=2r+2+t$ with $t\in\mathbb N$ by Proposition \ref{PR3}(2). We estimate
$$\wp_{r,|K|^2-\lceil\frac{|K|}{\sqrt{2}}\rceil}\le\Bigl\lfloor\frac{|K|^2-3r-s-2r-2-t}{2r}\Bigr\rfloor+2=\Bigl\lfloor\frac{|K|^2-r-s-2-t}{2r}\Bigr\rfloor.$$

We check that the inequality holds $r+s+t\ge 5$. As for $|K|\ge 23$ we have $r\ge 5$. we can assume that $|K|\in\{13,16,17,19\}$. If $|K|=13$, then $\lceil\frac{|K|}{\sqrt{2}}\rceil=10$, $r=3$, $s=1$ and, by Proposition \ref{PR3}(2), $t=1$ as $|K|^2-\lceil\frac{|K|}{\sqrt{2}}\rceil=159=4\cdot 39+r$ with $39=3\cdot 13$; so $r+s+t=5$. For $|K|=16$ we have $\lceil\frac{|K|}{\sqrt{2}}\rceil=12$, so $r=4$ and $s=0$, and $t=2r=8$ as $4\mid |K|^2-\lceil\frac{|K|}{\sqrt{2}}\rceil=244$; so $r+s+t=12$. If $|K|=17$, then $\lceil\frac{|K|}{\sqrt{2}}\rceil=13$ and thus $s=1$; so $r+s+t\ge 5$. If $|K|=19$, then $\lceil\frac{|K|}{\sqrt{2}}\rceil=14$ and thus $s=2$; so $r+s+t\ge 6$. So the inequality holds. 

Therefore we have
$$\wp_{r,|K|^2-\lceil\frac{|K|}{\sqrt{2}}\rceil}\le\Bigl\lfloor\frac{|K|^2-7}{2r}\Bigr\rfloor.$$

From this and $r>\frac{|K|-\sqrt{8}}{3\sqrt{2}}$, we get that 
$$10\Bigl[\wp_{r,|K|^2-\lceil\frac{|K|}{\sqrt{2}}\rceil}+\frac{2}{5}\Bigr]<10\Bigl[\frac{3\sqrt{2}(|K|^2-7)}{2(|K|-\sqrt{8})}+\frac{2}{5}\Bigr]=15\sqrt{2}|K|+64+\frac{15\sqrt{2}}{|K|-2\sqrt{2}}.$$
Hence
\begin{equation*}
\ln N_{\sigma}< \Bigl(15\sqrt{2}|K|+64+\frac{15\sqrt{2}}{|K|-2\sqrt{2}}\Bigr)\ln0.8326(|K|-1)+\ln 0.5533
\end{equation*}
and this implies that part (3) holds.

For $|K|=13$ we estimate $15\sqrt{2}|K|+64+\frac{15\sqrt{2}}{|K|-2\sqrt{2}}<299.44<2(15\sqrt{m})$ as $m\ge 122$. For $|K|\ge 16$ we have $\frac{15\sqrt{2}}{|K|-2\sqrt{2}}<1.6106$ and thus 
$$15\sqrt{2}|K|+64+\frac{15\sqrt{2}}{|K|-2\sqrt{2}}<2\Bigl(15\frac{|K|}{\sqrt{2}}+32.8053\Bigr)\le 2(15\sqrt{m}+32.8053).$$ 
as $|K|\le\sqrt{2m}$. Therefore $\ln N_{\sigma}<(15\sqrt{m}+32.8053)\ln (0.8326\sqrt{2m})^2$.
So, again as $2\cdot 0.8326^2<1.3865$, we get that $\ln N_{\sigma}<(15\sqrt{m}+32.8053)\ln (1.3865m)$. So part (2) holds and hence part (1) holds in this case as well.
\end{proof}

\section{Subexponential upper bounds for $n\ge 3$ and all $m$}\label{S24}

In this section we combine the results on medium $m$ with Propositions \ref{PR2} and \ref{PR19}(1.a) to get upper bounds for the finite $\pi_{n,m}(K)$s when $K$ is finite and the integer $m\in \{|K|^n-2,|K|^n\}$ takes the largest possible value for $K$ and $n\ge 3$. The bounds are either as in Section \ref{S23} or of the form 
$$\pi_{n,m}(K)\le^{C_1(2m)^{\frac{1}{n}}(C_2\ln 2m)^2 \bigl[\ln\bigl(\ln(2m)\bigr)+C_3\bigr]}$$ 
with $C_1$ and $C_2$ universal positive constants and $C_3$ a universal constant. 

\begin{theorem}\label{T13} Suppose that $n\in\mathbb N^{\ast}\setminus\{1,2\}$, $K$ is a finite field, and $m\in\llbracket2, |K|^n\rrbracket$. If $4\mid |K|$, then we also assume that $m\le |K|^n-2$. Then the following properties hold.

\medskip
{\bf (1)} If $|K|\ge 3$ and $m\ge 123$ (resp.\ and $m\in\llbracket121,|K|^n-2\rrbracket$), then we have $\pi_{n,m}(K)\le (8.4m)^{16(\sqrt{m}+3.75)}$ (resp.\ $\pi^{\S}_{n,m}(K)\le (8.4m)^{16(\sqrt{m}+3.75)}$).

\smallskip
{\bf (2)} If $|K|=2$, $n\le 9$, and $m\ge 5$, then $\pi_{n,m}(K)< (m-3)^{12\sqrt{m}+30}$.

\smallskip
{\bf (3)} If $|K|=2$ and $n\ge 10$, then $\pi_{n,m}(K)<m^{18.6545(\sqrt{m}+3.49394)}$.
\end{theorem}

\begin{proof}
Let $k:=|K|^n$. For part (1) we have $k\ge 125$.

Based on Corollary \ref{C19} we can assume that $\frac{\sqrt{k}}{\sqrt{2}}<m$. Based on Theorem \ref{T7}(2) we can also assume that $n<m-1$; so $m\ge n+2\ge 5$.

In this paragraph we assume that $|K|=2$. As $\kappa_{n,K}=n-1$ by Proposition \ref{PR16}(4), from Lemma \ref{L1}(1) we get that $\pi_{n,m}(K)\le (n-1)^m\le (m-3)^m$ and from Lemmas \ref{L20}(1) and \ref{P1}(2) we get that 
$$\pi_{n,m}(K)\le\pi_{n,2^n}(K)\le (n-1)^{2^n-n-1}.$$ 
If $m\le 199$, then $m<12\sqrt{m}+30$ and hence parts (2) and (3) hold as we have $\pi_{n,m}(K)\le (n-1)^m\le (m-3)^{12\sqrt{m}+30}$. If $m\ge 200$ and $n\le 9$, then part (2) holds as we have
$$\pi_{n,m}(K)\le\pi_{9,m}(K)\le 8^{502}<8^{505}<8^{\frac{199\ln(196)}{\ln(8)}}=196^{199}<(m-3)^{12(\sqrt{m}+9)}.$$ 
If $m\le 469$, then $m<18.65(\sqrt{m}+3.493)$, and therefore part (3) holds as we have inequalities $\pi_{n,m}(K)\le (n-1)^m\le (m-3)^{m}<m^{18.65(\sqrt{m}+3.493)}$.

We are left to prove part (1) for $m\ge 121$ (resp.\ $m\ge 119$) and part (3) for $m\ge 470$.

Let $(\underline{P},\underline{Q})=\bigl((P_1,\ldots,P_m),(Q_1,\ldots,Q_m)\bigr)\in\mathbb D_{n,m}(K)^2$. Let $\sigma\in\Perm(K^n)$ be such that $\sigma(P_i)=Q_i$ for all $i\in \llbracket1,m\rrbracket$; it is not unique and we define $l\in \llbracket1,k\rrbracket$ and $n_0\in\{1,n-1\}$ and we choose $\sigma$ depending on $m$ and $|K|$ as follows.

Assume first that $m\ge \frac{k}{2}$.  We choose $l=l_1\in\{k-n-1,k-n\}$ such that we have $l=k-n$ iff $|K|$ is odd and we are bounding $\pi^{\S}_{n,m}(K)$. We choose $\sigma$ such that there exists a triple $(\theta,a,b)\in\Alt(K^n)\times\AGL_n(K)\times\STGA_n(K)[n_0]$ with $\sigma=\theta a(K)b(K)$ and $\n(\theta)\le l$ as follows. If $|K|$ is either odd or such that $4\mid |K|$, then we can assume that $b=1_{\mathbb A^n_K}$, $n_0=1$, and, if $K$ is odd and we are bounding $\pi^{\S}_{n,m}(K)$, $a\in\ASL_n(K)$ by Lemma \ref{L20}(1), (2), and (5). If $|K|=2$, then we can assume that $n_0=n-1$ by Lemma \ref{L20}(4).

If $m<\frac{k}{2}$, then we choose the values $n_0:=1$, $l=l_2:=2m-6$ for part (1), and $l=l_2:=2m-2\lfloor\frac{n+1}{2}\rfloor\le 2m-10$ for part (3) and we can assume that $\sigma$ is even with $\n(\sigma)\le l$ by Proposition \ref{PR25}(4) applied to $q=2$ for part (1) and by Proposition \ref{PR25}(1) applied to $q=\lfloor\frac{n-1}{2}\rfloor$ for part (3); the values of $q$ are such that $\frac{\sqrt{k}}{\sqrt{2}}<m$ implies $|K|^{q-1}+1<m$.\footnote{If $l_1<l_2$, then one can work with the smaller $l$ as in the previous paragraph and this would provide a small improvement in many cases such as when $|K|\ge 3$.}

Let $r:=\Bigl\lfloor \frac{\lfloor\frac{\sqrt{k}}{\sqrt{2}}\rfloor}{3}\Bigr\rfloor=\lfloor \frac{\sqrt{k}}{3\sqrt{2}}\rfloor$; so $r > \frac{\sqrt{k}}{3\sqrt{2}}-1=\frac{\sqrt{k}-3\sqrt{2}}{3\sqrt{2}}$ and for $k=125$ we have $r=2$ and for $k>125$, (i.e., $k\ge 243$) we have $r\ge 4$.

Let $D_{\sigma}\in [1,\infty)$ be such that there exists $c\in\STGA_n(K)[\lfloor e^{D_{\sigma}}\rfloor]$ with $\sigma=c(K)$. 

Let $s:=\lfloor\frac{n}{2}\rfloor$; if $n$ is even, then $s\in\mathbb N^{\ast}\setminus\{1\}$ and if $n$ is odd, then $s\in\mathbb N^{\ast}$. If $|K|=2$, let $\epsilon:=0.3$. If $|K|\ge 3$, let $\epsilon:=0$. If $|K|=3$ let $\varepsilon:=0$ and if $|K|\ge 4$ let $\varepsilon:=\frac{1}{5}$. We have inequalities $s(|K|-1)+\lfloor\sqrt{\frac{|K|}{2}}\rfloor\le (|K|-1)\bigl(s+\frac{1}{2}\bigr)$ and $s(|K|-1)+\lfloor\frac{|K|}{\sqrt{2}}\rfloor\le (|K|-1)(s+1-\varepsilon)$. Thus for $|K|\ge 3$ we have an inequality
$$\Bigl[s(|K|-1)+\Bigl\lfloor\sqrt{\frac{|K|}{2}}\Bigr\rfloor\Bigr]^2\Bigl[s(|K|-1)+\Bigl\lfloor\frac{|K|}{\sqrt{2}}\Bigr\rfloor\Bigr]^4\le (|K|-1)^6\bigl(s+\frac{1}{2}\bigr)^2(s+1-\varepsilon)^4.$$

For $n=2s+1$ odd, we write $4s^6=\frac{(n-1)^6}{16}$ and
$$|K|^2(|K|-1)^6\bigl(s+\frac{1}{2}\bigr)^2(s+1-\varepsilon)^4=\frac{1}{64}|K|^2(|K|-1)^6(n)^2(n+1-2\varepsilon)^4.$$ 

Let $N:=\frac{13}{16}$ if $|K|\ge 3$ and $N:=1$ if $|K|=2$. We have an inequality $E_{n-2,|K|-1}\le N(n-1)^{1+\epsilon}(|K|-1)^4$ by Lemma \ref{L18}(1) and (2) applied to $l=n-1$. 

For $|K|\ge 3$ and $d\in \llbracket1,k\rrbracket$ we define
$$D_{n,|K|,d}:=\Bigl\lfloor\frac{d+2r-2}{2r}\Bigr\rfloor\ln \Bigl[\frac{N}{2^{12}}|K|^4(|K|-1)^{20}(n-1)^{1+\epsilon}(n)^4(n+1-2\varepsilon)^8\Bigr]$$
and 
\begin{equation*}
\begin{aligned}
D=D_{n,|K|}:&=D_{n,|K|,k-3r}+\ln\Bigl[\frac{1}{64}|K|^2(|K|-1)^{8}(n)^{2}(n+1-2\varepsilon)^4\Bigr]\\
&=D_{n,|K|,k-r}-\frac{1}{2}\ln[N(n-1)^{1+\epsilon}(|K|-1)^4].
\end{aligned}
\end{equation*}

For $|K|=2$ and $d\in \llbracket1,k\rrbracket$ we also define
$$D_{n,2,d}:=\Bigl\lfloor\frac{d+2r-2}{2r}\Bigr\rfloor\ln \Bigl[\frac{N}{256}(n-1)^{13+\epsilon}\Bigr]$$
and 
\begin{equation*}
\begin{aligned}
D=D_{n,2}:&=D_{n,|K|,k-3r}+\ln\Bigl[\frac{1}{16}(n-1)^6\Bigr]\\
&=D_{n,2,k-r}-\frac{1}{2}\ln[N(n-1)^{1+\epsilon}].
\end{aligned}
\end{equation*}

If $|K|=2$ we also define 
$$D'_{n,2,d}:=2\Bigl\lceil\frac{d}{2^{n-2}}\Bigr\rceil[2D^{-}_{n,1}+\ln 2]+\ln(n-1)-\ln(2)$$
and 
$$D':=5\ln(2)+\ln(n-1)+13\ln(D^{-}_{n,1}).$$

We apply Corollary \ref{C22}(2) to $\sigma$ if $m<\frac{k}{2}$ and to $\theta$ if $m\ge\frac{k}{2}$; for $s\ge 2$, the case $n=2s+1$ gives a greater $D_{\sigma}$ than the case $n=2s$. So for $|K|\ge 3$ we take 
$$D_{\sigma}:=\min(D_{n,|K|,l},D).$$
For $|K|=2$ we also apply Corollary \ref{C21}(1) and (2) to $\sigma$ if $m<\frac{k}{2}$ and Corollary \ref{C21}(1) to $\theta$ if $m\ge\frac{k}{2}$; so we can take
$$D_{\sigma}:=\min\bigl(D_{n,2,l}+\ln(n_0),D'_{n,2,l},D+\ln(n_0),D'\bigr)$$
and we recall that $\ln(n_0)=0$ if $m< \frac{k}{2}$.

We have reduced the problem to bounding $D_{\sigma}$ from above in terms of $m$ when $r=\lfloor \frac{\sqrt{k}}{3\sqrt{2}}\rfloor$ and $m\in \llbracket\lfloor\frac{\sqrt{k}}{\sqrt{2}} \rfloor+1, k\rrbracket$, with $n$ odd if either $|K|\ge 3$ or $|K|=2$ and for $n$ odd and $m<\frac{k}{2}$ (resp.\ $m\ge\frac{k}{2}$) we redefine $l=l_2$ (resp.\ $l=l_1$) as $2m-2\lfloor\frac{n}{2}\rfloor$ (resp.\ $k-n$).

First, we consider for $|K|\ge 3$ and hence $\epsilon=0$ the `$\ln$ expression' 
$$\frac{N}{2^{12}}|K|^4(|K|-1)^{20}(n-1)(n)^4(n+1-2\varepsilon)^8$$ 
that depends on $|K|$ and $n$, for any given $k$. We keep in mind that $n\geq 3$ and $|K|\geq 3,$ but we ignore the fact that $|K|$ and $n$ are integers. Rewriting the expression in terms of $|K|$ and $\ln k$, i.e., by eliminating $n=\frac{\ln(k)}{\ln(|K|)}$ we get a function of $x=|K|$, $\mathfrak f:[3,\infty)\rightarrow\mathbb R$, defined by the rule
$$\mathfrak f(x):=\frac{N}{2^{12}}x^4(x-1)^{20}\Bigl(\frac{\ln k}{\ln x}-1\Bigr)\Bigl(\frac{\ln k}{\ln x}\Bigr)^{4}\Bigl(\frac{\ln k}{\ln x}+1-2\varepsilon\Bigr)^8.$$ 

Then
$$\frac{\mathfrak f'(x)}{\mathfrak f(x)}=\frac{4}{x}+\frac{20}{x-1}-\frac{\frac{\ln k}{x(\ln x)^2}}{\frac{\ln k}{\ln x}-1}-\frac{\frac{4\ln k}{x(\ln x)^2}}{\frac{\ln k}{\ln x}}-\frac{\frac{8\ln k}{x(\ln x)^2}}
{\frac{\ln k}{\ln x}+1-2\varepsilon}.$$
By choosing the smaller positive and the larger negative terms, we get that $\frac{\mathfrak f'(x)}{\mathfrak f(x)}$ is strictly greater than
$$\frac{24}{x}-\frac{13\ln k}{x\ln x [\ln k-\ln x]}= \frac{24[\ln k-\ln x]\ln x -13\ln k}{x\ln x[\ln k-\ln x]}.$$

This bound for $\frac{\mathfrak f'(x)}{\mathfrak f(x)}$ is positive iff $24[\ln k-\ln x]\ln x-13\ln k>0$. As $\ln x <\ln k,$ this difference is the smallest when either $\ln x$ is the smallest, $\ln x=3$, or the largest, $\ln x=\frac{\ln(k)}{3}$. Substituting $\ln x =\frac{\ln k}{3}$ in $24[\ln k-\ln x]\ln x -13\ln k$ we obtain $\frac{\ln k}{3}(-39+\frac{48}{3}\ln k)>0$ based on the estimates $\ln(k)>\ln(12)>\frac{39\cdot 3}{48}=2.4375$. Also, $\frac{24\ln(3)^2}{24\ln(3)-13}<2.168$ gives $24[\ln k-\ln 4]\ln 4-13\ln k>0$. 

From the last two paragraphs we get that $\mathfrak f(x)$ is the largest when $x$ is the largest possible as a function of (fixed) $k$, which means $n=3$ and $|K|=k^{\frac{1}{3}}$.

Therefore, as for $n=3$ we have $n+1-2\varepsilon\le 4$ and we ignore that for $|K|\ge 4$ we can replace $4$ by $3.6$, for $|K|\ge 3$ we have 
$$D_{n,|K|,l}\le\Bigl(\Bigl\lfloor\frac{l-2r-2}{2r}\Bigr\rfloor+2\Bigr)\left[\ln N+5\ln 2+4\ln 3+4\ln |K|+20 \ln (|K|-1)\right]$$
$$=\Bigl(\Bigl\lfloor\frac{l-2}{2r}\Bigr\rfloor+1\Bigr)\bigl[\ln N +5\ln 2+4\ln3+\frac{4}{3}\ln k+20\ln (k^{\frac{1}{3}}-1)\bigr].$$

For $|K|\ge 3$, an entirely similar argument gives first, when the role of the pair (24,13) is replaced by $(20,6)$, that 
$$\ln\Bigl[\frac{1}{64}|K|^2(|K|-1)^{8}(n)^{2}(n+1-2\varepsilon)^4\Bigr]<2\ln6+\frac{2}{3}\ln k+8\ln (k^{\frac{1}{3}}-1)$$
and second that $D$ is less than 
$$\Bigl(\frac{k-3r-2}{2r}+2\Bigr)\Bigl[\ln N +5\ln 2+4\ln3+\frac{4}{3}\ln k+12 \ln (k^{\frac{1}{3}}-1)\Bigr]-\frac{1}{2}\ln 2-2\ln (k^{\frac{1}{3}}-1)$$
$$=\Bigl(\frac{k-2}{2r}+0.5\Bigr)\Bigl[\ln N +5\ln 2+4\ln3+\frac{4}{3}\ln k+12\ln (k^{\frac{1}{3}}-1)\Bigr]-\frac{1}{2}\ln 2-2\ln (k^{\frac{1}{3}}-1).$$

As in the proof of Theorem \ref{T12}, we consider two disjoint cases.

{\bf Case 1: $\frac{k}{2} \le m \le k.$} We first consider the subcase $|K|\ge 3$. We only use $D$ and not $D_{n,|K|,l}$ so the actual value of $l\in\{k-n-1,k-n\}$ is not used.

With $r> \frac{\sqrt{k}-3\sqrt{2}}{3\sqrt{2}}$ and $k\le 2m$, we get that $D$ is less than
$$\left(\frac{\frac{3}{2}(k-2)}{\sqrt{\frac{k}{2}}-3}+0.5\right)\Bigl[\ln N+5\ln 2+4\ln 3+\frac{4}{3}\ln (2m)+12\ln \bigl((2m)^{\frac{1}{3}}-1\bigr)\Bigr]-\frac{2}{3}\ln (2m-1).$$
As $k\le 2m$ and $\frac{\frac{3}{2}(k-2)}{\sqrt{\frac{k}{2}}-3}$ is an increasing function on $k$, we get that $D$ is less than
\begin{equation}\label{EQ37}
\Bigl(\frac{3(m-1)}{\sqrt{m}-3}+0.5\Bigr)\Bigl[\ln N+5\ln 2+ \frac{4}{3}\ln (2m)+12\ln \bigl((2m)^{\frac{1}{3}}-1\bigr)\Bigr]-\frac{2}{3}\ln (2m-1).
\end{equation}

We rewrite $\frac{3(m-1)}{\sqrt{m}-3}+0.5$ as $3\sqrt{m}+9.5+\frac{24}{\sqrt{m}-3}.$ For $m\geq 121$, this does not exceed $3\sqrt{m}+12.5$. We have $5\ln(2)+4\ln(3)+\ln(N)<7.6526$ and thus the factor within brackets of Display (\ref{EQ37}) is less than 
$$\frac{16}{3}\ln (2m)+7.6526<\frac{16}{3}[\ln (2m)+1.4349]<\frac{16}{3}[\ln(2m)+\ln(4.2)]=\frac{16}{3}\ln(8.4m).$$
We get that
$$D< 16\Bigr(\sqrt{m}+\frac{12.5}{3}\Bigl)\ln(8.4m)<16(\sqrt{m}+3.75)\ln(8.4m)-\frac{2}{3}\ln (2m-1).$$

Next we consider the subcase $|K|=2$. We have $l\le k-11$ regardless if we restrict to an odd $n\ge 11$ or we work with $n\ge 10$; so $l-2\le k-13$. We have
$$D+\ln(n_0)<\left[13.3\left(\Bigl\lfloor\frac{k-13}{2r}\Bigr\rfloor+0.5\right)+0.7\right]\ln(n-1);$$
thus
\begin{equation}\label{EQ38}
D+\ln(n_0)<\Bigl(13.3\Bigl\lfloor\frac{k-13}{2r}\Bigr\rfloor+7.35\Bigr)\ln(n-1)=\Bigl(13.3\Bigl\lfloor\frac{k-13}{2r}\Bigr\rfloor+7.35\Bigr)\ln\left(\frac{\ln(k)}{\ln 2}-1\right).\end{equation}

Inequality (\ref{EQ38}) and the inequalities $r > \frac{\sqrt{k}-3\sqrt{2}}{3\sqrt{2}}$ and $k\le 2m$, give, as in the previous subcase, that $D$ is less than 
$$\left(\frac{39.9(m-6.5)}{\sqrt{m}-3}+7.35\right)\ln\left(\frac{\ln(m)}{\ln 2}\right)=\Bigl(39.9\sqrt{m}+127.5+\frac{99.75}{\sqrt{m}-3}\Bigr)\ln\left(\frac{\ln(m)}{\ln 2}\right).$$
For $m\ge 470$, this is less than $(39.9\sqrt{m}+132.8401)\ln\left(\frac{\ln(m)}{\ln 2}\right)$ and hence than
$$39.9(\sqrt{m}+3.32933)\ln\left(\frac{\ln(m)}{\ln 2}\right)<14.1595(\sqrt{m}+3.32933)\ln\left(\frac{\ln(m)}{\ln 2}\right)^{2.8179}.$$ 
As for $m\ge 470$ we have $\left(\frac{\ln(m)}{\ln 2}\right)^{2.8179}<m$, we conclude that for $m\ge 470$ and $|K|=2$ we have $D<14.1595(\sqrt{m}+3.32933)\ln(m)$.

{\bf Case 2: $\lfloor\frac{\sqrt{k}}{\sqrt{2}}\rfloor <m < \frac{k}{2}.$} Then $l\le 2m-6$ for $|K|\ge 3$ and $l\le 2m-10$ if $|K|=2$ regardless if we restrict to an odd $n\ge 11$ or we work with $n\ge 10$. We first consider the subcase $|K|\ge 3$. 

As $D_{n,|K|,l}<\left(\frac{3\sqrt{2}(m-4)}{\sqrt{k}-3\sqrt{2}}+1\right)\bigl[\ln N+5\ln 2+4\ln(3)+ \frac{4}{3}\ln k+12\ln (k^{\frac{1}{3}}-1)\bigr]$, we also have $\left(\frac{3\sqrt{2}(m-4)}{\sqrt{k}-3\sqrt{2}}+1\right)(\frac{16}{3}\ln k +7.6526)$. Denoting $y:=\sqrt{k}$, this can be rewritten as $(\frac{3\sqrt{2}(m-4)}{y-3\sqrt{2}}+1)(\frac{32}{3}\ln y + 7.6526)$. The inequalities $\lfloor\frac{\sqrt{k}}{\sqrt{2}} \rfloor <m < \frac{k}{2}$ translate for $y$ into $\sqrt{2m}<y<\sqrt{2}m.$ 

We show that for every fixed $m\geq 28$, the maximum of the real-valued function $\mathfrak f:[\sqrt{2m}, \sqrt{2}m]\rightarrow\mathbb R$ defined by the rule $$\mathfrak f(y):=\Bigl[\frac{3\sqrt{2}(m-4)}{y-3\sqrt{2}}+1\Bigr]\Bigl(\frac{32}{3}\ln y + 7.6526\Bigr)$$ is obtained at $y=\sqrt{2m}$. It suffices to show that $\mathfrak f'(y)<0$ for $y\in [\sqrt{2m}, \sqrt{2}m]$. This is equivalent to
$$3\sqrt{2}(m-4)\Bigl(\frac{32}{3}\ln y + 7.6526\Bigr)> \frac{32}{3y}\bigl[3\sqrt{2}m+y-15\sqrt{2}\bigr](y-3\sqrt{2}).$$
As $y\leq \sqrt{2}m$, the right hand side is bounded from above by $\frac{32}{3}\cdot 4\sqrt{2}(m-3.75)$ and thus by $44\sqrt{2}(m-4)$. As $y\geq \sqrt{2m}>e^2,$ the left hand side is clearly larger. Thus 
$$D_{n,|K|,l}\le \Bigl(\frac{3\sqrt{2}(m-4)}{\sqrt{2m}-3\sqrt{2}}+1\Bigr)\Bigl(\frac{16}{3}\ln \sqrt{2}m +7.6526\Bigr).$$

As above we conclude that
$$D_{\sigma}\le D_{n,|K|,l}<(3\sqrt{m}+11.875)\cdot \frac{16}{3}\ln(8.4m)=16(\sqrt{m}+3.9584)\ln(5.94m).$$

Next we consider the subcase $|K|=2$. As $l\le 2m-10$ we have
$$D_{n,2,l}\le\Bigl[13.3\Bigl(\Bigl\lfloor\frac{2m-12}{2r}\Bigr\rfloor+1\Bigr)\Bigr]\ln(n-1).$$
So, as $r> \frac{\sqrt{\frac{k}{2}}-3}{3}$, $n+1=\frac{\ln(2k)}{\ln 2}$, and $\frac{\sqrt{k}}{\sqrt{2}}\le m$ implies $2^{n-1}=\frac{k}{2}<m^2$, we get
$$D_{\sigma}<13.3\Bigl(\Bigl\lfloor\frac{m-6}{r}\Bigr\rfloor+1\Bigr)\ln(n-1)<13.3\left(\frac{3(m-6)}{\sqrt{m}-3}+1\right)\ln\left(\frac{2\ln m}{\ln 2}\right)$$
$$=13.3\Bigl(3\sqrt{m}+10+\frac{9}{\sqrt{m}-3}\Bigr)\ln\left(\frac{2\ln m}{\ln 2}\right).$$
For $m\ge 470$, this is less than
$13.3(3\sqrt{m}+10.48182)\ln\left(\frac{2\ln m}{\ln 2}\right)$ and hence than
$$39.9(\sqrt{m}+3.49394)\ln\left(\frac{2\ln m}{\ln 2}\right)<18.6545(\sqrt{m}+3.49394)\ln\left(\frac{2\ln m}{\ln 2}\right)^{2.1389}.$$ 
As for $m\ge 470$ we have $\left(\frac{2\ln m}{\ln 2}\right)^{2.1389}<m$, we conclude that for $m\ge 470$ and $|K|=2$ we have $D_{\sigma}<18.6545(\sqrt{m}+3.49394)\ln(m)$.

Based on Cases 1 and 2, we get that parts (1) and (3) hold.\end{proof}

\begin{theorem}\label{T14} Suppose that $n\in\mathbb N^{\ast}\setminus\{1,2\}$, $K$ is a finite field, $N\in\{513,350,626\}$ with $N=626$ iff $|K|\ge 4$ and $N=350$ iff $|K|=3$, and $m\in\llbracket N,|K|^n\rrbracket$. If $4\mid |K|$, then we also assume that $m\le |K|^n-2$. Then the following properties hold.

\medskip
{\bf (1)} If $|K|\geq 4$, then
$$\ln\bigl(\pi_{n,m}(K)\bigr)\le
9.0169(2m)^{\frac{1}{n}}[\ln(2m)]^2[0.7214\ln\bigl(\ln(2m)\bigr)+1.278]$$
and, if $m\le |K|^n-2$, 
$$\ln\bigl(\pi^{\S}_{n,m}(K)\bigr)\le
9.0169(2m)^{\frac{1}{n}}[\ln(2m)]^2[0.7214\ln\bigl(\ln(2m)\bigr)+1.278].$$

{\bf (2)} If $|K|=3$, then
$\ln\bigl(\pi_{n,m}(K)\bigr)\le
10.495(2m)^{\frac{1}{n}}[\ln(2m)]^2[\ln\bigl(\ln(2m)\bigr)+0.7922]$
and,  if $m\le |K|^n-2$, $\ln\bigl(\pi^{\S}_{n,m}(K)\bigr)\le
10.495(2m)^{\frac{1}{n}}[\ln(2m)]^2[\ln\bigl(\ln(2m)\bigr)+0.7922]$.

\smallskip
{\bf (3)} If $|K|=2$, then
$\ln\bigl(\pi_{n,m}(K)\bigr)\le
27.0578(\frac{4m}{3})^{\frac{1}{n}}[\ln(\frac{4m}{3})]^2[\ln\bigl(\ln(\frac{4m}{3})\bigr)-0.673].$
\end{theorem}

\begin{proof}
Based on Corollary \ref{C19} we can assume that $\frac{|K|^{\frac{n}{2}}}{\sqrt{2}}<m$. Thus $n<\frac{\ln(2m^2)}{\ln |K|}$. Also, we can assume that $m\le |K|^n-1$ by Equation (\ref{EQ21}). If $|K|=2$, then the hypotheses $m\ge N=513$ gives $n\ge 10$.

We define the triple $(C_1,C_2,C_3)$ to be $(9.0169,0.7214,1.278)$ if $|K|\ge 4$, to be $(10.495,1,0.7922)$ if $|K|=3$, and to be $(27.0578,1,-0.673)$ if $|K|=2$. 

We consider three disjoint cases as follows.

{\bf Case 1: $\frac{|K|^{\frac{n}{2}}}{\sqrt{2}}<m \le 2\lceil
\frac{|K|^{n-1}}{4}\rceil.$} Then by Corollary \ref{C20},
$\ln\bigl(\pi^{\S}_{n,m}(K)\bigr)$ is bounded from above by $\ln D_{n,k}$ and hence by $\ln C_{n,k}^+$ of Lemma \ref{L16}(1). Thus, if $|K|\geq 3,$ then we have
$$\ln\bigl(\pi^{\S}_{n,m}(K)\bigr)<(2n^2+4n-9)\ln (|K|) + (2n^2+n-2)\ln (n-1)+\min\Bigl(2\ln(n-2)+\frac{1}{3},n-1\Bigr).$$
As $\ln|K| >1$, $n\geq 3$, $\ln n\ge 1$, and $m\ge 3$, we
have $n<\ln(2m^2)$ and
$$\ln\bigl(\pi^{\S}_{n,m}(K)\bigr) < \frac{10}{3}n^2\ln|K|+\frac{8}{3}n^2 \ln n< \frac{10}{3}\bigl[\ln
(2m^2)\bigr]^2\bigl[1 +\ln\bigl(\ln (2m^2)\bigr)\bigr]$$
and hence $\ln\bigl(\pi^{\S}_{n,m}(K)\bigr)<\frac{10}{3}[\ln(2m^2)]^3$.

Similarly, if $|K|=2$, then based on Lemma \ref{L16}(2) and $n\ge 10$ we get
$$\ln\bigl(\pi_{n,m}(K)\bigr)< (2n^2+n-2)\ln(n)+2n<(n^2\ln
n)\Bigl(\frac{21}{20}+\frac{2}{n\ln n}\Bigr)$$
$$< 2.9402[\ln
(2m^2)]^2\ln\bigl(\ln(2m^2)\bigr).$$

It is easy to see that, regardless of what $|K|$ is, these upper bounds are smaller than the required ones. So parts (1) to (3) hold in this case.

{\bf Case 2: $|K|=2$ and $2^{n-2}< m\le 2^n-1.$} Recall that $\ln D_{n,1}$ is, by Lemma \ref{L16}(2) and its definition in Corollary \ref{C20}, less than or equal to 
$$E(n,2):=n+\frac{1}{3}+(n^2+n-2)\ln(n-2)+(n^2+4)\ln(n-1).$$

Let $(\underline{P},\underline{Q})=\bigl((P_1,\ldots,P_m),(Q_1,\ldots,Q_m)\bigr)\in\mathbb D_{n,m}(K)$. Let $a\in\TGA_n(K)[D_{n,1}]$ be such that $a(P_i)=Q_i$ for each $i\in \llbracket1,2^{n-2}\rrbracket$ by Corollary \ref{C20}. Let $\sigma\in\perm(K^n)$ be such that $\sigma\bigl(a(P_i)\bigr)=Q_i$ for each $i\in\llbracket1,m\rrbracket$. So $\supp(\sigma)\subset K^n\setminus\{Q_1,\ldots,Q_{2^{n-2}}\}$ and based on Lemma \ref{L1}(1) we can assume that $\nu_{2,1}(\sigma)\le m-2^{n-2}\le 3\cdot 2^{n-2}-1$. If $m$ is even with $m\le 2^n-2$, then we assume that $\sigma$ is even by Lemma \ref{L1}(2). 

Let $b\in\TGA_n(K)[\ell_{n,K}(\sigma)]$ be such that $b(K)=\sigma$. For the automorphism $c:=ba\in\TGA_n(K)[D_{n,1}\ell_{n,K}(\sigma)]$ we have $c(\underline{P})=\underline{Q}$. 

Let $\mathfrak f:[\frac{2^n+4}{3},2^n+\frac{4}{3}]\rightarrow\mathbb R$ be the real-valued function defined by the rule 
$$\mathfrak f(y):=C_1(y)^{\frac{1}{n}}[\ln(y)]^2\bigl[C_2\ln\bigl(\ln(y)\bigr)+C_3\bigr],$$
with $y$ thought as $\frac{4m}{3}$.

If $m\in\llbracket2^{n-2}+1,2^{n-1}\rrbracket$, then either $m\le 2^{n-1}-1$ or $m=2^{n-1}$ and $\sigma$ is even; thus $\ell_{n,\sigma}(K)\le 2D_{n,1}^4$ by Corollary \ref{C21}(1) and (2) and therefore we have $c\in\TGA_n(K)[2(n-1)D_{n,1}^5]$. For $m\in\llbracket2^{n-2}+1,2^{n-1}\rrbracket$ we want to prove that the following inequality $\ln(2)+\ln(n-1)+5\ln(D_{n,1})\le\mathfrak f(\frac{2^n+4}{3})$ holds.

If $m\in\llbracket2^{n-1}+1,3\cdot 2^{n-2}\rrbracket$, then $\ell_{n,\sigma}(K)\le 8(n-1)D_{n,1}^8$ by Corollary \ref{C21}(1) and (2) and thus $c\in\TGA_n(K)[8(n-1)D_{n,1}^9]$; we want to prove that the following inequality 
$3\ln(2)+\ln(n-1)+9\ln(D_{n,1})\le\mathfrak f(\frac{2^{n+1}+4}{3})$ holds.

If $m\in\llbracket3\cdot 2^{n-2}+1,2^{n}-1\rrbracket$, then we have $\ell_{n,\sigma}(K)\le 32(n-1)D_{n,1}^{12}$ by Corollary \ref{C21}(1) and (2). Thus $c\in\TGA_n(K)[32(n-1)D_{n,1}^{13}]$ and we want to prove that the following inequality 
$5\ln(2)+\ln(n-1)+13\ln(D_{n,1})\le\mathfrak f(2^n+\frac{4}{3})$ holds. 

Let $\mathfrak f_1:[\frac{2^n+4}{3},2^n+\frac{4}{3}]\rightarrow\mathbb R$ be the real-valued function defined by the rule 
$$\mathfrak f_1(y):=\ln\frac{(n-1)E(n,2)}{2}+\frac{y-\frac{4}{3}}{\frac{2^{n-2}}{3}}\ln\bigl(\sqrt{2}E(n,2)\bigr)$$
$$=\Bigl(\frac{3y-4}{2^{n-2}}-1\Bigr)\ln(2)+\ln(n-1)+\Bigl(1+\frac{3y-4}{2^{n-2}}\Bigr)\ln\bigl(E(n,2)\bigr).$$ 
Let $d_n:=\frac{\frac{2^{n-2}}{3}\ln\frac{(n-1)E(n,2)}{2}}{\ln\bigl(\sqrt{2}E(n,2)\bigr)}-\frac{4}{3}$. The inequality $d_n< \frac{2^{n}-4}{3}$ is equivalent to the inequality $n-1<8E(n,2)^3$ and thus it holds as $E(n,2)\ge 2n$.

As $D_{n,1}\le E(n,2)$, from the second expression of $\mathfrak f_1$ we get that we have inequalities $\ln(2n-2)+5\ln(D_{n,1})\le\mathfrak f_1(\frac{2^n+4}{3})$, $\ln(8n-8)+9\ln(D_{n,1})\le\mathfrak f_1(\frac{2^{n+1}+4}{3})$, and $\ln(32n-32)+13\ln(D_{n,1})\le\mathfrak f_1(2^n+\frac{4}{3})$. Thus to prove the desired three inequalities above it suffices to prove that we have $\ln\bigl(\mathfrak f_1(y)\bigr)-\ln\bigl(\mathfrak f(y)\bigr)\le 0$ for each $y\in [\frac{2^n+4}{3},2^n+\frac{4}{3}]$. 

\phantomsection{We use the substitution $z:=\ln(y)$. Taking the derivative of the difference $\ln\bigl(\mathfrak f_1(y)\bigr)-\ln\bigl(\mathfrak f(y)\bigr)$ of the two functions we get that}\label{PH91}
$$\frac{e^z}{e^z+d_n}-\Bigl[\frac{1}{n}+\frac{2}{z}+\frac{1}{z(\ln z
+C_3)}\Bigr]$$
and we check that this expression is positive. We have $e^z=y\ge \frac{2^n+4}{3}> d_n$ and thus $\frac{e^z}{e^z+d_n}<\frac{1}{2}$. If $m\ge 1024$, then $\frac{4m}{3}>1365$; so $z\ge \ln(978)>7.2$ which implies that $\frac{2}{z}+\frac{1}{z(\ln z+C_3)}<0.39$. Therefore, if $m\ge 734$, the expression is strictly greater than $\frac{1}{2}-\frac{1}{10}-0.39=0.01$. If $m\in\llbracket513,1024\rrbracket$, then we can assume that $n=10$, so $e^z=y\ge \frac{2^{11}+4}{3}>2d_n$ and hence $\frac{e^z}{e^z+d_n}\ge \frac{2}{3}$. Thus, if $m\in\llbracket513,1024\rrbracket$, we have $\frac{4m}{3}=684$, $z\ge\ln(684)>6.5$, so the expression is strictly greater than $\frac{2}{3}-\frac{1}{10}-0.44>0.12$.

So to prove that $\ln\bigl(\mathfrak f_1(y)\bigr)-\ln\bigl(\mathfrak f(y)\bigr)\le 0$ for each $y\in [\frac{2^n+4}{3},2^n+\frac{4}{3}]$ it suffices to prove that the inequality
$$13\ln\bigl(E(n,2)\bigr)+\ln(32)+\ln(n-1)\le C_1\Bigl(\frac{4}{3}m\Bigr)^{\frac{1}{n}}\Bigl[C_2\ln(\frac{4}{3}m)\Bigr]^2\Bigl[\ln\bigl(\ln(\frac{4}{3}m)\bigr)+C_3\Bigr]$$
holds for one value of $m$ which is less than or equal to $3\cdot 2^{n-2}+1$, and hence we can assume that $m=3\cdot 2^{n-2}$. 

We have inequalities
$$13n+\frac{13}{3}+13(n^2+n-2)\ln\Bigl(1-\frac{2}{n}\Bigr)+13(n^2+4)\ln\Bigl(1-\frac{1}{n}\Bigr)+\ln(32n-32)$$
$$<13n+\frac{13}{3}-26(n^2+n-2)\Bigl(\frac{1}{n}+\frac{1}{n^2}\Bigr)-\frac{13(n^2+4)}{n}-\frac{13(n^2+4)}{2n^2}+\ln(32n-32)$$
$$<-52-26n-6.5+3.5+4.34+\ln(n-1)+\frac{26}{n^2}<-50.66-26n+\frac{26}{n^2}+\ln(n-1)$$
$$<-133.85\ln(n),$$
where for the last inequality we used that $n\ge 10$. 
Thus, as $\frac{26-80.4}{26}<-4.14$, it suffices to show that 
\begin{equation}\label{EQ38.5}
26n^2\Bigl(1+\frac{1}{2n}-\frac{4.14}{n^2}\Bigr)\ln(n)\le C_1\Bigl(\frac{4}{3}m\Bigr)^{\frac{1}{n}}\Bigl[C_2\ln(\frac{4}{3}m)\Bigr]^2\Bigl[\ln\bigl(\ln(\frac{4}{3}m)\bigr)+C_3\Bigr],
\end{equation}
which follows from the following three identities $\frac{4}{3}m=2^n$, $\ln(\frac{4}{3}m)=n\ln(2)$, and $\ln(n)= -\ln(2)+\ln\bigl(\ln(\frac{4}{3}m)\bigr)$, the inequality $\frac{13}{(\ln 2)^2}<27.0578$, and the fact that for $n\ge 10$ we have $-\ln (2)+(\frac{1}{2n}-\frac{4.14}{n^2})\ln(n)<-0.673=C_3$. So part (3) holds.

{\bf Case 3: $|K|\ge 3$ and $2\lceil \frac{|K|^{n-1}}{4}\rceil< m\le |K|^n-1.$} Let 
$$C(m):=C_1(2m)^{\frac{1}{n}}\bigl[\ln(2m)]^2[C_2\ln\bigl(\ln(2m)\bigr)+C_3\bigr].$$
We use a similar approach to the one of Case 2 that reduces the problem to proving that for each $\sigma\in\Alt(K^n)$ that satisfies $\n(\sigma)\le \min(2m,|K|^n)-2\lceil \frac{|K|^{n-1}}{4}\rceil$ we have
$$\ln\bigl(\ell^{\S}_{n,K}({\sigma})\bigr)+\ln\bigl(\pi^{\S}_{n,2\lceil \frac{|K|^{n-1}}{4}\rceil}(K)\bigr)\le C(m)$$
 (cf.\ Lemmas \ref{L1}(3) and \ref{L9}(1)).\footnote{Strictly speaking, the simplified expression $\min(2m,|K|^n)-2\bigl\lceil \frac{|K|^{n-1}}{4}\bigr\rceil$ can be replaced by $\min\bigl(2m-4\bigl\lceil \frac{|K|^{n-1}}{4}\bigr\rceil,|K|^n-2\bigl\lceil \frac{|K|^{n-1}}{4}\bigr\rceil\bigr)$ and the below simplified `calibration point' of $\frac{|K|^n}{2}$ can be replaced by $\frac{|K|^n}{2}+\bigl\lceil \frac{|K|^{n-1}}{4}\bigr\rceil$. The rational number $\frac{|K|^n}{\frac{|K|^n}{2}+\bigl\lceil \frac{|K|^{n-1}}{4}\bigr\rceil}$, which is the analog of $\frac{4}{3}$ in Case 2, would complicate the notation quite a bit and this explains the simplifications adopted despite the fact that they lead to weaker results.}

\phantomsection{Let $r:=\lfloor 2\frac{\lceil\frac{|K|^{n-1}}{4}\rceil}{3}\rfloor$; we have $r\ge 4$ and one can check that, by considering the restrictions made when $n=3$, in fact we have $r\ge n+1\ge 4$. For $|K|\ge 3$ let}\label{PH99}
$$E(n,|K|):=\ln(C_{|K|-1})+\epsilon_{n,|K|-1}+(4n^2+8n-14)\ln(|K|)+(4n^2+2n-3)\ln(n-1),$$
where $C_k:=\frac{k^3+k^2+1}{2k^3}\le \frac{13}{16}$ for $k\in\mathbb N^{\ast}\setminus\{1\}$ is neglected in what follows and where $\epsilon_{n,|K|-1}$ is a negative real number estimated in Lemma \ref{L19}.

We have $\pi^{\S}_{n,2\lceil \frac{|K|^{n-1}}{4}\rceil}(K)\le D_{n,|K|-1}\le E(n,|K|)$ by Corollary \ref{C20} and 
$$\ln\bigl(\ell^{\S}_{n,K}({\sigma})\bigr)\le \wp_{r,\min(2m,|K|^n)-2\lceil\frac{|K|^{n-1}}{4}\rceil}E(n,|K|)$$ 
by Corollaries \ref{C22}(1) and \ref{C23}(1) and (3). As
 Corollary \ref{C23} gives the inequality $\ln(E_{n-2,|K|-1})+2\ln\bigl(D_{n,|K|-1})\bigr)<E(n,|K|)$, after using the inequality (see Lemma \ref{L18}(1)) $E_{n-2,|K|-1}\le C_{|K|}(n-1)(|K|-1)^4$ we get that it suffices to prove that
$$-\frac{1}{2}\ln\bigl( C_{|K|}(n-1)(|K|-1)^4\bigr)+\bigl(\wp_{r,\min(2m,|K|^n)-2\lceil\frac{|K|^{n-1}}{4}\rceil}+0.5\bigr)E(n,|K|)$$
is less than or equal to $C(m)$.
Therefore, based on Proposition \ref{PR3}(4), it suffices to prove that
\begin{equation}\label{EQ39}
E(n,|K|)\Bigl(\Bigl\lfloor\frac{\min(2m,|K|^n)-2\lceil \frac{|K|^{n-1}}{4}\rceil-2}{2r}\Bigr\rfloor+1.5\Bigr)\le C(m).
\end{equation}

It is convenient to work with all real numbers (not only integers) $m$ in the open-closed interval
$\Bigl(\max\bigl(N,2\lceil\frac{|K|^{n-1}}{4}\rceil+1\bigr),|K|^n\Bigr]$.

First, we reduce our task to proving Inequality (\ref{EQ39})
for $m=\max\bigl(N,\frac{|K|^n}{2}\bigr)$. To do this, we note that for $m\ge\max\left(N,\frac{|K|^n}{2}\right)$
the left hand side of Inequality (\ref{EQ39}) is fixed. Thus, it is enough to restrict to $m\in
\Bigl(\max\bigl(N,2\lceil\frac{|K|^{n-1}}{4}\rceil+1\bigr),\max\bigl(N,\frac{|K|^n}{2}\bigr)\Bigr]$.

We consider two disjoint subcases as follows.

{\bf Subcase 3.1: $m\in \llbracket N,675\rrbracket$.} We only have the following four possibilities for $\max(N,2\lceil \frac{|K|^{n-1}}{4}\rceil+1)\le m\le 675$: $(n,|K|)\in\{(6,3),(5,4),(3,9),(3,11)\}$.

If $(n,|K|)=(6,3)$ and $m\in \llbracket350,675\rrbracket$, then we have
$$10.495(700)^{\frac{1}{6}}[\ln(700)]^2\bigl[\ln\bigl(\ln(700)\bigr)+0.7922\bigr]>3585.9,$$ 
while on the hand we have $E(6,3)<-68.49+178\ln(3)+153\ln(5)<373.31$ (as $\epsilon_{6,2}<-68.49$ by Lemma \ref{L19}(4)) and $7.5\cdot 372.65<2795<3585.9$.

If $(n,|K|)=(5,4)$, then $E(5,4)<-41.51+126\ln(4)+107\ln(4)<281.5$ as we have $\epsilon_{5,3}<-41.51$; so 
$$11.5\cdot 281.5=3237.25<3778<9.0169(700)^{\frac{1}{5}}[\ln(700)]^2[0.7214\ln\bigl(\ln(7000)\bigr)+1.278].$$

If $(n,|K|)=(3,11)$, then $E(3,11)<46\ln(11)+39\ln(2)<137.34$ as we have $\epsilon_{3,10}<0$; so 
$$63.5\cdot 137.34<8722<9050<9.0169(700)^{\frac{1}{5}}[\ln(700)]^2[0.7214\ln\bigl(\ln(7000)\bigr)+1.278].$$

If $(n,|K|)=(3,9)$, then $E(3,9)<46\ln(9)+39\ln(2)<128.11$ as we have $\epsilon_{3,8}<0$; so 
$57.5\cdot 128.11<7367<9050.$

{\bf Subcase 3.2: $m\ge 676$.} For $x:=\ln(2m)$, we have
$x\ge\ln(1352)>7.2$. Let $s\in\{0,1,2\}$ be such that $2\lceil\frac{|K|^{n-1}}{4}\rceil=3r+s$. Thus
\begin{equation}\label{EQ21.9}
\frac{2m-2\lceil\frac{|K|^{n-1}}{4}\rceil-2}{2r}+1.5=\frac{2m-s-2}{2r}.
\end{equation}
So to prove Inequality (\ref{EQ39}) it suffices to show that for $m\le\frac{|K|^n}{2}$ we have
\begin{equation}\label{EQ40}
\frac{2m-s-2}{2r}E(n,|K|)\le C(m).
\end{equation}

We prove that the difference of the derivatives with respect to $x$ of the logarithms of the two sides of
Inequality (\ref{EQ40}) is positive, for all $x$ in the range of
interest:
\begin{equation}\label{EQ41}
\frac{e^x}{e^x-s-2}-\left[ \frac{1}{n}+\frac{2}{x}+\frac{1}{x(\ln x
+2)}\right] >0.
\end{equation}

As $\frac{e^x}{e^x-s-2}>1$, $n\ge 3$, and for all $x\geq
7.2$ we have $\frac{2}{x}+\frac{1}{x(\ln x+2)}<0.32$, it follows that Inequality (\ref{EQ41}) holds. 

So to prove Inequality (\ref{EQ40}) we can assume that $m=\max\bigl(N,\frac{|K|^n}{2}\bigr)=\frac{|K|^n}{2}$.

We have $\wp_{r,|K|^n-2\lceil\frac{|K|^{n-1}}{4}\rceil}\le 3|K|$ by Equation (\ref{EQ21.9}) and Lemma \ref{L17}(3) if $|K|$ is odd and Lemma \ref{L17}(1) if $4\mid |K|$ as $|K|^{n-2}\ge 4^3=64>24$. So to prove Inequality (\ref{EQ39}) it suffices to show that for $2m=|K|^n$ we have
\begin{equation}\label{EQ42}
E(n,|K|)(3|K|+0.5)\le C(m).
\end{equation}

Next we estimate $E(n,|K|)$ by something simpler. It seems impossible to do
this well for all $n$ and $|K|,$ so our bounds are not optimal.

First we deal with the term $(4n^2+2n-3)\ln (n-1)$ in $E(n,|K|)$ which is the fastest
growing term when $|K|\ge 3$ is fixed. As $\ln (n-1) -\ln n= \ln
(1-\frac{1}{n})=-\sum_{i=1}^{\infty}\frac{1}{in^i}$ and $n\ge 3$, we get that
$$(4n^2+2n-3) \ln (n-1)<(4n^2+2n-3)\ln n-4n-4+\frac{2}{3n}\le (4n^2+2n-3)\ln n-4n-\frac{34}{9}.$$
The function $\mathfrak f_2: [1,\infty)\rightarrow\mathbb R$ defined by the rule
$$\mathfrak f_2(y):=\frac{(2y-3)\ln y-4y-\frac{34}{9}}{4y^2}$$ 
has absolute maximum of less than $\frac{1}{50}$ at a point in the interval $[30.3,30.4]$; in fact for each $y\in\mathbb N^{\ast}$ we have $\mathfrak f_2(y)\le\mathfrak f_2(30)=0.019469...$ and for each $y\in\llbracket1,11\rrbracket$ we have $\mathfrak f_2(y)<0$. We conclude that for $n\ge 3$ we have 
$$(4n^2+2n-3)\ln (n-1)< 4n^2\bigl(\ln n +\mathfrak f_2(n)\bigr)<\begin{cases}4n^2(\ln n +\frac{1}{50})\quad\textup{if}\quad n\ge 12
\\
4n^2\ln n\qquad\qquad\textup{if}\quad n\in\llbracket3,11\rrbracket.\end{cases}$$ 

Denoting $\varepsilon_{n,|K|}:=\mathfrak f_2(n)+\frac{\epsilon_{n,|K|-1}}{4n^2}$, it follows that for $|K|\ge 3$ we have
$$E(n,|K|) \le (4n^2+8n-14)\ln |K| + 4n^2\bigl(\ln n +\varepsilon_{n,|K|}\bigr).$$
Therefore we need to estimate
\begin{equation}\label{EQ43-}
(3|K|+0.5)\left[(4n^2+8n-14)\ln (|K|) + 4n^2\bigl(\ln n +\varepsilon_{n,|K|}\bigr)\right].
\end{equation}

We compute $4n^2=\frac{4[\ln (2m)]^2}{(\ln |K|)^2}$, $\ln(n)=\ln\bigl(\ln(2m)\bigr)-\ln\bigl(\ln(|K|)\bigr)$, 
$$(4n^2+8n-14)\ln|K|=4n^2\Bigl(1+\frac{2}{n}-\frac{7}{2n^2}\Bigr)\ln(|K|)=\frac{4[\ln(2m)]^2}{\ln
|K|}\Bigl(1+\frac{2}{n}-\frac{7}{2n^2}\Bigr),$$
and $3|K|+0.5=3(|K|+\frac{1}{6})=3\bigl[(2m)^{\frac{1}{n}}+\frac{1}{6}\bigr]=3(2m)^{\frac{1}{n}}\bigl(1+\frac{1}{6|K|}\bigr)$.

Putting all these together, for $|K|\ge 3$ we get the following upper bound for the expression of Display (\ref{EQ43-}) and hence also of the left hand side of Inequality (\ref{EQ40}):
\begin{equation}\label{EQ43}
\frac{12}{\ln |K|}\Bigl(1+\frac{1}{6|K|}\Bigr)(2m)^{\frac{1}{n}}[\ln(2m)]^2\Bigl[\frac{\ln\bigl(\ln(2m)\bigr)}{\ln |K|}+1+\frac{2}{n}-\frac{7}{2n^2}+\frac{\varepsilon_{n,|K|}-\ln\bigl(\ln(|K|)\bigr)}{\ln|K|}\Bigr].
\end{equation}

For $|K|\ge 4$, we have $\ln(4)^{-1}<0.7214$, $\varepsilon_{n,|K|}-\ln\bigl(\ln(|K|)\bigr)<\frac{1}{50}-\ln\bigl(\ln(4)\bigr)<0$, $\frac{12}{\ln |K|}\left(1+\frac{1}{6|K|}\right) \le \frac{25}{2\ln 4}<9.0169$, and, as $n\ge 3$, $1+\frac{2}{n}-\frac{7}{2n^2}\le\frac{23}{18}<1.278$. So the Display (\ref{EQ43}) is strictly less than
$9.0169(2m)^{\frac{1}{n}}[\ln(2m)]^2\left[0.7214\ln\bigl(\ln(2m)\bigr)+1.278\right]$ and thus Inequality (\ref{EQ42}) holds for $|K|\ge 4$. So part (1) holds.

For $|K|=3$, we have $n\ge 7$ and the expression of Display (\ref{EQ43}) is equal to 
\begin{equation}\label{EQ44}
\frac{12}{(\ln 3)^2}\Bigl(1+\frac{1}{18}\Bigr)(2m)^{\frac{1}{n}}[\ln(2m)]^2\Bigl[\ln\bigl(\ln(2m)\bigr)+\Bigl(1+\frac{2}{n}-\frac{7}{2n^2}\Bigr)\ln(3)+\varepsilon_{n,3}-\ln\bigl(\ln(3)\bigr)\Bigr].
\end{equation}
We have $\frac{\epsilon_{n,2}}{4n^2}<-\frac{7}{18}<-0.3888$ for $n\ge 12$ and $\frac{\epsilon_{n,2}}{4n^2}<-0.4528$ for $n\in\llbracket7,11\rrbracket$ based on Lemma \ref{L19}(4). Thus for $n\ge 12$, $\bigl(1+\frac{2}{n}-\frac{7}{2n^2}\bigr)\ln(3)+\varepsilon_{n,3}-\ln\bigl(\ln(3)\bigr)$ is strictly less than $\bigl(1+\frac{1}{6}-\frac{7}{288}\bigr)\ln(3)-0.3888+0.02-\ln\bigl(\ln(3)\bigr)<0.7922$, and for $n\in\llbracket7,11\rrbracket$, as $\mathfrak f_2(n)<0$, $\bigl(1+\frac{2}{n}-\frac{7}{2n^2}\bigr)\ln(3)+\varepsilon_{n,3}-\ln\bigl(\ln(3)\bigr)$ is strictly less than $\frac{17}{14}\ln(3)-0.4528-\ln\bigl(\ln(3)\bigr)<0.7872$. From all these and the inequality $\frac{12}{(\ln 3)^2}(1+\frac{1}{18})=\frac{38}{3(\ln 3)^2}<10.495$ we get that Inequality (\ref{EQ42}) holds for $|K|=3$. So part (2) holds.\end{proof}
 
\begin{corollary}\label{C24} Suppose that $n\in\mathbb N^{\ast}\setminus\{1,2\}$, $K$ is a finite field, $N\in\{513,350,626\}$ with $N=626$ iff $|K|\ge 3$ and $N=350$ iff $|K|=3$, and $m\in\llbracket N,|K|^n\rrbracket$. If $4\mid |K|$, then we also assume that $m\le |K|^n-2$. Then the following properties hold.

\medskip
{\bf (1)} If $|K|\geq 4$, then 
$$\ln\bigl(\pi_{n,m}(K)\bigr)< 11.3606|K|[\ln(2m)]^2[0.7214\ln\bigl(\ln(2m)\bigr)+1.278]$$
and, if $m\le |K|^n-2$, 
$$\ln\bigl(\pi^{\S}_{n,m}(K)\bigr)< 11.3606|K|[\ln(2m)]^2[0.7214\ln\bigl(\ln(2m)\bigr)+1.278].$$

{\bf (2)} If $|K|=3$, then $\ln\bigl(\pi_{n,m}(K)\bigr)< 35.3408[\ln(2m)]^2[\ln\bigl(\ln(2m)\bigr)+0.7922]$ and, if $m\le |K|^n-2$, $\ln\bigl(\pi^{\S}_{n,m}(K)\bigr)< 35.3408[\ln(2m)]^2[\ln\bigl(\ln(2m)\bigr)+0.7922]$.

\smallskip
{\bf (3)} If $|K|=2$, then $\ln\bigl(\pi_{n,m}(K)\bigr)< 55.696[\ln(\frac{4m}{3})]^2[\ln\bigl(\ln(\frac{4m}{3})\bigr)-0.673]$.
\end{corollary}

\begin{proof} 
As $m\le |K|^n,$ we have $m^{\frac{1}{n}}\le |K|.$ For part (1), we note that $2^{\frac{1}{n}}\le 2^{\frac{1}{3}}$ as $n\ge 3$ and $9.0169\cdot 2^{1/3}<11.3606$. For part (2) we note that $2^{\frac{1}{n}}\le 2^{\frac{1}{6}}$ as $n\ge 6$ and $10.495\cdot 2^{1/6}\cdot 3<35.3408$. For part (3) we note that $\bigl(\frac{4}{3}\bigr)^{\frac{1}{n}}\le \bigl(\frac{4}{3}\bigr)^{\frac{1}{10}}$ as $n\ge 10$ and $27.0578\cdot \bigl(\frac{4}{3}\bigr)^{1/10}\cdot 2 <55.696$. Based on all these, parts (1) to (3) follow from Theorem \ref{T14}(1) to (3) (respectively).
\end{proof}

\begin{remark}\normalfont\label{R9}
{\bf (1)} Using the same ideas for each fixed finite field $K$, we can get better upper bounds for $\ln\bigl(\pi_{n,m}(K)\bigr)$, especially for larger $|K|$ and $m$: if $\varepsilon>0,$ then $\ln\bigl(\pi_{n,m}(K)\bigr)$ with $n\ge \frac{\ln m}{\ln |K|}$ is bounded from above by $\bigl[\frac{12|K|+2}{(\ln |K|)^2}+\varepsilon\bigr](\ln m)^2 \ln (\ln m)$ for $m \gg 1$ by Display (\ref{EQ43}) for $|K|\ge 4$, by Display (\ref{EQ44}) for $|K|=3$, and by Inequalities (\ref{EQ38.5}) for $|K|=2$.

\smallskip
{\bf (2)} Similarly, for $n\in\mathbb N^{\ast}\setminus\{1,2\}$ fixed we can get an upper bound for $\ln\bigl(\pi_{n,m}(K)\bigr)$ that is asymptotically $C(n)m^{\frac{1}{n}}\ln(m)$, with $C(n)\in (0,\infty)$, an improvement of Theorem \ref{T14} by a multiple of $\ln(m)\ln\bigl(\ln(m)\bigr)$. 

\smallskip
{\bf (3)} Often one gets better results for $s\in \llbracket5,2p-1\rrbracket$ odd than for $s=3$ in Corollaries \ref{C17}(1) and \ref{C22} and in some cases the same is true for the upper bounds for the $N_{\sigma}$s that involve the $\n(\sigma)$s. For instance, if $K=\mathbb F_7$ and we consider the version of Corollary \ref{C17}(1) for $r=1$, then $r=1$ and hence $s\in\{3,5\}$ and $t\ge 1$; the upper bounds for $N_{\sigma}$ we would get from Corollary \ref{C17}(1) for $\varepsilon=0$ and $N=4$, say working with $\n(\sigma)\ge 15$, is $1728^{2\lfloor\frac{\n(\sigma)}{2}\rfloor}$ for $s=3$ based on Lemma \ref{P3}(1) and is $1440^{4\lfloor\frac{\n(\sigma)}{5}\rfloor+8}$ for $s=5$ based on Theorem \ref{T2}(7), the latter being smaller if $\lfloor\frac{\n(\sigma)}{2}\rfloor\ge 2\lfloor\frac{\n(\sigma)}{5}\rfloor+4$ and hence if $\n(\sigma)\ge 36$.

\smallskip
{\bf (4)} For $s=5$ and $p\ge 3$ (resp.\ $s\in \llbracket7,2p-1\rrbracket$ odd with $p\ge 5$), if any stronger version of Theorem \ref{T2}(7) or (9) (resp.\ \ref{T2}(8) or (9)) of the form 
$$\nu_{5,r}(\sigma)\le C_1\Bigl\lfloor\frac{\n(\sigma)}{5r}\Bigr\rfloor+C_2$$ 
(resp.\ $\nu_{s,r}(\sigma)\le C_1\lfloor\frac{\n(\sigma)}{sr}\rfloor+C_2$) with $(C_1,C_2)\in (0,\frac{3}{2})\times (0,\infty)$ holds for a suitable class of $\sigma\in A_l$, then one implicitly gets stronger versions of Corollary \ref{C17}(2) for such $\sigma$s. Unfortunately, the best constant $C_1$ one could hope to get for all $\sigma$s in the $A_l$s with $l\ge 5$ (resp.\ $l\ge 7$) when $l$ goes to infinity is $\frac{5}{3}$ (resp.\ $\frac{3}{2-\frac{\epsilon(s)}{s}}$ where $\epsilon(s)\in\{0,1,2\}$ is such that $3\mid 2s-\epsilon(s)$) by \cite{HR2}, Thm.\ 6; cf.\ Theorem \ref{P4}(2) and also Theorem \ref{T2}(9) which shows that the constant $\frac{15}{8}$ (resp.\ $\frac{3}{2-\frac{2}{s}}$) is actually achievable modulo either a correction term or the addition of an arbitrary positive constant. This explains why Theorems \ref{T12}, \ref{T13}, and \ref{T14} and Corollaries \ref{C22}(3) and \ref{C24} are worked out via $3$-cycles despite the fact that $E_{n-2,x,y}$ is a decreasing function on $y\ge 0$.

\smallskip
{\bf (5)} Variants of Corollaries \ref{C17} and \ref{C22} (and hence also of Theorems \ref{T12}, \ref{T13}, and \ref{T14})) can be obtained if in their proofs one only uses the union of conjugacy classes $\cup_{i=1}^{\frac{r}{j}} \mathcal C_{ij}$ with $j$ a fixed divisor of $r$.

\smallskip
{\bf (6)} Case 2 of the proof of Theorem \ref{T14} can be also worked out using $3$-cycles in a way similar to Case 3 but the resulting bounds are slightly weaker as $C_3$ is larger and the `calibration point' is larger than $\frac{4}{3}$ and depends on $n$.
\end{remark}

\section{Upper bounds via stabilizers}\label{S25}

In this section we study the rational numbers $\frac{\pi_{n,m+s}(K)}{\pi_{n,m}(K)}$ (resp.\ $\frac{\pi^{\S}_{n,m+s}(K)}{\pi^{\S}_{n,m}(K)}$) with $(n,m,s)\in (\mathbb N^{\ast})^3$ such that $\sqrt[n]{m+s}\le |K|$ (resp.\ $\sqrt[n]{m+s+2}\le |K|$). We begin with the following definition that continues Definition \ref{D2.5}.
 
\begin{definition}\label{D21}
Let $n\in\mathbb N^{\ast}\setminus\{1\}$ and $K$ a field. Let $s\in\mathbb N^{\ast}$ and a finite non-empty subset $Y$ of $K^n$ be such that $s+|Y|\le |K|^n$. 

\medskip
{\bf (1)} Let $m\in\mathbb N^{\ast}$ be such that $\sqrt[n]{m+s}\le |K|$. If $K$ is a finite field with $4\mid |K|$, then we assume that $m\le |K|^n-s-2$. Let 
$$\Xi_{n,m,s}(K):=\max\bigl(\pi_{Y,s}(K)|Y\subset K^n, |Y|=m\bigr)\in \llbracket1,\pi_{n,m+s}(K)\rrbracket.\index{$\Xi_{n,m,s}(K)$ invariant}$$

\smallskip
{\bf (2)} Let $m\in\mathbb N^{\ast}$ be such that $\sqrt[n]{m+s+2}\le |K|$. Let 
$$\Xi^{\S}_{n,m,s}(K):=\max\bigl(\pi^{\S}_{Y,s}(K)|Y\subset K^n, |Y|=m\bigr)\in \llbracket1,\pi^{\S}_{n,m+s}(K)\rrbracket.\index{$\Xi^{\S}_{n,m,s}(K)$ invariant}$$
\end{definition}

If in Definitions \ref{D2.5} and \ref{D21} we have $|Y|=1$, then $\pi_{Y,s}=\pi_{n,s+1}(K)=\Xi_{n,1,s}(K)$ and $\pi^{\S}_{Y,s}=\pi^{\S}_{n,s+1}(K)=\Xi^{\S}_{n,1,s}(K)$.

Let $\kappa_{n,K}\in\mathbb N^{\ast}$ be as in Proposition \ref{PR16}. Let $\mathcal Y_2$ and $\mathcal Y_3$ be as in Notation \ref{N4}. 

\begin{proposition}\label{PR26}
Let $K$ be a finite field. Let $(n,m,s)\in (\mathbb N^{\ast}\setminus\{1\})^2\times\mathbb N^{\ast}$ be such that $m+s\le |K|^n$. If $4\mid |K|$, then we assume that $m+s\le |K|^n-2$. Let $Y$ be a subset of $K^n$ with $m$ elements. Then the following properties hold.

\medskip
{\bf (1)} If $4\nmid |K|$, then we have inequalities
$$\pi_{Y,s}(K)\le\Xi_{n,m,s}(K)\le\Omega(\mathcal Y_2^{\min(s,|K|^n-m-1)})\le\kappa_{n,K}^{\min(s,|K|^n-m-1)}.$$

{\bf (2)} If $m+s\le |K|^n-2$, then we have inequalities
$$\pi_{Y,s}(K)\le\Xi_{n,m,s}(K)\le\Omega(\mathcal Y_3)^{\min(s,\lfloor\frac{|K|^n-m}{2}\rfloor)}$$
and
$$\pi^{\S}_{Y,s}(K)\le\Xi^{\S}_{n,m,s}(K)\le\Omega^{\S}(\mathcal Y_3)^{\min(s,\lfloor\frac{|K|^n-m}{2}\rfloor)}.$$

{\bf (3)} We have inequalities 
$$\max\bigl(\pi_{n,m}(K),\Xi_{n,m,s}(K)\bigr)\le\pi_{n,m+s}(K)\le \pi_{n,m}(K)\Xi_{n,m,s}(K).$$
If moreover $m+s\le |K|^n-2$, then we also have inequalities
$$\max\bigl(\pi^{\S}_{n,m}(K),\Xi^{\S}_{n,m,s}(K)\bigr)\le\pi^{\S}_{n,m+s}(K)\le \pi^{\S}_{n,m}(K)\Xi^{\S}_{n,m,s}(K).$$

\smallskip
{\bf (4)} If $s=1$ and $\d_Y\le n-1$, then $\pi^{\S}_{Y,1}(K)\le\min\bigl(m+1-\d_Y,\d_Y(|K|-1)\bigr)^2$.
\end{proposition}

\begin{proof}
Let $(\underline{P},\underline{Q})=\bigl((P_1,\ldots,P_{m+s}),(Q_1,\ldots,Q_{m+s})\bigr)\in\mathbb D_{n,m+s}(K)^2$. Let the permutation $\sigma\in\perm(K^n)$ be such that $\sigma(\underline{P})=\underline{Q}$.

We can assume that $Y=\{P_i|i\in \llbracket1,m\rrbracket\}$. We emphasize that for parts (1) and (2) we also assume that $Q_i=P_i$ for each $i\in \llbracket1,m\rrbracket$. 

There exists $l\in \llbracket0,\min(s,|K|^n-m-1)\rrbracket$ and $l$ transpositions $\tau_1,\ldots,\tau_l$ in $\perm(K^n)$ such that for $\theta:=\prod_{i=1}^l \tau_i$, we have $\theta(\underline{P})=\underline{Q}$ by Lemma \ref{L1}(1) applied to $l=|K|^n-m$. If $a\in\TGA_n(K)[\Omega(\mathcal Y_2^l)]$ is such that $a(K)=\theta$, then $a(\underline{P})=\underline{Q}$. This implies that the second inequality of part (1) holds. As the first and the third inequalities of part (1) are clear, part (1) holds.

For part (2), based on Lemma \ref{L1}(2) applied to $l=|K|^n-m$ we can assume that $\sigma$ is even and that $\nu_{3,1}(\sigma)\le\min\bigl(s,\lfloor\frac{|K|^n-m}{2}\rfloor\bigr)$. We consider $3$-cycles $\vartheta_1,\ldots,\vartheta_{\nu_{3,1}(\sigma)}$ in $\Alt(K^n)$ such that $\sigma=\prod_{i=1}^{\nu_{3,1}(\sigma)} \vartheta_i$. For each $i\in \llbracket1,l\rrbracket$ let $a_i\in \TGA_n(K)[\Omega(\mathcal Y_3)]$ and $b_i\in \STGA_n(K)[\Omega^{\S}(\mathcal Y_3)]$ be such that we have $a_i(K)=b_i(K)=\vartheta_i$. So, by defining automorphisms $a:=\prod_{i=1}^{\nu_{3,1}(\sigma)} a_i\in\TGA_n(K)[\Omega(\mathcal Y_3)^{\min(s,\lfloor\frac{|K|^n-m}{2}\rfloor)}]$ and $b:=\prod_{i=1}^{\nu_{3,1}(\sigma)} b_i\in\STGA_n(K)[\Omega^{\S}(\mathcal Y_3)^{\min(s,\lfloor\frac{|K|^n-m}{2}\rfloor}]$, we have $a(K)=b(K)=\sigma$ and thus $a(\underline{P})=b(\underline{P})=\underline{Q}$. This implies that part (2) holds.

To prove part (3), let $c_0\in\TGA_n(K)[\pi_{n,m}(K)]$ be such that $c_0(Q_i)=P_i$ for each $i\in \llbracket1,m\rrbracket$. Let $\sigma_1:=c_0(K)\sigma\in\perm(K^n)$. As $\sigma_1(P_i)=P_i$ for each $i\in \llbracket1,m\rrbracket$, we have $\sigma_1\in\Fix_{\TGA_n(K)}(Y)$. Let $c_1\in\Fix_{\TGA_n(K)}(Y)[\pi_{Y,s}(K)]$ be such that $c_1(P_i)=\sigma_1(P_i)$ for each $i\in \llbracket m+1,m+s\rrbracket$. Hence for $c:=c_0^{-1}c_1$ we have $c(\underline{P})=\sigma(\underline{P})=\underline{Q}$ and 
$$\ell(c)\le\ell(c_0)\ell(c_1)\le \pi_{n,m}(K)\pi_{Y,s}(K)\le \pi_{n,m}(K)\Xi_{n,m,s}(K).$$ 
This implies that $\pi_{n,m+s}(K)\le \pi_{n,m}(K)\Xi_{n,m,s}(K)$. 

Clearly, $\pi_{n,m}(K)\le \pi_{n,m+s}(K)$ and $\pi^{\S}_{n,m}(K)\le \pi^{\S}_{n,m+s}(K)$. To prove the inequality $\Xi_{n,m,s}(K)\le\pi_{n,m+s}(K)$, we can assume that $\Xi_{n,m,s}(K)=\pi_{Y,s}$ and that the pair $\bigl((P_{m+1},\ldots,P_{m+s}),(Q_{m+1},\ldots,Q_{m+s})\bigr)\in\mathbb D_s(\mathbb A^n_K\setminus Y)^2$ is such that for each automorphism $b_Y\in\Fix_{\TGA_n(K)}(Y)$ with $b_Y(P_i)=Q_i$ for every $i\in \llbracket m+1,m+s\rrbracket$ we have $\ell(b_Y)\ge\pi_{Y,s}=\Xi_{n,m,s}(K)$. Taking $Q_i=P_i$ for each $i\in \llbracket1,m\rrbracket$, for each $b_0\in\TGA_n(K)$ with $b_0(\underline{P})=\underline{Q}$ we have $b_0\in\Fix_{\TGA_n(K)}(Y)$ and $b_0(P_i)=Q_i$ for each $i\in \llbracket m+1,m+s\rrbracket$, so $\ell(b_0)\ge\Xi_{n,m,s}(K)$ and thus $\pi_{n,m}(K)\ge\Xi_{n,m,s}(K)$. The inequality $\Xi^{\S}_{n,m,s}(K)\le\pi^{\S}_{n,m+s}(K)$ is proved similarly. So part (3) holds.

For part (4), up to linear automorphisms we can assume that $\langle Y\rangle_{\aff}$ is contained in the zero locus $x_n=0$. For each $O\in\{P_{m+1},Q_{m+1}\}$, we check that there exists $f_O\in K[x_1,\ldots,x_{n-1}]$ such that $\deg(f_O)\le\min\bigl(m+1-\d_Y,\d_Y(|K|-1)\bigr)$, $f_O(Y)=\{0\}$, and for $$c_O:=\e\bigl(x_1,\ldots,x_{n-1},x_n+f_O(x_1,\ldots,x_{n-1})\bigr)\in\Fix_{\STGA_n(K)}(Y)$$ 
we have $c_O(O)\notin \langle Y\rangle_{\aff}$. For this we can assume that $O\in\langle Y\rangle_{\aff}$, and in this case the existence of $f_O$ follows from Theorem \ref{T5}(1) applied to the set $Y\cup\{O\}$ with $m+1$ elements. 

We can make the choices in such that a way that there exists $c_2\in\ASL_n(K)$ which acts identically on $\langle Y\rangle_{\aff}$ and moreover we have $c_2\bigl(c_{P_{m+1}}(P_{m+1})\bigr)=c_{Q_{m+1}}(Q_{m+1})$. Then for 
$$c_3:=c_{Q_{m+1}}^{-1}c_2c_{P_{m+1}}\in \Fix_{\STGA_n(K)}(Y)\bigl[\min\bigl(m+1-\d_Y,\d_Y(|K|-1)\bigr)^2\bigr]$$ 
we have $c_3(P_{m+1})=Q_{m+1}$. So part (4) holds.
\end{proof}

The following consequence supplements Corollary \ref{C21}.\footnote{For larger values of $m$ (resp.\ of $m$ and $s$) in their ranges of integers, the bounds in Corollary \ref{C21} are better than the ones of Corollary \ref{C25}(1) to (3) (resp.\ Corollary \ref{C25}(4)) and the opposite holds if ``larger'' is replaced by smaller.}

\begin{corollary}\label{C25}
For $n\in\mathbb N^{\ast}\setminus\{1\}$ let $\mathfrak e_n:=\frac{\varepsilon_{n,1,2}}{2}=\min\bigl(n+4\ln(n-2)+\frac{1}{3},2n-\frac{4}{3}\bigr)$. Let $l:=\lfloor\frac{n}{2}\rfloor$, $r\in\{\lceil\frac{2^n}{\sqrt{2}}\rceil+1,\lceil\frac{2^n}{\sqrt{2}}\rceil+2\}$ the largest such that $2^n>4\Bigl\lfloor\frac{\frac{r(r-1)}{2}}{2^n+1}\Bigr\rfloor-4$, and $r_0=r_0(2^l)\in\{\lceil2^{\frac{l-1}{2}}\rceil,\lceil2^{\frac{l-1}{2}}\rceil+1\}$ the smallest such that $2^{l-1}+2\le r_0^2-r_0$. Then the following properties hold.

\medskip
{\bf (1)} Suppose that $n=2l\ge 4$ and $m\in \llbracket r,2^n-1\rrbracket$. Then 
$$\pi_{n,m}(\mathbb F_2)\le (l-1)^4(2n-3)^{\lfloor\frac{m-r}{2}\rfloor}(n-1)^{m-r-2\lfloor\frac{m-r}{2}\rfloor}.$$

{\bf (2)} Suppose that $n=2l\ge 4$ and $m\in \llbracket 2^l,2^n-1\rrbracket$. Then 
$$\pi_{n,m}(\mathbb F_2)\le (2^l-r_0)^2(2^l-2)^2(2^l-1)^2(2n-3)^{\lfloor\frac{m-2^l}{2}\rfloor}(n-1)^{m-2^l-2\lfloor\frac{m-2^l}{2}\rfloor}.$$

{\bf (3)} Suppose that $n=2l+1\ge 7$. Let $m_1\in\llbracket 2^l,2r\rrbracket$ be the largest such that $m_1^2-m_1<2^{2l}+4$. If $m\in \llbracket2^{l-1},2^n-1\rrbracket$, then we have
$$\pi_{n,m}(\mathbb F_2)\le 4l^6(2n-3)^{\lfloor\frac{\max(m-m_1,0)}{2}\rfloor}(n-1)^{\max(m-m_1,0)-2\lfloor\frac{\max(m-m_1,0)}{2}\rfloor}.$$

{\bf (4)} If $n\ge 4$ and $m\in \llbracket2^{n-2},2^n-1\rrbracket$, then 
$$\pi_{n,m}(\mathbb F_2)\le (2n-3)^{\lfloor\frac{m-2^{n-2}}{2}\rfloor}(n-1)^{m-2^{n-2}-2\lfloor\frac{m-2^{n-2}}{2}\rfloor+n^2+4}(n-2)^{n^2+n-6}e^{\mathfrak e_n}.$$

{\bf (5)} If $n\ge 4$, $m\in \llbracket2^{n-2},2^n-3\rrbracket$, and $s\in \llbracket1,2^n-2-m\rrbracket$, then 
$$\pi_{n,m+s}(\mathbb F_2)\le(2n-3)^{\lfloor\frac{m-2^{n-2}}{2}\rfloor}(n-1)^{m-2^{n-2}-2\lfloor\frac{m-2^{n-2}}{2}\rfloor+2n^2+n-2} e^{\mathfrak e_n}N_{n,s,|K|},$$
where $N_{n,s,|K|}:=\min\bigl((2n-3)^{\lfloor\frac{s}{2}\rfloor}(n-1)^{s-2\lfloor\frac{s}{2}\rfloor},(2n-3)^{\min(s,\lfloor\frac{2^n-m}{2}\rfloor)}\bigr)$.
\end{corollary}

\begin{proof}
For part (1) (resp.\ (2) or (3) or (4)), let $t:=m-r\in \llbracket0,2^n-1-r\rrbracket$ (resp.\ $t:=m-2^l\in \llbracket0,2^n-1-2^l\rrbracket$ or $t:=\max(m-m_1,0)\in\llbracket0,2^n-1-m_1\rrbracket$ or $t:=m-2^{n-2}\in \llbracket0,2^n-1-2^{n-2}\rrbracket$). 

For $t=0$, part (1) (resp.\ (2) or (3) or (4)) is a particular case of Proposition \ref{PR21}(1.b) (resp.\ Proposition \ref{PR21}(2) or Proposition \ref{PR23}(1) or Corollary \ref{C20} combined with Lemma \ref{L16}(2)).

If $t>0$, then $\pi_{n,m}(\mathbb F_2)\le\pi_{n,m-t}(\mathbb F_2)\Omega(\mathcal Y_2^t)$ by Proposition \ref{PR26}(1) and (3); so part (1) (resp.\ (2) or (3) or (4)) follows from Proposition \ref{PR16}(4) and (5) and the prior paragraph.

Part (5) follows from part (4), the inequality $\Omega(\mathcal Y_3)\le 2n-3$ (see Proposition \ref{PR16}(5)), and Proposition \ref{PR26}(2) and (3).\end{proof}

\begin{definition}\label{D22}
Let $K$ be a field and $(m,l)\in (\mathbb N^{\ast})^2$ with $1\le m\le |K|^n$. Let $\natural_{n,m,l}(K)\in \llbracket1,m\rrbracket$\index{$\natural_{n,m,l}(K)$ invariant} (resp.\ $\natural^{\S}_{n,m,l}(K)\in \llbracket1,m\rrbracket$\index{$\natural^{\S}_{n,m,l}(K)$ invariant}) be the largest such that for every pair $\bigl((P_1,\ldots,P_m),(Q_1,\ldots,Q_m)\bigr)\in\mathbb D_{n,m}(K)^2$ there exists an automorphism $a\in\TGA_n(K)[l]$ (resp.\ $a\in\STGA_n(K)[l]$) such that we have the inequality $\break|\{i\in \llbracket1,m\rrbracket|a(P_i)=Q_i\}|\ge\natural_{n,m,l}(K)$ (resp.\ $|\{i\in \llbracket1,m\rrbracket|a(P_i)=Q_i\}|\ge\natural^{\S}_{n,m,l}(K)$). 
\end{definition}

Clearly, $\natural_{n,m,l}(K)\le \natural^{\S}_{n,m,l}(K)$. The sequence $\bigl(\natural_{n,m,l}(K)\bigr)_{l\ge 1}$ is non-decreasing, if $\pi_{n,m}(K)\in\mathbb N^{\ast}$ then the sequence $\bigl(\natural_{n,m,l}(K)\bigr)_{l\ge\pi_{n,m}(K)}$ is constant of constant value $m$, and if $\pi_{n,m}(K)\in\mathbb N^{\ast}\setminus\{1\}$, then $\natural_{n,m,\pi_{n,m}(K)-1}(K)<m$. Similarly, the sequence $\bigl(\natural^{\S}_{n,m,l}(K)\bigr)_{l\ge 1}$ is non-decreasing, if $\pi^{\S}_{n,m}(K)\in\mathbb N^{\ast}$ then the sequence $\bigl(\natural^{\S}_{n,m,l}(K)\bigr)_{l\ge\pi^{\S}_{n,m}(K)}$ is constant of constant value $m$, and if $\pi_{n,m}^{\S}(K)\in\mathbb N^{\ast}\setminus\{1\}$, then $\natural^{\S}_{n,m,\pi^{\S}_{n,m}(K)-1}(K)<m$.

The following applications of the $\natural_{n,m,l}(K)$s are modeled on Proposition \ref{PR25}.

\begin{proposition}\label{PR27}
Let $n\in\mathbb N^{\ast}\setminus\{1\}$ and $s\in \llbracket1,n\rrbracket$. Let $K$ be a finite field and $m\in \llbracket|K|^{s-1}+1,|K|^s\rrbracket$. Let $\epsilon\in\{0,1\}$ be such that we have $\epsilon=0$ iff $|K|=2$ and $m\in \llbracket2^s-s+1,2^s\rrbracket$. Then the following properties hold.

\medskip
{\bf (1)} We have $\natural_{n,m,1}(K)\ge s+\epsilon$.

\smallskip
{\bf (2)} Let $(\underline{P},\underline{Q})\in\mathbb D_{n,m}(K)^2$. Then for every $l\in\mathbb N^{\ast}$ there exists a pair 
$$(a_l,\sigma_l)\in\TGA_n(K)[l]\times\perm(K^n)$$ 
such that $\sigma_l\bigl(a_l(P_i)\bigr)=Q_i$ for each $i\in \llbracket1,m\rrbracket$ and $\n(\sigma_l)\in \llbracket0,2m-2\natural_{n,m,l}(K)\rrbracket$.

\smallskip
{\bf (3)} If $4\mid |K|$ we assume that $m\le |K|^n-2$. Then for each $l\in \llbracket1,\pi_{n,m}(K)-1\rrbracket$ and $t\in \llbracket1,\natural_{n,m,l}\rrbracket$, we have an inequality $\pi_{n,m}(K)\le l\Xi_{n,t,m-t}(K)$.

\smallskip
{\bf (4)} If $4\mid |K|$ we assume that $m\le |K|^n-2$. Then we have an inequality $\pi_{n,m}(K)\le\Xi_{n,s+\epsilon,m-s-\epsilon}(K)$.

\smallskip
{\bf (5)} If $s<n$, then $\natural^{\S}_{n,m,1}(K)\ge s+\epsilon$.

\smallskip
{\bf (6)} Let $(\underline{P},\underline{Q})\in\mathbb D_{n,m}(K)^2$. Then for every $l\in\mathbb N^{\ast}$ there exists a pair 
$$(a_l,\sigma_l)\in\STGA_n(K)[l]\times\perm(K^n)$$ 
such that $\sigma_l\bigl(a_l(P_i)\bigr)=Q_i$ for each $i\in \llbracket1,m\rrbracket$ and $\n(\sigma_l)\in \llbracket0,2m-2\natural^{\S}_{n,m,l}(K)\rrbracket$. 

\smallskip
{\bf (7)} Suppose that $m\le |K|^n-2$. Then for each $l\in \llbracket1,\pi^{\S}_{n,m}(K)-1\rrbracket$ and every $t\in \llbracket1,\natural^{\S}_{n,m,l}\rrbracket$, we have an inequality $\pi^{\S}_{n,m}(K)\le l\Xi^{\S}_{n,t,m-t}(K)$.

\smallskip
{\bf (8)} Suppose that $m\le |K|^n-2$ and $s\le n-1$. Then we have an inequality $\pi^{\S}_{n,m}(K)\le\Xi^{\S}_{n,s+\epsilon,m-s-\epsilon}(K)$.
\end{proposition}

\begin{proof}
Part (1) is only a translation of the proof of Proposition \ref{PR25}(1) and (2).

For part (2) we first take $a_l$ such that the set $I:=\{i\in \llbracket1,m\rrbracket|a_l(P_i)=Q_i\}$ has cardinality $\natural_{n,m,l}(K)$. Next we take $\sigma_l\in\perm(K^n)$ such that we have an identity $\supp(\sigma_l)=\{a_l(P_i),Q_i|i\in \llbracket1,m\rrbracket\setminus I\}$ and $\sigma_l(P_i)=Q_i$ for each $\llbracket1,m\rrbracket\setminus I$. Clearly, $\sigma_l\bigl(a_l(P_i)\bigr)=Q_i$ for each $i\in \llbracket1,m\rrbracket$ and $\n(\sigma_l)\in \llbracket0,2m-2\natural_{n,m,l}(K)\rrbracket$, so part (2) holds.

\phantomsection{Let $J$ be a subset of $I$ with $t$ elements. Let $Y:=\{Q_i|i\in J\}$; so $|Y|=|J|=t$. If $b_t\in\Fix_{\TGA_n(K)}(Y)[\pi_{Y,m-t}]$ is such that $b_t\bigl(a_l(P_i)\bigr)=Q_i$ for each $i\in \llbracket1,m\rrbracket\setminus J$ (equivalently, for each $i\in \llbracket1,m\rrbracket$), then for $c:=b_ta_l\in\TGA_n(K)[l\pi_{Y,m-t}]$ we have $c(\underline{P})=\underline{Q}$. As $\pi_{Y,m-t}\le\Xi_{n,t,m-t}(K)$, it follows that part (3) holds.}\label{PH96}

Part (4) follows from parts (1) and (3) applied to $l=1$ and $t=s+\epsilon$.

Parts (5) to (8) are proved in the same way as parts (1) to (4) (respectively).
\end{proof}

\section{Lower bounds for small $m$}\label{S26}

For lower bounds for the $\pi^{\S}_{n,m}(K)$s when $|K|\ge m$ we need the following notions and basic results on them.

\begin{definition}\label{D23}
Let $m\in\mathbb N^{\ast}$. Let $J$ be a subset of $\llbracket0,m-1\rrbracket$. Let $K$ be a field with $|K|\ge m$.

\medskip
{\bf (1)} We say that $\bigl((\alpha_1,\ldots,\alpha_m),(\beta_1,\ldots,\beta_m)\bigr)\in\mathbb D_{1,m}(K)^2$ is $J$-generic\index{$J$-generic} if for each pair $(f,g)\in K[x]^2$ with $1\le\deg(f)+\deg(g)\le m-1$ and $\deg(g)\in J$ there exists $i\in \llbracket1,m\rrbracket$ such that $f(\alpha_i)+g(\beta_i)\neq 0$.

\smallskip
{\bf (2)} Suppose $J\subset \llbracket0,m-2\rrbracket$. We say that $\bigl((\alpha_1,\ldots,\alpha_m),(\beta_1,\ldots,\beta_m)\bigr)\in\mathbb D_{1,m}(K)^2$ is weakly $J$-generic\index{$J$-generic!weakly $J$-generic} if for each pair $(f,g)\in K[x]^2$ with $1\le\deg(f)+\deg(g)\le m-2$ and $\deg(g)\in J$ there exists $i\in \llbracket1,m\rrbracket$ such that $f(\alpha_i)+g(\beta_i)\neq 0$.
\end{definition}

If $g(x)\in K^{\ast}\subset K[x]$, then for each $f\in K[x]$ such that $f(\alpha_i)+g(\beta_i)=0$ for all $i\in \llbracket1,m\rrbracket$ we have $f(x)+g(0)\in \prod_{i=1}^m (x-\alpha_i)K[x]$, hence either $\deg(f)=0$ or $\deg(f)\ge m$. Thus $J$-generic and $J\setminus\{0\}$-generic are equivalent notions.

\begin{lemma}\label{F12}
Let $m\in\mathbb N^{\ast}$. Let $K$ be a field with $|K|\ge m$ and let $J$ be a subset of $\llbracket0,m-1\rrbracket$. Let $\bigl((\alpha_1,\ldots,\alpha_m),(\beta_1,\ldots,\beta_m)\bigr)\in\mathbb D_{1,m}(K)^2$ be $J$-generic. Then the following properties hold.

\medskip
{\bf (1)} For each subset $I\subset J$, $\bigl((\alpha_1,\ldots,\alpha_m),(\beta_1,\ldots,\beta_m)\bigr)$ is $I$-generic.

\smallskip
{\bf (2)} If $J^*_m:=\{0\}\cup\{i\in \llbracket1,m-1\rrbracket|\llbracket1,m-i-1\rrbracket\subset J\}$, then the interchanged pair $\bigl((\beta_1,\ldots,\beta_m),(\alpha_1,\ldots,\alpha_m)\bigr)\in\mathbb D_{1,m}(K)^2$ is $J^*_m$-generic.

\smallskip
{\bf (3)} If $|K|\ge m+1$, $\alpha_{m+1}\in K\setminus\{\alpha_i|i\in \llbracket1,m\rrbracket\}$, and $\beta_{m+1}\in K\setminus\{\beta_i|i\in \llbracket1,m\rrbracket\}$, then $\bigl((\alpha_1,\ldots,\alpha_{m+1}),(\beta_1,\ldots,\beta_{m+1})\bigr)\in\mathbb D_{1,m+1}(K)^2$ is weakly $J$-generic. 

\smallskip
{\bf (4)} Let $h\in K[x]$ be the Lagrange interpolation polynomial of degree at most $m-1$ such that $h(\alpha_i)=\beta_i$ for all $i\in \llbracket1,m\rrbracket$. If $1\in J$ and $m\ge 3$, then $\deg(h)=m-1$ and $|K|>m$.
\end{lemma}

\begin{proof}
Parts (1) to (3) follow directly from definitions. 

We prove part (4). If $|K|=m\ge 3$, then $\deg(h)\le m-2$ by Hermite's criterion of \cite{Di}, Sect.\ 1, Subsect.\ 11. So for part (4) it suffices to show that the assumption $\deg(h)\le m-2$ leads to a contradiction. We have $1\le\deg\bigl(h(x)\bigr)+\deg(-x)\le m-1$ and $\deg(-x)=1\in J$. Thus the $J$-generic hypothesis gives that there exists $i\in \llbracket1,m\rrbracket$ such that $h(\alpha_i)-\beta_i\neq 0$, a contradiction.\end{proof}

\begin{lemma}\label{L22}
Let $m\in\mathbb N^{\ast}\setminus\{1\}$. Let $J$ be a non-empty subset of $\llbracket1,m-1\rrbracket$ and let $K$ be a field. For $\bigl((\alpha_1,\ldots,\alpha_m),(\beta_1,\ldots,\beta_m)\bigr)\in\mathbb D_{1,m}(K)^2$ we consider the following statements.

\medskip
{\bf (1)} The pair of $m$-tuples $\bigl((\alpha_1,\ldots,\alpha_m),(\beta_1,\ldots,\beta_m)\bigr)$ is $J$-generic.

\smallskip
{\bf (2)} For each $r\in J$, the determinant of the $m\times m$ matrix 
$$\Delta_r:=\begin{bmatrix} 
1 & \alpha_1 & \cdots & \alpha_1^{m-r-1} & \beta_1& \cdots & \beta_1^r\\
1 & \alpha_2 & \cdots & \alpha_2^{m-r-1} & \beta_2& \cdots & \beta_2^r\\
\cdots & \cdots & \cdots & \cdots & \cdots & \cdots & \cdots\\
1 & \alpha_m & \cdots & \alpha_m^{m-r-1} & \beta_m& \cdots & \beta_m^r\\
\end{bmatrix}$$
with entries in $K$ is non-zero.

\medskip
Then $(2)\Rightarrow (1)$. Moreover, if $J=\llbracket1,\max(J)\rrbracket$ then $(1)\Leftrightarrow (2)$.
\end{lemma}

\begin{proof}
Assume statement (2) holds. Let $r\in J$. We show that the assumption that there exists $(f,g)\in\ (K[x])^2$ such that $\deg(f)\le m-r-1$, $\deg(g)=r$, and $f(\alpha_i)+g(\beta_i)=0$ for each $i\in \llbracket1,m\rrbracket$ leads to a contradiction. This assumption implies that there exists $(\delta_0,\ldots,\delta_{m-1})\in K^{m-1}\times K^{\ast}$ such that we have 
\begin{equation}\label{EQ45}
\sum_{j=0}^{m-r-1} \delta_j\alpha_i^j+\sum_{j=1}^{r} \delta_{m-r-1+j}\beta_i^j=0
\end{equation}
for each $i\in \llbracket1,m\rrbracket$. Thus $\det(\Delta_r)=0$, a contradiction. Hence $(2)\Rightarrow (1)$.

\phantomsection{Assume statement (1) holds and $J=\llbracket1,\max(J)\rrbracket$. If there exists $r\in J$ such that $\Delta_r=0$, then there exists a non-zero $m$-tuple $(\delta_0,\ldots,\delta_{m-1})\in K^m$ such that Equation (\ref{EQ45}) holds for each $i\in \llbracket1,m\rrbracket$. 

If $\delta_{\m-r}=\cdots=\delta_{m-1}=0$, then the matrix of size $(m-r)\times (m-r)$ obtained from $\Delta_r$ by removing its last $r$ rows and columns has a zero determinant which is a non-zero Vandermonde determinant $W(\alpha_1,\ldots,\alpha_{m-r})$, a contradiction. So for the largest $r_0\in \llbracket0,m-1\rrbracket$ with $\delta_{r_0}\neq 0$ we have $r_0\in \llbracket m-r,m-1\rrbracket$.}\label{PH97}

Let $f(x):=\sum_{j=0}^{m-r-1} \delta_jx^j\in K[x]$ and $g(x):=\sum_{j=1}^{r_0+r+1-m} \delta_{m-r-1+j}x^j\in K[x]$. We have $\deg(g)=r_0+r+1-m\in \llbracket1,r\rrbracket\subset J$, $1\le\deg(f)+\deg(g)\le r_0\le m-1$, and $f(\alpha_i)+g(\beta_i)=0$ for each $i\in \llbracket1,m\rrbracket$, which contradicts (1). Hence $(1)\Rightarrow (2)$ if $J=\llbracket1,\max(J)\rrbracket$.
\end{proof}

\begin{proposition}\label{PR28}
Let $m\in \mathbb N^{\ast}\setminus\{1\}$. Let $K$ be a field with $|K|\ge m$. Then the following properties hold.

\medskip
{\bf (1)} Let $\bigl((\alpha_1,\ldots,\alpha_m),(\beta_1,\ldots,\beta_{m-1})\bigr)\in\mathbb D_{1,m}(K)\times \mathbb D_{1,m-1}(K)$. If $|K|\ge m+1$, then there exist $\beta_m\in K\setminus\{\beta_i|i\in \llbracket1,m-1\rrbracket\}$ and $f\in K[x]$ such that $\deg(f)=m-1$ and $f(\alpha_i)=\beta_i$ for each $i\in \llbracket1,m\rrbracket$.

\smallskip
{\bf (2)} If $|K|\ge m+1$, then there exists $\bigl((\alpha_1,\ldots,\alpha_m),(\beta_1,\ldots,\beta_m)\bigr)\in\mathbb D_{1,m}(K^{\ast})^2$ which is $\llbracket1,m-1\rrbracket$-generic.

\smallskip
{\bf (3)} Suppose that $|K|=m\ge 3$. Let $(\alpha_1,\ldots,\alpha_m)\in\mathbb D_{1,m}(K)$. Let $h\in K[x]$ be the Lagrange interpolating polynomial such that $\deg(h)\le m-1$, $h(\alpha_i)=\alpha_i$ for each $i\in \llbracket1,m-2\rrbracket$, $h(\alpha_{m-1})=\alpha_m$, and $h(\alpha_m)=\alpha_{m-1}$. Then $\deg(h)=m-2$.

\smallskip
{\bf (4)} If $|K|=m\ge 3$, then there exists $\bigl((\alpha_1,\ldots,\alpha_m),(\beta_1,\ldots,\beta_m)\bigr)\in\mathbb D_{1,m}(K)^2$ which is weakly $\llbracket1,|K|-2\rrbracket$-generic.

\smallskip
{\bf (5)} If $m\ge 3$, then there exists $\bigl((\alpha_1,\ldots,\alpha_m),(\beta_1,\ldots,\beta_m)\bigr)\in\mathbb D_{1,m}(K)^2$ which is $\llbracket1,|K|-1\rrbracket$-generic iff $|K|>m$.
\end{proposition}

\begin{proof}
For part (1) let $h(x)\in K[x]$ be the Lagrange interpolating polynomial of degree at most $m-2$ such that $h(\alpha_i)=\beta_i$ for each $i\in \llbracket1,m-1\rrbracket$. Let $\gamma\in K$ be such that for $f(x):=h(x)+\gamma\prod_{i=1}^{m-1}(x-\alpha_i)$ we have $f(\alpha_m)\in K\setminus\{\beta_i|i\in \llbracket1,m-1\rrbracket\}$; there exist $|K|-m+1\ge 2$ possibilities for $\gamma$ and thus we can choose $\gamma\in K^{\ast}$; so $\deg(f)=m-1$ and part (1) holds.

For part (2), let $\delta\in K^{\ast}$ be an element whose multiplicative order is at least $m$. E.g., if $K$ is finite then we can take $\delta$ to be a generator of the cyclic multiplicative group $K^{\ast}$. For $i\in \llbracket1,m\rrbracket$ let $\alpha_i:=\delta^i$ and $\beta_i:=\delta^{-i}$. For $r\in \llbracket1,m-1\rrbracket$, let the matrix $\Delta_r$ be as in Lemma \ref{L22}. As $\det(\Delta_r)$ is $\delta^{\frac{(m-r-1)(m-r)-r(r+1)}{2}}$ times the non-zero Vandermonde determinant $W(1,\delta^1,\ldots,\delta^{m-r-1},\delta^{-1},\ldots,\delta^{-r})$, it follows that $\det(\Delta_r)\neq 0$. So part (2) follows from Lemma \ref{L22}.

For part (3) we have $h(x)=x+h_0(x)\prod_{i=1}^{m-2}(x-\alpha_i)$ for a suitable $h_0(x)\in K[x]$ of degree $0$ or $1$ and hence of degree $0$ by Hermite's Criterion, see \cite{Di}, Sect.\ 1, Subsect.\ 11. So $\deg(h)=m-2$ and part (3) holds. 

Part (4) follows from part (2) applied to $m-1$ and Lemma \ref{F12}(3).

Part (5) follows from part (2) and Lemma \ref{F12}(4).\end{proof}

We have the following direct consequence of Lemmas \ref{F11.5}(1) and \ref{F12}(4).

\begin{corollary}\label{L23}
Let $m\in \mathbb N^{\ast}\setminus\{1,2\}$. Let $K$ be a field with $|K|> m$. Let the pair $\bigl((\alpha_1,\ldots,\alpha_m),(\beta_1,\ldots,\beta_m)\bigr)\in\mathbb D_{1,m}(K)^2$ be $\llbracket1,m-1\rrbracket$-generic. We consider the pair $(\underline{P},\underline{Q})=\bigl((P_1,\ldots,P_m),(Q_1,\ldots,Q_m)\bigr)\in\mathbb D_{2,m}(K)^2$ defined by $P_i:=(\alpha_i,0)$ and $Q_i:=(\beta_i,0)$ for each $i\in\llbracket1,m\rrbracket$. Then $\pi_{\underline{P},\underline{Q}}^{\S,\le 2}=(m-1)^2$.
\end{corollary}

We recall \cite{AM}, Main Thm.\ (1.1) and an application of the proofs of \cite{vdK}, Thms.\ 1 and 2 in the following form.

\begin{theorem}\label{T14.5}
Let $(g_1,g_2)\in K[x]^2$ be such that we have $K[x]=K[g_1,g_2]$ and let $s_1:=\deg(g_1)\in\mathbb N^{\ast}$ and $s_2:=\deg(g_2)\in\mathbb N^{\ast}$. Then $s_1\mid s_2$ or $s_2\mid s_1$ provided one of the following conditions holds.

\medskip
{\bf (1) (Abhyankar--Moh)} The greatest common divisor of $s_1$ and $s_2$ is not divisible by $\char(K)$.

\smallskip
{\bf (2)} There exists $(f_1,f_2)\in K[x_1,x_2]^2$ such that $a:=\e(f_1,f_2)$ is in $\GA_2(K)$ and we have $g_1(x)=f_1(x,0)$ and $g_2=f_2(x,0)$.
\end{theorem}

\begin{proof}
If condition (1) holds, then this is a restatement of \cite{AM}, Main Thm.\ (1.1). 

If condition (2) holds, then this is also well-known being the essence of the proofs of \cite{vdK}, Thms.\ 1 and 2 which is captured, for instance, by the algorithm in \cite{F1}, Sect.\ 2, and which gives that $(g_1,g_2)$ is obtained from a pair $(\alpha x+\beta,\gamma x+\delta)$ with $(\alpha,\beta,\gamma,\delta)\in K^4$ such that $\alpha\gamma\neq 0$ by performing a finite number of one of the following types (i) to (iii) of replacements, the first type (i) being performed at the most one time as the very last replacement. The three types are: replace a pair $(h_1,h_2)\in K[x]^2$ by (i) $(h_1+\alpha_1,h_2)$ with $\alpha_1\in K^{\ast}$, or by (ii) $(h_2,h_1)$, or by (iii) $(h_1,h_2+h\circ h_1)$ with $h\in K[x]$ of degree at least $2$. We recall the standard proof.  

We can assume that $s_1>s_2$. By replacing $(g_1,g_2)$ with $(g_1+\alpha_2, g_2)$ where $\alpha_2\in K$, we can assume that $(g_1,g_2)$ is obtained from a pair $(\alpha x+\beta,\gamma x+\delta)$ by performing alternatively replacements of only types (ii) and (iii). For each pair $(\hbar_1,\hbar_2)\in K[x]^2$ with $\bigl(\deg(\hbar_1),\deg(\hbar_2)\bigr)\in (\mathbb N^{\ast})^2\setminus\{(1,1)\}$ and obtainable from a pair $(\alpha x+\beta,\gamma x+\delta)$ by performing alternatively replacements of only types (ii) and (iii), let $l_{\hbar_1,\hbar_2}\in\mathbb N^{\ast}$ be the smallest number of replacements required to obtain it.

We show by induction on $l\in\mathbb N^{\ast}$ that for each such pair $(\hbar_1,\hbar_2)$ with $l_{\hbar_1,\hbar_2}=l$ we have $\frac{\deg(\hbar_1)}{\deg(\hbar_2)}\in\mathbb N^{\ast}\setminus\{1\}$ if $l$ is even and $\frac{\deg(\hbar_2)}{\deg(\hbar_1)}\in\mathbb N^{\ast}\setminus\{1\}$ if $l$ is odd. 

By the minimality of $l$, the first replacement is of type (iii) and involves an $h$ of degree at least $2$. Thus if $l=1$ and $(\hbar_1,\hbar_2)=\bigl(\alpha x+\beta,\gamma x+\delta+h(\alpha x+\beta)\bigr)$ with $h\in K[x]$, then $\frac{\deg(\hbar_2)}{\deg(\hbar_1)}=\deg(h)\in\mathbb N^{\ast}\setminus\{1\}$. So the base of the induction holds. For $l\ge 2$, the inductive passage from $l-1$ to $l$ goes as follows. If $l$ is even, then $\ell_{l_2,l_1}=l-1$ is odd and $\frac{\deg(\hbar_1)}{\deg(\hbar_2)}\in\mathbb N^{\ast}\setminus\{1\}$ by the inductive assumption applied to $(\hbar_2,\hbar_1)$. If $l$ is odd and $(\hbar_1,\hbar_2)=(\hbar_3,\hbar_4+h\circ \hbar_3)$ with $h\in K[x]$ of degree at least $2$ and $l_{\hbar_3,\hbar_4}=l-1$, then $l-1$ is even and by the inductive assumption applied to $(\hbar_3,\hbar_4)$ we have $\frac{\deg(\hbar_3)}{\deg(\hbar_4)}\in\mathbb N^{\ast}\setminus\{1\}$, hence $\frac{\deg(\hbar_2)}{\deg(\hbar_1)}=\deg(h)\frac{\deg(\hbar_3)}{\deg(\hbar_4)}\in\mathbb N^{\ast}\setminus\{1\}$. This ends the inductive step, the induction, and the proof.\end{proof}

\begin{theorem}\label{T15.6}
Let $m\in \mathbb N^{\ast}\setminus\{1\}$. Let $K$ be a field with $|K|> m\ge 3$. We consider a pair $\bigl((\alpha_1,\ldots,\alpha_m),(\beta_1,\ldots,\beta_m)\bigr)\in\mathbb D_{1,m}(K)^2$ which is $\llbracket1,\lfloor\frac{2m-1}{3}\rfloor\rrbracket$-generic. Let $F(x):=\prod_{i=1}^m (x-\alpha_i)$. Let $h(x)\in K[x]$ be the Lagrange interpolation polynomial such that $h(\alpha_i)=\beta_i$ for all $i\in \llbracket1,m\rrbracket$. Let $(F_1,F_2)\in K[x]^2$ be such that for the pair $(f,g):=(h+F_1F,F_2F)$ we have $\deg(f)<\deg(g)$ and $K[x]=K[f,g]$. If $\char(K)\mid \deg(g)$, then we assume that there exists $(f_1,f_2)\in K[x_1,x_2]^2$ such that $a:=\e(f_1,f_2)$ is in $\GA_2(K)$ and we have $f(x)=f_1(x,0)$ and $g(x)=f_2(x,0)$. Writing $m=3r+\epsilon$ with $(r,\epsilon)\in\mathbb N^{\ast}\times\{0,1,2\}$, the following properties hold.

\medskip
{\bf (1)} Suppose that $\epsilon=0$. Then we have $\deg(g)\ge 4r^2=\frac{4m^2}{9}$. If $\deg(g)=4r^2$, then either $r=1$ and $\deg(f)=2$ or $r\ge 2$ and $2r\mid\deg(f)\mid 4r^2$ with $\deg(f)\notin\{2r,4r^2\}$.

\smallskip
{\bf (2)} Suppose that $\epsilon=1$. Then we have $\deg(g)\ge 4r^2+4r>\frac{4m^2}{9}$. If $r=1$ and $\deg(g)=8$, then $\deg(f)=4$ and if $r\ge 2$ and $\deg(g)= 4r^2+4r$, then either $2r+2\mid\deg(f)\mid 4r^2$ with $\deg(f)\notin\{2r+2,4r^2+4r\}$ or $4r\mid\deg(f)\mid 4r^2+4r$ with $\deg(f)\neq 4r^2+4r$.

\smallskip
{\bf (3)} Suppose that $\epsilon=2$. Then we have $\deg(g)\ge 4r^2+6r+2>\frac{4m^2}{9}$. If $r=1$ and $\deg(g)=12$, then $\deg(f)\in\{4,6\}$ and if $r\ge 2$ and $\deg(g)=4r^2+6r+2$, then either $2r+2\mid\deg(f)\mid 4r^2+6r+2$ with $\deg(f)\notin\{2r+2,4r^2+6r+2\}$ (so $2r+1$ is not a prime) or $4r+2\mid\deg(f)\mid 4r^2+6r+2$ with $\deg(f)\neq 4r^2+6r+2$.

\smallskip
{\bf (4)} Suppose that $K$ has a primitive $m$-th root of unity $\gamma$ and for each $\in\llbracket1,m\rrbracket$ we have $\alpha_i=\gamma^i$ and $\beta_i=\gamma^{-i}$. Then the following properties hold.

\medskip\noindent
{\bf (4.a)} We have an inequality
$$\deg(g)\ge\bigl\lfloor\frac{m+2}{2}\bigr\rfloor(m-1)=\begin{cases} \frac{m^2-1}{2}\quad\quad\quad\quad {\rm if}\; m\;\textup{is odd}\\
\frac{m^2+m-2}{2}\quad\quad\;\;\, {\rm if}\; m\;\textup{is even}.\end{cases}$$

\noindent
{\bf (4.b)} If $m=5$, then $\deg(g)\ge 16$.

\smallskip\noindent
{\bf (4.c)} If $m=6$, then $\deg(g)\ge 24$. 

\smallskip\noindent
{\bf (4.d)} If $m=6$ and there exists $(f_1,f_2)\in K[x_1,x_2]^2$ such that $a:=\e(f_1,f_2)$ is in $\GA_2(K)$ and we have $f(x)=f_1(x,0)$ and $g(x)=f_2(x,0)$, then $\deg(g)\ge 25$. 
\end{theorem}

\begin{proof}
We have $\deg(h)=m-1$ by Lemma \ref{F12}(4). By denoting $n:=\deg(g)$, we have $n>\deg(f)\ge m-1=\deg(h)$. Thus $\frac{n}{\deg(f)}\in\mathbb N^{\ast}\setminus\{1\}$ by Theorem \ref{T14.5}. 

Let $(h_0,\ldots,h_{n-1})\in K[x]^n$ be an $n$-tuple of polynomials constructed from $(f,g)$ as follows. If $\char(K)\nmid n$, then they are constructed from $(f,g)$ in \cite{Ri}, (*) of Sect.\ 2, cf.\ also \cite{Kan}, Thm.\ R. If $\char(K)\mid n$, then the fact that they can be constructed in the same manner follows from the proof of \cite{Kan}, Prop.\ 1 as in the mentioned proof one can replace the usage of Theorem \ref{T14.5}(1) by the usage of Theorem \ref{T14.5}(2) which applies here by our hypotheses.

Recall that $h_0:=g$, we have 
\begin{equation}\label{EQ45.1}
h_i-f^i\in gK[g]+\sum_{j=1}^{i-1}f^jK[g]\subset K[x]
\end{equation}
for each $i\in \llbracket1,n-1\rrbracket$, and the set $\{\deg(h_i)|i\in \llbracket0,n-1\rrbracket\}$ is a complete residue system modulo $n$; thus $h_i=f^i$ for each $i\in\bigl\llbracket1,\frac{n}{\deg(f)}-1\bigr\rrbracket$.\footnote{For $i\in \llbracket1,n-1\rrbracket$, the $h_i$ in \cite{Kan}, Thm.\ R is uniquely determined up to a constant in $K$. Equation (\ref{EQ45.1}) determines $h_i$ uniquely by using $gK[g]$ instead of $K[g]$.}

\phantomsection{There exists $q\in\mathbb N^{\ast}$, a strictly increasing sequence $\bigl(\eta(i)\bigr)_{i\in \llbracket0,q+1\rrbracket}$ in $\mathbb N$, and $(l_1,\ldots,l_q)\in (\mathbb N^{\ast}\setminus\{1\})^q$, such that by defining $s_i:=\deg(h_{\eta(i)})$ for each $i\in \llbracket0,q\rrbracket$, the following properties hold by \cite{Kan}, Sect.\ 2, Steps 2 and 3 and Sect.\ 3, Prop.\ 8 and its proof (if $\char(K)\mid n$, cf. proof of \cite{Kan}, Prop.\ 1).}\label{PH98}

\medskip
{\bf (i)} We have $\eta(0)=0$, $\eta(1)=1$, $\eta(q+1)=n$, and $\eta(i+1)=l_i\eta(i)$ for $i\in \llbracket1,q\rrbracket$.

\smallskip
{\bf (ii)} We have $s_i=\prod_{j=i+1}^q l_i$ (equivalently, $s_i=\frac{n}{\eta(i+1)}$) for $i\in \llbracket0,q\rrbracket$.

\medskip
Thus $$s_0=\prod_{j=1}^q l_i=\eta(q+1)=n\ge m\ge 3,$$ 
$$s_1=\deg(h_1)=\deg(f)\ge\deg(h)=m-1\ge\Bigl\lceil\frac{2m}{3}\Bigr\rceil\ge\Bigr\lfloor\frac{2m}{3}\Bigr\rfloor\ge 2,$$ and for the strictly decreasing sequence $(s_i)_{i\in \llbracket0,q\rrbracket}$ we have the following divisibilities $1=s_q\mid s_{q-1}\mid\cdots\mid s_1\mid s_0=n$. Thus there exists a unique $i_0\in \llbracket1,q-1\rrbracket$ such that $s_{i_0}\ge \frac{2m}{3}$ and $s_{i_0+1}<\frac{2m}{3}$.\footnote{The function $\mathfrak f:[0,1]\rightarrow\mathbb R$, $\mathfrak f(x):=\min\bigl(x,(2-x)(1-x),x(1-\frac{x}{2})\bigr)$, attains its maximum $\frac{4}{9}$ only at $\frac{2}{3}$.} We have $s_{i_0}\mid s_1=\deg(f)$.

Moreover, for each $i\in\llbracket1,q\rrbracket$ we have an identity
\begin{equation}\label{EQ46.4}
K[x]=K[f,g]=K[h_{\eta(i-1)},h_{\eta(i)}]
\end{equation}
(the case $i=1$ follows directly from the construction of $h_{\eta(1)}$, for $i=2$ see \cite{Kan}, p.\ 420 if $\char(K)\nmid n$ and see the proof of  \cite{Kan}, Prop.\ 1 if $\char(K)\mid n$; the case $i\ge 3$ follows from the mentioned case by replacing the pair $(f,g)$ by $\bigl(h_{\eta(i-3)},h_{\eta(i-2)}\bigr)$).

In this paragraph we assume that $\eta(i_0+1)=\frac{n}{s_{i_0}}\ge \lfloor\frac{2m+2}{3}\rfloor$. Then for $n=s_{i_0}\frac{n}{s_{i_0}}$ we have $n\ge 4r^2$ if $\epsilon=0$, $n\ge 4r^2+4r+1$ if $\epsilon=1$, and $n\ge 4r^2+8r+4$ if $\epsilon=2$; so parts (2) and (3) hold. If $\epsilon=0$ and $n=4r^2$, then $s_{i_0}=\frac{n}{s_{i_0}}=2r$ and thus for $r=1$ we get that $\deg(f)=2$ and for $r\ge 2$ we get that $2r\mid\deg(f)\mid 4r^2$ and $\deg(f)\notin\{2r,4r^2\}$, hence part (1) also holds.

Therefore we can assume that $\eta(i_0+1)=\frac{n}{s_{i_0}}\le \lfloor\frac{2m-1}{3}\rfloor$; so
\begin{equation}\label{EQ46}
\frac{n}{s_{i_0}}\le 2r+\Bigl\lfloor\frac{2\epsilon-1}{3}\Bigr\rfloor=2r+\epsilon-1.
\end{equation}

Let $\bigl(\delta_1,\ldots,\delta_{\eta(i_0+1)-1}\bigr)\in K^{\eta(i_0+1)-1}$ be such that 
$$h_{\eta(i_0+1)}-f^{\eta(i_0+1)}-\sum_{j=1}^{\eta(i_0+1)-1} \delta_j f^j\in gK[x].$$ For $i\in [1,m]$, as $g(\alpha_i)=0$, by evaluating at $\alpha_i$ we get an identity 
\begin{equation}\label{EQ46.5}
h_{\eta(i_0+1)}(\alpha_i)-\beta_i^{\eta(i_0+1)}-\sum_{j=1}^{\eta(i_0+1)-1} \delta_j\beta_i^j=0.
\end{equation}
From this and the $\llbracket1,\lfloor\frac{2m-1}{3}\rfloor\rrbracket$-generic hypothesis we get that $\eta(i_0+1)+s_{i_0+1}\ge m$.\footnote{If the pair $\bigl((\alpha_1,\ldots,\alpha_m),(\beta_1,\ldots,\beta_m)\bigr)\in\mathbb D_{1,m}(K)^2$ is $\llbracket1,m-1\rrbracket$-generic, then we have $s_i+\eta(i)\ge m$ for each $i\in\llbracket1,q\rrbracket$.} Thus 
\begin{equation}\label{EQ47}
n= s_{i_0}\eta(i_0+1)\ge s_{i_0}(m-s_{i_0+1})=s_{i_0}\Bigl(m-\frac{s_{i_0}}{l_{i_0+1}}\Bigr)\ge s_{i_0}\left(m-\frac{s_{i_0}}{2}\right)
\end{equation} 
by properties (i) and (ii). To show that Relations (\ref{EQ47}) implies parts (1) to (3) we consider two cases as follows.

{\bf Case 1: we have $s_{i_0}>\frac{4m}{3}$.}\footnote{In this case we have $l_{i_0+1}\ge 3$.} From this, the inequality $m-s_{i_0+1}>m-\frac{2m}{3}$, and Relations (\ref{EQ47}) we get that $n\ge\lceil\frac{4m+1}{3}\rceil\lceil\frac{m+1}{3}\rceil$. Thus $n\ge 4r^2+5r+1$ if $\epsilon=0$, $n\ge 4r^2+6r+2$ if $\epsilon=1$, and $n\ge 4r^2+7r+3$ if $\epsilon=2$. Hence parts (1) to (3) hold.

{\bf Case 2: we have $s_{i_0}\le\frac{4m}{3}$.} We consider the function $\mathfrak f_{\epsilon}:\mathbb R\rightarrow\mathbb R$ defined by the rule $\mathfrak f_{\epsilon}(x):=x(3r+\epsilon-\frac{x}{2})$. Relations (\ref{EQ47}) gives $n\ge\mathfrak f_{\epsilon}(s_{i_0})$. 

Assume that $\epsilon=0$. Then $s_{i_0}\in \llbracket2r,4r\rrbracket$. The minimum value of $\mathfrak f_0$ on the interval $\llbracket2r,4r\rrbracket$ is $\mathfrak f_0(2r)=\mathfrak f_0(4r)=4r^2$, hence $n\ge 4r^2$. If $n=4r^2$, then Relations (\ref{EQ46}) gives $s_{i_0}\ge \lceil\frac{4r^2}{2r-1}\rceil=2r+2$, thus $s_{i_0}=4r$ and this is possible only when $r>1$ as $i_0\ge 1$ implies $s_{i_0}<n$. If $r\ge 2$ and we have $n=4r^2$ and $s_{i_0}=4r$, then $4r\mid\deg(f)\mid 4r^2$ with $\deg(f)\neq 4r^2$. Hence part (1) holds. 

Assume that $\epsilon=1$. Then $s_{i_0}\in \llbracket2r+1,4r+1\rrbracket$. The minimum value of $\mathfrak f_1$ on the interval $\llbracket2r+2,4r\rrbracket$ is $\mathfrak f_1(2r+2)=\mathfrak f_1(4r)=4r^2+4r$. If $s_{i_0}=2r+1$ (resp.\ $4r+1$), then $2\nmid s_{i_0}$ implies that $l_{i_0+1}\ge 3$ and Relations (\ref{EQ47}) gives 
$$n\ge (2r+1)\Bigl(3r+1-\Bigl\lfloor\frac{2r+1}{3}\Bigr\rfloor\Bigr)\ge (2r+1)(3r+1-r)=4r^2+4r+1$$ 
(resp.\ $n\ge (4r+1)(3r+1-\lfloor\frac{4r+1}{3}\rfloor)\ge (4r+1)(3r+1-2r+1)=4r^2+9r+2$); so $n>4r^2+4r$. If $n=4r^2+4r$, then Relations (\ref{EQ46}) gives $s_{i_0}\ge \lceil\frac{4r^2+4r}{2r}\rceil=2r+2$, hence $s_{i_0}\in\{2r+2,4r\}$; if $r=1$, then $s_{i_0}=4=\deg(f)$ and $n=8$ and if $r\ge 2$ then either $s_{i_0}=2r+2\mid\deg(f)\mid 4r^2+4r$ with $\deg(f)\notin\{2r+2,4r^2+4r\}$ or $s_{i_0}=4r\mid\deg(f)\mid 4r^2+4r$ with $\deg(f)\neq 4r^2+4r$. So part (2) holds. 

If $\epsilon=2$, then $s_{i_0}\in \llbracket2r+2,4r+2\rrbracket$. The minimum value of $\mathfrak f_2$ on the interval $[2r+2,4r+2]$ is $\mathfrak f_2(2r+2)=\mathfrak f_2(4r+2)=4r^2+6r+2$. If $n=4r^2+6r+2$, then Relations (\ref{EQ46}) gives $s_{i_0}\ge \lceil\frac{4r^2+6r+2}{2r+1}\rceil=2r+2$, hence $s_{i_0}\in\{2r+2,4r+2\}$; if $r=1$, then $s_{i_0}=\deg(f)\in\{4,6\}$ and $n=12$ and if $r\ge 2$ then we have either $s_{i_0}=2r+2\mid\deg(f)\mid 4r^2+6r+2$ with $\deg(f)\notin\{2r+2,4r^2+6r+2\}$ or $s_{i_0}=4r+2\mid\deg(f)\mid 4r^2+6r+2$ with $\deg(f)\neq 4r^2+6r+2$. So part (3) holds.

For part (4), we have $h(x)=x^{m-1}$ and $F(x)=x^m-1$. The proof of Proposition \ref{PR28}(2) gives that $\bigl((\alpha_1,\ldots,\alpha_m),(\beta_1,\ldots,\beta_m)\bigr)\in\mathbb D_{1,m}(K)^2$ is $\llbracket1,m-1\rrbracket$-generic. Let $i_1\in\llbracket1,m\rrbracket$ be the unique integer such that $s_{i_1}\ge m-1$ and $s_{i_1+1}\le m-2$. If $\eta(i_1+1)\ge m-1$, then $n=s_{i_1}\eta(i_1+1)\ge (m-1)^2$. Thus we can assume that $\eta(i_1+1)\le m-2$; hence $\frac{s_0}{s_1}=\eta(2)\le \eta(i_1+1)\le m-2$. 

For each $\iota\in\llbracket2,q\rrbracket$ such that $\eta(\iota)\le m-2$, an argument similar to the one of Equation (\ref{EQ46.5}) gives that there exists $\bigl(\delta_{\iota,1},\ldots,\delta_{\iota,\eta(i_0+1)-1}\bigr)\in K^{\eta(i_0+1)-1}$ such that for each $i\in\llbracket1,m\rrbracket$ we have 
\begin{equation*}\label{EQ47.3}
h_{\eta(\iota)}(\alpha_i)=\beta_i^{\eta(\iota)}+\sum_{j=1}^{\eta(\iota)-1} \delta_{\iota,j}\beta_i^j=\alpha_i^{m-\eta(\iota)}+\sum_{j=1}^{\eta(\iota)-1} \delta_{\iota,j}\alpha_i^{m-j}.
\end{equation*}
Hence there exists a unique polynomial $g_{\eta(\iota)}\in 1+xK[x]$ such that $\deg(g_{\iota})\le\eta(\iota)-1$ and $h_{\eta(\iota)}(x)-x^{m-\eta(\iota)}g_{\eta(\iota)}(x)$ is divisible by $x^m-1$. In particular, if $s_{\iota}\le m-1$, then $h_{\eta(\iota)}(x)=x^{m-\eta(\iota)}g_{\eta(\iota)}(x)$ and $\deg(g_{\eta(\iota)})=s_{\iota}+\eta(\iota)-m$ and thus $s_{\iota}+\eta(\iota)\ge m$.

For part (4.a), we show that the assumption $\eta(i_1+2)\le m-2$ leads to a contradiction. 
From this, by taking $\iota\in\{i_1+1,i_1+2\}$, the fact that $s_{\iota}\le m-2$, gives that
$$h_{\eta(\iota)}(x)=x^{m-\eta(\iota)}g_{\eta(\iota)}(x)$$
is divisible by $x^2$. From this and Equation (\ref{EQ46.4}) applied to $i=i_1+2$ we get that the two polynomials $h_{\eta(i_1+1)}$ and $h_{\eta(i_1+2)}$ in $x^2K[x]$ generate $K[x]$, a contradiction.

Therefore we have $\eta(i_1+2)\ge m-1$. 

We show that the inequality $n\ge\lfloor\frac{m+2}{2}\rfloor(m-1)$ holds. If $s_{i_1+1}\ge\frac{m+1}{2}$, then we have $n=s_{i_1+1}\eta(i_1+2)\ge \lceil\frac{m+1}{2}\rceil(m-1)$. Similarly if $\eta(i_1+1)\ge\frac{m+1}{2}$, then $n=s_{i_1}\eta(i_1+1)\ge \lceil\frac{m+1}{2}\rceil(m-1)$. 

Thus we can assume that $s_{i_1+1}<\frac{m+1}{2}$ and $\eta(i_1+1)<\frac{m+1}{2}$. From this and the inequality $s_{i_1+1}+\eta(i_1+1)\ge m$ we get that $s_{i_1+1}=\eta(i_1+1)=\frac{m}{2}$. Hence $m$ is even, $\deg(g_{\eta(i_1+1)})=0$, and $h_{\eta(i_1+1)}=x^{\frac{m}{2}}$. 

If $s_{i_1}\ge m+1$ then $n=s_{i_1}\eta(i_1+1)\ge  \frac{m(m+1)}{2}\ge\lfloor\frac{m+2}{2}\rfloor(m-1)$. Thus we can assume that $s_{i_1}\le m$. As $s_{i_1}\ge 2s_{i_1+1}=m$,  we conclude that $s_{i_1}=m$. Hence $h_{\eta(i_1)}(x)=\gamma_{\eta(i_1)}(x^m-1)+x^{m-\eta(\iota)}g_{\eta(\iota)}(x)$ with $\gamma_{\eta(i_1)}\in K^{\ast}$. From this and Equation (\ref{EQ46.4}) applied to $i=i_1+1$ we get that two polynomials $h_{\eta(i_1)}$ and $h_{\eta(i_1+1)}$ in $K[x]$ generate $K[x]$, thus the two polynomials $x^{m-\eta(\iota)}g_{\eta(\iota)}(x)$ and $x^{\frac{m}{2}}$ in $x^2K[x]$ generate $K[x]$, a contradiction.

We conclude that the inequality $n\ge\lfloor\frac{m+2}{2}\rfloor(m-1)$ holds. So part (4.a) holds.

For part (4.b), if $i_1\ge 2$, then $s_0\ge 2s_1\ge 4s_2\ge 4(m-1)=16$. Thus we can assume that $i_1=1$; so $s_2\le m-2$ and $\eta(2)\le m-2=3$. From the inequalities $2\le m-\eta(2)\le s_2\le m-2=3$ we get that $s_2\in\{2,3\}$. 

If $s_2=2$, then $h_{\eta(2)}(x)=x^2$ and from $K[x]=K[x^4+(x^5-1)F_1(x),x^2]$ we get that $s_1\ge 6$; as $3=m-s_2\le\eta(2)\le 3$, we get that $\eta(2)=3$ and hence $n=s_1\eta(2)\ge 18$.

Assume now that $s_2=3$; so $s_1=6$ or $s_1\ge 9$. If $s_1\ge 9$, then $s_0\ge 18$. Thus we can assume that $s_1=6$; so $\deg(F_1)=1$. If $\eta(2)\ge 3$, then $s_0=s_1\eta(2)\ge 18$. Thus we can assume that $\eta(2)=2$. Hence $h_{\eta(2)}(x)=x^3$. As we have an identity $K[x]=K[x^4+(x^5-1)F_1(x),x^3]$ with $\deg(F_1)=1$, we reached a contradiction. So part (4.b) holds.

For part (4.c), so $m=6$, we know that $n\ge 20$ by part (4.a) and we show that the assumption that $n\in\{20,21,22,23\}$ leads to a contradiction. 

As $5\le s_1\mid s_0=n$, it follows that $(s_0,s_1)\in\{(20,5),(20,10),(21,7),(22,11)\}$. If $s_1$ is a prime, i.e., $s_1\in\{5,7,11\}$, then $s_2=1$ and hence $\eta(2)\ge 6-s_2=5$, and it follows that $n=s_0=s_1\eta(2)\ge 5s_1\ge 25$, a contradiction. Hence we can assume that $(s_0,s_1)=(20,10)$. Thus $\eta(2)=2$ and hence $s_2\ge m-\eta(2)=4$. As $s_2$ divides $s_1$, we get that $s_2=5$ and therefore $\eta(3)=4$ and $s_3=1$, so $\eta(3)+s_3=5<6=m$, a contradiction. So part (4.c) holds.

Based on part (4.c), to prove part (4.d) it suffices to show that the assumption that $s_0=24$ leads to a contradiction. As $m-1=5\le s_1\mid s_0$, we have $s_1\in\{6,8,12\}$. 

We assume that $s_1=6$. Thus $F_1=\gamma_1\in K^{\ast}$ and $\eta(2)=4$. So $s_2\ge 6-\eta(2)=2$. Hence the $K$-algebra $K[x]=K[h_{\eta(1)},h_{\eta(2)}]=K[x^5+\gamma_1(x^6-1),x^2g_{\eta(2)}(x)]$ is equal to $K[x^5(\gamma_1x-1),x^2g_{\eta(2)}(x)]$, a contradiction.

We assume that $s_1=8$. Thus $\eta(2)=3$ and hence $s_2\ge 6-\eta(2)=3$. As $s_2$ divides $s_1$, it follows that $s_2=4$. Thus $h_{\eta(2)}(x)=x^3(1+\gamma_3x)$ for some $\gamma_1\in K^{\ast}$. There exists $(\gamma_0,\gamma_1,\gamma_2)\in K^2\times K^{\ast}$ such that $F_1(x)=\gamma_2x^2+\gamma_1x+\gamma_0$. As we have $f(x)=x^5+(x^6-1)(\gamma_2x^2+\gamma_1x+\gamma_0)$, for all $(\delta_0,\delta_1,\delta_2)\in K^2\times K^{\ast}$, the coefficient of $x^5$ in $f-\delta_2h_{\eta(2)}^2-\delta_1h_{\eta(2)}-\delta_0$ is non-zero and this contradicts the fact that the pair $\bigl(f,h_{\eta(2)}\bigr)$ is obtained from a pair of polynomials in $K[x]$ of degree at most $1$ using the replacements of type (ii) and (iii) used in the proof of Theorem \ref{T14.5}.

We assume that $s_1=12$. Thus $\eta(2)=2$. Hence $s_2\ge 6-\eta(2)=4$. As $s_2$ divides $s_1$, it follows that $s_2\in\{4,6\}$. So there exists a pair $(\delta_{2,1},\gamma_7)\in K^2$ such that $h_{\eta(2)}(x)=x^4+\delta_{2,1}x^5+\gamma_7(x^6-1)$ and we have $(\delta_{2,1},\gamma_7)=(0,0)$ if $s_2=4$ and $(\delta_{2,1},\gamma_7)\in K\times K^{\ast}$ if $s_2=6$. There exists $(\gamma_0,\ldots,\gamma_6)\in K^6\times K^{\ast}$ such that $F_1=\sum_{j=0}^6 \gamma_jx^j$. If $s_2=4$ (resp.\ $s_2=6)$, for all $(\delta_0,\delta_1,\delta_2,\delta_3)\in K^3\times K^{\ast}$ (resp.\ $(\delta_0,\delta_1,\delta_2)\in K^2\times K^{\ast}$), the coefficients of $x^7$ and $x$ in $f-\sum_{j=0}^{\frac{12}{s_4}} \delta_jh_{\eta(2)}^j$ are $\delta_1$ and $-\delta_1$ (respectively) and this contradicts the fact that the pair $(f,h_{\eta(2)})$ is obtained from a pair of polynomials in $K[x]$ of degree at most $1$ using the replacements of type (ii) and (iii) used in the proof of Theorem \ref{T14.5}. So part (4.d) holds.\end{proof}

Next example shows that Theorem \ref{T15.6}(4.a) cannot be refined in a way similar to Theorem \ref{T15.6}(4.b) and (4.c) for $m\ge 7$.

\begin{example}\normalfont\label{EX18.1}
Let $m=7$, $K$ a field that has a $7$-th primitive root of unity, and $h(x):=x^6$. For each $\delta\in K^{\ast}$ we consider polynomials $F_1(x):=\delta^2x+2\delta$, $f(x):=x^6+(x^7-1)F_1(x)$, $g(x):=-f(x)^4-\delta f(x)^3+\delta x^4+x^3$, and $F_3(x):=\delta x^4+x^3$ in $K[x]$. Clearly, $\deg(g)=4\deg(f)=32<36=(7-1)^2$. As 
$$g(x)\equiv -(x^{24}-x^3)-\delta(x^{18}-x^4)\equiv 0\pmod{(x^7-1)k[x]},$$
there exists $F_2(x)\in k[x]$ such that $g(x)=(x^7-1)F_2(x)$. We have 
$$K[f,g]=K[f,F_3]=K[f-F_3^2,F_3]=K[-\delta^2x-2\delta,F_3]=K[x].$$
\end{example}

Our lower bounds for the $\pi_{n,m}(K)$s when $|K|\ge m$ are grouped as follows. 

\begin{theorem}\label{T15}
Suppose that $(n,m)\in (\mathbb N^\ast\setminus\{1\})\times (\mathbb N^\ast\setminus\{1,2\})$. Let $K$ be a field with $|K|\ge m$. Then the following properties hold.

\medskip
{\bf (1)} If $|K|>m$, then $\pi_{n,m}(K)\ge m+1$.

\smallskip
{\bf (2)} If $m$ is odd and $|K|=m\ge 5$, then $\pi_{n,m}(K)\ge m+1$. 

\smallskip
{\bf (3)} If $4\mid |K|=m$, then $\pi_{n,m}(K)\ge m$.

\smallskip
{\bf (4)} Suppose that $|K|> m\ge 4$. Writing $m=3r+\epsilon$ with $(r,\epsilon)\in\mathbb N^{\ast}\times\{0,1,2\}$, let $\phi(r,\epsilon):=4r^2+\epsilon(5-\epsilon)r+\epsilon(\epsilon-1)$. Then $\pi_{2,m}(K)\ge\phi(r,\epsilon)\ge \frac{4m^2}{9}$.

\smallskip
{\bf (5)} If $m\ge 4$ and $K$ has a primitive $m$-th root of unity, then $\pi_{2,m}(K)\ge (m-1)^2$ for $m\in\{4,5,6\}$ and $\pi_{2,m}(K)\ge \lfloor\frac{m+2}{2}\rfloor(m-1)$ for $m\ge 7$.\end{theorem}

\begin{proof}
Let $(\underline{P},\underline{Q})\in\mathbb D_{n,m}(K)^2$. To prove parts (1) to (3) it suffices to show that we can choose $\underline{P}=(P_1,\ldots,P_m)$ and $\underline{Q}=(Q_1,\ldots,Q_m)$ such that for each $a\in\GA_n(K)$ with $a(\underline{P})=\underline{Q}$ we have $\ell(a)\ge m+1$ in parts (1) and (2) and $\ell(a)\ge m$ in part (3). 

If $|K|>m$ (resp.\ $|K|=m$) let $\bigl((\alpha_1,\ldots,\alpha_m),(\beta_1,\ldots,\beta_m)\bigr)\in\mathbb D_{1,m}(K)^2$ be such that there exists $f\in K[x]$ of degree $m-1$ (resp.\ $m-2$) with $f(\alpha_i)=\beta_i$ for each $i\in \llbracket1,m\rrbracket$ by Proposition \ref{PR28}(1) (resp.\ \ref{PR28}(3)). Let $F(x):=\prod_{i=1}^n (x-\alpha_i)\in K[x]$.

We choose $P_i:=(\alpha_i,0,\ldots,0)$ and $Q_i:=\bigl(\beta_i,0,\ldots,0\bigr)$ for $i\in \llbracket1,m\rrbracket$. Clearly, $(\underline{P},\underline{Q})\in\mathbb D_{n,m}(K)^2$. Let $a=\e(f_1,\ldots,f_n)\in\GA_n(K)$ be such that $a(\underline{P})=\underline{Q}$. There exists an $n$-tuple $(F_1,\ldots,F_n)\in K[x]^n$ such that the set 
$$\{f_1-f(x_1)-F_1(x_1)F(x_1)\}\cup\{f_i-F_i(x_1)F(x_1)|i\in \llbracket2,m\rrbracket\}$$ is contained in the ideal $(x_2,\ldots,x_n)$ of $R$. Thus
\begin{equation}\label{EQ48}
\ell(a)\ge\pi(a)\ge\pi\bigl(f(x_1)+F_1(x_1)F(x_1),F_2(x_1)F(x_1),\ldots,F_n(x_1)F(x_1)\bigr).
\end{equation}
As $a\in\GA_n(K)$, the $K$-algebra $K[x]$ is generated by $\{f+F_1F\}\cup\{F_iF|i\in \llbracket2,n\rrbracket\}$. 

As either $m\ge 4$ or $|K|>m=3$, $\deg(f)\ge 2$; so $f(x)$ does not generate $K[x]$. So there exists $i\in \llbracket1,m\rrbracket$ such that $F_i\neq 0$. Thus the right hand side of Inequalities (\ref{EQ48}) is at least $m=\deg(F)$. Hence part (3) holds. 

To prove parts (1) and (2), we show that the assumption that $F_i\in K$ for all $i\in \llbracket1,m\rrbracket$ leads to a contradiction. This assumption implies that the $K$-algebra $K[x]$ is generated by $f(x)$ and $F(x)$. Let $\varepsilon:=m-1-\deg(f)\in\{0,1\}$; so, if $\varepsilon=1$, then $|K|=m\ge 5$ is odd. As $\textup{g.c.d.}\bigl(\deg(f),\deg(F)\bigr)=\textup{g.c.d.}(m-1-\varepsilon,m)=1$, $\char(K)$ does not divide it and hence $m-1-\varepsilon\mid m$ by \cite{AM}, Main Thm.\ (1.1). It follows that $m-1-\varepsilon=1$, hence $m\le 3$, a contradiction. 

As there exists $i\in \llbracket1,m\rrbracket$ such that $F_i\notin K$, the right hand side of Inequalities (\ref{EQ48}) is at least $m+1$. So parts (1) and (2) hold.

For parts (4) and (5) we use the above notation for $n=2$. 

For part (4), we take the pair $\bigl((\alpha_1,\ldots,\alpha_m),(\beta_1,\ldots,\beta_m)\bigr)\in\mathbb D_{1,m}(K)^2$ to be $\llbracket1,m-1\rrbracket$-generic by Proposition \ref{PR28}(2).  Among all automorphisms $a$, we choose the one for which the pair 
$$(s_0,s_1):=\bigl(\deg(F_2F),\deg(f+F_1F)\bigr)\in\mathbb N^2$$ 
is minimal with respect to the total lexicographic order on $\mathbb N^2$. 

If $s_0\le s_1$, then $s_0\mid s_1$ by Theorem \ref{T14.5}(2) and given any $g\in xK[x]$, for $b:=\e\bigl(x_1+g(x_2),x_2\bigr)\in\SGA_2(K)$ we have $(ba)(\underline{P})=\underline{Q}$; by replacing $a$ with $ba$, $(f+F_1F,F_2F)$ gets replaced with $\bigl(f+F_1F+g\circ (F_2F),F_2F\bigr)$ and $g$ can be chosen so that $\deg\bigl(f+F_1F+g(F_2F)\bigr)<s_1$ by the proof of Theorem \ref{T15.6}. In the replacement of $a$ by $ba$ we might have $\ell(a)<\ell(ba)$ but we are only interested in the lower bounds of $\ell(a)$ and $\ell(ba)$ that are given by Inequalities (\ref{EQ48}) and clearly the one for $ba$ is less. Thus we can assume that $s_1<s_0$. In such a case we have $\deg(F_2F)\ge \phi(r,\epsilon)\ge \frac{4m^2}{9}$ by Theorem \ref{T15.6}(1) to (3) and hence part (4) holds.

For part (5), we choose the pair $\bigl((\alpha_1,\ldots,\alpha_m),(\beta_1,\ldots,\beta_m)\bigr)$ as in Theorem \ref{T15.6}(4). As above we argue that we can assume that $s_1<s_0$. In such a case we have $\deg(F_2F)\ge 16$ if $m=5$ and $\deg(F_2F)\ge \lfloor\frac{m+2}{2}\rfloor(m-1)$ by Theorem \ref{T15.6}(4.a) and hence part (5) holds.
\end{proof}

\section{Lower bounds for $n=2$ via a count of permutations}\label{S27}

For a finite field $K$ and a pair $(n,N)\in (\mathbb N^{\ast}\setminus\{1\})^2$, one would like to compute the cardinalities
$$|\{a\in\TGA_n(K)|\ell(a)=N\}|=|\TGA_n(K)[N]|-|\TGA_n(K)[N-1]|,$$ 
$$|\{a\in\STGA_n(K)|\ell(a)=N\}|=|\STGA_n(K)[N]|-|\STGA_n(K)[N-1]|,$$
$$|\{\sigma\in\Perm(K^n)|\ell_{n,K}(\sigma)=N\}|,\;\;\;\text{and}\;\;\;|\{\sigma\in\Alt(K^n)|\ell^{\S}_{n,K}(\sigma)=N\}|.$$
For $n\ge 3$ we do not know results on these cardinalities. But for $n=2$, based on van der Kulk's results (\cite{vdK}, Thms.\ 1 and 2) the generating series for the numbers 
$$|\{a\in\GA_2(K)|\ell(a)=N\}|=(|K|-1)|\{a\in\SGA_2(K)|\ell(a)=N\}|$$
with $N\ge 2$ and $|\AGL_2(K)|=(|K|-1)|\ASL_2(K)|$ is known based on \cite{DV}, Thms.\ 2.1 and one gets bounds of the following type (see \cite{DV}, Thm.\ 2.3)
$$1\le\frac{|\{a\in\GA_2(K)|\ell(a)=N\}|}{(|K|-1)^3(|K|+1)^2|K|^{N+1}}\le 1+\frac{\log_2(N)^{\log_2(N)}}{(|K|-1)^3(|K|+1)^2|K|^{\frac{N}{2}-7}}.$$

For applications to the $\pi_{n,m}(K)$s and the $\pi^{\S}_{n,m}(K)$s it is more convenient to compute or estimate the following cardinalities involving permutations.

\phantomsection{Let $\varrho_{n,K,N}\in\mathbb N^{\ast}$ (resp.\ $\varrho^{\S}_{n,K,N}\in\mathbb N^{\ast}$) be the cardinality of the set}\label{PH89} 
$$\varrho_{n,K}\bigl(\TGA_n(K)[N]\bigr)=\{\sigma\in\Perm(K^n)|\ell_{n,K}(\sigma)\le N\}$$
(resp.\ $\varrho_{n,K}\bigl(\STGA_n(K)[N]\bigr)=\{\sigma\in\Alt(K^n)|\ell^{\S}_{n,K}(\sigma)\le N\}$).
Clearly, $$\varrho_{n,K,N}=|\AGL_n(K)|+\sum_{i=2}^N |\{\sigma\in\Perm(K^n)|\ell_{n,K}(\sigma)=i\}|$$
and $\varrho^{\S}_{n,K,N}=|\ASL_n(K)|+\sum_{i=2}^N |\{\sigma\in\Alt(K^n)|\ell^{\S}_{n,K}(\sigma)=i\}|$. 

In this section, we use van der Kulk's results as in Section \ref{S10} to get upper bounds for $\varrho_{2,K,N}$ and apply these bounds to obtain a lower bound for $\pi_{2,|K|^2}(K)$.
 
We begin by recalling some values related to the zeta function and properties of an arithmetic function $\rho:\mathbb N^{\ast}\rightarrow\mathbb N^{\ast}$. 

\phantomsection{The zeta function $\zeta:(1,\infty)\rightarrow (1,\infty)$ defined by the rule $\zeta(s):=\sum_{i=1}^{\infty} \frac{1}{i^s}$ is a decreasing bijection. So there exists a unique $\epsilon_{\z}\in (1,\infty)$ such that $\zeta(\epsilon_{\z})=2$. We have $\epsilon_{\z}=1.72864...<1.73$ and $\C_2:=\frac{-1}{\epsilon_{\z}\zeta^{\prime}(\epsilon_{\z})}=0.31817...<0.32$. Let}\label{PH85} 
$$\C_1:=2^{1-\epsilon_{\z}}=\frac{2}{2^{\epsilon_{\z}}}=0.60346...<0.61 .$$

\begin{lemma}\label{L24} For $N\in \mathbb N^{\ast}$ let $\rho(N)$ be the number of finite sequences $(d_i)_{i\in \llbracket1,j\rrbracket}$ in $\mathbb N^{\ast}\setminus\{1\}$ with $\prod_{i=1}^j d_i\le N$.\footnote{The empty sequence with product $1$ is included and hence we have $\rho(1)=1$. The more traditional function on $N\in\mathbb N^{\ast}$ in the literature is $\rho(N)-1.$} Then the following properties hold.

\medskip
{\bf (1)} If $N\ge 2$, then $\rho(N)\le\C_1N^{\epsilon_{\z}}$ and the equality only holds for $N=2$.

\smallskip
{\bf (2)} If $C\in (0,\C_2)$, then $\rho(N)>CN^{\epsilon_{\z}}$ for $N\gg 1$.
\end{lemma}

\begin{proof}
By considering the possible values for $d_1$, so $d_1\in \llbracket1,N\rrbracket$, and the resulting inequality $\prod_{i=2}^j d_i\le \lfloor\frac{N}{d_1}\rfloor$, we get the recursive inequality 
\begin{equation}\label{EQ49}
\rho(N)\le 1+\sum_{i=2}^N \rho\Bigl(\Bigl\lfloor\frac{N}{i}\Bigr\rfloor\Bigr)
\end{equation}
(the term $1$ corresponds to the empty sequence).

Clearly, $\rho(2)=2=\C_12^{\epsilon_{\z}}$. Let $\mathfrak S(N)$ be the statement that $\rho(N)< N^{\epsilon_{\z}}$. We prove by induction on $N\ge 3$ that $\mathfrak S(N)$ is true. As $\rho(3)=3$, $\rho(4)=5$, $\rho(5)=6$, $\rho(6)=9$, $\rho(7)=10$, and $\rho(8)=14,$ it is easy to see that $\mathfrak S(N)$ holds for $N\in \llbracket2,8\rrbracket$; so the base of the induction holds. Assume that $N\ge 9$ and $\mathfrak S(i)$ is true for $i\in \llbracket2,N-1\rrbracket$. We have $1+\rho\bigl(\bigl\lfloor\frac{N}{\lfloor\frac{N}{3}\rfloor}\bigr\rfloor\bigr)=1+\rho(3)=4<0.6\cdot 3^{1.728}<\C_1\bigl(\frac{N}{\lfloor\frac{N}{3}\rfloor}\bigr)^{\epsilon_{\z}}$ as $\frac{N}{\lfloor\frac{N}{3}\rfloor}\in [3,4)$. Thus, starting by rearranging Inequality (\ref{EQ49}), we estimate
$$\rho(N)\le 1+\rho\Bigl(\Bigl\lfloor\frac{N}{\lfloor\frac{N}{3}\rfloor}\Bigr\rfloor\Bigr)+\sum_{i\in \llbracket2,N\rrbracket\setminus\{\lfloor\frac{N}{3}\rfloor\}} \rho\bigl(\bigl\lfloor\frac{N}{i}\bigr\rfloor\bigr)$$
$$<\C_1\Bigl(\frac{N}{\lfloor\frac{N}{3}\rfloor}\Bigr)^{\epsilon_{\z}}+\C_1\sum_{i\in \llbracket2,N\rrbracket\setminus\{\lfloor\frac{N}{3}\rfloor\}}\left(\frac{N}{i}\right)^{\epsilon_{\z}}=\C_1\sum_{i=2}^N \left(\frac{N}{i}\right)^{\epsilon_{\z}}$$
$$<\C_1N^{\epsilon_{\z}}[\zeta(\epsilon_{\z})-1]=\C_1N^{\epsilon_{\z}},$$
and hence $\mathfrak S(N)$ is true. This ends the induction. So part (1) holds. 

Part (2) follows directly from an asymptotic result of Kalm\'ar (see \cite{Kalm}, Eq.\ (27) or \cite{Hw}, Cor.\ 2).\end{proof}

\begin{lemma}\label{L25}
Let $(j,k,s)\in (\mathbb N^{\ast})^2\times\mathbb N$ with $k\ge 4$. We consider the set
$$B=B_{j,k,s}:=\Bigl\{\prod_{i=1}^j d_i|(d_1,\ldots,d_j)\in \llbracket2,k-1\rrbracket^j,\,\sum_{i=1}^j d_i=2j+s\Bigr\}$$
and we write $s=q(k-3)+r$ with $(q,r)\in\mathbb N\times\llbracket0,k-4\rrbracket$.
Then the following properties hold.

\medskip
{\bf (1)} We have $\min(B)=(r+2)2^{j-q-1}(k-1)^{q}=(r+2)2^{j-1}\bigl(\frac{k-1}{2}\bigr)^q$.

\smallskip
{\bf (2)} Then we have 
$$\frac{s+2j}{\ln\bigl(\min(B)\bigr)}=\frac{q(k-3)+r+2j}{q\ln(\frac{k-1}{2})+\ln(r+2)+(j-1)\ln(2)}\le\begin{cases} \frac{k-1}{\ln(k-1)}\quad\quad\;\, {\rm if}\;k\ge 5\\
\frac{2}{\ln{2}}\quad\quad\quad\quad {\rm if}\;k=4.\end{cases}$$ 
\end{lemma}

\begin{proof}
For part (1), let $(d_1,\ldots,d_j)\in \llbracket2,k-1\rrbracket^j$ be such that $\sum_{i=1}^j d_i=s+2j$ and $\prod_{i=1}^j d_i=\min(B)$. For $l\in\llbracket2,k-1\rrbracket$ let $\c_l:=|\{i\in\llbracket1,j\rrbracket|d_i=l\}|$. Therefore $\min(B)=\prod_{l=2}^{k-1} l^{\c_l}$. If $\sum_{l=3}^{k-2} \c_l\ge 2$, by choosing $(i_1,i_2)\in\llbracket1,j\rrbracket^2$ with $i_1\neq i_2$ and $3\le d_{i_1}\le d_{i_2}\le k-2$ and by replacing $(d_{i_1},d_{i_2})$ with $(d_{i_1}-1,d_{i_2}+1)$, as $(d_{i_1}-1)(d_{i_2}+1)=d_{i_1}d_{i_2}-d_{i_2}+d_{i_1}-1\le d_{i_1}d_{i_2}-1$, we reach a contradiction to $\prod_{i=1}^j d_i=\min(B)$. Thus $\sum_{l=3}^{k-2} \c_l\le 1$. We consider two cases as follows.

{\bf Case 1: $\c_3=\cdots=\c_{k-2}=0$.} Thus $j=\c_2+\c_{k-1}$, $s+2j=2\c_2+(k-1)\c_{k-1}$, and $\min(B)=2^{\c_2}(k-1)^{\c_{k-1}}$. Hence $\c_{k-1}=\frac{s}{k-3}$ and $\c_2=j-\frac{s}{k-3}$. This case is possible iff $r=0$. We have $s=q(k-3)$, $\c_{k-1}=q$, and $\min(B)=2^{j-q}(k-1)^{q}$. Hence $\min(B)=(r+2)2^{j-q-1}(k-1)^{q}$ as $r=0$.

{\bf Case 2: $k\ge 5$ and $\sum_{l=3}^{k-2} \c_l=1$.} Let $l\in\llbracket3,k-2\rrbracket$ be the unique integer such that $\c_l=1$. We have $j=\c_2+1+\c_{k-1}$, $s+2j=2\c_2+l+(k-1)\c_{k-1}$, and $\min(B)=l2^{\c_2}(k-1)^{\c_{k-1}}$. Thus $\c_{k-1}=j-1-\c_2$ and $$l=s+2j-2\c_2-(k-1)(j-1-\c_2)=s+k-1+(\c_2-j)(k-3).$$
As $3\le l\le k-2$, we get that 
$$s+1=q(k-3)+r+1\le (j-\c_2)(k-3)\le s+k-4=q(k-3)+k-4+r.$$ 
Hence $j-\c_2=q+1$. Thus $\c_2=j-q-1$, $l=r+2$, and $\c_{k-1}=q$. 
Therefore this case is possible iff $r\in\llbracket1,k-4\rrbracket$. We have $\min(B)=(r+2)2^{j-q-1}(k-1)^{q}$.

Part (1) follows from the two cases above.

For part (2), we first note that for $r=0$ we are in Case 1 and $j-q=\c_2\ge 0$ and for $r\in\llbracket1,k-4\rrbracket$ we are in Case 2 and $j-q=\c_1+1\ge 1$. We consider the function $\mathfrak f:[\max(q,1),\infty)\rightarrow\mathbb R$ defined by the rule 
$$\mathfrak f(y):=\frac{q(k-3)+r+2y}{q\ln(\frac{k-1}{2})+\ln(r+2)+(y-1)\ln(2)}.$$ 
We have $\frac{s+2j}{\ln\bigl(\min(B)\bigr)}=\mathfrak f(j)$ and 
\begin{equation*}
\begin{split}
\mathfrak f'(y)&=\frac{q[2\ln\bigl(\frac{k-1}{2}\bigr)-(k-3)\ln 2]+[2\ln\frac{r+2}{2}-r\ln(2)]}{[q\ln(\frac{k-1}{2})+\ln(r+2)+(y-1)\ln(2)]^2}\\
&=\frac{2q(k-1)[2\frac{\ln(k-1)}{k-1}-\frac{\ln 2}{2}]+2(r+2)[\frac{\ln(r+2)}{r+2}-\frac{\ln 2}{2}]}{[q\ln(\frac{k-1}{2})+\ln(r+2)+(y-1)\ln(2)]^2}.\end{split}
\end{equation*}

The function $\mathfrak f_1:[3,\infty)\rightarrow\mathbb R$ defined by $\mathfrak f_1(z):=\frac{z-1}{\ln(z-1)}$ is increasing for $z\ge 5$ with $\mathfrak f_1(3)=\mathfrak f_1(5)>\mathfrak f_1(4)$. Thus for $k\neq 4$ we have $\frac{k-1}{\ln(k-1)}\ge\frac{2}{\ln 2}$. If $k=4$, then $r=0$ and $f$ is constant for $q=0$ and $f$ is strictly decreasing for $q>0$. So, if $\mathfrak f$ is non-decreasing, then $\mathfrak f(j)\le\lim \limits_{y\rightarrow\infty} \mathfrak f(y)=\frac{2}{\ln{2}}$ and part (2) holds.

So we can assume that $\mathfrak f$ is strictly decreasing. For $r=q=0$ we have $\mathfrak f(j)=\frac{2}{\ln{2}}$. For $r=0$ and $q\ge 1$ we have $\mathfrak f(j)\le \mathfrak f(q)=\frac{k-1}{\ln(k-1)}$. If $r\ge 1$, then we have 
$$\mathfrak f(j)\le \mathfrak f(q+1)=\frac{q(k-1)+r+2}{q\ln(k-1)+\ln(r+2)}\le\frac{k-1}{\ln(k-1)},$$ as $\mathfrak f_1(r+3)\le \mathfrak f_1(k)$. So part (2) holds.\end{proof}

\phantomsection{For $k\in\mathbb N^{\ast}\setminus\{1,2\}$ we define}\label{PH101}
$$\mathfrak c_k:=\begin{cases}-1+\epsilon_{\z}+\frac{(k-1)\ln(1+\frac{1}{k-1})}{\ln(k-1)}\in (0.72864,1.29507)\quad\;\;\,\,\, {\rm if}\;k\ge 6\\
-0.61+\epsilon_{\z}+\frac{4\ln(1+\frac{1}{4})}{\ln 4}=1.7624...\quad\quad\quad\quad\quad\quad\quad\,\; {\rm if}\;k=5\\
\epsilon_{\z}=1.72864...\quad\quad\quad\quad\quad\quad\quad\quad\quad\quad\quad\quad\quad\quad\quad\;\,\, {\rm if}\;k=4\\
\frac{2\ln 3}{\ln 2}-3+\epsilon_{\z}=1.89856...\;\quad\quad\quad\quad\quad\quad\quad\quad\quad\quad\;\, {\rm if}\;k=3.\end{cases}$$

\begin{proposition}\label{PR30}
Let $K$ be a finite field with $|K|\ge 3$ and $N\in\mathbb N^{\ast}\setminus\{1\}$. Then the following properties hold.

\medskip
{\bf (1)} We have inequalities
$$\varrho_{2,K,N}\le (|K|^2-1)^2N^{|K|+\mathfrak c_{|K|}-\epsilon_{\z}}\rho(N)\le (|K|^2-1)^2N^{|K|+\mathfrak c_{|K|}}.$$

{\bf (2)} We have inequalities
$$\varrho^{\S}_{2,K,N}\le(|K|^2-1)(|K|+1)N^{|K|+\mathfrak c_{|K|}-\epsilon_{\z}}\rho(N)\le (|K|^2-1)(|K|+1)N^{|K|+\mathfrak c_{|K|}}.$$
\end{proposition}

\begin{proof}
For $\sigma\in\varrho_{n,K}\bigl(\GA_n(K)[N]\bigr)$, let $a\in\GA_2(K)[N]$ be such that $a(K)$ is a minimal representation of $\sigma$ by Proposition \ref{PR14}(1). Let $j=\j_a\in\mathbb N$ be the Furter's length of $a$. We consider a reduced standard product decomposition $a=a_1\cdots a_jb$ (see Definition \ref{D13}(2) and (3)). For $i\in \llbracket1,j\rrbracket$ let $d_i:=\ell(a_i)$. Then we have $\ell(a)=\prod_{i=1}^j d_i\le N$ by Equation (\ref{EQ14}) and $d_i\in \llbracket2,|K|-1\rrbracket$ for each $i\in \llbracket1,j\rrbracket$ by Definitions \ref{D13}(1) and (3) and \ref{D15}(2). 

Each such particular type of a reduced standard product decomposition leads to a family of automorphisms $a$. We bound the total number of permutations $a(K)$ that can be obtained in this way, i.e., we bound $\bigl|\varrho_{2,K}\bigl(\GA_2(K)[N]\bigr)\bigr|$ from above following in essence the proof of \cite{DV}, Thm.\ 2.1 for the next four paragraphs.\footnote{For $N\ge |K|$, part (1) is a stronger result than what one would get by applying \cite{DV}, Thm.\ 2.3 directly.}

We fix an arbitrary sequence $(d_i)_{i\in \llbracket1,j\rrbracket}$ in $\llbracket2,|K|-1\rrbracket$. 

The total number of possible $b\in\AGL_2(K)$ is $|K|^2(|K|^2-1)(|K|^2-|K|)$. For $j\ge 1$ the total number of possible $a_jb$ is at most the product of the total number of possible $a_j$ and the number $(|K|^2-|K|)(|K|^2-1)$, i.e., the number of representatives of a quotient of $\AGL_2(K)$ by a subgroup isomorphic to $K^2$, by Remark \ref{R3}(3).

The total number of each reduced non-affine triangular automorphism $a_i$ of degree length $d_i$ is $(|K|+1)(|K|-1)|K|^{d_i}$, where $|K|+1$ is the total number of directions\footnote{Equivalently, it is the index $[\GL_2(K):B_2(K)]$ by (and with the notation of) Remark \ref{R3}(2).}, and we have $d_i-1$ coefficients of a polynomial $f_i\in x^2K[x]$ that defines $a_i$, the leading one being non-zero. Thus the total number of possible permutations $a_i(K)$ given by such $a_i$s is bounded from above for $i=1$ by $(|K|+1)(|K|-1)|K|^{d_1-2}$ and for each $i\in \llbracket2,j\rrbracket$, as we have $\dir(a_i)\neq\dir(a_{i-1})$, by $|K|(|K|-1)|K|^{d_i-2}$.

We get that the total number of permutations $a(K)$ with $a$ in the family of minimal representations that correspond to $(d_i)_{i\in \llbracket1,j\rrbracket}$ (i.e., with $\ell(a_i)=d_i$ for each $i\in \llbracket1,j\rrbracket$) is bounded from above for $j=0$ by 
$$B_{\emptyset}:=|K|^2(|K|^2-1)(|K|^2-|K|)$$ 
and for $j\ge 1$ by $(|K|^2-1)(|K|^2-|K|)(|K|+1)|K|^{j-1}\prod_{i=1}^j
(|K|-1)|K|^{d_i-2}$ and thus by
\begin{equation}\label{EQ37.9} B_{(d_i)_{i\in \llbracket1,j\rrbracket}}:=(|K|^2-1)^2(1-|K|^{-1})^j|K|^{\sum_{i=1}^j d_i}.
\end{equation}

Let $s:=-2j+\sum_{i=1}^j d_i$. We have $s\in\mathbb N$ and $s+2j=\sum_{i=1}^j d_i$. 

We write $s=q(|K|-3)+r$ with $(q,r)\in\mathbb N\times\llbracket0,|K|-4\rrbracket$. From Lemma \ref{L25}(1) we get that $(r+2)2^{j-1}\bigl(\frac{|K|-1}{2}\bigr)^{q}\le N$ if $|K|\ge 4$. 
So for $|K|=4$ we have
\begin{equation}\label{EQ50}
\frac{\ln |K|^{s+2j}}{\ln N}\le \ln(|K|)\frac{q(|K|-3)+r+2j}{q\ln(\frac{|K|-1}{2})+\ln(r+2)+(j-1)\ln(2)}\le \frac{2\ln |K|}{\ln 2}=4
\end{equation}
and for $|K|\ge 5$ we have 
\begin{equation}\label{EQ51}
\frac{\ln |K|^{s+2j}}{\ln N}\le \ln(|K|)\frac{q(|K|-3)+r+2j}{q\ln(\frac{|K|-1}{2})+\ln(r+2)+(j-1)\ln(2)}\le \frac{(|K|-1)\ln |K|}{\ln(|K|-1)}\end{equation}
by Lemma \ref{L25}(2). 
If $|K|=3$, then $d_1=\cdots=d_j=2$, $s=0$, $N\ge 2^j$, and $\frac{\ln |K|^{2j}}{\ln (N)}\le \frac{2\ln(3)}{\ln(2)}.$

From the prior paragraph and the definition of $\mathfrak c_{|K|}$ we get that for all $K$ we have 
$\frac{\ln|K|^{\sum_{i=1}^j d_i}}{\ln (N)}\le |K|+\mathfrak c_{|K|}-\epsilon_{\z}$. Hence 
$$|K|^{\sum_{i=1}^j d_i}\le N^{|K|+\mathfrak c_{|K|}-\epsilon_{\z}}.$$
As for all $K$ we have
$$B_{\emptyset}=\frac{|K|^3}{|K|+1}(|K|^2-1)^2<(|K|^2-1)(|K|^2-1)^2\le 2^{|K|}(|K|^2-1)^2\le N^{|K|}(|K|^2-1)^2,$$ 
for $|K|=5$ the inequality $\frac{|K|^3}{|K|+1}\le 2^{|K|-0.61}$ gives $B_{\emptyset}\le N^{|K|-0.61}(|K|^2-1)^2$, and for $|K|\ge 7$ the inequality $(|K|^2-1)<2^{|K|-1}$ gives $B_{\emptyset}\le N^{|K|-1}(|K|^2-1)^2$, we conclude that the total number of permutations $a(K)$ that can be given by minimal representations $a\in\GA_2(K)[N]$ is strictly less than 
$$(|K|^2-1)^2N^{|K|+\mathfrak c_{|K|}-\epsilon_{\z}}\rho_{|K|-1}(N)\le (|K|^2-1)^2N^{|K|+\mathfrak c_{|K|}-\epsilon_{\z}}\rho(N),$$ 
where $\rho_{|K|-1}(N)\in \llbracket1,\rho(N)\rrbracket$ is the number of finite sequences $(d_i)_{i\in \llbracket1,j\rrbracket}$ in the set $\llbracket2,|K|-1\rrbracket$ with $\prod_{i=1}^j d_i\le N$. As $\rho(N)\le N^{\epsilon_{\z}}$ by Lemma \ref{L24}(1), part (1) holds.

Part (2) is proved similarly as we have $a=a_1\cdots a_jb\in\SGA_2(K)$ iff $b$ has Jacobian determinant $1$. So the only change required for proving part (2) relates to the fact that $|\ASL_2(K)|=|K|^3(|K|^2-1)$.\end{proof}

\begin{theorem}\label{T16} Let $K$ be a finite field. Then the following properties hold.

\medskip
{\bf (1)} We have a strict inequality
\begin{equation}\label{EQ52}
\ln\bigl(\pi_{2,|K|^2}(K)\bigr) > \frac{|K|^2[\ln (|K|^2-2)-1]+3-\ln [(1-|K|^{-2})(|K|^2-2)^2]}{|K|+\mathfrak c_{|K|}}.
\end{equation}

{\bf (2)} We have a strict inequality
$$\ln\bigl(\pi^{\S}_{2,|K|^2-2}(K)\bigr) > \frac{|K|^2[\ln (|K|^2-2)-1]+3-\ln [(2|K|^{-1}+2|K|^{-2})(|K|^2-2)^2]}{|K|+\mathfrak c_{|K|}}.$$
\end{theorem}

\begin{proof}
For part (1), for $|K|=2$ our lower bound is negative and for $|K|=3$ it is less than $1$. If $4\mid |K|$, then $\pi_{2,|K|^2}(K)=\infty$. So we can assume that $|K|$ is odd and at least $5$. 

Let $N:=\pi_{2,|K|^2}(K)$; so $N\ge 2$ and $\varrho_{n,K}$ restricts to a surjective function $\GA_2(K)[N]\rightarrow\perm(K^2)$ defined by the rule $a\mapsto a(K)$. Thus $(|K|^2)!\le\varrho_{2,K,N}$. From this and Proposition \ref{PR30}(1) we get by transitivity the following inequality $(|K|^2)!<(|K|^2-1)^2 N^{|K|+\mathfrak c_{|K|}}$. Hence
\begin{equation}\label{EQ53}
(|K|^2-2)!< (1-|K|^{-2}) N^{|K|+\mathfrak c_{|K|}}.
\end{equation}

Recall that for $m\in\mathbb N^{\ast}$ we have
$$\ln (m!)=\sum_{i=1}^m \ln i > \int \limits_1^m \ln t \ dt = m\ln m-m+1.$$
So for $m=|K|^2-2$, by taking $\ln$ of Inequality (\ref{EQ53}) we get the inequality
$$\ln(1-|K|^{-2})+(\ln N)(|K|+\mathfrak c_{|K|}) > (|K|^2-2) \bigl(\ln (|K|^2-2)\bigr) -|K|^2 +3$$
$$=|K|^2[\ln (|K|^2-2)-1]+3-2\ln (|K|^2-2),$$
which implies that Inequality (\ref{EQ52}) holds.

Part (2) is proved similarly based on Proposition \ref{PR30}(2), the only change being that the expression $1-|K|^{-2}$ gets replaced by $2|K|^{-1}+2|K|^{-2}$, the factor $2$ coming from the replacement of $\perm(K^2)$ by $\Alt(K^2)$.
\end{proof}

\begin{theorem}\label{T17} Let $K$ be a finite field with $|K|\geq 3$. Let $m\in \llbracket3,|K|^2-2\rrbracket$. Then the following properties hold.

\medskip
{\bf (1)} We have a strict inequality
\begin{equation}\label{EQ54}
\ln \bigl(\pi_{2,m}(K) \bigr) >\frac{(m-2)\ln\left(\frac{|K|^2-m}{e}\right)-(|K|^2-2)\ln\left(\frac{K|^2-m}{|K|^2-2}\right)-\ln(1-|K|^{-2})}{|K|+\mathfrak c_{|K|}}.
\end{equation}
In particular, by taking $m=|K|^2-2$ we have
\begin{equation}\label{EQ55}
\ln\bigl(\pi_{2,|K|^2-2}(K)\bigr) > \frac{(|K|^2-2)\ln\left(\frac{|K|^2-2}{e}\right)+2[1-\ln(2)]-\ln(1-|K|^{-2})}{|K|+\mathfrak c_{|K|}}.
\end{equation}

{\bf (2)} We have a strict inequality
$$\ln \bigl(\pi^{\S}_{2,m}(K) \bigr) >\frac{(m-2)\ln\left(\frac{|K|^2-m}{e}\right)-(|K|^2-2)\ln\left(\frac{K|^2-m}{|K|^2-2}\right)-\ln(|K|^{-1}+|K|^{-2})}{|K|+\mathfrak c_{|K|}}.$$
In particular, by taking $m=|K|^2-2$ we have
$$\ln\bigl(\pi^{\S}_{2,|K|^2-2}(K)\bigr) > \frac{(|K|^2-2)\ln\left(\frac{|K|^2-2}{e}\right)+2[1-\ln(2)]-\ln(|K|^{-1}+|K|^{-2})}{|K|+\mathfrak c_{|K|}}.$$
\end{theorem}

\begin{proof} We first prove part (1). Let $N:=\pi_{2,m}(K)$; we have $N\ge 2$. The total number of injective functions $Y\rightarrow K^2$, where $Y\subset K^2$ with $|Y|=m$, is $m!{{|K|^2}\choose{m}}^2$. The number of permutations $a(K)$ with $a\in\GA_2(K)[N]$ is less than $(|K|^2-1)^2N^{|K|+\mathfrak c_{|K|}}$ by Proposition \ref{PR30}(1). Each such permutation $a(K)$ gives ${|K|^2}\choose{m}$ functions. So, as the $a(K)$s give all possible injective functions $Y\rightarrow K^2$, we have an inequality
$$(|K|^2-1)^2 N^{|K|+\mathfrak c_{|K|}}>\frac{m!{{|K|^2}\choose{m}}^2}{{{|K|^2}\choose{m}}}=\frac{(|K|^2)!}{(|K|^2-m)!}= \prod_{i=|K|^2-m+1}^{|K|^2} i.$$
Thus
\begin{equation}\label{EQ56}
(1-|K|^{-2}) N^{|K|+\mathfrak c_{|K|}}>\prod_{i=|K|^2-m+1}^{|K|^2-2} i.
\end{equation}

We have
$$\ln\Bigl(\prod_{i=|K|^2-m+1}^{|K|^2-2} i\Bigr)=\sum_{i=|K|^2-m+1}^{|K|^2-2} \ln i \ \ > \int \limits_{|K|^2-m}^{|K|^2-2} \ln t \ dt =x(\ln x -1) \biggl|_{|K|^2-m}^{|K|^2-2}. $$
This evaluates and simplifies to $(m-2)[-1+\ln (|K|^2-m)]-(|K|^2-2)\ln\left(\frac{|K|^2-m}{|K|^2-2}\right).$ 
We get that 
$$(|K|+\mathfrak c_{|K|})\ln N +\ln(1-|K|^{-2})>(m-2)\ln\left(\frac{|K|^2-m}{e}\right)-(|K|^2-2)\ln\left(\frac{|K|^2-m}{|K|^2-2}\right)$$
from which Inequality (\ref{EQ54}) follows.

Part (2) is proved similarly to part (1) based on Proposition \ref{PR30}(2), the only change being that the expression $1-|K|^{-2}$ gets replaced by $|K|^{-1}+|K|^{-2}$.
\end{proof}

We have the following estimates of Furter's lengths of finite fields.

\begin{corollary}\label{C26} Let $K$ be a finite field with $|K|\ge 3$. If $4\nmid |K|$ (resp.\ $4\mid |K|$), then
$$\j_K\ge \frac{|K|^2\ln\left(\frac{|K|^2-2}{e}\right)+3-\ln [(1-|K|^{-2})(|K|^2-2)^2]}{(|K|+\mathfrak c_{|K|})\ln(|K|-1)}$$
(resp.\ $\j_K\ge \frac{(|K|^2-2)\ln\bigl(\frac{|K|^2-2}{e}\bigr)+2[1-\ln(2)]-\ln(2|K|^{-1}+2|K|^{-2})}{(|K|+\mathfrak c_{|K|})\ln(|K|-1)}$).
\end{corollary}

\begin{proof}
This follows from Theorem \ref{T16}(1) (resp.\ \ref{T17}(2)) and Corollary \ref{C12}.
\end{proof}

We have the following asymptotic ranges for $n=2$ and finite fields.

\begin{corollary}\label{C27}
For $\varepsilon\in (0,2)$ the following properties hold.

\medskip
{\bf (1)} There exists $C(\varepsilon)\in\mathbb N^{\ast}$ such that for each finite filed $K$ with $|K|\ge C(\varepsilon)$ and $|K|$ odd (resp.\ and $|K|$ an odd prime or a power of $3$) we have
$$(2-\varepsilon) |K|\ln |K|<\ln\bigl(\pi_{2,|K|^2}(K)\bigr)<18|K|\ln |K|$$
(resp.\ $(2-\varepsilon) |K|\ln |K|<\ln\bigl(\pi_{2,|K|^2}(K)\bigr)<(13.5+\epsilon)|K|\ln |K|$).

\smallskip
{\bf (2)} Let $p\ge 5$ be a prime. There exists $C_p(\varepsilon)\in\mathbb N^{\ast}$ such that for each finite field $K$ with $|K|\ge C_p(\varepsilon)$ and $\char(K)=p$ we have
$$(2-\varepsilon) |K|\ln |K|<\ln\bigl(\pi_{2,|K|^2-2}(K)\bigr)\le \ln\bigl(\pi^{\S}_{2,|K|^2-2}(K)\bigr)<\Bigl(\frac{27p}{2(p-1)}+\varepsilon\Bigr)|K|\ln |K|.$$

{\bf (3)} There exists $E(\varepsilon)\in\mathbb N^{\ast}$ such that for each finite field $K$ with $|K|\ge E(\varepsilon)$ and $4\mid |K|$ we have
$$(2-\varepsilon) |K|\ln |K|<\ln\bigl(\pi_{2,|K|^2-2}(K)\bigr)\le \ln\bigl(\pi^{\S}_{2,|K|^2-2}(K)\bigr)<18|K|\ln |K|.$$

{\bf (4)} There exists $J(\varepsilon)\in\mathbb N^{\ast}$ such that for each finite field $K$ with $|K|\ge J(\varepsilon)$ (resp.\ with $|K|\ge J(\varepsilon)$ and $|K|$ an odd prime or a power of $3$) we have
$$(2-\varepsilon) |K|<\j_K\le 18|K|\;\;\;(\textup{resp.}\;\;\;(2-\varepsilon) |K|<\j_K\le 13.5|K|\Bigl(1+\frac{3}{|K|}\Bigr)$$
and $\frac{3}{|K|}$ can be replaced by $\frac{1}{|K|}$ if $p=3$).

{\bf (5)} Let $p\ge 5$ be a prime. There exists $J_p(\varepsilon)\in\mathbb N^{\ast}$ such that for each finite field $K$ with $|K|\ge J_p(\varepsilon)$ and $\char(K)=p$ we have
$$(2-\varepsilon) |K|<\j_K\le \Bigl(\frac{27p}{2(p-1)}+\varepsilon\Bigr)|K|.$$
\end{corollary}

\begin{proof}
The lower bound of part (1) follows from Theorem \ref{T16}(1), as the right hand side of Inequality (\ref{EQ52}) is asymptotic to $2|K|\ln|K|$. 

The lower bounds of parts (2) and (3) follow from Theorem \ref{T17}(1), as the right hand side of Inequality (\ref{EQ55}) is asymptotic to $2|K|\ln|K|$.

The upper bounds of part (1) and (3) with $18$ follow from Corollary \ref{C18}(1) and (2).

The upper bound of part (1) with $13.5+\epsilon$ follows from Corollaries \ref{C18}(1) and \ref{C4}(2).

The upper bound of part (2) follows from Corollary \ref{C18}(1) and the fact that we have $\limsup \limits_{q\to\infty} \nu^+_{p,p^{q-1}}(p^{2q})\le\frac{3p}{2(p-1)}$ by Theorem \ref{T2}(10).

The lower bound of parts (4) and (5) follows from Corollary \ref{C26}, as the right hand side of its displays are asymptotic to $2|K|$. 

The upper bounds of part (4) follow from Corollaries \ref{C18}(1) and \ref{C4}(2).

The upper bounds of part (5) follow from Corollaries \ref{C18}(1) and \ref{C4}(3).\end{proof}

\begin{remark}\normalfont\label{R10}
{\bf (1)} Suppose that $K$ is finite with $|K|\ge 3$. For $j\in\llbracket 1,\j_K\rrbracket$, let $\j(K,j)$ be the number of permutations $\sigma\in\Perm(K^2)$ of the form $\sigma=\prod_{i=1}^j a_i(K)$ with $a_i$ reduced non-affine triangular for each $i\in\llbracket1,j\rrbracket$ and $\dir(a_i)\neq\dir(a_{i+1})$ for each $i\in\llbracket1,j-1\rrbracket$. Let 
$$\N_{|K|}:=(|K|-1)\Bigl(\sum_{i=0}^{|K|-3} |K|^i\Bigr)=|K|^{|K|-2}-1<|K|^{|K|-2}.$$

We have the recursive inequalities $$\j(K,1)\le (|K|+1)\N_{|K|}=(|K|+1)(|K|^{|K|-2}-1)$$ (cf.\ Remark \ref{R3}(2) and Equation (\ref{EQ37.9}) applied to $j=1$ and $d_1\in\llbracket2,|K|-1\rrbracket$) and 
$$\j(K,j)\le|K|\N_{|K|}\j(K,j-1)$$
for $j\in\llbracket2,\j_K\rrbracket$.
By considering the cardinality of the quotient set $\Alt(|K|^2)/\ASL_2(K)$ we get that 
$$-1+\frac{(|K|^2)!}{2|K|^2(|K|+1)(|K|^2-|K|)}\le\sum_{j=1}^{\j_K} \j(K,j)\le (|K|+1)\N_{|K|}\frac{(|K|\N_{|K|})^{\j_K}-1}{|K|\N_{|K|}-1}.$$
In this way one can get an improvement of Corollary \ref{C26} which itself can be improved if one considers better bounds for the $\j(K,j)$s.
For instance, from Example \ref{EX12} we get that if $|K|=p^q$ with $p$ a prime and $q\in\mathbb N^{\ast}\setminus\{1\}$, then we have
$$\j(K,2)\le |K|(|K|+1)(|K|^{|K|-2}-1)^2- \frac{|K|(|K|+1)}{2}\Bigl[\sum_{i=1}^{q-1} \binom{n}{i}_p(p^i-1)\Bigr]^2,$$
where $\binom{n}{i}_p$ is the $p$-binomial coefficient that computes the number of $d$-dimensional $\mathbb F_p$-vector subspaces of $\mathbb F_p^q$. 

\smallskip
{\bf (2)} Based on Remark \ref{R5} and the general result \cite{BBS}, Thm.\ 1.1 on diameters of alternating groups (see also \cite{BGHHSS}, Thm.\ 1.2) one gets that for $|K|\ge 4$, $\j_K$ has an upper bound of the form $|K|^{14}\log(|K|^2)^{O(1)}$.
\end{remark}

\begin{theorem}\label{T18} Let $\bigl((K_i,m_i)\bigr)_{i\in\mathbb N^{\ast}}$ be a sequence of pairs with each $K_i$ a finite field and $m_i\in \llbracket2,|K_i|^2\rrbracket$. If $\lim \limits_{i\to +\infty}|K_i|=+\infty$, then the following statements are equivalent.

\medskip
{\bf (1)} There exists $C\in (0,\infty)$ such that $\pi^{\S}_{2,m_i}(K_i) <m_i^C$ for each $i\in\mathbb N^{\ast}$.

\smallskip
{\bf (2)} There exists $C\in (0,\infty)$ such that $\pi_{2,m_i}(K_i) <m_i^C$ for each $i\in\mathbb N^{\ast}$.

\smallskip
{\bf (3)} The set $\{\frac{m_i}{|K_i|}|i\in\mathbb N^{\ast}\}$ of rational numbers is bounded.
\end{theorem}

\begin{proof}
Clearly, $(1)\Rightarrow (2)$.

{$\pmb{(2)\Rightarrow (3)}$.} We prove the contrapositive $\neg (3)\Rightarrow\neg (2)$; thus we assume that the set $\{\frac{m_i}{|K_i|}|i\in\mathbb N^{\ast}\}$ is unbounded. By taking a subsequence, we can assume that $\lim \limits_{\i\to +\infty}\frac{m_i}{|K_i|}=+\infty$. By replacing $m_i$ by $\min(m_i,|K_i|^2-2)$ we can also assume that $m_i\leq |K_i|^2-2$. As $[(m_i-2)-(|K_i|^2-2)]\ln\bigl(1-\frac{m_i-2}{|K_i|^2-2}\bigr)\ge 0$ we get that
$$(m_i-2)\ln (|K_i|^2-m_i)-(|K_i|^2-2)\ln\Bigl(1-\frac{m_i-2}{|K_i|^2-2}\Bigr)\ge (m_i-2)\ln (|K_i|^2-2).$$ Thus, by Theorem \ref{T17}(1) we have
$$\ln \bigl(\pi_{2,m_i}(K_i) \bigr) >\frac{(m_i-2)\ln (|K_i|^2-2)-m_i+2-\ln(1-|K_i|^{-1})}{|K_i|+\mathfrak c_{|K_i|}}$$
$$=\ln|K_i|\frac{\frac{m_i-2}{|K_i|}\bigl(\frac{\ln (|K_i|^2-2)}{\ln|K_i|}-\frac{1}{\ln|K_i|}\bigr)-\frac{\ln(1-|K_i|^{-2})}{|K_i|\ln |K_i|}}{1+\frac{\mathfrak c_{|K_i|}}{|K_i|}}.$$
As $\lim \limits_{i\to +\infty }|K_i|=\lim \limits_{i\to +\infty }\frac{m_i-2}{|K_i|}=+\infty$ implies
$\lim \limits_{i\to +\infty}\frac{\ln \bigl(\pi_{2,m_i}(K_i) \bigr)}{\ln |K_i|}=+\infty$, $\neg (2)$ holds.

{$\pmb{(3)\Rightarrow (1)}$.} Suppose there exists $D\in (0,\infty)$ such that $m_i <D|K_i|$ for each $i\in\mathbb N^{\ast}$. Taking $C\ge 4$, for each $i\in\mathbb N^{\ast}$ such that $m_i < \frac{|K|}{\sqrt{2}}$ we have $$\pi^{\S}_{2,m_i}(K_i) <m_i^4\le m_i^C$$ by the identity $\pi_{2,2}(K_i)=1$ and Theorem \ref{T8}(1.b) applied when $m_i\ge 3$. For $i\in\mathbb N^{\ast}$ such that $\frac{|K_i|}{\sqrt{2}}\le m_i\le D|K_i|$ we use the proof of Theorem \ref{T12}, specifically the first inequality of Inequalities (\ref{EQ36}). With $r_i:=\bigl\lfloor \frac{\lceil\frac{|K_i|}{\sqrt{2}}\rceil}{3}\bigr\rfloor$, if $|K_i|\ge 11$ we have $r_i\ge 3$ and $r_i>\frac{|K_i|-\sqrt{8}}{3\sqrt{2}}$ and hence we get that
$$\ln \bigl(\pi^{\S}_{2,m_i}(K_i)\bigr) \le 12\ln(|K_i|-1)\Bigl(\frac{m_i-4}{r_i}+2\Bigr)<12\Bigl(D\frac{3\sqrt{2}|K_i|}{|K_i|-2\sqrt{2}}+2\Bigr) \ln (|K_i|).$$ It follows that there exists a constant $E\in (0,\infty)$ that does not depend on $i$ and such that $\ln \bigl(\pi^{\S}_{2,m_i}(K_i)\bigr) <E\ln |K_i|\le E (\ln m_i + \ln \sqrt{2}) <(E+1)\ln m_i$. So statement (1) holds for $C:=\max(4,E+1)$.
\end{proof}

The next remarks point out that Proposition \ref{PR30} can be improved in multiple small ways. 

\begin{remark}\normalfont\label{R11} 
For each sequence $\underline{d}=(d_i)_{i\in\llbracket1,j\rrbracket}$ in $\mathbb N^{\ast}\setminus\{1\}$, we consider the set $\mathcal F_j=\mathcal F_j(\underline{d},|K|)$ of automorphisms $a\in\GA_2(K)$ that admit reduced standard product decompositions $a=a_1\cdots a_jb$ with $a_i$ a reduced non-triangular automorphism with $\ell(a_i)=d_i$ for every $i\in\llbracket1,j\rrbracket$. 

\medskip
{\bf (1)} For $|K|\ge 5$ we have $\frac{|K|^3}{|K|+1}\le 2^{|K|}\cdot\frac{|K|-1}{|K|}$ and hence the factor $(|K|^2-1)^2$ in Proposition \ref{PR30})(1) can be replaced by $(|K|^2-1)^2\frac{|K|-1}{|K|}$ as in the proof of Proposition \ref{PR30} the factor $\bigl(1-\frac{1}{|K|}\bigr)^j$ was rounded up to $1$. 

\smallskip
{\bf (2)} Improvements can be made by utilizing Remark \ref{R10} or the fact that $\frac{d_i}{\ln (d_i)}$ is much smaller than $\frac{|K|}{\ln |K|}$ for small $d_i\in\llbracket1,|K|-1\rrbracket$, among which the latter extrapolates into the fact that the first inequalities in Lemma \ref{L25}(2) and Inequalities (\ref{EQ50}) and (\ref{EQ51}) are often strict.

\smallskip
{\bf (3)} Such improvements would lead to slightly better bounds in Theorems \ref{T16} and \ref{T17}. However, they would not affect the results or simplify the proofs of the qualitative results Corollary \ref{C27} and Theorem \ref{T18}, for the following reason. For $|K|>2$ and $N=(|K|-1)^j$ with $j\in\llbracket1,\j_K\rrbracket$ we can consider the finite sequence $\underline{d}:=(d_i)_{i\in \llbracket1,j\rrbracket}$ with $d_i=|K|-1$ for each $i\in \llbracket1,j\rrbracket$. The contribution of this one sequence to the upper bound for $\varrho_{2,K, N}$ is
$$B_{(d_i)_{i\in \llbracket1,j\rrbracket}}=(|K|^2-1)^2(1-|K|^{-1})^j|K|^{\sum_{i=1}^j d_i}$$
$$=|K|^{j(|K|-2)}(|K|^2-1)^2 (|K|-1)^j>(|K|^2-1)^2(|K|-1)^{j(|K|-1)}=(|K|^2-1)^2N^{|K|-1}.$$

Due to this, to get an asymptotically better power of $N$ in the estimate in Proposition \ref{PR30} one would have to prove that many of the automorphisms in $\mathcal F_j$ produce the same permutations. 

\smallskip
{\bf (4)} If $|K|$ is not a prime we have $\j_K\le 18|K|$, if $|K|$ is a prime different from $3$ we have $\j_K\le 13.5|K|\bigl(1+\frac{3}{|K|}\bigr)$, and if $\char(K)=3$ with $|K|\ge 27$ we have $\j_K\le 13.5|K|\bigl(1+\frac{1}{|K|}\bigr)$ by Corollary \ref{C18}(4) and (5). For each $a\in\mathcal F_j$, the number of automorphisms $b\in\mathcal F_j$ with $b(K)=a(K)$ is equal to the number of $K$-valued points of a quasi-affine variety $\mathcal V_{a(K),j}$ over $\Spec K$ in $\mathbb A^{j(|K|-2)+4}\sqcup \mathbb A^{j(|K|-2)+4}$. Here the two copies relate to the fact that for the $|K|+1$ possible directions, it is more convenient to use $2$ parameters, one of them taking only $2$ values. If $j\in\llbracket |K|,\j_K\rrbracket$, then $\mathcal V_{a(K),j}$ has at least one irreducible component $\mathcal V_{a(K),j}^0$ which is of dimension at least $|K|(j-|K|)+2$ and which has $K$-valued points. Without a more explicit description of the $\mathcal V_a$s, it is hard to estimate in a refined way the number of $b$s. 

\smallskip
{\bf (5)} Note that $\mathcal F_j\subset \mathcal F_{j+2}$ if $j\le \j_K-2$. The set $\{l\in\llbracket0,\j_K-j\rrbracket|a\in\mathcal F_{j+l}\}$ contains $\llbracket0,\j_K-j\rrbracket\cap 2\mathbb N$ and $|\mathcal V_{a,j+l}^0(K)|$ is an increasing function on such $l$s.

To give examples, let $\sigma_0$ be the identity element of $\perm(K^2)$ and we write $|K|=p^q$ with $p$ a prime and $q\in\mathbb N^{\ast}$. 

If $l\in\{1,\ldots,j-1\}$ is such that $a\in\mathcal F_{j-l}$, then $|\mathcal V_{a(K),j}(K)|\ge \frac{|K|}{|K|+1}|\mathcal V_{\sigma_0,l}(K)|$ as for $b\in\mathcal F_{j-l}$ that corresponds to a point in $\mathcal V_{\sigma_0,j-l}(K)$, if $b=b_1\cdots b_{l}c_b$ and $\break a=a_1\cdots a_{j-l}c_a$ are standard product decompositions, then $ba\in\mathcal F_j$ provided $\dir(b_l)\neq\dir(a_1)$ and the number of such $b$s with $\dir(b_l)\neq\dir(a_1)$ is $\frac{|K|}{|K|+1}|\mathcal V_{\sigma_0,l}(K)|$.

If $\j_K\ge 4$ and $q\ge 2$, then from Example \ref{EX12} we get that 
\begin{equation}\label{EQ57}
|\mathcal V_{\sigma_0,4}(K)|\ge |K|(|K|+1)\Bigl[\sum_{i=1}^{q-1} \binom{n}{i}_p(p^i-1)\Bigr]^2,
\end{equation}
where $\binom{n}{i}_p$ is as Remark \ref{R10}(1). Thus if we have $5\le j\le\j_K$, then for each $a\in\mathcal F_{j-4}$ we have an inequality $|\mathcal V_{a(K),j}(K)|\ge |K|^2[\sum_{i=1}^{q-1} \binom{n}{i}_p(p^i-1)^2]$. 

If $\j_K\ge 6$, then from Example \ref{EX14}(2) we get that 
$$|\mathcal V_{\sigma_0,6}(K)|\ge |K|(|K|+1)(|K|-1)^2.$$ 
Thus if we have $7\le j\le\j_K$, then for each $a\in\mathcal F_{j-6}$ we have an inequality $|\mathcal V_{a(K),j}(K)|\ge |K|^2(|K|-1)^2$. 

From Example \ref{EX14}(1) and Inequality (\ref{EQ57}) we get that for $\j_K\ge 8$ we have $|\mathcal V_{\sigma_0,8}(K)|\ge 3!\binom{|K|^2}{3}$ and for $\j_K\ge 8$ and $q\ge 2$ we have
$$|\mathcal V_{\sigma_0,8}(K)|\ge 3!\binom{|K|^2}{3}+|K|^3(|K|+1)\Bigl[\sum_{i=1}^{q-1} \binom{n}{i}_p(p^i-1)\Bigr]^4.$$ 
\end{remark}

Directly from Theorem \ref{T18} we get the following result.

\begin{corollary}\label{C28}
Let $\varepsilon\in (0,1)$. Then for each pair $(C,N)\in (0,\infty)\times\mathbb N^{\ast}$, there exists $D\in\mathbb N^{\ast}$ such that for every finite field $K$ with $|K|\ge D$ we have $\pi_{2,\lfloor|K|^{1+\varepsilon}\rfloor}(K)>C|K|^{N+N\varepsilon}$.
\end{corollary}

\section{Values and estimates for small $m$}\label{S28}

In this section we use prior bounds to either provide ranges for or compute the $\pi_{n,m}(K)$s when $m\in\llbracket3,7\rrbracket$. 

For a field $K$ and $(n,m)\in (\mathbb N^{\ast})^2$, let 
$$\mathbb D_{n,m}^{\sum=0}(K):=\bigl\{(P_1,\ldots,P_m)\in\mathbb D_{n,m}(K)|\sum_{i=1}^m P_i=0\bigr\}$$
and $\mathbb D_{n,m}^{\sum\neq 0}(K):=\mathbb D_{n,m}(K)\setminus\mathbb D_{n,m}^{\sum=0}(K)$. If $K$ is finite, then for $|K|\neq 2$ we have $\mathbb D_{n,|K|}^{\d=1}\subset \mathbb D_{n,K|}^{\sum=0}$ and for $|K|=2$ we have $\mathbb D_{n,|K|}^{\d=1}\subset \mathbb D_{n,K|}^{\sum\neq0}$.

We begin with a general lemma for fields $K$ of characteristic $2$.

\begin{lemma}\label{F13}
Let $K$ be a field of characteristic $2$. Let $m\in 2\mathbb N^{\ast}$ be such that we have $\sqrt{m}\le |K|$. Let $\Gamma_2(K)[2]$ be the subgroup of $\GA_2(K)$ generated by $\GA_2(K)[2]$. Then the following properties hold.

\medskip
{\bf (1)} The action $\mathbb T_{2,m}$ restricts to actions $\Gamma_2(K)[2]\times\mathbb D_{2,m}^{\sum=0}(K)\rightarrow \mathbb D_{2,m}^{\sum=0}(K)$ and $\Gamma_2(K)[2]\times\mathbb D_{2,m}^{\sum\neq 0}(K)\rightarrow \mathbb D_{2,m}^{\sum\neq 0}(K)$.

\smallskip
{\bf (2)} Suppose that both sets $\mathbb D_{2,m}^{\sum= 0}(K)$ and $\mathbb D_{2,m}^{\sum\neq 0}(K)$ are non-empty. If $K$ is a finite field with $4\mid  |K|$, then we assume that $m\le |K|-2$. Then we have $\Gamma_2(K)[2]\neq\GA_2(K)$.
\end{lemma}

\begin{proof}
For part (1), it suffices to show that for each $a\in\AGL_2(K)$ we have $a\bigl(\mathbb D_{2,m}^{\sum=0}(K)\bigr)\subset\mathbb D_{2,m}^{\sum=0}(K)$ and for every $f\in K[x]$ with $\deg(f)\le 2$, by denoting $b:=\e\bigl(x_1,x_2+f(x_1)\bigr)$ we have $b\bigl(\mathbb D_{2,m}^{\sum=0}(K)\bigr)\subset\mathbb D_{2,m}^{\sum=0}(K)$. 

To check these two properties, let $\underline{P}:=(P_1,\ldots,P_m)\in\mathbb D_{2,m}^{\sum=0}(K)$. Writing $P_i=(\alpha_i,\beta_i)\in K^2$ for $i\in\llbracket1,m\rrbracket$, we have $\sum_{i=1}^m \alpha_i=\sum_{i=1}^m \beta_i=0$. 

As $m$ is even and $\char(K)=2$, the equation $\sum_{i=1}^m P_i=0$ is equivalent to $\sum_{i=1}^m (P_i-P_1)=0$ and thus it is invariant under affine automorphisms. Therefore we have $a(\underline{P})\in \mathbb D_{2,m}^{\sum=0}(K)$. 

Writing $f(x)=\gamma_2 x^2+\gamma_1 x+\gamma_0$ with $(\gamma_2,\gamma_1,\gamma_0)\in K^3$, we compute that 
$$\sum_{i=1}^m \beta_i+f(\alpha_i)=\Bigl(\sum_{i=1}^m\beta_i\Bigr)+\gamma_0\Bigl(\sum_{i=1}^m \alpha_i\Bigr)^2+\gamma_1\Bigl(\sum_{i=1}^m \alpha_i\Bigr)+m\gamma_0=0.$$
Thus $b(\underline{P})\in \mathbb D_{2,m}^{\sum=0}(K)$. So part (1) holds. 

Part (2) follows from part (1) and the fact that $\GA_2(K)$ acts transitively on $\mathbb D_{2,m}(K)$ by Theorem \ref{T1}(3).
\end{proof}

We have the following application of Proposition \ref{PR16.5}(6) for $m=4$.

\begin{lemma}\label{L26}
Let $(\underline{P},\underline{Q})=\bigl((P_1,P_2,P_3,P_4),(Q_1,Q_2,Q_3,Q_4)\bigr)\in\mathbb D_{2,4}(\mathbb F_4)$. Let $Y_1:=\{P_i|i\in\llbracket1,4\rrbracket\}$ and $Y_2:=\{Q_i|i\in\llbracket1,4\rrbracket\}$. Then the following properties hold.

\medskip
{\bf (1)} If $Y_1$ and $Y_2$ are both non-collinear, then $\pi^{\le 2}_{\underline{P},\underline{Q}}\le 6$.

\smallskip
{\bf (2)} If only one of the sets $Y_1$ and $Y_2$ is collinear, then $\pi^{\le 2}_{\underline{P},\underline{Q}}\le 6$.
\end{lemma}

\begin{proof}
Let $\gamma\in\mathbb F_4\setminus\mathbb F_2$. Thus $\mathbb F_4=\{0,1,\gamma,\gamma+1\}$. 

For part (1), for $\iota\in\{1,2\}$ let $B_{\iota}$ be the subset of $\llbracket1,4\rrbracket$ defined as follows. If $\c_{Y_{\iota}}=2$, then $B_{\iota}:=\emptyset$. If $\c_{Y_{\iota}}=3$ and $\iota=1$ (resp.\ $\iota=2$), let $B_{\iota}:=\{l\}$ with $l$ the unique element such that $Y_1\setminus\{P_l\}$ (resp.\ $Y_2\setminus\{Q_l\}$) is collinear. For every subset $B$ of $\llbracket1,4\rrbracket$ such that $|B|=3$ and $B_1\cup B_2\subset B$, the sets $\{P_i|i\in B\}$ and $\{Q_i|i\in B\}$ are affinely independent. Up to a reindexing, we can assume that $B=\{1,2,3\}$. Up to affine automorphisms we can assume that $P_1=Q_1=(0,0)$, $P_2=Q_2=(1,0)$, and $P_3=Q_3=(0,1)$. Let $Z:=\{P_1,P_2,P_3\}$. We have $\{P_4,Q_4\}\subset\mathbb F_4^2\setminus Z$. 

For each transposition $\tau\in\perm(\mathbb F_4^2)$, let $\tau_{\s}\in\perm(\mathbb F_4^2)$ be the transposition with $\supp(\tau_{\s})=\lambda_{\tau}\setminus\supp(\tau)$, where $\lambda_{\tau}$ is the line generated by $\supp(\tau)$, and let $\tau^+:=\tau\tau_{\s}$. So $o(\tau^+)=2$, $\supp(\tau^+)=\lambda_{\tau}$, and $\n(\tau^+)=4$. Let $a_{\tau^+}\in\SGA_2(\mathbb F_4)[3]$ be such that $a_{\tau^+}(\mathbb F_4)=\tau_+$ and $\j_{a_{\tau^+}}=1$. All such $a_{\tau^+}$ are of the form $b\shift_{V\oplus W}^{v_0+w_0}b^{-1}$ with $\dim_{\mathbb F_4}(V)=\dim_{\mathbb F_4}(W)=1$ and $b\in\AGL_2(\mathbb F_4)$ such that $\lambda_{\tau}=b(V+w_0)$. 

As in the case $Q_4=P_4$ we can take $a=1_{\mathbb A^2_K}$, we can assume that $Q_4\neq P_4$. To show that there exists $a\in\Fix_{\GA_2(\mathbb F_4)}(Z)[6]$ such that $\j_a\le 2$ and $a(P_4)=Q_4$, based on the existence of the automorphism of $\mathbb F_4$ that maps $\gamma$ to $\gamma+1$ and on the fact that the subgroup of $\AGL_2(\mathbb F_2)$ that stabilizes the set $Z$ (equivalently, fixes $(1,1)$) has $6$ elements, we can also assume that $P_4\in\{(1,1),(\gamma,1),(\gamma,0)\}$. 

{\bf Case 1: $P_4=(1,1)$.} By the same reason we can assume that $Q_4\in\{(\gamma,0),(\gamma,1)\}$. Let $a_1:=a_{\bigl(P_4\,(\gamma,0)\bigr)^+}\in\SGA_2(\mathbb F_4)[3]$. If $Q_4=(\gamma,0)$, then we can take $a:=a_1$. If $Q_4=(\gamma,1)$, then we can take $a:=\e(x_1,x_2+x_1^2+x_1)a_1\in\SGA_2(\mathbb F_4)[6]$. 

{\bf Case 2: $P_4=(\gamma,1)$.} If $Q_4$ does not belong to a line that passes through $P_4$ and a point in $Z$, then we can take $a:=a_{(P_4\,Q_4)^+}\in\SGA_2(\mathbb F_4)[3]$. Based on Case 1 we can assume that $Q_4\neq (1,1)$. Based on the last two sentences we can assume that $Q_4\in\{(\gamma+1,1),(\gamma+1,\gamma),(1,\gamma+1),(0,\gamma),(\gamma+1,\gamma+1)\}$. If the first coordinate of $Q_4$ is $\gamma+1$, then we can take $a:=\e(x_1,x_2+\beta x_1^2+\beta x_1)a_{\bigl(P_4\,(\gamma+1,0)\bigr)^+}$ in $\SGA_2(\mathbb F_4)[6]$, where $\beta\in\mathbb F_4$ is the second coordinate of $Q_4$. If $Q_4=(\epsilon,\gamma+\epsilon)$ with $\epsilon\in\{0,1\}$, then we can similarly take $a:=\e(x_1+x_2^2+x_2,x_1)a_{\bigl(P_4\,(1+\epsilon,\gamma+\epsilon)\bigr)^+}\in\SGA_2(\mathbb F_4)[6]$. Note that if $Q_4\neq (1,1)$, then all automorphisms of this case fix $(1,1)$ and hence the case holds if $P_1=Q_1$ is $(1,1)$ and not $(0,0)$.

{\bf Case 3: $P_4=(\gamma,0)$.} Via the translation with $(0,1)$ and the reindexing of the points, this case gets reduced to Case 2. 

As in all three cases above we have $\j_a\le 2$, we get that part (1) holds.

For part (2), we can assume that $\c_{Y_1}=4$ and, up to affine automorphisms and the automorphism of $\mathbb F_4$ that interchanges $\gamma$ and $\gamma+1$, that $P_1=Q_1=(0,0)$, $P_2=Q_2=(1,0)$, $P_3=(\gamma+\varepsilon,0)$ and $P_4=(\gamma+1+\varepsilon,0)$ with $\varepsilon\in\{0,1\}$, $Q_3=(0,1)$, and $Q_4\in\{(1,1),(1,\gamma),(0,\gamma)\}$. If $Q_4=(1,1)$, then we can take $a:=\e(x_1+(\gamma+\varepsilon)x_2,x_2)\e(x_1,x_2+x_1^2+x_1)\in\SGA_2(\mathbb F_4)[2]$. If $Q_4\neq (1,1)$, then we can take $a:=\e\bigl(x_1+g(x_2),x_2\bigr)\e\bigl(x_1,x_2+f(x_1)\bigr)\in\SGA_2(\mathbb F_4)[6]$ with $(f,g)\in\mathbb F_4[x]^2$ such that $\deg(g)\le 2$, $\deg(f)\le 3$, $f(0)=f(1)=g(0)=0$, $f(\gamma+\varepsilon)=1$, $f(\gamma+1+\varepsilon)=\gamma$, $g(1)=\gamma+\varepsilon$, and $g(\gamma)+\gamma+1+\varepsilon$ is the first coordinate of $Q_4$. As $\j_a\le 2$ in all cases, part (2) holds.
\end{proof}

\begin{theorem}\label{T19}
Let $(n,m)\in (\mathbb N^\ast\setminus\{1\})^2$ and a field $K$ be such that $\sqrt[n]{m}\le |K|$. Then the following properties hold.

\medskip
{\bf (1)} We have $\pi_{n,3}(\mathbb F_3)=\pi^{\S}_{n,3}(\mathbb F_3)=2$ and $\pi_{n,3}(K)=\pi^{\S}_{n,3}(K)=4$ for $|K|\ge 4$. 

\smallskip
{\bf (2)} For $|K|=5$ we have $9=\pi_{2,4}(K)\le\pi^{\S}_{2,4}(K)\le 36$ and for $|K|\ge 7$ we have identities $\pi_{2,4}(K)=\pi^{\S}_{2,4}(K)=9$. Moreover, for $|K|\ge 5$ and $\n\ge 3$ we have $5\le \pi_{n,4}(K)\le\pi^{\S}_{n,4}(K)=9$.

\smallskip
{\bf (3)} We have $\pi_{2,4}(\mathbb F_4)=6$ and $4\le\pi_{n,4}(\mathbb F_4)\le\pi^{\S}_{n,4}(\mathbb F_4)\le 6$ for $n\ge 3$.

\smallskip
{\bf (4)} If $n\ge 3$, then $\pi_{n,4}(\mathbb F_2)=2$.

\smallskip
{\bf (5)} We have $\pi_{2,4}(\mathbb F_3)=4$ and $2\le\pi_{n,4}(\mathbb F_3)\le \pi^{\S}_{n,4}(\mathbb F_3)=4$ for $n\ge 3$.

\smallskip
{\bf (6)} We have $2\le \pi_{3,5}(\mathbb F_2)\le\pi_{3,6}(\mathbb F_2)\le 3$ and $2\le\pi_{3,7}(\mathbb F_2)=\pi_{3,8}(\mathbb F_2)\le 6$.

\smallskip
{\bf (7)} We have the following inequalities $6\le \pi_{n,5}(\mathbb F_5)\le \pi^{\S}_{n,5}(\mathbb F_5)\le 12$ for $n\ge 4$, $6\le \pi_{3,5}(\mathbb F_5)\le 16$, $8\le \pi_{2,5}(\mathbb F_5)\le 16$, and $8\le \pi^{\S}_{2,5}(\mathbb F_5)\le 36$. In particular, $\pi_{2,5}(\mathbb F_5)\in\{8,9,12,16\}$.

\smallskip
{\bf (8)} If $K$ contains a primitive $5$-th root of unity, then $\pi_{2,5}(K)=16$.

\smallskip
{\bf (9)} We have the following inequalities $7\le \pi_{n,6}(\mathbb F_7)\le \pi^{\S}_{n,6}(\mathbb F_7)\le 80$ for $n\ge 3$, $25\le \pi_{2,6}(\mathbb F_7)\le 80$, and $25\le\pi^{\S}_{2,6}(\mathbb F_7)\le 400$.

\smallskip
{\bf (10)} If $|K|\ge 8$, then we have $7\le \pi_{n,6}(K)\le 25$ for $n\ge 3$, $16\le \pi_{2,6}(K)\le 25$ if $|K|\not\equiv 1\pmod{6}$, and $\pi_{2,6}(K)=25$ if $|K|\equiv 1\pmod{6}$. In particular, for $|K|\in\{8,9\}$ we have $\pi_{2,6}(K)\in\{16,18,20,21,24,25\}$.

\smallskip
{\bf (11)} We have the following inequalities $8\le \pi_{n,7}(\mathbb F_7)\le \pi^{\S}_{n,7}(\mathbb F_7)\le 30$ if $n\ge 6$, $8\le \pi_{n,7}(\mathbb F_7)\le \pi^{\S}_{n,7}(\mathbb F_7)\le 150$ if $n\in\{3,4,5\}$, and $20\le \pi_{2,7}(\mathbb F_7)\le \pi^{\S}_{2,7}(\mathbb F_7)\le 900$.
\end{theorem}

\begin{proof}
For part (1), we have $\pi_{n,3}(\mathbb F_3)\le\pi^{\S}_{n,3}(\mathbb F_3)\le 2$ by Theorem \ref{T7}(3). As $\AGL_2(K)$ acts on $\mathbb D^{\d=1}_{n,3}(\mathbb F_3)$ and $\mathbb D^{\d=2}_{n,3}(\mathbb F_3)$, we have $2\le\pi_{n,3}(\mathbb F_3)$. So part (1) holds for $\mathbb F_3$. As $(3-1)^2=3+1$, for $|K|\ge 4$ part (1) follows from Theorems \ref{T7}(2) and \ref{T15}(1). So part (1) holds.

Part (2) holds for $n\ge 3$ as $\pi^{\S}_{n,4}(K)\le 9$ by Theorem \ref{T7}(2) and $5\le \pi_{n,4}(K)$ by Theorem \ref{T15}(1). We check that part (2) holds for $n=2$. 

We have $8\le\pi_{2,4}(K)$ by Theorem \ref{T15}(4). For $|K|\ge 7$ we have $\pi^{\S}_{2,4}(K)\le 9$ by Theorem \ref{T7}(1). For $|K|=5$ we have $\pi_{2,4}(K)\le 9$ by Theorem \ref{T8}(4) and $\pi^{\S}_{2,4}(K)\le 36$ by Theorem \ref{T8}(1.b). If $|K|\equiv 1\pmod{4}$, then $K$ has a primitive $4$-th root of unity and thus $\pi_{2,4}(K)\ge 9$ by Theorem \ref{T15}(5).

Thus to complete the proof of part (2) we are left to show that for $|K|\ge 7$ with $|K|\not\equiv 1\pmod{4}$, the assumption that $\pi_{2,4}(K)=8$ leads to a contradiction. 

Let $(\underline{P},\underline{Q})=\bigl((P_1,P_2,P_3,P_4),(Q_1,Q_2,Q_3,Q_4)\bigr)\in\mathbb D_{2,4}(K)^2$.

Let $\bigl((\alpha_1,\alpha_2,\alpha_3,\alpha_4),(\beta_1,\beta_2,\beta_3,\beta_4)\bigr)\in\mathbb D_{1,4}(K)^2$ be $\{1,2,3\}$-generic by Proposition \ref{PR28}(2) or by the concrete examples below. We take 
$$(\underline{P},\underline{Q})=\bigl((P_1,P_2,P_3,P_4),(Q_1,Q_2,Q_3,Q_4)\bigr)\in\mathbb D_{2,4}(K)^2$$ 
such that $P_i:=(\alpha_i,0)$ and $Q_i=(\beta_i,0)$ for each $i\in\llbracket1,4\rrbracket$. We have $\pi_{\underline{P},\underline{Q}}^{\S,\le 2}=9$ by Corollary \ref{L23}. From this and our assumption we get that there exists a special automorphism $a\in\SGA_2(K)[8]$ such that $\j_a\ge 3$ and $a(\underline{P})=\underline{Q}$. 

As $\ell(a)\le 8$ and $\j_a\ge 3$ and as $8=2\cdot 2\cdot 2$ is the only way to write a number in $\llbracket1,8\rrbracket$ as a product of three numbers in $\mathbb N^{\ast}\setminus\{1\}$, from Proposition \ref{PR11}(4.a) we get that $\ell(a)=8$, $\j_a=3$, and we have a product decomposition $a=a_1a_2a_3$ with each $a_i$ a non-affine triangular automorphism in $\GA_2(K)[2]$. 

To reach a contradiction we take $\alpha_1=\beta_1=0$, $\alpha_2=\beta_2=1$. If $\char(K)=2$, then $8\mid |K|$ and we choose first $\alpha\in K\setminus\{0,1\}$ and second $\beta\in K\setminus\{0,1\}$ such that $\beta^2+\beta+1\neq 0$ and $(\beta^4+\alpha+1)(\beta+\alpha)\neq (\beta^2+\alpha)(\beta^2+\alpha+1)$; as $0$ and $1$ are solutions in $K$ of the equation $(x^4+\alpha+1)(x+\alpha)-(x^2+\alpha)(x^2+\alpha+1)=0$, such a $\beta$ exists as we are excluding at most $7$ values in $K$. If $\char(K)\neq 2$, then $|K|\ge 7$ and we choose $\alpha\in K\setminus\{0,1,-1,-2,-2^{-1}\}$ and $\beta:=\alpha+1$. If $\char(K)=2$ we take $\alpha_3=\alpha$, $\alpha_4=\alpha+1$, $\beta_3=\beta$, and $\beta_4=\beta^2$ and if $\char(K)\neq 2$ we take $\alpha_3=\beta_3=\alpha$ and $\alpha_4=-\beta_4=\beta=\alpha+1$. The assumptions on $\alpha$ and $\beta$ are such that Lemma \ref{L22} gives that the pair $\bigl((\alpha_1,\alpha_2,\alpha_3,\alpha_4),(\beta_1,\beta_2,\beta_3,\beta_4)\bigr)$ is $\{1,2,3\}$-generic. 

If, $\char(K)=2$, then we have 
$$(\underline{P},\underline{Q})=\bigl((P_1,P_2,P_3,P_4),(Q_1,Q_2,Q_3,Q_4)\bigr)\in\mathbb D^{\sum=0}_{2,4}(\mathbb F_4)\times \mathbb D^{\sum\neq 0}_{2,4}(\mathbb F_4),$$
and from this and the identity $a(\underline{P})=\underline{Q}$ we get a contradiction to Lemma \ref{F13}(1).

Thus we can assume that $\char(K)\neq 2$. Up to translations we can also assume that each $a_i$ fixes $P_1$. As $GL_2(K)$ acts $2$-transitively on the lines passing through the origin, we can assume that there exist two triples $(f,g,h)\in (xK[x])^3$ and $(c_1,c_2,c_3)\in\GL_2(K)^3$ such that $\deg(f)=\deg(g)=\deg(h)=2$ and we have $a=c_3^{-1}\e\bigl(x_1,x_2-h(x_1)\bigr)c_2^{-1}e\bigl(x_1+g(x_2),x_2)\bigr)e\bigl(x_1,x_2+f(x_1)\bigr)c_1$. Thus
$$[e\bigl(x_1+g(x_2),x_2)\bigr)e\bigl(x_1,x_2+f(x_1)\bigr)c_1](\underline{P})=[c_2\e\bigl(x_1,x_2+h(x_1)\bigr)c_3](\underline{Q}).$$
Let $\bigl((\gamma_1,\gamma_2),(\delta_1,\delta_2)\bigr)\in (K^2\setminus\{(0,0)\})^2$ be such that $c_1\bigl((1,0)\bigr)=(\gamma_1,\gamma_2)$ and $c_3\bigl((1,0)\bigr)=(\delta_1,\delta_2)$. We write $c_2(x_1,x_2)=(\delta_{11}x_1+\delta_{12}x_2,\delta_{21}x_1+\delta_{22}x_2)$.

We get that for $i\in\llbracket1,4\rrbracket$, $\bigl(\alpha_i\gamma_1+g(\alpha_i\gamma_2+f(\alpha_i\gamma_1)),\alpha_i\gamma_2+f(\alpha_i\gamma_1)\bigr)$ is equal to 
$$\bigl(\delta_{11}\beta_i\delta_1+\delta_{12}\beta_i\delta_2+\delta_{12}h(\beta_i\delta_1),\delta_{21}\beta_i\delta_1+\delta_{12}\beta_i\delta_2+\delta_{22}h(\beta_i\delta_1)\bigr).$$
If $\gamma_1=0$, then by replacing $f$ with $f(0)$ we get that $\pi_{\underline{P},\underline{Q}}^{\S,\le 2}\le 4$, a contradiction to $\pi_{\underline{P},\underline{Q}}^{\S,\le 2}=9$. Thus $\gamma_1\neq 0$. 

We have $\alpha_i\gamma_2+f(\alpha_i\gamma_1)=\delta_{21}\beta_i\delta_1+\delta_{12}\beta_i\delta_2+\delta_{22}h(\beta_i\delta_1)$ for each $i\in\llbracket1,4\rrbracket$. From the $\{1,2,3\}$-generic property we get that the only pairs $(f_1,g_1)\in (xK[x])^2$ such that $\deg(f_1)=\deg(g_1)=2$ and $f_1(\alpha_i)=g_1(\beta_i)$ for each $i\in\llbracket1,4\rrbracket$ are of the form $(\gamma x^2,\gamma x^2)$ with $\gamma\in K^{\ast}$. Thus there exists $\gamma\in K^{\ast}$ such that 
$$\gamma_2x+f(\gamma_1x)=(\delta_{21}\delta_1+\delta_{12}\delta_2)x+\delta_{22}h(x\delta_1)=\gamma x^2.$$
It follows that for $f_0(x):=\gamma_1x+g(\gamma x^2)$ and $g_0(x):=(\delta_{11}\delta_1+\delta_{12}\delta_2)x+\delta_{12}h(x\delta_1)$ we have $f_0(\alpha_i)=g_0(\beta_i)$ for each $i\in\llbracket1,4\rrbracket$.

Based on this and the fact that $\alpha_i^2=\beta_i^2$ for each $i\in\llbracket1,4\rrbracket$ we get that there exists a triple $(\gamma_3,\gamma_4,\gamma_5)\in K\times K^{\ast}\times K$ such that for $f_1(x):=\gamma_1x+\gamma_3x^2+\gamma_4x^4$ we have $f_1(\alpha_i)=\gamma_5\beta_i$ for each $i\in\llbracket1,4\rrbracket$.

If $\gamma_5=0$, then $f_1(x)$ is a multiple of $\prod_{i=1}^4 (x-\alpha_i)=x(x-1)(x-\alpha)(x-\beta)$ by an element in $K^{\ast}$ and we get that $\sum_{i=1}^4\alpha_i=1+\alpha+\beta=0$, a contradiction. Thus we can assume that $\gamma_5\neq 0$. 

For the two polynomials $f_2(x):=x$ and $g_2(x):=\gamma_1^{-1}\gamma_5x-\gamma_1^{-1}\gamma_3x^2-\gamma_1^{-1}\gamma_4x^4$ in $K[x]$ we have $\alpha_i=f_2(\alpha_i)=g_2(\beta_i)$ for each $i\in\llbracket1,4\rrbracket$ and $\deg(g_2)=4$. Thus $$g_2(x)=h_3(x)-\gamma_1^{-1}\gamma_4\prod_{i=1}^4 (x-\beta_i),$$ 
where $h_3(x)\in K[x]$ is the Lagrange interpolation polynomial of degree at most $3$ defined by $\alpha_i=h_3(\beta_i)$ for each $i\in\llbracket1,4\rrbracket$. As $\deg(h_3)=3$ by Lemma \ref{F12}(2) and (4) (or a direct computation) and $\sum_{i=1}^4 \beta_i=0$, the coefficient of $x^3$ in $g_2$ is equal to the one in $h_3$ and thus it is non-zero, a contradiction. 
So part (2) holds.

For part (3), for $n\ge 2$ we have $4\le\pi_{n,4}(\mathbb F_4)$ by Theorem \ref{T15}(3) and for $n\ge 3$ we have $\pi_{n,4}(\mathbb F_4)\le 6$ by Theorem \ref{T7}(3). Thus part (3) holds if $n\ge 3$.

From Proposition \ref{PR16.5}(6), Lemma \ref{F11.5}(4), and Lemma \ref{L26}(1) and (2) we get that $\pi_{2,4}(\mathbb F_4)\le 6$. Thus $\pi_{2,4}(\mathbb F_4)\in\{4,5,6\}$. As $\pi_{2,4}(\mathbb F_4)\neq 5$ by Corollary \ref{C11}, we get that $\pi_{2,4}(\mathbb F_4)\in\{4,6\}$. 

We show that the assumption $\pi_{2,4}(\mathbb F_4)=4$ leads to a contradiction. This assumption gives that for each $(\underline{P},\underline{Q})=\bigl((P_1,P_2,P_3,P_4),(Q_1,Q_2,Q_3,Q_4)\bigr)\in\mathbb D_{2,4}(\mathbb F_4)^2$ there exists a permutation $\sigma\in\Alt(\mathbb F_4^2)$ (see Theorem \ref{T6}(3)) such that $\sigma(\underline{P})=\underline{Q}$ and $\ell_{2,\mathbb F_4}(\sigma)\le 4$. Let $a\in\GA_2(K)$ be a minimal representation of $\sigma$ such that $\ell(a)=\ell_{2,\mathbb F_4}(\sigma)$ by Proposition \ref{PR14}(1.b). Thus $\ell(a)\le 4$ and hence from Proposition \ref{PR11}(4.a) we get that either $\j_a=2$ and we can write $a=bc$ with $(b,c)\in\GA_2(\mathbb F_4)[2]^2$ and $j_b=\j_c=1$ or $\j_a=1$. 

We now take $\underline{P}\in\mathbb D_{2,4}(\mathbb F_4)^{\d=1}\subset\mathbb D_{2,4}(\mathbb F_4)^{\sum=0}$ and $\underline{Q}\in\mathbb D_{2,4}(\mathbb F_4)^{\sum\neq 0}$ such that for $Y_2:=\{Q_1,Q_2,Q_3,Q_4\}$ we have $\t_{Y_2}=2$. We first show that we can assume that $a$ belongs to the subgroup $\Gamma_2(\mathbb F_4)[2]$ of $\GA_2(\mathbb F_4)$ generated by $\GA_2(K)[2]$. 

Clearly, if $\j_a=2$, then $(b,c)\in\Gamma_2(\mathbb F_4)[2]^2$ and thus $a\in\Gamma_2(\mathbb F_4)[2]$.

If $\j_a=1$, then there exists $f\in K[x]$ and a pair $(c_1,c_2)\in\AGL_2(K)^2\subset\Gamma_2(\mathbb F_4)[2]^2$ such that $a=c_2^{-1}\e\bigl(x_1+x_2+f(x_1)\bigr)c_1$. Hence $\e\bigl(x_1,x_2+f(x_1)\bigr)\bigl(c_1(\underline{P})\bigr)=c_2(\underline{Q})$. Not to introduce extra notation, based on Lemma \ref{F13}(1) and the fact that tame numbers are invariant under affine automorphisms, we can assume that we have $c_1=c_2=1_{\mathbb A^2_{\mathbb F_4}}$. Thus $\e\bigl(x_1,x_2+f(x_1)\bigr)(\underline{P})=\underline{Q}$. We write $P_i=(\alpha_i,\beta_i)\in\mathbb F_4$ for $i\in\llbracket1,4\rrbracket$; so for $i\in\llbracket1,4\rrbracket$ we have $Q_i=\bigl(\alpha_i,\beta_i+f(\alpha_i)\bigr)$. As $\t_{Y_2}=2$, for the subset $W:=\{\alpha_i|i\in\llbracket1,4\rrbracket\}$ of $\mathbb F_4$ we have $|W|\neq 4$ and thus, as $\sum_{i=1}^4 \alpha_i=0$, we have $|W|\in\{1,2\}$. Hence we can assume that $\deg(f)\le 1$. So $a\in\AGL_2(K)\Gamma_2(\mathbb F_4)[2]$.

As $a\in \Gamma_2(\mathbb F_4)[2]$ and $a(\underline{P})=\underline{Q}$, we reached a contradiction to Lemma \ref{F13}(1). We conclude that $\pi_{2,4}(\mathbb F_4)\neq 4$. Hence $\pi_{2,4}(\mathbb F_4)=6$. So part (3) holds.

For part (4), $\AGL_n(\mathbb F_2)$ acts transitively on $\mathbb D^{\d=3}_{n,4}(\mathbb F_2)$ and $\mathbb D^{\d=2}_{n,4}(\mathbb F_2)$ by Lemma \ref{F7}(1) and (3.c); also, $\e(x_1,x_2,x_3+x_1^2+x_1x_2+x_2^2,x_4,\ldots,x_n)\in\TGA_n(\mathbb F_2)[2]$ maps $\mathbb F_2^2\times\{0\}$ onto an affinely independent set in $\mathbb F_2^n$. So $\pi_{n,4}(\mathbb F_2)=2$ and part (4) holds.

For part (5), we have $\pi_{2,4}(\mathbb F_3)\le 4$ by Theorem \ref{T8}(4). For $n\ge 3$ we have $\pi^{\S}_{n,4}(\mathbb F_3)\le 4$ by Theorem \ref{T7}(2). For $n\ge 2$ we have $2\le\pi_{n,4}(\mathbb F_3)$ by either part (1) or Lemma \ref{L11}. As $\pi_{2,4}(\mathbb F_3)\neq 3$ by Corollary \ref{C11}, we get that $\pi_{2,4}(\mathbb F_3)\in\{2,4\}$. We are left to show that the assumption that $\pi_{2,4}(\mathbb F_3)\neq 4$ leads to a contradiction. 

Let the pair $(\underline{P},\underline{Q})=\bigl((P_1,P_2,P_3,P_4),(Q_1,Q_2,Q_3,Q_4)\bigr)\in\mathbb D_{2,4}(\mathbb F_3)$ be such that $P_1=Q_1=(0,0)$, $P_2=Q_2=(1,0)$, $P_3=Q_3=(0,1)$, $P_4=(1,1)$, and $Q_4=(1,2)$. The assumption implies that $\pi_{2,4}(\mathbb F_3)=2$. So there exists $a\in\GA_2(\mathbb F_3)[2]$ with $a(\underline{P})=\underline{Q}$. As $P_1+P_4=P_2+P_3$ but $Q_1+Q_4\neq Q_2+Q_3$, we have $a\notin\AGL_2(\mathbb F_3)$. Thus $\ell(a)=2$. From this and Proposition \ref{PR11}(1) and (4.a) we get that $\j_a=1$, i.e., $a$ is non-affine triangular. Therefore we can write $a=c\e\bigl(x_1,x_2+f(x_1)\bigr)b$ with $f\in\mathbb F_3[x_1]$ of degree $2$ and $(c,b)\in\GL_2(\mathbb F_3)\times\AGL_2(\mathbb F_3)$. 

We write $b(\underline{P})=(O_1,O_2,O_3,O_4)$ and $c^{-1}(\underline{Q})=(Q_5,Q_6,Q_7,Q_8)$. We get that the sets $\{O_1,O_2,O_3\}$ and $\{Q_5,Q_6,Q_7\}$ are linearly independent, we have two identities $O_1+O_4=O_2+O_3$ and $Q_6+2Q_7=2Q_5+Q_8$ and a non-identity $Q_5+Q_8\neq Q_6+Q_7$. We write $O_i=(\alpha_i,\beta_i)$ for $i\in \llbracket1,4\rrbracket$ and $Q_j=(\alpha_j,\beta_j)$ for $j\in \llbracket5,8\rrbracket$. Let $Z:=\{\alpha_i|i\in \llbracket1,4\rrbracket\}\subset\mathbb F_3$. As $\{O_1,O_2,O_3\}$ is linearly independent, we have $|Z|\ge 2$ and hence $|Z|\in\{2,3\}$. We consider two cases as follows.

{\bf Case 1: $|Z|=2$.} As $O_1+O_4=O_2+O_3$ and $2$ is invertible in $\mathbb F_3$, up to a reindexing of the set $\llbracket1,4\rrbracket$, we can assume that $\alpha_1=\alpha_2$ and $\alpha_3=\alpha_4$. For $j\in \llbracket5,8\rrbracket$ we have $\alpha_j=\alpha_{j-4}$ and $\beta_j=\alpha_{j-4}+f(\alpha_{j-4})$. Hence $Q_5+Q_8=Q_6+Q_7$, a contradiction.

{\bf Case 2: $|Z|=3$.} Let $\gamma\in\{0,2\}$. If $\gamma=2$ let $\delta\in\{0,1\}$ and if $\gamma=0$ let $\delta\in\{1,2\}$; so $\delta\neq\gamma$. As $O_1+O_4=O_2+O_3$, we can assume that $(\alpha_1,\alpha_2,\alpha_3,\alpha_4)$ is either $(\delta,\gamma,\gamma,2\gamma-\delta)$ or $(\gamma,2\gamma-\delta,\delta,\gamma)$. As for $j\in \llbracket5,8\rrbracket$ we have $\alpha_j=\alpha_{j-4}$, we compute that $\alpha_6+2\alpha_7-2\alpha_5-\alpha_8$ is $-2\gamma-\delta$ or $2\gamma+\delta$ (respectively), and hence it is non-zero, a contradiction to $Q_6+2Q_7=2Q_5+Q_8$.

From the two cases we conclude that part (5) holds.

For part (6), the lower bounds follow from part (4) and the fact that the sequence $\bigl(\pi_{3,m}(\mathbb F_2)\bigr)_{m\in \llbracket4,7\rrbracket}$ is non-decreasing. If $m\in\{5,6\}$ (resp.\ $m\in\{7,8\}$), then to prove that $\pi_{3,6}(\mathbb F_2)\le 3$ (resp.\ $\pi_{3,7}(\mathbb F_2)\le 6$), based on Proposition \ref{PR25}(1) (resp.\ \ref{L20}(4)) it suffices to show that each $\sigma\in\perm(\mathbb F_2^3)$ which is a product of at most $2$ transpositions (resp.\ each $\sigma\in\Alt(\mathbb F_2^3)$ with $\n(\sigma)\le 4$) is $a(\mathbb F_2)$ for an automorphism $a\in\STGA_3(\mathbb F_2)[3]$. We can assume that $\sigma$ is not the identity permutation. So $\sigma$ is a transposition or a $3$-cycle or a product of two disjoint transpositions whose support is affinely dependent or a product of two disjoint transpositions whose support is affinely independent; as $\AGL_3(\mathbb F_2)$ acts transitively on these types of permutations, the existence of $a\in\STGA_3(\mathbb F_2)[3]$ follows from Example \ref{EX7}. Thus part (6) holds.

The upper (resp.\ lower) bounds of the first sentence of part (7) follow from Theorem \ref{T7}(3) for $n\ge 4$ and Theorem \ref{T8}(4) if $n\in\{2,3\}$ (resp.\ Theorem \ref{T15}(5) for $n=2$ applied to only $4$ points and Theorem \ref{T15}(2) for $\n\ge 3$). The second sentence of part (7) follows from the first one and Corollary \ref{C11}.

For part (8) we have $|K|\ge 11$ and thus $\pi_{2,5}(K)\le 16$ by Theorem \ref{T7}(1). As $\pi_{2,5}(K)\ge 16$ by Theorem \ref{T15}(5), we get that part (8) holds. 

The upper (resp.\ lower) bounds of part (9) follow from Theorem \ref{T8}(2) for $n\ge 3$ and Theorem \ref{T8}(1.b) and (3) for $n=2$ (resp.\ Theorem \ref{T15}(1) for $\n\ge 3$ and Theorem \ref{T15}(5) for $n=2$ as $\mathbb F_7$ has a $6$-th primitive root of unity).

For part (10), the upper bounds of the first sentence follow from Theorem \ref{T8}(4) for $n\ge 3$ and Theorem \ref{T8}(3) for $n=2$ (resp.\ Theorem \ref{T15}(1) for $\n\ge 3$, the inequality $\pi_{2,5}(K)\le\pi_{2,6}(K)$ for $n=2$, and Theorem \ref{T15}(5) for $n=2$ and $|K|\equiv 1\pmod{6}$); the last sentence follows from the first one and Corollary \ref{C11}.

The upper (resp.\ lower) bounds of part (11) follow from Theorem \ref{T7}(3) for $n\ge 6$, Theorem \ref{T8}(2) for $n\in\{3,4,5\}$, and  Theorem \ref{T8}(1.b) for $n=2$ (resp.\ from $\pi_{2,6}(\mathbb F_7)\le\pi_{2,7}(\mathbb F_7)$ for $n=2$ and Theorem \ref{T15}(2) for $\n\ge 3$).\end{proof} 

\section{Basic properties of the $\upsilon_{n,m}(K)$s}\label{S29}

This section presents basic results on the $\upsilon^{\S}_{n,m}(K)$s and $\upsilon_{n,m}(K)$s introduced in Corollary \ref{C1}(2) and (3).

Throughout this and the next section let $n\in\mathbb N^{\ast}\setminus\{1\}$, $K$ a field, and $m\in\mathbb N$ be such that $\sqrt[n]{m}\le |K|$. We have $\upsilon_{n,m}(K)\le\pi_{n,m}(K)$ and $\upsilon^{\S}_{n,m}(K)\le\pi^{\S}_{n,m}(K)$. So each upper bound for $\pi_{n,m}(K)$ is also an upper bound for $\upsilon_{n,m}(K)$ and each upper bound for $\pi^{\S}_{n,m}(K)$ is also an upper bound for $\upsilon^{\S}_{n,m}(K)$. Based on this, we concentrate more on formulas of, identities and inequalities between, and lower bounds for the $\upsilon_{n,m}(K)$s and $\upsilon^{\S}_{n,m}(K)$s.

\begin{lemma}\label{F14}
The following properties hold.

\medskip
{\bf (1)} For a permutation $\sigma:W\rightarrow W$ of a set $W$ and subsets $Y$ and $Z$ of $W$ we have $\sigma(Y)=Z$ iff $\sigma(W\setminus Y)=W\setminus Z$. 

\smallskip
{\bf (2)} If $K$ is a finite field, then for each $m\in \llbracket0,|K|^n\rrbracket$ we have the symmetric properties $\upsilon_{n,m}(K)=\upsilon_{n,|K|^n-m}(K)$ and $\upsilon^{\S}_{n,m}(K)=\upsilon^{\S}_{n,|K|^n-m}(K)$.
\end{lemma}

\begin{proof}
Part (1) is clear and part (2) follows directly from part (1) applied to permutations $\sigma$ in $\Perm(K^n)$ of $\Alt(K^n)$.
\end{proof}

The following result is the analog of the statement that if a field $K$ and $m\in\mathbb N^{\ast}$ are such that $m\le |K|$, then $\pi_{1,m}(K)=1$ iff $m\in\{1,2\}$ and $\pi^{\S}_{1,m}(K)=1$ iff either $m=1$ or $m=|K|=2$.

\begin{lemma}\label{L27}
Let $K$ be a field. Let $m\in \llbracket1,|K|\rrbracket$. Let $\mathcal P_m(K)$ be the set of subsets of $K$ of cardinality $m$. Then the following properties hold.

\medskip
{\bf (1)} The group $\AGL_1(K)$ acts transitively on $\mathcal P_m(K)$ iff one of the following conditions holds.

\medskip\noindent
{\bf (1.a)} We have $m\in\{0,1,2\}$.

\smallskip\noindent
{\bf (1.b)} The field $K$ is finite with $|K|\ge 5$ and $m\in \{|K|-2,|K|-1,|K|\}$.

\smallskip\noindent
{\bf (1.c)} We have $K\cong\mathbb F_8$ and $m\in\{3,5\}$.

\medskip
{\bf (2)} The group $\AGL_1(K)$ acts simply transitively on $\mathcal P_m(K)$ iff the pair $(|K|,m)$ belongs to the set $\{(2,1),(8,3),(8,5)\}$.

\smallskip
{\bf (3)} The subgroup of $\AGL_1(K)$ formed by translations acts transitively on $\mathcal P_m(K)$ iff one of the following conditions holds.

\medskip\noindent
{\bf (3.a)} We have $m\in\{0,1\}$.

\smallskip\noindent
{\bf (3.b)} The field $K$ is finite and $m\in \{|K|-1,|K|\}$.

\medskip
{\bf (4)} The subgroup of $\AGL_1(K)$ formed by translations acts simply transitively on $\mathcal P_m(K)$ iff $m=1$ or $K$ is finite with $m=|K|-1$.
\end{lemma}

\begin{proof}
For the `if' of parts (1) and (2) and the `only if' of part (2), based on Lemma \ref{F14}(1) we can assume that either (1.a) holds or (1.c) holds and $m=3$. 

If (1.a) holds, then it is well-known that $\AGL_1(K)$ acts transitively on $\mathcal P_m(K)$. The action is simply transitive iff we have $\binom{|K|}{m}=|K|(|K|-1),$ which only holds when $(|K|,m)=(2,1)$. 

If (1.c) holds and $m=3$, then $|\AGL_1(\mathbb F_8)|=8\cdot 7=56=\binom{8}{3}$ and each non-identity element $a\in\AGL_1(\mathbb F_8)$ has either order $2$ or $7$, and it is easy to see that for every $Y\in \mathcal P_3(\mathbb F_8)$ we have $a(Y)\neq Y$. Hence $\AGL_1(\mathbb F_8)$ acts simply transitively on $\mathcal P_3(\mathbb F_8)$.

For the `only if' of part (1) we have to show that if conditions (1.a) to (1.c) do not hold, then $\AGL_1(K)$ does not act transitively on $\mathcal P_m(K)$. As the case $K$ is infinite is well-known, we can assume that $K$ is finite with $3\le m\le |K|-3$. Based on Lemma \ref{F14}(1) we can assume that $2m\le|K|$. We show that the assumption that $\AGL_1(K)$ acts transitively on $\mathcal P_m(K)$ leads to a contradiction. This assumption implies that $\binom{|K|}{m}$ divides $|AGL_1(K)|=|K|(|K|-1).$ As $\binom{|K|}{m}\geq \binom{|K|}{3}=\frac{|K|(|K|-1)(|K|-2)}{6},$ this implies $|K|\leq 8.$ From $|K|\ge 2m\ge 6$ and the fact that we are assuming that (1.c) does not hold we get that $(|K|,m)$ is either $(7,3)$ or $(8,4)$ which is impossible as $35\nmid 42$ and $70\nmid 56$.

Parts (3) and (4) follow easily from the fact that if $K$ is finite then we have $\binom{|K|}{m}>|K|$ if $m\in\llbracket2,m-2\rrbracket$.\end{proof}

We have the following variant of Lemma \ref{L11}.

\begin{lemma}\label{L28}
Let $n\in\mathbb N^\ast\setminus\{1\}$, $K$ a field, and $m\in\mathbb N$ with $\sqrt[n]{m}\le|K|$. Then the following statements are equivalent.

\medskip
{\bf (1)} We have $\upsilon^{\S}_{n,m}(K)=1$.

\smallskip
{\bf (2)} We have $\upsilon_{n,m}(K)=1$.

\smallskip
{\bf (3)} One of the following conditions holds.

\medskip\noindent
{\bf (3.a)} We have $m\in\{0,1,2\}$.

\smallskip\noindent
{\bf (3.b)} The field $K$ is finite and $m\in \{|K|^n-2,|K|^n-1,|K|^n\}$.

\smallskip\noindent
{\bf (3.c)} We have $K\cong\mathbb F_2$, $n\ge 3$, and $m\in\{3,|K|^n-3\}$.
\end{lemma}

\begin{proof}
Clearly, $(1)\Rightarrow (2)$. 

To prove that $(2)\Rightarrow (3)$, we can assume that $m\ge 3$ and that $\sqrt[n]{2m}\le |K|$ by Lemma \ref{F14}(2) in case $K$ is finite. It suffices to prove that if $|K|>2$ or $m>3$, then there exist subsets $Y$ and $Z$ of $K^n$ with $m$ elements such 
that $a(Y)\neq Z$ for each $a\in\AGL_n(K)$. We consider three disjoint cases as follows.

{\bf Case 1: $|K|=\infty$.} We can take $Y$ and $Z$ such that $\d_Y=1$ and $\d_Z=2$; clearly $a(Y)\neq Z$ for each $a\in\AGL_n(K)$.

{\bf Case 2: $|K|=2$.} Let $q\in \llbracket1,n-2\rrbracket$ be such that $m\in \llbracket2^q+1,2^{q+1}\rrbracket$. If $q=1$, then $m>3$ and hence $m=4=q+3$. If $q\ge 2$, then $m\ge 2^q+1\ge q+3$. So we can choose $Z$ such that $\d_Z=q+2$. Taking $Y$ such that $\d_Y=q+1$, we have $a(Y)\neq Z$ for each $a\in\AGL_n(K)$.

{\bf Case 3: $2<|K|<\infty$.} Let $q\in \llbracket0,n-1\rrbracket$ be such that $m\in \llbracket|K|^q+1,|K|^{q+1}\rrbracket$. As $m\ge 3$, we have $m\ge q+3$ if $q=0$. If $q\in \llbracket1,n-2\rrbracket$, then $m\ge |K|^q+1\ge q+3$ as $|K|\ge 3$. So for $q\le n-2$, by choosing $Z$ such that $\d_Z=q+2\le n$ and by taking $Y$ such that $\d_Y=q+1$, we have $a(Y)\neq Z$ for each $a\in\AGL_n(K)$. Thus we can assume that $q=n-1$. Hence $m\in \llbracket|K|^{n-1}+1,\lfloor\frac{|K|^n}{2}\rfloor\rrbracket$. 

To show the existence of $Y$ and $Z$ it suffices to show that the number $\binom{|K|^n}{m}$ of subsets of $K^n$ of cardinality $m$ does not divide $|AGL_n(K)|$. As 
$$|AGL_n(K)|=|K|^n\prod_{i=0}^{n-1} (|K|^n-|K|^i)<|K|^{n(n+1)},$$
to prove this, recall that $\binom{|K|^n}{s}$ is an increasing function on $s\in \llbracket|K|^{n-1},\lfloor\frac{|K|^n}{2}\rfloor\rrbracket$. So if we have $\binom{|K|^n}{|K|^{n-1}}\ge |K|^{n(n+1)}$, then we also have $\binom{|K|^n}{m}>|\AGL_n(K)|$. Using the standard inequality $\binom{|K|^n}{|K|^{n-1}}\ge \Bigl(\frac{|K|^n}{|K|^{n-1}}\Bigr)^{|K|^{n-1}}=|K|^{|K|^{n-1}}$ and the fact that $|K|^{n-1}>n(n+1)$ except when $(n,|K|)\in\{(2,3),(2,4),(2,5),(3,3)\}$, it follows that for $(n,|K|)\notin\{(2,3),(2,4),(2,5),(3,3)\}$ we have $\binom{|K|^n}{m}>|\AGL_n(K)|$. 

If $(n,|K|)$ is $(2,3)$ (resp.\ $(2,4)$, $(2,5)$, or $(3,3)$), then $\binom{|K|^n}{|K|^{n-1}+1}$ is $126$ (resp.\ $4368$, $177100$, or $8436285$), $|\AGL_n(K)|$ is $432$ (resp.\ $2880$, $12000$, or $303264$), and $\binom{|K|^n}{|K|^{n-1}+1}>|\AGL_n(K)|$ except when $(n,|K|)=(2,3)$. If $(n,|K|)=(2,3)$, then $m=4$ and $\binom{|K|^n}{m}=126$ does not divide $432=|\AGL_2(\mathbb F_3)|$. 

Thus $(2)\Rightarrow (3)$.

The implication $(3)\Rightarrow (1)$ is clear. For instance, if condition (3.c) holds, then we can assume that $m=3$ by Lemma \ref{F14}, hence $\pi^{\S}_{n,m}(K)=1$ by Lemma \ref{L11} and therefore $\upsilon^{\S}_{n,m}(K)=1$. 
\end{proof}

As for a finite field $K$ we have $\upsilon^{\S}_{n,|K|^n-i}(K)=1$ for $i\in\{0,1,2\}$ by Lemma \ref{L28}, in general the inequalities $\upsilon_{n,m}(K)\le\pi_{n,m}(K)$ and $\upsilon^{\S}_{n,m}(K)\le\pi^{\S}_{n,m}(K)$ are strict. This can happen even if $6\le 2m\le|K|^n$ as one can see based on the following result and Theorem \ref{T19}(1) (e.g., $\upsilon_{n,3}(K)=2<4=\pi_{n,3}(K)$ if $n\ge 2$ and $|K|\in\{4,5,8\}$).

\begin{proposition}\label{PR31}
If $K\in\{\mathbb F_3,\mathbb F_4,\mathbb F_5,\mathbb F_8\}$, then $\upsilon_{n,3}(K)=\upsilon^{\S}_{n,3}(K)=2$ for each $n\in \mathbb N^{\ast}\setminus\{1\}$.
\end{proposition}

\begin{proof}
Let $Y=\{P_1,P_2,P_3\}$ and $Z=\{Q_1,Q_2,Q_3\}$ be subsets of $K^n$ with $3$ elements. If $Y$ is collinear and $Z$ is non-collinear, then there exists no $c\in\AGL_n(K)$ such that $c(Y)=Z$. Hence $2\le\upsilon_{n,3}(K)$. As $\upsilon_{n,3}(K)\le\upsilon^{\S}_{n,3}(K)$, to complete the proof it suffices to show that there exists $a\in\STGA_n(K)[2]$ such that $a(Y)=Z$.

If $Y$ and $Z$ are collinear, then up to special affine automorphisms we can assume that $Y\cup Z$ is contained in the zero locus $x_2=\cdots=x_n=0$; so there exists $(\alpha,\beta)\in K^{\ast}\times K$ such that for $a:=\e(\alpha x_1+\beta,\alpha^{-1}x_2,\ldots,x_n)\in\ASL_n(K)$ we have $a(Y)=Z$ by Lemma \ref{L27}(1).

Thus we can assume that $\d_Y=2$. As $\d_Z\le 2$, up to special affine automorphisms we can assume that the zero locus $x_3=\cdots=x_n=0$ contains $Y\cup Z$. Thus we can assume that $n=2$. The proof of Theorem \ref{T7}(2) shows the existence of $a\in\SGA_2(K)[2]$ such that $a(P_i)=Q_i$ for each $i\in\{1,2,3\}$.\end{proof}

We have the following consequence of several results above.

\begin{corollary}\label{C29}
Let $n\in\mathbb N^{\ast}\setminus\{1\}$, $K$ a finite field, and $m\in \llbracket3,\lfloor\frac{|K|^n}{2}\rfloor\rrbracket$. Then the following properties hold.

\medskip
{\bf (1)} If $\upsilon_{n,m}(K)=\pi_{n,m}(K)$, then $\upsilon_{n,m-1}(K)\le\upsilon_{n,m}(K)$.

\smallskip
{\bf (2)} If $\pi_{n,m}(K)=2$, then $\upsilon_{n,m-1}(K)\le\upsilon_{n,m}(K)=2$.

\smallskip
{\bf (3)} If $\upsilon^{\S}_{n,m}(K)=\pi^{\S}_{n,m}(K)$, then $\upsilon^{\S}_{n,m-1}(K)\le\upsilon^{\S}_{n,m}(K)$.

\smallskip
{\bf (4)} If $\pi^{\S}_{n,m}(K)=2$, then $\upsilon^{\S}_{n,m-1}(K)\le\upsilon^{\S}_{n,m}(K)=2$.

\smallskip
{\bf (5)} If $n\ge 3$, then $\upsilon_{n,4}(\mathbb F_2)=2$.

\smallskip
{\bf (6)} If $n\ge 4$, then $\upsilon_{n,5}(\mathbb F_2)=2$.
\end{corollary}

\begin{proof}
Part (1) holds as we have $\upsilon_{n,m-1}(K)\le\pi_{n,m-1}(K)\le\pi_{n,m}(K)=\upsilon_{n,m}(K)$.

For part (2), we have $(n,m,|K|)\neq (3,3,2)$ by Lemma \ref{L11}(2) and $|K|^n\ge 8$. From this and Lemma \ref{L28} we get that $1\neq\upsilon_{n,m}(K)$. Therefore we have relations $2\le\upsilon_{n,m}(K)\le\pi_{n,m}(K)=2$ which imply that $\upsilon_{n,m}(K)=\pi_{n,m}(K)=2$. From this and part (1) we get that part (2) holds.

Parts (3) and (4) are proved similarly to parts (1) and (2) (respectively).

As $\upsilon_{n,4}(\mathbb F_2)\le\pi_{n,4}(\mathbb F_2)$, part (5) follows from the identity $\pi_{n,4}(\mathbb F_2)=2$ of Theorem \ref{T19}(4) and the relation $\upsilon_{n,4}(\mathbb F_2)\neq 1$ by Lemma \ref{L28}.

For part (6), let $Y$ and $Z$ be two subsets of $\mathbb F_2^n$ with $5$ elements. We consider the subset $W:=\{O,O_1,O_2,O_3,O_1+O_2\}$ of $\mathbb F_2^n$, where $O:=(0,\ldots,0)$, $O_1:=(1,0,\ldots,0)$, $O_2:=(0,1,0,\ldots,0)$, and $O_3:=(0,0,1,0,\ldots,0)$. Taking $b:=\e(x_1,x_2,x_3,x_4+x_1x_2,x_5,\ldots,x_n)\in\TGA_n(\mathbb F_2)[2]$, the set $b(W)$ is linearly independent. From this and Lemma \ref{F7}(3.d), as in the proof of Proposition \ref{PR31} we argue with ``collinear'' and ``non-collinear'' replaced by `linearly dependent' and `linearly independent' (respectively) that there exists $c\in\TGA_n(\mathbb F_2)[2]$ with $c(Y)=Z$. So part (6) holds.
\end{proof}

To facilitate applications of prior sections, we note the following obvious fact.

\begin{fact}\label{F15}
The number $\upsilon_{n,m}(K)\in\mathbb N^{\ast}$ (resp.\ $\upsilon^{\S}_{n,m}(K)\in\mathbb N^{\ast}$) is the smallest with the property that for every pair $\bigl((P_1,\ldots,P_m),(Q_1,\ldots,Q_m)\bigr)\in\mathbb D_{n,m}(K)^2$ there exists $a\in\TGA_n(K)[\upsilon_{n,m}(K)]$ (resp.\ $a\in\STGA_n(K)[\upsilon^{\S}_{n,m}(K)]$) and $\sigma\in S_m$ such that $a(P_i)=Q_{\sigma(i)}$ for each $i\in \llbracket1,m\rrbracket$. 
\end{fact}

\section{Applications to the $\upsilon_{n,m}(K)$s}\label{S30}

In this section we apply the previous ones to study the $\upsilon_{n,m}(K)$s and $\upsilon^{\s}_{n,m}(K)$s. 

We have the following variant of Definitions \ref{D2.5} and \ref{D21}.

\begin{definition}\label{D24}
Let $n\in\mathbb N^{\ast}\setminus\{1\}$, $K$ a field, $s\in\mathbb N^{\ast}$, and a finite non-empty subset $Y$ of $K^n$ be such that $s+|Y|\le |K|^n$. 

\medskip
{\bf (1)} For a subgroup $H$ of $\GA_n(K)$, let $\Stab_H(Y)$ be the subgroup of $H$ that maps $Y$ onto $Y$ (thus $\Stab_H(Y)=H\cap\Stab_{\GA_n(K)}(Y)$). 

\smallskip
{\bf (2)} If $K$ is a finite field with $4\mid |K|$, then we assume that the inequality $\max(s,|K|^n-2-|Y|-s)\ge 2$ holds. Let $\upsilon_{Y,s}(K)\in\mathbb N^{\ast}$\index{$\upsilon_{Y,s}(K)$ invariant} (resp.\ $\upsilon^{\S}_{Y,s}(K)\in\mathbb N^{\ast}$\index{$\upsilon^{\S}_{Y,s}(K)$ invariant}) be the smallest such that for every two subsets $Y_0$ and $Z_0$ of $K^n\setminus Y$ with $s$ elements there exists $a\in\Stab_{\TGA_n(K)}(Y)[\upsilon_{Y,s}(K)]$ (resp.\ $a\in\Stab_{\STGA_n(K)}(Y)[\upsilon_{Y,s}(K)]$) with $a(Y_0)=Z_0$.

\smallskip
{\bf (3)} Let $m\in \llbracket1,|K|^n-s\rrbracket$. If $K$ is a finite field with $4\mid |K|$, then we assume that $\max(s,|K|^n-s-m)\ge 2$. Let 
$$\xi_{n,m,s}(K):=\max\bigl(\upsilon_{Z,s}(K)|Z\subset K^n, |Z|=m\bigr)\in \llbracket1,\pi_{n,m+s}(K)\rrbracket\index{$\xi_{n,m,s}(K)$ invariant}$$
and
$$\xi^{\S}_{n,m,s}(K):=\max\bigl(\upsilon^{\S}_{Z,s}(K)|Z\subset K^n, |Z|=m\bigr)\in \llbracket1,\pi^{\S}_{n,m+s}(K)\rrbracket.\index{$\xi^{\S}_{n,m,s}(K)$ invariant}$$
\end{definition}

Let $q\in\mathbb N^{\ast}$ with $\sqrt[n]{1+q}\le |K|$. If $|Y|=1$, then $\upsilon_{Y,q}(K)=\upsilon_{n,q+1}(K)$ and $\upsilon^{\S}_{Y,q}(K)=\upsilon^{\S}_{n,q+1}(K)$, and thus $\xi_{n,1,s}(K)=\Xi_{n,1,q}(K)$ and $\xi^{\S}_{n,m,q}(K)=\Xi^{\S}_{n,1,q}(K)$.

We have the following variant of Proposition \ref{PR26}(1) to (3).

\begin{proposition}\label{PR32}
Let $K$ be a finite field. Let $(n,m,s)\in (\mathbb N^{\ast}\setminus\{1\})^2\times\mathbb N^{\ast}$ be such that $m+s\le |K|^n$. If $4\mid |K|$, then we assume that $\max(s,|K|^n-2-m-s)\ge 2$. Let $Y$ be a subset of $K^n$ with $m$ elements. Then the following properties hold.

\medskip
{\bf (1)} We have inequalities
$$\upsilon_{Y,s}(K)\le\min\bigl(\pi_{Y,s}(K),\xi_{n,m,s}(K)\bigr)\le\max\bigl(\pi_{Y,s}(K),\xi_{n,m,s}(K)\bigr)\le\Xi_{n,m,s}(K).$$

{\bf (2)} We have an inequality $\upsilon_{n,m+s}(K)\le \upsilon_{n,m}(K)\xi_{n,m,s}(K)$.

\smallskip
{\bf (3)} If $s\ge 2$ we assume that moreover $m+2s\le |K|^n+2$. Then we have an inequality $\upsilon_{n,s}(K)\le\xi_{n,m,s}(K)$.

\smallskip
{\bf (4)} Parts (1) to (3) also hold for $(\upsilon,\xi,\pi)$ replaced by $(\upsilon^{\S},\xi^{\S},\pi^{\S})$.
\end{proposition}

\begin{proof}
Part (1) follows directly from the definitions. 

To prove part (2), let $Y_1=\{P_i|i\in \llbracket1,m+s\rrbracket\}$ and $Z_1=\{Q_i|i\in \llbracket1,m+s\rrbracket\}$ be two subsets of $K^n$ with $m+s$ elements. We consider two subsets $Y_0\subset Y_1$ and $Z_0\subset Z_1$ with $m$ elements. Reindexing, we can assume that $Y_0=\{P_i|\i\in \llbracket1,m\rrbracket\}$ and $Z_0=\{Q_i|i\in \llbracket1,m\rrbracket\}$. 

First we choose $a\in\TGA_n(K)[\upsilon_{n,m}(K)]$ such that $a(Y_0)=Z_0$. Next we choose $b\in\TGA_n(K)[\upsilon_{Z_0,s}(K)]$ such that $b\in\Stab_{\TGA_n(K)}(Z_0)$ and $b\bigl(a(Y_1\setminus Y_0)\bigr)=Z_1\setminus Z_0$. For $c:=ba$ we have $c(Y_0)=Z_0$. From this and the inequalities 
$$\ell(c)\le\ell(b)\ell(a)\le \upsilon_{n,m}(K)\upsilon_{Z_0,s}(K)\le\upsilon_{n,m}(K)\xi_{n,m,s}(K)$$
we get that part (2) holds.

For part (3), as $\upsilon_{n,1}(K)=1$ we can assume that $s\ge 2$. We consider two subsets $Y_0$ and $Z_0$ of $K^n$ such that $|Y_0|=|Z_0|=s$ and for each $a\in\TGA_n(K)$ with $a(Y_0)=Z_0$ we have $\ell(a)\ge\upsilon_{n,s}(K)$. Up to a special affine automorphism we can assume that $|Y_0\cap Z_0|\ge 2$. Hence $|Y_0\cup Z_0|\le 2s-2$. From this and the inequality $m+2s\le |K|^n+2$ we get that there exists a subset $Y\subset K^n\setminus (Y_0\cup Z_0)$ with $|Y|=m$. For each $b\in\TGA_n(K)$ with $b(Y)=Y$ and $b(Y_0)=Z_0$ we have $\ell(b)\ge\upsilon_{n,s}(K)$, hence $\upsilon_{Y,s}(K)\ge\upsilon_{n,s}(K)$. From this and the inequality $\xi_{n,m,s}(K)\ge\upsilon_{Y,s}(K)$ we get that part (3) holds.

Part (4) is clear as the same proofs of parts (1) to (3) apply to $(\upsilon^{\S},\xi^{\S},\pi^{\S})$.
\end{proof}

We have the following variant of Definition \ref{D22}.

\begin{definition}\label{D25}
Let a field $K$ and a triple $(n,m,l)\in (\mathbb N^{\ast})^3$ be such that $n\ge 2$, $m\ge 3$, and $\sqrt[n]{m}\le |K|$. Let $\flat_{n,m,l}(K)\in \llbracket1,m\rrbracket$\index{$\flat_{n,m,l}(K)$ invariant} (resp.\ $\flat^{\S}_{n,m,l}(K)\in \llbracket1,m\rrbracket$\index{$\flat^{\S}_{n,m,l}(K)$ invariant}) be the largest such that for every two subsets $Y$ and $Z$ of $K^n$ with $m$ elements there exists $a\in\TGA_n(K)[l]$ (resp.\ $a\in\STGA_n(K)[l]$) such that $|a(Y)\cap Z|\ge\flat_{n,m,l}(K)$ (resp.\ $|a(Y)\cap Z|\ge\flat^{\S}_{n,m,l}(K)$). 
\end{definition}

The sequence $\bigl(\flat_{n,m,l}(K)\bigr)_{l\ge 1}$ is non-decreasing with $\flat_{n,m,\upsilon_{n,m}(K)+l}(K)=m$ for each $l\in\mathbb N$, if $\upsilon_{n,m}(K)>1$ then we have an inequality $\flat_{n,m,\upsilon_{n,m}(K)-1}(K)<m$, and clearly $\flat_{n,m,l}(K)\ge\natural_{n,m,l}(K)$. The same holds with $(\flat,\natural)$ replaced by $(\flat^{\S},\natural^{\s})$. 

Our applications of the $\flat_{n,m,l}(K)$s and $\flat^{\S}_{n,m,l}(K)$s are modeled on Propositions \ref{PR25} and \ref{PR27}.

\begin{proposition}\label{PR33}
Let $n\in\mathbb N^{\ast}\setminus\{1\}$ and $K$ a finite field. If $|K|=2$ let $s\in \llbracket1,n-1\rrbracket$ and if $|K|\ge 3$ let $s\in \llbracket1,n\rrbracket$. Let $m\in\bigl\llbracket|K|^{s-1}+1,|K|^s\bigr\rrbracket$ if $s\le n-1$ and $m\in\bigl\llbracket|K|^{s-1}+1,\lfloor\frac{|K|^s}{2}\rfloor\bigr\rrbracket$ if $s=n$ (thus $2m\le |K|^n$). Let $Y$ and $Z$ be subsets of $K^n$ with $m$ elements. For each $i\in\bigl\llbracket1,\lfloor\frac{|K|^n}{2}\rfloor-s-1\bigr\rrbracket$ let $\mathcal C_i$ be a conjugacy class of $\perm(K^n)$ with $\n(\mathcal C_i)=2i$ and $\c_{\textup{odd}}(\mathcal C_i)=0$. Let $\kappa$ be $2n-3$ if $|K|=2$, be $4n$ if $|K|=3$, and be $E_{n-2,|K|-1}$ of Equation (\ref{EQ13}) if $|K|\ge 4$. Then the following properties hold.

\medskip
{\bf (1)} We have $\flat_{n,m,1}(K)\ge s+1$. If moreover, $s\le n-1$, then $\flat^{\S}_{n,m,1}(K)\ge s+1$.

\smallskip
{\bf (2)} For each $l\in\mathbb N^{\ast}$ there exists $(a_l,\sigma_l)\in\TGA_n(K)[l]\times\perm(K^n)$ (resp.\ $(a_l,\sigma_l)\in\STGA_n(K)[l]\times\perm(K^n)$) such that $\sigma_l\bigl(a_l(Y)\bigr)=Z$ and the cyclic type of $\sigma_l$ is arbitrary subject to the two restrictions $\n(\sigma_l)\in 2\llbracket0,m-\flat_{n,m,l}(K)\rrbracket$ (resp.\ $\n(\sigma_l)\in 2\llbracket0,m-\flat^{\S}_{n,m,l}(K)\rrbracket$) and $\c_{\textup{odd}}(\sigma_l)=0$. Also, if $\sigma_l$ is a transposition, then there exists a $3$-cycle $\theta_l\in\Alt(K^n)$ of non-collinear points such that $\theta_l\bigl(a_l(Y)\bigr)=Z$.

\smallskip
{\bf (3)} If $4\nmid |K|$, then for each $l\in \llbracket1,\upsilon_{n,m}(K)-1\rrbracket$ we have an inequality
$$\upsilon_{n,m}(K)\le l\max\bigl(\Omega(\mathcal C_i)|i\in \llbracket1,m-\flat_{n,m,l}(K)\rrbracket\bigr).$$
In particular, by taking $l=1$ we get the inequality 
$$\upsilon_{n,m}(K)\le \max\bigl(\Omega(\mathcal C_i)|i\in \llbracket1,m-s-1\rrbracket\bigr).$$
Moreover, in the last two sentences we can replace $\Omega(\mathcal C_1)=\kappa_{n,K}$ by $\kappa$.

\smallskip
{\bf (4)} If $\cup_{i=2}^{m-s-1} \mathcal C_i\subset\Alt(K^n)$, then for every $l\in \llbracket1,\upsilon^{\S}_{n,m}(K)-1\rrbracket$ we have an inequality $\upsilon^{\S}_{n,m}(K)\le l\max\bigl(\kappa,\Omega^{\S}(\mathcal C_i)|i\in \llbracket2,m-\flat^{\S}_{n,m,l}(K)\rrbracket\bigr)$.
In particular, by taking $l=1$ we get the inequality $\upsilon^{\S}_{n,m}(K)\le \max\bigl(\kappa,\Omega^{\S}(\mathcal C_i)|i\in \llbracket2,m-s-1\rrbracket\bigr).$\end{proposition}

\begin{proof}
Let $Y_0\subset Y$ and $Z_0\subset Z$ be linearly independent subsets with $s+1$ elements. If $a\in\AGL_n(K)$ is such that $a(Y_0)=Z_0$, then $a(Y)\cap Z\supset Z_0$; so $|a(Y)\cap Z|\ge s+1$. If $s\le n-1$, then we can assume that $a\in\ASL_n(K)$. Based on the last two sentences and very definitions we get that part (1) holds.

For part (2), let $a_l\in\TGA_n(K)[l]$ (resp.\ $a_l\in\STGA_n(K)[l]$) be such that we have $|a_l(Y)\cap Z|\ge\flat_{n,m,l}(K)]$ (resp.\ $|a_l(Y)\cap Z|\ge\flat^{\S}_{n,m,l}(K)]$). Thus $a_l(Y)\setminus Z$ and $Z\setminus a_l(Y)$ have the same cardinality $s_l\in \llbracket0,m-\flat_{n,m,l}(K)\rrbracket$ (resp.\ $s_l\in \llbracket0,m-\flat^{\S}_{n,m,l}(K)\rrbracket$). 

Let 
$$B_l:=\{\sigma\in\perm(K^n)|\sigma\bigl(a_l(Y)\setminus Z\bigr)=Z\setminus a_l(Y)\;\textup{and}\;\n(\sigma)=2s_l\}\subset\perm(K^n).$$
For every $\sigma\in B_l$, we have $\supp(\sigma)=[a_l(Y)\setminus Z]\cup [Z\setminus a_l(Y)]$, $\c_{\textup{odd}}(\sigma)=0$, $\n(\sigma)\in 2\llbracket0,m-\flat_{n,m,l}(K)\rrbracket$ (resp.\ $\n(\sigma)\in 2\llbracket0,m-\flat^{\S}_{n,m,l}(K)\rrbracket$), and $\sigma\bigl(a_l(Y)\bigr)=Z$. Thus we can take $\sigma_l\in B_l$ such that the first sentence of part (2) holds. We have $s_l=1$ iff $B_l$ consists of a single transposition and iff $B_l$ contains a transposition. 

If $s_l=1$ and $B_l=\{(P\;Q)\}$, then we can take $\theta=(P\;Q\; O)$ for an arbitrary point $O\in K^n\setminus (Y\cup Z\cup\lambda)$, where $\lambda$ is the line in $K^n$ that passes through $P$ and $Q$; such a point exists as $|K|^n>|K|+\bigl|\lfloor\frac{|K|^n}{2}\rfloor\bigr|-1$. So part (2) holds.

Let $\mathcal C_0$ be the conjugacy class in $\perm(K^n)$ of the identity element; we have $\Omega(\mathcal C_0)=1$. The first sentence of part (3) follows from part (2) by taking $\sigma_l$ in $B_l\cap \mathcal C_{s_l}$ as for every $b_l\in\TGA_n(K)[\Omega(\mathcal C_{s_l})]$ such that $b_l(K)=\sigma_l$, by defining $c_l:=b_la_l\in\TGA_n(K)$ we have $\ell(c_l)\le l\Omega(\mathcal C_l)$ and $c_l(Y)=Z$; thus 
$$\upsilon_{n,m}(K)\le l\Omega(\mathcal C_{s_l})\le l\max\bigl(\Omega(\mathcal C_i)|i\in \llbracket1,m-\flat_{n,m,l}(K)\rrbracket\bigr).$$ 
The second sentence of part (3) follows from the first sentence and part (1). The last sentence of part (3) is proved similarly but, in case $\sigma_l$ is a transposition, we assume that $b_l\in\STGA_n(K)[\kappa]$ is such that $b_l(K)=\theta_l$ by Example \ref{EX7} and $\pi_{2,2}(\mathbb F_2)=1$ if $|K|=2$ and by Proposition \ref{PR10}(2) and (3) if $|K|\ge 3$.

The proof of part (4) is the same as the proof of part (3) with $\mathcal C_1$ replaced by $\mathcal Y_3$ in all cases, i.e., if $s_l=1$, we use $\theta_l\in\mathcal Y_3^{\d=2}\subset\mathcal Y_3$ instead of $\sigma_l$.
\end{proof}

\begin{corollary}\label{C30} Let $K$ be a finite field with $|K|\ge 4$. Let $r$ and $r_0$ be as in Theorem \ref{T8}. Let $m\in\bigl\llbracket|K|+1,\bigl\lfloor\frac{|K|^2}{2}\bigr\rfloor\bigr\rrbracket$. Then 
the following properties hold.

\medskip
{\bf (1)} Suppose that $|K|$ is odd. Let $\epsilon\in\{5,10\}$ be such that it is $5$ iff $|K|$ is a prime. Let $\Omega:=\Omega^{\S}(\mathcal Y_{|K|})$. Then we have inequalities
$$\upsilon^{\S}_{2,m}(K)\le \max\Bigl(4(|K|-1)^4,\Omega^{\lfloor\frac{2m+3|K|-10}{|K|-1}\rfloor}\Bigr)\le [(|K|-r_0)^2(|K|-2)^2(|K|-1)^{\epsilon}]^{|K|+3}.$$
If $|K|\in\{5,7\}$, then we can replace the exponent $|K|+3$ by $|K|+2$.

{\bf (2)} Suppose that $p:=\char(K)$ is odd. Let $q\in\mathbb N^{\ast}$ be such that $|K|=p^q$. Let $\Omega:=\Omega^{\S}(p^{q-1}\mathcal Y_p)$. Then we have inequalities
$$\upsilon^{\S}_{2,m}(K)\le \max\Bigl(4(|K|-1)^4,\Omega^{\nu^+_{p,p^{q-1}}(2m-6-\cycle)}\Bigr)$$
$$\le [(|K|-r_0)^2(|K|-2)^2(|K|-1)^5]^{\nu^+_{p,p^{q-1}}(|K|^2-7-\cycle)}.$$

{\bf (3)} Suppose that $|K|=2^q$ with $q\in\mathbb N^{\ast}\setminus\{1\}$. Let $\Omega:=\Omega^{\S}(2^{q-1}\mathcal Y_2)$. Then we have inequalities
$$\upsilon_{2,m}(K)\le \max\Bigl(4(|K|-1)^4,\Omega^{\lfloor\frac{m-3}{2^{q-1}}\rfloor+2}\Bigr)\le [(2^q-r_0)^2(2^q-2)^2(2^q-1)^5]^{2^q+1}.$$
\end{corollary}

\begin{proof}
We have $m\ge 5$. We apply Proposition \ref{PR33} to $n=s=2$ and $l=1$. For parts (1) and (2), for $i\in \llbracket2,m-3\rrbracket$ we take $\mathcal C_i$ to be the conjugacy class in $\perm(K^2)$ formed by products $\sigma_i=\tau\theta_i$ with $\tau$ a transposition, $\theta_i$ an $(2i-2)$-cycle, and $\supp(\tau)\cap\supp(\theta_i)=\emptyset$. We have $\n(\sigma_i)=2i\le 2(m-3)$ and $\c(\sigma_i)=2$. Taking $\mathcal C_1:=\mathcal Y_3$, we have $\Omega^{\S}(\mathcal C_1)\le 4(|K|-1)^4$ by Proposition \ref{PR10}(3) and Example \ref{EX8} and this explains the $4(|K|-1)^4$ terms in the maxima.

For part (1) we use $|K|$-cycles. We have $\nu_{|K|,1}(\sigma_i)\le\bigl\lfloor\frac{2i+3|K|-4}{|K|-1}\bigr\rfloor\le\bigl\lfloor\frac{2m+3|K|-10}{|K|-1}\bigr\rfloor$ by Theorem \ref{P4}(1). As for $i\in \llbracket2,m-3\rrbracket$ we have $\mathcal C_i\subset\Alt(K^2)$, it follows that $\Omega^{\S}(\mathcal C_i)\le\Omega^{\lfloor\frac{2m+3|K|-10}{|K|-1}\rfloor}$. From this, the inequality $\Omega\le (|K|-r_0)^2(|K|-2)^2(K|-1)^{\epsilon}$ (see Corollary \ref{C18}(1)), the inequalities $\bigl\lfloor\frac{2m+3|K|-10}{|K|-1}\bigr\rfloor\le\bigl\lfloor\frac{|K|^2+3|K|-11}{|K|-1}\bigr\rfloor\le|K|+3$ if $|K|\ge 9$ and $\bigl\lfloor\frac{|K|^2+3|K|-11}{|K|-1}\bigr\rfloor\le|K|+2$ if $|K|\in\{5,7\}$, and Proposition \ref{PR33}(3), we get that part (1) holds. 

For part (2) we use $p^{q-1}\mathcal Y_{p}$. We have inequalities $2(m-3)\le |K|^2-7$ and $\Omega^{\S}(\mathcal C_i)\le\Omega^{\nu^+_{p,p^{q-1}}(2m-6-\cycle)}\le\Omega^{\nu^+_{p,p^{q-1}}(|K|^2-7-\cycle)}$. From these, Proposition \ref{PR33}(3), and by Corollary \ref{C18}(1) the inequality $\Omega\le (|K|-r_0)^2(|K|-2)^2(K|-1)^5$, we get that part (2) holds. 

For part (3) we use products of $2^{q-1}$ disjoint transpositions and for $i\in \llbracket2,m-3\rrbracket$ even (resp.\ odd) we take $\mathcal C_i$ to be $i\mathcal Y_2$ (resp.\ to be formed by products $\sigma_i=\vartheta\vartheta_i$ with $\vartheta\in\mathcal Y_4$, $\vartheta_i\in (i-2)\mathcal Y_2$, and $\supp(\vartheta)\cap\supp(\vartheta_i)=\emptyset$); clearly, $\mathcal C_i\subset\Alt(K^2)$. If $i$ is even (resp.\ odd), then each $\sigma\in\mathcal C_i$ is a product $\sigma=\prod_{j=1}^{\nu_i} \theta_j$ with each $\theta_j\in 2\mathcal Y_{2^{q-1}}$ and $\nu_i:=\lfloor\frac{i}{2^{q-1}}\rfloor$ if $2^{q-1}\mid i$ and $\nu_i:=\lfloor\frac{i}{2^{q-1}}\rfloor+2$ if $2^{q-1}\nmid i$ (resp.\ $\nu_i:=\lfloor\frac{i-1}{2^{q-1}}\rfloor+2$) as one can see by an argument similar to the one used in the proof of Lemma \ref{F2}. From these, Proposition \ref{PR33}(4), and $\Omega\le (2^q-r_0)^2(2^q-2)^2(2^q-1)^5$ by Corollary \ref{C18}(2), we get that part (2) holds. 
\end{proof}

We have the following variant of Proposition \ref{PR26}.

\begin{proposition}\label{PR34}
Let $(n,m,s)\in\mathbb (N^{\ast} \setminus\{1\})^2\times\mathbb N^{\ast}$ and $K$ a finite field be such that $m\in \llbracket s,|K|^n\rrbracket$. For $i\in \llbracket2,s\rrbracket$ let $\mathcal C_i$ be a conjugacy class of $\perm(K^n)$ with $\n(\mathcal C_i)=2i$ and $\c_{\textup{odd}}(\mathcal C_i)=0$. If $4\mid |K|$ we assume that $\cup_{i=2}^s \mathcal C_i\subset\Alt(K^n)$ and define $\mathcal C_1:=\mathcal Y_3$. If $4\nmid |K|$ let $\mathcal C_1\in\{\mathcal Y_2,\mathcal Y_3\}$. Then the following properties hold.

\medskip
{\bf (1)} Let $\Omega_s:=\max\bigl(\Omega(\mathcal C_i)|i\in \llbracket1,s\rrbracket\bigr)$. Then we have inequalities
$$\frac{1}{\Omega_s}\le \frac{\upsilon_{n,m}(K)}{\upsilon_{n,m-s}(K)}\le\Omega_s.$$

{\bf (2)} Suppose that $\cup_{i=2}^s \mathcal C_i\subset\Alt(K^n)$ and $\mathcal C_1=\mathcal Y_3$. Let 
$$\Omega^{\S}_s:=\max\bigl(\Omega^{\S}(\mathcal C_i)|i\in \llbracket1,s\rrbracket\bigr).$$ Then we have inequalities
$$\frac{1}{\Omega^{\S}_s}\le \frac{\upsilon^{\S}_{n,m}(K)}{\upsilon^{\S}_{n,m-s}(K)}\le\Omega^{\S}_s.$$
\end{proposition}

\begin{proof}
To prove that $\upsilon_{n,m}(K)\le\Omega_s\upsilon_{n,m-s}(K)$, we consider two subsets $Y$ and $Z$ of $K^n$ with $m$ elements. Let $Y_0\subset Y$ and $Z_0\subset Z$ be two subsets with $m-s$ elements. Based on Lemma \ref{L28}, if $s=1$ we can assume that $m\le |K|^n-2$. 

Let $a\in\TGA_n(K)[\upsilon_{n,m-s}(K)]$ be such that $a(Y_0)=Z_0$. We can also assume that $a(Y)\neq Z$. Then there exists $\sigma\in\cup_{i=1}^s \mathcal C_i$ such that $\sigma\bigl(a(Y)\bigr)=Z$. We detail on the only case in which $\supp(\sigma)$ must intersect non-trivially $a(Y)\cap Z$. If $s=1$, $|a(Y)\cap Z|=m-1\le |K|^n-3$, and $\mathcal C_1=\mathcal Y_3$, then for the only points $P\in\ a(Y)\setminus Z$ and $Q\in Z\setminus a(Y)$ we can take $\sigma=(Q\;P\;O)$ with $O\in K^n\setminus \bigl(a(Y)\cup Z\bigr)$ arbitrary.

For $b\in\TGA_n(K)[\Omega_s]$ with $b(K)=\sigma(K)$ let $c:=ba\in\TGA_n(K)[\Omega_s\upsilon_{n,m-s}(K)]$. We have $c(Y)=Z$. Thus $\upsilon_{n,m}(K)\le\Omega_s\upsilon_{n,m-s}(K)$.

The inequality $\upsilon_{n,m-s}(K)\le\Omega_s\upsilon_{n,m}(K)$ follows from the inequality of the previous paragraph applied to $m$ equal to $|K|^n-m+s\in \llbracket s,|K|^n\rrbracket$ and Lemma \ref{F14}(2). So part (1) holds.

Part (2) is proved entirely in the same manner as part (1).\end{proof}

Based on Fact \ref{F15}, we have the following permutation variant of Definition \ref{D23}.

\begin{definition}\label{D26}
Let $m\in \mathbb N^{\ast}\setminus\{1\}$ and $s\in \llbracket1,m-1\rrbracket$. 

\medskip
{\bf (1)} Let $\mathfrak g(m,s)\in\mathbb N^{\ast}\setminus\{1\}$ be the smallest with the property that for each field $K$ with $|K|\ge\mathfrak g(m,s)$, there exists $\bigl((\alpha_1,\ldots,\alpha_m),(\beta_1,\ldots,\beta_m)\bigr)\in\mathbb D_{1,m}(K)^2$ such that for every permutation $\sigma\in S_m$, the pair $\bigl((\alpha_1,\ldots,\alpha_m),(\beta_{\sigma(1)},\ldots,\beta_{\sigma(m)})\bigr)$ in $\mathbb D_{1,m}(K)^2$ is $\llbracket1,s\rrbracket$-generic.

\smallskip
{\bf (2)} Suppose that $m\ge 4$ and $s\in \llbracket2,m-1\rrbracket$. Let $\mathfrak h(m,s)\in\mathbb N^{\ast}\setminus\{1\}$ be the smallest with the property that for each field $K$ with $|K|\ge\mathfrak h(m,s)$, there exists a pair $\bigl((\alpha_1,\ldots,\alpha_m),(\beta_1,\ldots,\beta_m)\bigr)\in\mathbb D_{1,m}(K)^2$ such that for every permutation $\sigma\in S_m$, each polynomial $h_{\sigma}\in K[x]$ with the property that $h_{\sigma}(\alpha_i)=\beta_{\sigma(i)}$ for every $i\in \llbracket1,m\rrbracket$ has degree at least $s$.
\end{definition}

The next two results justify Definition \ref{D26}. First, we have the following variant of Proposition \ref{PR28}.

\begin{proposition}\label{PR35}
Let $m\in \mathbb N^{\ast}\setminus\{1\}$ and $s\in \llbracket1,m-1\rrbracket$. Then the following properties hold.

\medskip
{\bf (1)} We have an inequality $\mathfrak g(m,s)\le m+m!\frac{s(s+1)}{2}$.

\smallskip
{\bf (2)} 
The set $\llbracket\mathfrak g(m,s),m\rrbracket$ contains no integral power of a prime.
\end{proposition}

\begin{proof}
Let $K$ be a field with $|K|\ge m+m!\frac{s(s+1)}{2}$. We fix an arbitrary $m$-tuple $(\alpha_1,\ldots,\alpha_m)\in\mathbb D_{1,m}(K^{\ast})$. For $(r,\sigma)\in \llbracket1,s\rrbracket\times S_m$ we consider the $m\times m$ matrix
$$\Delta_{r,\sigma}:=\begin{bmatrix} 
1 & \alpha_1 & \cdots & \alpha_1^{m-r-1} & x_{\sigma(1)}& \cdots & x_{\sigma(1)}^r\\
1 & \alpha_2 & \cdots & \alpha_2^{m-r-1} & x_{\sigma(2)}& \cdots & x_{\sigma(2)}^r\\
\cdots & \cdots & \cdots & \cdots & \cdots & \cdots & \cdots\\
1 & \alpha_m & \cdots & \alpha_m^{m-r-1} & x_{\sigma(m)}& \cdots & x_{\sigma(m)}^r\\
\end{bmatrix}$$
with entries in $K[x_1,\ldots,x_m]$. Endowing $\mathbb N^m$ with the lexicographic order, the greatest $(d_1,\ldots,d_m)\in\mathbb N^m$ with the property that the coefficient of $\prod_{i=1}^m x_i^{d_i}$ in $\det(\Delta_{r,\sigma})$ is non-zero is $(r,r-1,\ldots,1,0,\ldots,0)$, the coefficient of $x_1^rx_2^{r-1}\cdots x_r^1$ being up to sign equal to the Vandermonde determinant $W(\alpha_{i_1},\ldots,\alpha_{i_{m-r}})$ with $\{i_j|j\in \llbracket1,m-r\rrbracket\}=\llbracket1,m\rrbracket\setminus\sigma^{-1}(\llbracket1,r\rrbracket)$. 
This implies that the product 
$$f:=\prod_{(r,\sigma)\in \llbracket1,s\rrbracket\times S_m} \det(\Delta_{r,\sigma})\in K[x_1,\ldots,x_m]$$ 
is not a constant. Clearly, the partial degree of $h:=f\prod_{(i,j)\in \llbracket1,m\rrbracket^2,i<j} (x_i-x_j)$ in each variable $x_i$ with $i\in \llbracket1,m\rrbracket$ is at most 
$$m-1+\sum_{(r,\sigma)\in \llbracket1,s\rrbracket\times S_m} r=m-1+|S_m|\sum_{r=1}^s r=m-1+m!\frac{s(s+1)}{2}<|K|.$$ 

Thus the function $K^m\rightarrow K$ defined by evaluating $h$ is not the zero function. Hence there exists $(\beta_1,\ldots,\beta_m)\in K^m$ such that $h(\beta_1,\ldots,\beta_m)\neq 0$. Therefore $(\beta_1,\ldots,\beta_m)\in\mathbb D_{1,m}(K)$ and for each $\sigma\in S_m$, as $\prod_{r=1}^s \det\bigl(\Delta_{r,\sigma}(\beta_1,\ldots,\beta_m)\bigr)\neq 0$, from Lemma \ref{L22}(2) we get that $\bigl((\alpha_1,\ldots,\alpha_m),(\beta_{\sigma(1)},\ldots,\beta_{\sigma(m)})\bigr)\in\mathbb D_{1,m}(K)^2$ is $\llbracket1,s\rrbracket$-generic. So part (1) holds.

For part (2), it suffices to show that the assumption that there exists a finite field $L$ with $|L|\in \llbracket\mathfrak g(m,s),m\rrbracket$ leads to a contradiction. As $|L|\ge \mathfrak g(m,s)$, there exists $\bigl((\alpha_1,\ldots,\alpha_m),(\beta_1,\ldots,\beta_m)\bigr)\in\mathbb D_{1,m}(K)^2$ such that for every permutation $\sigma\in S_m$, the pair $\bigl((\alpha_1,\ldots,\alpha_m),(\beta_{\sigma(1)},\ldots,\beta_{\sigma(m)})\bigr)$ in $\mathbb D_{1,m}(K)^2$ is $\llbracket1,s\rrbracket$-generic. So $|L|\ge m$ and thus $|L|=m$. But $|L|>m+1$ by Lemma \ref{F12}(4), a contradiction.
\end{proof}

\begin{lemma}\label{L29}
For $m\in\mathbb N^{\ast}\setminus\{1,2,3\}$ the following properties hold.

\medskip
{\bf (1)} The sequence $\bigl(\mathfrak h(m,s)\bigr)_{s\in \llbracket2,m-1\rrbracket}$ is non-decreasing.

\smallskip
{\bf (2)} We have $\mathfrak h(m,m-1)\le\mathfrak g(m,1)$.

\smallskip
{\bf (3)} The set $\llbracket\mathfrak h(m,2),m+2\rrbracket$ contains no integral power of a prime.

\smallskip
{\bf (4)} We have $\mathfrak h(5,2)=9$ and for $m\neq 5$, $\mathfrak h(m,2)=p^q+1$ where $q\in\mathbb N^{\ast}$, $p$ is a prime, and $p^q$ is the largest positive integer power of a prime less than $m+3$.
\end{lemma}

\begin{proof}
Part (1) is clear by the very definition.

To prove part (2) it suffices to show that for a finite field $K$ with $|K|\ge \mathfrak g(m,1)$ there exists $\bigl((\alpha_1,\ldots,\alpha_m),(\beta_1,\ldots,\beta_m)\bigr)\in\mathbb D_{1,m}(K)^2$ such that for every permutation $\sigma\in S_m$, the Lagrange interpolating polynomial $h_{\sigma}\in K[x]$ defined by the properties that $\deg(h_{\sigma})\le m-1$ and $h_{\sigma}(\alpha_i)=\beta_{\sigma(i)}$ for each $i\in \llbracket1,m\rrbracket$ has degree $m-1$. We choose $\bigl((\alpha_1,\ldots,\alpha_m),(\beta_1,\ldots,\beta_m)\bigr)\in\mathbb D_{1,m}(K)^2$ such that $\bigl((\alpha_1,\ldots,\alpha_m),(\beta_{\sigma(1)},\ldots,\beta_{\sigma(m)})\bigr)\in\mathbb D_{1,m}(K)^2$ is $\{1\}$-generic for each $\sigma\in S_m$ by Definition \ref{D26}(1). As $\deg(h_{\sigma})=m-1$ by Lemma \ref{F12}(4), part (2) holds.

We begin the proofs of parts (3) and (4) by showing that $\mathfrak h(m,2)\le m+3$ if $m\neq 5$ and $\mathfrak h(5,2)\le 9$. Let $K$ be a field with $|K|\ge m+3$. It suffices to show that if $m\neq 5$ or $|K|\ge 9$, then there exists a pair $\bigl((\alpha_1,\ldots,\alpha_m),(\beta_1,\ldots,\beta_m)\bigr)\in\mathbb D_{1,m}(K)^2$ such that for every permutation $\sigma\in S_m$, each polynomial $h_{\sigma}\in K[x]$ with the property that $h_{\sigma}(\alpha_i)=\beta_{\sigma(i)}$ for every $i\in \llbracket1,m\rrbracket$ has degree at least $2$. As we have $\mathfrak h(m,2)\le\mathfrak h(m,m-1)\le\mathfrak g(m,1)$ by parts (1) and (2), we can assume that $|K|\le \mathfrak g(m,1)$ by Lemma \ref{F12}(4). 

We first show that the number $\binom{|K|}{m}$ of subsets $Y$ of $K$ with $m$ elements is greater than $|K|(|K|-1)$ except when $(m,|K|)\in\{(4,7),(5,8)\}$. This is equivalent to $\prod_{i=2}^{m-1} (|K|-i)>m!$. To check this, as the left hand side is an increasing function of $|K|$, we can assume that $|K|=m+3$ if $m\ge 6$ and $|K|=m+4$ if $m\in\{4,5\}$. So we need to check that $\frac{(m+1)!}{3!}>m!$ if $m\ge 6$ and $\frac{(m+2)!}{4!}>m!$ if $m\in\{4,5\}$ which are true as they are equivalent to $m+1>6$ if $m\ge 6$ and to $(m+1)(m+2)>24$ if $m\in\{4,5\}$. As $|\AGL_1(K)|=|K|(|K|-1)$, it follows that if $(m,|K|)\notin\{(4,7),(5,8)\}$ then there exist subsets $Y$ and $Z$ of $K$ with $m$ elements such that $Z\neq a(Y)$ for each $a=\e(f_a)\in\AGL_1(K)$. This also holds if $m=4$ and $K=\mathbb F_7$ as, by identifying $\mathbb F_7=\llbracket0,6\rrbracket$, we can take $Y=K\setminus \{0,1,2\}$ and $Z=K\setminus\{0,1,3\}$, but it does not hold if $m=5$ and $|K|=8$ by Lemma \ref{L27}(1).

So if $(m,|K|)\neq (5,8)$ and $\bigl((\alpha_1,\ldots,\alpha_m),(\beta_1,\ldots,\beta_m)\bigr)\in\mathbb D_{1,m}(L)^2$ is such that $Y=\{\alpha_i|i\in \llbracket1,m\rrbracket\}$ and $Z=\{\beta_i|i\in \llbracket1,m\rrbracket\}$, then for each $\sigma\in S_m$ and polynomial $h_{\sigma}\in K[x]$ with the property that $h_{\sigma}(\alpha_i)=\beta_{\sigma(i)}$ for each $i\in \llbracket1,m\rrbracket$ we have $h_{\sigma}\notin\{f_a|a\in\AGL_1(K)\}$ and hence $\deg(h_{\sigma})\ge 2$. 

We conclude that $\mathfrak h(m,2)\le m+3$ if $m\neq 5$ and $\mathfrak h(5,2)\le 9$.

For part (3), we first show that the assumption that $\llbracket\mathfrak h(m,2),m+2\rrbracket$ contains $|L|$ with $L$ a finite field leads to a contradiction. As $\mathbb D_{1,m}(L)$ is non-empty, we have $|L|\in\{m,m+1,m+2\}$. For $\bigl((\alpha_1,\ldots,\alpha_m),(\beta_1,\ldots,\beta_m)\bigr)\in\mathbb D_{1,m}(L)^2$, let $Y_0:=L\setminus \{\alpha_i|i\in \llbracket1,m\rrbracket\}$ and $Z_0:=L\setminus\{\beta_i|i\in \llbracket1,m\rrbracket\}$. As $|Y_0|=|Z_0|\le 2$, there exists $(\gamma,\delta)\in L^{\ast}\times L$ such that for $f(x):=\gamma x+\delta\in L[x]$ we have $f(Y_0)=Z_0$, and hence $f(\{\alpha_i|i\in \llbracket1,m\rrbracket\})=\{\beta_i|i\in \llbracket1,m\rrbracket\}$. Thus, if $\sigma\in S_m$ is such that $f(\alpha_i)=\beta_{\sigma(i)}$ for every $i\in \llbracket1,m\rrbracket$, as $\deg(f)=1$, we have $\mathfrak h(m,2)>|L|$, a contradiction. Thus part (3) holds. 

As in the prior paragraph we argue that $\mathfrak h(5,2)\neq 8$.

Part (4) follows from the last three paragraphs.\end{proof}

We have the following variant of Theorem \ref{T15}.

\begin{theorem}\label{T20}
Suppose that $(n,m)\in (\mathbb N^\ast\setminus\{1\})^2$. Let $K$ be a field. Then the following properties hold.

\medskip
{\bf (1)} If $m\ge 3$ and $|K|\ge\mathfrak g(m,1)$ (e.g., if $|K|\ge m+m!\ge 9$ by Proposition \ref{PR35}(1)), then $\upsilon^{\S}_{n,m}(K)\ge \upsilon_{n,m}(K)\ge m+1$.

\smallskip
{\bf (2)} Suppose that $m\ge 4$. If either $m\neq 5$ and $|K|\ge m+3$ or $m=5$ and $|K|\ge 9$, then $\upsilon^{\S}_{n,m}(K)\ge \upsilon_{n,m}(K)\ge m$.

\smallskip
{\bf (3)} Suppose that $m\ge 4$ and $|K|\ge\mathfrak g(m,\lfloor\frac{2m-1}{3}\rfloor)$ (e.g., this holds if we have $|K|\ge m+m!\frac{\lfloor\frac{2m-1}{3}\rfloor^2+\lfloor\frac{2m-1}{3}\rfloor}{2}\ge 76$ by Proposition \ref{PR35}(1)). Writing $m=3r+\epsilon$ with $(r,\epsilon)\in\mathbb N^{\ast}\times\{0,1,2\}$, let $\phi(r,\epsilon):=4r^2+\epsilon(5-\epsilon)r+\epsilon(\epsilon-1)$. Then 
$$\upsilon^{\S}_{2,m}(K)\ge\upsilon_{2,m}(K)\ge\phi(r,\epsilon)\ge \frac{4m^2}{9}.$$
\end{theorem}

\begin{proof}
For all parts we can assume that $|K|\ge \sqrt[n]{m}$. Based on Fact \ref{F15} it suffices to show that there exist subsets $Y=\{P_1,\ldots,P_m\}$ and $Z=\{Q_1,\ldots,Q_m\}$ of $K^n$ with $m$ elements such that for every $(a_{\sigma},\sigma)\in\GA_n(K)\times S_m$ with $a_{\sigma}(P_i)=Q_{\sigma(i)}$ for each $i\in \llbracket1,m\rrbracket$ we have $\ell(a_{\sigma})\ge N$, where $N$ is $m+1$ in part (1), is $m$ in part (2), and is $\phi(r,\epsilon)$ in part (3). 

For all parts we choose $P_i:=(\alpha_i,0,\ldots,0)$ and $Q_i:=\bigl(f(\alpha_i),0,\ldots,0\bigr)$ for each $i\in \llbracket1,m\rrbracket$ for a suitable $(\alpha_1,\ldots,\alpha_m,f)\in \mathbb D_{1,m}(K)\times K[x]$. Let $\underline{P}:=(P_1,\ldots,P_m)$. For $\sigma\in S_m$, let $\underline{Q}_{\sigma}:=\bigl(Q_{\sigma(1)},\ldots,Q_{\sigma(m)}\bigr)$. So, $(\underline{P},\underline{Q}_{\sigma})\in\mathbb D_{1,m}(K)^2$. 

For part (1), as $|K|\ge\mathfrak g(m,1),$ we can choose $(\alpha_1,\ldots,\alpha_m,f)$ such that the pair $(\underline{P},\underline{Q}_{\sigma})\in\mathbb D_{n,m}(K)^2$ is $\{1\}$-generic for each $\sigma\in S_m$. From the proof of Theorem \ref{T15}(1) applied to $(\underline{P},\underline{Q}_{\sigma})$ we get that $\ell(a_{\sigma})\ge N$. So part (1) holds.

For part (2), we have $|K|\ge\mathfrak h(m,2)$ by Lemma \ref{L29}(4). So we can choose $(\alpha_1,\ldots,\alpha_m,f)$ such that for each $\sigma\in S_m$ the Lagrange polynomial $f_{\sigma}\in K[x]$ defined by $\deg(f_{\sigma})\le m-1$ and $f_{\sigma}(\alpha_i):=f(\alpha_{\sigma(i)})$ for each $i\in \llbracket1,m\rrbracket$ is not linear. As $\deg(f_{\sigma})\ge 2$, as in the proof of Theorem \ref{T15}(3) applied to $(\underline{P},\underline{Q}_{\sigma})$ we argue that $\ell(a_{\sigma})\ge N$. So part (2) holds.

For part (3), as $|K|\ge\mathfrak g(m,\lfloor\frac{2m-1}{3}\rfloor),$ we can choose $(\alpha_1,\ldots,\alpha_m,f)$ such that $(\underline{P},\underline{Q}_{\sigma})\in\mathbb D_{n,m}(K)^2$ is $\llbracket1,\lfloor\frac{2m-1}{3}\rfloor\rrbracket$-generic for each $\sigma\in S_m$. From the proof of Theorem \ref{T15}(4) applied to $(\underline{P},\underline{Q}_{\sigma})$ we get that $\ell(a_{\sigma})\ge N$. So part (3) holds.
\end{proof}

We have the following variant of Theorem \ref{T17}.

\begin{theorem}\label{T21} Let $K$ be a finite field with $|K|\geq 3$ and $m\in\bigl\llbracket3,\bigl\lfloor\frac{|K|^2}{2}\bigr\rfloor\bigr\rrbracket$. Then the following properties hold.

\medskip
{\bf (1)} The number $\ln \bigl(\upsilon_{2,m}(K) \bigr)$ is strictly greater than
$$\frac{(m-2)[\ln (|K|^2-m)-\frac{m+1}{m-2}\ln(m+1)]-(|K|^2-2)\ln\left(\frac{|K|^2-m}{|K|^2-2}\right)+\ln \frac{4e|K|^2}{|K|^2-1}}{|K|+\mathfrak c_{|K|}}.$$

{\bf (2)} The number $\ln \bigl(\upsilon^{\S}_{2,m}(K) \bigr)$ is strictly greater than
$$\frac{(m-2)[\ln (|K|^2-m)-\frac{m+1}{m-2}\ln(m+1)]-(|K|^2-2)\ln\left(\frac{|K|^2-m}{|K|^2-2}\right)+\ln \frac{4e|K|^2}{|K|+1}}{|K|+\mathfrak c_{|K|}}.$$
\end{theorem}

\begin{proof} We prove part (1). Let $N:=\upsilon_{2,m}(K)$; we have $N\ge 2$ by Lemma \ref{L28}. Due to the existence of $\sigma\in S_m$ in Fact \ref{F15}, Inequality (\ref{EQ56}) becomes
\begin{equation}\label{EQ58}
(1-|K|^{-2}) N^{|K|+\mathfrak c_{|K|}}>\frac{m!{{|K|^2}\choose{m}}^2}{|K|^2(|K|^2-1)m!{{|K|^2}\choose{m}}}=\left(\prod_{i=|K|^2-m+1}^{|K|^2-2} i\right)\left(\prod_{i=2}^m i^{-1}\right).
\end{equation}

We have
$$\ln\left(\prod_{i=2}^m i\right)=\sum_{i=2}^{m} \ln i \ \ < \int \limits_{2}^{m+1} \ln t \ dt =x(\ln x -1) \biggl|^{m+1}_{2}=-m+2-D_m,$$
where
$$D_m:=\ln(4e)-(m+1)\ln(m+1).$$
So based on this and the last part of the proof of Theorem \ref{T17}(1), Equation (\ref{EQ58}) gives that
$$(|K|+\mathfrak c_{|K|})\ln N +\ln (1-|K|^{-2})$$
is strictly greater than
$$(m-2)[\ln (|K|^2-m)]-(|K|^2-2)\ln\left(\frac{|K|^2-m}{|K|^2-2}\right)+D_m$$
and hence part (1) holds.

Part (2) is proved similarly based on Proposition \ref{PR30}(2), the only change being that the expression $1-|K|^{-2}$ gets replaced by $|K|^{-1}+|K|^{-2}$.\end{proof}

We have the following asymptotic ranges for $n=2$ and finite fields which are modeled on Corollary \ref{C27}.

\begin{corollary}\label{C31}
Let $\varepsilon\in (0,2)$. Let $\upsilon_{2,2}:=9$, $\upsilon_{2,3}:=13.5$, and $\upsilon_{2,p}:=\frac{27p}{2(p-1)}$ for a prime $p\ge 5$. Then the following properties hold.

\medskip
{\bf (1)} There exists $C(\varepsilon)\in\mathbb N^{\ast}$ such that for every finite field $K$ with $|K|\ge C(\varepsilon)$ we have
$$(2-\varepsilon) |K|\ln |K|<\ln\bigl(\upsilon_{2,\lfloor\frac{|K|^2}{2}\rfloor}(K)\bigr)\le \ln\bigl(\upsilon^{\S}_{2,\lfloor\frac{|K|^2}{2}\rfloor}(K)\bigr)<(14+\varepsilon)|K|\ln |K|.$$

{\bf (2)} Let $p$ be a prime. There exists $D_p(\varepsilon)\in\mathbb N^{\ast}$ such that for every finite field $K$ with $|K|\ge D_p(\varepsilon)$ and $\char(K)=p$ we have inequalities\footnote{We have $\upsilon_{2,p}<14$ iff $p\ge 29$.}
$$(2-\varepsilon) |K|\ln |K|<\ln\bigl(\upsilon_{2,\lfloor\frac{|K|^2}{2}\rfloor}(K)\bigr)\le \ln\bigl(\upsilon^{\S}_{2,\lfloor\frac{|K|^2}{2}\rfloor}(K)\bigr)<[\min(14,\upsilon_{2,p})+\varepsilon]|K|\ln |K|.$$

{\bf (3)} There exists $E(\varepsilon)\in\mathbb N^{\ast}$ such that for every finite field $K$ with $|K|$ a prime greater or equal to $E(\varepsilon)$ we have
$$(2-\varepsilon) |K|\ln |K|<\ln\bigl(\upsilon_{2,\lfloor\frac{|K|^2}{2}\rfloor}(K)\bigr)\le \ln\bigl(\upsilon^{\S}_{2,\lfloor\frac{|K|^2}{2}\rfloor}(K)\bigr)<(9+\varepsilon)|K|\ln |K|.$$
\end{corollary}

\begin{proof}
The lower bounds of all parts follow from Theorem \ref{T21}(1), as its expression is asymptotic to $2|K|\ln|K|$. 

The upper bounds with $14$ for parts (1) and (2) follow from Corollary \ref{C30}(1) applied to $\epsilon=10$ as the $\ln$ of its bound is asymptotic to $14|K|\ln|K|$.

We now check the upper bound of part (2) that involves $\upsilon_{2,p}$. For $p=2$, the upper bound follows from Corollary \ref{C30}(3) as the $\ln$ of its bound is asymptotic to $9|K|\ln|K|$. For $p\ge 5$ (resp.\ $p=3$), the upper bound follows from Corollaries \ref{C4}(3) (resp.\ \ref{C4}(2)) and \ref{C30}(2) as the $\ln$ of its bound is asymptotic to $\upsilon_{2,p}|K|\ln|K|$.

The upper bound for part (3) follows from Corollary \ref{C30}(1) applied to $\epsilon=5$ as the $\ln$ of its bound is asymptotic to $9|K|\ln|K|$.\end{proof}

\section{Polynomial bounds for sequences}\label{S31}

In this section we prove the following strengthening of Theorem \ref{T18}.

\begin{theorem}\label{T22} Let $\bigl((K_i,m_i,q_i)\bigr)_{i\in\mathbb N^{\ast}}$ be a sequence of triples with each $K_i$ a finite field, $m_i\in \llbracket2,\lfloor\frac{|K_i|^2}{2}\rfloor\rrbracket$, and $q_i\in \llbracket1,m_i-1\rrbracket$. We assume that $\lim \limits_{i\to +\infty}|K_i|=+\infty$. We also assume that $2q_i\le m_i$ for each $i\in\mathbb N^{\ast}$ or the set $\{\frac{q_i}{|K_i|}|i\in\mathbb N^{\ast}\}$ of rational numbers is bounded. Then the following statements are equivalent.

\medskip
{\bf (1)} There exists $C\in (0,\infty)$ such that $\upsilon^{\S}_{2,m_i}(K_i) <m_i^C$ for each $i\in\mathbb N^{\ast}$.

\medskip
{\bf (2)} There exists $C\in (0,\infty)$ such that $\upsilon_{2,m_i}(K_i) <m_i^C$ for each $i\in\mathbb N^{\ast}$.

\medskip
{\bf (3)} There exists $C\in (0,\infty)$ such that $\pi^{\S}_{2,m_i}(K_i) <m_i^C$ for each $i\in\mathbb N^{\ast}$.

\medskip
{\bf (4)} There exists $C\in (0,\infty)$ such that $\pi_{2,m_i}(K_i) <m_i^C$ for each $i\in\mathbb N^{\ast}$.

\smallskip
{\bf (5)} The set $\{\frac{m_i}{|K_i|}|i\in\mathbb N^{\ast}\}$ of rational numbers is bounded.

\smallskip
{\bf (6)} There exists $C\in (0,\infty)$ such that $\Xi^{\S}_{2,q_i,m_i-q_i}(K_i) <m_i^C$ for each $i\in\mathbb N^{\ast}$.

\smallskip
{\bf (7)} There exists $C\in (0,\infty)$ such that $\Xi_{2,q_i,m_i-q_i}(K_i) <m_i^C$ for each $i\in\mathbb N^{\ast}$.

\smallskip
{\bf (8)} There exists $C\in (0,\infty)$ such that $\xi^{\S}_{2,q_i,m_i-q_i}(K_i) <m_i^C$ for each $i\in\mathbb N^{\ast}$.

\smallskip
{\bf (9)} There exists $C\in (0,\infty)$ such that $\xi_{2,q_i,m_i-q_i}(K_i) <m_i^C$ for each $i\in\mathbb N^{\ast}$.
\end{theorem}

\begin{proof} By Theorem \ref{T18}, (3) to (5) are equivalent. We have $(3)\Rightarrow (1)\Rightarrow (2)$ as $\upsilon_{2,m_i}(K_i)\le\upsilon^{\S}_{2,m_i}(K_i)\le \pi^{\S}_{2,m_i}(K_i)$. Also, we have $(3)\Rightarrow (6)\Rightarrow (7)\wedge (8)$ and $(7)\vee (8)\Rightarrow (9)$ as $\xi_{2,q_i,m_i-q_i}(K_i)\le \Xi_{2,q_i,m_i-q_i}(K_i)\le \Xi^{\S}_{2,q_i,m_i-q_i}(K_i)\le\pi^{\S}_{2,m_i}(K_i)$ and $\xi_{2,q_i,m_i-q_i}(K_i)\le\xi^{\S}_{2,q_i,m_i-q_i}(K_i)$ by Definitions \ref{D21}(1) and (2) and \ref{D24}(3) and Proposition \ref{PR32}(1) and (4). Hence it suffices to prove the two implications $(9)\Rightarrow (2)\Rightarrow (5)$. 

{$\pmb{(2)\Rightarrow (5)}$.} As for a fixed $m\in\mathbb N^{\ast}$, the set of all values $\pi_{2,m}(K)$ with $K$ a finite field is finite by Theorem \ref{T7}(1), we can assume that $m_i\ge 5$ for each $i\in\mathbb N^{\ast}$.

As $|K_i|$ approaches $+\infty,$ the denominator of the lower bound in Theorem \ref{T21}(1) is asymptotically equivalent to $|K_i|$ and $\lim \limits_{i\to +\infty } \frac{4e(\ln |K_i|^2)-\ln(|K_i|^2-1)}{|K_i|\ln |K_i|}=0$. Moreover, we have $\ln (m_i)\le 2\ln |K_i|-\ln 2$. The last two sentences, Theorem \ref{T21}(1), and statement (2) imply that for some constant $C_1\in (0,\infty)$ that does not depend on $i$, the expression
$$(m_i-2)\left[\ln (|K_i|^2-m_i)-\frac{m_i+1}{m_i-2}\ln(m_i+1)\right]-(|K_i|^2-2)\ln\left(1-\frac{m_i-2}{|K_i|^2-2}\right)$$
is strictly less than $C_1 |K_i| \ln (|K_i|)$.

Denoting $t_i:=\frac{m_i-2}{|K_i|^2-2}$ and $s_i:=\frac{m_i}{|K|^2_i}$, we have $\frac{3}{|K_i|^2-2} \leq t_i<s_i\le\frac{1}{2}$. Rewriting $\frac{1}{|K_i|^2-2}$ times the above inequality using $t_i$ instead of $m_i$, we get
$$t_i\left[\ln\left((|K_i|^2-2)(1-t_i)\right)-\left(1+\frac{3}{t_i (|K_i|^2-2)}\right)\ln\bigl(t_i (|K_i|^2-2)+3\bigr)\right]-\ln(1-t_i)$$ is strictly less than $C_1 \frac{\ln (|K_i|)}{|K_i|-\frac{2}{|K_i|}}$.

Because $t_i (1+\frac{3}{t_i (|K_i|^2-2)}) \ln(1+\frac{3}{t_i (|K_i|^2-2)}) \leq 2t_i \ln\left(1+\frac{3}{t_i (|K_i|^2-2)}\right) < \frac{6}{(|K_i|^2-2)},$ for a different constant $C_2\in (0,\infty)$ we get that
$$ t_i\left[\ln (|K_i|^2-2) + \ln (1-t_i) - \left(1+\frac{3}{t_i (|K_i|^2-2)}\right) \ln\bigl(t_i (|K_i|^2-2)\bigr)\right]-\ln(1-t_i)$$
is strictly less than $C_2 \frac{\ln (|K_i|)}{|K_i|}$.

Finally, with another constant $C_3\in (0,\infty)$, we can simplify this to
$$ -t_i\ln (t_i) -(1-t_i)\ln (1-t_i) < C_3 \frac{\ln (|K_i|)}{|K_i|}.$$

This implies that $-t_i\ln (t_i) < C_3 \frac{\ln (|K_i|)}{|K_i|}.$ As $t_i\leq \frac{1}{2},$ we get that $t_i < \frac{C_3}{\ln 2}\frac{\ln (|K_i|)}{|K_i|}.$ For all large enough $i$ this is less than $|K_i|^{-\frac{1}{2}},$ so $-\ln(t_i) > \frac{1}{2} \ln (|K_i|)$; thus $\frac{t_i}{2} \ln (|K_i|)< C_3 \frac{\ln (|K_i|)}{|K_i|}$ by transitivity, hence $\frac{m_i-2}{|K_i|}< t_i|K_i|<2C_3$. So $(2)\Rightarrow (5)$.

It follows that statements (1) to (5) are equivalent. 

{$\pmb{(9)\Rightarrow (2)}$.} We have $2(m_i-q_i)+q_i=2m_i-q_i\le 2m_i\le |K_i|^2$. From this, Proposition \ref{PR32}(3) applied to $(m,s)=(q_i,m_i-q_i)$, and statement (9) we get that for each $i\in\mathbb N^{\ast}$ we have inequalities $\upsilon_{2,m_i-q_i}(K_i)\le \xi_{2,q_i,m_i-q_i}(K_i)\le m_i^C$. 

We assume that $2q_i\le m_i$ for each $i\in\mathbb N^{\ast}$. If $m_i\ge 4$, then $\bigl(\frac{m_i}{2}\bigr)^{2C}\ge m_i^C$. Thus, as $2(m_i-q_i)\ge m_i$ and as $\upsilon_{2,m_i-q_i}(K_i)=1$ if $m_i-q_i\le 2$ and therefore if $m_i<4$, it follows that 
$\upsilon_{2,m_i-q_i}(K_i)\le (m_i-q_i)^{2C}$. From this and the implication $(2)\Rightarrow (5)$ applied to the sequence of triples $\bigl((K_i,m_i-q_i,1)\bigr)_{i\in\mathbb N^{\ast}}$ we get that the set $\{\frac{m_i-q_i}{|K_i|}|i\in\mathbb N^{\ast}\}$ is bounded. From this and the inequality $2(m_i-q_i)\ge m_i$ we get that statement (5) holds. As $(5)\Leftrightarrow (2)$, we get that $(9)\Rightarrow (2)$ holds if $2q_i\le m_i$ for each $i\in\mathbb N^{\ast}$. 

We now assume that the set $\{\frac{q_i}{|K_i|}|i\in\mathbb N^{\ast}\}$ of rational numbers is bounded. From this and the implication $(5)\Rightarrow (2)$ applied to the sequence of triples $\bigl((K_i,q_i,1)\bigr)_{i\in\mathbb N^{\ast}}$ we get that there exists $D\in (0,\infty)$ such that $\upsilon_{2,q_i}(K_i)\le q_i^D\le m_i^D$ for each $i\in\mathbb N^{\ast}$. From this and Proposition \ref{PR32}(2) applied again to $(m,s)=(q_i,m_i-q_i)$ we get that for each $i\in\mathbb N^{\ast}$ we have inequalities 
$$\upsilon_{2,m_i}(K_i)\le\upsilon_{2,q_i}(K_i)\xi_{2,q_i,m_i-q_i}\le m_i^Dm_i^C=m_i^{C+D}.$$ Thus $(9)\Rightarrow (2)$ holds if the set $\{\frac{q_i}{|K_i|}|i\in\mathbb N^{\ast}\}$ of rational numbers is bounded.

We conclude that $(9)\Rightarrow (2)$ and that the theorem holds.
\end{proof}

\section{Open problems and questions}\label{S32}

The open problems and questions are grouped on topics as follows.

\subsection{On odd permutations}\label{S32.1}
For a field $K$, let $O_n(K):=\GA_n(K)/SGA_n(K)$. Let $O_{n,1}(K)$ be the kernel of the surjective Jacobian determinant homomorphism $O_n(K)\rightarrow K^{\ast}$. If $K$ is finite, let $O_{n,2}(K)$ be the kernel of the homomorphism $O_n(K)\rightarrow\perm(K^n)/\Alt(K^n)$ defined by the rule $a\mapsto a(K)\Alt(K^n)$ (cf.\ Theorem \ref{T3}); we have $[O_n(K):O_{n,1}(K)]=|K|-1$ and $[O_n(K):O_{n,2}(K)]\in\{1,2\}$. We identify $O_{n,1}(K)=\GA_n^{\det=1}(K)/SGA_n(K)$, where $\GA_n^{\det=1}(K)$ is the normal subgroup of $\GA_n(K)$ formed by automorphisms of Jacobian determinant $1$. 

\smallskip\noindent
{\bf Question.} Is it true that if $n\ge 3$ and $K$ is finite with $|K|\ge 3$, we have an inclusion $O_{n,1}(K)\subset O_{n,2}(K)$ (equivalently, there exists no automorphism $a\in\GA_n(K)$ such that $a(K)\notin\Alt(K^n)$ and the Jacobian determinant of $a$ is $1$)?

\smallskip
By viewing $\GA_n^{\det=1}(K)$ as the group of $K$-valued points of a rationally connected infinite dimensional reduced group scheme over $\Spec K$, Theorem \ref{T3} would favor a positive answer.

\smallskip\noindent
{\bf Problem.} If $K$ is a finite field with $|K|$ odd, identify pairs $(l,s)\in (\mathbb N^{\ast}\setminus\{1\})^2$ such that there exists $a\in\TGA_n(K)[l]$ with $a(K)$ odd and $\n\bigl(a(K)\bigr)\le s$.

\subsection{On large fields}\label{S32.2} 
Let $K$ be a field and $m\in\mathbb N^{\ast}\setminus\{1\}$. 

\smallskip\noindent
{\bf Question.} Is it true that there exists a smallest $N_m\in\mathbb N^{\ast}$ (resp.\ $N^{\S}_m\in\mathbb N^{\ast}$) which is a power of a prime and such that for $|K|\ge N_m$ (resp.\ $|K|\ge N^{\S}_m$) the value of $\pi_{2,m}(K)$ (resp.\ $\pi^{\S}_{2,m}(K)$) does not depend on $K$? 

\smallskip
For instance, we have $N_2=N_2^{\S}=1$, $N_3=N_3^{\S}=4$ by Theorem \ref{T19}(1), and $N_4=5$ and $N_4^{\S}\in\{5,7\}$ by Theorem \ref{T19}(2).

\smallskip\noindent
{\bf Problem.} If $N_m$ (resp.\ $N^{\S}_m$) exists, compute the value $\pi_m$ (resp. $\pi^{\S}_m$) of $\pi_{2,m}(K)$ (resp.\ $\pi^{\S}_{2,m}(K)$) one gets for $|K|\ge N_m$ (resp.\ for $|K|\ge N^{\S}_m$).

\smallskip
We have $\pi_2=\pi_2^{\S}=1$, $\pi_3=\pi_3^{\S}=4$ by Theorem \ref{T19}(1), and $\pi_4=\pi_4^{\S}=9$ by Theorem \ref{T19}(2). In general, if $N_m$ (resp.\ $N_m^{\S}$) exists, then by and with the notation of Corollary \ref{C2} we have
$\lfloor\frac{m+2}{2}\rfloor(m-1)\le\pi_m\le \pi^{\S}_m\le (m-1)^2$ for $m\ge 7$ and $\pi_m=\pi^{\S}_m=(m-1)^2$ for $m\in\{5,6\}$.

\subsection{On $n=2$ for finite fields}\label{S32.3} Let $K$ be a finite field with $|K|\ge 3$.

\smallskip\noindent
{\bf Problem 1.} Classify all $\sigma\in\Perm(K^2)$ that have strict minimal representations.

\smallskip\noindent
{\bf Problem 2.} Compute the Furter's lengths of finite fields.

\smallskip\noindent
{\bf Problem 3.} Decide which of the following possible ten limits 
$$\lim \limits_{|K|\in 2\mathbb N+1,|K|\to \infty} \frac{\ln \pi_{2,|K|^2}(K)}{|K|\ln |K|},\;\;\;\;\lim \limits_{n\to \infty} \frac{\ln \pi_{2,2^{2n}-2}(\mathbb F_{2^n})}{n2^n\ln 2},$$
$$\lim \limits_{|K|\in 2\mathbb N+1,|K|\to \infty} \frac{\ln \pi^{\S}_{2,|K|^2-2}(K)}{|K|\ln |K|},\;\;\;\;\lim \limits_{n\to \infty} \frac{\ln \pi^{\S}_{2,2^{2n}-2}(\mathbb F_{2^n})}{n2^n\ln 2},$$
$$\lim \limits_{|K|\in 2\mathbb N+1,|K|\to \infty} \frac{\ln \upsilon_{2,\frac{|K|^2-1}{2}}(K)}{|K|\ln |K|},\;\;\;\;\lim \limits_{n\to \infty} \frac{\ln \upsilon_{2,2^{2n-1}}(\mathbb F_{2^n})}{n2^n\ln 2},$$
$$\lim \limits_{|K|\in 2\mathbb N+1,|K|\to \infty} \frac{\ln \upsilon^{\S}_{2,\frac{|K|^2-1}{2}}(K)}{|K|\ln |K|},\;\;\textup{and}\;\;\lim \limits_{n\to \infty} \frac{\ln \upsilon^{\S}_{2,2^{2n-1}}(\mathbb F_{2^n})}{n2^n\ln 2},$$
$$\lim \limits_{|K|\in 2\mathbb N+1,|K|\to \infty} \frac{\j_K}{|K|},\;\;\textup{and}\;\;\lim \limits_{n\to \infty} \frac{\j_{\mathbb F_{2^n}}}{2^n},$$
exists and compute all those that exist.

\medskip
If the first (resp.\ second, third, fourth, one among the fifth to eight, ninth, or tenth) limit exists, then the limit is at least $2$ and at most $13.5$ (resp.\ $18$, $13.5$, $18$, $9$, $13.5$, or $18$) by Corollaries \ref{C18}, \ref{C27}, and \ref{C31}(2).

\subsection{On special values of $\ell_{n,K}$}\label{S32.4} Let $n\in\mathbb N^{\ast}\setminus\{1\}$ and $K$ a finite field. 

\smallskip\noindent
{\bf Problem 1.} If $|K|$ is odd, compute the constant $\kappa_{n,K}$ of Proposition \ref{PR16}(1). 

\smallskip\noindent
{\bf Problem 2.} If $s\in \llbracket3,n+1\rrbracket$, with $s$ odd if $4\mid |K|$, compute the constants $\kappa^{\d=s-1}_{n,K;s}$ and $\kappa^{\S,\d=s-1}_{n,K;s}$ of Lemma \ref{F9}(1).

\smallskip\noindent
{\bf Problem 3.} Let $s\in\llbracket3,|K|\rrbracket$. If $s$ is even and $4\nmid |K|$ (resp.\ if $s$ is odd) compute the set $\Pi_{n,s\textup{-cycle}}(K)$ (resp.\ the sets $\Pi_{n,s\textup{-cycle}}(K)$ and $\Pi^{\S}_{n,s\textup{-cycle}}(K)$).

\subsection{On atoms and regularity}\label{S32.5} Let $K$ be a field. For atoms of semigauges we refer to \cite{Pro}, Sect.\ 2.
For instance, each $a\in\GA_n(K)$ with $\ell(a)$ a prime number is an atom of the semigauge $\ln \ell_{\GA_n(K)}$.

\smallskip\noindent
{\bf Problem 1.} Describe all atoms of the semigauge $\ln \ell_{\GA_n(K)}$. 

\smallskip\noindent
{\bf Problem 2.} If $K$ is finite, describe all atoms of the semigauge $\ln\ell_{n,K}$ and all atoms of the semigauge $\ln\ell^{-}_{n,K}$. 

\smallskip\noindent
{\bf Problem 3.} Classify all fields $K$ for which $\ell_{\GA_2(K)}:\GA_2(K)\rightarrow\mathbb N^{\ast}$ is L-regular. 

\smallskip\noindent
{\bf Problem 4.} If either $n=2$ with $|K|\notin\{2,4\}$ or $n\ge 3$ and $|K|=2$, classify all cases when $\ell_{n,K}$ is regular (or L-regular).

\subsection{On invariants of finite subsets}\label{S32.6} Let $K$ be a finite field and $n\in\mathbb N^{\ast}\setminus\{1\}$. Let $\mathcal P_{\ge 2}(K^n)$ be the set formed by all subsets of $K^n$ with at least $2$ elements. We consider the function 
$$\mathbb I_{n,K}:\mathcal P_{\ge 2}(K^n)\rightarrow \mathbb N^{14}$$ 
defined by $\mathbb I_{n,K}(Y):=(|Y|,\d_Y,\s_Y,\overline{\s}_Y,\q_Y,\mq_Y,\L_Y,\l_Y,\c_Y,\r_Y,\dir_Y,\w_Y,\v_Y,\t_Y)$, with the convention that $\L_Y:=0$ if $Y$ contains the zero vector of $K^n$. 

\smallskip\noindent
{\bf Problem.} Describe the image of $\mathbb I_{n,K}$ and the orbits of the action of $\AGL_n(K)$ on the fibers of $\mathbb I_{n,K}$. In particular, study the invariant $\min(\l_Y,\overline{\s}_Y)$ and the set of all possible quintuples $(|Y|,\d_Y,\s_Y,\w_Y,\v_Y)$.

\subsection{On algebraic groups of automorphisms}\label{S32.7} Suppose that $K$ is a finite field. Let $X$ be a reduced scheme of finite type over $\Spec K$. Let $G_1$, $G_2$, and $G_3$ be linear algebraic groups over $\Spec K$ that act faithfully on $X$. Suppose that for each finite field extension $K\rightarrow K^{\prime}$, the subgroup of $\perm\bigl(X(K^{\prime})\bigr)$ generated by $G_1(K^{\prime})$ and $G_2(K^{\prime})$ contains $G_3(K^{\prime})$. 

\smallskip\noindent
{\bf Question.} Is it true that there exists $l\in\mathbb N^{\ast}$ such that for each finite field extension $K\rightarrow K^{\prime}$, every $\sigma\in G_3(K^{\prime})$ is a product $\prod_{i=1}^l \theta_{i,1}\theta_{i,2}$ with $\theta_{i,j}\in G_j(K^{\prime})$ for each $(i,j)\in \llbracket1,l\rrbracket\times\{1,2\}$?

\subsection{On more variants of the $\pi_{n,m}(K)$s}\label{S32.8} Let $(n,m)\in (\mathbb N^{\ast}\setminus\{1,2\})^2$. Let $K$ be a field with $|K|\ge\sqrt[n]{m}$ (resp.\ $|K|\ge\sqrt[n]{m+2}$). Besides $\pi_{n,m}^{\S}(K)$ which was extensively used and, when $K$ is finite, $\pi_{n,m}^{\E}(K)$ which was only occasionally used and has a natural variant $\upsilon_{n,m}^{\E}(K)$, we have the following extra variants of $\pi_{n,m}(K)$.

Let $\pi^{\GA}_{n,m}(K)$ and $\upsilon^{\GA}_{n,m}(K)$ (resp.\ $\pi^{\SGA}_{n,m}(K)$ and $\upsilon^{\SGA}_{n,m}(K)$) be defined similarly to $\pi_{n,m}(K)$ and $\upsilon_{n,m}(K)$ (resp.\ $\pi^{\S}_{n,m}(K)$ and $\upsilon^{\S}_{n,m}(K)$) but using $\GA_n(K)$ instead of $\TGA_n(K)$ (resp.\ using $\SGA_n(K)$ instead of $\STGA_n(K)$). 

Let $\pi^-_{n,m}(K)\in\mathbb N^{\ast}$\index{$\pi^-_{n,m}(K)$ invariant} (resp.\ $\pi^{-,\S}_{n,m}(K)\in\mathbb N^{\ast}$\index{$\pi^{-,\S}_{n,m}(K)$ invariant}) be the smallest with the property that for each $(\underline{P},\underline{Q})\in \mathbb D_{n,m}(K)^2$, there exists $a\in\TGA_n(K)$ (resp.\ $a\in\STGA_n(K)$) such that we have $a(\underline{P})=\underline{Q}$ and $\pi(a)\le\pi^-_{n,m}(K)$ (resp.\ $\pi(a)\le\pi^{-,\S}_{n,m}(K)$). 

Let $\upsilon^-_{n,m}(K)\in\mathbb N^{\ast}$\index{$\upsilon^-_{n,m}(K)$ invariant} (resp.\ $\upsilon^{-,\S}_{n,m}(K)\in\mathbb N^{\ast}$\index{$\upsilon^{-,\S}_{n,m}(K)$ invariant}) be the smallest with the property that for all subsets $Y$ and $Z$ of $K^n$ of cardinality $m$, there exists $a\in\TGA_n(K)$ (resp.\ $a\in\STGA_n(K)$) with $\pi(a)\le \upsilon^-_{n,m}(K)$ (resp.\ $\pi(a)\le \upsilon^{-,\S}_{n,m}(K)$) and $a(Y)=Z$. 

As the notation suggests, as $\pi(a)\le\ell(a)$ for each $a\in\GA_n(K)$, we have the following inequalities $\pi^{-}_{n,m}(K)\le\pi_{n,m}(K)$, $\pi^{-,\S}_{n,m}(K)\le\pi^{\S}_{n,m}(K)$, $\upsilon^-_{n,m}(K)\le \upsilon_{n,m}(K)$, and $\upsilon^{-,\S}_{n,m}(K)\le \upsilon^{\S}_{n,m}(K)$. 

\smallskip\noindent
{\bf Problem.} Classify all cases when $\pi^{-}_{n,m}(K)<\pi_{n,m}(K)$ (or $\pi^{-,\S}_{n,m}(K)<\pi^{-,\S}_{n,m}(K)$ or $\upsilon^-_{n,m}(K)<\upsilon_{n,m}(K)$ or $\upsilon^{-,\S}_{n,m}(K)<\upsilon^{\S}_{n,m}(K)$) and estimate their difference.

\subsection{On special finite fields}\label{S32.9} We consider $\star\in\{\pi,\upsilon,\pi^\E,\upsilon^\E,\pi^\GA,\upsilon^\GA\}$ (resp.\ $\star\in\{\pi^{\S},\upsilon^{\S},\pi^\SGA,\upsilon^{\S},\upsilon^\SGA\}$). 

\smallskip\noindent
{\bf Problem.} Find all finite fields $K$ with the property that there exists a constant $C_K\in (0,\infty)$ such that we have $\star_{n,m}(K)< |K|^{C_K}$ for each $n\in\mathbb N^{\ast}\setminus\{1\}$ and every $m\in \llbracket1,|K|^n-1\rrbracket$ (resp.\ $m\in \llbracket1,|K|^n-2\rrbracket$). 

\subsection{On tame automorphisms}\label{S32.10} Suppose that $n\ge 3$ and $p$ is a prime. Let $K$ be a finite field with $\char(K)=p$. Let 
$$\GA^{\triv}_n(K):=\{a\in\GA_n(K)|a(K)\,\textup{is the identity permutation}\}$$ 
be the kernel of $\varrho_{n,K}$. We denote by $\overline{L}$ an algebraic closure of a field $L$. 

\smallskip\noindent
{\bf Question 1.} Is $\GA^{\triv}_n(K)$ a subgroup of $\TGA_n(K)$ or $\TGA_n(\overline{K})$? 

\smallskip\noindent
{\bf Question 2.} Suppose that there exist infinitely many primes $p$ such that we have $\TGA_n(\overline{K})=\GA_n(\overline{K})$. Is it true that $\TGA_n(L)=\GA_n(L)$ for each algebraically closed field $L$ of characteristic $0$?

\smallskip\noindent
{\bf Question 3.} Is it true that there exists no $a\in\TGA_3(K)$ such that $a^p=1_{\mathbb A^3_K}$ and the set of fixed points of the permutation $a(\overline{K})$ of $\overline{K}^3$ is for $p>2$ the zero locus $x_1x_3+x_2^2=0$ and for $p=2$ the zero locus $x_3(x_1x_3+x_2^2)=0$?

\smallskip
It is a standard piece of algebraic geometry that if the answer to Question 3 is positive, then Nagata's Conjecture made in \cite{N}, Conj.\ 3.1 is true over all fields and one obtains a new proof of \cite{SU}, Cor.\ 9.

\smallskip\noindent
{\bf Problem 1.} Find all pairs $(n,p)$ such that $\TGA_n(\overline{K})=\GA_n(\overline{K})$.

\smallskip\noindent
{\bf Problem 2.} Find all fields $L$ such that the identity $\TGA_n(L)=\GA_n(L)$ holds iff the identity $\TGA_n(\overline{L})=\GA_n(\overline{L})$ holds.

\smallskip\noindent
{\bf Problem 3.} Compute the invariants as in Section \ref{S14} for each $\mathcal W_{f,g,h,d}\in\GA_3(K)$ as in Remark \ref{R12.1}.

\smallskip
For instance, if $\deg(f)\ge 1$ we have 
$$\ell(\mathcal W_{f,g,h,d})=\deg(f)\bigl[\max\bigl(\deg(f),2\bigr)\deg(g)+d-\deg(g)+\deg(h)\bigr]-1.$$

\subsection{On fixed loci}\label{S32.11} 
Let $K$ be an arbitrary algebraically closed field. The following is an extrapolation of Question 3 of Subsection \ref{S32.10}.

\smallskip\noindent
{\bf Problem.} Classify the set $\mathcal V_{n,K,\textup{tame}}$ (resp.\ $\mathcal V_{n,K,\textup{wild}}$) formed by all reduced closed subschemes of $\mathbb A^n_K$ whose $K$-values points are the fixed points of a permutation of $K^n$ defined by $a(K)$ with $a\in\TGA_n(K)$ (resp.\ $a\in\GA_n(K)\setminus\TGA_n(K)$).

\subsection{On varying $n$}\label{S32.12} Let $K$ be a field and $m\in\mathbb N^{\ast}\setminus\{1,2\}$. As in the proof of Theorem \ref{T7}(2) we argue that the sequences $\bigl(\upsilon_{m-1+i,m}(K)\bigr)_{i\in\mathbb N}$ and $\bigl(\upsilon^{\S}_{m-1+i,m}(K)\bigr)_{i\in\mathbb N}$ are non-increasing. 

\smallskip\noindent
{\bf Problem 1.} Classify all triples $(n,m,K)$ for which $\pi_{n,m}(K)\neq\pi_{n+1,m}(K)$ or $\pi^{\S}_{n,m}(K)\neq\pi^{\S}_{n+1,m}(K)$ or $\upsilon_{n,m}(K)\neq\upsilon_{n+1,m}(K)$ or $\upsilon^{\S}_{n,m}(K)\neq\upsilon^{\S}_{n+1,m}(K)$. In particular, classify all non-constant sequences $\bigl(\pi_{m-1+i,m}(K)\bigr)_{i\in\mathbb N}$, $\bigl(\pi^{\S}_{m-1+i,m}(K)\bigr)_{i\in\mathbb N}$, $\bigl(\upsilon_{m-1+i,m}(K)\bigr)_{i\in\mathbb N}$, and $\bigl(\upsilon^{\S}_{m-1+i,m}(K)\bigr)_{i\in\mathbb N}$.

\smallskip\noindent
{\bf Problem 2.} Suppose that $K$ is finite. Let $\sigma\in\Perm(K^n)$ with $n\in\mathbb N^{\ast}\setminus\{1\}$. For $i\in\mathbb N$, by identifying $K^n$ with the zero locus $x_{n+1}=\cdots=x_{n+i}=0$ in $K^{n+i}$, let $\sigma_i\in\Perm(K^{n+i})$ be such that it extends $\sigma$ and we have $\supp(\sigma_i)=\supp(\sigma)$. Classify all the $\sigma$s for which the sequence $\bigl(\ell_{n+i,K}(\sigma_i)\bigr)_{i\in\mathbb N}$ is non-increasing.

\subsection{On counting permutations}\label{S32.13} Let $K$ be a finite field.

\smallskip\noindent
{\bf Problem.} For a pair $(n,N)\in (\mathbb N^{\ast}\setminus\{1,2\})\times (\mathbb N^{\ast}\setminus\{1\})$ and a finite field $K$, estimate the numbers $\varrho_{n,K,N}$ introduced in the beginning of Section \ref{S27}.

\subsection{On polynomial bounds}\label{S32.14} Let $K$ be a finite field.

\smallskip\noindent
{\bf Problem 1.} Find quadruples $(N,\epsilon,C,D)\in\mathbb N\times (0,1]\times (0,\infty)^2$ such that for a pair $(n,m)\in (\mathbb N^{\ast})^2$ the two inequalities $\pi_{n,m}(K)<\infty$ (resp.\ $\pi^{\S}_{n,m}(K)<\infty$) and $m\le C|K|^{n\epsilon}$ imply that $\pi_{n,m}(K)\le Dm^N$ (resp.\ $\pi^{\S}_{n,m}(K)\le Dm^N$).

\smallskip\noindent
{\bf Problem 2.} Is it true that for each finite field $K$ there exist a pair $(C,N)\in (\mathbb N^{\ast})^2$ with the property that for every $(\underline{P},\underline{Q})\in\mathbb D_{n,m}(K)^2$ with $(n,m)\in (\mathbb N^{\ast})^2$, if there exists an automorphism $a\in\GA_n(K)$ with $a(\underline{P})=\underline{Q}$, then we can choose $a$ such that $\ell(a)<Cm^N$? A weaker version of this: is this true for $|K|=2$ at least?

\smallskip
Note that we do not require $a$ to be tame and it is conceivable that some robust enough families of wild automorphisms can yield improvements to Corollary \ref{C24} and Remark \ref{R9}(1).

\subsection{On special invariants}\label{S32.15} Among the many possible problems, only two are mentioned here. 

\smallskip\noindent
{\bf Problem 1.} For $(m,s)\in\mathbb N^{\ast}\times\mathbb N$ with $s\le m-1$ (resp.\ $s\in \llbracket3,m-2\rrbracket$), compute $\mathfrak g(m,s)$ (resp.\ $\mathfrak h(m,s)$) of Definition \ref{D26}(1) (resp.\ (2)).

\smallskip\noindent
{\bf Problem 2.} For $n\ge 3$, get $n$-dimensional variants of the $2$-dimensional results Lemma \ref{L5}, Proposition \ref{PR8} and Corollary \ref{C9}.

\subsection{On fractional unital growth}\label{S32.16} Let $K$ a finite field and $n\in\mathbb N^{\ast}\setminus\{1\}$.

\smallskip\noindent
{\bf Problem.} Compute the sets of rational numbers $\bigl\{\frac{\upsilon_{n,m}(K)}{\upsilon_{n,m-1}(K)}\bigl|m\in \bigl\llbracket1,\lfloor\frac{|K|^n}{2}\rfloor\bigr\rrbracket\bigr\}$ and $\bigl\{\frac{\upsilon^{\S}_{n,m}(K)}{\upsilon^{\S}_{n,m-1}(K)}\bigl|m\in\bigl\llbracket1,\lfloor\frac{|K|^n}{2}\rfloor\bigr\rrbracket\bigr\}$ and their intersections with the interval $(0,1)$.

\subsection{On sequences of invariants}\label{S32.17} Let $K$ be a field and $(n,m)\in (\mathbb N^{\ast}\setminus\{1\})^2$ with $\sqrt[n]{m}\le |K|$. 

\smallskip\noindent
{\bf Problem.} If $\pi_{n,m}(K)\in\mathbb N^{\ast}\setminus\{1\}$ (resp.\ $\pi^{\S}_{n,m}(K)\in\mathbb N^{\ast}\setminus\{1\}$ or $\upsilon_{n,m}(K)\in\mathbb N^{\ast}\setminus\{1\}$ or $\upsilon^{\S}_{n,m}(K)\in\mathbb N^{\ast}\setminus\{1\}$), determine ranges for and inequalities between elements of the sequence $(\natural_{n,m,l})_{l=1}^{\pi_{n,m}(K)-1}$ (resp.\ $(\natural^{\S}_{n,m,l})_{l=1}^{\pi^{\S}_{n,m}(K)-1}$ or $(\flat_{n,m,l})_{l=1}^{\upsilon_{n,m}(K)-1}$ or $(\flat^{\S}_{n,m,l})_{l=1}^{\upsilon_{n,m}(K)-1}$).

\subsection{On algebraic entropies}\label{S32.18} Let $K$ be a field and $(n,s,j)\in (\mathbb N^{\ast}\setminus\{1\})^3$. Recall that for each $e\in\End_n(K)$ and every $(l,t)\in (\mathbb N^{\ast})^2$ we have $\pi(e^{l+t})\le\pi(e^l)\pi(e^t)$ and hence the limit, called the algebraic entropy\index{algebraic entropy} of $e$ (cf.\ \cite{BV}, Sect.\ 2; see also, for instance \cite{G-POS}, Subsect.\ 1.2), $\mathcal E(e):=\lim \limits_{l\to \infty} \frac{\ln\pi(e^l)}{l}$ always exists. 

By similar reasons the limits of the following definition always exist.

\begin{definition}\label{D27} Let $e\in\End_n(K)$ with $\pi(e)\ge 2$ and $a\in\GA_n(K)$. 

\medskip
{\bf (1)} By the normalized algebraic entropy of $e$ we mean the limit 
$$\mathcal E^{\n}(e):=\lim \limits_{l\to \infty} \frac{\ln\pi(e^l)}{l\ln\bigl(\pi(e)\bigr)}\in [0,1].$$ 

{\bf (2)} By the symmetric algebraic entropy\index{algebraic entropy!symmetric algebraic entropy} of $a$ we mean 
$$\mathcal E^{\symm}(a):=\lim \limits_{l\to \infty} \frac{\ln\ell(a^l)}{l}\in [0,\infty).$$

{\bf (3)} Suppose that $a\notin\AGL_n(K)$. By the normalized symmetric algebraic entropy\index{algebraic entropy!normalized symmetric algebraic entropy} of $a$ we mean the limit $\mathcal E^{\nsymm}(e):=\lim \limits_{l\to \infty} \frac{\ln\ell(a^l)}{l\ln\bigl(\ell(a)\bigr)}\in [0,1]$.

\smallskip
{\bf (4)} Suppose that $n=2$. By the minimal algebraic entropy\index{algebraic entropy!minimal symmetric algebraic entropy} of $a$ we mean the limit $\mathcal J(a):=\lim \limits_{l\to \infty} \frac{\j_{a^l}}{l}\in [0,\j_a]$. 

\smallskip
{\bf (5)} Suppose that $n=2$. If $a\notin\AGL_2(K)$, then by the normalized minimal algebraic entropy\index{algebraic entropy!normalized minimal symmetric algebraic entropy} of $a$ we mean the limit $\mathcal J^{\n}(a):=\lim \limits_{l\to \infty} \frac{\j_{a^l}}{l\j_a}\in [0,1]$.
\end{definition}

\smallskip\noindent
{\bf Problem 1.} Compute the subsets $\mathbb E_{n,s}(K):=\{\mathcal E^{\nsymm}(a)|a\in\TGA_n(K),\ell(a)=s\}$ and $\mathbb E_{2,s,j}(K):=\{\mathcal E^{\nsymm}(a)|a\in\TGA_n(K),\,\ell(a)=s,\,\j_a=j\}$ of the interval $[0,1]$.

\begin{example}\normalfont\label{EX20}
{\bf (1)} For each $c\in\AGL_n(K)$ we have $\mathcal E^{\nsymm}(cac^{-1})=\mathcal E^{\nsymm}(a)$.

\smallskip
{\bf (2)} For $a:=\e(x_2+x_1^s,x_1)\in\GA_2(K)$ we have $\ell(a)=s$ and $\mathcal E^{\nsymm}(a)=1$. Thus $1\in \mathbb E_{2,s,1}(K)$. 

\smallskip
{\bf (3)} Suppose that $s$ is not a prime number. From Proposition \ref{PR11}(4.b) applied to $j=2$ and $a=a_1a_2$ with $\j_{a_1}=\j_{a_2}=1$, $\dir(a_1)\neq\dir(a_2)$, and $\ell(a)=\ell(a_1)\ell(a_2)=s$ we get that $\ell(a^l)=\ell(a)^l$ for each $l\in\mathbb N^{\ast}$ and thus $\mathcal E^{\nsymm}(a)=1$. Thus $1\in \mathbb E_{2,s,2}(K)$. 

\smallskip
{\bf (4)} It is easy to see that $\mathbb E_{2,s}(K)\subset [0,\frac{1}{2}]\cup\{1\}$.
\end{example}

{\bf Problem 2.} Find the subsets $\mathbb J_s(K):=\{\mathcal J^{\n}(a)|a\in\GA_2(K),\ell(a)=s\}\subset [0,1]$.

One expects to solve Problem 1 for $n=2$ and Problem 2 based solely on \cite{F3}, Props.\ 4 and 6.

{\bf Problem 3.} Classify all $a\in\TGA_n(K)$ with $\ell(a)=s$ and $\mathcal E^{\nsymm}(a)=1$ and all $b\in\GA_2(K)$ with $\ell(b)=s$ and $\mathcal J^{\n}(b)=1$.

\smallskip
The next problem is a natural extrapolation of \cite{BV}, Sect.\ 5, Conj.\ 1.

\smallskip
{\bf Problem 4.} Classify all $a\in\TGA_n(K)$ with $\ell(a)=s$ and the property that the series $\sum_{l=0}^{\infty} \mathcal E^{\symm}(a^l)x^l\in\mathbb Z[[x]]$ is a rational function with integer coefficients.

\smallskip
For $n=2$ and $K=\mathbb C$, Problem 4 is solved in \cite{F3}, Subsect.\ 2.4.

The last problem is a natural extrapolation of \cite{G-POS}, Subsect.\ 3.2.

\smallskip
{\bf Problem 5.} Classify all endomorphisms $\e\in\End_n(K)$ with the property that for each $a\in\TGA_n(K)$ with $\ell(a)=s$ we have $\mathcal E(aea^{-1})=1$.

\medskip\smallskip
\hbox{Alexander Borisov,\;\;Email: aborisov@binghamton.edu}
\hbox{Address: Department of Mathematics and Statistics, Binghamton University,}
\hbox{Binghamton, P.\ O.\ Box 6000, New York 13902-6000, U.S.A.}

\medskip\smallskip
\hbox{Ofer Gabber,\;\;\;E-mail: gabber@ihes.fr}
\hbox{Address: IH\'ES, Le Bois-Marie, 35 Route de Chartres,}
\hbox{F-91440 Bures-sur-Yvette, France.}

\medskip\smallskip
\hbox{Adrian Vasiu,\;\;\;Email: avasiu@binghamton.edu}
\hbox{Address: Department of Mathematics and Statistics, Binghamton University,}
\hbox{Binghamton, P.\ O.\ Box 6000, New York 13902-6000, U.S.A.}

\newpage\section*{Notation}\label{N6}

\subsubsection*{Capital letters}

\begin{itemize}\footnotesize

\item
\hyperref[PH1]{$A$} a set

\item 
\hyperref[N1]{${A_l}$} the alternating subgroup of $S_l$

\item 
\hyperref[PH16]{$B$} a Borel subgroup of a reductive group over $\Spec K$, with an upper right index

\item 
\hyperref[L16]{$C$}, \hyperref[L16]{$C^+$}, \hyperref[PH94]{$\overline{C}$}, \hyperref[L16]{$\overline{C}^+$}, \hyperref[EXT5]{$D$}, \hyperref[C20]{$D^+$}, \hyperref[C20]{$D^-$}, and \hyperref[C20]{$\overline{D}^+$} positive constants, often with a lower right index

\item 
\hyperref[EQ12]{$E_{l,x,y}$}, \hyperref[EQ13]{$E_{l,x}$}, and \hyperref[PH99]{$E(n,|K|)$} real numbers

\item 
\hyperref[PH16m]{$F$} a polynomial, often with a lower right index

\item 
\hyperref[PH15]{$G$} a linear algebraic group over $\Spec K$, with $K$ a perfect field

\item 
\hyperref[PH15]{$G^0$} the connected component of the identity element of $G$ 

\item 
\hyperref[PH15]{$H$} a reductive group over $\Spec K$

\item 
\hyperref[PH15j]{$H^{\der}$} the derived group of $H$

\item 
\hyperref[PH15j]{$H^{\sc}$} the simply connected semisimple group cover of $H^{\der}$ 

\item 
\hyperref[D7]{$I$} a set, often with a lower right index

\item 
\hyperref[PH96]{$J$} a finite subset of either $\mathbb N$ or $\mathbb N^2$

\item 
\hyperref[PH2]{$K$} the main field

\item 
\hyperref[PH2]{$K^{\ast}$} the set or the multiplicative group of non-zero elements of a field $K$

\item 
\hyperref[PR5]{$L$} an additional field related to $K$

\item 
\hyperref[PH15-]{$M$} a smooth monoid or a group scheme over $\Spec K$

\item 
\hyperref[PH16n]{$N$} a natural number, often with a lower right index

\item 
\hyperref[PH2]{$O$}, \hyperref[PH2]{$P$}, and \hyperref[PH2]{$Q$} points in $K^n$, often with a lower right index

\item 
\hyperref[PH5]{$R$} the polynomial $K$-algebra $K[x_1,\ldots,x_n]$

\item 
\hyperref[PH15]{$R_u(G)$} the unipotent radical of a connected linear algebraic group $G$ over $\Spec K$

\item 
\hyperref[D2]{$ST_{n,m}$} a class of fields and a subclass of the class of fields $T_{n,m}$

\item 
\hyperref[N1]{$S_l$} the group of permutations of $\llbracket1,l\rrbracket$

\item 
\hyperref[D2]{$T_{n,m}$} a class of fields

\item 
\hyperref[PH16]{$U$} the unipotent radical of a Borel subgroup $B$, with an upper right index

\item 
\hyperref[PH13j]{$V$} and \hyperref[PH13j]{$W$} vector spaces or affine sets over $K$ 

\item 
\hyperref[PH97]{$W(\alpha_1,\ldots,\alpha_n)$} a Vandermonde determinant with $(\alpha_1,\ldots,\alpha_n)\in K^n$

\item 
\hyperref[PH3]{$X$} a scheme of finite type over $\Spec K$, often assumed to be reduced

\item 
\hyperref[C1]{$Y$} a finite subset of $K^n$, often with a lower right index

\item 
\hyperref[C1]{$Z$} a finite subset of $K^n$, often with a lower right index

\item 
\hyperref[PH15j]{$Z_H$} the kernel of the central isogeny $H^{\sc}\rightarrow H^{\der}$

\end{itemize}

\subsubsection*{Capital letters symbols} 

\begin{itemize}\footnotesize

\item 
\hyperref[N5]{$(P)$} a $K$-valued point of $\mathbb P^{n-1}_K$, often with a lower right index

\item 
\hyperref[PH1]{$|A|$} the number of elements of a set $A$; so $|A|\in\mathbb N\cup\{\infty\}$

\item 
\hyperref[PH9]{$\AGL_n(K)$} the subgroup of $\GA_n(K)$ of affine automorphisms

\item 
\hyperref[PH10]{$\ASL_n(K)$} the subgroup of $\GA_n(K)$ or $\AGL_n(K)$ of special affine automorphisms

\item 
\hyperref[PH1]{$\Alt(A)$} the subgroup of $\perm(K^n)$ formed by even permutations of a finite set $A$

\item 
\hyperref[PH85]{$\C_1$} the real number $2^{1-\epsilon_{\z}}=0.60346...$

\item 
\hyperref[PH85]{$\C_2$} the real number $0.31817...$ equal to $\frac{-1}{\epsilon_{\z}\zeta^{\prime}(\epsilon_{\z})}$

\item 
\hyperref[PH5]{$\End_n(K)$} the multiplicative monoid of endomorphisms of $\mathbb A^n_K$

\item 
\hyperref[D2.5]{$\Fix_H(Y)$} the subgroup of the group of automorphism $H$ that fixes each point of $Y$

\item 
\hyperref[PH2]{$\GA_n(K)$} or $\GA(\mathbb A^n_K)$ the group of automorphisms of $\mathbb A^n_K$

\item 
\hyperref[PH4a]{$\GGA_n(K)$} or $\GGA(\mathbb A^n_K)$ the normal subgroup of $\GA_n(K)$ that contains $\SGA_n(K)$

\item 
\hyperref[PH10a]{$\GL_n(K)$} the general linear group of invertible $n\times n$ matrices with entries in $K$ identified with the subgroup of $\GA_n(K)$ or $\AGL_n(K)$ formed by linear automorphisms 

\item 
\hyperref[D4]{$\L_Y$} an invariant of a subset $Y$ of $K^n$ that does not contain the zero vector

\item 
\hyperref[D6]{$\L_k$} an arithmetic function with $k\in\mathbb N^{\ast}\setminus\{1\}$

\item 
\hyperref[PH1]{$\perm(A)$} the group of permutations of a set $A$

\item 
\hyperref[S18]{$\R_{L/K}$} the Weil restriction functor, with $L$ a finite field extension of $K$

\item 
\hyperref[PH4a]{$\SGA_n(K)$} or $\SGA(\mathbb A^n_K)$ the group of special automorphisms of $\mathbb A^n_K$

\item 
\hyperref[PH10a]{$\SL_n(K)$} the special linear group of $n\times n$ matrices of determinant $1$ with entries in $K$ identified with the subgroup of $\GA_n(K)$ formed by linear automorphisms of Jacobian determinant $1$

\item 
\hyperref[PH11]{$\STGA_n(K)$} the intersection $\TGA_n(K)\cap\SGA_n(K)$ or the group of special tame automorphisms of $\mathbb A^n_K$ or the group of tame autommorphisms of $\mathbb A^n_K$ of Jacobian determinant $1$

\item 
\hyperref[D24]{$\Stab_{\GA_n(K)}(Y)$} and $\Stab_{\TGA_n(K)}(Y)$ stabilizers subgroups of a non-empty finite subset $Y$ of $K^n$

\item 
\hyperref[PH10b]{$\TGA_n(K)$} the subgroup of $\GA_n(K)$ of tame automorphisms

\item 
\hyperref[D10.2]{$\Perm(K^n)$} the group of tame permutations of $K^n$, with $K$ a finite field

\item 
\hyperref[D4]{$\langle Y_{\aff}\rangle$} the affine span of $Y$ over $K$

\item
\hyperref[PH92]{$\underline{O}$}, \hyperref[PH3+]{$\underline{P}$}, and \hyperref[PH3]{$\underline{Q}$} tuples of points in $K^n$

\end{itemize}

\subsubsection*{Small letters}

\begin{itemize}\footnotesize

\item 
\hyperref[EQ3]{$a$}, \hyperref[EQ3]{$b$}, and \hyperref[PH12b]{$c$} automorphisms, often with a lower right index

\item 
\hyperref[D6]{$d$} a natural number for affine dimensions, often with a lower right index 

\item 
\hyperref[PH6]{$e$} an endomorphism

\item 
\hyperref[PH6]{$f$} and \hyperref[PH16z]{$g$} polynomials in some $R=K[x_1,\ldots,x_n]$ 

\item 
\hyperref[PR9]{$h$} a function from a set $Y$ to $K$, polynomial in some $R$, or valued point of some monoid

\item 
\hyperref[PH3]{$i$} and \hyperref[PH3]{$j$} indices of natural numbers

\item 
\hyperref[L3]{$k$} a natural number related to the cardinality of a finite field $K$

\item 
\hyperref[N1]{$l$} a natural numbers, often the index for $S_l$ and $A_l$

\item 
\hyperref[PH3]{$m$} the usual notation for the number of points

\item 
\hyperref[D5]{$m_{i,Y,a}\in\mathbb N$} an invariant related to $Y\subset K^n$, $a\in\AGL_n(K)$, and $i\in \llbracket0,n\rrbracket$

\item 
\hyperref[PH4]{$n$} the usual dimension of the affine spaces

\item 
\hyperref[PH16]{$o(H)$} the order of the kernel of the central isogeny $H^{\sc}\rightarrow H^{\der}$

\item 
\hyperref[N1]{$o(\sigma)$} the order of a permutation $\sigma$

\item 
\hyperref[PH12a]{$p$} a prime in $\mathbb N$

\item 
\hyperref[PH12a]{$q$} a natural number used for exponents such as for finite fields with $p^q$ elements, for integral parts of quotients, and, often with an index, as cardinalities of certain fibers

\item 
\hyperref[EQ26]{$q_k(n)$} a natural number with $(q,n)\in (\mathbb N^{\ast}\setminus\{1\})^2$

\item 
\hyperref[PH1]{$r$} a natural number, often for remainders in divisions and for numbers of disjoint $s$-cycles

\item 
\hyperref[PH1]{$s$} a natural number for the length of cycles or exponents or the second coordinate of the tame types, often with a lower right index

\item 
\hyperref[F1]{$t$} a natural number index

\item 
\hyperref[D10]{$v$} and \hyperref[D10]{$w$} vectors

\item 
\hyperref[PH5]{$x$} an indeterminate for polynomials, often with a lower right index

\item 
\hyperref[PH16y]{$y$} and \hyperref[PH91]{$z$} indeterminates for real-valued functions

\end{itemize}

\subsubsection*{Small letters symbols} 

\begin{itemize}\footnotesize

\item 
\hyperref[N1]{$\c(\sigma)$} the number of cycles of a permutation $\sigma$

\item 
\hyperref[D4]{$\c_Y$} the maximum number of collinear points in $Y$

\item 
\hyperref[N1]{$\c_s(\sigma)$} the number of $s$-cycles of a permutation $\sigma$

\item 
\hyperref[N1]{$\c_{\textup{even}}(\sigma)$} the number of cycles of even length of a permutation $\sigma$

\item 
\hyperref[N1]{$\c_{\textup{odd}}(\sigma)$} the number of cycles of odd length of a permutation $\sigma$

\item 
\hyperref[S4]{$\char(K)$} the characteristic of the field $K$ 

\item 
\hyperref[D4]{$\d_Y$} the affine dimension of the affine span of $Y\subset K^n$

\item 
\hyperref[PH5]{$\deg(f)$} the degree of a polynomial $f\in R$

\item 
\hyperref[PH90]{$\det$} the determinant function for square matrices and for the Jacobian determinant surjective homomorphism $\det:\GA_n(K)\rightarrow K^{\ast}$

\item 
\hyperref[D12]{$\dir(a)\in\mathbb P_K^1(K)$} the direction of a non-affine triangular automorphism $a\in\GA_2(K)$

\item 
\hyperref[D19]{$\dir_Y$} an invariant of the finite subset $Y$ with at least $2$ elements

\item  
\hyperref[PH6]{$\e(f_1,\ldots,f_n)$} the element of $\End_n(K)$ with $e^{\#}(x_i)=f_i$ for each $i\in\llbracket1,n\rrbracket$

\item 
\hyperref[D13]{$\j_a$} the Furter's length of $a\in\GA_2(K)$

\item 
\hyperref[D15]{$\j_{K}$} the Furter's length of a finite field $K$

\item 
\hyperref[D15]{$\j_{\sigma}$} the Furter's length of $\sigma\in\Perm(K^2)$

\item 
\hyperref[L5]{$\l=\l_{r,s}^A$} a function from $\mathcal F_{r,s}(A)$ to $\llbracket0,s\rrbracket$

\item 
\hyperref[D6]{$\l_K$}, \hyperref[D6]{$\l_K^{[d]}$}, and \hyperref[D6]{$\l_K^{[\le d]}$} arithmetic functions with $K$ a field

\item 
\hyperref[D4]{$\l_Y$} an invariant of a subset $Y$ of $K^n$

\item 
\hyperref[D5]{$\mq_Y$} the multisecant number of $Y$

\item 
\hyperref[D5]{$\mq_Y^P$} the multisecant number of $Y$ at $P$

\item 
\hyperref[D4]{$\overline{m}_{i,Y,P,a}\in\mathbb N$} an invariant related to $P\in Y\subset K^n$, $a\in\AGL_n(K)$, and $i\in \llbracket0,n\rrbracket$

\item 
\hyperref[N1]{$\n(\sigma)$} the cardinality of the support of a permutation $\sigma$

\item 
\hyperref[D5]{$\q_Y$} the secant number of $Y$

\item 
\hyperref[D5]{$\q_Y^P$} the secant number of $Y$ at $P$

\item 
\hyperref[D8]{$\r_Y$} the representation degree of $Y$

\item 
\hyperref[D8]{$\r_Y^P$} the representation degree of $Y$ at $P$

 \item 
\hyperref[D5]{$\s_Y$} the capacity of $Y$, i.e., the minimum of all the $\s_a(Y)$s
 
\item 
\hyperref[D5]{$\s_a(Y)$} an invariant related to $a\in\AGL_n(K)$ and $Y\subset K^n$

\item 
\hyperref[D4]{$\overline{\s}_Y$} the strict capacity of $Y$

\item 
\hyperref[D4]{$\overline{\s}_Y^P$} the strict capacity of $Y$ at $P\in Y$

\item 
\hyperref[D5]{$\overline{\s}_a^P(Y)$} an invariant related to $a\in\AGL_n(K)$ and $P\in Y\subset K^n$

\item 
\hyperref[D9]{$\overline{\s}_k$} the modulo $k$ strict capacity function

\item 
\hyperref[D10]{$\shift_{V\oplus W}^{v_0+W_0}$} the selective shift permutation in $\perm(K^n)$ with $K$ a finite field, a direct sum decomposition $K^n=V\oplus W$, $v_0\in V\setminus\{0\}$, and a non-empty subset $W_0$ of $W$

\item 
\hyperref[D10]{$\shift_{V\oplus W}^{v_0+w_0}$} the selective shift permutation when $W_0=\{w_0\}$
 
\item 
\hyperref[N1]{$\supp(\sigma)$} the support of a permutation $\sigma$

\item 
\hyperref[D20]{$\t_Y$} an invariant of the finite subset $Y$ of $K^n$

\item 
\hyperref[D20]{$\v_Y$} an invariant of the finite subset $Y$ of $K^n$

\item 
\hyperref[D20]{$\v_{m,|K|}$} the invariant of $K$ and $m\in\mathbb N^{\ast}\setminus\{1\}$ which is $1$ if $K$ is infinite and it is the smallest $\v$ such that $|K|^{\v}\ge m$ if $K$ is finite

\item 
\hyperref[D19]{$\w_Y$} an invariant of a finite subset $Y$ of $K^n$ with at least $2$ elements

\end{itemize}

\subsubsection*{Capital Greek letters}

\begin{itemize}\footnotesize

\item 
\hyperref[EXT4]{$\Delta$} the polynomial in three variables related to Nagata's automorphism

\item 
\hyperref[L22]{$\Delta_r$} and \hyperref[PR35]{$\Delta_{r,\sigma}$} square matrices of size $r\times r$ with entries in $K$

\item 
\hyperref[PH8a]{$\Gamma$} and \hyperref[PH8a]{$\Gamma[l]$} subsets of $\GA_n(K)$ or abstract groups

\item 
\hyperref[S4]{$\Lambda$} a left action of $G$ on $X$

\item 
\hyperref[N4]{$\Omega(\nabla)$} the maximum of the set $\{\ell_{n,K}(\sigma)|\sigma\in\nabla\}$ 

\item 
\hyperref[N4]{$\Omega^{\S}(\nabla)$} the maximum of the set $\{\ell^{\S}_{n,K}(\sigma)|\sigma\in\nabla\}$ 

\item 
\hyperref[PH14]{$\Pi_{n,s\textup{-cycle}}(K)$} a finite subset of $\mathbb N$ with $K$ a finite field and $s\in \mathbb N^{\ast}\setminus\{1\}$, $s$ being odd if $4\mid |K|$

\item 
\hyperref[N4]{$\Pi^{\S}_{n,s\textup{-cycle}}(K)$} a finite subset of $\mathbb N$ with $K$ a finite field and $s\in\mathbb N^{\ast}\setminus\{1\}$ odd

\item 
\hyperref[S4]{$\Theta$} an action of $M$ on $X$

\item 
\hyperref[D21]{$\Xi_{n,m,q}(K)$} an invariant of $K$

\item 
\hyperref[D21]{$\Xi^{\S}_{n,m,q}(K)$} an invariant of $K$

\end{itemize}

\subsubsection*{Small Greek letters}

\begin{itemize}\footnotesize 

\item 
\hyperref[EX3]{$\alpha$} an element of $K$, often with a lower right index

\item 
\hyperref[EX3]{$\beta$} an element of $K$, often with a lower right index

\item 
\hyperref[EX5]{$\chi^P_Y$} a polynomial in $K[x_1,\ldots,x_n]$ with $P\in Y\subset K^n$

\item 
\hyperref[EX6]{$\delta$} an element of $K$, often with a lower right index

\item 
\hyperref[D3]{$\digamma_{r,\sigma}$} a natural number defined for $r\in\mathbb N^{\ast}\setminus\{1\}$ and $\sigma\in S_l$

\item 
\hyperref[D3]{$\digamma_{r,l}$} a natural number defined for $(r,l)\in (\mathbb N^{\ast}\setminus\{1\})\times\mathbb N^{\ast}$

\item 
\hyperref[C2]{$\epsilon$} a small non-negative integer 

\item 
\hyperref[PH85]{$\epsilon_{\z}$} the real number $1.72864...$ such that $\zeta(\epsilon_{\z})=2$

\item 
\hyperref[L19]{$\epsilon_{n,k}$} a real number with $(n,k)\in (\mathbb N^{\ast}\setminus\{1\})^2$

\item 
\hyperref[PH98]{$\eta(i)$} a natural number which is the degree of a polynomial

\item 
\hyperref[EX3]{$\gamma$} an element of $K$, often with a lower right index

\item 
\hyperref[PH14e]{$\iota$} a small positive integer

\item 
\hyperref[PR16]{$\kappa_{n,K}=\kappa_{n,K;2}^{\d=1}$} the degree length of a (each) transposition in $\perm(K^n)$ when $4\nmid |K|$

\item 
\hyperref[F9]{$\kappa_{n,K;3}^{\d=1}$} the degree length of any $3$-cycle in $\perm(K^n)$ with collinear support when $|K|=3$

\item 
\hyperref[F9]{$\kappa_{n,K;4}^{\d=2}$} the degree length of any $4$-cycle in $\perm(K^n)$ with affinely dependent support when $|K|=2$ and $n\ge 2$

\item 
\hyperref[F9]{$\kappa_{n,K;s}^{\d=s-1}$} the degree length of any $s$-cycle in $\perm(K^n)$ with affinely independent support, with $K$ a finite field and $s\in\llbracket1,n+1\rrbracket$, $s$ being odd if $4\mid |K|$

\item 
\hyperref[PH14j]{$\lambda$} a line in $K^n$ or a $K$-valued point of $\mathbb P^{n-1}_K$

\item 
\hyperref[N1]{$\nu^+_{s,r}(\sigma)$} the smallest number of factors required to write $\sigma$ as a product of permutations that are products of precisely $r$ disjoint $s$-cycles, when $\sigma$ is even, $s$ is odd, and $r\ge 2$

\item 
\hyperref[L3]{$\nu_{k,1}(k^2)$} the maximum of  the set $\{\nu_{k,1}(\sigma)|\sigma\in A_{k^2}\}$, with $k\in 1+2\mathbb N^{\ast}$

\item 
\hyperref[N1]{$\nu_{s,r}(\sigma)$} the smallest number of factors required to write $\sigma$ as a product of permutations that are products of at most $r$ disjoint $s$-cycles,  with $(r,s)\in\mathbb N^{\ast}\times(\mathbb N^{\ast}\setminus\{1\})$ and $\sigma$ even if $s$ is odd

\item 
\hyperref[C4]{$\nu^+_{s,r}(l)$} the maximum of the $\nu^+_{s,r}(\sigma)$s with $\sigma\in A_l$

\item 
\hyperref[C4]{$\nu^+_{s,r}(l-\cycle)$} the maximum of the $\nu^+_{s,r}(\sigma)$s with $\sigma$ a disjoint product of a transposition and a cycle of even length 

\item 
\hyperref[N1]{$\nu^{\textup{even}}_{s,r}(\sigma)$} the smallest number of factors required to write $\sigma$ as a product of permutations in $A_l$ that are products of at most $r$ disjoint $s$-cycles, with $(r,s)\in\mathbb N^{\ast}\times (2\mathbb N^{\ast})$ and $\sigma\in A_l$

\item 
\hyperref[N2]{$\overline{\nu}^{(l)}_{s,r}(\underline{\c})$} the number of products of at most $r$ disjoint $s$-cycles in $S_l$ required to write a $\theta$ with $\c_i(\theta)=\c_i$ if $i\in\llbracket2,s-1\rrbracket$ and $\c_i(\theta)=0$ if $i\ge s$, with $\underline{\c}=(\c_2,\ldots,c_{s-1})\in\mathbb N^{s-2}$, and with $l\ge\n(\theta)=\sum_{i=2}^{s-1} i\c_i$

\item 
\hyperref[N2]{$\overline{\nu}_{s,r}(\sigma)$} an invariant of a permutation $\sigma$ used to estimate from above $\nu_{s,r}(\sigma)$

\item 
\hyperref[N4]{$\omega(\nabla)$} the minimum of the set $\{\ell_{n,K}(\sigma)|\sigma\in\nabla\}$ 

\item 
\hyperref[N4]{$\omega^{\S}(\nabla)$} the minimum of the set $\{\ell^{\S}_{n,K}(\sigma)|\sigma\in\nabla\}$

\item 
\hyperref[D5]{$\pi$} a linear projection between affine spaces, often with a lower right index

\item 
\hyperref[PH6]{$\pi(a)$} the degree of an automorphism $a$ of $\mathbb A^n_K$

\item 
\hyperref[PH7]{$\pi(e)$} or \hyperref[PH7]{$\pi(f_1,\ldots,f_n)$} the maximum degree of the polynomials (that define $e$)

\item 
\hyperref[D2]{$\pi_{\underline{P},\underline{Q}}$} the degree length required to move $\underline{P}$ to $\underline{Q}$ via tame automorphisms

\item 
\hyperref[D14]{$\pi^{\E}_{\underline{P},\underline{Q}}$} the degree length required to move $\underline{P}$ to $\underline{Q}$ via tame automorphisms that represent even permutations of $\perm(K^n)$, with $K$ a finite field and $(\underline{P},\underline{Q})\in\mathbb D_{n,m}(K)^2$

\item 
\hyperref[D2]{$\pi^{\S}_{\underline{P},\underline{Q}}$} the degree length required to move $\underline{P}$ to $\underline{Q}$ via special tame automorphisms

\item 
\hyperref[D17]{$\pi^{\S,\le j}_{\underline{P},\underline{Q}}$} the degree length required to move $\underline{P}$ to $\underline{Q}$ via special tame automorphisms $a$ of $\mathbb A^2_K$ with $\j_a\le j$

\item 
\hyperref[D17]{$\pi^{\le j}_{\underline{P},\underline{Q}}$} the degree length required to move $\underline{P}$ to $\underline{Q}$ via tame automorphisms $a$ of $\mathbb A^2_K$  with $\j_a\le j$

\item 
\hyperref[D2.5]{$\pi_{Y,s}(K)$} an invariant of $Y$

\item 
\hyperref[D2.5]{$\pi^{\S}_{Y,s}(K)$} an invariant of $Y$

\item 
\hyperref[D2]{$\pi_{n,m}(K)$} a main invariant of $K$

\item 
\hyperref[D14]{$\pi^{\E}_{n,m}(K)$} an invariant of $K$, the supremum of the set formed by the positive integers $\pi^{\E}_{\underline{P},\underline{Q}}$ with $(\underline{P},\underline{Q})\in\mathbb D_{n,m}(K)^2$

\item
\hyperref[D2]{$\pi^{\S}_{n,m}(K)$} a main invariant of $K$

\item 
\hyperref[S32.8]{$\pi_{n,m}^{-,\S}(K)$} an invariant of $K$

\item 
\hyperref[S32.8]{$\pi_{n,m}^{-}(K)$} an invariant of $K$

\item 
\hyperref[D19]{$\psi_Y$} a function from $\{(i,j)\in\llbracket1,m\rrbracket|i<j\}$ to $\mathbb P_K^{n-1}(K)$ attached to a finite set $Y$ with at least $2$ elements

\item 
\hyperref[N1]{$\sigma$} a permutation, often with a lower right index

\item 
\hyperref[N2]{$\sigma_s^{(t)}$} a permutation associated recursively to $\sigma$, an integer $s\ge 4$, and $t\in\mathbb N$, with $\sigma_s^{(0)}=\sigma$

\item 
\hyperref[PH14f]{$\tau$} a transposition, often with a lower right index

\item 
\hyperref[N1]{$\theta$} a permutation, often with a lower right index

\item 
\hyperref[D24]{$\upsilon_{Y,s}(K)$} an invariant of the finite non-empty subset $Y$ of $K^n$ and $s\in\mathbb N^{\ast}$ when we have $\sqrt[n]{|Y|+s}\le |K|$ and, if $K$ is finite with $4\mid |K|$, $\sqrt[n]{|Y|+s+2}\le |K|$

\item 
\hyperref[D24]{$\upsilon^{\S}_{Y,s}(K)$} an invariant of the finite non-empty subset $Y$ of $K^n$ and $s\in\mathbb N^{\ast}$ when we have $\sqrt[n]{|Y|+s+2}\le |K|$

\item 
\hyperref[C1]{$\upsilon_{n,m}(K)$} a main invariant of $K$

\item 
\hyperref[C1]{$\upsilon^{\S}_{n,m}(K)$} a main invariant of $K$

\item 
\hyperref[S32.8]{$\upsilon_{n,m}^{-,\S}(K)$} an invariant of $K$

\item 
\hyperref[S32.8]{$\upsilon_{n,m}^{-}(K)$} an invariant of $K$

\item 
\hyperref[EXT3]{$\varepsilon$} a small real number

\item 
\hyperref[L19]{$\varepsilon_{n,k,1}$} and \hyperref[L16]{$\varepsilon_{n,k,2}$} real numbers indexed by $(n,k)\in (\mathbb N^{\ast}\setminus\{1\})^2$

\item 
\hyperref[S4]{$\varrho_{\Theta}$} the monoid homomorphism associated to an action $\Theta$

\item 
\hyperref[PH89]{$\varrho_{n,K,N}$} the cardinality of the image of $\TGA_n(K)[N]$ under $\varrho_{n,K}$ when $K$ is a finite field and $N\in\mathbb N^{\ast}\setminus\{1\}$

\item 
\hyperref[PH89]{$\varrho^{\S}_{n,K,N}$} the cardinality of the image of $\STGA_n(K)[N]$ under $\varrho_{n,K}$ when $K$ is a finite field and $N\in\mathbb N^{\ast}\setminus\{1\}$

\item 
\hyperref[PH12]{$\varrho_{n,K}$} the permutation representation of $\GA_n(K)$ on $K^n$

\item 
\hyperref[F1]{$\varsigma$} a permutation, often with a lower right index

\item 
\hyperref[PH14i]{$\vartheta$} a permutation, often with a lower right index

\item 
\hyperref[D24]{$\xi_{n,m,q}(K)$} an invariant of $K$

\item 
\hyperref[D24]{$\xi^{\S}_{n,m,q}(K)$} an invariant of $K$

\end{itemize}

\subsubsection*{Blackboard bold letters}

\begin{itemize}\footnotesize 

\item 
\hyperref[PH2]{$\mathbb A^r_K$} the affine space of dimension $r\in\mathbb N$ over $\Spec K$

\item 
\hyperref[PH2]{$\mathbb D_{1,m}(K^{\ast})$} or $\mathbb D_m(\mathbb G_{\m,K})$ the set of $m$-tuples of distinct points in $K^{\ast}$

\item 
\hyperref[PH2]{$\mathbb D_{n,m}(K)$} or $\mathbb D_m(\mathbb A^n_K)\subset (K^n)^m$ the set of $m$-tuples of distinct points in $K^n=\mathbb A^n_K(K)$

\item 
\hyperref[N3]{$\mathbb D^{\d=i}_{n,m}(K)$} the set of $m$-tuples of distinct points in $K^n$ whose affine dimension is $i$, with $i\in\llbracket1,\min(n,m-1)\rrbracket$

\item 
\hyperref[N3]{$\mathbb D^{\d\ge i}_{n,m}(K)$} the set of $m$-tuples of distinct points in $K^n$ whose affine dimension is at least $i$, with $i\in\llbracket1,\min(n,m-1)-1\rrbracket$

\item 
\hyperref[N3]{$\mathbb D^{\d\le i}_{n,m}(K)$} the set of $m$-tuples of distinct points in $K^n$ whose affine dimension is at most $i$, with $i\in\llbracket2,\min(n,m-1)\rrbracket$

\item 
\hyperref[S28]{$\mathbb D_{n,m}^{\sum=0}(K)$} the set of $m$-tuples of distinct points in $K^n$ that add up to $0$

\item 
\hyperref[S28]{$\mathbb D_{n,m}^{\sum\neq 0}(K)$} the set of $m$-tuples of distinct points in $K^n$ that add up to a non-zero point

\item 
\hyperref[PH12a]{$\mathbb F_{p^q}$} a finite field with $p^q$ elements

\item 
\hyperref[PH3a]{$\mathbb G_{\a,K}$} the split connected unipotent group over $\Spec K$ of dimension $1$

\item 
\hyperref[PH3a]{$\mathbb G_{\m,K}$} the split multiplicative group over $\Spec K$ of dimension $1$

\item 
\hyperref[S32.6]{$\mathbb I_{n,K}$} a function from $\mathcal P_{\ge 2}(K^n)$ to $\mathbb N^{14}$ that encodes the fourteen invariants introduced

\item 
\hyperref[PH4]{$\mathbb T_{n,m}(K)$} or $\mathbb T_m(\mathbb A^n_K)$ the left action of $\GA_n(K)$ on $\mathbb D_{n,m}(K)$

\end{itemize}

\subsubsection*{Calligraphic letters}

\begin{itemize}\footnotesize 

\item 
\hyperref[EXT2]{$\mathcal C$} a conjugacy class in either $A_l$ or $S_l$

\item 
\hyperref[S32.18]{$\mathcal E(e)$} the algebraic entropy of $e\in\End_n(K)$

\item 
\hyperref[D27]{$\mathcal E^{\nsymm}(a)$} the normalized symmetric algebraic entropy of $a\in\GA_n(K)\setminus\AGL_n(K)$

\item 
\hyperref[D27]{$\mathcal E^{\n}(e)$} the normalized algebraic entropy of $e\in\End_n(K)$ with $\pi(e)\ge 2$

\item 
\hyperref[D27]{$\mathcal E^{\symm}(a)$} the symmetric algebraic entropy of $a\in\GA_n(K)$

\item 
\hyperref[D7]{$\mathcal F_{r,s}(A)$} the set whose elements are $r$ pairwise disjoint finite subsets of the set $A$ whose cardinalities add up to $s$

\item 
\hyperref[D7]{$\mathcal I$} and \hyperref[D7]{$\mathcal J$} sets of subsets, often with an upper right index

\item 
\hyperref[S14]{$\mathcal N=\mathcal N_K$} the Nagata automorphism of $\mathbb A^3_K$

\item 
\hyperref[L27]{$\mathcal P_m(A)$} the set of subsets of $A$ of cardinality $m$

\item 
\hyperref[S32.6]{$\mathcal P_{\ge 2}(K)$} the set of finite subsets of $K$ of cardinality at least $2$

\item 
\hyperref[D5]{$\mathcal Q^P_Y$} the set of minimal $n$-tuple secants of $Y$ at $P$

\item 
\hyperref[EXT1]{$\mathcal S$} the solution set of a system of equations, often with a lower right index

\item 
\hyperref[P12]{$\mathcal U$} a union of conjugacy classes of $\perm(K^n)$ or $\Perm(K^n)$

\item 
\hyperref[N4]{$\mathcal Y_s$} the conjugacy class $1\mathcal Y_s$ in $\perm(K^n)$; it has subsets $\mathcal Y_s^{\d\le s-2}$, $\mathcal Y_s^{\d\ge 2}$, $\mathcal Y_s^{\d=s-1}$ if $s\le n+1$, and $\mathcal Y_s^{\d=1}$ if $s\le |K|$, with $s$ and $|K|$ subject to certain restrictions (inequalities) 

\item 
\hyperref[N4]{$t\mathcal Y_s$} the conjugacy class in $\perm(K^n)$ of products of $t$ disjoint $s$-cycles when $K$ is finite, $s\in\mathbb N^{\ast}\setminus\{1\}$, and $t\in\mathbb N^{\ast}$ with $st\le |K|^n$

\end{itemize}

\subsubsection*{Special characters}

\begin{itemize}\footnotesize 

\item 
\hyperref[PH8]{$1_{\mathbb A^n_K}$} the identity automorphism of $\mathbb A^n_K$

\item 
\hyperref[EQ1]{$\ell(a)$} the degree length of an automorphism $a$ of $\mathbb A^n_K$

\item 
\hyperref[PH8]{$\ell_{\GA_n(K)}$} the length function on $\GA_n(K)$ with values in $\mathbb N^{\ast}$

\item 
\hyperref[PH16d]{$\ell_{n,K}$} the length function on $\varrho_{n,K}\bigl(\TGA_n(K)\bigr)$ with values in $\mathbb N^{\ast}$

\item
\hyperref[D18]{$\eth(\nabla)$} a positive real number associated to a non-empty subset $\nabla$ of $\Perm(K^n)$

\item
\hyperref[D18]{$\eth^{\S}(\nabla)$} a positive real number associated to a non-empty subset $\nabla$ of $\Alt(K^n)$

\item
\hyperref[D25]{$\flat^{\S}_{n,m,l}(K)$} an invariant of $K$

\item
\hyperref[D25]{$\flat_{n,m,l}(K)$} an invariant of $K$

\item
\hyperref[PH14e]{$\hbar$} a function between finite subsets of $\mathbb N$

\item 
\hyperref[PH1]{$\llbracket r,s\rrbracket$} the set of integers $i$ with $r\le i\le s$

\item 
\hyperref[PH86]{$\mathfrak C$} the subset $\{\e(\alpha x_1+x_2,x_1)|\alpha\in K\}$ of $\GA_2(K)$

\item 
\hyperref[PH17]{$\mathfrak I$} the ideal $(x_1^{|K|-1}-x_1,\ldots,x_n^{|K|-1}-x_n)$ of $K[x_1,\ldots,x_n]$ when $K$ is finite

\item 
\hyperref[PH101]{$\mathfrak c_k$} real constants indexed by $k\in\mathbb N$ with $k\ge 4$

\item 
\hyperref[PH102]{$\mathfrak f$} a real-valued function with domain contained in $\mathbb R$, often with a lower right index

\item 
\hyperref[D26]{$\mathfrak g(m,q)$} a natural number greater than or equal to $2$

\item 
\hyperref[D26]{$\mathfrak h(m,q)$} a natural number greater than or equal to $2$

\item 
\hyperref[N1]{$\nabla$} a set of permutations

\item 
\hyperref[N1]{$\nabla^n$} the set of products of $n$ permutations in $\nabla$ with $n\in\mathbb N$

\item 
\hyperref[D22]{$\natural_{n,m,l}(K)$} an invariant of $K$

\item 
\hyperref[D22]{$\natural^{\S}_{n,m,l}(K)$} an invariant of $K$

\item 
\hyperref[D3]{$\wp_{r,l}$} the natural number (3-deviation exponent) defined for $(r,l)\in (\mathbb N^{\ast}\setminus\{1\})\times\mathbb N^{\ast}$

\end{itemize}

\printindex


\newpage\begin{thebibliography}{MMMMM}

\bibitem[AM]{AM} S.\ S.\ Abhyankar and T.\ T.\ Moh,
\sl Embeddings of the line in the plane,
\rm J.\ Reine Angew.\ Math.\ 276 (1975), 148--166

\bibitem[AFKKZ]{AFKKZ} I.\ Arzhantsev, H.\ Flenner, S.\ Kaliman, F.\ Kutzschebauch, M.\ Zaidenberg,
\sl Flexible varieties and automorphism groups,
\rm Duke Math.\ J.\ {\bf 162} (2013), no.\ 4, 767--823

\bibitem[BBS]{BBS} L.\ Babai, R.\ Beals, and Á.\ Seress,
\sl On the diameter of the symmetric group: polynomial bounds,
\rm Proceedings of the Fifteenth Annual ACM-SIAM Symposium on Discrete Algorithms, Association for Computing Machinery (ACM), New York, 2004, 1108--1112

\bibitem[BGHHSS]{BGHHSS} J.\ Bamberg, N.\ Gill, T.\ P.\ Hayes, H.\ A.\ Helfgott, Á.\ Seress, and P.\ Spiga,
\sl Bounds on the diameter of Cayley graphs of the symmetric group,
\rm J.\ Algebraic Combin.\ {\bf 40} (2014), no.\ 1, 1--22

\bibitem[Bar]{Bar} V.\ G.\ Bardakov,
\sl Factorization of even permutations into two factors of given cyclic structure,
\rm Diskret.\ Mat.\ {\bf 5} (1993), no.\ 1, 70--90; translation in Discrete Math.\ Appl.\ {\bf 3} (1993), no.\ 4, 385--406

\bibitem[BCW]{BCW} H.\ Bass, E.\ H.\ Connell, and D.\ Wright,
\sl The Jacobian conjecture: reduction of degree and formal expansion of the inverse,
\rm Bull.\ Amer.\ Math.\ Soc.\ (N.S.) {\bf 7} (1982), no.\ 2, 287--330

\bibitem[BV]{BV} M.\ P.\ Bellon and C.\ M.\ Viallet,
\sl Algebraic entropy,
\rm Comm.\ Math.\ Phys.\ {\bf 204} (1999), no.\ 2, 425--437

\bibitem[Ber]{Ber} E.\ Bertram,
\sl Even permutations as a product of two conjugate cycles,
\rm J.\ Combinatorial Theory Ser.\ A {\bf 12} (1972), 368--380

\bibitem[BH]{BH} E.\ Bertram and M.\ Herzog,
\sl Powers of cycle-classes in symmetric groups,
\rm J.\ Combin.\ Theory Ser.\ A {\bf 94} (2001), no.\ 1, 87--99

\bibitem[Bor]{Bor} A.\ Borel,
\sl Linear algebraic groups. Second edition,
\rm Grad.\ Texts in Math., Vol.\ {\bf 126}, Springer-Verlag, New York, 1991

\bibitem[BT]{BT} A.\ Borel and J.\ Tits,
\sl Groupes réductifs,
\rm Inst.\ Hautes Études Sci.\ Publ.\ Math.\ No.\ {\bf 27} (1965), 55--150

\bibitem[Bre]{Bre} J.\ L.\ Brenner,
\sl Covering theorems for FINASIGs. VIII--Almost all conjugacy classes in $\mathcal A_n$ have exponent $\le 4$,
\rm J.\ Austral.\ Math.\ Soc.\ Ser.\ A {\bf 25} (1978), no.\ 2, 210--214

\bibitem[BE]{BE} J.\ L.\ Brenner and R.\ J.\ Evans,
\sl Even permutations as a product of two elements of order five,
\rm J.\ Combin.\ Theory Ser.\ A {\bf 45} (1987), no.\ 2, 196--206

\bibitem[BR1]{BR1} J.\ L.\ Brenner and J.\ Riddell,
\sl Covering theorems for finite nonabelian simple groups. VII. Asymptotics in the alternating groups,
\rm Ars Combin.\ {\bf 1} (1976), no.\ 1, 77--108

\bibitem[BR2]{BR2} J.\ L.\ Brenner and J.\ Riddell,
\sl Research Problems: Noncanonical factorization of a permutation, 
\rm Amer.\ Math.\ Monthly {\bf 84} (1977), no.\ 1, 39--40

\bibitem[Ch]{Ch} I.\ M.\ Chiswell,
\sl Abstract length functions in groups,
\rm Math.\ Proc.\ Cambridge Philos.\ Soc.\ {\bf 80} (1976), no.\ 3, 451--463

\bibitem[Co]{Co} B.\ Conrad,
\sl The structure of solvable groups over general fields,
\rm Panor.\ Synthèses {\bf 46} [Panoramas and Syntheses], Société Mathématique de France, Paris, 2015, 159--192

\bibitem[CGP]{CGP} B.\ Conrad, O.\ Gabber, and G.\ Prasad,
\sl Pseudo-reductive groups. Second edition,
\rm New Math.\ Monogr., {\bf 26}, Cambridge University Press, Cambridge, 2015

\bibitem[DG]{DG} M.\ Demazure and P.\ Gabriel, 
\sl Groupes alg\'ebriques. Tome I: G\'eom\'etrie alg\'ebrique, g\'en\'eralit\'es, groupes commutatifs. Avec un appendice {\it Corps de classes local} par Michiel Hazewinkel,
\rm Masson \& Cie, \'Editeur, Paris; North-Holland Publishing Co., Amsterdam, 1970

\bibitem[Di]{Di} L.\ E.\ Dickson,
\sl The analytic representation of substitutions on a power of a prime number of letters with a discussion of the linear group,
\rm Ann.\ of Math.\ {\bf 11} (1896/97), no.\ 1--6, 65--120

\bibitem[DV]{DV} V.\ Drensky and J.-T.\ Yu,
\sl Automorphisms of polynomial algebras and Dirichlet series,
\rm J.\ Algebra {\bf 321} (2009), no.\ 1, 292--302

\bibitem[Dv]{Dv} Y.\ Dvir,
\sl Covering properties of permutation groups,
\rm Lecture Notes in Math.\ {\bf 1112}, Springer-Verlag, Berlin, 1985, 197--221

\bibitem[E]{E} J.\ Edmonds,
\sl Minimum partition of a matroid into independent subsets,
\rm J.\ Res.\ Nat.\ Bur.\ Standards Sect.\ B {\bf 69B} (1965), 67--72

\bibitem[FM]{FM} S.\ Friedland and J.\ Milnor,
\sl Dynamical properties of plane polynomial automorphisms,
\rm Ergodic Theory Dynam.\ Systems {\bf 9} (1989), no.\ 1, 67--99

\bibitem[F1]{F1} J.-P.\ Furter,
\sl On the variety of automorphisms of the affine plane,
\rm J.\ Algebra {\bf 195} (1997), no.\ 2, 604--623

\bibitem[F2]{F2} J.-P.\ Furter,
\sl On the length of polynomial automorphisms of the affine plane,
\rm Math.\ Ann.\ {\bf 322} (2002), no.\ 2, 401--411

\bibitem[F3]{F3} J.-P.\ Furter,
\sl On the degree of iterates of automorphisms of the affine plane,
\rm Manuscripta Math.\ {\bf 98} (1999), no.\ 2, 183--193

\bibitem[FP]{FP} J.-P.\ Furter and P.-M.\ Poloni,
\sl On the maximality of the triangular subgroup,
\rm Ann.\ Inst.\ Fourier (Grenoble) {\bf 68} (2018), no.\ 1, 393--421

\bibitem[G-POS]{G-POS} D.\ Gómez-Pérez, A.\ Ostafe, and I.\ Shparlinski,
\sl Algebraic entropy, automorphisms and sparsity of algebraic dynamical systems and pseudorandom number generators,
\rm Math.\ Comp.\ {\bf 83} (2014), no.\ 287, 1535--1550

\bibitem[Ha]{Ha} N.\ Harrison,
\sl Real length functions in groups,
\rm Trans.\ Amer.\ Math.\ Soc.\ {\bf 174} (1972), 77--106

\bibitem[HKL1]{HKL1} M.\ Herzog, G.\ Kaplan, and A.\ Lev,
\sl Representation of permutations as products of two cycles,
\rm Discrete Math.\ {\bf 285} (2004), no.\ 1-3, 323--327

\bibitem[HKL2]{HKL2} M.\ Herzog, G.\ Kaplan, and A.\ Lev,
\sl Covering the alternating groups by products of cycle classes,
\rm J.\ Combin.\ Theory Ser.\ A {\bf 115} (2008), no.\ 7, 1235--1245

\bibitem[HR1]{HR1} M.\ Herzog and K.\ B.\ Reid,
\sl Representation of permutations as products of cycles of fixed length,
\rm J.\ Austral.\ Math.\ Soc.\ Ser.\ A {\bf 22} (1976), no.\ 3, 321--331

\bibitem[HR2]{HR2} M.\ Herzog and K.\ B.\ Reid,
\sl Number of factors in $k$-cycle decompositions of permutations,
\rm Lecture Notes in Math., Vol.\ {\bf 560}, Springer-Verlag, Berlin-New York, 1976, pp.\ 123--131

\bibitem[Hw]{Hw} H.-K.\ Hwang,
\sl Distribution of the number of factors in random ordered factorizations of integers,
\rm J.\ Number Theory {\bf 81} (2000), no.\ 1, 61--92

\bibitem[Je]{Je} Z.\ Jelonek,
\sl The extension of regular and rational embeddings,
\rm Math.\ Ann.\ {\bf 277} (1987), no.\ 1, 113--120

\bibitem[Jo]{Jo} G.\ A.\ Jones,
\sl Primitive permutation groups containing a cycle,
\rm Bull.\ Aust.\ Math.\ Soc.\ {\bf 89} (2014), no.\ 1, 159--165

\bibitem[Kali]{Kali} S.\ Kaliman,
\sl Extensions of isomorphisms between affine algebraic subvarieties of $k^n$ to automorphisms of $k^n$, 
\rm Proc.\ Amer.\ Math.\ Soc.\ {\bf 113} (1991), no.\ 2, 325--334

\bibitem[KZ]{KZ} S.\ Kaliman and M.\ Zaidenberg,
\sl Affine modifications and affine hypersurfaces with a very transitive automorphism group,
\rm Transform.\ Groups {\bf 4} (1999), no.\ 1, 53--95

\bibitem[Kalm]{Kalm} L.\ Kalm\'ar,
\sl \"Uber die mittler Anzahl der Produktdarstellungen der Zahlen, erste Mitteilung,
\rm Acta Litterarum ac Scientiarum, Szeged {\bf 5} (1930--32), 95--107

\bibitem[Kan]{Kan} M.\ C.\ Kang,
\sl On Abhyankar--Moh's epimorphism theorem,
\rm Amer.\ J.\ Math.\ {\bf 113} (1991), no.\ 3, 399--421

\bibitem[Kur]{Kur} S.\ Kuroda,
\sl Elementary reducibility of automorphisms of a vector group,
\rm Saitama Math.\ J.\ {\bf 29} (2012), 79--87

\bibitem[LPS]{LPS} M.\ W.\ Liebeck, C.\ E.\ Praeger, and J.\ Saxl,
\sl A classification of the maximal subgroups of the finite alternating and symmetric groups,
\rm J.\ Algebra {\bf 111} (1987), no.\ 2, 365--383

\bibitem[L]{L} R.\ C.\ Lyndon,
\sl Length functions in groups,
\rm Math.\ Scand.\ {\bf 12} (1963), 209--234

\bibitem[M]{M} T.\ Moh,
\sl An application of algebraic geometry to encryption: tame transformation method,
\rm Rev.\ Mat.\ Iberoamericana {\bf 19} (2003), no.\ 2, 667--685

\bibitem[N]{N} M.\ Nagata,
\sl On automorphism group of $k[x,y]$,
\rm Department of Mathematics, Kyoto University, Lectures in Mathematics, No.\ {\bf 5}, Kinokuniya Book Store Co., Ltd., Tokyo, 1972

\bibitem[Poo]{Poo} B.\ Poonen,
\sl Bertini theorems over finite fields,
\rm Ann.\ of Math.\ (2) {\bf 160} (2004), no.\ 3, 1099--1127

\bibitem[Pro]{Pro} D.\ Promislow,
\sl Equivalence classes of length functions on groups,
\rm Proc.\ London Math.\ Soc.\ (3) {\bf 51} (1985), no.\ 3, 449--477

\bibitem[Ri]{Ri} D.\ R.\ Richman,
\sl On the computation of minimal polynomials,
\rm J.\ Algebra {\bf 103} (1986), no.\ 1, 1--17

\bibitem[Ro]{Ro} M.\ Rosenlicht,
\sl Questions of rationality for solvable groups over nonperfect fields,
\rm Ann.\ Mat.\ Pura Appl.\ (4) {\bf 61} (1963), 97--120

\bibitem[SGA3-1]{SGA3-1}
\sl Schémas en groupes (SGA 3). Tome I. Propriétés générales des schémas en groupes [Group schemes (SGA 3). Vol.\ I. General properties of group schemes],
\rm Doc.\ Math.\ (Paris), {\bf 7} [Mathematical Documents (Paris)], Société Mathématique de France, Paris, 2011

\bibitem[SU]{SU} I.\ P.\ Shestakov and U.\ Umirbaev,
\sl The tame and the wild automorphisms of polynomial rings in three variables,
\rm J.\ Amer.\ Math.\ Soc.\ {\bf 17} (2004), no.\ 1, 197--227

\bibitem[Sp]{Sp} T.\ A.\ Springer,
\sl Linear algebraic groups. Second edition,
\rm Progr.\ Math., {\bf 9}, Birkh\"auser Boston, Inc., Boston, MA, 1998

\bibitem[Sr]{Sr} V.\ Srinivas,
\sl On the embedding dimension of an affine variety,
\rm Math.\ Ann.\ {\bf 289} (1991), no.\ 1, 125--132

\bibitem[St]{St} R.\ Steinberg,
\sl Lectures on Chevalley groups,
\rm Univ.\ Lecture Ser., {\bf 66}, American Mathematical Society, Providence, RI, 2016

\bibitem[TK]{TK} T.\ Tanaka and H.\ Kaneta,
\sl Automorphism groups of vector groups over a field of positive characteristic,
\rm Hiroshima Math.\ J.\ {\bf 38} (2008), no.\ 3, 437--446

\bibitem[T]{T} J.\ Tits,
\sl Groupes de Whitehead de groupes algébriques simples sur un corps (d'après V.\ P.\ Platonov et al.),
\rm Lecture Notes in Math., {\bf 677}, Springer, Berlin, 1978, Exp.\ No.\ 505, pp.\ 218--236

\bibitem[UY]{UY} U.\ Umirbaev and J.-T.\ Yu,
\sl The strong Nagata conjecture,
\rm Proc.\ Natl.\ Acad.\ Sci.\ USA {\bf 101} (2004), no.\ 13, 4352--4355

\bibitem[vdK]{vdK} W.\ van der Kulk,
\sl On polynomial rings in two variables,
\rm Nieuw Arch.\ Wisk.\ (3) {\bf 1} (1953), 33--41

\bibitem[W]{W} D.\ Wright,
\sl Abelian subgroups of $\Aut_k(k[X,Y])$ and applications to actions on the affine plane,
\rm Illinois J.\ Math.\ {\bf 23} (1979), no.\ 4, 579--634
\end{thebibliography}
\end{document}